\documentclass[a4paper]{amsart}

\RequirePackage{amsmath} 
\RequirePackage{amssymb}
\usepackage{amscd,latexsym,amsthm,amsfonts,amssymb,amsmath,amsxtra}
\usepackage[colorlinks=true,urlcolor=blue,citecolor=blue]{hyperref}
\usepackage{color}
\usepackage[all]{xy}
\usepackage{bm}
\usepackage{mathtools}
\usepackage{ mathrsfs }
\usepackage{xcolor}
\usepackage{comment}
\usepackage{marginnote}
\usepackage{enumitem}
\usepackage{thmtools, thm-restate}

\let\Re\undefined
\let\Im\undefined

\DeclareMathOperator{\Re}{Re}
\DeclareMathOperator{\Im}{Im}

\DeclareMathOperator{\Tr}{Tr}

\DeclareMathOperator{\supp}{supp}

\DeclareMathOperator{\GL}{GL}

\newcommand{\floor}[1]{{\left\lfloor#1\right\rfloor}}

	\newcommand{\Res}{\operatorname{Res}}
	
	\newcommand{\Eis}{\operatorname{Eis}}

	\newcommand{\sgn}{\operatorname{sgn}}
	
	\newcommand{\Ad}{\operatorname{Ad}}

	\newcommand{\Reg}{\operatorname{Reg}}
	\newcommand{\fin}{\operatorname{fin}}
	\newcommand{\du}{\operatorname{du}}
	\newcommand{\tw}{\operatorname{tw}}
	\newcommand{\gen}{\operatorname{gen}}
	\newcommand{\const}{\operatorname{const}}
	\newcommand{\diag}{\operatorname{diag}}
	\newcommand{\cusp}{\operatorname{cusp}}

	\newcommand{\Vol}{\operatorname{Vol}}

	\newcommand{\Ind}{\operatorname{Ind}}
	
	\newcommand{\RNum}[1]{\uppercase\expandafter{\romannumeral #1\relax}}

\begin{document}
\theoremstyle{plain}
\newtheorem{thm}{Theorem}[section]
	
\newtheorem{cor}[thm]{Corollary}
\newtheorem{thmy}{Theorem}
\renewcommand{\thethmy}{\Alph{thmy}}
\newenvironment{thmx}{\stepcounter{thm}\begin{thmy}}{\end{thmy}}
\newtheorem{cory}{Corollary}
\renewcommand{\thecory}{\Alph{cory}}
\newenvironment{corx}{\stepcounter{thm}\begin{cory}}{\end{cory}}
\newtheorem{hy}[thm]{Hypothesis}
\newtheorem*{thma}{Theorem A}
\newtheorem*{corb}{Corollary B}
\newtheorem*{thmc}{Theorem C}
\newtheorem{lemma}[thm]{Lemma}  
\newtheorem{prop}[thm]{Proposition}
\newtheorem{conj}[thm]{Conjecture}
\newtheorem{fact}[thm]{Fact}
\newtheorem{claim}[thm]{Claim}
	
\theoremstyle{definition}
\newtheorem{defn}[thm]{Definition}
\newtheorem{example}[thm]{Example}
\theoremstyle{remark}
	
\newtheorem{remark}[thm]{Remark}	
\numberwithin{equation}{section}

\title[]{Rankin--Selberg Subconvexity via Spectral Reciprocity}% 
\author{Peter Humphries and Liyang Yang}
	
\address{Department of Mathematics, University of Virginia, Charlottesville, VA 22904, USA}
\email{pclhumphries@gmail.com}

\address{253-37 Caltech, Pasadena\\
CA 91125, USA}
\email{lyyang@caltech.edu}

\begin{abstract}
We establish explicit subconvex bounds for central values of Rankin--Selberg
$L$-functions $L(1/2,\pi\times\pi')$ associated with pairs of unitary cuspidal
automorphic representations of $\mathrm{GL}_2$ over a number field.
Building on the spectral reciprocity framework of Michel and Venkatesh,
we develop a refined, fully explicit form of spectral reciprocity that
allows for precise control of conductors and local test vectors. As a consequence, we obtain an explicit subconvex bound, which, even over
$F=\mathbb{Q}$, improves all previously known results.

We further apply these bounds to several arithmetic problems.
These include effective equidistribution of CM suborbits on quaternionic
Shimura varieties, quantitative equidistribution of totally geodesic
submanifolds, a uniform quantitative form of dihedral quantum unique
ergodicity over number fields, and an application to distinguishing
cuspidal automorphic representations.
\end{abstract}
	
\date{\today}
\maketitle
\tableofcontents

\section{Introduction}

Let $\pi$ and $\pi'$ be unitary automorphic representations of $\GL_2$ over a
number field $F$, and let $L(s,\pi\times\pi')$ denote the associated
Rankin--Selberg $L$-function.
The corresponding \emph{subconvexity problem} \textbf{(ScP)} asks whether there
exists $\delta>0$ such that 
\begin{equation}\label{1.1}
L(1/2,\pi\times\pi')
\ll C(\pi\times\pi')^{\frac{1}{4}-\delta},
\end{equation}
where $C(\pi\times\pi')$ is the analytic conductor.

The problem \textbf{(ScP)} is one of the central and challenging questions
in analytic number theory.
When $F=\mathbb{Q}$ and $\pi'$ is fixed and cuspidal, estimates of the form
\eqref{1.1} with explicit values of $\delta$ have been established in a variety
of aspects, leading to numerous applications in arithmetic and geometry; see,
for instance,
\cite{Sar01, KMV02, LY02, Mic04, HM06, Har11, AKKLLY24}.

For general number fields, a fundamental breakthrough was achieved by Michel and Venkatesh \cite[Theorem~1.2]{MV10}. They proved that for unitary automorphic representations $\pi$ and $\pi'$ of $\GL_2/F$, there exist positive constants $A$ and $\delta$ such that
\begin{equation}\label{e1.2}
L(1/2,\pi\times\pi')
\ll_F C(\pi')^{A} C(\pi)^{\frac{1}{2}-\delta}.
\end{equation}
Although the exponents $A$ and $\delta$ in \eqref{e1.2} are \textit{implicit}, this result has already had far-reaching consequences, including applications to equidistribution and arithmetic problems in great generality; see, for example,
\cite{Zha05, MV06, Mic07, MW09}.

For $F\neq \mathbb{Q}$ and general cuspidal $\pi$ and $\pi'$, no explicit value of $\delta$ has previously been established. 
Nevertheless, both from the perspective of methodology and in view of refined arithmetic applications, it is highly desirable to obtain \emph{explicit} exponents $A$ and $\delta$ in \eqref{e1.2}. 
The main goal of this paper is to establish a refined spectral reciprocity formula, building on the framework of \cite{MV10}, and to exploit it to derive explicit subconvex bounds with concrete values of $(A,\delta)$, improving upon all previously known results in the case $F=\mathbb{Q}$.

\subsection{Explicit Subconvexity for $\GL_2\times\GL_2$}

We now state our main result, which provides the explicit exponents $(A,\delta)=(5/2,\,1/64)$ in \eqref{e1.2}, thereby refining the theorem of Michel--Venkatesh \cite{MV10}.

\begin{thmx}\label{A}
Let $F$ be a number field and let $0<\varepsilon<10^{-3}$. Let $\pi$ and $\pi'$ be unitary cuspidal automorphic representations of 
$\mathrm{GL}_2/F$, with analytic conductors $C(\pi)$ and $C(\pi')$, respectively. Then
\begin{equation}\label{e1.3}
L(1/2,\pi\times\pi')\ll C(\pi')^{\frac{5}{2}+\varepsilon}C(\pi)^{\frac{1}{2}-\frac{1}{64}+\varepsilon}, 
\end{equation}
where the implied constant depends only on $F$ and $\varepsilon$.
\end{thmx}

Theorem \ref{A} yields a substantial explicit improvement over previously known
quantitative bounds.
In particular, the strongest earlier result in the level aspect attained the
exponent $\delta = 1/524$ for prime levels \cite[Theorem 1.1]{AKKLLY24}.
A non-cuspidal analogue of Theorem \ref{A} is established in \cite{Yan26b}.
The strength and uniformity of the exponent in \eqref{e1.3}, together with its
non-cuspidal counterpart, make these bounds particularly well suited for
arithmetic applications requiring effective control of central values of
Rankin--Selberg $L$-functions.

\subsection{Arithmetic Applications}
In this subsection, we present several arithmetic applications of Theorem \ref{A} and its non-cuspidal analogue. 
These include effective equidistribution results for CM suborbits and totally geodesic submanifolds, quantitative forms of dihedral quantum unique ergodicity, and an application to distinguishing cuspidal automorphic representations.

\subsubsection{Equidistribution of CM-suborbits}\label{sec1.1.1}
Let $F$ be a totally real field and $B$ be a quaternion algebra over $F$. Fix a maximal order $\mathcal{O}_B$ of $B$. Let $G=\mathrm{Res}_{F/\mathbb{Q}}(B^{\times}/F^{\times})$. Then there exists $0\leq d\leq [F:\mathbb{Q}]$ such that 
\begin{align*}
G(\mathbb{R})=\prod_{v\mid\infty}(B\otimes F_v)^{\times}/F_v^{\times}\simeq \mathrm{PGL}_2(\mathbb{R})^d\times \mathrm{SO}_3^{[F:\mathbb{Q}]-d}.
\end{align*}
Via such an isomorphism $G(\mathbb{R})$ acts on $\mathbb{X}:=(\mathbb{C}\setminus \mathbb{R})^d$. 

Let $U$ be an open compact subgroup of $G(\mathbb{A}_{\mathbb{Q},\fin})$. 
Define the analytic variety 
\begin{align*}
M_U(\mathbb{C}):=G(\mathbb{Q})\backslash\mathbb{X}\times G(\mathbb{A}_{\mathbb{Q},\fin})/U.
\end{align*} 
A point $\boldsymbol{z}=[(x,g)]\in M_U(\mathbb{C})$ is called a \textit{CM-point} if and only if 
it is represented by $(x,g)\in \mathbb{X}\times G(\mathbb{A}_{\mathbb{Q},\fin})$ 
such that the stabilizer of $\boldsymbol{z}$ in $G(\mathbb{Q})$ is a torus 
$T:= K^{\times}/F^{\times}$, 
where $E/F$ is a quadratic CM extension embedded into $B$. Let 
\begin{align*}
\mathcal{O}_{\boldsymbol{z}}:=E\cap g\mathcal{O}_Bg^{-1}
\end{align*}
be the associated order. Denote by $\mathrm{disc}(\boldsymbol{z}):=\mathrm{disc}(\mathcal{O}_{\boldsymbol{z}})$ its discriminant. 

For a CM-point $\boldsymbol{z}=[(x,g)]$ and an open subgroup $[T']$ of 
$T(\mathbb{Q})\backslash T(\mathbb{A}_{\mathbb{Q},\fin})$, define the corresponding 
CM-suborbit by 
\begin{align*}
\mathrm{Orb}_{T'}(\boldsymbol{z}):=\big\{[(x,t'g)] :\ t'\in [T']\big\}.
\end{align*}

We will simply write $\mathrm{Orb}(\boldsymbol{z})$ for $\mathrm{Orb}_{T'}(\boldsymbol{z})$ when the 
chosen subgroup $[T']$ is implicit.  
Note that $\mathrm{Orb}(\boldsymbol{z})$ is finite; let $\#\mathrm{Orb}(\boldsymbol{z})$ denote its cardinality. 

Combining Theorem \ref{A} and its non-cuspidal analogue \cite{Yan26b} with 
\cite[Corollary 3.3]{Zha05}, we obtain the following result. 
\begin{thmx}\label{B}
Let $0<\varepsilon<10^{-3}$.  Let $M$ be a fixed connected component of $M_U$, equipped with the
canonical probability measure $\mu$, and let $f \in C^\infty_c(M)$.  
Then for every CM-point $\boldsymbol{z}\in M$ with CM-suborbit $\mathrm{Orb}(\boldsymbol{z})$, one has
\begin{equation}\label{eq1.1}
\bigg|
\frac{1}{\# \mathrm{Orb}(\boldsymbol{z})}
\sum_{y\in \mathrm{Orb}(\boldsymbol{z})} f(y) -\int_M f  d\mu
\bigg|\ll
\frac{\mathrm{disc}(\boldsymbol{z})^{\frac{1}{2}- \frac{1}{128} + \varepsilon}}
     {\#\mathrm{Orb}(\boldsymbol{z})},
\end{equation}
where the implied constant depends only on $f$, $F$, and $\varepsilon$. In particular, as 
$\mathrm{disc}(\boldsymbol{z})\to\infty$, the CM-suborbit $\mathrm{Orb}(\boldsymbol{z})$ becomes 
equidistributed in $M$ provided that $\#\mathrm{Orb}(\boldsymbol{z})\gg \mathrm{disc}(\boldsymbol{z})^{\frac{1}{2}- \frac{1}{128} + \varepsilon}$. 
\end{thmx}

%\begin{remark}
%Since the subconvexity Theorem \ref{A} is uniform, the above estimate \eqref{eq1.1} holds uniformly for $U$. 
%\end{remark}

\begin{remark}
When $F=\mathbb{Q}$ and $U=\widehat{\mathcal{O}_B^{\times}}$ is \textit{fixed}, Theorem \ref{B} simultaneously improves the quantitative sparse 
equidistribution results of Michel \cite{Mic04} and Harcos--Michel \cite{HM06}: 
\begin{itemize}
\item if $B/\mathbb{Q}$ is definite (so $d=0$), the estimate \eqref{eq1.1} 
improves the exponent in Michel \cite[Theorem 10]{Mic04} from $1/2115$ to 
$1/128-\varepsilon$;

\item if $B/\mathbb{Q}$ is indefinite (so $d=1$, and the Galois orbit of $\boldsymbol{z}$ amounts to the CM-orbit $\mathrm{Orb}_T(\boldsymbol{z})$), the estimate \eqref{eq1.1} 
improves the exponent in Harcos--Michel 
\cite[Theorem 6]{HM06} from $1/5298$ to $1/128-\varepsilon$.
\end{itemize}
\end{remark}

According to \cite[Corollary 3.8]{Zha05}, the above comequidistribution holds for Galois CM-points with a fixed maximal
Mumford-Tate group $H=T$. In particular, any infinite set of such CM-points are Zariski
dense.

\subsubsection{Equidistribution of Totally Geodesic Submanifolds}
Let $F$ be a number field and $E/F$ be quadratic. Let $\mathcal{O}$ be an order in $E$, and write $K = K_{\infty}K_{\fin}$ with 
$K_{\infty}\subset \mathrm{GL}_2(F_{\infty})$ a maximal compact subgroup and 
$K_{\fin}\subset \mathrm{GL}_2(\mathbb{A}_{F,\fin})$ an open compact subgroup.
 Set $X=Z(\mathbb{A}_F)\mathrm{GL}_2(F)\backslash\mathrm{GL}_2(\mathbb{A}_F)/K$, equipped with the
canonical probability measure $\mu$.

Let $T$ be a standard torus obtained by an optimal embedding of $E^{\times}\hookrightarrow \mathrm{GL}_2/F$ relative to $K_{\fin}$.  
Then there exists $g_{\infty}\in \mathrm{GL}_2(F_{\infty})\hookrightarrow \mathrm{GL}_2(\mathbb{A}_F)$ such that
\begin{align*}
K_{T,\infty}:=g_{\infty}K_{\infty}g_{\infty}^{-1}\cap T(F_{\infty})
\end{align*}
is the maximal compact subgroup of $T(F_{\infty})$. Then
\begin{align*}
\bigl(T(F)\backslash T(\mathbb{A}_F)/(K_{T,\infty}\widehat{\mathcal{O}}^{\times})\bigr) g_{\infty}
=
\bigcup_{\mathfrak{a}\in \mathrm{Pic}(\mathcal{O})}\gamma_{\mathfrak{a}},
\end{align*}
where $\gamma_{\mathfrak{a}}$ denotes the fiber above $\mathfrak{a}$ under the projection
\begin{align*}
T(F)\backslash T(\mathbb{A}_F)/(K_{T,\infty}\widehat{\mathcal{O}}^{\times})
\longrightarrow
T(F)\backslash T(\mathbb{A}_{F,\fin})/\widehat{\mathcal{O}}^{\times}\simeq \mathrm{Pic}(\mathcal{O}).
\end{align*}

By combining Theorem \ref{A} and its non-cuspidal analogue \cite{Yan26b} with 
\cite[\textsection 4--\textsection 5]{MW09}, together with the argument underlying the proof of \cite[Theorem 6]{HM06},
we obtain the following quantitative equidistribution result.
\begin{thmx}
Let $0<\varepsilon<10^{-3}$, and let $G\subseteq \mathrm{Pic}(\mathcal{O})$ be a subgroup.
Let $f \in C^\infty_c(X)$. Then 
\begin{equation}\label{1.2}
\bigg|\frac{1}{\#G}\sum_{\mathfrak{a}\in G}\frac{1}{\Vol(\gamma_{\mathfrak{a}})}\int_{\gamma_{\mathfrak{a}}}fd\mu-\int_{X}fd\mu\bigg|\ll [\mathrm{Pic}(\mathcal{O}): G]\cdot \mathrm{disc}(\mathcal{O})^{-\frac{1}{128}+\varepsilon},  
\end{equation}	
where the implied constant depends only on $f$, $F$ and $\varepsilon$.  
\end{thmx}

By taking $G$ to be trivial in \eqref{1.2}, we obtain the equidistribution of each individual $\gamma_{\mathfrak{a}}$ in $X$: 
\begin{cor}
Suppose that $\mathrm{disc}(\mathcal{O}) \to \infty$ and that $\#\mathrm{Pic}(\mathcal{O})\ll \mathrm{disc}(\mathcal{O})^{\frac{1}{128}-\varepsilon}$. Then, for every $\mathfrak{a}\in \mathrm{Pic}(\mathcal{O})$,  $\gamma_{\mathfrak{a}}$ becomes equidistributed in $X$ as $\mathrm{disc}(\mathcal{O})\to\infty$. 
\end{cor}
\begin{remark}
If $T(F_v)$ is noncompact for every Archimedean place $v\mid\infty$, then each
$\gamma_{\mathfrak{a}}$ is a totally geodesic submanifold of $X$. More generally, when
$T(F_v)$ is noncompact for some Archimedean places $v\mid\infty$ and compact for the remaining ones,
$\gamma_{\mathfrak{a}}$ is still a totally geodesic submanifold of $X$, of lower dimension,
given by a product of geodesic directions at the noncompact places and fixed points at the compact places.
In particular, when $F=\mathbb{Q}$ and $E/F$ is real quadratic, \eqref{1.2} improves the results of
Duke \cite{Duk88} and Popa \cite{Pop06}, and when $F=\mathbb{Q}(i)$,
\eqref{1.2} improves the results of Martin--Whitehouse \cite{MW09}.
\end{remark}

\subsubsection{Quantitative Dihedral Quantum Unique Ergodicity}
Following Sarnak's original observation \cite{Sar01}, extended to the Maass form 
setting by \cite{LY02}, our Theorem \ref{A}, together with its non-cuspidal analogue \cite{Yan26b}, yields a quantitative form of dihedral Quantum Unique Ergodicity (QUE) in the sense of \cite{RS94}.

In addition, Kowalski--Michel--VanderKam \cite[Conjecture 1.5]{KMV02} proposed a 
QUE conjecture in the level aspect. For holomorphic forms of fixed weight $k$ 
and varying squarefree level $q$, this was established by Nelson 
\cite[Theorem 1.3]{Nel11}. Nelson--Pitale--Saha \cite{NPS14} subsequently extended 
this to the uniform setting as $qk \to \infty$ for general levels $q$, obtaining 
in particular power savings in the depth aspect. 

In this subsection, we establish a quantitative, uniform QUE result over number 
fields in the dihedral setting, refining the results mentioned above.

Let $F$ be a number field. Set $X := Z(\mathbb{A}_F)\mathrm{GL}_2(F)\backslash \mathrm{GL}_2(\mathbb{A}_F)/K$, where $K$ is a maximal compact subgroup of $\mathrm{GL}_2(\mathbb{A}_F)$. We then obtain the following quantitative QUE result. 

\begin{thmx}\label{thmD}
Suppose that $\pi$ is a dihedral representation of $\mathrm{GL}_2/F$. 
Let $\phi \in \pi$ be a \emph{nice} cusp form (see Definition \ref{defn1.7} in 
\textsection\ref{sec15.2}). 
Then for every $f \in C_c(X)$,
\begin{equation}\label{1.3}
\bigg|
\int_{Z(\mathbb{A}_F)G(F)\backslash G(\mathbb{A}_F)}
    f(x)|\phi(x)|^2dx
-
\int_X f(x)dx
\bigg|
\ll
C(\pi)^{-\frac{1}{128}+\varepsilon}
C_{\fin}(\pi)^{-\frac{1}{16}},
\end{equation}
where the implied constant depends only on $f$, $F$, and $\varepsilon$.
\end{thmx}

\begin{remark}
When $F=\mathbb{Q}$, \eqref{1.3} yields a direct strengthening of the results of
\cite{Sar01} and \cite{LY02} in the spectral parameter aspect, and provides an improved error term--namely, a
\textit{uniform power saving} in the analytic conductor--compared with
\cite[Theorem 1.2]{NPS14}, in the dihedral case.
\end{remark}

\begin{remark}
Let $E/F$ be a fixed quadratic extension and $\mathcal{O}$ be an order in $E$. Let $\chi$ be a character of $\mathrm{Pic}(\mathcal{O})$, and let $\pi=\mathrm{AI}_{E/F}(\chi)$ be the automorphic induction of $\chi$. Then $C(\pi\times\widetilde{\pi})\asymp C(\pi)\asymp \mathrm{disc}(\mathcal{O})$, where the implied constant depends on $E/F$. Similar to \cite{NPS14}, we obtain directly from an estimate of the corresponding trilinear forms (see Proposition \ref{prop15.2}) and the convexity bound that
\begin{align*}
\bigg|
\int_{Z(\mathbb{A}_F)G(F)\backslash G(\mathbb{A}_F)}
    f(x)|\phi(x)|^2dx-\int_X f(x)dx\bigg|
\ll
C(\pi)^{-\frac{1}{4}+\varepsilon}.
\end{align*}

However, for a general dihedral $\pi$, the power saving in \eqref{1.3} is genuinely
nontrivial and relies essentially on our explicit subconvexity bounds, together
with the Burgess bound established in \cite{Yan26}. 
\end{remark}

\subsubsection{Distinguishing Cuspidal Representations}
Following the arguments in the proof of \cite[Corollary 1.3]{KMV02}, we obtain the following generalization.  
\begin{thmx}
Let $F$ be a number field and let $0<\varepsilon<10^{-3}$. 
Let $\pi$ and $\pi'$ be unitary cuspidal automorphic representations of 
$\mathrm{GL}_2/F$. Then there exists an integral ideal $\mathfrak{n}\subseteq\mathcal{O}_F$ with
\begin{equation}\label{1.5}
N_F(\mathfrak{n})\ll \min\Big\{C(\pi')^{5+\varepsilon}C(\pi)^{1-\frac{1}{32}+\varepsilon},C(\pi)^{5+\varepsilon}C(\pi')^{1-\frac{1}{32}+\varepsilon}\Big\},
\end{equation}
such that $\lambda_{\pi}(\mathfrak{n})\neq \lambda_{\pi'}(\mathfrak{n})$, where $\lambda_{\pi}(\mathfrak{n})$ and $\lambda_{\pi'}(\mathfrak{n})$ are the normalized Hecke eigenvalues. Here, the implied constant in \eqref{1.5} depends only on $F$ and $\varepsilon$. 
\end{thmx}
\begin{remark}
When $F=\mathbb{Q}$, $\pi'$ is fixed, and the central character is \textit{trivial}, the bound \eqref{1.5} improves the exponent in \cite[Corollary~1.3]{KMV02} from $1/40$ to $1/32$.
\end{remark}

\subsection{Discussion of the Proof}
\subsubsection{The Spectral Reciprocity}
Let $\mathfrak{L}$, $\mathfrak{M}$, and $\mathfrak{N}$ be integral ideals, and let
$\omega$ and $\omega'$ be unitary Hecke characters.
Set $\mathbf{C}_{\infty} := \prod_{v\mid\infty} \mathbf{C}_v$, where
$\mathbf{C}_v \geq C_v(\pi')$ for each Archimedean place $v$.
Let $\mathcal{F}(\mathfrak{N},\omega)$ denote the set of automorphic
representations of $\GL_2/F$ with conductor ideal $\mathfrak{N}$ and central
character $\omega$.

Throughout this paper, the principal analytic inputs are the establishment of a
spectral reciprocity formula (see Theorem \ref{thmA} in \textsection\ref{sec3.5}) and the
local and global estimates of the relevant period integrals for suitably chosen
test vectors (see \textsection\ref{sec8}--\textsection\ref{sect13}).
At a heuristic level, the resulting spectral reciprocity takes the following
form:
\begin{multline}\label{1.8}
\int_{\substack{\pi\in \mathcal{F}(\mathfrak{N},\omega)\\ C_v(\pi)\ll \mathbf{C}_v,\ v\mid\infty}}\lambda_{\pi}(\mathfrak{L})|L(1/2,\pi\times\pi')|^2d\mu_{\pi}\rightsquigarrow C_{\infty}(\pi')^{9/4+\varepsilon}C_{\infty}(\omega\overline{\omega}')^{-\vartheta}\\
\mathbf{C}_{\infty}^{1/2+\vartheta+\varepsilon}[N_F(\mathfrak{N}),N_F(\mathfrak{M})]^{1/2+\varepsilon}N_F(\mathfrak{L})^{-1/2+\varepsilon}\sum_{\mathfrak{a}\supseteq \mathfrak{L}}|\lambda_{\pi'}(\mathfrak{a})|\int_{\substack{\sigma\in \mathcal{F}(\mathfrak{L}\mathfrak{M};\mathbf{1})\\ C_v(\sigma)\ll C_v(\pi')^{2+\varepsilon},\ v\mid\infty}}\\
L(1/2,\Ad\pi'\times\sigma)^{1/2}|L(1/2,\sigma)|^{3/2}|L(1/2,\sigma \times  \omega\overline{\omega}')|d\mu_{\sigma},
\end{multline}
where $\vartheta$ denotes any admissible exponent toward the Ramanujan
conjecture. 

In the special case where $\omega=\omega'=\mathbf{1}$, the level $\mathfrak{N}$
is squarefree, $\mathfrak{M}=\mathcal{O}_F$, and $\pi'$ is tempered,
a reciprocity formula of the form \eqref{1.8} was established by
Zacharias \cite{Zac20}, following the approach of Michel and Venkatesh
\cite{MV10}.

The most delicate situation arises when the central characters and the level
$\mathfrak{N}$ are allowed to be general, and when $\pi'$ is arbitrary.
To treat this case, we introduce an analytic regularization that differs from
that of \cite{MV10}, together with a careful local analysis of the trilinear forms on the dual side, which constitutes the bulk of this paper.
These ingredients allow us to establish the reciprocity formula
\eqref{1.8} in full generality.

\subsubsection{Moment Estimates and Hybrid Subconvexity}
By H\"{o}lder's inequality, the integral
\begin{align*}
\int_{\substack{\sigma\in \mathcal{F}(\mathfrak{L}\mathfrak{M};\mathbf{1})\\ C_v(\sigma)\ll C_v(\pi')^{2+\varepsilon},\ v\mid\infty}}
L(1/2,\Ad\pi'\times\sigma)^{1/2}|L(1/2,\sigma)|^{3/2}|L(1/2,\sigma \times  \omega\overline{\omega}')|d\mu_{\sigma}
\end{align*}
is bounded by
\begin{multline*}
\max_{\substack{\sigma\in \mathcal{F}(\mathfrak{L}\mathfrak{M};\mathbf{1})\\ C_v(\sigma)\ll C_v(\pi')^{2+\varepsilon},\ v\mid\infty}}|L(1/2,\sigma \times  \omega\overline{\omega}')|
\bigg[\int_{\substack{\sigma\in \mathcal{F}(\mathfrak{L}\mathfrak{M};\mathbf{1})\\ C_v(\sigma)\ll C_v(\pi')^{2+\varepsilon},\ v\mid\infty}}
|L(1/2,\sigma)|^{4}d\mu_{\sigma}\bigg]^{\frac{3}{8}}\\
\bigg[\int_{\substack{\sigma\in \mathcal{F}(\mathfrak{L}\mathfrak{M};\mathbf{1})\\ C_v(\sigma)\ll C_v(\pi')^{2+\varepsilon},\ v\mid\infty}}
L(1/2,\Ad\pi'\times\sigma)d\mu_{\sigma}\bigg]^{\frac{1}{2}}\bigg[\int_{\substack{\sigma\in \mathcal{F}(\mathfrak{L}\mathfrak{M};\mathbf{1})\\ C_v(\sigma)\ll C_v(\pi')^{2+\varepsilon},\ v\mid\infty}}
1d\mu_{\sigma}\bigg]^{\frac{1}{8}}.
\end{multline*}

The individual $L$-values and moments appearing above are bounded as follows:
\begin{itemize}
\item
The first moment of $L(1/2,\Ad\pi'\times\sigma)$ is bounded using
\cite[Theorem C]{Yan25c}; see Theorem \ref{thm14.2} in
\textsection\ref{sec14.2.1}.

\item
The fourth moment of $L(1/2,\sigma)$ is bounded by the result established
in \cite{Yan26b}; see Theorem \ref{thm14.4} in
\textsection\ref{sec14.2.3}.

\item
The factor $\bigl|L(1/2,\sigma \times \omega\overline{\omega}')\bigr|$ is
bounded uniformly in terms of $C(\pi')$ and $N_F(\mathfrak{L}\mathfrak{M})$
by the hybrid bound of \cite[Corollary 1.10]{Yan26}, which yields a
Burgess-type bound in the $\omega\overline{\omega}'$-aspect; see
Theorem \ref{thm14.3} in \textsection\ref{14.2.4}. 
\end{itemize}

\subsubsection{Removing the Temperedness Assumption}  
In order to control the factor $|\lambda_{\pi'}(\mathfrak{a})|$ in
\eqref{1.8}, the temperedness of
$\pi'$ is assumed in \cite{Zac20}.
We show that this assumption can be removed by combining the
Rankin--Selberg method with bounds arising from symmetric power lifts of
$\GL_2$, building on an observation of \cite{Yan21}.

\subsubsection{Amplification}
Combining the estimates obtained above with \eqref{1.8}, and incorporating
amplification in each term of the regularized spectral reciprocity,
we deduce that for $L\gg 1$,
\begin{multline*}
|L(1/2,\overline{\pi}\times\pi')|^2\ll C_{\infty}(\pi')^{\frac{19}{4}+\varepsilon}\mathbf{C}_{\infty}^{1+\varepsilon}N_F(\mathfrak{M})^{\frac{3}{2}+\varepsilon}[N_F(\mathfrak{N}),N_F(\mathfrak{M})]^{1+\varepsilon}L^{-1+\varepsilon}\\
+\mathbf{C}_{\infty}^{\frac{1}{2}+\vartheta+\varepsilon}C_{\infty}(\omega\overline{\omega}')^{-\vartheta}
[N_F(\mathfrak{N}),N_F(\mathfrak{M})]^{\frac{1}{2}+\varepsilon}N_F(\mathfrak{M})C_{\infty}(\pi')^{\frac{43}{8}+\varepsilon}C_{\fin}(\Ad\pi')^{\frac{1}{2}}\\
\Big[L^{4+\varepsilon}N_F(\mathfrak{M})^{\frac{1}{2}}C(\omega\overline{\omega}')^{\frac{1}{4}}C_{\infty}(\pi')^{\frac{1}{8}}
+L^{3+\varepsilon}N_F(\mathfrak{M})^{\frac{1}{4}}C(\omega\overline{\omega}')^{\frac{3}{8}}\Big].
\end{multline*}

Optimizing the parameter $L$ in the above inequality yields
Theorem \ref{A}.

\subsection{Outline of the Paper}

The structure of this paper is as follows.
\begin{itemize}
\item
In \textsection\ref{sec2}, we introduce notation and collect several
preliminary results.

\item
In \textsection\ref{sec3}, we establish a regularized spectral reciprocity
formula for general test vectors.

\item
In \textsection\ref{sec4}, we specify the choice of test vectors and derive an
explicit form of the spectral reciprocity.

\item
In \textsection\ref{sec5}, we compute the twisted moment and derive lower
bounds for it; along the way, we refine the Archimedean newvector theory
developed in \cite{MV10}.

\item
In \textsection\ref{sec6}, we establish a sharp upper bound for the amplified
constant term on the dual side.

\item
In \textsection\ref{sec8}--\textsection\ref{sect13}, we compute and estimate
the local period integrals appearing on the dual side.

\item
In \textsection\ref{sec14}, we assemble the local analysis of the period
integrals to obtain global bounds for the dual side in terms of moments of
$L$-functions, and we also control the contribution of the amplified residual
terms.

\item
Finally, in \textsection\ref{sec15}, we combine the preceding global estimates
and optimize the amplification parameter to prove Theorem~\ref{A}.
We also incorporate our local analysis of the period integrals into the
framework of Theorem~\ref{thmD} and deduce Theorem~\ref{thmD} from
Theorem~\ref{A} and its non-cuspidal analogue.
\end{itemize}

\textbf{Acknowledgements}
We are especially grateful to Gergely Harcos for helpful discussions and valuable suggestions. We would also like to thank Rizwanur Khan, Riad Masri, Philippe Michel, and Matthew Young for their helpful comments and suggestions.

\section{Some Preliminaries}\label{sec2}
\subsection{Notation}
\subsubsection{Fields and measures}\label{2.1.1}
Let $F$ be a number field with ring of integers $\mathcal{O}_F$. Let $[F:\mathbb{Q}]$ be the degree. Let $N_F$ be the absolute norm. Let $\mathfrak{O}_F$ be the different of $F$. Let $\mathbb{A}_F$ be the adele group of $F$. Let $\Sigma_F$ be the set of places of $F$. Denote by $\Sigma_{F,\fin}$ (resp. $\Sigma_{F,\infty}$) the set of nonarchimedean (resp. archimedean) places. For $v\in \Sigma_F$, we denote by $F_v$ the corresponding local field and $\mathcal{O}_{v}$ its ring of integers. For a nonarchimedean place $v$, let  $\mathfrak{p}_v$ be the maximal prime ideal in $\mathcal{O}_{v}$. Given an integral ideal $\mathcal{I}$, we say $v\mid \mathcal{I}$ if $\mathcal{I}\subseteq \mathfrak{p}_v$. Fix a uniformizer $\varpi_{v}\in\mathfrak{p}_v$. Denote by $e_v(\cdot)$ the evaluation relative to $\varpi_v$ normalized as $e_v(\varpi_v)=1$. Let $q_v$ be the cardinality of the residue field $\mathcal{O}_v/\mathfrak{p}_v$. We use $v\mid\infty$ to indicate an archimedean place $v$ and write $v<\infty$ if $v$ is nonarchimedean. Let $|\cdot|_v$ be the norm in $F_v$. Put $|\cdot|_{\infty}=\prod_{v\mid\infty}|\cdot|_v$ and $|\cdot|_{\fin}=\prod_{v<\infty}|\cdot|_v$. Let $|\cdot|_{\mathbb{A}_F}=|\cdot|_{\infty}\otimes|\cdot|_{\fin}$. We will simply write $|\cdot|$ for $|\cdot|_{\mathbb{A}_F}$ in calculation over $\mathbb{A}_F^{\times}$ or its quotient by $F^{\times}$.   

Denote by $\Tr_F$ the trace map, extended to $\mathbb{A}_F\rightarrow \mathbb{A}_{\mathbb{Q}}$. Let $\psi_{\mathbb{Q}}$ be the additive character on $\mathbb{Q}\backslash \mathbb{A}_{\mathbb{Q}}$ such that $\psi_{\mathbb{Q}}(t_{\infty})=\exp(2\pi it_{\infty})$, for $t_{\infty}\in \mathbb{R}\hookrightarrow\mathbb{A}_{\mathbb{Q}}$. Let $\psi=\psi_{\mathbb{Q}}\circ \Tr_F$. Then $\psi(t)=\prod_{v\in\Sigma_F}\psi_v(t_v)$ for $t=(t_v)_v\in\mathbb{A}_F$. For $v\in \Sigma_F$, let $dt_v$ be the additive Haar measure on $F_v$, self-dual relative to $\psi_v$. Then $dt=\prod_{v\in\Sigma_F}dt_v$ is the standard Tamagawa measure on $\mathbb{A}_F$. Let $d^{\times}t_v=\zeta_v(1)dt_v/|t_v|_v$, where $\zeta_v(\cdot)$ is the local Dedekind zeta factor. In particular, $\Vol(\mathcal{O}_{v}^{\times},d^{\times}t_v)=\Vol(\mathcal{O}_{v},dt_v)=|\mathfrak{O}_{v}|_v^{-1/2}$ for all finite place $v$, where $|\mathfrak{O}_{v}|_v$ denotes the local norm of the local different $\mathfrak{O}_v$. Moreover, $\Vol(F\backslash\mathbb{A}_F; dt)=1$ and $\Vol(F\backslash\mathbb{A}_F^{(1)},d^{\times}t)=\underset{s=1}{\Res}\ \zeta_F(s)$, where $\mathbb{A}_F^{(1)}$ is the subgroup of ideles $\mathbb{A}_F^{\times}$ with norm $1$, and $\zeta_F(s)=\prod_{v<\infty}\zeta_v(s)$ is the finite Dedekind zeta function.

\subsubsection{Algebraic groups and measures}
For an algebraic group $H$ over $F$, we will denote by $[H]:=H(F)\backslash H(\mathbb{A}_F)$. We equip measures on $H(\mathbb{A}_F)$ as follows: for each unipotent group $U$ of $H$, we equip $U(\mathbb{A}_F)$ with the Haar measure such that, $U(F)$ being equipped with the counting measure and the measure of $[U]$ is $1$. We equip the maximal compact subgroup $K$ of $H(\mathbb{A}_F)$ with the Haar measure such that $K$ has total mass $1$. When $H$ is split, we also equip the maximal split torus of $H$ with Tamagawa measure induced from that of $\mathbb{A}_F^{\times}$. 

Let $G=\mathrm{GL}(2)$. Let $B$ be the Borel subgroup of $G$ consisting of upper triangle matrices. Let $T$ (resp. $N$) be the Levi component (resp. unipotent radical) of $B$. Let $\overline{G}=Z\backslash G=\mathrm{PGL}(2)$, where $Z$ is the center of $G$. Let $B_0$ be the mirabolic subgroup of $G$, which is isomorphic to $Z\backslash B$. Let $A$ be the Levi component of $B_0$. 

\subsubsection{Level structures}\label{sec2.1.3}
Let $K_{\infty}=\otimes_{v\mid\infty}K_v$, where $K_v=\mathrm{O}(2)$ if $F_v\simeq \mathbb{R}$, and $K_v=\mathrm{U}(2)$ if $F_v\simeq\mathbb{C}$. Let $K_{\fin}=\otimes_{v\in \Sigma_{\fin}}K_v$, where   $K_v=G(\mathcal{O}_v)$. Define $K:=K_{\infty}\otimes K_{\fin}$, which is  a maximal compact subgroup of $G(\mathbb{A}_F)$. 

Let $v\in \Sigma_{\fin}$ and  $n\geq 0$ be an integer. Define
\begin{align*}
&K_v[n]:=\Big\{\begin{pmatrix}
k_{11}& k_{12}\\
k_{21}& k_{22}	
\end{pmatrix}\in K_v:\ k_{21}\in \varpi_v^n\mathcal{O}_v
\Big\}.
\end{align*} 
In particular, $K_v[0]=K_v$. Let $\mathfrak{L}\subset\mathcal{O}_F$ be an integral ideal. Define $K_{\fin}[\mathfrak{L}]:=\otimes_{v\in \Sigma_{\fin}}K_v[e_v(\mathfrak{L})]$,  where $e_v(\cdot)$ is the standard evaluation of $F_v$.     

\subsubsection{Additive Characters}\label{sec2.1.4}
Let $\psi_{\mathbb{Q}}$ be the additive character on $\mathbb{Q}\backslash \mathbb{A}_{\mathbb{Q}}$ such that $\psi_{\mathbb{Q}}(t_{\infty})=\exp(2\pi it_{\infty})$, for $t_{\infty}\in \mathbb{R}\hookrightarrow\mathbb{A}_{\mathbb{Q}}$. Let $\psi=\psi_{\mathbb{Q}}\circ \Tr_F$, where $\Tr_F$ is the trace map, extended to $\mathbb{A}_F\rightarrow \mathbb{A}_{\mathbb{Q}}$. Then $\psi(t)=\prod_{v\in\Sigma_F}\psi_v(t_v)$ for $t=(t_v)_v\in\mathbb{A}_F$. Let $\theta$: $[N]\longrightarrow \mathbb{C}^{\times}$ be the generic character induced from $\psi$ via $\theta\left(\begin{pmatrix}
	1& b\\
	&1
\end{pmatrix}\right)=\psi(b)$.

\subsubsection{Automorphic Representations}
Let $\omega$ be a unitary Hecke character over $F$. We denote by $\mathcal{A}([G],\omega)$ the set of non-equivalent unitary automorphic representations of $G(\mathbb{A}_F)$ with central character $\omega$. Let $\mathcal{A}_0([G],\omega)$ be the subset of $\mathcal{A}([G],\omega)$ consisting of cuspidal representations. For $\sigma\in \mathcal{A}_0([G],\omega)$, we denote by $\mathfrak{B}(\sigma)$ an orthonormal basis of $\sigma$ consisting of $K$-isotypical pure tensors.

For an integral ideal $\mathfrak{q}$, we denote by $\mathcal{F}(\mathfrak{q};\omega)$ the set of unitary generic 
automorphic representations of $\mathrm{GL}_2/F$ with central 
character $\omega$ and conductor dividing $\mathfrak{q}$. Let $\mathcal{F}_0(\mathfrak{q};\omega):=\mathcal{F}(\mathfrak{q};\omega)\cap \mathcal{A}_0([G],\omega)$.  

\subsubsection{Other Conventions}
%For a function $h$ on $G(\mathbb{A}_F)$, we define $h^*$ by assigning $h^*(g)=\overline{h({g}^{-1})}$, $g\in G(\mathbb{A}_F)$. Let $F_1(s), F_2(s)$ be two meromorphic functions. Write $F_1(s)\propto F_2(s)$ if there exists an \textit{entire} function $E(s)$ such that $F_1(s)=E(s)F_2(s)$. Denote by $\alpha\asymp \beta$ for $\alpha, \beta \in\mathbb{R}$ if there are absolute constants $c$ and $C$ such that $c\beta\leq \alpha\leq C\beta$.

Let $\mathcal{S}(\mathbb{A}_F^2)$ denote the space of Schwartz-Bruhat functions on $\mathbb{A}_F^2$. Let $s\in \mathbb{C}$. We use the notation $\Re(s)\gg 1$ to indicate that $\Re(s)$ is large enough, say, $\Re(s)>10$. Let $F_1(s), F_2(s)$ be two meromorphic functions. Write $F_1(s)\propto F_2(s)$ if there exists an \textit{entire} function $E(s)$ such that $F_1(s)=E(s)F_2(s)$, in which case we also say that $F_1$ is a \textit{holomorphic multiple} of $F_2$. 

For an automorphic representation $\Pi$, we denote the complete $L$-function by $\Lambda(s,\Pi)$ and the finite part, with the archimedean factors removed, by $L(s,\Pi)$.

Throughout, we adhere to the $\varepsilon$-convention, where $\varepsilon$ is always a positive number that can be arbitrarily small, although its magnitude may vary between different instances.

\subsection{Hecke operators}\label{sec2.2}
Let $\sigma=\otimes_{v\in\Sigma}\sigma_v$ be an automorphic representation of $G(\mathbb{A}_F)$, with central character $\omega$. Let $\mathcal{P}$ be a finite subset of primes in $\Sigma_{\fin}$. Let $\mathfrak{L}=\prod_{v\in \mathcal{P}}\mathfrak{p}_v^{\ell_v}$, $\ell_v\in \mathbb{Z}_{\geq 0}$. The $\mathfrak{L}$-th Hecke operator is defined by 
\begin{equation}\label{hecke}
T_{\mathfrak{L}}(\varphi):=\prod_{v\in \mathcal{P}}q_v^{-\frac{\ell_v}{2}}\int_{K_v\diag(\varpi_v^{\ell_v},1)K_v} \sigma_v(y_v)\varphi dy_v,\ \ \varphi\in \sigma.
\end{equation} 

When $\varphi$ is right invariant under $\prod_{v\in \mathcal{P}}K_v$, it becomes an eigenform for the operator $T_{\mathfrak{L}}=\prod_{v\in \mathcal{P}}T_{\mathfrak{p}_v^{\ell_v}}$. We denote the eigenvalue of $T_{\mathfrak{p}_v^{\ell_v}}$ as $\lambda_{\sigma}(\mathfrak{p}_v^{\ell_v})$. 

\subsubsection{Hecke relations}
Let $v<\infty$. We have the Hecke relation
\begin{equation}\label{hecke.}
\lambda_{\sigma}(\mathfrak{p}_v^{\ell_v})^2=\sum_{n_v=0}^{\ell_v}\omega_v(\varpi_v^{n_v})\lambda_{\sigma}(\mathfrak{p}_v^{2\ell_v-2n_v}).
\end{equation}

Let $\omega_v^{-1/2}$ be a branch of the square root of $\omega$. It is well known that the operator $\omega_v^{-1/2}T_{\mathfrak{p}_v^{\ell_v}}$ is Hermitian. Hence, 
\begin{equation}\label{eq2.3}
\overline{\lambda_{\sigma}(\mathfrak{p}_v^{\ell_v})}=\omega_v^{-1}(\varpi_v^{\ell_v})\cdot \lambda_{\sigma}(\mathfrak{p}_v^{\ell_v}).
\end{equation}

\subsubsection{The Ramanujan bounds}\label{sec2.2.2}
Let $\vartheta$ be a parameter towards the Ramanujan conjecture for unitary cuspidal automorphic representations of $\mathrm{GL}(2)/F$, namely, for any $v<\infty$ and any unitary cuspidal automorphic representation $\sigma=\otimes_v\sigma_v$ with $v<\infty$ and $L_v(s,\sigma_v)=(1-\alpha_vq_v^{-s})^{-1}(1-\beta_vq_v^{-s})^{-1}$, we have  
\begin{equation}\label{eq6.6}
\max\{|\alpha_v|,|\beta_v|\}<q_v^{\vartheta}.
\end{equation}

In particular, by \cite{KS03} and \cite{BB11} we can take $\vartheta= 7/64$.

\subsection{Eisenstein Series and Spectral Decomposition}
\subsubsection{Eisenstein series}\label{sec2.3.1}
Let $\xi_1, \xi_2\in \widehat{F^{\times}\backslash \mathbb{A}_F^{(1)}}$ and $s\in \mathbb{C}$. Let $H(\xi_1,\xi_2,s)$ be the space of functions $f$ on $G(\mathbb{A}_F)$ satisfying $f\in L^2(K)$ and 
\begin{align*}
f\left(\begin{pmatrix}
a&b\\
&d
\end{pmatrix}g\right)=\xi_1(a)\xi_2(d)|a/d|^{s+1/2}f(g),\ \ a, d\in \mathbb{A}_F^{\times}, \ b\in \mathbb{A}_F,\ g\in G(\mathbb{A}_F). 
\end{align*}
We denote by $\pi_{\xi_1,\xi_2,s}$ the representation of  $G(\mathbb{A}_F)$ given by  the  right translation on $H(\xi_1,\xi_2,s)$. 

Let $H(\xi_1,\xi_2)$ be the space defined by 
\begin{align*}
\Big\{f\in L^2(K):\ f\left(\begin{pmatrix}
a& b\\
&d
\end{pmatrix}k\right)=\xi_1(a)\xi_2(d)f(k),\ \forall\ k\in K,\ \begin{pmatrix}
a& b\\
&d
\end{pmatrix}\in B(\mathbb{A}_F)\cap K
\Big\}.
\end{align*}

Given $f\in H(\xi_1,\xi_2)$ and $s\in \mathbb{C}$, we define 
\begin{equation}\label{eq2.2}
S(f)\left(\begin{pmatrix}
a& b\\
&d
\end{pmatrix}k,s\right):=\xi_1(a)\xi_2(d)\Big|\frac{a}{d}\Big|^{s+1/2}f(k),
\end{equation}
where $a,d\in \mathbb{A}_F^{\times}$, $b\in \mathbb{A}_F$, and $k\in K$. Then 
\begin{align*}
S:\ H(\xi_1,\xi_2)\longrightarrow H(\xi_1,\xi_2,s), \ \ f\mapsto S(f)(\cdot,s)
\end{align*}
is an isomorphism, giving the trivialization of the fiber bundle $\bigcup_{s\in \mathbb{C}}H(\xi_1,\xi_2,s)$ over $\mathbb{C}$. As a consequence, 
\begin{align*}
H(\xi_1,\xi_2,s)=\Big\{S(f)(\cdot,s):\ f\in H(\xi_1,\xi_2)\Big\}.
\end{align*}
Through the transformation $S$ we may regard $\pi_{\xi_1,\xi_2,s}$ as a representation of $G(\mathbb{A}_F)$ on $H(\xi_1,\xi_2)$, and thus 
\begin{equation}\label{2.2}
\pi_{\xi_1,\xi_2,s}(g)f(k)=S(f)(kg,s),\ \ f\in H(\xi_1,\xi_2),\ k\in K,\ g\in G(\mathbb{A}_F). 
\end{equation} 

Define the associated Eisenstein series by
\begin{align*}
E(g,f,s):=\sum_{\gamma\in B(F)\backslash G(F)}S(f)(\gamma g,s),\ \ \Re(s)>1/2.
\end{align*}

It is known that $E(g,f,s)$ converges absolutely in $\Re(s)>1/2$, and admits a meromorphic continuation to $\mathbb{C}$, with only one possible simple pole at $s=1/2$ in $\Re(s)\geq 0$. Henceforth, we regard $E(g,f,s)$ as a meromorphic function. 

For later purpose, we fix once and for all an orthonormal basis of $K$-finite functions $\mathfrak{B}(\xi_1,\xi_2)$ of the Hilbert space $H(\xi_1,\xi_2)$.

\begin{comment}
by the right translation on the space $H(\xi_1,\xi_2,s)$ of classes of functions on $G(\mathbb{A}_F)$ satisfying $\int_K|f(k)|^2dk<\infty$, and 
\begin{align*}
f\left(\begin{pmatrix}
a& b\\
&d
\end{pmatrix}g\right)=\xi_1(a)\xi_2(d)|a/d|^{s+1/2}f(g)
\end{align*} 
for all $a,d\in \mathbb{A}_F^{\times}$, $b\in \mathbb{A}_F$, and $g\in G(\mathbb{A}_F)$.

\bigskip 
\bigskip 
\bigskip

For $a,d\in\mathbb{A}_F^{\times}$, define $H(\diag(a,d))=\log|ad^{-1}|$. Let $g\in G(\mathbb{A}_F)$. By Iwasawa decomposition $g=umk$, where $u\in N(\mathbb{A}_F)$, $m\in T(\mathbb{A}_F)$, and $k\in K$. Define the height function as $H(g):=H(m)$.

Let $\xi_1, \xi_2\in \widehat{F^{\times}\backslash \mathbb{A}_F^{(1)}}$. Let $\varphi\in \Ind_B^G(\xi_1,\xi_2)$, and $\mu\in \mathbb{C}$ with $\Re(\mu)>1/2$. Define the Eisenstein series by
\begin{equation}\label{2.1}
E(x,\varphi,\mu):=\sum_{\delta\in B(F)\backslash G(F)}\varphi(\delta g)e^{(\mu+\frac{1}{2})H(\delta x)}.
\end{equation}

It is well known that $E(x,\varphi,\mu)$ converges absolutely in $\Re(\mu)>1/2$, and admits a meromorphic continuation to $\mathbb{C}$, with only one possible simple pole at $s=1/2$ in $\Re(s)\geq 0$. Henceforth, we regard $E(x,\varphi,\mu)$ as a meromorphic function. 
\end{comment}

\subsubsection{The spectral decomposition}\label{sec2.3.2}
Denote by $\langle\cdot,\cdot\rangle$ the Petersson inner product equipped in the space $L^2\left([G],\omega\right)$, where $\omega$ is a unitary Hecke character. 
Let $\varphi\in L^2\left([G],\omega\right)$ be a smooth vector. Utilizing the spectral decomposition (see \cite{GJ79} or \cite[I.1.4]{Cog07}), we express $\varphi$ as in \cite[Theorem 2.16]{Wu14}: 
\begin{equation}\label{2.6}
\varphi(g)=\frac{1}{\Vol([\overline{G}])}\sum_{\substack{\eta\in \widehat{F^{\times}\backslash\mathbb{A}_F^{(1)}},\ \eta^2=\omega}}\langle \varphi,\eta\circ\det\rangle \eta(\det g)+\varphi_{\cusp}(g)+\varphi_{\Eis}(g),
\end{equation}
which converges absolutely and uniformly on any compact subset. Here
\begin{align*}
&\varphi_{\cusp}(g):=\sum_{\sigma\in \mathcal{A}_0([G],\omega)}\sum_{\phi\in\mathfrak{B}(\sigma)}\langle \varphi, \phi\rangle \phi(g),\\
&\varphi_{\Eis}(g):=\sum_{\substack{\xi\in \widehat{F^{\times}\backslash\mathbb{A}_F^{(1)}}}} \frac{1}{4\pi i}\int_{i\mathbb{R}}\sum_{f\in \mathfrak{B}(\xi,\xi^{-1}\omega)} \langle \varphi,E(\cdot,f,s)\rangle E(g,f,s)ds.
\end{align*}
If $\varphi$ is $K$-finite, then the sums over $\mathfrak{B}(\sigma)$ and $\mathfrak{B}(\xi,\xi^{-1}\omega)$ in the above definitions are finite, depending on the $K$-type of $\varphi$. 

\subsection{Automorphic Forms and Rankin--Selberg Convolution}
\subsubsection{Eisenstein series}\label{sec3.1.2}
Let $\chi_1$ and $\chi_2$ be  unitary Hecke characters over $F$. Let $\Phi=\otimes_{v\in \Sigma}\Phi_v$ be a Schwartz-Bruhat function on $\mathbb{A}_F^2$. For $\Re(s)>1/2$, we define 
\begin{equation}\label{2.3}
h_{\Phi,\chi_1,\chi_2}(g,s):=\chi_1(\det g)|\det g|^s\int_{\mathbb{A}_F^{\times}}\Phi((0,t)g)\chi_1\chi_2^{-1}(t)|t|^{2s}d^{\times}t,
\end{equation}
where $g\in G(\mathbb{A}_F)$. Notice that $h_{\Phi,\chi_1,\chi_2}(\cdot,s)\in H(\chi_1,\chi_2,s-1/2)$, as
\begin{equation}\label{eq2.6}
h_{\Phi,\chi_1,\chi_2}\left(\begin{pmatrix}
a& b\\
&d
\end{pmatrix}g,s\right)=\chi_1(a)\chi_2(d)|a|^s|d|^{-s}h_{\Phi,\chi_1,\chi_2}(g,s).
\end{equation}

For $\Re(s)>1$, we define the Eisenstein series by 
\begin{equation}\label{2.6.}
E_{\Phi,\chi_1,\chi_2}(g,s):=\sum_{\delta\in B(F)\backslash G(F)}h_{\Phi,\chi_1,\chi_2}(\delta g,s).
\end{equation}

By Poisson summation, $E(g,s;\Phi,\chi)$ admits a meromorphic continuation to $s\in\mathbb{C}$, with a functional equation 
\begin{align*}
E_{\Phi,\chi_1,\chi_2}(g,s)=E_{\widehat{\Phi},\chi_1^{-1},\chi_2^{-1}}(g^{\iota},1-s),
\end{align*}
where $g^{\iota}$ is the transverse inverse of $g$, and $\widehat{\Phi}$ is the Fourier transform of $\Phi$.

\subsubsection{Rankin--Selberg periods}\label{sec2.4.2}
Let $\varphi$ be an automorphic form on $[G]$. The associated Whittaker function with respect to $\theta$ is defined by 
\begin{equation}\label{2.4}
W_{\varphi}(g;\theta):=\int_{[N]}\varphi(ug)\overline{\theta(u)}du.
\end{equation} 

Substituting the definition of $h_{\Phi,\chi_1,\chi_2}(g,s)$ into the above integral, we derive that $W_{E_{\Phi,\chi_1,\chi_2}(\cdot,s)}(g)$ is equal to 
\begin{equation}\label{f2.3}
\chi_1(\det g)|\det g|^{s}\int_{\mathbb{A}_F^{\times}}\int_{\mathbb{A}_F}\Phi((t,bt)g)\overline{\psi}_F(b)db\chi_1\chi_2^{-1}(t)|t|^{2s}d^{\times}t,
\end{equation}
which admits a meromorphic continuation to $s\in \mathbb{C}$.

Let $\Phi\in \mathcal{S}(\mathbb{A}_F^2)$ be a  Schwartz-Bruhat function on $\mathbb{A}_F^2$. Let $\varphi_1$ and $\varphi_2$ be generic automorphic forms on $[G]$ with central characters $\omega_1$ and $\omega_2$, respectively. Let $\chi_1$ and $\chi_2$ be unitary Hecke characters with $\omega_1\omega_2\chi_1\chi_2=\mathbf{1}$. Let 
\begin{equation}\label{2.12}
\Psi(s,\varphi_1,\varphi_2,h_{\Phi,\chi_1,\chi_2}):=\int_{N(\mathbb{A}_F)\backslash \overline{G}(\mathbb{A}_F)}W_{\varphi_1}(x;\theta)W_{\varphi_2}(x;\overline{\theta})h_{\Phi,\chi_1,\chi_2}(x,s)dx
\end{equation}
be the Rankin--Selberg convolution in the Whittaker form, where $h_{\Phi,\chi_1,\chi_2}(x,s)$ is defined by \eqref{2.3}. It converges absolutely in $\Re(s)\gg 1$ and  admits a meromorphic continuation to $s\in \mathbb{C}$. 

Let $\Re(s)\gg 1$. Utilizing the condition $\omega_1\omega_2\chi_1\chi_2=\mathbf{1}$, we obtain  
\begin{align*}
\Psi(s,\varphi_1,\varphi_2,h_{\Phi,\chi_1,\chi_2})=\int W_{\varphi_1}(x;\theta)W_{\varphi_2}(x;\overline{\theta})\chi_1(\det x)\Phi((0,1)x)|\det x|^sdx,
\end{align*}
where the integral is over $x\in N(\mathbb{A}_F)\backslash G(\mathbb{A}_F)$.

\begin{comment}
\begin{align*}
\Psi(s,\varphi_1,\varphi_2,E_{\Phi,\chi_1,\chi_2}(\cdot,s)):=&\int_{N(\mathbb{A}_F)\backslash \overline{G}(\mathbb{A}_F)}W_{\varphi_1}(x;\theta)W_{\varphi_2}(x;\overline{\theta})\\
&\chi_1(\det x)|\det x|^s\int_{\mathbb{A}_F^{\times}}\Phi((0,t)x)\chi_1\chi_2^{-1}(t)|t|^{2s}d^{\times}tdx\\
=&\int_{\mathbb{A}_F^{\times}}\int_{N(\mathbb{A}_F)\backslash \overline{G}(\mathbb{A}_F)}W_{\varphi_1}(tx;\theta)W_{\varphi_2}(tx;\overline{\theta})\overline{\Omega_2(t)\omega_2(t)}\\
&\chi_1(\det tx)|\det tx|^s\Phi((0,1)tx)\chi_1^{-1}\chi_2^{-1}(t)dxd^{\times}t\\
=&\int_{N(\mathbb{A}_F)\backslash G(\mathbb{A}_F)}W_{\varphi_1}(x;\theta)W_{\varphi_2}(x;\overline{\theta})\chi_1(\det x)\Phi((0,1)x)|\det x|^sdx.
\end{align*}

Notice that 
\begin{align*}
\Psi(s,\varphi_1,\varphi_2,E_{\Phi,\chi_1,\chi_2}(\cdot,s))=\int_{N(F)\backslash \overline{G}(\mathbb{A}_F)}W_{\varphi_1}(x;\theta)\varphi_2(x)h_{\Phi,\chi_1,\chi_2}(x,s)dx.
\end{align*}

\subsubsection{Meromorphic Continuation}
Let $\Psi(s,\varphi_1,\varphi_2,\Phi)$  be defined as in \textsection\ref{sec2.4.2}. Suppose that $\varphi_1\in \pi_{\xi_1,\xi_2,s}$ and $\varphi_2\in \pi_{\xi_1',\xi_2',s'}$. Then $\Omega_2=\xi_1\xi_2$ and $\omega_2=\xi_1'\xi_2'$. 
\end{comment}

\section{Regularized Spectral Reciprocity}\label{sec3}
Let $\omega, \omega', \chi_1, \chi_2, \chi_1', \chi_2'$ be unitary Hecke characters over $F$. Let $\pi'=\otimes_v'\pi_v'\in L_0^2([G],\omega')$. Let $\phi_1$ and $\phi_2$ be two $K$-finite vectors in $\pi'$. Let $\Phi_1,\Phi_2\in \mathcal{S}(\mathbb{A}_F^2)$. Denote by $E_1(g,s):=E_{\Phi_1,\chi_1,\chi_2}(g,s)$, and $E_2(g,s):=E_{\Phi_2,\chi_1',\chi_2'}(g,s)$. For simplicity we write $h_2(\cdot,s)=h_{\Phi_2,\chi_1',\chi_2'}(\cdot,s)$, which is defined by \eqref{2.3}. Suppose that $\chi_1\chi_2=\chi_1'\chi_2'$. 

In this section we aim to establish a general regularized spectral reciprocity formula based on two expansions of the inner product $\langle \phi_1E_1(\cdot,s_1),\phi_2E_2(\cdot,\overline{s_2})\rangle$. A detailed statement is given in Theorem \ref{thmA} in \textsection\ref{sec3.5}.

\subsection{Decomposition of the Inner Product}\label{sec3.1}
\subsubsection{The twisted moment}\label{sec3.1.1}
Let $\Re(s_1)\gg 1$ and $\Re(s_2)\gg 1$. Following the notation in \textsection\ref{sec2.3.2} we define
\begin{align*}
\mathcal{M}_{\cusp}^{\tw}(s_1,s_2):=\sum_{\sigma\in \mathcal{A}_0([G],\omega)}\sum_{\varphi\in\mathfrak{B}(\sigma)}\langle \phi_1E_1(\cdot,s_1),\varphi\rangle\langle \varphi,\phi_2E_2(\cdot,\overline{s_2})\rangle;
\end{align*}
and define the function $\mathcal{M}_{\Eis}^{\tw}(s_1,s_2)$ by 
\begin{align*}
\sum_{\substack{\xi\in \widehat{F^{\times}\backslash\mathbb{A}_F^{(1)}}}} \frac{1}{4\pi i}\int_{i\mathbb{R}}\sum_{f\in \mathfrak{B}(\xi,\xi^{-1})} \langle \phi_1E_1(\cdot,s_1),E(\cdot,f,s)\rangle \langle E(\cdot,f,s),\phi_2E_2(\cdot,\overline{s_2})\rangle ds.
\end{align*}
We call $\mathcal{M}_{\cusp}^{\tw}(s_1,s_2)$ (resp. $\mathcal{M}_{\Eis}^{\tw}(s_1,s_2)$) the cuspidal part (resp. the continuous part) of the twisted moment. 

\subsubsection{The dual moment}\label{sec3.1.2.}
Let $\Re(s_1)\gg 1$ and $\Re(s_2)\gg 1$. Write  $h_2:=h_2(\cdot,\overline{s_2})$ for simplicity. We define 
\begin{equation}\label{f3.1}
\mathcal{M}_{\cusp}^{\du}(s_1,s_2):=\sum_{\sigma\in \mathcal{A}_0([G],\mathbf{1})}\sum_{\varphi\in\mathfrak{B}(\sigma)}\langle \phi_1\overline{\phi_2},\varphi\rangle\Psi(s_2,\varphi,E_1(\cdot,s_1),\overline{h_2});
\end{equation}
and define the function $\mathcal{M}_{\Eis}^{\du}(s_1,s_2)$ by 
\begin{align*}
\sum_{\substack{\xi\in \widehat{F^{\times}\backslash\mathbb{A}_F^{(1)}}}} \frac{1}{4\pi i}\int_{i\mathbb{R}}\sum_{f\in \mathfrak{B}(\xi,\xi^{-1})} \langle \phi_1\overline{\phi_2},E(\cdot,f,s)\rangle \Psi(s_2,E(\cdot,f,s),E_1(\cdot,s_1),\overline{h_2})ds.
\end{align*}
Here $\Psi(s_2,E(\cdot,f,s),E_1(\cdot,s_1),\overline{h_2})$ is defined as in \textsection\ref{sec2.4.2}. Define $\mathcal{M}_{\const}^{\du}(s_1,s_2)$ by 
\begin{align*}
\int_{A(F)N(\mathbb{A}_F)\backslash \overline{G}(\mathbb{A}_F)}\int_{[N]}\phi_1(u_1x)\overline{\phi_2(u_1x)}du_1\int_{[N]}E_1(u_2x,s_1)du_2\overline{h_2(x,\overline{s_2})}dx.
\end{align*}

We call $\mathcal{M}_{\cusp}^{\du}(s_1,s_2)$ (resp. $\mathcal{M}_{\Eis}^{\du}(s_1,s_2)$) the cuspidal part (resp. the continuous part) of the dual moment. Call $\mathcal{M}_{\const}^{\du}(s_1,s_2)$ the constant part of the dual moment. 

\subsubsection{A coarse spectral reciprocity}
\begin{prop}\label{prop3.1}
Let notation be as before. Let $\Re(s_1)\gg 1$ and $\Re(s_2)\gg 1$. Then 
\begin{equation}\label{3.1}
\mathcal{M}_{\cusp}^{\tw}(s_1,s_2)+\mathcal{M}_{\Eis}^{\tw}(s_1,s_2)=\mathcal{M}_{\const}^{\du}(s_1,s_2)+\mathcal{M}_{\cusp}^{\du}(s_1,s_2)+\mathcal{M}_{\Eis}^{\du}(s_1,s_2).
\end{equation}	
\end{prop}
\begin{proof}
By definition, we have
\begin{align*}
\langle \phi_1E_1(\cdot,s_1),\phi_2E_2(\cdot,\overline{s_2})\rangle=\int_{[\overline{G}]}\phi_1(x)\overline{\phi_2(x)}E_1(x,s_1)\overline{E_2(x,\overline{s_2})}dx.
\end{align*}
Write $h_2(x,s)$ for $h(g,s;\Phi_2,\overline{\omega\omega'})$. By unfolding $E_2(x,\overline{s_2})$,
\begin{align*}
\langle \phi_1E_1(\cdot,s_1),\phi_2E_2(\cdot,\overline{s_2})\rangle=\int_{B(F)Z(\mathbb{A}_F)\backslash G(\mathbb{A}_F)}\phi_1(x)\overline{\phi_2(x)}E_1(x,s_1)\overline{h_2(x,\overline{s_2})}dx.
\end{align*}

Substituting the Fourier expansion 
\begin{align*}
\phi_1(x)\overline{\phi_2(x)}=\int_{[N]}\phi_1(ux)\overline{\phi_2(ux)}du+\sum_{\delta\in A(F)}\int_{[N]}\phi_1(u\delta x)\overline{\phi_2(u\delta x)}\overline{\theta(u)}du
\end{align*}
into the above identity we obtain 
\begin{equation}\label{3.2}
\langle \phi_1E_1(\cdot,s_1),\phi_2E_2(\cdot,\overline{s_2})\rangle=\mathcal{M}_{\const}^{\du}(s_1,s_2)+\mathcal{M}_{\gen}^{\du}(s_1,s_2),
\end{equation}
where
\begin{align*}
\mathcal{M}_{\gen}^{\du}(s_1,s_2):=\int_{N(F)\backslash \overline{G}(\mathbb{A}_F)}\int_{[N]}\phi_1(ux)\overline{\phi_2(ux)}\overline{\theta(u)}duE_1(x,s_1)\overline{h_2(x,\overline{s_2})}dx.
\end{align*}

By a change of variable and the fact that $h_2(ux,\overline{s_2})=h_2(x,\overline{s_2})$ for all $u\in N(\mathbb{A}_F)$, we derive that 
\begin{align*}
\mathcal{M}_{\gen}^{\du}(s_1,s_2)=\int_{N(\mathbb{A}_F)Z(\mathbb{A}_F)\backslash G(\mathbb{A}_F)}W_{\phi_1\overline{\phi_2}}(x;\theta)W_{E_1(\cdot,s_1)}(x;\overline{\theta})\overline{h_2(x,\overline{s_2})}dx,
\end{align*}
where $W_{\phi_1\overline{\phi_2}}(x;\theta)$ and $W_{E_1(\cdot,s_1)}(x;\overline{\theta})$ are the Whittaker functions defined by \eqref{2.4}. 

Since $\phi_1\overline{\phi_2}$ decays rapidly in $[G]$, $\mathcal{M}_{\gen}^{\du}(s_1,s_2)$ converges absolutely. We may execute the spectral decomposition (see \textsection\ref{sec2.3.2}) to the function $\phi_1\overline{\phi_2}$, obtaining
\begin{equation}\label{3.3}
\mathcal{M}_{\gen}^{\du}(s_1,s_2)=\mathcal{M}_{\cusp}^{\du}(s_1,s_2)+\mathcal{M}_{\Eis}^{\du}(s_1,s_2).
\end{equation}
Notice that the residual spectrum does not contribute to $\mathcal{M}_{\gen}^{\du}(s_1,s_2)$. 

Since $\Re(s_1)\gg 1$ and $\Re(s_2)\gg 1$, along with the fact that Whittaker functions decay rapidly in vertical strips, both $\mathcal{M}_{\cusp}^{\du}(s_1,s_2)$ and $\mathcal{M}_{\Eis}^{\du}(s_1,s_2)$ converges absolutely. Therefore, by \eqref{3.2} and \eqref{3.3}, we have
\begin{equation}\label{eq3.4}
\langle \phi_1E_1(\cdot,s_1),\phi_2E_2(\cdot,\overline{s_2})\rangle=\mathcal{M}_{\const}^{\du}(s_1,s_2)+\mathcal{M}_{\cusp}^{\du}(s_1,s_2)+\mathcal{M}_{\Eis}^{\du}(s_1,s_2).
\end{equation}

The formula \eqref{3.1} follows from \eqref{eq3.4} and the spectral decomposition of $\phi_1E_1(\cdot,s_1)$ given in \eqref{2.6}. Since $\phi_1E_1(\cdot,s_1)$ only project to the generic representations, the residual spectrum does not contribute to the decomposition. 
\end{proof}

\subsection{Meromorphic Continuation of $\mathcal{M}_{\const}^{\du}(s_1,s_2)$}\label{sec3.2}
Recall that the constant part of the dual moment  $\mathcal{M}_{\const}^{\du}(s_1,s_2)$ is defined by 
\begin{align*}
\int_{A(F)N(\mathbb{A}_F)\backslash \overline{G}(\mathbb{A}_F)}\int_{[N]}\phi_1(u_1x)\overline{\phi_2(u_1x)}du_1\int_{[N]}E_1(u_2x,s_1)du_2\overline{h_2(x,\overline{s_2})}dx.
\end{align*}

Let $\Re(s_1)>1$. By Bruhat decomposition, we have
\begin{align*}
\int_{[N]}E_1(u_2x,s_1)du_2=h_1(x,s_1)+\int_{N(\mathbb{A}_F)}h_1(wux,s_1)du,
\end{align*}
where $h_1(\cdot,s_1):=h_{\Phi_1,\chi_1,\chi_2}(\cdot,s_1)$. Consequently,
\begin{equation}\label{3.5..}
\mathcal{M}_{\const}^{\du}(s_1,s_2)=\mathcal{M}_{\const,1}^{\du}(s_1,s_2)+\mathcal{M}_{\const,2}^{\du}(s_1,s_2),
\end{equation}
where
\begin{align*}
\mathcal{M}_{\const,1}^{\du}(s_1,s_2):=&\int_{A(F)N(\mathbb{A}_F)\backslash \overline{G}(\mathbb{A}_F)}\int_{[N]}\phi_1(u_1x)\overline{\phi_2(u_1x)}du_1h_1(x,s_1)\overline{h_2(x,\overline{s_2})}dx,
\end{align*}
and $\mathcal{M}_{\const,2}^{\du}(s_1,s_2)$ is defined by 
\begin{align*}
\int_{A(F)N(\mathbb{A}_F)\backslash \overline{G}(\mathbb{A}_F)}\int_{[N]}\phi_1(u_1x)\overline{\phi_2(u_1x)}du_1\int_{N(\mathbb{A}_F)}h_1(wux,s_1)du\overline{h_2(x,\overline{s_2})}dx.
\end{align*}

\begin{lemma}\label{lem3.6.}
Let notation be as before. Then $\mathcal{M}_{\const,1}^{\du}(s_1,s_2)$ converges absolutely in $\big\{(s_1,s_2)\in \mathbb{C}^2:\ \Re(s_1)\gg 1,\ \Re(s_2)\gg 1\big\}$, and it admits a meromorphic continuation $\mathcal{M}_{\const,1}^{\du,\heartsuit}(s_1,s_2)$ to $\mathbb{C}^2$, satisfying
\begin{equation}\label{eq3.5.}
\mathcal{M}_{\const,1}^{\du,\heartsuit}(s_1,s_2)\propto \frac{\Lambda(s_1+s_2,\pi'\times\overline{\pi}'\otimes\chi_1\overline{\chi_1'})\Lambda(2s_1,\chi_1\overline{\chi_2})\Lambda(2s_2,\overline{\chi_1'}\chi_2')}{\Lambda(2s_1+2s_2,\chi_1^2\chi_1'^{-2})}.
\end{equation}
\end{lemma}
\begin{proof}
Since $h_1(\cdot,s_1)$ and $h_2(\cdot,\overline{s_2})$ are invariant under the left multiplication by $N(\mathbb{A}_F)$, then 
\begin{align*}
\mathcal{M}_{\const,1}^{\du}(s_1,s_2)=\int_{B(F)Z(\mathbb{A}_F)\backslash G(\mathbb{A}_F)}\phi_1(x)\overline{\phi_2(x)}h_1(x,s_1)\overline{h_2(x,\overline{s_2})}dx,
\end{align*}
which converges absolutely in $\Re(s_1+s_2)\gg 1$.

By the Fourier expansion of $\phi_1$ we obtain 
\begin{align*}
\mathcal{M}_{\const,1}^{\du}(s_1,s_2)=\int_{N(\mathbb{A}_F)Z(\mathbb{A}_F)\backslash G(\mathbb{A}_F)}W_{\phi_1}(x;\theta)\overline{W_{\phi_2}(x;\theta)}h_1(x,s_1)\overline{h_2(x,\overline{s_2})}dx.
\end{align*}

Using the Iwasawa decomposition, 
\begin{align*}
\mathcal{M}_{\const,1}^{\du}(s_1,s_2)=&\int_K\int_{A(\mathbb{A}_F)}W_{\phi_1}(ak;\theta)\overline{W_{\phi_2}(ak;\theta)}\chi_1\chi_1'^{-1}(\det a)\\
&\qquad |\det a|^{s_1+s_2-1}h_1(k,s_1)\overline{h_2(k,\overline{s_2})}d^{\times}adk,
\end{align*}
which converges absolutely if $\Re(s_1)\gg 1$, and $\Re(s_2)\gg 1$.

By the definition \eqref{2.3}, $h_1(k,s_1)\propto \Lambda(2s_1,\chi_1\chi_2^{-1})$, $\overline{h_2(x,\overline{s_2})}\propto \Lambda(2s_2,\chi_1'^{-1}\chi_2')$. Applying the Casselman-Shalika formula, the inner integral over $a$ is a holomorphic multiple of the quotient $\Lambda(s_1+s_2,\pi'\times\overline{\pi}'\otimes\chi_1\chi_1'^{-1})/\Lambda(2s_1+2s_2,\chi_1^2\chi_1'^{-2})$. Therefore, $\mathcal{M}_{\const,1}^{\du}(s_1,s_2)$ converges absolutely in $\big\{(s_1,s_2)\in \mathbb{C}^2:\ \Re(s_1)\gg 1,\ \Re(s_2)\gg 1\big\}$, and admits a meromorphic continuation $\mathcal{M}_{\const,1}^{\du,\heartsuit}(s_1,s_2)$ to $\mathbb{C}^2$, satisfying \eqref{eq3.5.}. So Lemma \ref{lem3.6.} holds. 
\end{proof}

\begin{lemma}\label{lem3.7}
Let notation be as before. Then $\mathcal{M}_{\const,2}^{\du}(s_1,s_2)$ converges absolutely in $\big\{(s_1,s_2)\in \mathbb{C}^2:\ \Re(s_1)\gg 1,\ \Re(s_2)\gg 1\big\}$, and it admits a meromorphic continuation $\mathcal{M}_{\const,2}^{\du,\heartsuit}(s_1,s_2)$ to $\mathbb{C}^2$, satisfying
\begin{equation}\label{eq3.5}
\mathcal{M}_{\const,2}^{\du,\heartsuit}(s_1,s_2)\propto \frac{\Lambda(1+s_2-s_1,\pi'\times\overline{\pi}'\otimes \chi_2')\Lambda(2s_1-1,\chi_1\overline{\chi_2})\Lambda(2s_2,\overline{\chi_1'}\chi_2')}{\Lambda(2+2s_2-2s_1,\chi_2'^{2})}.
\end{equation}
\end{lemma}
\begin{proof}
By definition, $\mathcal{M}_{\const,1}^{\du}(s_1,s_2)$ is equal to 
\begin{align*}
\int_{B(F)\backslash \overline{G}(\mathbb{A}_F)}\phi_1(x)\overline{\phi_2(x)}\int_{N(\mathbb{A}_F)}h_1(wux,s_1)du\overline{h_2(x,\overline{s_2})}dx.
\end{align*}

By the Fourier expansion of $\phi_1$ the function $\mathcal{M}_{\const,1}^{\du}(s_1,s_2)$ transfers to  
\begin{equation}\label{3.15}
\int_{N(\mathbb{A}_F)\backslash \overline{G}(\mathbb{A}_F)}W_{\phi_1}(x;\theta)\overline{W_{\phi_2}(x;\theta)}\int_{N(\mathbb{A}_F)}h_1(wux,s_1)du\overline{h_2(x,\overline{s_2})}dx.
\end{equation}

Let $a\in \mathbb{A}_F^{\times}$. By \eqref{eq2.6}, we have
\begin{equation}\label{3.16}
\int_{N(\mathbb{A}_F)}h_1\left(
wu\begin{pmatrix}
a&\\
&1
\end{pmatrix}x,s_1\right)du=\chi_2(a)|a|^{1-s_1}\int_{N(\mathbb{A}_F)}h_1(wux,s_1)du.
\end{equation}

Using the Iwasawa decomposition in \eqref{3.15}, along with \eqref{3.16}, we obtain 
\begin{align*}
\mathcal{M}_{\const,2}^{\du}(s_1,s_2)=\int_K &\int_{\mathbb{A}_F^{\times}}W_{\phi_1}\left(\begin{pmatrix}
a\\
&1
\end{pmatrix}k
;\theta\right)\overline{W_{\phi_2}\left(\begin{pmatrix}
a\\
&1
\end{pmatrix}k
;\theta\right)}\\
&\int_{N(\mathbb{A}_F)}h_1(wuk,s_1)du\overline{h_2(k,\overline{s_2})}\chi_1\chi_1'^{-1}\chi_2(a)|a|^{s_2-s_1}d^{\times}adk,
\end{align*}
which converges absolutely if $\Re(s_1)\gg 1$, $\Re(s_2)\gg 1$, and $\Re(s_2-s_1)\gg 1$.

By Casselman-Shalika formula, the integral relative to $a\in \mathbb{A}_F^{\times}$ in the above inner integral is a holomorphic multiple of $\Lambda(1+s_2-s_1,\pi'\times\overline{\pi}'\otimes \chi_1\chi_1'^{-1}\chi_2)/\Lambda(2+2s_2,\chi_1^2\chi_1'^{-2}\chi_2^2)$. Moreover, 
by \eqref{2.3}, $\overline{h_2(x,\overline{s_2})}\propto \Lambda(2s_2,\chi_1'^{-1}\chi_2')$, and 
\begin{align*}
h_1(wuk,s_1)=\chi_1(\det wk)\int_{\mathbb{A}_F^{\times}}\Phi_1((0,t)wuk)\chi_1\chi_2^{-1}(t)|t|^{2s_1}d^{\times}t.
\end{align*}

Let $u=\begin{pmatrix}
1& b\\
&1
\end{pmatrix}$. Since $(0,t)wuk=(bt,t)wk$, then by a change of variable, 
\begin{align*}
\int_{N(\mathbb{A}_F)}h_1(wuk,s_1)du=\chi_1(\det wk)\int_{\mathbb{A}_F^{\times}}\int_{\mathbb{A}_F}\Phi_1((b,t)wk)\chi_1\chi_2^{-1}(t)|t|^{2s_1-1}dbd^{\times}t,
\end{align*}
which is a Tate integral for $\Lambda(2s_1-1,\chi_1\chi_2^{-1})$. 

Therefore, $\mathcal{M}_{\const,2}^{\du}(s_1,s_2)$ admits a meromorphic continuation $\mathcal{M}_{\const,2}^{\du,\heartsuit}(s_1,s_2)$ to $\mathbb{C}^2$, satisfying 
\begin{align*}
\mathcal{M}_{\const,2}^{\du,\heartsuit}(s_1,s_2)\propto \frac{\Lambda(1+s_2-s_1,\pi'\times\overline{\pi}'\otimes \chi_1\overline{\chi_1'}\chi_2)\Lambda(2s_1-1,\chi_1\overline{\chi_2})\Lambda(2s_2,\overline{\chi_1'}\chi_2')}{\Lambda(2+2s_2-2s_1,\chi_1^2\chi_1'^{-2}\chi_2^2)}.
\end{align*}
Hence, \eqref{eq3.5} follows from the assumption that 
$\chi_1\chi_2=\chi_1'\chi_2'$.
\end{proof}

We conclude from \eqref{3.5..}, and Lemmas \ref{lem3.6.} and \ref{lem3.7} the following.
\begin{prop}\label{prop3.4}
Let notation be as before. Then $\mathcal{M}_{\const}^{\du}(s_1,s_2)$ converges absolutely in $\big\{(s_1,s_2)\in \mathbb{C}^2:\ \Re(s_1)\gg 1,\ \Re(s_2)\gg 1\big\}$, and it admits a meromorphic continuation
\begin{align*}
\mathcal{M}_{\const}^{\du,\heartsuit}(s_1,s_2)=\mathcal{M}_{\const,1}^{\du,\heartsuit}(s_1,s_2)+\mathcal{M}_{\const,2}^{\du,\heartsuit}(s_1,s_2)
\end{align*}
to $(s_1,s_2)\in\mathbb{C}^2$, where $\mathcal{M}_{\const,1}^{\du,\heartsuit}(s_1,s_2)$ and $\mathcal{M}_{\const,2}^{\du,\heartsuit}(s_1,s_2)$ are the meromorphic functions defined in Lemma \ref{lem3.6.} and Lemma \ref{lem3.7}, respectively.  
\end{prop}

\subsection{Analytic Behaviors of $\Psi(s_2,E(\cdot,f,s),E_1(\cdot,s_1),\overline{h_2})$}
\subsubsection{A Hecke integral for $\mathrm{GL}_2\times\mathrm{GL}_1$}\label{sec3.2.1}
Let $f'\in H(\xi_1,\xi_2)$. Let $s,s_1, s_2\in \mathbb{C}$. Define (at least formally) 
\begin{align*}
\mathcal{H}_{E(\cdot,f',s)}(s_1,s_2):=\int_{\mathbb{A}_F^{\times}}W_{E(\cdot,f',s)}\left(\begin{pmatrix}
t\\
&1
\end{pmatrix};\theta\right)\chi_1'^{-1}\chi_2(t)|t|^{s_2-s_1}d^{\times}t.
\end{align*}
\begin{lemma}\label{lem3.2}
Let notation be as before. Then $\mathcal{H}_{E(\cdot,f',s)}(s_1,s_2)$ converges absolutely in the region $\Re(s_2-s_1)-|\Re(s)|>1/2$. Moreover, $\mathcal{H}_{E(\cdot,f',s)}(s_1,s_2)$ admits a meromorphic continuation to $(s,s_1,s_2)\in \mathbb{C}^3$ satisfying 
\begin{align*}
\mathcal{H}_{E(\cdot,f',s)}(s_1,s_2)\propto \frac{\Lambda(1/2+s_2-s_1+s,\chi_1'^{-1}\chi_2\xi_1)\Lambda(1/2+s_2-s_1-s,\chi_1'^{-1}\chi_2\xi_2)}{\Lambda(1+2s,\xi_1\xi_2^{-1})}.
\end{align*}
\end{lemma}
\begin{proof}
By \cite[Lemma on p.6]{JZ87},  we can write the flat section $S(f')(g,s)$ as a finite linear combination
\begin{equation}\label{eq3.10}
S(f')(g,s)=\frac{1}{\Lambda(1+2s,\xi_1\xi_2^{-1})}\sum_{j}P_j(s)h_{\Phi_j^{\circ},\xi_1,\xi_2}(g,s+1/2),
\end{equation}
where $\Phi_j^{\circ}\in \mathcal{S}(\mathbb{A}_F^2)$ is of Gauss type in the archimedean places, and each $P_j(s)$ is the reciprocal of a polynomial in $s$ and in $q_v^{-s}$ for finitely many places $v$ which has no zeros in the half-plane $\Re(s)>-1/2$. Then 
\begin{align*}
E(g,f',s)=\frac{1}{\Lambda(1+2s,\xi_1\xi_2^{-1})}\sum_{j}P_j(s)E_{\Phi_j^{\circ},\xi_1,\xi_2}(g,s+1/2).
\end{align*}
As a consequence, we have
\begin{equation}\label{3.4.}
\mathcal{H}_{E(\cdot,f',s)}(s_1,s_2)=\frac{1}{\Lambda(1+2s,\xi_1\xi_2^{-1})}\sum_{j}P_j(s)\mathcal{H}_{E_{\Phi_j^{\circ},\xi_1,\xi_2}(\cdot,s+1/2)}(s_1,s_2),
\end{equation}
where $\mathcal{H}_{E_{\Phi_j^{\circ},\xi_1,\xi_2}(\cdot,s+1/2)}(s_1,s_2)$ is defined by 
\begin{align*}
\int_{\mathbb{A}_F^{\times}}W_{E_{\Phi_j^{\circ},\xi_1,\xi_2}(\cdot,s+1/2)}\left(\begin{pmatrix}
t\\
&1
\end{pmatrix};\theta\right)\chi_1'^{-1}\chi_2(t)|t|^{s_2-s_1}d^{\times}t.
\end{align*}

By the meromorphicity of $E_{\Phi_j^{\circ},\xi_1,\xi_2}(ug,s+1/2)$, 
the Whittaker function
\begin{equation}\label{3.5.}
W_{E_{\Phi_j^{\circ},\xi_1,\xi_2}(\cdot,s+1/2)}(g;\theta):=\int_{[N]}E_{\Phi_j^{\circ},\xi_1,\xi_2}(ug,s+1/2)\overline{\theta(u)}du
\end{equation}
is a meromorphic function of $s$. When $\Re(s)>1/2$, we may plug the series expression of $E_{\Phi_j^{\circ},\xi_1,\xi_2}(g,s+1/2)$ (see \eqref{2.6.}) into \eqref{3.5.}, obtaining 
\begin{align*}
W_{E_{\Phi_j^{\circ},\xi_1,\xi_2}(\cdot,s+1/2)}(g;\theta)=\int_{N(\mathbb{A}_F)}h_{\Phi_j^{\circ},\xi_1,\xi_2}(wug,s+1/2)\overline{\theta(u)}du.
\end{align*}
Here $w=\begin{pmatrix}
&1\\
1
\end{pmatrix}$ is the Weyl element.

Substituting the definition of $h_{\Phi_j^{\circ},\xi_1,\xi_2}(g,s)$ in \eqref{2.3} into the above integral yields
\begin{equation}\label{3.4}
W_{*}(g;\theta)=\xi_1(\det g)|\det g|^{s+\frac{1}{2}}\int_{\mathbb{A}_F^{\times}}\widehat{g\Phi_j^{\circ}}(t',t'^{-1};\psi)\xi_1\xi_2^{-1}(t')|t'|^{2s}d^{\times}t',
\end{equation}
where $W_*(g;\theta):=W_{E_{\Phi_j^{\circ},\xi_1,\xi_2}(\cdot,s+1/2)}(g;\theta)$,  $\Re(s)>1/2$, and \begin{equation}\label{3.5}
\widehat{g\Phi_j^{\circ}}(t_1,t_2;\psi):=\int_{\mathbb{A}}\Phi_j^{\circ}((t_1,u)g)\overline{\psi(t_2u)}du
\end{equation}
is the partial Fourier transform relative to the second entry.

Since $\Phi_j^{\circ}\in \mathcal{S}(\mathbb{A}_F^2)$, as a function of $s$, the integral $W_{E_{\Phi_j^{\circ},\xi_1,\xi_2}(\cdot,s+1/2)}(g;\theta)$ converges absolute everywhere. Thus, \eqref{3.4} holds for all $s\in \mathbb{C}$. 

Taking advantage of \eqref{3.4},  $\mathcal{H}_{E_{\Phi_j^{\circ},\xi_1,\xi_2}(\cdot,s+1/2)}(s_1,s_2)$ simplifies to  
\begin{equation}\label{3.6}
\int_{\mathbb{A}_F^{\times}}\xi_1(t)|t|^{s+\frac{1}{2}}\int_{\mathbb{A}_F^{\times}}\widehat{\Phi_j^{\circ}}(tt',t'^{-1};\psi)\xi_1\xi_2^{-1}(t')|t'|^{2s}d^{\times}t'\chi_1'^{-1}\chi_2(t)|t|^{s_2-s_1}d^{\times}t.
\end{equation}

Changing variable $t\mapsto tt'^{-1}$ and $t'\mapsto t'^{-1}$ in \eqref{3.6} leads to 
\begin{align*}
\int_{\mathbb{A}_F^{\times}}\int_{\mathbb{A}_F^{\times}}\widehat{\Phi_j^{\circ}}(t,t';\psi)\chi_1'^{-1}\chi_2\xi_2(t')|t'|^{1/2+s_2-s_1-s}d^{\times}t'\chi_1'^{-1}\chi_2\xi_1(t)|t|^{1/2+s_2-s_1+s}d^{\times}t,
\end{align*}
which is a Tate integral for $\Lambda(1/2+s_2-s_1-s,\chi_1'^{-1}\chi_2\xi_2)\Lambda(1/2+s_2-s_1+s,\chi_1'^{-1}\chi_2\xi_1)$, converges absolutely in 
\begin{align*}
\begin{cases}
1/2+\Re(s_2-s_1-s)>1\\
1/2+\Re(s_2-s_1+s)>1,
\end{cases}
\end{align*}
which amounts to $\Re(s_2-s_1)-|\Re(s)|>1/2$. 

Hence, $\mathcal{H}_{E_{\Phi_j^{\circ},\xi_1,\xi_2}(\cdot,s+1/2)}(s_1,s_2)$ converges absolutely in $\Re(s_2-s_1)-|\Re(s)|>1/2$, and admits a meromorphic continuation to $(s,s_1,s_2)\in \mathbb{C}^3$ with 
\begin{align*}
\frac{\mathcal{H}_{E_{\Phi_j^{\circ},\xi_1,\xi_2}(\cdot,s+1/2)}(s_1,s_2)}{\Lambda(1/2+s_2-s_1-s,\chi_1'^{-1}\chi_2\xi_2)\Lambda(1/2+s_2-s_1+s,\chi_1'^{-1}\chi_2\xi_1)}
\end{align*}
being entire. Therefore, Lemma \ref{lem3.2} follows from \eqref{3.4.}.
\end{proof}

\subsubsection{The Eisenstein period $\Psi(s_2,E(\cdot,f,s),E_1(\cdot,s_1),\overline{h_2})$}
\begin{lemma}\label{lem3.3}
Let notation be as before. Suppose that $\xi_1\xi_2\chi_1\chi_2\chi_1'^{-1}\chi_2'^{-1}=\mathbf{1}$. 
Then the integral $\Psi(s_2,E(\cdot,f,s),E_1(\cdot,s_1),\overline{h_2})$ converges absolutely in 
\begin{align*}
\begin{cases}
\Re(s_2-s_1)-|\Re(s)|>1/2\\
\Re(s_2+s_1)-|\Re(s)|>3/2,
\end{cases}
\end{align*}
and it admits a meromorphic continuation to $(s,s_1,s_2)\in \mathbb{C}^3$, satisfying 
\begin{equation}\label{3.8.}
\Psi(s_2,E(\cdot,f,s),E_1(\cdot,s_1),\overline{h_2})\propto \frac{\Lambda(s_2,\pi_{\xi_1,\xi_2,s}\times\pi_{\chi_1,\chi_2,s_1-1/2}\times\chi_1'^{-1})}{\Lambda(1+2s,\xi_1\xi_2^{-1})}. 
\end{equation}
\end{lemma}
\begin{remark}
The Rankin--Selberg $L$-function  $\Lambda(s_2,\pi_{\xi_1,\xi_2,s}\times\pi_{\chi_1,\chi_2,s_1-1/2}\times\chi_1'^{-1})$ fatorizes as $\Lambda(1/2+s_2-s_1+s,\xi_1\chi_1'^{-1}\chi_2)\Lambda(1/2+s_2-s_1-s,\xi_2\chi_1'^{-1}\chi_2)\Lambda(s_1+s_2+s-1/2,\xi_1\chi_1'^{-1}\chi_1)\Lambda(s_1+s_2-s-1/2,\xi_2\chi_1'^{-1}\chi_1).
$
\end{remark}

\begin{proof}
Similar to the calculation \eqref{3.4} we have, for arbitrary $s_1\in \mathbb{C}$, that 
\begin{align*}
W_{E_1(\cdot,s_1)}(g;\overline{\theta})=|\det g|^{s_1}\int_{\mathbb{A}_F^{\times}}\widehat{g\Phi_1}(t,t^{-1};\overline{\psi})\omega\omega'(t)|t|^{2s_1-1}d^{\times}t.
\end{align*}
where $\widehat{g\Phi_1}(t,t^{-1};\overline{\psi})$ is defined analogously to \eqref{3.5}. Substituting this into the  definition in \textsection\ref{sec2.4.2} , the integral $\Psi(s_2,E(\cdot,f,s),E_1(\cdot,s_1),\overline{h_2})$ becomes
\begin{align*}
\int_{N(\mathbb{A}_F)\backslash \overline{G}(\mathbb{A}_F)}&W_{E(\cdot,f,s)}(x;\theta)\int_{\mathbb{A}_F^{\times}}\int_{\mathbb{A}_F}\Phi_1((t,ut)x)\psi(u)du\chi_1\chi_2^{-1}(t)|t|^{2s_1}d^{\times}t\\
&\int_{\mathbb{A}_F^{\times}}\overline{\Phi_2}((0,t')x)\chi_1'^{-1}\chi_2'(t')|t'|^{2s_2}d^{\times}t'\chi_1\chi_1'^{-1}(\det x)|\det x|^{s_1+s_2}dx.
\end{align*}
	
Performing the change of variables $t\mapsto tt'$ and $x\mapsto t'x$, the above integral transform to 
\begin{align*}
\int_{N(\mathbb{A}_F)\backslash G(\mathbb{A}_F)}&W_{E(\cdot,f,s)}(x;\theta)\int_{\mathbb{A}_F^{\times}}\int_{\mathbb{A}_F}\Phi_1((t,ut)x)\psi(u)du\chi_1\chi_2^{-1}(t)|t|^{2s_1}d^{\times}t\\
&\int_{\mathbb{A}_F^{\times}}\overline{\Phi_2}((0,1)x)\chi_1\chi_1'^{-1}(\det x)|\det x|^{s_1+s_2}dx,
\end{align*}
which, after the change of variable $x\mapsto \begin{pmatrix}
1& u\\
&1
\end{pmatrix}x$, simplifies to 
\begin{align*}
\int_{G(\mathbb{A}_F)}&W_{E(\cdot,f,s)}(x;\theta)\int_{\mathbb{A}_F^{\times}}\Phi_1((t,0)x)\chi_1\chi_2^{-1}(t)|t|^{2s_1}d^{\times}t\\
&\int_{\mathbb{A}_F^{\times}}\overline{\Phi_2}((0,1)x)\chi_1\chi_1'^{-1}(\det x)|\det x|^{s_1+s_2}dx.
\end{align*}
	
With a further change of variable $t\mapsto t^{-1}$ and then $x\mapsto   \begin{pmatrix}
t\\
&1
\end{pmatrix}x$, the above integral, which represents $\Psi(s_2,E(\cdot,f,s),E_1(\cdot,s_1),\overline{h_2})$, is equal to 
\begin{equation}\label{eq3.8}
\int_{G(\mathbb{A}_F)}\mathcal{H}_{E(\cdot,f,s)}(x,s_1,s_2)\Phi_1((1,0)
x)\overline{\Phi_2}((0,1)x)\overline{\chi_1'}\chi_1(\det x)|\det x|^{s_1+s_2}dx,
\end{equation}
where 
\begin{align*}
\mathcal{H}_{E(\cdot,f,s)}(x,s_1,s_2):=\int_{\mathbb{A}_F^{\times}}W_{E(\cdot,f,s)}\left(\begin{pmatrix}
t\\
&1
\end{pmatrix}x;\theta\right)\chi_1'^{-1}\chi_2(t)|t|^{s_2-s_1}d^{\times}t.
\end{align*}

Notice that for each $x$, $\mathcal{H}_{E(\cdot,f,s)}(x,s_1,s_2)$ is a Hecke integral investigated in \textsection\ref{sec3.2.1}. To separate the translation, we consider the spectral decomposition. Let $g\in G(\mathbb{A}_F)$. By \eqref{2.2} and the spectral decomposition, we have 
\begin{align*}
S(f)(wugx,s)=\sum_{f'\in \mathfrak{B}(\xi_1,\xi_2)}\langle \pi_{\xi_1,\xi_2,s}(wugx)f,f'\rangle f'(1).
\end{align*}
Along with a change of basis $f'\mapsto \pi_{\xi_1,\xi_2,s}(wug)f'$ the above identity becomes  
\begin{align*}
S(f)(wugx,s)=\sum_{f'\in \mathfrak{B}(\xi_1,\xi_2)}\langle \pi_{\xi_1,\xi_2,s}(x)f,f'\rangle \pi_{\xi_1,\xi_2,s}(wug)f'(1).
\end{align*}
Together with one more application of \eqref{2.2} we obtain 
\begin{equation}\label{3.8}
S(f)(wugx,s)=\sum_{f'\in \mathfrak{B}(\xi_1,\xi_2)}\langle \pi_{\xi_1,\xi_2,s}(x)f,f'\rangle S(f')(wug).
\end{equation}

Substituting \eqref{3.8} into the Whittaker function 
\begin{align*}
W_{E(\cdot,f,s)}(gx;\theta)=\int_{N(\mathbb{A}_F)}S(f)(wugx,s)\overline{\theta(u)}du
\end{align*}
we then derive that
\begin{align*}
W_{E(\cdot,f,s)}(gx;\theta)=\sum_{f'\in \mathfrak{B}(\xi_1,\xi_2)}\langle \pi_{\xi_1,\xi_2,s}(x)f,f'\rangle W_{E(\cdot,f',s)}(g;\theta),\ \ g\in G(\mathbb{A}_F). 
\end{align*}

Therefore, the integral \eqref{eq3.8} boils down to 
\begin{equation}\label{3.10}
\Psi(s_2,E(\cdot,f,s),E_1(\cdot,s_1),\overline{h_2})=\sum_{f'\in \mathfrak{B}(\xi_1,\xi_2)}\mathcal{H}_{E(\cdot,f',s)}(s_1,s_2)\mathcal{G}_{f'}(s,s_1,s_2),
\end{equation}
where we set $\Phi(x):=\Phi_1((1,0)
x)\overline{\Phi_2}((0,1)x)$, and 
\begin{align*}
\mathcal{G}_{f'}(s,s_1,s_2):=\int_{G(\mathbb{A}_F)}\langle \pi_{\xi_1,\xi_2,s}(x)f,f'\rangle \Phi(x)\chi_1'^{-1}\chi_1(\det x)|\det x|^{s_1+s_2}dx.
\end{align*}

The analytic behavior of $\mathcal{H}_{E(\cdot,f',s)}(s_1,s_2)$ has been investigated in Lemma \ref{lem3.2}. Moving forward, we shift our focus to $\mathcal{G}_{f'}(s,s_1,s_2)$. By \eqref{2.2} we obtain 
\begin{align*}
\mathcal{G}_{f'}(s,s_1,s_2)=\int_{G(\mathbb{A}_F)}\int_K S(f)(kx,s)\overline{f'(k)}dk \Phi(x)\chi_1'^{-1}\chi_1(\det x)|\det x|^{s_1+s_2}dx.
\end{align*}

Utilizing the Iwasawa coordinates $x=\begin{pmatrix}
a& b\\
&c
\end{pmatrix}k'$ and the definition of $S(f)(\cdot,s)$ in \eqref{eq2.2}, the integral $\mathcal{G}_{f'}(s,s_1,s_2)$ becomes  
\begin{align*}
&\int_{\mathbb{A}_F^{\times}}\int_{\mathbb{A}_F^{\times}}\int_{\mathbb{A}_F}\int_{K}\int_{K} \Psi\left(k^{-1}\begin{pmatrix}
a& b\\
&c\end{pmatrix}k'\right)f(k)\overline{f'(k)}\chi_1'^{-1}\chi_1(\det k^{-1}k')dkdk'db\\
&\qquad\qquad|a|^{s_1+s_2+s-1/2}\xi_1(a)\chi_1'^{-1}\chi_1(a)|c|^{s_1+s_2-s-1/2}\xi_2(c)\chi_1'^{-1}\chi_1(c)d^{\times}ad^{\times}c,
\end{align*}
which is a Tate integral for 
\begin{align*}
\Lambda(s_1+s_2+s-1/2,\xi_1\chi_1'^{-1}\chi_1)\Lambda(s_1+s_2-s-1/2,\xi_2\chi_1'^{-1}\chi_1).
\end{align*}
Hence, $\mathcal{G}_{f'}(s,s_1,s_2)$ converges absolutely in 
\begin{align*}
\begin{cases}
\Re(s_1+s_2+s-1/2)>1\\
\Re(s_1+s_2-s-1/2)>1,
\end{cases}
\end{align*}
which is equivalent to $\Re(s_1+s_2)-|\Re(s)|>3/2$. Moreover, $\mathcal{G}_{f'}(s,s_1,s_2)$ admits a meromorphic continuation to $(s,s_1,s_2)\in \mathbb{C}^3$, satisfying 
\begin{align*}
\mathcal{G}_{f'}(s,s_1,s_2)\propto \Lambda(s_1+s_2+s-1/2,\xi_1\chi_1'^{-1}\chi_1)\Lambda(s_1+s_2-s-1/2,\xi_2\chi_1'^{-1}\chi_1).
\end{align*}

Therefore, Lemma \ref{lem3.3} follows from Lemma \ref{lem3.2} and \eqref{3.10}.
\end{proof}
\begin{remark}
We adopt the notation $\mathcal{G}$ in \eqref{3.10} because $\mathcal{G}_{f'}(s,s_1,s_2)$ emerges as a Godement integral (for Eisenstein series).	
\end{remark}

\subsection{Zero-free Regions of Hecke $L$-functions}\label{sec3.3}
Let $\chi$ be a unitary Hecke character over $F$. There is an effectively computable absolute constant $c>0$ and a constant $c_{\chi}>0$ depending on $\chi$ (and $F$), such that Hecke $L$-function $L(s, \chi)$ has no zero in the region 
$$\big\{s=\sigma+it:\ \sigma> 1-\frac{c}{c_{\chi}+\log (t^2+1)}\big\}.$$ 
Then $\Lambda(1+2s,\chi)\neq 0$ in $s\in \mathcal{D}_{\chi}(\varepsilon)$, where $0<\varepsilon\leq 1$ and 
\begin{equation}\label{eq3.12}
\mathcal{D}_{\chi}(\varepsilon):=\Big\{s=\sigma+it:\ -\frac{c\varepsilon}{10[c_{\chi}+\log (t^2+1)]}< \sigma< \frac{c\varepsilon}{10[c_{\chi}+\log (t^2+1)]}\Big\}.
\end{equation}    
Upon replacing $c$ with a smaller constant, we may assume that if $s\in\mathcal{D}_{\chi}(1)$, then $|\Re(s)|<10^{-1}$.  

Let $\mathcal{D}_{\chi}^+(\varepsilon):=\mathcal{D}_{\chi}(\varepsilon)\cap \{s\in \mathbb{C}:\ \Re(s)> 0\}$, and $\mathcal{D}_{\chi}^-(\varepsilon):=\mathcal{D}_{\chi}(\varepsilon)\cap \{s\in \mathbb{C}:\ \Re(s)< 0\}$. We denote the boundaries of $\mathcal{D}_{\chi}(\varepsilon)$ by 
\begin{align*}
&\mathcal{C}_{\chi}^+(\varepsilon):=\Big\{s=\sigma+it:\ \sigma=\frac{c}{4[c_{\chi}+\log (t^2+1)]}\Big\},\\
&\mathcal{C}_{\chi}^-(\varepsilon):=\Big\{s=\sigma+it:\ \sigma=-\frac{c}{4[c_{\chi}+\log (t^2+1)]}\Big\}.
\end{align*}
For simplicity, we write $\mathcal{C}_{\chi}^{+}$ (resp. $\mathcal{C}_{\chi}^{-}$) for $\mathcal{C}_{\chi}^{+}(1)$ (resp. $\mathcal{C}_{\chi}^{-}(1)$).

\subsection{Meromorphic Continuation of $\mathcal{M}_{\Eis}^{\du}(s_1,s_2)$}

Let $\Re(s_2-s_1)>1/2$ and $
\Re(s_2+s_1)>3/2$. Let $\xi\in \widehat{F^{\times}\backslash\mathbb{A}_F^{(1)}}$. Define (at least formally)
\begin{align*}
\mathcal{M}_{\Eis}^{\du}(s_1,s_2;\xi):= \frac{1}{4\pi i}\int_{i\mathbb{R}}\sum_{f} \langle \phi_1\overline{\phi_2},E(\cdot,f,s)\rangle \Psi(s_2,E(\cdot,f,s),E_1(\cdot,s_1),\overline{h_2})ds,
\end{align*}
where $f\in \mathfrak{B}(\xi,\xi^{-1})$. Noticing that $s=-\overline{s}$ for $s\in i\mathbb{R}$, we have 
\begin{equation}\label{3.12}
\mathcal{M}_{\Eis}^{\du}(s_1,s_2;\xi)=\frac{1}{4\pi i}\int_{i\mathbb{R}}\sum_{f\in \mathfrak{B}(\xi,\xi^{-1})} \langle \phi_1\overline{\phi_2},E(\cdot,f,-\overline{s})\rangle \Psi(\cdots)ds,
\end{equation}
where $\Psi(\cdots):=\Psi(s_2,E(\cdot,f,s),E_1(\cdot,s_1),\overline{h_2})$.

Due to the rapid decay of $\phi_1\overline{\phi_2}$, the inner product $\langle \phi_1\overline{\phi_2},E(\cdot,f,-\overline{s})\rangle$ converges absolutely for $s\in \mathbb{C}$, defining an entire function. Since $\phi_1$ and $\phi_2$ are $K$-finite, $\mathfrak{B}(\xi,\xi^{-1})$ is a finite set. With Lemma \ref{lem3.3} and the fact that Rankin--Selberg periods decreasing rapidly in a fixed  vertical strip, the function $\mathcal{M}_{\Eis}^{\du}(s_1,s_2;\xi)$ converges absolutely in the region $\Re(s_2-s_1)>1/2$ and $
\Re(s_2+s_1)>3/2$. In particular, we may consider 
\begin{align*}
\mathcal{M}_{\Eis}^{\du}(s_1,s_2;\xi)=\frac{1}{2\pi i}\int_{i\mathbb{R}} \mathfrak{F}_{\xi}(s,s_1,s_2)ds,
\end{align*}
where 
\begin{equation}\label{f}
\mathfrak{F}_{\xi}(s,s_1,s_2):=\frac{1}{2}\sum_{f\in \mathfrak{B}(\xi,\xi^{-1})}\langle \phi_1\overline{\phi_2},E(\cdot,f,-\overline{s})\rangle \Psi(s_2,E(\cdot,f,s),E_1(\cdot,s_1),\overline{h_2}).
\end{equation}

\begin{lemma}\label{lem3.6}
Let notation be as before.  Let $f\in \mathfrak{B}(\xi_1,\xi_2)$, where $\xi_1, \xi_2\in \widehat{F^{\times}\backslash\mathbb{A}_F^{(1)}}$. Then the function $\langle \phi_1\overline{\phi_2},E(\cdot,f,-\overline{s})\rangle \Psi(s_2,E(\cdot,f,s),E_1(\cdot,s_1),\overline{h_2})
$ admits a meromorphic continuation to $\mathbb{C}^3$, and is a holomorphic multiple of  
\begin{equation}\label{3.17}
\frac{\Lambda(1/2-s,\pi'\times\overline{\pi}'\otimes \overline{\xi}_1)\Lambda(s_2,\pi_{\xi_1,\xi_2,s}\times\pi_{\chi_1,\chi_2,-1/2}\times\chi_1'^{-1})}{\Lambda(1-2s,\xi_1^{-1}\xi_2)\Lambda(1+2s,\xi_1\xi_2^{-1})}.
\end{equation}
\end{lemma}
\begin{proof}
Let $\Re(s)\leq -2$, By unfolding and the Fourier expansion of cusp forms, 
\begin{align*}
\langle \phi_1\overline{\phi_2},E(\cdot,f,-\overline{s})\rangle=\int_{N(\mathbb{A}_F)Z(\mathbb{A}_F)\backslash G(\mathbb{A}_F)}W_{\phi_1}(x;\theta)\overline{W_{\phi_2}(x;\theta)}\overline{S(f)(x,-\overline{s})}dx.
\end{align*}

Using the Iwasawa decomposition and \eqref{eq3.10} we can write $\langle \phi_1\overline{\phi_2},E(\cdot,f,-\overline{s})\rangle$ as 
\begin{align*}
\sum_{j}\frac{\overline{P_j(-\overline{s})}}{\Lambda(1-2s,\xi_1^{-1}\xi_2)}\int_{N(\mathbb{A}_F)Z(\mathbb{A}_F)\backslash G(\mathbb{A}_F)}W_{\phi_1}(x;\theta)\overline{W_{\phi_2}(x;\theta)}\overline{h_{\Phi_j^{\circ},\xi_1,\xi_2}(x,-\overline{s}+1/2)}dx
\end{align*}
which converges absolutely in $\Re(s)\leq -2$, and is the Rankin--Selberg integral for $\Lambda(1/2-s,\pi'\times\overline{\pi}'\otimes \overline{\xi}_1)/\Lambda(1-2s,\xi_1^{-1}\xi_2)$. 

Hence, $\langle \phi_1\overline{\phi_2},E(\cdot,f,-\overline{s})\rangle$ admits a meromorphic continuation to $s\in \mathbb{C}$, with
\begin{equation}\label{3.13}
\langle \phi_1\overline{\phi_2},E(\cdot,f,-\overline{s})\rangle \propto \frac{\Lambda(1/2-s,\pi'\times\overline{\pi}'\otimes \overline{\xi}_1)}{\Lambda(1-2s,\xi_1^{-1}\xi_2)}.
\end{equation}

Then Lemma \ref{lem3.6} follows from Lemma \ref{lem3.3}, and \eqref{3.17} is a consequence of \eqref{3.8.} and \eqref{3.13}.\end{proof}

\subsubsection{Meromorphic continuation of $\mathcal{M}_{\Eis}^{\du}(s_1,s_2;\xi)$} 
Define the region 
\begin{equation}\label{3.21.}
\mathcal{D}:=\big\{(s_1,s_2)\in \mathbb{C}: \ -1/2<\Re(s_2-s_1)<1/2,\ 1/2<\Re(s_2+s_1)<3/2\big\}.
\end{equation}

By Lemma \ref{lem3.6}, the function $\mathfrak{F}_{\xi}(s,s_1,s_2)$ admits a 
meromorphic continuation to $(s,s_1,s_2)\in\mathbb{C}^3$. 
For simplicity, we continue to write $\mathfrak{F}_{\xi}(s,s_1,s_2)$ for its 
meromorphic continuation. Let $(s_1,s_2)\in \mathcal{D}$. We define 
\begin{multline}\label{3.27}
\mathcal{M}_{\Res}^{\du,\heartsuit}(s_1,s_2):=\underset{s=s_1-s_2+\frac{1}{2}}{\Res}\  \mathfrak{F}_{\chi_1\overline{\chi_2}}(s,s_1,s_2)-\underset{s=s_2-s_1-\frac{1}{2}}{\Res}\  \mathfrak{F}_{\overline{\chi_1}\chi_2}(s,s_1,s_2)\\
\underset{s=\frac{3}{2}-s_1-s_2}{\Res}\  \mathfrak{F}_{\mathbf{1}}(s,s_1,s_2)-\underset{s=s_1+s_2-\frac{3}{2}}{\Res}\  \mathfrak{F}_{\mathbf{1}}(s,s_1,s_2).
\end{multline}

\begin{prop}\label{prop3.7}
Let notation be as before. Assume $\chi_1=\chi_1'$ and $\chi_2=\chi_2'$. The function $\mathcal{M}_{\Eis}^{\du}(s_1,s_2;\xi)$ admits a meromorphic continuation $\widetilde{\mathcal{M}}_{\Eis}^{\du}(s_1,s_2;\xi)$ to $\mathcal{D}$. 
In particular, for $(s_1,s_2)\in \mathcal{D}$, we have
\begin{align*}
\widetilde{\mathcal{M}}_{\Eis}^{\du}(s_1,s_2;\xi)=\frac{1}{2\pi i}\int_{i\mathbb{R}}\mathfrak{F}_{\xi}(s,s_1,s_2)ds+\mathcal{M}_{\Res}^{\du,\heartsuit}(s_1,s_2;\xi).
\end{align*}
\end{prop}
\begin{proof}
Let $\chi=\xi^2$. Let $s\in \mathcal{D}_{\chi}$, which is defined in \eqref{eq3.12}. By Lemma \ref{lem3.6} and the global functional equation, 
\begin{multline*}
\mathfrak{F}_{\xi}(s,s_1,s_2)\propto \Lambda(1-2s,\xi^{-2})^{-1}\Lambda(1+2s,\xi^2)^{-1}\Lambda(1/2+s,\pi'\times\overline{\pi}'\otimes \xi)\\
\Lambda(1/2+s_2-s_1+s,\overline{\chi_1}\chi_2\xi)\Lambda(1/2+s_2-s_1-s,\overline{\chi_1}\chi_2\overline{\xi})\\
\Lambda(s_1+s_2+s-1/2,\xi)\Lambda(s_1+s_2-s-1/2,\overline{\xi}).
\end{multline*}

We will make use of this analytic behavior to produce the meromorphic continuation of the function $\big\{(s_1,s_2)\in \mathbb{C}^2:\ \Re(s_1)<\min\{1/2+\Re(s_2), 1/2\}\big\}$. 

Recall that the function $\mathcal{M}_{\Eis}^{\du}(s_1,s_2;\xi)$ converges absolutely in the region $\Re(s_2-s_1)>1/2$ and $
\Re(s_2+s_1)>3/2$. To start with, we take $s_1\in \mathbb{C}$ with $\Re(s_1)>1/2$. Then $\mathcal{M}_{\Eis}^{\du}(s_1,s_2;\xi)$ converges absolutely in the region $
\Re(s_2)>1/2+\Re(s_1)$. We shall do the meromorphic continuation strip by strip from the right to the left of the $x$-axis. 

\begin{itemize}
\item Let $0<\varepsilon<10^{-3}$. Suppose $\Re(s_1)>1/2+\varepsilon$. Let $\mathcal{D}_{\chi}^+(\varepsilon)$ and $\mathcal{C}_{\chi}^{+}(\varepsilon)$ be defined as in \textsection\ref{sec3.3}. Consider $s_2\in 1/2+s_1+\mathcal{D}_{\chi}^+(\varepsilon)$. By Lemma \ref{lem3.6} we may shift the contour and apply Cauchy's formula to obtain 
\begin{equation}\label{3.21}
\frac{1}{2\pi i}\int_{i\mathbb{R}}\mathfrak{F}_{\xi}(s,s_1,s_2)ds=\frac{1}{2\pi i}\int_{\mathcal{C}_{\chi}^{+}(\varepsilon)}\mathfrak{F}_{\xi}(s,s_1,s_2)ds-\underset{s=s_2-s_1-\frac{1}{2}}{\Res}\  \mathfrak{F}_{\xi}(s,s_1,s_2).
\end{equation}
By Lemma \ref{lem3.6} we have $\underset{s=s_2-s_1-1/2}{\Res}\  \mathfrak{F}_{\xi}(s,s_1,s_2)\equiv 0$ unless $\xi=\overline{\chi_1}\chi_2$. Assume $\xi=\overline{\chi_1}\chi_2$. Utilizing \eqref{3.17} we obtain   
\begin{align*}
\underset{s=s_2-s_1-\frac{1}{2}}{\Res}\  \mathfrak{F}_{\xi}(s,s_1,s_2)\propto \frac{\Lambda(s_2-s_1,\pi'\times\overline{\pi}'\otimes \overline{\chi_1}\chi_2))\Lambda(2s_2-1,\overline{\chi_1}\chi_2)\Lambda(2s_1,\chi_1\overline{\chi_2})}{\Lambda(2-2s_2+2s_1,\chi_1^2\overline{\chi_2}^2)}.
\end{align*}

Notice that $\frac{1}{2\pi i}\int_{\mathcal{C}_{\chi}^{+}(\varepsilon)}\mathfrak{F}_{\xi}(s,s_1,s_2)ds$ is holomorphic in 
\begin{align*}
\mathcal{D}_1:=\big\{(s_1,s_2)\in \mathbb{C}: \ s_2\in 1/2+s_1+\mathcal{D}_{\chi}(\varepsilon),\ \Re(s_1)>1/2+\varepsilon\big\} 
\end{align*}
and $\underset{s=s_2-s_1-1/2}{\Res}\  \mathfrak{F}_{\xi}(s,s_1,s_2)$ is meromorphic in $\mathbb{C}$. We then obtain a meromorphic continuation $\widetilde{\mathcal{M}}_{\Eis}^{\du}(s_1,s_2;\xi)$ of $\mathcal{M}_{\Eis}^{\du}(s_1,s_2;\xi)$ to the region   
\begin{align*}
\mathcal{D}_2:=\mathcal{D}_1\bigcup \big\{(s_1,s_2)\in \mathbb{C}: \ \Re(s_2)\geq  1/2+\Re(s_1),\ \Re(s_1)>1/2+\varepsilon\big\}. 
\end{align*}

\item Consider $s_2\in 1/2+s_1+\mathcal{D}_{\chi}^{-}(\varepsilon)$ and $\Re(s_1)>1/2+\varepsilon$. Shifting contour and apply Cauchy's formula according to Lemma \ref{lem3.6}, we derive  
\begin{equation}\label{3.22}
\frac{1}{2\pi i}\int_{\mathcal{C}_{\chi}^{+}(\varepsilon)}\mathfrak{F}_{\xi}(s,s_1,s_2)ds=\frac{1}{2\pi i}\int_{i\mathbb{R}}\mathfrak{F}_{\xi}(s,s_1,s_2)ds+\underset{s=s_1-s_2+\frac{1}{2}}{\Res}\  \mathfrak{F}_{\xi}(s,s_1,s_2),
\end{equation}
where $\underset{s=s_1-s_2+1/2}{\Res}\  \mathfrak{F}_{\xi}(s,s_1,s_2)\equiv 0$ unless $\xi=\chi_1\overline{\chi_2}$. 

Assume $\xi=\chi_1\overline{\chi_2}$. Then $\underset{s=s_1-s_2+1/2}{\Res}\  \mathfrak{F}_{\xi}(s,s_1,s_2)$ is a holomorphic multiple of the meromorphic function 
\begin{align*}
\frac{\Lambda(1+s_1-s_2,\pi'\times\overline{\pi}'\otimes \chi_1\overline{\chi_2})\Lambda(2s_1,\chi_1\overline{\chi_2})\Lambda(2s_2-1,\overline{\chi_1}\chi_2)}{\Lambda(2+2s_1-2s_2,\chi_1^2\overline{\chi_2}^2)}.
\end{align*}

Notice that $\frac{1}{2\pi i}\int_{i\mathbb{R}}\mathfrak{F}_{\xi}(s,s_1,s_2)ds$ is holomorphic in the region 
\begin{align*}
\mathcal{D}_3:=\big\{(s_1,s_2)\in \mathbb{C}: \ 3/2-\Re(s_1)<\Re(s_2)<1/2+\Re(s_1),\ \Re(s_1)>1/2+\varepsilon\big\}.
\end{align*}
This leads to the meromorphic continuation of $\widetilde{\mathcal{M}}_{\Eis}^{\du}(s_1,s_2;\xi)$ to $\mathcal{D}_3$.

\item Let $s_2\in 3/2-s_1+\mathcal{D}_{\chi}^{+}$ and $\Re(s_1)>1/2+\varepsilon$. We can shift the contour, utilizing \eqref{3.17}, to derive  
\begin{equation}\label{3.18}
\frac{1}{2\pi i}\int_{i\mathbb{R}}\mathfrak{F}_{\xi}(s,s_1,s_2)ds=\frac{1}{2\pi i}\int_{\mathcal{C}_{\chi}^{+}}\mathfrak{F}_{\xi}(s,s_1,s_2)ds-\underset{s=s_1+s_2-3/2}{\Res}\  \mathfrak{F}_{\xi}(s,s_1,s_2).
\end{equation}

By the analytic behavior of $\mathfrak{F}_{\xi}(s,s_1,s_2)$ (see Lemma \ref{lem3.6}), we have
\begin{align*}
\underset{s=s_1+s_2-3/2}{\Res}\  \mathfrak{F}_{\xi}(s,s_1,s_2)\equiv 0
\end{align*}
unless $\xi=\mathbf{1}$. With $\xi=\mathbf{1}$, it follows from \eqref{3.17} that 
\begin{align*}
\underset{s=s_1+s_2-3/2}{\Res}\  \mathfrak{F}_{\xi}(s,s_1,s_2)\propto \frac{\Lambda(s_1+s_2-1,\pi'\times\overline{\pi}')\Lambda(2s_2-1,\overline{\chi_1}\chi_2)\Lambda(2-2s_1,\overline{\chi_1}\chi_2)}{\Lambda(4-2s_1-2s_2)}.
\end{align*}

We notice that the function $\frac{1}{2\pi i}\int_{\mathcal{C}_{\chi}^{+}}\mathfrak{F}_{\xi}(s,s_1,s_2)ds$ 
is holomorphic in 
\begin{align*}
\big\{(s_1,s_2)\in \mathbb{C}: \ s_2\in 3/2-s_1+\mathcal{D}_{\chi},\ \Re(s_1)>1/2+\varepsilon\big\}.
\end{align*}

\item Let $s_2\in 3/2-s_1+\mathcal{D}_{\chi}^-$ and $\Re(s_1)>1/2+\varepsilon$. By \eqref{3.17} and contour shifting, 
\begin{equation}\label{3.19}
\frac{1}{2\pi i}\int_{\mathcal{C}_{\chi}^{+}}\mathfrak{F}_{\xi}(s,s_1,s_2)ds=\frac{1}{4\pi i}\int_{i\mathbb{R}}\mathfrak{F}_{\xi}(s,s_1,s_2)ds+\underset{s=3/2-s_1-s_2}{\Res}\  \mathfrak{F}_{\xi}(s,s_1,s_2).
\end{equation}

Similar to the previous analysis, $\underset{s=3/2-s_1-s_2}{\Res}\  \mathfrak{F}_{\xi}(s,s_1,s_2)\equiv 0$ unless $\xi=\mathbf{1}$. Under the assumption that $\xi=\mathbf{1}$, we have 
\begin{align*}
\underset{s=3/2-s_1-s_2}{\Res}\  \mathfrak{F}_{\xi}(s,s_1,s_2)\propto \frac{\Lambda(2-s_1-s_2,\pi'\times\overline{\pi}')\Lambda(2-2s_1,\overline{\chi_1}\chi_2)\Lambda(2s_2-1,\overline{\chi_1}\chi_2)}{\Lambda(4-2s_1-2s_2)}.
\end{align*}

Gathering together  \eqref{3.21}, \eqref{3.22}, \eqref{3.18} and \eqref{3.19}, we have extended $\mathcal{M}_{\Eis}^{\du}(s_1,s_2;\xi)$ to a meromorphic function $\widetilde{\mathcal{M}}_{\Eis}^{\du}(s_1,s_2;\xi)$ to the region 
\begin{align*}
\big\{(s_1,s_2)\in\mathbb{C}^2:\ s_2+s_1\in 3/2+\mathcal{D}_{\chi},\ \Re(s_1)>1/2+\varepsilon\big\}.
\end{align*} 
In particular, when $s_2\in 3/2-s_1+\mathcal{D}_{\chi}^-$ and $\Re(s_1)>1/2+\varepsilon$, the function $\widetilde{\mathcal{M}}_{\Eis}^{\du}(s_1,s_2;\xi)$ is explicitly defined by 
\begin{align*}
&\frac{1}{2\pi i}\int_{i\mathbb{R}}\mathfrak{F}_{\xi}(s,s_1,s_2)ds-\mathbf{1}_{\xi=\mathbf{1}}\underset{s=s_1+s_2-\frac{3}{2}}{\Res}\  \mathfrak{F}_{\xi}(s,s_1,s_2)+\mathbf{1}_{\xi=\mathbf{1}}\underset{s=\frac{3}{2}-s_1-s_2}{\Res}\  \mathfrak{F}_{\xi}(s,s_1,s_2)\\
&\qquad -\mathbf{1}_{\xi=\overline{\chi_1}\chi_2}\underset{s=s_2-s_1-1/2}{\Res}\  \mathfrak{F}_{\xi}(s,s_1,s_2)+\mathbf{1}_{\xi=\chi_1\overline{\chi_2}}\underset{s=s_1-s_2+1/2}{\Res}\  \mathfrak{F}_{\xi}(s,s_1,s_2).
\end{align*}

Notice that the function $\frac{1}{2\pi i}\int_{i\mathbb{R}}\mathfrak{F}_{\xi}(s,s_1,s_2)ds$ is holomorphic in  
\begin{align*}
\mathcal{D}=\big\{(s_1,s_2)\in \mathbb{C}: \ -1/2<\Re(s_2-s_1)<1/2,\ 1/2<\Re(s_2+s_1)<3/2\big\}.
\end{align*}
Therefore, the function $\widetilde{\mathcal{M}}_{\Eis}^{\du}(s_1,s_2;\xi)$ further extends to a meromorphic function in $\mathcal{D}$, which is defined by \eqref{3.21.}.
\end{itemize}

As a consequence of the above discussion, Proposition \ref{prop3.7} holds. 
\end{proof}
\begin{remark}
With further effort we may derive a meromorphic continuation of $\mathcal{M}_{\Eis}^{\du}(s_1,s_2;\xi)$ to $(s_1,s_2)\in \mathbb{C}^2$. However, for the purpose we only need to focus on the continuation to a small neighborhood of $(1/2,1/2)$. Note that when $\Re(s_1)< 1/2$, we need to start with $\Re(s_2)>3/2-\Re(s_1)$  instead of $\Re(s_2)>1/2+\Re(s_1)$ in the process of contour shifting. 
\end{remark}

\begin{remark}
Each $\underset{s=*}{\Res}\  \mathfrak{F}_{\xi}(s,s_1,s_2)$, as a meromorphic function, may have a pole at $(s_1,s_2)=(1/2,1/2)$ of order at most $3$. However, the singularity of the linear combination of them should be matched with that of $\widetilde{\mathcal{M}}_{\Eis}^{\du}(s_1,s_2)$ near $(1/2,1/2)$. 
\end{remark}

\subsubsection{Meromorphic continuation of $\mathcal{M}_{\Eis}^{\du}(s_1,s_2)$} 
Let $(s_1,s_2)\in\mathcal{D}$. We define
\begin{align*}
\mathcal{M}_{\Eis}^{\du,\heartsuit}(s_1,s_2):=\sum_{\xi\in \widehat{F^{\times}\backslash\mathbb{A}_F^{(1)}}}\frac{1}{2\pi i}\int_{i\mathbb{R}}\mathfrak{F}_{\xi}(s,s_1,s_2)ds.
\end{align*}

\begin{prop}\label{prop3.10}
Let notation be as before. Suppose that $\chi_1=\chi_1'$ and $\chi_2=\chi_2'$. Then $\mathcal{M}_{\Eis}^{\du}(s_1,s_2)$ admits a meromorphic continuation $\widetilde{\mathcal{M}}_{\Eis}^{\du}(s_1,s_2)$ to $\mathcal{D}$. Moreover, for $(s_1,s_2)\in\mathcal{D}$, we have
\begin{align*}
\widetilde{\mathcal{M}}_{\Eis}^{\du}(s_1,s_2)=\mathcal{M}_{\Eis}^{\du,\heartsuit}(s_1,s_2)+\mathcal{M}_{\Res}^{\du,\heartsuit}(s_1,s_2).
\end{align*}
Moreover, we have the following analytic behaviors:
\begin{enumerate}
\item[(a).]  The function $\underset{s=s_2-s_1-1/2}{\Res}\  \mathfrak{F}_{\overline{\chi_1}\chi_2}(s,s_1,s_2)$ is a holomorphic multiple of 
\begin{align*}
\frac{\Lambda(s_2-s_1,\pi'\times\overline{\pi}'\otimes \overline{\chi_1}\chi_2))\Lambda(2s_2-1,\overline{\chi_1}\chi_2)\Lambda(2s_1,\chi_1\overline{\chi_2})}{\Lambda(2-2s_2+2s_1,\chi_1^2\overline{\chi_2}^2)};
\end{align*}

\item[(b).]   The function $\underset{s=s_1-s_2+1/2}{\Res}\  \mathfrak{F}_{\chi_1\overline{\chi_2}}(s,s_1,s_2)$ is a holomorphic multiple of 
\begin{align*}
\frac{\Lambda(1+s_1-s_2,\pi'\times\overline{\pi}'\otimes \chi_1\overline{\chi_2})\Lambda(2s_1,\chi_1\overline{\chi_2})\Lambda(2s_2-1,\overline{\chi_1}\chi_2)}{\Lambda(2+2s_1-2s_2,\chi_1^2\overline{\chi_2}^2)};
\end{align*}
\item[(c).]   The function $\underset{s=s_1+s_2-3/2}{\Res}\  \mathfrak{F}_{\mathbf{1}}(s,s_1,s_2)$ is a holomorphic multiple of 
\begin{align*}
\frac{\Lambda(s_1+s_2-1,\pi'\times\overline{\pi}')\Lambda(2s_2-1,\overline{\chi_1}\chi_2)\Lambda(2-2s_1,\overline{\chi_1}\chi_2)}{\Lambda(4-2s_1-2s_2)};
\end{align*}

\item[(d).]   The function $\underset{s=3/2-s_1-s_2}{\Res}\  \mathfrak{F}_{\mathbf{1}}(s,s_1,s_2)$ is a holomorphic multiple of 
\begin{align*}
\frac{\Lambda(2-s_1-s_2,\pi'\times\overline{\pi}')\Lambda(2-2s_1,\overline{\chi_1}\chi_2)\Lambda(2s_2-1,\overline{\chi_1}\chi_2)}{\Lambda(4-2s_1-2s_2)}.
\end{align*}
\end{enumerate}
\end{prop}
\begin{proof}

By Lemma \ref{lem3.6}, we have 
\begin{align*}
\mathfrak{F}_{\xi}(s,s_1,s_2)\propto \frac{\Lambda(1/2+s,\pi'\times\overline{\pi}'\otimes \xi)\Lambda(s_2,\pi_{\xi,\xi^{-1},s}\times\pi_{\chi_1,\chi_2,-1/2}\times\chi_1'^{-1})}{\Lambda(1-2s,\xi^{-2})\Lambda(1+2s,\xi^2)}.
\end{align*}

Therefore, the analytic behaviors (a)--(d) follows. 
\end{proof}

\subsection{The Regularized Spectral Reciprocity}\label{sec3.5} 
Since $\phi_1$ and $\phi_2$ decay rapidly in $[G]$, the function $\langle \phi_1E_1(\cdot,s_1),\phi_2E_2(\cdot,\overline{s_2})\rangle$ converges absolutely for all $(s_1,s_2)\in \mathbb{C}$. Hence, $\langle \phi_1E_1(\cdot,s_1),\phi_2E_2(\cdot,\overline{s_2})\rangle$ extends to an entire function in $\mathbb{C}^2$. 

\begin{thm}\label{thmA}
Let notation be as before. Suppose that $\chi_1=\chi_1'$ and $\chi_2=\chi_2'$. Then 
\begin{multline*}
\mathcal{M}_{\cusp}^{\tw}(s_1,s_2)+\mathcal{M}_{\Eis}^{\tw}(s_1,s_2)=\mathcal{M}_{\const}^{\du,\heartsuit}(s_1,s_2)+\mathcal{M}_{\cusp}^{\du}(s_1,s_2)\\
+\mathcal{M}_{\Eis}^{\du,\heartsuit}(s_1,s_2)+\mathcal{M}_{\Res}^{\du,\heartsuit}(s_1,s_2)
\end{multline*}
holds as an identity of meromorphic functions in 
\begin{align*}
\mathcal{D}=\big\{(s_1,s_2)\in \mathbb{C}: \ -1/2<\Re(s_2-s_1)<1/2,\ 1/2<\Re(s_2+s_1)<3/2\big\}.
\end{align*}
Here the meromorphic functions $\mathcal{M}_{\const}^{\du,\heartsuit}(s_1,s_2)$,  $\mathcal{M}_{\Eis}^{\du,\heartsuit}(s_1,s_2)$ and $\mathcal{M}_{\Res}^{\du,\heartsuit}(s_1,s_2)$ are defined in Proposition \ref{prop3.4} and Proposition \ref{prop3.10}, respectively. 
\end{thm}

%\subsection{The Amplified Regularized Spectral Reciprocity}

\section{Explicit Choice of the Automorphic Data}\label{sec4}
Let $\omega=\otimes_v'\omega_v$ and $\omega'=\otimes_v'\omega_v'$ be unitary Hecke characters on $F^{\times}\backslash \mathbb{A}_F^{\times}$, which are trivial on $\mathbb{R}_+^{\times}$. Let $\pi=\otimes_v'\pi_v\in \mathcal{A}([G],\omega)$ and $\pi'=\otimes_v'\pi_v'\in \mathcal{A}_0([G],\omega')$. 
\subsection{Conductors}\label{sec4.1}
Let $C(\pi):=\otimes_vC_v(\pi)$ be the analytic conductor of $\pi$, where each local conductor $C_v(\pi)$ is defined as follows.
\begin{itemize}
\item For $v<\infty$, we denote by  $r_{\pi_v}\geq 0$ the conductor exponent of $\pi_v$, which is the unique integer such that $\pi_v$ has a vector which is $K_{v}[r_{\pi_v}]$-invariant but is not $K_{v}[r_{\pi_v}-1]$-invariant. Here $K_v[n]:=\Big\{\begin{pmatrix}
k_{11}& k_{12}\\
k_{21}& k_{22}	
\end{pmatrix}\in K_v:\ k_{21}\in \varpi_v^n\mathcal{O}_v
\Big\}$, $n\in \mathbb{Z}$. Let $C_{v}(\pi):=q_v^{r_{\pi_v}}$ be the local conductor of $\pi_v$. 

\item For $v\mid\infty$, the local $L$-function of $\pi_v$ is a product of shifted Gamma factors of the shape $\ell_v(s,\pi_v)=\Gamma_v(s+\beta_{1,v})\Gamma_v(s+\beta_{2,v})$. Here $\Gamma_v$ is the Gamma function over $F_v$, and $\beta_{1,v}$ and $\beta_{2,v}$ are complex numbers. Let
\begin{equation}\label{cond}
C_v(\pi):=[(1+|\beta_{1,v}|)(1+|\beta_{2,v}|)]^{[F_v:\mathbb{R}]}.
\end{equation}
\end{itemize}  

Let $\chi_v$ be a character of $F_v^{\times}$, where $v<\infty$. We define $r_{\chi_v}\geq 0$ to be the least nonnegative integer such that $\chi_v$ is invariant on $1+\mathfrak{p}_v^{r_{\chi_v}}$ but is not invariant on $1+\mathfrak{p}_v^{r_{\chi_v}-1}$. This integer $r_{\chi_v}$ is called the conductor exponent of $\chi_v$.   

Let $\mathfrak{N}:=\prod_{v<\infty}\mathfrak{p}_v^{r_{\pi_v}}$ and $\mathfrak{M}:=\prod_{v<\infty}\mathfrak{p}_v^{r_{\pi_v'}}$. Let $m_v:=\max\{r_{\pi_v},r_{\pi_v'}\}$. Denote by $C_{\fin}(\pi,\pi'):=\prod_{v<\infty}q_v^{m_v}$.

For $v\mid\mathfrak{M}\mathfrak{N}$, we denote by $m_v':=r_{\omega_v\omega_v'^{-1}}=r_{\omega_v^{-1}\omega_v'}$. Then $0\leq m_v'\leq m_v$. 

%Let $C_{\fin}(\pi)=\prod_{v<\infty}C_v(\pi)$ be the arithmetic conductor of $\pi$ and let $C_{\infty}(\pi)=\prod_{v\mid\infty}C_v(\pi)$ be the archimedean conductor of $\pi$. 

%Let $\chi=\otimes_v\chi_v$ be a Hecke character over $F$. For $v<\infty$, we denote by $\varpi_v^{r_{\chi_v}}$ the smallest nonnegative integer such that $\chi_v$ is trivial over $1+\varpi_v^{r_{\chi_v}}\mathcal{O}_v$ but not over $1+\varpi_v^{r_{\chi_v}-1}\mathcal{O}_v$. 

\begin{comment}
For $F_v\simeq \mathbb{R}$, $\chi_v=\sgn^{n_v'}|\cdot|^{i\kappa_v}$, $n_v'\in\{0,1\}$, we define 
\begin{align*}
C_v(\chi)=1+\Big|\frac{n_v'+i\kappa_v}{2}\Big|;
\end{align*}
if $F_v\simeq \mathbb{C}$, and $\chi_v(a)=(a|a|^{-1})^{n_v'}|a|^{2i\kappa_v}$, $a\in F_v^{\times}$, we define
\begin{align*}
C_v(\chi):=(1+|i\kappa_v+|n_v'|/2|)^2;
\end{align*}

Let $C_v(\chi)=q_v^{r_{\chi_v}}$. Denote by $C_{\infty}(\chi):=\otimes_{v\mid\infty}C_v(\chi)$ and $C_{\fin}(\chi):=\otimes_{v<\infty}C_v(\chi)$. Let $C(\chi):=C_{\infty}(\chi)C_{\fin}(\chi)$ be the \textit{analytic} conductor of $\chi$.
\end{comment}

\subsection{The Subconvexity Data}\label{sec4.2.}
Let $\pi$ and $\pi'$ be as before. Let $t\in \mathbb{R}$. The major goal of this paper is to prove a subconvex bound of the form
\begin{align*}
L(1/2+it,\pi\times\pi')\ll_{\pi'} C(\pi\otimes |\cdot|^{it})^{1/2-\delta}
\end{align*}
for some $\delta>0$. Notations $\pi$, $\pi'$ and $t$ will be frequently used henceforth.

\subsection{The Amplification Data}\label{sec4.2}
Let $L>10^{3}$ be a sufficiently large constant, and define 
\begin{equation}\label{equ4.2}
\mathcal{P}:=\big\{v\in \Sigma_{F,\fin}:\ v\nmid \mathfrak{D}_F\mathfrak{M}\mathfrak{N},\ L<q_v\leq 2L,\ |\lambda_{\sigma'}(\mathfrak{p}_v)|\leq q_v^{\varepsilon}\big\}
\end{equation}
which is a finite set of non-Archimedean places.
Here $\mathfrak{D}_F$ denotes the discriminant of $F$.
By the Rankin--Selberg method together with bounds arising from the symmetric
power lifts of $\pi'_v$ (see \cite[\textsection 4]{Yan21}), we obtain  
\begin{equation}\label{equ4.3}
L^{1-10\varepsilon}\ll \#\mathcal{P}\ll L. 
\end{equation}

Let $\boldsymbol{\ell}=(\ell_v)_{v\in \mathcal{P}}\in \mathbb{Z}_{\geq 0}^{\#\mathcal{P}}$ be a sequence of positive  integers. Let $\boldsymbol{\alpha}=(\alpha_v)_{v\in \mathcal{P}}\in \mathbb{C}^{\#\mathcal{P}}$ be a sequence of complex numbers. 

Let $\mathfrak{L}=\prod_{v\in \mathcal{P}}\mathfrak{p}_v^{\ell_v'}$ for some $0\leq \ell_v'\leq 2\ell_v$, $v\in \mathcal{P}$. Denote by $\|\mathfrak{L}\|=\prod_{v\in \mathcal{P}}q_v^{\ell_v'}$. We say $v\mid\mathfrak{L}$ if $v\in \mathcal{P}$ and $\ell_v'\geq 1$.  

\subsection{The Regularization Data}\label{sect4.4}
Let $t\in\mathbb{R}$ be as in \textsection\ref{sec4.2.}, and $L$ be defined as in \textsection\ref{sec4.2}. Let $v\mid\infty$. Let $\gamma_v\in \mathbb{C}$ be such that $L_v(s,\omega_v^{-1}\omega_v')=\Gamma_v(s+\gamma_v)$. By Pigeonhole principle, there exists some $(\log C_{\infty}(\pi')+\log C_{\fin}(\pi,\pi')+\log L)^{-10}\leq \varepsilon_0\leq 2(\log C_{\infty}(\pi')+\log C_{\fin}(\pi,\pi')+\log L)^{-10}$ such that for all $|s|=\varepsilon_0$, we have 
\begin{equation}\label{equa4.2}
\min_{v\mid\infty}|it+\gamma_v/2-s|\geq \frac{1}{4[F:\mathbb{Q}](\log C_{\infty}(\pi')+\log C_{\fin}(\pi,\pi')+\log L)^{10}},
\end{equation}
where $[F:\mathbb{Q}]$ refers to the degree of the number field $F.$

\subsection{Choice of Cusp Forms in $\pi'$}\label{sec4.3}
\begin{itemize}
\item For $v<\infty$, we take $W_{\pi_v'}$ to be the local new vector in the Kirillov model of $\pi_v'$, normalized by $\langle W_{\pi_v'},W_{\pi_v'}\rangle_v=1$.
\item For $v\mid\infty$, we fix $\varphi_v^{\circ}\in C_c^{\infty}((0,+\infty))$ with $\varphi_v^{\circ}(1)=1$ and $\varphi_v^{\circ}(t_v)=0$ if $|t_v-1|_v>1$. For $t_v\in F_v^{\times}$, we define $\varphi_v(t_v)=\varphi_v^{\circ}(|t_v|_v^{1/[F_v:\mathbb{R}]})$.

Take $W_{\pi_v'}$ to be the vector in the Kirillov model satisfying that 
\begin{align*}
W_{\pi_v'}\left(\begin{pmatrix}
y_v\\
& 1
\end{pmatrix}\right)=\varphi_v(y_v)\ \ \text{for all $y_v\in F_v^{\times}$}.
\end{align*}
\end{itemize}

Let $\phi\in \pi'$ be the cusp form corresponding to $\otimes_vW_{\pi_v'}$. Let $\boldsymbol{d}_{\mathfrak{L}}=\otimes_{v\in \mathcal{P}}\ \boldsymbol{d}_v\in G(\mathbb{A}_F)$, where $\boldsymbol{d}_v:=\begin{pmatrix}
1&\\
&\varpi_v^{\ell_v'}
\end{pmatrix}$. Define  
\begin{equation}\label{eq4.1}
\phi_{\mathfrak{L}}:=\prod_{v\in\mathcal{P}}\omega_v^{-1}(\varpi_v^{l_v'})(q_v^{\frac{\ell_v'}{2}}+\ell_v'q_v^{-\frac{\ell_v'}{2}}|\mathbb{F}_v^{\times}|)\cdot \pi'(\boldsymbol{d}_{\mathfrak{L}})\phi.
\end{equation}

\subsection{Construction of Eisenstein Series}\label{sec4.4} 
Let $0<\varepsilon^*<10^{-3}$. Let $t\in \mathbb{R}$ be as in \textsection\ref{sec4.2.}. Let $\mathcal{P}$ be the finite set defined by \eqref{equ4.2} in \textsection\ref{sec4.2}. For $v\mid\infty$, let $h_v$ be a fixed nonnegative smooth  function on $\mathbb{R}$, supported in the interval $|t|\leq \varepsilon^*$, satisfying $h_v(t_v)\equiv 1$ when $|t|\leq \varepsilon^*/2$. 
\begin{itemize}
\item For $v\mid\infty$, we fix a tiny constant $\varepsilon_1>0$ and let 
\begin{equation}\label{f4.5}
\textbf{C}_v(\pi,t,\pi'):=\max\big\{C_v(\pi\otimes|\cdot|^{it}),C_v(\pi')\big\}^{1+\varepsilon_1},
\end{equation}
where $C_v(\cdot)$ is the local archimedean conductor defined by \eqref{cond}.
Define
\begin{equation}\label{eq6.28}
\Phi_v(t_{1,v},t_{2,v}):=\textbf{C}_v^{\, 1/2}h_v(\textbf{C}_v|t_{1,v}|_v)h_v(|t_{2,v}|_v-1),
\end{equation}
where $\textbf{C}_v=\textbf{C}_v(\pi,t,\pi')$.

\item For $v\mid\mathfrak{M}\mathfrak{N}$, namely,  $m_v:=\max\{r_{\pi_v},r_{\pi_v'}\}\geq 1$ (see \textsection\ref{sec4.1}), we set
\begin{equation}\label{equa4.5}
\Phi_v(t_{1,v},t_{2,v}):=\Vol(K_v[m_v])^{-1}\omega_v\omega_v'^{-1}(t_{2,v})\mathbf{1}_{\mathfrak{p}_v^{m_v}}(t_{1,v})\mathbf{1}_{\mathcal{O}_v^{\times}}(t_{2,v}). 
\end{equation}
\item For the remaining $v$'s, i.e., $v<\infty$, and $v\nmid\mathfrak{M}\mathfrak{N}$, we set
\begin{equation}\label{equ4.9}
\Phi_v(t_{1,v},t_{2,v}):=\mathbf{1}_{\mathcal{O}_v}(t_{1,v})\mathbf{1}_{\mathcal{O}_v}(t_{2,v}).
\end{equation}
\end{itemize}
Let $\Phi(t_1,t_2):=\otimes_v\Phi_v(t_{1,v},t_{2,v})\in \mathcal{S}(\mathbb{A}_F^2)$, $t_i=\otimes_vt_{i,v}\in \mathbb{A}_F$, $1\leq i\leq 2$. Notice that $\Phi(\cdot,\cdot)$ also depends on the parameter $t\in \mathbb{R}$. 

Following the notation in \textsection\ref{sec3.1.2}, for $g\in G(\mathbb{A}_F)$, we define
\begin{equation}\label{4.1}
\begin{cases}
E(g,s):=E_{\Phi,\mathbf{1},\omega\omega'^{-1}}(g,s),\\ 
E_{\mathfrak{L}}(g,s):=E_{\Phi,\mathbf{1},\omega\omega'^{-1}}(g\boldsymbol{d}_{\mathfrak{L}},s).
\end{cases}
\end{equation}
For simplicity we write $h_2(g,s)=h_{\Phi,\mathbf{1},\omega\omega'^{-1}}(g\boldsymbol{d}_{\mathfrak{L}},s)$, which is defined by \eqref{2.3}. It is noteworthy that \eqref{4.1} aligns with the notation in \textsection\ref{sec3} by taking $\Phi_1=\Phi$, $\Phi_2((\cdot,\cdot))=|\det \boldsymbol{d}_{\mathfrak{L}}|^{s}\Phi((\cdot,\cdot)\boldsymbol{d}_{\mathfrak{L}})$, $\chi_1=\chi_1'=\mathbf{1}$ $\chi_2=\chi_2'=\omega\omega'^{-1}$, $E(g,s)=E_1(\cdot,s)$, and $E_{\mathfrak{L}}(g,s)=E_2(\cdot,s)$.

By definition, the $v$-th component of $h_2(g,s)=h_{\Phi,\mathbf{1},\omega\omega'^{-1}}(g\boldsymbol{d}_{\mathfrak{L}},s)$ is 
\begin{equation}\label{eq4.10}
\overline{h_{2,v}(g_v\boldsymbol{d}_v,\overline{s})}=|\det g_v\boldsymbol{d}_v|_v^{s}\int_{F_v^{\times}}\overline{\Phi_v((0,t_v)g_v\boldsymbol{d}_v)}\omega_v\omega_v'^{-1}(t_v)|t_v|_v^{2s}d^{\times}t_v,
\end{equation}
for all $g_v\in G(F_v)$. Hence, a straightforward calculation yields the following. 
\begin{lemma}\label{lem4.1}
Let notation be as before. Let $v\mid\mathfrak{L}$, i.e., $\ell_v'\geq 1$. Let $k_v\in K_v$. Then 
\begin{multline}\label{4.10}
|\det \boldsymbol{d}_v|_v^{-s}\overline{h_{2,v}(k_v\boldsymbol{d}_v,\overline{s})}=L_v(2s,\omega_v\omega_v'^{-1})\omega_v^{-1}\omega_v'(\varpi_v^{l_v'})q_v^{2sl_v'}\mathbf{1}_{k_v\in K_v[l_v']}\\
+L_v(2s,\omega_v\omega_v'^{-1})\sum_{j=0}^{\ell_v'-1}\omega_v^{-1}\omega_v'(\varpi_v^j)q_v^{2sj}\mathbf{1}_{k_v\in K_v[j]-K_v[j+1]}.
\end{multline}
\end{lemma}

\subsection{Notation}\label{sec4.5}
Let $\mathfrak{L}$ and $\phi$ be defined \textsection\ref{sec4.2} and \textsection\ref{sec4.3}. 
Henceforth, we will take $\phi_1=\phi$, $\phi_2=\phi_{\mathfrak{L}}$, $E_1(\cdot,s)=E(g,s)$, and $E_2(\cdot,s)=E_{\mathfrak{L}}(g,s)$. The integrals $\mathcal{M}_{\cusp}^{\tw}(s_1,s_2)$, $\mathcal{M}_{\Eis}^{\tw}(s_1,s_2)$ defined in \textsection\ref{sec3.1.1}, and $\mathcal{M}_{\cusp}^{\du}(s_1,s_2)$, $\mathcal{M}_{\Eis}^{\du,\heartsuit}(s_1,s_2)$ and $\mathcal{M}_{\const}^{\du,\heartsuit}(s_1,s_2)$ defined in \textsection\ref{sec3.1.2.} will be denoted as $\mathcal{M}_{\cusp}^{\tw}(s_1,s_2;\mathfrak{L})$, $\mathcal{M}_{\Eis}^{\tw}(s_1,s_2;\mathfrak{L})$, $\mathcal{M}_{\cusp}^{\du}(s_1,s_2;\mathfrak{L})$, $\mathcal{M}_{\Eis}^{\du,\heartsuit}(s_1,s_2;\mathfrak{L})$ and $\mathcal{M}_{\const}^{\du,\heartsuit}(s_1,s_2;\mathfrak{L})$, respectively. Under this notation, Theorem \ref{thmA} boils down to the following. 
\begin{cor}\label{cor4.1}
Let notation be as before. Then 
\begin{multline*}
\mathcal{M}_{\cusp}^{\tw}(s_1,s_2;\mathfrak{L})+\mathcal{M}_{\Eis}^{\tw}(s_1,s_2;\mathfrak{L})=\mathcal{M}_{\const}^{\du,\heartsuit}(s_1,s_2;\mathfrak{L})\\
+\mathcal{M}_{\cusp}^{\du}(s_1,s_2;\mathfrak{L})+\mathcal{M}_{\Eis}^{\du,\heartsuit}(s_1,s_2;\mathfrak{L})+\mathcal{M}_{\Res}^{\du,\heartsuit}(s_1,s_2;\mathfrak{L})
\end{multline*}
holds as an identity of meromorphic functions in 
\begin{align*}
\mathcal{D}=\big\{(s_1,s_2)\in \mathbb{C}: \ -1/2<\Re(s_2-s_1)<1/2,\ 1/2<\Re(s_2+s_1)<3/2\big\}.
\end{align*}
\end{cor}

\subsection{The Amplified Regularized Spectral Reciprocity}\label{sec4.6}
Let $\boldsymbol{\alpha}=(\alpha_v)_{v\in \mathcal{P}}\in \mathbb{C}^{\#\mathcal{P}}$. For a function $\mathcal{M}_*^*(s_1,s_2;\mathfrak{L})$, we define
\begin{multline}\label{equ4.13}
\mathcal{M}_*^*(s_1,s_2;\boldsymbol{\alpha},\boldsymbol{\ell}):=
\sum_{\substack{v,v'\in \mathcal{P}\\ v\neq v'}}\alpha_v\overline{\alpha_{v'}}\omega_{v'}^{-1}(\varpi_{v'}^{\ell_{v'}})\mathcal{M}_*^*(s_1,s_2;\mathfrak{p}_v^{\ell_v}\mathfrak{p}_{v'}^{\ell_{v'}})\\
+\sum_{v\in \mathcal{P}}\sum_{n_v=0}^{\ell_v}|\alpha_v|_v^2\omega_v^{-1}(\varpi_v^{\ell_v-n_v})\mathcal{M}_*^*(s_1,s_2;\mathfrak{p}_v^{2\ell_v-2n_v}).
\end{multline}
Here $\mathcal{M}_*^*\in \big\{\mathcal{M}_{\cusp}^{\tw},\ \mathcal{M}_{\Eis}^{\tw},\ \mathcal{M}_{\const}^{\du,\heartsuit},\ \mathcal{M}_{\cusp}^{\du},\ \mathcal{M}_{\Eis}^{\du,\heartsuit},\ \mathcal{M}_{\Res}^{\du,\heartsuit}\big\}$. Then Corollary \eqref{cor4.1} becomes
\begin{cor}\label{cor4.2}
Let notation be as before. Then 
\begin{multline*}
\mathcal{M}_{\cusp}^{\tw}(s_1,s_2;\boldsymbol{\alpha},\boldsymbol{\ell})+\mathcal{M}_{\Eis}^{\tw}(s_1,s_2;\boldsymbol{\alpha},\boldsymbol{\ell})=\mathcal{M}_{\const}^{\du,\heartsuit}(s_1,s_2;\boldsymbol{\alpha},\boldsymbol{\ell})\\
+\mathcal{M}_{\cusp}^{\du}(s_1,s_2;\boldsymbol{\alpha},\boldsymbol{\ell})+\mathcal{M}_{\Eis}^{\du,\heartsuit}(s_1,s_2;\boldsymbol{\alpha},\boldsymbol{\ell})+\mathcal{M}_{\Res}^{\du,\heartsuit}(s_1,s_2;\boldsymbol{\alpha},\boldsymbol{\ell})
\end{multline*}
holds as an identity of meromorphic functions in 
\begin{align*}
\mathcal{D}=\big\{(s_1,s_2)\in \mathbb{C}: \ -1/2<\Re(s_2-s_1)<1/2,\ 1/2<\Re(s_2+s_1)<3/2\big\}.
\end{align*}
\end{cor}

\subsection{Evaluation at $s_1=1/2+it$ and $s_2=1/2-it$}
Notice that $\mathcal{M}_{\const}^{\du,\heartsuit}(s_1,s_2;\boldsymbol{\alpha},\boldsymbol{\ell})$ and $\mathcal{M}_{\Eis}^{\du,\heartsuit}(s_1,s_2;\boldsymbol{\alpha},\boldsymbol{\ell})$ are just meromorphic at $s_1=1/2+it$ and $s_2=1/2-it$, with poles of order at most $3$. Based on Corollary \ref{cor4.2} and Cauchy's formula, we have the following. 
\begin{thm}\label{thm4.4}
Let notation be as before. Let $C_{\fin}(\pi,\pi')$ be the conductor defined in \textsection\ref{sec4.1},  $t\in \mathbb{R}$ be the subconvexity parameter in \textsection\ref{sec4.2.}, and $L>10^{3}$ be the amplification parameter in \textsection\ref{sec4.2}. Let $\varepsilon_0$ be the parameter defined in \textsection\ref{sect4.4}. Let $s_1=1/2+it$ and $s_2=1/2-it$. Then 
\begin{multline*}
\mathcal{M}_{\cusp}^{\tw}(s_1,s_2;\boldsymbol{\alpha},\boldsymbol{\ell})+\mathcal{M}_{\Eis}^{\tw}(s_1,s_2;\boldsymbol{\alpha},\boldsymbol{\ell})=\mathcal{M}_{\const}^{\du,\heartsuit,\Reg}(s_1,s_2;\boldsymbol{\alpha},\boldsymbol{\ell})\\
+\mathcal{M}_{\cusp}^{\du}(s_1,s_2;\boldsymbol{\alpha},\boldsymbol{\ell})+\mathcal{M}_{\Eis}^{\du,\heartsuit}(s_1,s_2;\boldsymbol{\alpha},\boldsymbol{\ell})+\mathcal{M}_{\Res}^{\du,\heartsuit,\Reg}(s_1,s_2;\boldsymbol{\alpha},\boldsymbol{\ell}),
\end{multline*}
where for $*\in \{\const,\Res\}$, we define 
\begin{equation}\label{eq4.5}
\mathcal{M}_{*}^{\du,\heartsuit,\Reg}(s_1,s_2;\boldsymbol{\alpha},\boldsymbol{\ell}):=\frac{1}{2\pi i}\int_{|s|=\varepsilon_0}\frac{\mathcal{M}_{*}^{\du,\heartsuit}(s_1-s,s_2+2s;\boldsymbol{\alpha},\boldsymbol{\ell})}{s}ds.
\end{equation}
\end{thm}

\begin{remark}
In the above definition, we utilize $s_2+2s$ instead of $s_2+s$ because $\mathcal{M}_{\Eis}^{\du,\heartsuit,\Reg}(s_1,s_2;\boldsymbol{\alpha},\boldsymbol{\ell})$ might exhibit a pole on the line $s_1=s_2$ as shown in Proposition \ref{prop3.10}. By using $2s$ instead of $s$, we aim to circumvent any potential poles along the contour.
\end{remark}

%The notations presented in this section will be used extensively throughout the remainder of this paper.

\section{The Twisted Moment}\label{sec5}
Let $\boldsymbol{\ell}=(\ell_v)_{v\in \mathcal{P}}\in \mathbb{Z}_{\geq 0}^{\#\mathcal{P}}$, and $\boldsymbol{\alpha}=(\alpha_v)_{v\in \mathcal{P}}\in \mathbb{C}^{\#\mathcal{P}}$. Let $\mathcal{M}_{\cusp}^{\tw}(s_1,s_2;\boldsymbol{\alpha},\boldsymbol{\ell})$ and $\mathcal{M}_{\Eis}^{\tw}(s_1,s_2;\boldsymbol{\alpha},\boldsymbol{\ell})$ be the 
the amplified twisted moments  defined in \textsection\ref{sec4.6}. The main result in this section is the following estimate.
\begin{prop}\label{prop5.3}
Let notation be as before. Let $\pi\in \mathcal{A}_0([G],\omega)$. Suppose that $\pi$ is right invariant under $K_{\fin}[\mathfrak{N}]$. Let $s_1=1/2+it$, where $t\in \mathbb{R}$ is the parameter in \textsection\ref{sec4.2.}. Let $L>1$ and $\mathcal{P}$ be defined in \eqref{equ4.2} in \textsection\ref{sec4.2}.  Then 
\begin{align*}
\mathcal{M}_{\cusp}^{\tw}(s_1,\overline{s_1};\boldsymbol{\alpha},\boldsymbol{\ell})+\mathcal{M}_{\Eis}^{\tw}(s_1,\overline{s_1};\boldsymbol{\alpha},\boldsymbol{\ell})\gg \frac{|L(1/2+it,\overline{\pi}\times\pi')|^2\big|\sum_{v\in \mathcal{P}}\alpha_v\lambda_{\pi}(\mathfrak{p}_v^{\ell_v})\big|^2}{\textbf{C}_{\infty}(\pi,t,\pi')^{1+\varepsilon}C_{\fin}(\pi,\pi')^{\varepsilon}\log^2 L},
\end{align*}  
where the implied constant depends on $F$ and $\varepsilon$. Here $C_{\fin}(\pi,\pi')$ is defined as in \textsection\ref{sec4.1}, and $\textbf{C}_{\infty}(\pi,t,\pi'):=\prod_{v\mid\infty}\textbf{C}_v(\pi,t,\pi')$ is the archimedean conductor defined as in \eqref{f4.5} in \textsection\ref{sec4.4}. 
\end{prop}

\subsection{The Amplified Twisted Moment}
\begin{lemma}
Let notation be as before. Let $\omega, \xi\in \widehat{F^{\times}\backslash\mathbb{A}_F^{(1)}}$, and $s\in \mathbb{C}$. Let $\sigma\in \mathcal{A}_0([G],\omega)$, or $\sigma=\pi_{\xi,\overline{\xi}\omega,s}$ (see \eqref{2.2}). Let $\varphi\in \sigma$. Suppose that $\langle \phi E(\cdot,s_1),\varphi\rangle\not\equiv 0$. Then $\sigma_{v}$ is spherical at $v\in \mathcal{P}$, and 
\begin{equation}\label{eq5.1}
\langle \varphi,\phi_{\mathfrak{L}}E_{\mathfrak{L}}(\cdot,\overline{s_2})\rangle=\lambda_{\sigma}(\mathfrak{L})\langle \varphi,\phi E(\cdot,\overline{s_2})\rangle.
\end{equation}
\end{lemma}
\begin{proof}
Analyzing the $K$-types, we have $\langle \phi E(\cdot,s_1),\varphi\rangle\equiv 0$ unless $\varphi$ is right invariant under $\prod_{v<\infty, v\nmid\mathfrak{MN}}K_v$. In particular, $\langle \phi_1E_1(\cdot,s_1),\varphi\rangle\equiv 0$ unless $\varphi$ is right $\prod_{v\in \mathcal{P}}K_v$-invariant.

Let $w=\begin{pmatrix}
&1\\
1
\end{pmatrix}$. Let $m\geq 0$. By Cramer's rule, 
\begin{align*}
K_v=\bigsqcup_{\alpha\in \mathcal{O}_v/\mathfrak{p}_v^m}\begin{pmatrix}
1&\alpha\\
&1
\end{pmatrix}wK_v[m]w\bigsqcup\bigsqcup_{\beta\in \mathfrak{p}_v/\mathfrak{p}_v^m}\begin{pmatrix}
&1\\
1 & \beta
\end{pmatrix}wK_v[m]w.
\end{align*}

Since $wK_v[m]w\begin{pmatrix}
\varpi_v^m\\
&1	
\end{pmatrix}=\begin{pmatrix}
\varpi_v^m\\
&1	
\end{pmatrix}K_v[m]$, then 
\begin{align*}
K_v\begin{pmatrix}
\varpi_v^m\\
&1	
\end{pmatrix}K_v=\bigsqcup_{\alpha\in \mathcal{O}_v/\mathfrak{p}_v^m}\begin{pmatrix}
\varpi_v^m&\alpha\\
&1	
\end{pmatrix}K_v\bigsqcup \bigsqcup_{j=1}^m\bigsqcup_{\beta\in \mathbb{F}_v^{\times}} \begin{pmatrix}
\varpi_v^{m-j} &\beta\\
& \varpi_v^j
\end{pmatrix}K_v,
\end{align*}
from which we deduce that 
\begin{equation}\label{5.1}
\Vol(K_v\begin{pmatrix}
\varpi_v^m\\
&1	
\end{pmatrix}K_v)=q_v^m+m|\mathbb{F}_v^{\times}|.
\end{equation}	
	
Recall the definition \eqref{hecke} of Hecke operators in \textsection\ref{sec2.2}: 
\begin{align*}
T_{\mathfrak{p}_v^m}(\varphi)=q_v^{-\frac{m}{2}}\int_{K_v\begin{pmatrix}
\varpi_v^m\\
&1	
\end{pmatrix}K_v}\sigma_v(y_v)\varphi dy_v.
\end{align*}

Substituting \eqref{5.1} into the above integral and change coordinates, we have 
\begin{align*}
T_{\mathfrak{p}_v^m}(\varphi)(g)=(q_v^{\frac{m}{2}}+mq_v^{-\frac{m}{2}}|\mathbb{F}_v^{\times}|)\int_{K_v}\varphi\left(gk_v\begin{pmatrix}
\varpi_v^m\\
&1
\end{pmatrix}\right)dk_v,
\end{align*}
for $g\in G(\mathbb{A}_F)$. Therefore, 
\begin{equation}\label{5.2}
(q_v^{\frac{m}{2}}+mq_v^{-\frac{m}{2}}|\mathbb{F}_v^{\times}|)\int_{K_v}\varphi\left(gk_v\begin{pmatrix}
\varpi_v^m\\
&1
\end{pmatrix}\right)dk_v=\lambda_{\sigma}(\mathfrak{p}_v^{m})\varphi(g),
\end{equation}
where $\lambda_{\sigma}(\mathfrak{p}_v^{m})$ is the Hecke eigenvalue. By definition \eqref{eq4.1}, 
\begin{align*}
\langle \varphi,\phi_{\mathfrak{L}}E_{\mathfrak{L}}(\cdot,\overline{s_2})\rangle=\prod_{v\in\mathcal{P}}\omega_v(\varpi_v^{l_v'})(q_v^{\frac{\ell_v'}{2}}+\ell_v'q_v^{-\frac{\ell_v'}{2}}|\mathbb{F}_v^{\times}|)\int_{[\overline{G}]}\varphi(g)\overline{\phi(g\boldsymbol{d})E(g\boldsymbol{d},\overline{s_2})}dg,
\end{align*}
where $\boldsymbol{d}=\otimes_{v\in \mathcal{P}}\ \begin{pmatrix}
1&\\
&\varpi_v^{\ell_v'}
\end{pmatrix}\in G(\mathbb{A}_F)$. Thus,
\begin{equation}\label{5.4}
\langle \varphi,\phi_{\mathfrak{L}}E_{\mathfrak{L}}(\cdot,\overline{s_2})\rangle=\prod_{}\omega_v(\varpi_v^{l_v'})(q_v^{\frac{\ell_v'}{2}}+\ell_v'q_v^{-\frac{\ell_v'}{2}}|\mathbb{F}_v^{\times}|)\int \varphi(g\boldsymbol{d}^{-1})\overline{\phi(g)E(g,\overline{s_2})}dg,
\end{equation}
where $v\in\mathcal{P}$, and the integral is over $g\in [\overline{G}]$. Since $\phi(\cdot)E(\cdot,\overline{s_2})$ is right $\otimes_{v\in\mathcal{P}}K_v$-invariant, then by a change of variable $g\mapsto g\otimes_{v\in\mathcal{P}}k_v$,
\begin{align*}
\int_{[\overline{G}]}\varphi(g\boldsymbol{d}^{-1})\overline{\phi(g)E(g,\overline{s_2})}dg=\int_{[\overline{G}]}\prod_{v\in\mathcal{P}}\int_{K_v}\varphi(gk_v\boldsymbol{d}^{-1})dk_v\overline{\phi(g)E(g,\overline{s_2})}dg
\end{align*}

By \eqref{5.2} we obtain 
\begin{equation}\label{5.5}
\prod_{v\in\mathcal{P}}\omega_v(\varpi_v^{l_v'})(q_v^{\frac{\ell_v'}{2}}+\ell_v'q_v^{-\frac{\ell_v'}{2}}|\mathbb{F}_v^{\times}|)\int_{K_v}\varphi(gk_v\boldsymbol{d}^{-1})dk_v=\prod_{v\in\mathcal{P}}\lambda_{\sigma}(\mathfrak{p}_v^{\ell_v'})\cdot \varphi(g),
\end{equation}
which is equal to $\lambda_{\sigma}(\mathfrak{L})\cdot \varphi(g)$. Consequently, \eqref{eq5.1} follows from \eqref{5.4} and  \eqref{5.5}.
\end{proof}

\begin{prop}\label{prop5.2}
Let notation be as before. Then 
\begin{equation}\label{eq5.6}
\mathcal{M}_{\cusp}^{\tw}(s_1,s_2;\boldsymbol{\alpha},\boldsymbol{\ell})=\sum_{\sigma}\sum_{\varphi\in\mathfrak{B}(\sigma)}\langle \phi E(\cdot,s_1),\varphi\rangle\langle \varphi,\phi E(\cdot,\overline{s_2})\rangle\Big|\sum_{v\in \mathcal{P}}\alpha_v\lambda_{\sigma}(\mathfrak{p}_v^{\ell_v})\Big|^2,
\end{equation}
where $\sigma\in \mathcal{A}_0([G],\omega)$; and $\mathcal{M}_{\cusp}^{\Eis}(s_1,s_2;\boldsymbol{\alpha},\boldsymbol{\ell})$ is equal to 
\begin{align*}
\sum_{\substack{\xi}} \frac{1}{4\pi i}\int_{i\mathbb{R}}\sum_{f} \langle \phi E(\cdot,s_1),E(\cdot,f,s)\rangle \langle E(\cdot,f,s),\phi E (\cdot,\overline{s_2})\rangle \Big|\sum_{v\in \mathcal{P}}\alpha_v\lambda_{\pi_{\xi,\overline{\xi}\omega,s}}(\mathfrak{p}_v^{\ell_v})\Big|^2ds,
\end{align*}
where $\xi\in \widehat{F^{\times}\backslash\mathbb{A}_F^{(1)}}$, and $f\in \mathfrak{B}(\xi,\xi^{-1})$. 
\end{prop}
\begin{proof}
Opening the square and using \eqref{eq2.3} we obtain 
\begin{align*}
\Big|\sum_{v\in \mathcal{P}}\alpha_v\lambda_{\sigma}(\mathfrak{p}_v^{\ell_v})\Big|^2=\sum_{\substack{v,v'\in \mathcal{P}\\ v\neq v'}}\overline{\alpha_{v'}}\alpha_v\overline{\omega}_{v'}(\varpi_{v'}^{\ell_{v'}})\lambda_{\sigma}(\mathfrak{p}_v^{\ell_v}\mathfrak{p}_{v'}^{\ell_{v'}})+\sum_{v\in \mathcal{P}}|\alpha_v|_v^2\overline{\omega}_v(\varpi_v^{\ell_v})\lambda_{\sigma}(\mathfrak{p}_v^{2\ell_v}).
\end{align*}
Utilizing the Hecke relation  \eqref{hecke.} we have
\begin{align*}
\sum_{v\in \mathcal{P}}|\alpha_v|_v^2\overline{\omega}_v(\varpi_v^{\ell_v})\lambda_{\sigma}(\mathfrak{p}_v^{2\ell_v})=\sum_{v\in \mathcal{P}}\sum_{n_v=0}^{\ell_v}|\alpha_v|_v^2\omega_v^{-1}(\varpi_v^{\ell_v-n_v})\lambda_{\sigma}(\mathfrak{p}_v^{2\ell_v-2n_v}).
\end{align*}

Therefore, Proposition \ref{prop5.2} follows from the above identities and the definition in \textsection\ref{sec4.6}.
\end{proof}

\subsection{A Lower Bound for Archimedean Period Integrals}
\subsubsection{Bessel functions}\label{sec5.2.1}
Let $j_{\pi_v'}$ be the Bessel function associated with $\pi_v'$, e.g., see \cite{Cog14}. Then for all $g_v\in G(F_v)$, we have 
\begin{equation}\label{eq5.7}
W_{\pi_v'}\left(\begin{pmatrix}
a_v\\
&1	
\end{pmatrix}wg_v\right)=\omega_v'(a_v)\int_{F_v^{\times}}j_{\pi_v'}(a_vy_v)W_{\pi_v'}\left(\begin{pmatrix}
y_v\\
&1
\end{pmatrix}g_v\right)d^{\times}y_v.
\end{equation}

Let $\chi_v$ be a general multiplicative character of $F_v^{\times}$. Upon replacing $W_v'(\cdot)$ with $W_v'(\cdot)\cdot \chi_v(\det(\cdot))$ in \eqref{eq5.7} we obtain 
\begin{equation}\label{5.8}
j_{\pi_v'\otimes\chi_v}(y_v)=\chi_v^{-1}(-y_v)j_{\pi_v'}(y_v).
\end{equation}

By functional equation, for $\Re(s)\ll 0$ we have 
\begin{equation}\label{f5.8}
\gamma_v(s,\pi_v',\psi_v)
=\int_{F_v^{\times}}j_{\pi_v'}(y_v)|y_v|_v^{1/2-s}d^{\times}y_v,
\end{equation}
where $\gamma_v(s,\pi_v',\psi_v)$ is the $\gamma$-factor associated with $\pi_v'$ relative to the unramified character $\psi_v$. Let $\chi_v$ be a character of $F_v^{\times}$. Combining \eqref{5.8} with \eqref{f5.8} yields 
\begin{equation}\label{twist}
\gamma_v(1/2-s,\pi_v'\otimes\chi_v^{-1},\psi_v)=\chi_v(-1)\int_{F_v^{\times}}\chi_v(y_v)j_{\pi_v'}(y_v)|y_v|_v^{s}d^{\times}y_v,
\end{equation}
where the integral converges in $-1/2+\vartheta\leq \Re(s)\leq -\vartheta$ if $\chi_v$ is unitary.

\begin{comment}
\begin{equation}\label{eq5.8}
\gamma_v(s,\pi_v'\otimes\chi_v,\psi_v)=\chi_v(-1)\int_{F_v^{\times}}\chi_v^{-1}(y_v)j_{\pi_v'}(y_v)|y_v|_v^{1/2-s}d^{\times}y_v,
\end{equation}
where $\gamma_v(s,\pi_v',\psi_v)$ is the $\gamma$-factor associated with $\pi_v'$ relative to the unramified character $\psi_v$.
\end{comment}

\subsubsection{Mellin transforms and inversions}\label{sec5.2.2}
Denote by $|\cdot|$ the usual norm in $\mathbb{C}$, i.e., for $z=x+iy\in\mathbb{C}$, with $x, y\in\mathbb{R}$, $|z|:=\sqrt{x^2+y^2}$. Let $v\mid\infty$. Then $F_v$ is isomorphic to $\mathbb{R}$ or $\mathbb{C}$. Recall the absolute values in $F_v$: if $F_v\simeq \mathbb{R}$, we set $|a_v|_v:=|a_v|$ for $a_v\in \mathbb{R}$; if $F_v\simeq\mathbb{C}$, set $|a_v|_v=|a_v|_v^2$ for $a_v\in \mathbb{C}$. 
\begin{itemize}
\item Let $\varphi$ be a Schwartz function on $\mathbb{R}^{\times}$. Let $\delta\in \mathbb{Z}/2\mathbb{Z}$. The Mellin transform $\mathcal{M}_{\delta}\varphi(s)$ of order $\delta$ is defined by 
\begin{align*}
\mathcal{M}_{\delta}\varphi(s):=\int_{\mathbb{R}^{\times}}\varphi(r)\sgn(r)^{\delta}|r|^sd^{\times}r. 
\end{align*}

Then we have the Mellin inversion:
\begin{equation}\label{m1}
\varphi(r)=\frac{1}{4\pi i}\sum_{\delta\in \mathbb{Z}/2\mathbb{Z}}\sgn(r)^{\delta}\int_{(\alpha)}\mathcal{M}_{\delta}\varphi(s)|r|^{-s}ds.
\end{equation}

Let $F_v\simeq \mathbb{R}$. Combining \eqref{twist} and \eqref{m1} we obtain 
\begin{equation}\label{5.17}
j_{\pi_v'}(y_v)=\frac{1}{4\pi i}\sum_{\delta\in \mathbb{Z}/2\mathbb{Z}}\sgn(-y_v)^{\delta}\int_{(\alpha)}\gamma_v(1/2-s,\pi_v'\otimes\sgn^{\delta},\psi_v)|y_v|_v^{-s}ds.
\end{equation}

\item Let $\varphi$ be a Schwartz function on $\mathbb{C}^{\times}$. Let $m\in \mathbb{Z}$. The Mellin transform $\mathcal{M}_{m}\varphi(s)$ of order $\delta$ is defined by 
\begin{align*}
\mathcal{M}_{m}\varphi(s):=\int_{\mathbb{R}^{\times}}\varphi(z)[z]^{m}|z|^sd^{\times}z=2\int_0^{\infty}\int_0^{2\pi}\varphi(re^{i\phi})e^{im\phi}r^s d\phi d^{\times}r. 
\end{align*}
Here, for $z=re^{i\phi}\in \mathbb{C}$ in the polar form, we denote $[z]=e^{i\phi}$ and $|z|=r$. 

We have the Mellin inversion:
\begin{equation}\label{m2}
\varphi(z)=\frac{1}{8\pi^2i}\sum_{m\in \mathbb{Z}}[z]^{-m}\int_{(\alpha)}\mathcal{M}_{m}\varphi(s)|z|^{-s}ds,
\end{equation}
or in the polar coordinates,
\begin{align*}
\varphi(re^{i\phi})=\frac{1}{8\pi^2i}\sum_{m\in \mathbb{Z}}e^{-im\phi}\int_{(\alpha)}\mathcal{M}_{m}\varphi(s)r^{-s}ds.
\end{align*}

Let $F_v\simeq \mathbb{C}$. Combining \eqref{twist} and \eqref{m2} we obtain 
\begin{equation}\label{5.17.}
j_{\pi_v'}(y_v)=\frac{1}{8\pi^2 i}\sum_{m\in \mathbb{Z}}[-y_v]^{-m}\int_{(\alpha)}\gamma_v(1/2-s,\pi_v'\otimes[\cdot]^{-m},\psi_v)|y_v|_v^{-s}ds.
\end{equation} 
\end{itemize}

\begin{comment}

Executing the Stirling formula leads to  
\begin{equation}\label{eq5.21}
|\gamma_v(1/2-s,\pi_v',\psi_v)|\asymp C_v(\pi'\otimes|\cdot|^{\Im(s)})^{\Re(s)}.
\end{equation} 
As a consequence, \eqref{5.17} and \eqref{5.17.} converges absolutely in $\alpha\leq -10$.

\begin{align*}
c=\begin{cases}
2\pi i, \ & \text{if $F_v\simeq \mathbb{R}$ and $\pi_v'$ }
\end{cases}
\end{align*}
\end{comment}

\subsubsection{Pointwise bounds for Whittaker functions}\label{sec5.2.2}
Let notation be as in \textsection\ref{sec5.2.2}. Let $a_v\in F_v^{\times}$. Consider the polar coordinates of $a_v$:
\begin{itemize}
\item If $F_v\simeq \mathbb{R}$, then $a_v=r\gamma$, where $r=|a_v|_v\in \mathbb{R}_+^{\times}$, and $\gamma\in \{\pm1\}$.
\item If $F_v\simeq \mathbb{C}$, then $a_v=re^{i\phi}$, where $r=|a_v|_v^{1/2}\in \mathbb{R}_+^{\times}$, and $\phi\in [0,2\pi)$. 
\end{itemize}
We call $r$ the modulus of $a_v$. By \eqref{eq5.7} (with $g_v=I_2$) and a change of variables, 
\begin{align*}
W_{\pi_v'}\left(\begin{pmatrix}
a_v\\
&1	
\end{pmatrix}w\right)=\omega_v'(a_v/r)W_{\pi_v'}\left(\begin{pmatrix}
r\\
&1	
\end{pmatrix}w\right).
\end{align*}
Thus, we may regard $W_{\pi_v'}\left(\begin{pmatrix}
a_v\\
&1	
\end{pmatrix}w\right)$ as a differential function of $r$ since $\omega_v'(r)$ is of the form $r^{\beta}$ for some constant $\beta\in i\mathbb{R}_+$. 

Let $l\geq 1$. It is known (see \cite[Proposition 3.2.3]{MV10}, \cite[(29)]{BH08}) that there exists $d$ (relying on $l$) such that 
\begin{equation}\label{5.15}
W_{\pi_v'}(\diag(a_v,1)w)\ll \mathcal{S}_d(W_{\pi_v'})\cdot \min\{|a_v|_v^{\frac{1}{2}-\vartheta}, |a_v|_v^{-m}\},
\end{equation}
where $\mathcal{S}_d(W_{\pi_v'})$ is the Sobolev norm, and  the implied constant depends on $l$. Nevertheless, for our purposes, we need an improvement of \eqref{5.15}, which is proved as follows. 
\begin{lemma}\label{lem5.4}
Let notation be as before. Let $\vartheta$ be a parameter towards the Ramanujan conjecture for $\pi_v'$. Let $a_v\in F_v$ with modulus $r\in\mathbb{R}_+^{\times}$. Let  $m\geq 1$ and $l\geq 0$. Then 
\begin{equation}\label{5.12}
r^l\frac{\partial^l}{\partial r^l}W_{\pi_v'}\left(\begin{pmatrix}
a_v\\
&1	
\end{pmatrix}w\right)\ll \min\{C_v(\pi')^{-\frac{1}{2}+\vartheta}|a_v|_v^{\frac{1}{2}-\vartheta}, C_v(\pi')^{m}|a_v|_v^{-m}\},
\end{equation}
where the implied constant depends only on $l$, $m$, $F_v$ and $\varphi_v$. Here $C_v(\pi')$ is the archimedean conductor of $\pi_v'$ as defined in \eqref{cond} in \textsection\ref{10.10}.
\end{lemma}
\begin{proof}
Recall the definition in \textsection\ref{sec4.3}: for $y_v\in F_v^{\times}$, $W_{\pi_v'}(\diag(y_v,1))=\varphi_v(y_v)$, where $\varphi_v\in C_c^{\infty}(F_v^{\times})$ with $\varphi_v(1)=1$ and $\varphi_v(t_v)=\varphi_v(|t_v|_v)$ for $t_v\in\mathcal{O}_v^{\times}$. Let 
\begin{align*}
\mathcal{M}\varphi_v(s):=\int_{F_v^{\times}}\varphi_v(y_v)|y_v|_v^{s}d^{\times}y_v
\end{align*}
be the Mellin transform of $\varphi_v$. By definition, for $s\not\in \mathbb{Z}_{\leq 0}$,
\begin{equation}\label{equ5.10}
\mathcal{M}\varphi_v(s)\ll \prod_{j=0}^m|s+j|^{-1},
\end{equation}
where the implied constant depends only on $m$, $\varphi_v$, and $F_v$. 

\begin{itemize}
\item Suppose $F_v\simeq \mathbb{R}$. Then 
\begin{align*}
\mathcal{M}\varphi_v(s)=2\int_{0}^{\infty}\varphi_v(r)r^{s}d^{\times}r
\end{align*}

Let $\alpha\in \mathbb{R}$. By Mellin inversion, we have
\begin{equation}\label{equ5.8.}
W_{\pi_v'}\left(\begin{pmatrix}
y_v\\
&1
\end{pmatrix}\right)=\frac{1}{4\pi i}\int_{(\alpha)}\mathcal{M}\varphi_v(s)r^{-s}ds=\frac{1}{4\pi i}\int_{(\alpha)}\mathcal{M}\varphi_v(s)|y_v|_v^{-s}ds.
\end{equation}

\item Suppose $F_v\simeq\mathbb{C}$. Then 
\begin{align*}
\mathcal{M}\varphi_v(s)=\int_{\mathbb{C}^{\times}}\varphi_v(z)|z|_v^{s}d^{\times}z
\end{align*}

Write $z=re^{i\phi}$. Since $d^{\times}z=2d^{\times}rd\phi$. So 
\begin{align*}
\mathcal{M}\varphi_v(s)=2\int_0^{2\pi}\int_{0}^{\infty}\varphi_v(re^{i\phi})r^{2s}d^{\times}rd\phi=4\pi\int_{0}^{\infty}\varphi_v(r)r^{2s}d^{\times}r.
\end{align*}
Let $\alpha\in \mathbb{R}$. By Mellin inversion, we have
\begin{equation}\label{equ5.8}
W_{\pi_v'}\left(\begin{pmatrix}
y_v\\
&1
\end{pmatrix}\right)=\varphi_v(r)=\frac{1}{4\pi^2 i}\int_{(\alpha)}\mathcal{M}\varphi_v(s)|y_v|_v^{-s}ds.
\end{equation}
 
\end{itemize}

Let $\alpha\ll 0$. Substituting \eqref{equ5.8.} and \eqref{equ5.8} into 
\begin{equation}\label{eq5.18}
W_{\pi_v'}\left(\begin{pmatrix}
a_v\\
&1	
\end{pmatrix}w\right)=\omega_v'(a_v)\int_{F_v^{\times}}j_{\pi_v'}(a_vy_v)W_{\pi_v'}\left(\begin{pmatrix}
y_v\\
&1
\end{pmatrix}\right)d^{\times}y_v,
\end{equation}
we deduce that
\begin{align*}
W_{\pi_v'}\left(\begin{pmatrix}
a_v\\
&1	
\end{pmatrix}w\right)=\frac{\omega_v'(a_v)}{4\pi^{[F_v:\mathbb{R}]}i}\int_{(\alpha)}\mathcal{M}\varphi_v(s)\int_{F_v^{\times}}j_{\pi_v'}(a_vy_v)|y_v|_v^{-s}d^{\times}y_vds.
\end{align*}
In conjunction with \eqref{f5.8}, the above expression becomes
\begin{equation}\label{5.11}
W_{\pi_v'}\left(\begin{pmatrix}
a_v\\
&1	
\end{pmatrix}w\right)=\frac{\omega_v'(a_v)}{4\pi^{[F_v:\mathbb{R}]}i}\int_{(\alpha)}\mathcal{M}\varphi_v(s)\gamma_v(1/2+s,\pi_v',\psi_v)|a_v|_v^sds.
\end{equation}

By definition of gamma factors, $\gamma_v(s+1/2,\pi_v',\psi_v)$ is holomorphic in the range $\Re(s)< 1/2-\vartheta$. Moreover, executing the Stirling formula leads to  
\begin{equation}\label{eq5.21}
|\gamma_v(1/2-s,\pi_v',\psi_v)|\asymp C_v(\pi'\otimes|\cdot|^{\Im(s)})^{\Re(s)}.
\end{equation} 

Denote by $W_{\pi_v'}^{(l)}:=\frac{\partial^l}{\partial r^l}W_{\pi_v'}.$  Shifting contour to $\alpha=1/2-\vartheta$ in \eqref{5.11}, along with \eqref{equ5.10} and \eqref{eq5.21}, we deduce that 
\begin{equation}\label{5.13}
r^lW_{\pi_v'}^{(l)}\left(\begin{pmatrix}
a_v\\
&1	
\end{pmatrix}w\right)\ll \frac{|\gamma_v(1-\vartheta,\pi_v',\psi_v)|}{|a_v|_v^{\vartheta-1/2}}\ll C_v(\pi')^{-\frac{1}{2}+\vartheta}|a_v|_v^{\frac{1}{2}-\vartheta},
\end{equation}
where the implied constant depends only on $l$, $F_v$ and $\varphi_v$. 

On the other hand, for $m\in 1/2+\mathbb{Z}_{\geq 0}$, by shifting contour in \eqref{5.11}, we obtain 
\begin{equation}\label{5.14}
r^lW_{\pi_v'}^{(l)}\left(\begin{pmatrix}
a_v\\
&1	
\end{pmatrix}w\right)\ll |\gamma_v(1/2-m,\pi_v',\psi_v)|\cdot |a_v|_v^{-m}\ll C_v(\pi')^{m}|a_v|_v^{-m},
\end{equation}
where the implied constant depends only on $l$, $F_v$ and $\varphi_v$. 

Therefore, \eqref{5.12} follows from \eqref{5.13} and \eqref{5.14}.
\end{proof}

\begin{comment}
\begin{align*}
|a_v|_v^{-s}\int_{(\alpha')}\mathcal{M}\varphi_v(s')\gamma_v(1/2-s',\pi_v',\psi_v)\int_{(\alpha)}\gamma_v(1/2-s,\pi_v'\otimes\sgn^{\delta},\psi_v)\int_{F_v^{\times}}\psi_v(y_vc_v)|y_v|_v^{-(s+s')}\Phi(y_v)d^{\times}y_vdsds'
\end{align*}

\end{comment}

\subsubsection{Approximation}\label{sec5.2.3}
Taking advantage of preparations in \textsection\ref{sec5.2.1} and \textsection\ref{sec5.2.2} we will establish the following explicit estimates for difference of Whittaker values.
%refining the arguments in \cite[\textsection 3.6.3--\textsection 3.6.5]{MV10}. 
\begin{lemma}\label{lem5.6}
Let notation be as before. Suppose $F_v\simeq \mathbb{R}$. Let $c_v\in F_v$. Then 
\begin{equation}\label{f5.16}
W_{\pi_v'}\left(\begin{pmatrix}
a_v\\
c_v&1	
\end{pmatrix}\right)-W_{\pi_v'}\left(\begin{pmatrix}
a_v\\
&1	
\end{pmatrix}\right)\ll |a_v|_v^{1/2-\vartheta}|c_v|_vC_v(\pi')^{1+\varepsilon}, 
\end{equation}
where the implied constant depends only on $\varepsilon$, $F_v$, and $\varphi_v$.
\end{lemma}
\begin{proof}
By \eqref{eq5.7}, we have
\begin{align*}
W_{\pi_v'}\left(\begin{pmatrix}
a_v\\
c_v&1	
\end{pmatrix}\right)=\omega_v'(a_v)\int_{F_v^{\times}}j_{\pi_v'}(a_vy_v)W_{\pi_v'}\left(\begin{pmatrix}
y_v\\
&1
\end{pmatrix}\begin{pmatrix}
1& c_v\\
&1	
\end{pmatrix}w\right)d^{\times}y_v.
\end{align*}

As a consequence, we can express the difference 
\begin{equation}\label{eq5.16}
\mathcal{D}':=W_{\pi_v'}\left(\begin{pmatrix}
a_v\\
c_v&1	
\end{pmatrix}\right)-W_{\pi_v'}\left(\begin{pmatrix}
a_v\\
&1	
\end{pmatrix}\right)
\end{equation}
as the integral 
\begin{equation}\label{5.16}
\omega_v'(a_v)\int_{F_v^{\times}}(\psi_v(y_vc_v)-1)j_{\pi_v'}(a_vy_v)W_{\pi_v'}\left(\begin{pmatrix}
y_v\\
&1
\end{pmatrix}w\right)d^{\times}y_v.
\end{equation}

\begin{comment}
\begin{align*}
\gamma_v(1/2-s,\pi_v',\psi_v)
=\int_{F_v^{\times}}j_{\pi_v'}(y_v)|y_v|_v^{s}d^{\times}y_v
\end{align*}
\end{comment}

Substituting \eqref{5.17} into \eqref{5.16} leads to 
\begin{equation}\label{5.19}
\mathcal{D}'=\frac{\omega_v'(a_v)}{4\pi i}\sum_{\delta\in \mathbb{Z}/2\mathbb{Z}}\int_{(\alpha)}\mathcal{I}_{\delta}(s)\cdot \gamma_v(1/2-s,\pi_v'\otimes\sgn^{\delta},\psi_v)|a_v|_v^{-s}ds,
\end{equation}
where $\mathcal{D}'$ is the difference defined by \eqref{eq5.16} and 
\begin{align*}
\mathcal{I}_{\delta}(s):=\int_{F_v^{\times}}(\psi_v(y_vc_v)-1)W_{\pi_v'}\left(\begin{pmatrix}
y_v\\
&1
\end{pmatrix}w\right)\sgn(-y_v)^{\delta}|y_v|_v^{-s}d^{\times}y_v.
\end{align*}

Let $\alpha:=-1/2+\vartheta$. We integrate by parts in the polar coordinates system, along with the estimate in Lemma \ref{lem5.4} (with $l=0$ therein), obtaining 
\begin{align*}
|\mathcal{I}_{\delta}(s)|\ll & |s|^{-1}\int_{F_v^{\times}}|\psi_v(y_vc_v)-1|\cdot \Big|\frac{\partial}{\partial y_v}W_{\pi_v'}\left(\begin{pmatrix}
y_v\\
&1
\end{pmatrix}w\right)\Big|\cdot |y_v|_v^{1-\alpha}d^{\times}y_v\\
&+|s|^{-1}\int_{F_v^{\times}}\Big|W_{\pi_v'}\left(\begin{pmatrix}
y_v\\
&1
\end{pmatrix}w\right)\Big|\cdot |c_v|_v|y_v|_v^{1-\alpha}d^{\times}y_v.
\end{align*}

Let $\varepsilon>0$ and $m=\floor{100\varepsilon^{-1}}$. By Lemma \ref{lem5.4}, 
\begin{align*}
|\mathcal{I}_{\delta}(s)|\ll & |s|^{-1}\int_{F_v^{\times}}|y_vc_v|_v\cdot (C_v(\pi')^{-1}|y_v|_v)^{\frac{1}{2}-\vartheta} \cdot |y_v|_v^{\frac{1}{2}-\vartheta}\mathbf{1}_{|y_v|_v\leq C_v(\pi')^{1+\varepsilon}}d^{\times}y_v\\
&+|s|^{-1}\int_{F_v^{\times}}|y_vc_v|_v\cdot (C_v(\pi')^{-1}|y_v|_v)^{-m} \cdot |y_v|_v^{\frac{1}{2}-\vartheta}\mathbf{1}_{|y_v|_v> C_v(\pi')^{1+\varepsilon}}d^{\times}y_v.
\end{align*}

A direct estimate leads to 
\begin{equation}\label{5.21}
|\mathcal{I}_{\delta}(s)|\ll |c_v|_v|s|^{-1}\cdot C_v(\pi')^{1+\varepsilon}\cdot C_v(\pi')^{(1+2\varepsilon)(1/2-\vartheta)},
\end{equation}
where the implied constant depends only on $\varepsilon$, $F_v$, and $\varphi_v$. Substituting \eqref{5.21} into \eqref{5.19} we obtain 
\begin{align*}
\mathcal{D}'\ll |a_v|_v^{1/2-\vartheta}|c_v|_v\cdot C_v(\pi')^{1+\varepsilon}\cdot C_v(\pi')^{1/2-\vartheta}\int_{(-1/2+\vartheta)}\frac{|\gamma_v(1/2-s,\pi_v',\psi_v)|}{|s|}ds.
\end{align*}

By \eqref{eq5.21} with $\Re(s)=1/2+\vartheta$, we obtain $|\gamma_v(1/2-s,\pi_v',\psi_v)|\asymp C_v(\pi'\otimes|\cdot|^{\Im(s)})^{-1/2+\vartheta}$. It then follows from the above estimate that  
\begin{equation}\label{equa5.16}
\mathcal{D}'\ll |a_v|_v^{1/2-\vartheta}|c_v|_vC_v(\pi')^{1+\varepsilon}, 
\end{equation}
where the implied constant depends only on $\varepsilon$, $F_v$, and $\varphi_v$.

Then \eqref{f5.16} follows from \eqref{equa5.16}.
\end{proof}

\begin{lemma}\label{lemma5.7}
Let notation be as before. Suppose $F_v\simeq \mathbb{R}$. Let $a_v\in F_v^{\times}$, and  $c_v\in F_v$. Let $m\geq 100$ and $\varepsilon>0$. Suppose $|a_v|_v\geq C_v(\pi')^{1+\varepsilon}$. Then 
\begin{equation}\label{f5.28}
W_{\pi_v'}\left(\begin{pmatrix}
a_v\\
c_v&1	
\end{pmatrix}\right)-W_{\pi_v'}\left(\begin{pmatrix}
a_v\\
&1	
\end{pmatrix}\right)\ll \frac{|c_v|_v|a_v|_v^{(1+\varepsilon)^2}+|c_v|_v^{m}|a_v|_v^{\varepsilon(1+\varepsilon)}}{|a_v|_v^{m-\varepsilon}C_v(\pi')^{-m-3}}, 
\end{equation}
where the implied constant depends only on $\varepsilon$, $m$, $F_v$, and $\varphi_v$.
\end{lemma}
\begin{proof}
Let $\mathrm{V}$ be a nonnegative smooth function on $\mathbb{R}_+^{\times}$ satisfying $V(r)\equiv 1$ for $r\leq |a_v|_v^{1+\varepsilon}$, $V(r)\equiv 0$ in $r\geq 2|a_v|_v^{1+\varepsilon}$, and $|\partial^m V(r)/\partial r^m|\leq 1$, $\forall$ $m\geq 1$. 

According to \eqref{5.16} we have the decomposition 
\begin{equation}\label{f5.29}
\mathcal{D}'=\mathcal{D}'_1+\mathcal{D}'_2,
\end{equation}
where 
\begin{align*}
\mathcal{D}'_1=&\omega_v'(a_v)\int_{F_v^{\times}}(\psi_v(y_vc_v)-1)j_{\pi_v'}(a_vy_v)W_{\pi_v'}\left(\begin{pmatrix}
y_v\\
&1
\end{pmatrix}w\right)V(|y_v|_v)d^{\times}y_v,\\
\mathcal{D}'_2=&\omega_v'(a_v)\int_{F_v^{\times}}(\psi_v(y_vc_v)-1)j_{\pi_v'}(a_vy_v)W_{\pi_v'}\left(\begin{pmatrix}
y_v\\
&1
\end{pmatrix}w\right)(1-V(|y_v|_v))d^{\times}y_v.
\end{align*}

By the explicit formula for $j_{\pi_v'}(\cdot)$ in \cite[Proposition 6.1]{CPS90}, and \cite[Proposition 18.3]{Qi20}, in conjunction with 
\begin{align*}
K_{\nu}(r)\propto \frac{\pi}{2r}e^{-r}\Big[1+\frac{4\nu^2-1}{8r}+\frac{(4\nu^2-1)(4\nu^2-9)}{2!(8r)^2}+\cdots\Big],
\end{align*}
and the bound in \cite[Theorem 2.1]{Ole06}
\begin{align*}
\sup_{r\geq 0}\sqrt{r}J_{\nu}(r)\leq \sqrt{\nu^{\frac{1}{3}}+2\nu^{-\frac{1}{3}}+\nu^{-1}},
\end{align*}
we obtain the coarse bound 
\begin{equation}\label{equ5.32}
|j_{\pi_v'}(a_vy_v)|\leq 10^{3}C_v(\pi')^{1/12}|a_vy_v|_v^{1/4}
\end{equation}
for $|a_vy_v|_v\geq C_v(\pi')^{2+2\varepsilon}$. By \eqref{equ5.32} and Lemma \ref{lem5.4} (with $l=0$) we obtain 
\begin{equation}\label{eq5.33}
\mathcal{D}'_2\ll \int_{|y_v|_v\geq |a_v|_v^{1+\varepsilon}}\frac{C_v(\pi')^{1/12}|c_vy_v|_v}{|a_vy_v|_v^{-1/4}}\cdot \frac{C_v(\pi')^{m+2}}{|y_v|_v^{m+2}}d^{\times}y_v\ll \frac{|c_v|_vC_v(\pi')^{m+3}}{|a_v|_v^{(1+\varepsilon)m}},
\end{equation}
where the implied constant depends only on $m$ and $\varphi_v$.
 
Now we proceed to handle $\mathcal{D}'_1$. Write $y_v=r\gamma_v$, where $r=|y_v|_v\in \mathbb{R}_+^{\times}$, and $\gamma_v\in \{\pm1\}$. Substituting the Mellin inversions \eqref{5.17} and \eqref{5.11} into \eqref{5.16} leads to  
\begin{align*}
\mathcal{D}'_1=&-\frac{\omega_v'(a_v)}{16\pi^2}\sum_{\substack{\delta\in \mathbb{Z}/2\mathbb{Z}\\
\gamma_v\in \{\pm 1\}}}\sgn(-a_v\gamma_v)^{\delta}\omega_v'(\gamma_v)\int_{(\alpha)}\gamma_v(1/2-s,\pi_v'\otimes\sgn^{\delta},\psi_v)|a_v|_v^{-s}\\
&\int_{(\alpha')}\mathcal{M}\varphi_v(s')\gamma_v(1/2+s',\pi_v',\psi_v)\int_{\mathbb{R}_+^{\times}}(\psi_v(r\gamma_vc_v)-1)V(r)r^{-s+s'}d^{\times}rds'ds.
\end{align*}

Let $\alpha'=-10$ and $\alpha=m-10-\varepsilon$. Integrating by parts, we obtain 
\begin{align*}
\mathcal{D}'_1\ll \int_{(\alpha)}&|\gamma_v(1/2-s,\pi_v'\otimes\sgn^{\delta},\psi_v)||a_v|_v^{-\alpha}\int_{(-10)}|\mathcal{M}\varphi_v(s')\gamma_v(1/2+s',\pi_v',\psi_v)|\\
&\prod_{j=1}^{m}\frac{1}{|j-s+s'-1|}\int_{F_v^{\times}}\bigg|\frac{\partial^m\big[(\psi_v(r\gamma_vc_v)-1)V(r)\big]}{\partial r^m}\bigg|r^{\varepsilon}d^{\times}rds'ds.
\end{align*}

Utilizing \eqref{eq5.21} we deduce from the above inequality that 
\begin{equation}\label{f5.34}
\mathcal{D}'_1\ll \frac{C_v(\pi')^{m-10-\varepsilon}}{|a_v|_v^{m-10-\varepsilon}}\cdot C_v(\pi')^{10}\cdot\bigg[|c_v|_v|a_v|_v^{(1+\varepsilon)^2}+|c_v|_v^{m}|a_v|_v^{\varepsilon(1+\varepsilon)}\bigg].
\end{equation}

Then \eqref{f5.28} follows from \eqref{f5.29}, \eqref{eq5.33}, and \eqref{f5.34}. 
\end{proof}

\begin{lemma}\label{lem5.7}
Let notation be as before. Suppose $F_v\simeq \mathbb{C}$. Let $c_v\in F_v$. Let $l\geq 100$. Let $\delta=1/2-\vartheta$. Then 
\begin{align*}
\mathcal{D}':=W_{\pi_v'}\left(\begin{pmatrix}
a_v\\
c_v&1	
\end{pmatrix}\right)-W_{\pi_v'}\left(\begin{pmatrix}
a_v\\
&1	
\end{pmatrix}\right)
\end{align*}
is majorized by 
\begin{equation}\label{5.30}
C_v(\pi')^{\varepsilon}\cdot \min\bigg\{|a_v|_v^{\delta}|c_v|_vC_v(\pi'), \frac{|c_v|_v+|c_v|_v^{l}}{|a_v|_v^{l-1}C_v(\pi')^{-l-3}}\bigg\}+|a_v|_v^{\delta}|c_v|_v^{\frac{1}{2}}C_v(\pi')^{\frac{1}{2}+\varepsilon}, 
\end{equation}
where the implied constant depends only on $\varepsilon$, $l$, $F_v$, and $\varphi_v$.
\end{lemma}
\begin{proof}
Parallel to \eqref{5.16} we have
\begin{equation}\label{5.31}
\mathcal{D}'=\omega_v'(a_v)\int_{F_v^{\times}}(\psi_v(y_vc_v)-1)j_{\pi_v'}(a_vy_v)W_{\pi_v'}\left(\begin{pmatrix}
y_v\\
&1
\end{pmatrix}w\right)d^{\times}y_v.
\end{equation}

Substituting \eqref{5.17.} into \eqref{5.31} leads to 
\begin{equation}\label{5.33}
\mathcal{D}'=\frac{\omega_v'(a_v)}{8\pi^2 i}\sum_{m\in\mathbb{Z}}\int_{(\alpha)}\mathcal{I}_m(s)\cdot \gamma_v(1/2-s,\pi_v'\otimes[\cdot]^{-m},\psi_v)|a_v|_v^{-s}ds,
\end{equation}
where
\begin{align*}
\mathcal{I}_m(s):=\int_{F_v^{\times}}(\psi_v(y_vc_v)-1)W_{\pi_v'}\left(\begin{pmatrix}
y_v\\
&1
\end{pmatrix}w\right)[-y_v]^{-m}|y_v|_v^{-s}d^{\times}y_v.
\end{align*}

Let $\omega_v'(y_v)=[y_v]^{m_0}|y_v|_v^{\beta}$ be the central character of $\pi_v'$, where $m_0\in \mathbb{Z}$ and $\beta\in i\mathbb{R}_+$ (since $\omega_v'$ is unitary). Write $y_v=re^{i\phi}$. By \eqref{eq5.18} we have 
\begin{align*}
W_{\pi_v'}\left(\begin{pmatrix}
y_v\\
&1
\end{pmatrix}w\right)=W_{\pi_v'}\left(\begin{pmatrix}
r\\
&1
\end{pmatrix}w\right)\omega_v'(e^{i\phi})=e^{im_0\phi}\cdot W_{\pi_v'}\left(\begin{pmatrix}
r\\
&1
\end{pmatrix}w\right).
\end{align*}
Hence, under polar coordinates  $y_v=re^{i\phi}$ and $c_v=r'e^{i\phi'}$, we obtain 
\begin{align*}
\mathcal{I}_m(s)=2\int_{0}^{\infty}\int_0^{2\pi}(e^{4\pi i rr'\cos(\phi+\phi')}-1)W_{\pi_v'}\left(\begin{pmatrix}
r\\
&1
\end{pmatrix}w\right)\frac{(-1)^md\phi d^{\times}r}{e^{i(m-m_0)\phi}r^{2s}}.
\end{align*}
In conjunction with \eqref{5.33}, along with the change of variable $m\mapsto m+m_0$, we derive the decomposition $\mathcal{D}'=\sum_{m\in \mathbb{Z}}\mathcal{D}'_m$, where 
\begin{equation}\label{eq5.34}
\mathcal{D}'_m=\frac{(-1)^{m+m_0}\omega_v'(a_v)}{4\pi^2 i}\int_{(\alpha)}\mathcal{J}_m(s)\cdot \gamma_v(1/2-s,\widetilde{\pi}_v'\otimes[\cdot]^{-m},\psi_v)|a_v|_v^{-s}ds,
\end{equation}
where $\widetilde{\pi}_v'=\pi_v'\otimes\omega_v'^{-1}$ is the contragredient of $\pi_v'$, and 
\begin{equation}\label{f5.39}
\mathcal{J}_m:=\int_{0}^{\infty}\int_0^{2\pi}(e^{4\pi i rr'\cos(\phi+\phi')}-1)e^{-im\phi}d\phi W_{\pi_v'}\left(\begin{pmatrix}
r\\
&1
\end{pmatrix}w\right)r^{-2s}d^{\times}r.
\end{equation}

Similar to the proof of Lemma \ref{lem5.6} and Lemma \ref{lemma5.7}, we have
\begin{equation}\label{eq5.32}
\mathcal{D}'_0\ll C_v(\pi')^{\varepsilon}\cdot \min\bigg\{|a_v|_v^{1/2-\vartheta}|c_v|_vC_v(\pi'), \frac{(|c_v|_v+|c_v|_v^{l})C_v(\pi')^{l+3}}{|a_v|_v^{l-1}}\bigg\},
\end{equation}
where the implied constant depends only on $\varepsilon$, $l\geq 100$, $F_v$, and $\varphi_v$.

Now we assume $|m|\geq 1$. Denote by 
\begin{align*}
\mathcal{R}_m(r):=&\int_0^{2\pi}(e^{4\pi i rr'\cos(\phi+\phi')}-1)e^{-im\phi}d\phi=\int_0^{2\pi}e^{4\pi i rr'\cos(\phi+\phi')}e^{-im\phi}d\phi\\
=&e^{im\phi'}\int_0^{2\pi}e^{4\pi i rr'\cos\phi}e^{im\phi}d\phi=2\pi i^me^{im\phi'}J_m(4\pi rr'),
\end{align*}
where $J_m$ is the Bessel $J$-function. 

Utilizing the uniform bound   $J_m(4\pi rr')\ll m^{1/6}(rr')^{-1/2}$ we conclude that \eqref{f5.39} converges absolutely in $\Re(s)<1/4-\vartheta$. Hence, $\mathcal{J}_m$ is the Hankel transform of order $m$ of the function $W_{\pi_v'}\left(\begin{pmatrix}
r\\
&1
\end{pmatrix}w\right)r^{-2s}$.

Let $l_1\geq 1$. Integrating by parts, we obtain 
\begin{equation}\label{eq5.36}
\mathcal{R}_m(r)=\frac{1}{(im)^{l_1}}\int_0^{2\pi}e^{-im\phi}\frac{\partial^{l_1}}{\partial\phi^{l_1}}e^{4\pi i rr'\cos(\phi+\phi')}d\phi.
\end{equation}
In particular, it follows from  \eqref{eq5.36} that 
\begin{equation}\label{eq5.37}
\mathcal{R}_m(r)\ll (rr')^{l_1}m^{-l_1},
\end{equation}
where the implied constant depends only on $l_1$.

Let $l_2\geq 0$. Integrating by parts $l_2$ times, we obtain 
\begin{equation}\label{5.36}
\mathcal{J}_m=\prod_{l'=1}^{l_2}\frac{1}{2s+l'-1}\int_{0}^{\infty}\frac{\partial^{l_2}}{\partial r^{l_2}}\bigg[\mathcal{R}_m(r) W_{\pi_v'}\left(\begin{pmatrix}
r\\
&1
\end{pmatrix}w\right)\bigg]\cdot r^{l_2-2s}d^{\times}r.
\end{equation}

Let $g=\begin{pmatrix}
r\\
&1
\end{pmatrix}w$. By Newton's binomial theorem,
\begin{equation}\label{5.39}
\frac{\partial^{l_2}}{\partial r^{l_2}}\bigg[\mathcal{R}_m(r) W_{\pi_v'}(g)\bigg]=\sum_{j=0}^{l_2}\binom{l_2}{j}\frac{\partial^{j}\mathcal{R}_m(r)}{\partial r^{j}} \cdot \frac{\partial^{l_2-j}}{\partial r^{l_2-j}}W_{\pi_v'}(g).
\end{equation}

Utilizing the expression  \eqref{eq5.36} for $\mathcal{R}_m(r)$ we have
\begin{align*}
\frac{\partial^{j}}{\partial r^{j}}\mathcal{R}_m(r)=\frac{(2\pi i r')^j}{(im)^{l_1}}\int_0^{2\pi}e^{-im\phi}\frac{\partial^{l_1}}{\partial\phi^{l_1}}\cos^j(\phi+\phi')e^{4\pi i rr'\cos(\phi+\phi')}d\phi,
\end{align*}
which, by Newton's binomial theorem, is equal to 
\begin{equation}\label{5.41}
\frac{(2\pi i r')^j}{(im)^{l_1}}\sum_{n=0}^{l_1}\binom{l_1}{n}\int_0^{2\pi}e^{-im\phi}\cdot \frac{\partial^{l_1-n}\cos^j(\phi+\phi')}{\partial\phi^{l_1-n}}\cdot \frac{\partial^{n}e^{4\pi i rr'\cos(\phi+\phi')}}{\partial\phi^{n}}d\phi.
\end{equation}

By a straight forward calculation for \eqref{5.41}, we obtain   
\begin{equation}\label{5.42}
r^j\frac{\partial^{j}}{\partial r^{j}}\mathcal{R}_m(r)\ll \frac{(rr')^j}{m^{l_1}}\cdot (1+(rr')^{l_1}).
\end{equation}

Combining \eqref{eq5.37} with \eqref{5.42} we obtain 
\begin{equation}\label{5.43}
r^j\frac{\partial^{j}}{\partial r^{j}}\mathcal{R}_m(r)\ll \frac{(rr')^j(1+(rr')^{l_1})\mathbf{1}_{j\geq 1}+(rr')^{l_1}\mathbf{1}_{j=0}}{m^{l_1}}.
\end{equation}

It the follows from \eqref{5.36}, \eqref{5.39}, and \eqref{5.43} that 
\begin{equation}\label{eq5.43}
\mathcal{J}_m\ll \frac{1}{m^{l_1}}\sum_{j=0}^{l_2}\int_{0}^{\infty}R_j(r)\cdot \bigg|r^{l_2-j}\frac{\partial^{l_2-j}}{\partial r^{l_2-j}}W_{\pi_v'}\left(\begin{pmatrix}
r\\
&1
\end{pmatrix}w\right)\bigg|\frac{d^{\times}r}{r^{2\alpha}},
\end{equation}
where 
\begin{align*}
R_j(r):=\big[(rr')^j(1+(rr')^{l_1})\mathbf{1}_{j\geq 1}+(rr')^{l_1}\mathbf{1}_{j=0}\big]\cdot \prod_{l'=1}^{l_2}|2s+l'-1|^{-1}.
\end{align*}

Let $l_1=l_2=1$, $\alpha=-1/2+\vartheta$. By \eqref{eq5.43} and Lemma \ref{lem5.4}, 
\begin{align*}
\mathcal{J}_m\ll \frac{1}{|s|}\int_{0}^{\infty}\frac{rr'+(rr')^{2}}{m}\cdot \min\bigg\{\frac{r^{1-2\vartheta}}{C_v(\pi')^{\frac{1}{2}-\vartheta}}, \frac{C_v(\pi')^{l_3}}{r^{2l_3}}\bigg\}\cdot r^{-2\alpha}d^{\times}r.
\end{align*}

Let $\varepsilon>0$, and $l_3=\floor{100\varepsilon^{-1}}$. Then $\mathcal{J}_m\ll \mathcal{J}_m^{(1)}+\mathcal{J}_m^{(2)}$, where 
\begin{align*}
\mathcal{J}_m^{(1)}:=&\frac{1}{|s|}\int_{0}^{C_v(\pi')^{1/2+\varepsilon}}\frac{rr'+(rr')^{2}}{m}\cdot \frac{r^{1-2\vartheta}}{C_v(\pi')^{\frac{1}{2}-\vartheta}}\cdot r^{1-2\vartheta}d^{\times}r,\\
\mathcal{J}_m^{(2)}:=&\frac{1}{|s|}\int_{C_v(\pi')^{1/2+\varepsilon}}^{\infty}\frac{rr'+(rr')^{2}}{m^{l_1}}\cdot \frac{C_v(\pi')^{l_3}}{r^{2l_3}}\cdot r^{1-2\vartheta}d^{\times}r.
\end{align*}

By a straightforward estimate we obtain 
\begin{equation}\label{5.44}
\mathcal{J}_m\ll \frac{C_v(\pi')^{\varepsilon}}{|s|}\cdot \frac{r'C_v(\pi')^{1/2+\varepsilon}+r'^{2}C_v(\pi')^{1+2\varepsilon}}{m}\cdot C_v(\pi')^{(1/2+\varepsilon)(1-2\vartheta)}.
\end{equation}

Substituting \eqref{5.44} into \eqref{eq5.34} yields 
\begin{equation}\label{5.45}
\mathcal{D}'_m\ll \frac{|a_v|_v^{1/2-\vartheta}}{m}\int_{(\alpha)} \frac{\big[r'C_v(\pi')^{1/2+\varepsilon}+r'^{2}C_v(\pi')^{1+2\varepsilon}\big]\cdot C_v(\pi')^{1/2-\vartheta}}{C_v(\widetilde{\pi}_v'\otimes[\cdot]^{-m}|\cdot|^{\Im(s)})^{1/2-\vartheta}}\cdot \frac{ds}{|s|},
\end{equation}
where $\alpha=-1/2+\vartheta$. By definition of $C_v(\widetilde{\pi}_v'\otimes[\cdot]^{-m}|\cdot|^{\Im(s)})$ and a brute force calculation, we obtain  
\begin{equation}\label{5.46}
\sum_{m\neq 0}\mathcal{D}'_m\ll |a_v|_v^{1/2-\vartheta}C_v(\pi')^{\varepsilon}\cdot \Big[|c_v|_v^{1/2}C_v(\pi')^{1/2}+|c_v|_vC_v(\pi')\Big].
\end{equation}
 
Therefore, \eqref{5.30} follows from \eqref{eq5.32} and \eqref{5.46}. 
\end{proof}

\begin{lemma}\label{lem5.9}
Let notation be as before. Suppose $F_v\simeq \mathbb{C}$. Let $a_v\in F_v^{\times}$, and  $c_v\in F_v$. Let $m\geq 100$ and $\varepsilon>0$. Suppose $|a_v|_v\geq C_v(\pi')^{1+\varepsilon}$. Then 
\begin{align*}
\mathcal{D}':=W_{\pi_v'}\left(\begin{pmatrix}
a_v\\
c_v&1	
\end{pmatrix}\right)-W_{\pi_v'}\left(\begin{pmatrix}
a_v\\
&1	
\end{pmatrix}\right)
\end{align*}
is majorized by 
\begin{equation}\label{5.52}
\frac{(|c_v|_v^{1/2}+|c_v|_v^{m/2})\cdot C_v(\pi')^{m-20}}{|a_v|_v^{m-30}}+\frac{|c_v|_vC_v(\pi')^{m+3}}{|a_v|_v^{(1+\varepsilon)m}}, 
\end{equation}
where the implied constant depends only on $\varepsilon$, $m$, $F_v$, and $\varphi_v$.
\end{lemma}
\begin{proof}
Let $\mathrm{V}$ be a nonnegative smooth function on $\mathbb{R}_+^{\times}$ satisfying $V(r)\equiv 1$ for $r\leq |a_v|_v^{(1+\varepsilon)/2}$, $V(r)\equiv 0$ in $r\geq 2|a_v|_v^{(1+\varepsilon)/2}$, and $|\partial^m V(r)/\partial r^m|\ll 1$, $\forall$ $m\geq 1$. 

Similar to \eqref{f5.29} we have $\mathcal{D}'=\mathcal{D}'_1+\mathcal{D}'_2$, where 
\begin{align*}
\mathcal{D}'_1=&\omega_v'(a_v)\int_{F_v^{\times}}(\psi_v(y_vc_v)-1)j_{\pi_v'}(a_vy_v)W_{\pi_v'}\left(\begin{pmatrix}
y_v\\
&1
\end{pmatrix}w\right)V(|y_v|_v^{1/2})d^{\times}y_v,\\
\mathcal{D}'_2=&\omega_v'(a_v)\int_{F_v^{\times}}(\psi_v(y_vc_v)-1)j_{\pi_v'}(a_vy_v)W_{\pi_v'}\left(\begin{pmatrix}
y_v\\
&1
\end{pmatrix}w\right)(1-V(|y_v|_v^{1/2}))d^{\times}y_v.
\end{align*}

By the explicit formula for $j_{\pi_v'}(\cdot)$ in \cite{BM03}, and \cite[Proposition 18.5]{Qi20}, in conjunction with the well known estimates for the classical Bessel $J$-functions and $K$-functions, we obtain the coarse bound 
\begin{equation}\label{5.53}
|j_{\pi_v'}(a_vy_v)|\leq 10^{3}C_v(\pi')|a_vy_v|_v
\end{equation}
for $|a_vy_v|_v\geq C_v(\pi')^{2+2\varepsilon}$. In conjunction with Lemma \ref{lem5.4} we obtain 
\begin{equation}\label{5.54}
\mathcal{D}'_2\ll \int_{|y_v|_v\geq |a_v|_v^{1+\varepsilon}}\frac{C_v(\pi')|c_vy_v|_v}{|a_vy_v|_v^{-1}}\cdot \frac{C_v(\pi')^{m+2}}{|y_v|_v^{m+2}}d^{\times}y_v\ll \frac{|c_v|_vC_v(\pi')^{m+3}}{|a_v|_v^{(1+\varepsilon)m}},
\end{equation}
where the implied constant depends only on $m$ and $\varphi_v$.

Now we proceed to handle $\mathcal{D}'_1$. Write $y_v=re^{i\phi}$, $c_v=r'e^{i\phi'}$, where $r, r'\in \mathbb{R}_+^{\times}$, and $0\leq\phi, \phi'<2\pi$. Let $\omega_v'=|\cdot|^{i\beta}[\cdot]^{m_0}$, where $m_0\in \mathbb{Z}$ and $\beta\in \mathbb{R}$. Utilizing the Mellin inversions \eqref{5.17.} and \eqref{5.11} we obtain 
\begin{align*}
\mathcal{D}'_1&=-\frac{\omega_v'(a_v)}{32\pi^4}\int_0^{2\pi}\sum_{m\in \mathbb{Z}}[-a_v]^{-m}e^{-i(m-m_0)\phi}\int_{(\alpha)}\gamma_v(1/2-s,\pi_v'\otimes[\cdot]^{-m},\psi_v)|a_v|_v^{-s}ds\\
& \int_{(\alpha')}\mathcal{M}\varphi_v(s')\gamma_v(1/2+s',\pi_v',\psi_v)ds'\int_{0}^{\infty}(\psi_v(re^{i\phi}c_v)-1)r^{2(s'-s+i\beta)}V(r)d^{\times}rd\phi.
\end{align*}
\begin{comment}
\begin{align*}
\mathcal{D}'_1=&-\frac{\omega_v'(a_v)}{16\pi^2}\sum_{\substack{\delta\in \mathbb{Z}/2\mathbb{Z}\\
\gamma_v\in \{\pm 1\}}}\sgn(-a_v\gamma_v)^{\delta}\omega_v'(\gamma_v)\int_{(\alpha)}\gamma_v(1/2-s,\pi_v'\otimes\sgn^{\delta},\psi_v)|a_v|_v^{-s}\\
&\int_{(\alpha')}\mathcal{M}\varphi_v(s')\gamma_v(1/2+s',\pi_v',\psi_v)\int_{\mathbb{R}_+^{\times}}(\psi_v(r\gamma_vc_v)-1)V(r)r^{-s+s'}d^{\times}rds'ds.
\end{align*}
\end{comment}

Let $\mathcal{R}_m(r)$ be defined as in \eqref{eq5.36}. We can rewrite $\mathcal{D}_1'$ as 
\begin{equation}\label{f5.55}
\mathcal{D}'_1=-\frac{\omega_v'(a_v)}{32\pi^4}\sum_{m\in \mathbb{Z}}[-a_v]^{-m-m_0}\int_{(\alpha)}|a_v|_v^{-s}\int_{(\alpha')}\Upsilon_v(s,s')\mathcal{J}_m(s,s')ds'ds.
\end{equation}
where $\Upsilon_v(s,s'):=\mathcal{M}\varphi_v(s')\gamma_v(1/2-s,\widetilde{\pi}_v'\otimes[\cdot]^{-m},\psi_v)\gamma_v(1/2+s',\pi_v',\psi_v)$, and 
\begin{align*}
\mathcal{J}_m(s,s'):=\int_{0}^{\infty}\mathcal{R}_m(r)V(r)r^{2(s'-s+i\beta)}d^{\times}r.
\end{align*}

Let $l_1\geq 100$ and $l_2\geq l_1+100$. Similar to \eqref{eq5.43} we have 
\begin{equation}\label{eq5.55}
\mathcal{J}_m(s,s')\ll \frac{1}{m^{l_1}}\sum_{j=0}^{l_2}\int_{0}^{\infty}R_j(r)\cdot r^{l_2-j}\frac{\partial^{l_2-j}}{\partial r^{l_2-j}}V(r)\cdot r^{2(\alpha'-\alpha)}d^{\times}r,
\end{equation}
where 
\begin{align*}
R_j(r):=\big[(rr')^j(1+(rr')^{l_1})\mathbf{1}_{j\geq 1}+(rr')^{l_1}\mathbf{1}_{j=0}\big]\cdot \prod_{l'=1}^{l_2}|2s+l'-1|^{-1}.
\end{align*}

By a direct calculation, along with $\frac{\partial^{l_2-j}}{\partial r^{l_2-j}}V(r)\ll 1$, we obtain from \eqref{eq5.55} that
\begin{equation}\label{5.56}
\mathcal{J}_m(s,s')\ll \frac{1}{m^{l_1}}\prod_{l'=1}^{l_2}|2s+l'-1|^{-1}\int_{0}^{|a_v|_v^{\frac{1+\varepsilon}{2}}}(rr'+(rr')^{l_1})\cdot r^{l_2} r^{2(\alpha'-\alpha)}d^{\times}r.
\end{equation}

Take $l_1\geq 100$, $\alpha=l_1-20$, $l_2=\alpha+10$, and $\alpha'=-10$. Combining \eqref{f5.55} with \eqref{5.56}, along with \eqref{eq5.21}, we derive that 
\begin{equation}\label{5.58}
\mathcal{D}'_1\ll |a_v|_v^{20-l_1+\frac{(1+\varepsilon)(10-l_1)}{2}}\Big[r'|a_v|_v^{\frac{1+\varepsilon}{2}}+r'^{l_1}|a_v|_v^{\frac{(1+\varepsilon)l_1}{2}}\Big]\cdot C_v(\pi')^{l_1-20},
\end{equation}
where the implied constant depends only on $\varepsilon$, $l_1$ and $\varphi_v$.

Then \eqref{5.52} follows from \eqref{5.54} and \eqref{5.58}. 
\end{proof}

\subsubsection{Archimedean period integrals: lower bounds}
We will refine the arguments in \cite[\textsection 3.6.3--\textsection 3.6.5]{MV10} by gathering the estimates in \textsection\ref{sec5.2.2} and \textsection\ref{sec5.2.3}. 

\begin{prop}\label{prop5.7}
Let notation be as before. Let $v\mid\infty$. Let $\varphi_v\in C_c^{\infty}(F_v^{\times})$ be defined as in \textsection\ref{sec4.3}. Let $t\in \mathbb{R}$. Let 
$W_{\pi_v,t}$ be the vector in the Kirillov model of $\pi_v\otimes|\cdot|_v^{-it}$  satisfying that $W_{\pi_v,t}(\diag(y_v,1))=\varphi_v(y_v)$ for all $y_v\in F_v^{\times}$. Let $W_{\pi_v}:=W_{\pi_v,t}\otimes|\cdot|^{it}$. Let $\Re(s)\gg 1$. Define 
\begin{align*}
\Psi_v(s,\overline{W_{\pi_v}},W_{\pi_v'}):=\int_{N(F_v)\backslash G(F_v)}\overline{W_{\pi_v}(x_v)}W_{\pi_v'}(x_v)\Phi_v(\eta x_v)|\det x_v|_v^{s}dx_v.
\end{align*}
Then $\Psi_v(s,\overline{W_{\pi_v}},W_{\pi_v'})$ admits a holomorphic continuation to $\Re(s)>2\vartheta$, and 
\begin{equation}\label{f5.26}
\Psi_v(1/2+it,\overline{W_{\pi_v}},W_{\pi_v'})\gg \textbf{C}_v(\pi,t,\pi')^{-1/2}\cdot |\overline{W_{\pi_v}}(I_2)W_{\pi_v'}(I_2)|,
\end{equation}	
where the implied constant depends only on $F_v$, $h_v$, and $\varphi_v$. Here $\varphi_v$ is defined as in \textsection\ref{sec4.3}, $h_v$ and $\textbf{C}_v(\pi,t,\pi')$ are defined as in \textsection\ref{sec4.4}.
\end{prop}
\begin{proof}
Using Bruhat decomposition, we can write $\Psi_v(s,\overline{W_{\pi_v}},W_{\pi_v'})$ as 
\begin{equation}\label{eq5.22}
\int_{F_v}\int_{F_v^{\times}}\overline{W_{\pi_v,t}\left(\begin{pmatrix}
a_v\\
c_v&1
\end{pmatrix}\right)}W_{\pi_v'}\left(\begin{pmatrix}
a_v\\
c_v&1
\end{pmatrix}\right)|a_v|_v^{s-it-1}d^{\times}a_v\mathcal{I}(c_v,s)dc_v,
\end{equation}
where
\begin{equation}\label{5.22}
\mathcal{I}(c_v,s):=\int_{F_v^{\times}}\Phi_v((z_vc_v,z_v))\omega_v^{-1}\omega_v'(z_v)|z_v|_v^{2s-2it}d^{\times}z_v.
\end{equation}

Substituting the definition \eqref{eq6.28} into \eqref{5.22} we obtain 
\begin{equation}\label{eq5.28}
\mathcal{I}(c_v,s)=\textbf{C}_v^{\, 1/2}\int_{F_v^{\times}}h_v(\textbf{C}_v|z_vc_v|_v)h_v(|z_v|_v-1)|z_v|_v^{2s-2it}d^{\times}z_v,
\end{equation}
where $h_v$ is a fixed nonnegative smooth function on $\mathbb{R}$, supported in the disc $|t_v|_v\leq \varepsilon^*$ and $h_v(t_v)\equiv 1$ when $|t_v|_v\leq \varepsilon^*/2$. Hence, $\mathcal{I}(c_v,s)$ defines an entire function in $s\in \mathbb{C}$, and a smooth function in $c_v$ supported in $|c_v|_v\leq \varepsilon^*\textbf{C}_v^{\, -1}$.

Following the strategy in \cite[\textsection 3.6.4]{MV10}, 
we start with $\mathcal{D}:=W_{\pi_v}\left(\begin{pmatrix}
a_v\\
c_v&1
\end{pmatrix}\right)-W_{\pi_v}\left(\begin{pmatrix}
a_v\\
&1
\end{pmatrix}\right)$ and $\mathcal{D}':=W_{\pi_v'}\left(\begin{pmatrix}
a_v\\
c_v&1
\end{pmatrix}\right)-W_{\pi_v'}\left(\begin{pmatrix}
a_v\\
&1
\end{pmatrix}\right)$. By Lemmas \ref{lem5.6}, \ref{lemma5.7}, \ref{lem5.7}, and \ref{lem5.9}, we obtain 
\begin{align*}
\int_{F_v^{\times}}\overline{\mathcal{D}}\cdot W_{\pi_v'}\left(\begin{pmatrix}
a_v\\
&1
\end{pmatrix}\right)|a_v|_v^{s-1-it}d^{\times}a_v\ll |c_v|_vC_v(\pi')^{1+\varepsilon},
\end{align*}
where the implied constant depends on $\varepsilon$, $F_v$ and $\varphi_v$. Hence, 
\begin{equation}\label{5.24}
\int_{F_v}\int_{F_v^{\times}}\overline{\mathcal{D}}\cdot W_{\pi_v'}\left(\begin{pmatrix}
a_v\\
&1
\end{pmatrix}\right)\frac{d^{\times}a_v\mathcal{I}(c_v,s)dc_v}{|a_v|_v^{1+it-s}}\ll 
C_v(\pi')^{1+\varepsilon}\textbf{C}_v^{\, -\frac{3}{2}}.
\end{equation}

Similarly, we have
\begin{equation}\label{5.25}
\int_{F_v}\int_{F_v^{\times}}\mathcal{D}' \cdot W_{\pi_v}\left(\begin{pmatrix}
a_v\\
&1
\end{pmatrix}\right)\frac{d^{\times}a_v\mathcal{I}(c_v,s)dc_v}{|a_v|_v^{1+it-s}}\ll 
C_v(\pi)^{1+\varepsilon}\textbf{C}_v^{\, -\frac{3}{2}}.
\end{equation}
 
To handle the remaining integral, we truncate $a_v$ into the ranges: $|a_v|_v\leq \textbf{C}_v^{\, 1+\varepsilon}$ and $|a_v|_v>\textbf{C}_v^{\, 1+\varepsilon}$. By Lemma \ref{lem5.6} and Lemma \ref{lem5.7},
\begin{align*}
\int_{|a_v|_v\leq \textbf{C}_v^{\, 1+\varepsilon}}\frac{\overline{\mathcal{D}}\cdot \mathcal{D}'}{|a_v|_v^{1+it-s}}d^{\times}a_v\ll |c_v|_v^2(C_v(\pi')C_v(\pi))^{1+\varepsilon}\int_{|a_v|_v\leq\textbf{C}_v^{\, 1+\varepsilon}}|a_v|_v^{s-2\vartheta-it}d^{\times}a_v,
\end{align*}
which converges absolutely in the range $\Re(s-2\vartheta-it)>0$, i.e., $\Re(s)>2\vartheta$. In particular, when $s=1/2+it$, we have
\begin{equation}\label{5.26}
\int_{F_v}\int_{|a_v|_v\leq \textbf{C}_v^{\, 1+\varepsilon}}\frac{\overline{\mathcal{D}}\cdot \mathcal{D}'\cdot \mathcal{I}(c_v,s)}{|a_v|_v^{1+it-s}}d^{\times}a_vdc_v\ll \frac{(C_v(\pi')C_v(\pi))^{1+\varepsilon}\textbf{C}_v^{\, 1/2-2\vartheta+\varepsilon}}{\textbf{C}_v^3},
\end{equation}
where the implied constant depends on $\varepsilon$, $F_v$ and $\varphi_v$. 

Suppose $|a_v|_v\geq \textbf{C}_v^{\, 1+\varepsilon}$. Substituting the constraint $|c_v|_v\leq \varepsilon^*\textbf{C}_v^{\, -1}$ into Lemma \ref{lemma5.7} and Lemma \ref{lem5.9}, we obtain 
\begin{align*}
\mathcal{D}\ll C_v(\pi)^{l+50}|a_v|_v^{-l},\ \ \ \mathcal{D}'\ll C_v(\pi')^{l+50}|a_v|_v^{-l},
\end{align*}
where the implied constant depends on $\varphi_v$, $F_v$ and $l$. Therefore,
\begin{align*}
\int_{|a_v|_v>\textbf{C}_v^{\, 1+\varepsilon}}\frac{\overline{\mathcal{D}}\cdot \mathcal{D}'}{|a_v|_v^{1+it-s}}d^{\times}a_v\ll (C_v(\pi')C_v(\pi))^{l+50}\int_{|a_v|_v>\textbf{C}_v^{\, 1+\varepsilon}}|a_v|_v^{\Re(s)-1-2l}d^{\times}a_v,
\end{align*}
which converges absolutely in $\Re(s)-1-2l<0$, i.e., $\Re(s)<1+2l$. By enlarging $l$ we obtain the meromorphic continuation of $\int_{|a_v|_v>\textbf{C}_v^{\, 1+\varepsilon}}\overline{\mathcal{D}}\cdot \mathcal{D}\cdot |a_v|_v^{s-1-it}d^{\times}a_v$ to $s\in \mathbb{C}$, with the bound (by taking $l>\floor{200\varepsilon^{-1}}$)
\begin{equation}\label{5.27}
\int_{|a_v|_v>\textbf{C}_v^{\, 1+\varepsilon}}\overline{\mathcal{D}}\cdot \mathcal{D}'\cdot |a_v|_v^{s-1-it}d^{\times}a_v\ll \textbf{C}_v^{\, -50}
\end{equation}
uniformly holding for $\Re(s)\geq -1$. 

Therefore, gathering together the above discussions, we conclude that the local integral $\Psi_v(s,\overline{W_{\pi_v}},W_{\pi_v'})$ admits a meromorphic continuation to $\Re(s)>2\vartheta$. Moreover, combining the estimates \eqref{5.24}, \eqref{5.25}, \eqref{5.26}, and \eqref{5.27}, we obtain 
\begin{equation}\label{5.32}
\Psi_v(1/2+it,\overline{W_{\pi_v}},W_{\pi_v'})=\mathcal{J}\cdot \int_{F_v}\mathcal{I}(c_v,1/2+it)dc_v,
\end{equation}
where the factor $\mathcal{J}$ is defined by 
\begin{align*}
\int_{F_v^{\times}}\overline{W_{\pi_v,t}\left(\begin{pmatrix}
a_v\\
&1
\end{pmatrix}\right)}W_{\pi_v'}\left(\begin{pmatrix}
a_v\\
&1
\end{pmatrix}\right)|a_v|_v^{-\frac{1}{2}}d^{\times}a_v=\int_{F_v^{\times}}\frac{|\varphi_v(a_v)|^2}{\sqrt{|a_v|_v}}d^{\times}a_v.
\end{align*}

By the definition of $\varphi_v$, we have $\mathcal{J}\gg 1$. Therefore, by \eqref{eq5.28} and \eqref{5.32}, 
\begin{equation}\label{5.34}
\Psi_v(1/2+it,\overline{W_{\pi_v}},W_{\pi_v'})\gg \int_{F_v}\int_{F_v}\frac{h_v(\textbf{C}_v|z_vc_v|)h_v(|z_v|-1)}{\textbf{C}_v^{\, -1/2}}dz_vdc_v,
\end{equation}
where the implied constant depends only on $\varepsilon$, $F_v$, and $\varphi_v$. According the construction of $h_v$ in \textsection\ref{sec4.4}, we have from \eqref{5.34} that $\Psi_v(1/2+it,\overline{W_{\pi_v}},W_{\pi_v'})\gg \textbf{C}_v^{\, -1/2}$, which is \eqref{f5.26}.
\end{proof}

\subsection{Proof of Proposition \ref{prop5.3}}
\begin{proof}[Proof of Proposition \ref{prop5.3}]
For $v<\infty$, we take $W_{\pi_v}\not \equiv 0$ to be a local new vector in the Whittaker model of $\pi_v$. For $v\mid\infty$, we take $W_{\pi_v}$ as defined in Proposition \ref{prop5.7}. Let $\varphi$ be the cusp form in $\pi$ corresponding to $W_{\pi}=\otimes_vW_{\pi_v}$. Let $\Re(s)\gg 1$. By unfolding the Eisenstein series $E(\cdot,s)$, 
\begin{align*}
\langle \phi E(\cdot,s),\varphi\rangle=\int_{N(\mathbb{A}_F)\backslash G(\mathbb{A}_F)}\overline{W_{\pi}(x)}W_{\pi'}(x)\Phi(\eta x)|\det x|^{s}dx,
\end{align*}
where $\eta=(0,1)$, and $\Phi=\otimes_v\Phi_v\in \mathcal{S}(\mathbb{A}_F)$ has been defined in \textsection\ref{sec4.4}. Therefore, 
\begin{equation}\label{5.5.}
\langle \phi E(\cdot,s),\varphi\rangle=\prod_{v\in\Sigma_F}\Psi_v(s,\overline{W_{\pi_v}},W_{\pi_v'}),
\end{equation}
where each factor $\Psi_v(s,\overline{W_{\pi_v}},W_{\pi_v'})$ is defined by 
\begin{align*}
\int_{N(F_v)\backslash G(F_v)}\overline{W_{\pi_v}(x_v)}W_{\pi_v'}(x_v)\Phi_v(\eta x_v)|\det x_v|_v^{s}dx_v.
\end{align*}

\begin{itemize}
\item Suppose that $v\mid\mathfrak{M}\mathfrak{N}$. Recall \eqref{equa4.5}:
\begin{align*}
\Phi_v(t_{1,v},t_{2,v}):=\Vol(K_v[m_v])^{-1}\omega_v\omega_v'^{-1}(t_{2,v})\mathbf{1}_{\mathfrak{p}_v^{m_v}}(t_{1,v})\mathbf{1}_{\mathcal{O}_v^{\times}}(t_{2,v}),
\end{align*}
where $m_v:=\max\{r_{\pi_v},r_{\pi_v'}\}$. Using the Iwasawa decomposition and the explicit formula for local new vectors we have 
\begin{equation}\label{eq5.9}
\Psi_v(s,\overline{W_{\pi_v}},W_{\pi_v'})=\frac{\overline{W_{\pi_v}}(I_2)L_v(s,\overline{\pi}_v\times\pi_v')}{L_v(2s,\omega_v^{-1}\omega_v')}.
\end{equation}

\item Suppose that $v<\infty$, and $v\nmid\mathfrak{M}\mathfrak{N}$. Recall \eqref{equ4.9}:
\begin{align*}
\Phi_v(t_{1,v},t_{2,v}):=\mathbf{1}_{\mathcal{O}_v}(t_{1,v})\mathbf{1}_{\mathcal{O}_v}(t_{2,v}).
\end{align*}
By Casselman-Shalika formula we have 
\begin{equation}\label{eq5.9.}
\Psi_v(s,\overline{W_{\pi_v}},W_{\pi_v'})=\overline{W_{\pi_v}}(I_2)L_v(s,\overline{\pi}_v\times\pi_v').
\end{equation}

\item Suppose that $v\mid\infty$. By Proposition \ref{prop5.7},
\begin{equation}\label{5.7}
\Psi_v(1/2+it,\overline{W_{\pi_v}},W_{\pi_v'})\gg \textbf{C}_{\infty}(\pi,t,\pi')^{-1/2}\cdot |\overline{W_{\pi_v}}(I_2)W_{\pi_v'}(I_2)|.
\end{equation}
\end{itemize}

By \eqref{5.5.}, \eqref{eq5.9} and \eqref{eq5.9.}, $\langle \phi E(\cdot,s),\varphi\rangle$ is equal to 
\begin{align*}
\frac{L(s,\overline{\pi}\times\pi')\prod_{v\mid\mathfrak{M}\mathfrak{N}}L_v(2s,\omega_v^{-1}\omega_v')^{-1}}{\prod_{v\in \mathcal{P}}L_v(2s,\omega_v^{-1}\omega_v')}\prod_{v\mid\infty}\Psi_v(s,\overline{W_{\pi_v}},W_{\pi_v'})\prod_{v<\infty}\overline{W_{\pi_v}(I_2)}.
\end{align*}
In conjunction with \eqref{5.7} we derive that 
\begin{equation}\label{5.9}
\big|\langle \phi E(\cdot,1/2+it),\varphi\rangle\big|\gg \frac{|L(1/2+it,\overline{\pi}\times\pi')||W_{\pi}(I_2)|}{\textbf{C}_{\infty}(\pi,t,\pi')^{1/2}\cdot\mathcal{E}},
\end{equation}
where 
\begin{align*}
\mathcal{E}:=\prod_{v\in \mathcal{P}}\Vol(K_v[0])^{-1}\prod_{v\in \mathcal{P}}|L_v(1+2it,\omega_v^{-1}\omega_v')|\prod_{v\mid\mathfrak{M}\mathfrak{N}}|L_v(1+2it,\omega_v^{-1}\omega_v')|.
\end{align*}

By triangle inequality and Mertens formula, 
\begin{equation}\label{eq5.12}
\prod_{\substack{v\in \mathcal{P}\\ \text{or $v\mid\mathfrak{M}\mathfrak{N}$}}}|L_v(1+2it,\omega_v^{-1}\omega_v')|\leq \prod_{q_v\ll LC_{\fin}(\pi,\pi')}(1-q_v^{-1})^{-1}\ll \log LC_{\fin}(\pi,\pi'),
\end{equation}
where the implied constant depends only on $F$. 

Let $\varphi^{\circ}:=\varphi/\langle\varphi,\varphi\rangle$. By \eqref{5.9}, \eqref{eq5.12}, and \cite[Corollary 0.3]{HL94} we obtain 
\begin{equation}\label{5.10}
\langle \phi E(\cdot,1/2+it),\varphi^{\circ}\rangle\gg \frac{|L(1/2+it,\overline{\pi}\times\pi')|}{\textbf{C}_{\infty}(\pi,t,\pi')^{1/2+\varepsilon}C_{\fin}(\pi,\pi')^{\varepsilon}\log L},
\end{equation}
where the implied constant depends on $F$ and $\varepsilon$.

Notice that $\pi\in \mathcal{A}_0([G],\omega)$ and $\varphi^{\circ}\in \mathfrak{B}(\pi)$ occurs in the cuspidal spectrum in \eqref{eq5.6}, and $\mathcal{M}_{\Eis}^{\tw}(1/2,1/2;\boldsymbol{\alpha},\boldsymbol{\ell})\geq 0$. Then Proposition \ref{prop5.3} follows from \eqref{5.10} and Proposition \ref{prop5.2}.
\end{proof}

\section{Estimates of the Dual Moment: the Constant Part}\label{sec6}
Let $\boldsymbol{\alpha}=(\alpha_v)_{v\in \mathcal{P}}$ and $\boldsymbol{\ell}=(\ell_v)_{v\in \mathcal{P}}$ be defined as in \textsection\ref{sec4.2}. Let $\varepsilon_0$ be the parameter defined in \textsection\ref{sect4.4}. 
Recall that the amplified constant part of the dual side is  
\begin{align*}
\mathcal{M}_{\const}^{\du,\heartsuit,\Reg}(s_1,s_2;\boldsymbol{\alpha},\boldsymbol{\ell}):=\frac{1}{2\pi i}\int_{|s|=\varepsilon_0}\frac{\mathcal{M}_{\const}^{\du,\heartsuit}(s_1-s,s_2+2s;\boldsymbol{\alpha},\boldsymbol{\ell})}{s}ds,\tag{\ref{eq4.5}}
\end{align*}
where $\mathcal{M}_{\const}^{\du,\heartsuit}(s_1-s,s_2+2s;\boldsymbol{\alpha},\boldsymbol{\ell})$ is defined as in \textsection\ref{sec4.6}. 

The main result in this section is the following estimates.

\begin{prop}\label{prop6.1}
Let notation be as before. Let $\varepsilon>0$. Then    
\begin{align*}
\mathcal{M}_{\const}^{\du,\heartsuit,\Reg}(1/2+it,1/2-it;\boldsymbol{\alpha},\boldsymbol{\ell})\ll 
C_{\fin}(\pi,\pi')\Delta\sum_{\substack{v\in \mathcal{P}}}|\alpha_v|_v^2,
\end{align*}
where $\Delta:=C_{\infty}(\pi\otimes|\cdot|^{it})^{\varepsilon}C_{\infty}(\pi')^{\varepsilon}C_{\infty}(\omega\omega'^{-1} |\cdot|^{2it})^{\varepsilon}C_{\fin}(\pi,\pi')^{\varepsilon}L^{\varepsilon}$. Here the implied constant depends only on $\varepsilon$ and the base field $F$.   
\end{prop}

The proof of Proposition \ref{prop6.1} will be provided in \textsection\ref{sec6.3} below, utilizing the estimates established in \textsection\ref{sec6.1}--\textsection\ref{sec6.2}.

\subsection{Decomposition of $\mathcal{M}_{\const}^{\du,\heartsuit,\Reg}(s_1,s_2;\mathfrak{L})$}
By definition, $\mathcal{M}_{\const}^{\du,\heartsuit,\Reg}(s_1,s_2;\boldsymbol{\alpha},\boldsymbol{\ell})$ is a linear combination of 
\begin{align*}
\mathcal{M}_{\const}^{\du,\heartsuit,\Reg}(s_1,s_2;\mathfrak{L}):=\frac{1}{2\pi i}\int_{|s|=\varepsilon_0}\frac{\mathcal{M}_{\const}^{\du,\heartsuit}(s_1-s,s_2+2s;\mathfrak{L})}{s}ds,
\end{align*}
where $\mathfrak{L}=\prod_{v\in \mathcal{P}}\mathfrak{p}_v^{\ell_v'}$ for some $0\leq \ell_v'\leq 2\ell_v$, $v\in \mathcal{P}$, and $\mathcal{M}_{\const}^{\du,\heartsuit}(s_1-s,s_2+2s;\mathfrak{L})$ is defined as in \textsection\ref{sec4.5}. By  \eqref{3.5..} in \textsection\ref{sec3.2}, we have
\begin{align*}
\mathcal{M}_{\const}^{\du,\heartsuit}(s_1,s_2;\mathfrak{L})=\mathcal{M}_{\const,1}^{\du,\heartsuit}(s_1,s_2;\mathfrak{L})+\mathcal{M}_{\const,2}^{\du,\heartsuit}(s_1,s_2;\mathfrak{L}).
\end{align*}
As a result, we can break down $\mathcal{M}_{\const}^{\du,\heartsuit,\Reg}(s_1,s_2;\mathfrak{L})$ as 
\begin{equation}\label{equa6.1}
\mathcal{M}_{\const}^{\du,\heartsuit,\Reg}(s_1,s_2;\mathfrak{L})=\mathcal{M}_{\const,1}^{\du,\heartsuit,\Reg}(s_1,s_2;\mathfrak{L})+\mathcal{M}_{\const,2}^{\du,\heartsuit,\Reg}(s_1,s_2;\mathfrak{L}),
\end{equation}
where for $j=1,2$, the function $\mathcal{M}_{\const,j}^{\du,\heartsuit,\Reg}(s_1,s_2;\mathfrak{L})$ is defined by
\begin{equation}\label{equ6.2}
\frac{1}{2\pi i}\int_{|s|=\varepsilon_0}\frac{\mathcal{M}_{\const,j}^{\du,\heartsuit}(s_1-s,s_2+2s;\mathfrak{L})}{s}ds.
\end{equation}

We shall investigate the integrals $\mathcal{M}_{\const,j}^{\du,\heartsuit,\Reg}(1/2+it,1/2-it;\mathfrak{L})$, $j=1,2$.
\subsection{The Constant Part: $\mathcal{M}_{\const,1}^{\du,\heartsuit,\Reg}(1/2+it,1/2-it;\mathfrak{L})$}\label{sec6.1}
In this section we aim to establish an upper bound for $\mathcal{M}_{\const,1}^{\du,\heartsuit,\Reg}(1/2+it,1/2-it;\mathfrak{L})$. See Proposition \ref{prop6.4} in \textsection\ref{sec6.1.3}.

\subsubsection{Shifted Rankin--Selberg convolutions}
Let $\mathcal{P}$ be defined by \eqref{equ4.2} in \textsection\ref{sec4.2}. Let $v\in \mathcal{P}$. Then $\pi_v'$ is unramified. 
Let $L_v(s,\pi_v')=(1-\alpha_v q_v^{-s})^{-1}(1-\beta_vq_v^{-s})^{-1}$ be the local $L$-factor of the standard $L$-function of $\pi_v'$. Then $\alpha_v\beta_v=\omega_v'(\varpi_v)$, and $\alpha_v+\beta_v=\lambda_{\pi'}(\mathfrak{p}_v)$. 

\begin{lemma}\label{lem6.2}
Let $v\in\mathcal{P}$. Let $\ell\geq 0$ and $i\geq 0$. Let $\Re(s)>-100^{-1}$. Then 
\begin{equation}\label{eq6.2}
\sum_{m\geq i}W_{\pi_v'}\left(\begin{pmatrix}
\varpi_v^{m+\ell}\\
&1
\end{pmatrix}\right)\overline{W_{\pi_v'}\left(\begin{pmatrix}
\varpi_v^{m}\\
&1
\end{pmatrix}\right)}q_v^{-ms}=\frac{H_vL_v(1+s,\pi_v'\times\overline{\pi}_v')}{q_v^{\ell/2+(1+s)i}},
\end{equation}
where $m$ ranges over nonnegative integers, and $H_v=H_v(s,\ell,i)$ is defined by   
\begin{equation}\label{eq6.3}
\lambda_{\pi'}(\mathfrak{p}_v^{\ell+i})\overline{\lambda_{\pi'}(\mathfrak{p}_v^i)}-c_1q_v^{-1-s}+c_2q_v^{-2-2s}-\lambda_{\pi'}(\mathfrak{p}_v^{\ell+i-1})\overline{\lambda_{\pi'}(\mathfrak{p}_v^{i-1})}q_v^{-3-3s},
\end{equation}
with 
\begin{align*}
c_1:=&\omega_v'^{-1}(\varpi_v)\lambda_{\pi'}(\mathfrak{p}_v^{\ell+i+1})\overline{\lambda_{\pi'}(\mathfrak{p}_v^{i-1})}+\omega_v'(\varpi_v)\lambda_{\pi'}(\mathfrak{p}_v^{\ell+i-1})\overline{\lambda_{\pi'}(\mathfrak{p}_v)}\overline{\lambda_{\pi'}(\mathfrak{p}_v^i)},\\
c_2:=&\lambda_{\pi'}(\mathfrak{p}_v^{\ell+i})\overline{\lambda_{\pi'}(\mathfrak{p}_v)}\overline{\lambda_{\pi'}(\mathfrak{p}_v^{i-1})}+\omega_v'(\varpi_v)\lambda_{\pi'}(\mathfrak{p}_v^{\ell+i-2})\overline{\lambda_{\pi'}(\mathfrak{p}_v^{i})}.
\end{align*} 
Here we set $\lambda_{\pi'}(\mathfrak{p}_v^{0})=1$, $\lambda_{\pi'}(\mathfrak{p}_v^{-1})=0$, and $\lambda_{\pi'}(\mathfrak{p}_v^{-2})=-\omega_v^{-1}(\varpi_v)$. 
\end{lemma}
\begin{proof}
For simplicity we write $\alpha=\alpha_v$ and $\beta=\beta_v$. By Casselman-Shalika formula, 
\begin{align*}
W_{\pi_v'}\left(\begin{pmatrix}
\varpi_v^n\\
&1
\end{pmatrix}\right)=\frac{\alpha^{n+1}-\beta^{n+1}}{\alpha-\beta}\cdot q_v^{-\frac{n}{2}}W_{\pi_v'}(I_2),\ \ n\geq 0.
\end{align*} 

Note that $W_{\pi_v'}(I_2)=1$. Therefore, the left hand side of \eqref{eq6.2} is equal to 
\begin{align*}
&q_v^{-\frac{\ell}{2}}\sum_{m\geq i}\frac{(\alpha^{m+\ell+1}-\beta^{m+\ell+1})(\overline{\alpha}^{m+1}-\overline{\beta}^{m+1})}{(\alpha-\beta)(\overline{\alpha}-\overline{\beta})}\cdot q_v^{-m-ms}\\
=&q_v^{-\frac{\ell}{2}}\sum_{m\geq i}\frac{(\alpha^{\ell+1}\overline{\alpha}(\alpha\overline{\alpha})^m+\beta^{\ell+1}\overline{\beta}(\beta\overline{\beta})^m-\alpha^{\ell+1}\overline{\beta}(\alpha\overline{\beta})^m-\beta^{\ell+1}\overline{\alpha}(\beta\overline{\alpha})^m)}{(\alpha-\beta)(\overline{\alpha}-\overline{\beta})q_v^{m(1+s)}}.
\end{align*}

Therefore,
\begin{equation}\label{6.5.}
\sum_{m\geq i}W_{\pi_v'}\left(\begin{pmatrix}
\varpi_v^{m+\ell}\\
&1
\end{pmatrix}\right)\overline{W_{\pi_v'}\left(\begin{pmatrix}
\varpi_v^{m}\\
&1
\end{pmatrix}\right)}q_v^{-ms}=q_v^{-\frac{\ell}{2}-(1+s)i}S(\alpha,\beta),
\end{equation}
where $S(\alpha,\beta)$ is defined by  
\begin{align*}
\frac{1}{(\alpha-\beta)(\overline{\alpha}-\overline{\beta})}\cdot \Bigg[\frac{\alpha^{\ell+1}\overline{\alpha}(\alpha\overline{\alpha})^i}{1-\alpha\overline{\alpha}q_v^{-1-s}}+\frac{\beta^{\ell+1}\overline{\beta}(\beta\overline{\beta})^i}{1-\beta\overline{\beta}q_v^{-1-s}}-\frac{\alpha^{\ell+1}\overline{\beta}(\alpha\overline{\beta})^i}{1-\alpha\overline{\beta}q_v^{-1-s}}-\frac{\beta^{\ell+1}\overline{\alpha}(\beta\overline{\alpha})^i}{1-\beta\overline{\alpha}q_v^{-1-s}}\Bigg].
\end{align*}

By a brute force calculation, the factor
\begin{align*}
\frac{\alpha^{\ell+1}\overline{\alpha}(\alpha\overline{\alpha})^i}{1-\alpha\overline{\alpha}q_v^{-1-s}}+\frac{\beta^{\ell+1}\overline{\beta}(\beta\overline{\beta})^i}{1-\beta\overline{\beta}q_v^{-1-s}}-\frac{\alpha^{\ell+1}\overline{\beta}(\alpha\overline{\beta})^i}{1-\alpha\overline{\beta}q_v^{-1-s}}-\frac{\beta^{\ell+1}\overline{\alpha}(\beta\overline{\alpha})^i}{1-\beta\overline{\alpha}q_v^{-1-s}}
\end{align*}
in $S(\alpha,\beta)$ is equal to the product of $L_v(1+s,\pi_v'\times\overline{\pi}_v')$ with 
\begin{align*}
(\alpha^{\ell+i+1}-\beta^{\ell+i+1})(\overline{\alpha}^{i+1}-\overline{\beta}^{i+1})-c_1'q_v^{-1-s}+c_2'q_v^{-2-2s}-c_3'q_v^{-3-3s},	
\end{align*}
where 
\begin{align*}
c_1':=&\overline{\alpha}\overline{\beta}(\alpha^{\ell+i+2}-\beta^{\ell+i+2})(\overline{\alpha}^{i}-\overline{\beta}^{i})+\alpha\beta(\overline{\alpha}+\overline{\beta})(\alpha^{\ell+i}-\beta^{\ell+i})(\overline{\alpha}^{i+1}-\overline{\beta}^{i+1}),\\
c_2':=&(\alpha^{\ell+i+1}-\beta^{\ell+i+1})(\overline{\alpha}^i-\overline{\beta}^i)(\overline{\alpha}+\overline{\beta})+\alpha\beta(\alpha^{\ell+i-1}-\beta^{\ell+i-1})(\overline{\alpha}^{i+1}-\overline{\beta}^{i+1}),
\end{align*}
and $c_3':=(\alpha^{\ell+i}-\beta^{\ell+i})(\overline{\alpha}^{i}-\overline{\beta}^{i})$. As a consequence, we obtain 
\begin{equation}\label{6.5}
S(\alpha,\beta)=\Bigg[\lambda_{\pi'}(\mathfrak{p}_v^{\ell+i})\overline{\lambda_{\pi'}(\mathfrak{p}_v^i)}-\frac{c_1'}{q_v^{1+s}}+\frac{c_2'}{q_v^{2+2s}}-\frac{\lambda_{\pi'}(\mathfrak{p}_v^{\ell+i-1})\overline{\lambda_{\pi'}(\mathfrak{p}_v^{i-1})}}{q_v^{3+3s}}\Bigg]\cdot *,
\end{equation}
where the factor $*$ represents the local $L$-factor $L_v(1+s,\pi_v'\times\overline{\pi}_v')$.

Therefore, \eqref{eq6.2} follows from \eqref{6.5.} and \eqref{6.5}. 
\end{proof}

\begin{comment}
\textcolor{blue}{By the Ramanujan bound  \eqref{eq6.6} in \textsection\ref{sec2.2.2},  $H_v(s,\ell,i)$ has the following asymptotic behavior}: 
\begin{equation}\label{H}
\begin{cases}
H_v(s,0)=\zeta_v(2+2s)^{-1}=1-q_v^{-2-2s},\\
H_v(s,1)=\lambda_{\pi'}(\mathfrak{p}_v)+O(q_v^{-1+\theta-\Re(s)}),\\
H_v(s,\ell)=\lambda_{\pi'}(\mathfrak{p}_v^{\ell})+O(q_v^{-1+\ell\theta-\Re(s)}),\ \ \ell\geq 2.
\end{cases}
\end{equation}
\end{comment}

\begin{lemma}\label{lemma6.3}
Let notation be as before. Let $v\in \mathcal{P}$. Let $\ell\geq 0$ and $i\geq 0$. Let $s_1, s_2\in \mathbb{C}$. Define 
\begin{align*}
\mathcal{RS}_1(s_1,s_2):=&-\Vol(K_v[\ell])\sum_{j=1}^{\ell}\omega_v^{-1}\omega_v'(\varpi_v^j)q_v^{2s_2j}\sum_{m\geq 0}W_{\pi_v'}\left(\begin{pmatrix}
	\varpi_v^m\\
	&1
\end{pmatrix}\right)\\
&\sum_{i=j}^{\ell}q_v^{\ell-i-1}\mathbf{1}_{m=i-1}\overline{W_{\pi_v'}\left(\begin{pmatrix}
	\varpi_v^{m+\ell-i}\\
& \varpi_v^{i}
\end{pmatrix}
\right)}q_v^{-(s_1+s_2-1)m}.
\end{align*}
Then 
\begin{align*}
\mathcal{RS}_1(s_1,s_2)=&-\Vol(K_v[\ell])q_v^{\ell-2}\omega_v^{-1}(\varpi_v)\omega_v'^{-1}(\varpi_v^{\ell-2})q_v^{-\frac{\ell-2}{2}}q_v^{2s_2}\sum_{0\leq m\leq \ell-2}q_v^{-(s_1+s_2)m}\\
&\sum_{0\leq j\leq \min\{m,\ell-2-m\}}\lambda_{\pi'}(\mathfrak{p}_v^{\ell-2-2j})\cdot \frac{1-\omega_v^{-1}\omega_v'(\varpi_v^{m+1})q_v^{2s_2(m+1)}}{1-\omega_v^{-1}\omega_v'(\varpi_v)q_v^{2s_2}}.
\end{align*}
\end{lemma}
\begin{proof}
By definition, 
\begin{align*}
\mathcal{RS}_1(s_1,s_2)=&-\Vol(K_v[\ell])q_v^{\ell-2}\sum_{j=1}^{\ell}\omega_v^{-1}\omega_v'(\varpi_v^j)q_v^{2s_2j}\sum_{j-1\leq m\leq \ell-1}W_{\pi_v'}\left(\begin{pmatrix}
	\varpi_v^m\\
	&1
\end{pmatrix}\right)\\
&\overline{W_{\pi_v'}\left(\begin{pmatrix}
	\varpi_v^{\ell-1}\\
& \varpi_v^{m+1}
\end{pmatrix}
\right)}q_v^{-(s_1+s_2)m}.
\end{align*}

Swapping the sums and computing the inner geometric series yields
\begin{align*}
\mathcal{RS}_1(s_1,s_2)=&-\Vol(K_v[\ell])q_v^{\ell-2}\omega_v^{-1}\omega_v'(\varpi_v)q_v^{2s_2}\sum_{0\leq m\leq \ell-2}W_{\pi_v'}\left(\begin{pmatrix}
	\varpi_v^m\\
	&1
\end{pmatrix}\right)\\
&\overline{W_{\pi_v'}\left(\begin{pmatrix}
	\varpi_v^{\ell-1}\\
& \varpi_v^{m+1}
\end{pmatrix}
\right)}q_v^{-(s_1+s_2)m}\cdot \frac{1-\omega_v^{-1}\omega_v'(\varpi_v^{m+1})q_v^{2s_2(m+1)}}{1-\omega_v^{-1}\omega_v'(\varpi_v)q_v^{2s_2}}.
\end{align*}

For $m\geq 0$,  utilizing the  Casselman-Shalika formula, we have  
\begin{align*}
W_{\pi_v'}\left(\begin{pmatrix}
	\varpi_v^m\\
	&1
\end{pmatrix}\right)\overline{W_{\pi_v'}\left(\begin{pmatrix}
	\varpi_v^{\ell-1}\\
& \varpi_v^{m+1}
\end{pmatrix}
\right)}=\frac{\omega_v'^{-1}(\varpi_v^{\ell-1})\lambda_{\pi'}(\mathfrak{p}_v^m)\lambda_{\pi'}(\mathfrak{p}_v^{\ell-2-m})}{q_v^{(\ell-2)/2}},
\end{align*}
which, in conjunction with the Hecke relations \eqref{hecke.} and \eqref{eq2.3}, is equal to 
\begin{align*}
\omega_v'^{-1}(\varpi_v^{\ell-1})q_v^{-\frac{\ell-2}{2}}\sum_{0\leq j\leq \min\{m,\ell-2-m\}}\lambda_{\pi'}(\mathfrak{p}_v^{\ell-2-2j}).
\end{align*}

Bringing together the above discussions yields Lemma \ref{lemma6.3}.
\end{proof}

\begin{lemma}\label{lemma6.4}
Let notation be as before. Let $v\in \mathcal{P}$. Let $\ell\geq 0$ and $i\geq 0$. Let $s_1, s_2\in \mathbb{C}$ with $\Re(s_1+s_2)>1$. Define 
\begin{align*}
\mathcal{RS}_2(s_1,s_2):=&|\mathbb{F}_v^{\times}|\Vol(K_v[\ell])\sum_{j=1}^{\ell}\omega_v^{-1}\omega_v'(\varpi_v^j)q_v^{2s_2j}\sum_{m\geq 0}W_{\pi_v'}\left(\begin{pmatrix}
	\varpi_v^m\\
	&1
\end{pmatrix}\right)\\
&\sum_{i=j}^{\ell}q_v^{\ell-i-1}\mathbf{1}_{m\geq i}\overline{W_{\pi_v'}\left(\begin{pmatrix}
	\varpi_v^{m+\ell-i}\\
& \varpi_v^{i}
\end{pmatrix}
\right)}q_v^{-(s_1+s_2-1)m}.
\end{align*}
Then $\mathcal{RS}_2(s_1,s_2)$ admits a meromorphic continuation to $(s_1,s_2)\in \mathbb{C}^2$, given by 
\begin{align*}
\mathcal{RS}_2(s_1,s_2)=&|\mathbb{F}_v^{\times}|\Vol(K_v[\ell])\sum_{i=1}^{\ell}\omega_v'^{-1}(\varpi_v^i)q_v^{\ell-i-1}\sum_{j=1}^i\omega_v^{-1}\omega_v'(\varpi_v^j)q_v^{2s_2j}\\
&\Big[H_v^-(s_1,s_2,\ell,i)\mathbf{1}_{2i\leq \ell}+H_v^+(s_1,s_2,\ell,i)\mathbf{1}_{2i>\ell}\Big]\cdot \frac{L_v(s_1+s_2,\pi_v'\times\overline{\pi}_v')}{q_v^{\ell/2+(s_1+s_2-1)i}}
,
\end{align*}
where 
\begin{align*}
H_v^-(s_1,s_2,\ell,i):=&\overline{\lambda_{\pi'}(\mathfrak{p}_v^{\ell-i})}\lambda_{\pi'}(\mathfrak{p}_v^i)-\frac{c_1^{i,\ell}}{q_v^{s_1+s_2}}+\frac{c_2^{i,\ell}}{q_v^{2(s_1+s_2)}}-\frac{\overline{\lambda_{\pi'}(\mathfrak{p}_v^{\ell-i-1})}\lambda_{\pi'}(\mathfrak{p}_v^{i-1})}{q_v^{3(s_1+s_2)}},\\
H_v^+(s_1,s_2,\ell,i):=&\lambda_{\pi'}(\mathfrak{p}_v^{\ell-i})\overline{\lambda_{\pi'}(\mathfrak{p}_v^i)}-\frac{c_3^{i,\ell}}{q_v^{s_1+s_2}}+\frac{c_4^{i,\ell}}{q_v^{2(s_1+s_2)}}-\frac{\lambda_{\pi'}(\mathfrak{p}_v^{i-1})\overline{\lambda_{\pi'}(\mathfrak{p}_v^{\ell-i-1})}}{q_v^{3(s_1+s_2)}}.
\end{align*}
Here the coefficients $c_j^{i,\ell}$, $1\leq j\leq 4$, are defined by  
\begin{align*}
c_1^{i,\ell}:=&\omega_v'(\varpi_v)\overline{\lambda_{\pi'}(\mathfrak{p}_v^{\ell-i+1})}\lambda_{\pi'}(\mathfrak{p}_v^{i-1})+\omega_v'^{-1}(\varpi_v)\overline{\lambda_{\pi'}(\mathfrak{p}_v^{\ell-i-1})}\lambda_{\pi'}(\mathfrak{p}_v)\lambda_{\pi'}(\mathfrak{p}_v^i),\\
c_2^{i,\ell}:=&\overline{\lambda_{\pi'}(\mathfrak{p}_v^{\ell-i})}\lambda_{\pi'}(\mathfrak{p}_v)\lambda_{\pi'}(\mathfrak{p}_v^{i-1})+\omega_v'^{-1}(\varpi_v)\overline{\lambda_{\pi'}(\mathfrak{p}_v^{\ell-i-2})}\lambda_{\pi'}(\mathfrak{p}_v^{i}),\\
c_3^{i,\ell}:=&\omega_v'^{-1}(\varpi_v)\lambda_{\pi'}(\mathfrak{p}_v^{i+1})\overline{\lambda_{\pi'}(\mathfrak{p}_v^{\ell-i-1})}+\omega_v'(\varpi_v)\lambda_{\pi'}(\mathfrak{p}_v^{i-1})\overline{\lambda_{\pi'}(\mathfrak{p}_v)\lambda_{\pi'}(\mathfrak{p}_v^{\ell-i})},\\
c_4^{i,\ell}:=&\lambda_{\pi'}(\mathfrak{p}_v^{i})\overline{\lambda_{\pi'}(\mathfrak{p}_v)}\overline{\lambda_{\pi'}(\mathfrak{p}_v^{\ell-i-1})}+\omega_v'(\varpi_v)\lambda_{\pi'}(\mathfrak{p}_v^{i-2})\overline{\lambda_{\pi'}(\mathfrak{p}_v^{\ell-i})}.
\end{align*} 
\end{lemma}
\begin{proof}
By definition, we have 
\begin{align*}
\mathcal{RS}_2(s_1,s_2)=&|\mathbb{F}_v^{\times}|\Vol(K_v[\ell])\sum_{i=1}^{\ell}\omega_v'^{-1}(\varpi_v^i)q_v^{\ell-i-1}\sum_{j=1}^i\omega_v^{-1}\omega_v'(\varpi_v^j)q_v^{2s_2j}\\
&\sum_{m\geq i}W_{\pi_v'}\left(\begin{pmatrix}
	\varpi_v^m\\
	&1
\end{pmatrix}\right)\overline{W_{\pi_v'}\left(\begin{pmatrix}
	\varpi_v^{m+\ell-2i}\\
& 1
\end{pmatrix}
\right)}q_v^{-(s_1+s_2-1)m},
\end{align*}
which converges absolutely in the region $\Re(s_1+s_2)>1$.

Utilizing Lemma \ref{lem6.2}, the sum
\begin{align*}
\sum_{m\geq i}W_{\pi_v'}\left(\begin{pmatrix}
	\varpi_v^m\\
	&1
\end{pmatrix}\right)\overline{W_{\pi_v'}\left(\begin{pmatrix}
	\varpi_v^{m+\ell-2i}\\
& 1
\end{pmatrix}
\right)}q_v^{-m(s_1+s_2-1)}
\end{align*}
is equal to 
\begin{align*}
\Big[H_v^-(s_1,s_2,\ell,i)\mathbf{1}_{2i\leq \ell}+H_v^+(s_1,s_2,\ell,i)\mathbf{1}_{2i>\ell}\Big]\cdot \frac{L_v(s_1+s_2,\pi_v'\times\overline{\pi}_v')}{q_v^{\ell/2+(s_1+s_2-1)i}}.
\end{align*}
Therefore, Lemma \ref{lemma6.4} follows. 
\end{proof}

Let $v\mid\mathfrak{L}$. Let $s_1, s_2\in \mathbb{C}$ with $\Re(s_1+s_2)>1$. Define 
\begin{align*}
\mathcal{RS}(s_1,s_2)=&\sum_{j=1}^{\ell_v'}\omega_v^{-1}\omega_v'(\varpi_v^j)q_v^{2s_2j}\int_{K_v[j]}\int_{F_v^{\times}}W_{\pi_v'}\left(\begin{pmatrix}
	a_v\\
	&1
\end{pmatrix}k_v\right)\\
&\overline{W_{\pi_v'}\left(\begin{pmatrix}
	a_v\\
	&1
\end{pmatrix}k_v\begin{pmatrix}
1\\
&\varpi_v^{\ell_v'}
\end{pmatrix}
\right)}|a_v|_v^{s_1+s_2-1}d^{\times}a_vdk_v.
\end{align*}

\begin{lemma}\label{lemm6.5}
Let notation be as before. Then $\mathcal{RS}(s_1,s_2)$ admits a meromorphic continuation to $(s_1,s_2)\in \mathbb{C}^2$, given by 
\begin{equation}\label{eq6.7}
\mathcal{RS}(s_1,s_2)=\mathcal{RS}_1(s_1,s_2)+\mathcal{RS}_2(s_1,s_2),
\end{equation}
where $\mathcal{RS}_1(s_1,s_2)$ has been defined in Lemma \ref{lemma6.3}, and $\mathcal{RS}_2(s_1,s_2)$ has been defined in Lemma \ref{lemma6.4}. In particular, when $s_1=1/2+it-s$, and $s_2=1/2-it+2s$, where $|s|=\varepsilon_0$ is defined in \textsection\ref{sect4.4}, we have 
\begin{equation}\label{eq6.8}
\frac{\mathcal{RS}(s_1,s_2)}{L_v(s_1+s_2,\pi_v'\times\overline{\pi}_v')}\ll \ell^2 q_v^{-\frac{\ell}{2}+5\ell\varepsilon_0}\sum_{0\leq i\leq \ell}|\lambda_{\pi'}(\mathfrak{p}_v^{i})|,
\end{equation}
where the implied constant is absolute. 
\end{lemma}
\begin{proof}
Write $\ell=\ell_v'$. By Cramer's rule, we have 
\begin{align*}
K_v[j]=\sum_{\tau\in \mathfrak{p}_v^j/\mathfrak{p}_v^{\ell}}\begin{pmatrix}
1&\\
\tau & 1
\end{pmatrix}K_v[\ell].
\end{align*}
For $\tau \in \mathfrak{p}_v^j/\mathfrak{p}_v^{\ell}$, there is a unique $j\leq i\leq \ell$ such that $\tau=\varpi_v^{i}c_v$ with $c_v\in (\mathcal{O}_v/\mathfrak{p}_v^{\ell-i})^{\times}$. As a result, we have  
\begin{align*}
&\int_{K_v[j]}\overline{W_{\pi_v'}\left(\begin{pmatrix}
	a_v\\
	&1
\end{pmatrix}k_v\begin{pmatrix}
1\\
&\varpi_v^{\ell}
\end{pmatrix}
\right)}dk_v\\
=&\Vol(K_v[\ell])\sum_{i=j}^{\ell}\sum_{c_v\in (\mathcal{O}_v/\mathfrak{p}_v^{\ell-i})^{\times}}\overline{W_{\pi_v'}\left(\begin{pmatrix}
	a_v\\
	&1
\end{pmatrix}\begin{pmatrix}
1&\\
\varpi_v^{i}c_v& \varpi_v^{\ell}
\end{pmatrix}
\right)}.
\end{align*}
Here we identify $(\mathcal{O}_v/\mathfrak{p}_v^{\ell-i})^{\times}$ with the singleton $\{I_2\}$. Notice that 
\begin{equation}\label{eq6.9}
\begin{pmatrix}
1&\\
\varpi_v^{i}c_v& \varpi_v^{\ell}
\end{pmatrix}=\begin{pmatrix}
&1\\
\varpi_v^{\ell}& \varpi_v^{i}c_v
\end{pmatrix}w=\begin{pmatrix}
\varpi_v^{\ell-i}c_v^{-1}&1\\
& \varpi_v^{i}c_v
\end{pmatrix}w\begin{pmatrix}
1&\varpi_v^{\ell-i}c_v^{-1}\\
& -1\end{pmatrix},
\end{equation}
where $w=\begin{pmatrix}
&1\\
1
\end{pmatrix}$. Therefore, 
\begin{align*}
W_{\pi_v'}\left(\begin{pmatrix}
	a_v\\
	&1
\end{pmatrix}\begin{pmatrix}
1&\\
\varpi_v^{i}c_v& \varpi_v^{\ell}
\end{pmatrix}
\right)=W_{\pi_v'}\left(\begin{pmatrix}
	a_v\\
	&1
\end{pmatrix}\begin{pmatrix}
\varpi_v^{\ell-i}&1\\
& \varpi_v^{i}c_v
\end{pmatrix}
\right),
\end{align*}
which is equal to $\overline{\psi_v(c_v^{-1}a_v\varpi_v^{-i})}W_{\pi_v'}\left(\begin{pmatrix}
	a_v\varpi_v^{\ell-i}\\
& \varpi_v^{i}
\end{pmatrix}
\right)$. Consequently, 
\begin{align*}
\mathcal{RS}(s_1,s_2)=&\Vol(K_v[\ell])\sum_{j=1}^{\ell}\omega_v^{-1}\omega_v'(\varpi_v^j)q_v^{2s_2j}\sum_{m\geq 0}W_{\pi_v'}\left(\begin{pmatrix}
	\varpi_v^m\\
	&1
\end{pmatrix}\right)\\
&\sum_{i=j}^{\ell}\sum_{c_v\in (\mathcal{O}_v/\mathfrak{p}_v^{\ell-i})^{\times}}\frac{\psi_v(c_v\varpi_v^{m-i})}{q_v^{(s_1+s_2-1)m}}\cdot \overline{W_{\pi_v'}\left(\begin{pmatrix}
	\varpi_v^{m+\ell-i}\\
& \varpi_v^{i}
\end{pmatrix}
\right)}.
\end{align*}

Taking advantage of the Ramanujan sums
\begin{align*}
\sum_{c_v\in (\mathcal{O}_v/\mathfrak{p}_v^{\ell-i})^{\times}}\psi_v(c_v\varpi_v^{m-i})=q_v^{\ell-i-1}\cdot \big[-\mathbf{1}_{m=i-1}+|\mathbb{F}_v^{\times}|\cdot \mathbf{1}_{m\geq i}\big],
\end{align*}
we thus have $\mathcal{RS}(s_1,s_2)=\mathcal{RS}_1(s_1,s_2)+\mathcal{RS}_2(s_1,s_2)$, obtaining \eqref{eq6.7}.  

Let $s_1=1/2+it-s$, and $s_2=1/2-it+2s$, where $|s|=\varepsilon_0$. By Lemma \ref{lemma6.3},
\begin{equation}\label{6.9}
\frac{\mathcal{RS}_1(s_1,s_2)}{L_v(s_1+s_2,\pi_v'\times\overline{\pi}_v')}\ll \ell^2 q_v^{-\frac{\ell}{2}+5\ell\varepsilon_0}\mathbf{1}_{\ell\geq 2}\sum_{0\leq j\leq \ell/2-1}|\lambda_{\pi'}(\mathfrak{p}_v^{\ell-2-2j})|,
\end{equation}
where the implied constant is absolute. Moreover, by Lemma \ref{lemma6.4}, we obtain 
\begin{align*}
\frac{\mathcal{RS}_2(s_1,s_2)}{L_v(s_1+s_2,\pi_v'\times\overline{\pi}_v')}\ll \sum_{i=1}^{\ell}\Big[|H_v^-(s_1,s_2,\ell,i)|\mathbf{1}_{2i\leq \ell}+|H_v^+(s_1,s_2,\ell,i)|\mathbf{1}_{2i>\ell}\Big]\cdot q_v^{-\ell/2+\varepsilon_0i}
,
\end{align*}
where the implied constant is absolute. Making use of Hecke relations we obtain 
\begin{align*}
|H_v^-(s_1,s_2,\ell,i)|\mathbf{1}_{2i\leq \ell}+|H_v^+(s_1,s_2,\ell,i)|\mathbf{1}_{2i>\ell}\ll \sum_{0\leq i\leq \ell}|\lambda_{\pi'}(\mathfrak{p}_v^{i})|.
\end{align*}
Therefore, we derive the majorization 
\begin{equation}\label{6.10}
\frac{\mathcal{RS}_2(s_1,s_2)}{L_v(s_1+s_2,\pi_v'\times\overline{\pi}_v')}\ll q_v^{-\ell/2+\ell\varepsilon_0}\sum_{0\leq i\leq \ell}|\lambda_{\pi'}(\mathfrak{p}_v^{i})|,
\end{equation}
where the implied constant is absolute. 

Hence, the estimate \eqref{eq6.8} follows from \eqref{eq6.7}, \eqref{6.9} and \eqref{6.10}. 
\end{proof}

\subsubsection{The Contribution from $\mathcal{M}_{\const,1}^{\du,\heartsuit}(s_1,s_2)$}
Let $h=h_{\Phi,\mathbf{1},\omega\omega'^{-1}}(\cdot,s)$ be defined by \eqref{2.3}, and $\phi_{\mathfrak{L}}$ be defined in \eqref{eq4.1}. By definition in \textsection\ref{sec3.2}, we have  
\begin{equation}\label{6.2}
\mathcal{M}_{\const,1}^{\du,\heartsuit}(s_1,s_2;\mathfrak{L})=\gamma\cdot\int_{B(F)Z(\mathbb{A}_F)\backslash G(\mathbb{A}_F)}\phi(x)\overline{\phi(x\boldsymbol{d}_{\mathfrak{L}})}h(x,s_1)\overline{h(x\boldsymbol{d}_{\mathfrak{L}},\overline{s_2})}dx,
\end{equation}
where $\gamma:=\prod_{v\in\mathcal{P}}\omega_v(\varpi_v^{\ell_v'})(q_v^{\ell_v'/2}+\ell_v'q_v^{-\ell_v'/2}|\mathbb{F}_v^{\times}|)$,  $\boldsymbol{d}_{\mathfrak{L}}=\otimes_{v\in \mathcal{P}}\ \boldsymbol{d}_v\in G(\mathbb{A}_F)$, with $\boldsymbol{d}_v=\begin{pmatrix}
1&\\
&\varpi_v^{\ell_v'}
\end{pmatrix}$. We may write $\boldsymbol{d}_{\mathfrak{L}}=\otimes_{v\in \Sigma_F}\ \boldsymbol{d}_v\in G(\mathbb{A}_F)$, where $\boldsymbol{d}_v=\begin{pmatrix}
1&\\
&\varpi_v^{\ell_v'}
\end{pmatrix}$ with $\ell_v':=0$ if $v\not\in \mathcal{P}$. %We start with the local integral at $v\mid\mathfrak{L}.$

%\begin{lemma}
%Let notation be as before. Then
%\begin{align*}

%\end{align*}
%\end{lemma}

The following is a more concise version of Lemma \ref{lem3.6.}. 
\begin{lemma}\label{lem6.3}
Let notation be as before. Then 
\begin{align*}
\mathcal{M}_{\const,1}^{\du,\heartsuit}(s_1,s_2;\mathfrak{L})=\frac{L(s_1+s_2,\pi'\times\overline{\pi}')L(2s_1,\omega^{-1}\omega')L(2s_2,\omega\omega'^{-1})}{\zeta_F(2s_1+2s_2)}\prod_{v\in\Sigma_F}I_v^{\sharp}
\end{align*}
where the local factors $I_v^{\sharp}:=I_v^{\sharp}(s_1,s_2;\mathfrak{L})$ are defined as follows:
\begin{itemize}
\item at $v\mid\mathfrak{L}$, we define $I_v^{\sharp}(s_1,s_2;\mathfrak{L})$ by
\begin{align*}
\frac{(q_v^{\ell_v'/2}+\ell_v'q_v^{-\ell_v'/2}|\mathbb{F}_v^{\times}|)\omega_v(\varpi_v^{\ell_v'})\zeta_v(2s_1+2s_2)\mathcal{RS}(s_1,s_2)}{q_v^{s_2\ell_v'}L_v(2s_2,\omega_v\omega_v'^{-1})L_v(s_1+s_2,\pi_v'\times\overline{\pi}_v')}+\frac{\overline{\lambda_{\pi'}(\mathfrak{p}_v^{\ell})}\omega_v(\varpi_v^{\ell_v'})}{q_v^{s_2\ell_v'}},
\end{align*}
where the function $\mathcal{RS}(s_1,s_2)$ is defined right before Lemma \ref{lemm6.5};
\item at $v\mid\mathfrak{M}\mathfrak{N}$, with $m_v:=\max\{r_{\pi_v},r_{\pi_v'}\}$, we define  
\begin{align*}
I_v^{\sharp}(s_1,s_2;\mathfrak{L}):=\prod_{v\mid\mathfrak{M}\mathfrak{N}}\Vol(K_v[m_v])^{-1}L_v(2s_1,\omega_v^{-1}\omega_v')^{-1}L_v(2s_2,\omega_v\omega_v'^{-1})^{-1};
\end{align*}
\item and at archimedean places $v\mid\infty$, 
\begin{align*}
I_v^{\sharp}(s_1,s_2;\mathfrak{L}):=&\int_{K_v}\int_{F_v^{\times}}W_{\pi_v'}\left(\begin{pmatrix}
	a_v\\
	&1
\end{pmatrix}k_v\right)\overline{W_{\pi_v'}\left(\begin{pmatrix}
	a_v\\
	&1
\end{pmatrix}k_v\right)}\\
&h_v(k_v,s_1)\overline{h_v(k_v,\overline{s_2})}|a_v|_v^{s_1+s_2-1}d^{\times}a_vdk_v;
\end{align*}
\item at the remaining places, i.e., $v<\Sigma_{F,\fin}$ and $v\nmid\mathfrak{M}\mathfrak{N}$, $I_v^{\sharp}(s_1,s_2;\mathfrak{L})\equiv 1$. 
\end{itemize}
\end{lemma}
\begin{proof}
Let $\Re(s_1)\gg 1$ and $\Re(s_2)
\gg 1$. Executing Fourier expansion in \eqref{6.2} and using the Iwasawa decomposition we have
\begin{align*}
\mathcal{M}_{\const,1}^{\du,\heartsuit}(s_1,s_2;\mathfrak{L})=\prod_{v\in \Sigma_F}\Psi_{1,v}(s_1,s_2),
\end{align*}
where for each $v\in \Sigma_F$,, the factor $\Psi_{1,v}(s_1,s_2)$ is defined by 
\begin{align*}
\Psi_{1,v}(s_1,s_2):=&\omega_v(\varpi_v^{\ell_v'})\int_{K_v}\int_{F_v^{\times}}W_{\pi_v'}\left(\begin{pmatrix}
	a_v\\
	&1
\end{pmatrix}k_v\right)\overline{W_{\pi_v'}\left(\begin{pmatrix}
	a_v\\
	&1
\end{pmatrix}k_v\boldsymbol{d}_v\right)}\\
&(q_v^{\ell_v'/2}+\ell_v'q_v^{-\ell_v'/2}|\mathbb{F}_v^{\times}|)h_v(k_v,s_1)\overline{h_v(k_v\boldsymbol{d}_v,\overline{s_2})}|a_v|_v^{s_1+s_2-1}d^{\times}a_vdk_v.
\end{align*}

By definition, 
\begin{equation}\label{6.7}
h_v(k_v,s_1)=\int_{F_v^{\times}}\Phi_v((0,t_v)k_v)\omega_v^{-1}\omega_v'(t_v)|t_v|_v^{2s_1}d^{\times}t_v
\end{equation}
and 
\begin{equation}\label{6.8}
\overline{h_v(k_v\boldsymbol{d}_v,\overline{s_2})}=|\det \boldsymbol{d}_v|_v^{s_2}\int_{F_v^{\times}}\overline{\Phi_v((0,t_v)k_v\boldsymbol{d}_v)}\omega_v\omega_v'^{-1}(t_v)|t_v|_v^{2s_2}d^{\times}t_v.
\end{equation}

Let $v<\infty$. Write $k_v=\begin{pmatrix}
k_{11}& k_{12}\\
k_{21}& k_{22}
\end{pmatrix}$, and $m_v:=\max\{r_{\pi_v},r_{\pi_v'}\}$.  Making use of the definition of $\Phi_v(\cdot,\cdot)$ in \textsection\ref{sec4.4}, we have the following cases.

\begin{comment}

Then
\begin{itemize}
\item Suppose that $v\mid\mathfrak{L}$, 

\begin{align*}
h_v(k_v,s_1)=\int_{F_v^{\times}}\Phi_v((0,t_v)k_v)\omega_v^{-1}\omega_v'(t_v)|t_v|_v^{2s_1}d^{\times}t_v=\Vol(K_v[m_v])^{-1}\omega_v\omega_v'^{-1}(k_{22})\mathbf{1}_{k_{22}\in \mathcal{O}_v^{\times}}\mathbf{1}_{k_v\in K_v[m_v]},
\end{align*} 
\end{itemize}

where $m_v:=\max\{r_{\pi_v},r_{\pi_v'}\}$. Likewise, 
\begin{align*}
\overline{h_v(k_v\boldsymbol{d}_v,\overline{s_2})}=|\det \boldsymbol{d}_v|_v^{s_2}\int_{F_v^{\times}}\overline{\Phi_v((0,t_v)k_v\boldsymbol{d}_v)}\omega_v\omega_v'^{-1}(t_v)|t_v|_v^{2s_2}d^{\times}t_v,
\end{align*}
which is equal to 
\begin{align*}
q_v^{s_2\ell_v'}\Vol(K_v[m_v])^{-1}\mathbf{1}_{k_v\in K_v[\ell_v']}\mathbf{1}_{k_{22}\in \mathcal{O}_v^{\times}}\mathbf{1}_{k_v\in K_v[m_v]}\omega_v^{-1}\omega_v'(k_{22})\omega_v^{-1}\omega_v'(\varpi_v^{\ell_v'}).
\end{align*}

\begin{align*}
q_v^{s_2\ell_v'}\Vol(K_v[m_v])^{-2}\omega_v^{-1}\omega_v'(\varpi_v^{\ell_v'})\int_{K_v}\int_{F_v^{\times}}&W_{\pi_v'}\left(\begin{pmatrix}
	a_v\\
	&1
\end{pmatrix}k_v\right)\overline{W_{\pi_v'}\left(\begin{pmatrix}
	a_v\\
	&1
\end{pmatrix}k_v\boldsymbol{d}_v\right)}\\
&\mathbf{1}_{k_v\in K_v[\ell_v']}\mathbf{1}_{k_{22}\in \mathcal{O}_v^{\times}}\mathbf{1}_{k_v\in K_v[m_v]}|a_v|_v^{s_1+s_2-1}d^{\times}a_vdk_v.
\end{align*}
\end{comment}

\begin{itemize}
\item Suppose that $v\mid\mathfrak{L}$, i.e.,   $\ell_v'\geq 1$. Then $m_v=0$. Thus by \eqref{4.10} we have
\begin{align*}
\Psi_{1,v}(s_1,s_2)=&\delta_v^{(1)}\mathcal{RS}(s_1,s_2)+\delta_v^{(1)} L_v(2s_2,\omega_v\omega_v'^{-1})\Psi_{1,v}^{(0)}(s_1+s_2),
\end{align*}
where $\delta_v^{(1)}:=q_v^{-s_2\ell_v'}(q_v^{\ell_v'/2}+\ell_v'q_v^{-\ell_v'/2}|\mathbb{F}_v^{\times}|)\omega_v(\varpi_v^{\ell_v'})L_v(2s_1,\omega_v^{-1}\omega_v')$;  and 
\begin{align*}
\Psi_{1,v}^{(0)}(s):=\int_{K_v}\int_{F_v^{\times}}W_{\pi_v'}\left(\begin{pmatrix}
	a_v\\
	&1
\end{pmatrix}k_v\right)\overline{W_{\pi_v'}\left(\begin{pmatrix}
	a_v\\
	&1
\end{pmatrix}k_v\boldsymbol{d}_v\right)}|a_v|_v^{s-1}d^{\times}a_vdk_v.
\end{align*}

Since $W_{\pi_v'}$ is spherical, then 
\begin{align*}
\Psi_{1,v}^{(0)}(s)=\int_{F_v^{\times}}W_{\pi_v'}\left(\begin{pmatrix}
	a_v\\
	&1
\end{pmatrix}\right)\int_{K_v}\overline{W_{\pi_v'}\left(\begin{pmatrix}
	a_v\\
	&1
\end{pmatrix}k_v\begin{pmatrix}
1\\
&\varpi_v^{\ell}
\end{pmatrix}\right)}dk_v|a_v|_v^{s-1}d^{\times}a_v.
\end{align*}

Taking advantage of \eqref{5.2} and the standard calculation of Rankin--Selberg convolutions, we then obtain that $\delta_v^{(1)}L_v(2s_2,\omega_v\omega_v'^{-1})\Psi_{1,v}^{(0)}(s_1+s_2)$ becomes
\begin{align*}
\frac{\overline{\lambda_{\pi'}(\mathfrak{p}_v^{\ell})}\omega_v(\varpi_v^{\ell_v'})L_v(2s_1,\omega_v^{-1}\omega_v')L_v(2s_2,\omega_v\omega_v'^{-1})L_v(s_1+s_2,\pi_v'\times\overline{\pi}_v')}{q_v^{s_2\ell_v'}\zeta_v(2s_1+2s_2)}.
\end{align*}

Therefore, by Lemma \ref{lemm6.5} we derive 
\begin{align*}
\Psi_{1,v}(s_1,s_2)=&q_v^{-s_2\ell_v'}(q_v^{\ell_v'/2}+\ell_v'q_v^{-\ell_v'/2}|\mathbb{F}_v^{\times}|)\omega_v(\varpi_v^{\ell_v'})L_v(2s_1,\omega_v^{-1}\omega_v')\mathcal{RS}(s_1,s_2)\\
&+\frac{\overline{\lambda_{\pi'}(\mathfrak{p}_v^{\ell})}L_v(2s_1,\omega_v^{-1}\omega_v')L_v(2s_2,\omega_v\omega_v'^{-1})L_v(s_1+s_2,\pi_v'\times\overline{\pi}_v')}{\omega_v^{-1}(\varpi_v^{\ell_v'})q_v^{s_2\ell_v'}\zeta_v(2s_1+2s_2)}.
\end{align*}

As a result, $I_v^{\sharp}(s_1,s_2;\mathfrak{L})$ is given by 
\begin{align*}
\frac{(q_v^{\ell_v'/2}+\ell_v'q_v^{-\ell_v'/2}|\mathbb{F}_v^{\times}|)\omega_v(\varpi_v^{\ell_v'})\zeta_v(2s_1+2s_2)\mathcal{RS}(s_1,s_2)}{q_v^{s_2\ell_v'}L_v(2s_2,\omega_v\omega_v'^{-1})L_v(s_1+s_2,\pi_v'\times\overline{\pi}_v')}+\frac{\overline{\lambda_{\pi'}(\mathfrak{p}_v^{\ell})}\omega_v(\varpi_v^{\ell_v'})}{q_v^{s_2\ell_v'}}.
\end{align*}

\item Suppose $v\mid\mathfrak{M}\mathfrak{N}$, i.e.,  $m_v:=\max\{r_{\pi_v},r_{\pi_v'}\}\geq 1$. Then 
\begin{align*}
\Psi_{1,v}(s_1,s_2)=\int_{K_v}\int_{F_v^{\times}}W_{\pi_v'}&\left(\begin{pmatrix}
	a_v\\
	&1
\end{pmatrix}k_v\right)\overline{W_{\pi_v'}\left(\begin{pmatrix}
	a_v\\
	&1
\end{pmatrix}k_v\right)}\\
&h_v(k_v,s_1)\overline{h_v(k_v,\overline{s_2})}|a_v|_v^{s_1+s_2-1}d^{\times}a_vdk_v,
\end{align*}
which, by \eqref{6.7} and \eqref{6.8}, reduces to 
\begin{align*}
\Psi_{1,v}(s_1,s_2)=&\int_{K_v}\int_{F_v^{\times}}W_{\pi_v'}\left(\begin{pmatrix}
	a_v\\
	&1
\end{pmatrix}k_v\right)\overline{W_{\pi_v'}\left(\begin{pmatrix}
	a_v\\
	&1
\end{pmatrix}k_v\right)}\\
&\Vol(K_v[m_v])^{-2}\mathbf{1}_{k_v\in K_v[m_v]}|a_v|_v^{s_1+s_2-1}d^{\times}a_vdk_v.
\end{align*}

Since $W_{\pi_v'}$ is right invariant by $K_v[m_v]$, by \cite[Theorem 4.1]{Miy14} we obtain that 
\begin{align*}
\Psi_{1,v}(s_1,s_2)=&\frac{L_v(s_1+s_2,\pi_v'\times\overline{\pi}_v')}{\Vol(K_v[m_v])\zeta_v(2s_1+2s_2)}.
\end{align*}

\item Suppose that $v<\infty$, and $v\nmid\mathfrak{M}\mathfrak{N}$. Then  
\begin{align*}
\Psi_{1,v}(s_1,s_2)=&\int_{K_v}\int_{F_v^{\times}}W_{\pi_v'}\left(\begin{pmatrix}
	a_v\\
	&1
\end{pmatrix}k_v\right)\overline{W_{\pi_v'}\left(\begin{pmatrix}
	a_v\\
	&1
\end{pmatrix}k_v\right)}\\
&L_v(2s_1,\omega_v^{-1}\omega_v')L_v(2s_2,\omega_v\omega_v'^{-1})|a_v|_v^{s_1+s_2-1}d^{\times}a_vdk_v.
\end{align*}

By Casselman-Shalika formula, we obtain 
\begin{align*}
\Psi_{1,v}(s_1,s_2)=&\frac{L_v(2s_1,\omega_v^{-1}\omega_v')L_v(2s_2,\omega_v\omega_v'^{-1})L_v(s_1+s_2,\pi_v'\times\overline{\pi}_v')}{\zeta_v(2s_1+2s_2)}.
\end{align*}

\end{itemize}
	
Combining the above discussions, Lemma \ref{lem6.3} follows.
\end{proof}

\subsubsection{An upper bound  of $\mathcal{M}_{\const,1}^{\du,\heartsuit,\Reg}(1/2+it,1/2-it;\mathfrak{L})$}\label{sec6.1.3}
Let $C_{\fin}(\pi,\pi')$ be the conductor defined in \textsection\ref{sec4.1},  $t\in \mathbb{R}$ be the subconvexity parameter in \textsection\ref{sec4.2.}, and $L>10^{3}$ be the amplification parameter in \textsection\ref{sec4.2}. Let $s_1=1/2+it-s$, and $s_2=1/2-it+2s$, where $|s|=\varepsilon_0$. 

\begin{prop}\label{prop6.4}
Let notation be as before. Let $\varepsilon>0$. Then   
\begin{align*}
\mathcal{M}_{\const,1}^{\du,\heartsuit,\Reg}(1/2+it,1/2-it;\mathfrak{L})\ll 
\frac{C_{\fin}(\pi,\pi')\Delta_1}{\sqrt{\|\mathfrak{L}\|}}\prod_{v\mid\mathfrak{L}}\Big[|\lambda_{\pi'}(\mathfrak{p}_v^{\ell_v'})|+1\Big],
\end{align*}
where the implied constant depends only on $\varepsilon$ and the base field $F$. Here $\Delta_1:=C_{\fin}(\pi,\pi')^{\varepsilon}C_{\infty}(\pi')^{\varepsilon}C_{\infty}(\omega\omega'^{-1} |\cdot|^{2it})^{\varepsilon}C_{\infty}(\pi'\times\overline{\pi}')^{\varepsilon}$, and $\|\mathfrak{L}\|:=\prod_{v\in\mathcal{P}}q_v^{\ell_v'}=\prod_{v\mid\mathfrak{L}}q_v^{\ell_v'}$ is defined as in \textsection\ref{sec4.2}.
\end{prop} 
\begin{proof} 
For $v\in \Sigma_F$, let $I_v^{\sharp}(s_1,s_2;\mathfrak{L})$ be the local factor defined in Lemma \ref{lem6.3}. 

Let $v\mid\mathfrak{L}$. Then $I_v^{\sharp}(1/2+it-s,1/2-it+2s;\mathfrak{L})$ is equal to 
\begin{align*}
\frac{\omega_v(\varpi_v^{\ell_v'})}{q_v^{(1/2-it+2s)\ell_v'}}\Bigg[\frac{(q_v^{\ell_v'/2}+\ell_v'q_v^{-\ell_v'/2}|\mathbb{F}_v^{\times}|)\mathcal{RS}(1/2+it-s,1/2-it+2s)}{\zeta_v(2+2s)^{-1}L_v(1-2it+4s,\omega_v\omega_v'^{-1})L_v(1+s,\pi_v'\times\overline{\pi}_v')}+\overline{\lambda_{\pi'}(\mathfrak{p}_v^{\ell})}\Bigg].
\end{align*}

By the estimate \eqref{eq6.8} in Lemma \ref{lemm6.5}, we obtain 
\begin{equation}\label{6.11}
I_v^{\sharp}(1/2+it-s,1/2-it+2s;\mathfrak{L})\ll \ell_v'^2 q_v^{-\frac{\ell_v'}{2}+7\ell_v'\varepsilon_0}\sum_{0\leq i\leq \ell_v'}|\lambda_{\pi'}(\mathfrak{p}_v^{i})|,
\end{equation}
where the implied is absolute. Thus, by \eqref{6.11}, and the inequality 
\begin{equation}\label{eq6.15}
\sum_{0\leq i\leq \ell_v'}|\lambda_{\pi'}(\mathfrak{p}_v^{i})|\leq 100\ell_v'\cdot (|\lambda_{\pi'}(\mathfrak{p}_v^{\ell_v'})|+1),
\end{equation}
we deduce the following estimate: 
\begin{equation}\label{6.13}
\prod_{v\mid\mathfrak{L}}I_v^{\sharp}(1/2+it-s,1/2-it+2s;\mathfrak{L})\ll \|\mathfrak{L}\|^{-\frac{1}{2}+\varepsilon}\prod_{v\mid\mathfrak{L}}\Big[|\lambda_{\pi'}(\mathfrak{p}_v^{\ell_v'})|+1\Big].
\end{equation}
Here the implied constant depends only on $\varepsilon.$ 

By definition, $\prod_{v\mid\mathfrak{M}\mathfrak{N}}I_v^{\sharp}(1/2+it-s,1/2-it+2s;\mathfrak{L})$ is equal to 
\begin{align*}
\prod_{v\mid\mathfrak{M}\mathfrak{N}}\frac{\Vol(K_v[m_v])^{-1}}{L_v(1+2it-2s,\omega_v^{-1}\omega_v')L_v(1-2it+4s,\omega_v\omega_v'^{-1})}\ll C_{\fin}(\pi,\pi')^{1+\varepsilon}.
\end{align*}

Let $v\mid\infty$. By definition (see Lemma \ref{lem6.3} ),  $I_v^{\sharp}(1/2+it-s,1/2-it+2s;\mathfrak{L})$ is 
\begin{align*}
\int_{F_v}h_v\left(\begin{pmatrix}
1&\\
c_v& 1	
\end{pmatrix}
,1/2+it-s\right)\overline{h_v\left(\begin{pmatrix}
1&\\
c_v& 1	
\end{pmatrix}
,1/2+it+2\overline{s}\right)}\mathcal{I}(c_v)dc_v,
\end{align*}
where 
\begin{align*}
\mathcal{I}(c_v):=\int_{F_v^{\times}}W_{\pi_v'}\left(\begin{pmatrix}
	a_v\\
	c_v&1
\end{pmatrix}\right)\overline{W_{\pi_v'}\left(\begin{pmatrix}
	a_v\\
	c_v&1
\end{pmatrix}\right)}|a_v|_v^{s}d^{\times}a_v.
\end{align*}
\begin{comment}
Recall that $W_{\pi_v'}\left(\begin{pmatrix}
	a_v\\
	&1
\end{pmatrix}\right)$ is a fixed smooth function of $a_v\in F_v^{\times}$ with compact support with $W_{\pi_v'}(I_2)=1$. Then 
\begin{align*}
\int_{F_v^{\times}}W_{\pi_v'}&\left(\begin{pmatrix}
	a_v\\
	&1
\end{pmatrix}\right)\overline{W_{\pi_v'}\left(\begin{pmatrix}
	a_v\\
	&1
\end{pmatrix}\right)}|a_v|_v^{|\Re(s)|}d^{\times}a_v\ll 1.
\end{align*}
\end{comment}

Recall that $|s|=\varepsilon_0\leq 2(\log C_{\infty}(\pi')+\log C_{\fin}(\pi,\pi')+\log L)^{-10}$, see \textsection\ref{sect4.4}. By Lemma \ref{lem5.6}, we have 
\begin{align*}
\int_{F_v^{\times}}\bigg|W_{\pi_v'}\left(\begin{pmatrix}
	a_v\\
	c_v&1
\end{pmatrix}\right)\bigg|^2|a_v|_v^{s}d^{\times}a_v\ll & C_v(\pi')^{\varepsilon_0}\int_{F_v^{\times}}\bigg|W_{\pi_v'}\left(\begin{pmatrix}
	a_v\\
	c_v&1
\end{pmatrix}\right)\bigg|^2d^{\times}a_v\\
&+C_v(\pi')^{-10}(|c_v|_v+|c_v|_v^l+|c_v|_vC_v(\pi')^{1+\varepsilon}).
\end{align*}
As a consequence, 
\begin{equation}\label{eq6.14}
\int_{F_v^{\times}}W_{\pi_v'}\left(\begin{pmatrix}
	a_v\\
	c_v&1
\end{pmatrix}\right)\overline{W_{\pi_v'}\left(\begin{pmatrix}
	a_v\\
	c_v&1
\end{pmatrix}\right)}|a_v|_v^{s}d^{\times}a_v\ll C_v(\pi')^{\varepsilon_0},
\end{equation}
where the implied constant depends only on $\varepsilon$, $F_v$, and $\varphi_v$. Here we use the fact that $\langle W_{\pi_v'},W_{\pi_v'}\rangle\ll 1$.

Utilizing \eqref{eq6.14}, in conjunction with Cauchy-Schwarz inequality, the function $|I_v^{\sharp}(1/2+it-s,1/2-it+2s;\mathfrak{L})|$ is majorized by 
\begin{align*}
\Bigg[\int_{K_v}|h_v(k_v,1/2+it-s)|^2dk_v\Bigg]^{\frac{1}{2}}\Bigg[\int_{K_v}|h_v(k_v,1/2+it+2\overline{s})|^2dk_v\Bigg]^{\frac{1}{2}}.
\end{align*}
Similar to the proof of \cite[Definition 3.6.4]{MV10} we obtain 
\begin{equation}\label{6.14}
|I_v^{\sharp}(1/2+it-s,1/2-it+2s;\mathfrak{L})|\ll 1,
\end{equation}
where the implied constant depends only on $F_v$ and the fixed function $h$ (see \textsection\ref{sec4.4}).

Combining \eqref{6.13}, \eqref{6.14} with Lemma \ref{lem6.3}, we thus obtain 
\begin{equation}\label{eq6.16}
\mathcal{M}_{\const,1}^{\du,\heartsuit}(\cdots)\ll \frac{|\mathcal{L}_1(s;\pi')|C_{\fin}(\pi,\pi')^{1+\varepsilon}}{\sqrt{\|\mathfrak{L}\|}}\prod_{v\mid\mathfrak{L}}\Big[\sum_{0\leq i\leq \ell_v'}|\lambda_{\pi'}(\mathfrak{p}_v^{i})|\Big],
\end{equation}
where $\mathcal{M}_{\const,1}^{\du,\heartsuit}(\cdots):=\mathcal{M}_{\const,1}^{\du,\heartsuit}(1/2+it-s,1/2-it+2s;\mathfrak{L})$, and 
\begin{align*}
\mathcal{L}_1(s;\pi'):=\frac{L(1+s,\pi'\times\overline{\pi}')L(1+2it-2s,\omega^{-1}\omega')L(1-2it+4s,\omega\omega'^{-1})}{\zeta_F(2+2s)}.
\end{align*}

Appealing the well known upper bound for $L$-functions in the zero-free region, 
\begin{equation}\label{eq6.17}
\mathcal{L}_1(s;\pi')\ll C_{\fin}(\pi,\pi')^{\varepsilon}C_{\infty}(\omega\omega'^{-1} |\cdot|^{2it})^{\varepsilon}C_{\infty}(\pi'\times\overline{\pi}')^{\varepsilon},
\end{equation}
where the implied constant depends on $\varepsilon$ and $F$.

By defintion \eqref{equ6.2} and the triangle inequality,
\begin{equation}\label{6.18}
\mathcal{M}_{\const,1}^{\du,\heartsuit,\Reg}(1/2+it,1/2-it;\mathfrak{L})\ll \max_{|s|=\varepsilon_0}|\mathcal{M}_{\const,1}^{\du,\heartsuit}(\cdots)|,
\end{equation}
where $\mathcal{M}_{\const,1}^{\du,\heartsuit}(\cdots):=\mathcal{M}_{\const,1}^{\du,\heartsuit}(1/2+it-s,1/2-it+2s;\mathfrak{L})$.

Therefore, Proposition \ref{prop6.4} follows from  \eqref{eq6.16}, \eqref{eq6.17} and \eqref{6.18}.
\end{proof}

\subsection{The Constant Part: $\mathcal{M}_{\const,2}^{\du,\heartsuit,\Reg}(1/2+it,1/2-it;\mathfrak{L})$}\label{sec6.2}
In this section we aim to establish an upper bound for $\mathcal{M}_{\const,2}^{\du,\heartsuit,\Reg}(1/2+it,1/2-it;\mathfrak{L})$.

\subsubsection{Meromorphic continuation of the archimedean integrals} 
Let $v\mid\infty$. Let $h_v$ be the fixed function in \textsection\ref{sec4.4}, and 
\begin{align*}
\Phi_v(t_{1,v},t_{2,v}):=\textbf{C}_v^{\, 1/2}h(\textbf{C}_v|t_{1,v}|)h(|t_{2,v}|-1)\omega_v\omega_v'^{-1}(t_{2,v}). \tag{\ref{eq6.28}}
\end{align*}
Then $\Phi_v$ is a Schwartz function in $\mathbb{F}_v^2$. For $t_v\in F_v$ and $k_v\in K_v$, we define 
\begin{equation}\label{beta}
\beta_v(t_v,k_v):=\int_{F_v}\Phi_v\left((t_v,b_v)
k_v\right)db_v.
\end{equation}
Let $\widehat{\beta_v}(t_v,k_v)$ be the Fourier transform relative to the first variable:
\begin{align*}
\widehat{\beta_v}(t_v,k_v):=\int_{F_v}\int_{F_v}\Phi_v\left((u_v,b_v)
k_v\right)\psi_v(u_vt_v)db_vdu_v.
\end{align*}

For $s\in \mathbb{C}$, we define 
\begin{align*}
Z_v(s,\beta_v(\cdot,k_v),\omega_v^{-1}\omega_v'):=\int_{F_v^{\times}}\beta_v(t_v,k_v)\omega_v^{-1}\omega_v'(t_v)|t_v|_v^{s}d^{\times}t_v,
\end{align*} 
which converges absolutely in $\Re(s)>0$, and is a local integral representing the archimedean $L$-factor $L_v(s,\omega_v^{-1}\omega_v')$. Let 
\begin{equation}\label{f6.26}
Z_v(s,\widehat{\beta_v}(\cdot,k_v),\omega_v\omega_v'^{-1}):=\int_{F_v^{\times}}\widehat{\beta_v}(t_v,k_v)\omega_v\omega_v'^{-1}(t_v)|t_v|_v^{s}d^{\times}t_v,
\end{equation}
which converges absolutely in $\Re(s)>0$. We have the functional equation 
\begin{equation}\label{eq6.26}
Z_v(1-s,\widehat{\beta_v}(\cdot,k_v),\omega_v\omega_v'^{-1})=\gamma_v(s,\omega_v^{-1}\omega_v',\psi_v)Z_v(s,\beta_v,\omega_v^{-1}\omega_v'),
\end{equation}
where $\gamma_v(s,\omega_v^{-1}\omega_v',\psi_v)$ is the gamma factor of $\omega_v^{-1}\omega_v'$ with respect to the fixed additive character $\psi_v$. (see \textsection\ref{sec2.1.4}). In particular, \eqref{eq6.26} gives a meromorphic continuation of $Z_v(s,\beta_v,\omega_v^{-1}\omega_v')$ from $\Re(s)>0$ to the region $\Re(s)<1.$

Let $v\mid\infty$. Let $0< \Re(s_1)<1$ and $s_2\in \mathbb{C}$. Define 
\begin{equation}\label{J_v}
\mathcal{J}_v^{\heartsuit}(s_1,s_2;k_v):=\frac{Z_v(2-2s_1,\widehat{\beta_v}(\cdot,k_v),\omega_v\overline{\omega_v'})L_v(2s_1-1,\overline{\omega_v}\omega_v')\overline{h_v(k_v,\overline{s_2})}}{\varepsilon_v(2s_1-1,\omega_v^{-1}\omega_v',\psi_v)L_v(2-2s_1,\omega_v\omega_v'^{-1})}.
\end{equation}
Here $\varepsilon_v(s,\omega_v^{-1}\omega_v',\psi_v)$ is the $\varepsilon$-factor of $\omega_v^{-1}\omega_v'$ with respect to the fixed additive character $\psi_v$. 

\begin{lemma}\label{lemma6.5}
Let notation be as before. Let $v\mid\infty$. Then $\mathcal{J}_v^{\heartsuit}(s_1,s_2;k_v)/L_v(2s_1-1,\overline{\omega_v}\omega_v')$ converges absolutely in the region $\{(s_1,s_2)\in \mathbb{C}^2:\ 0<\Re(s_1)<1\}$. Moreover, for $0<\Re(s_1)<1$ and $s_2\in \mathbb{C}$, we have 
\begin{equation}\label{f6.29}
\frac{\mathcal{J}_v^{\heartsuit}(s_1,s_2;k_v)}{L_v(2s_1-1,\overline{\omega_v}\omega_v')}\ll \frac{4^{\Re(s_2)}\textbf{C}_v^{\, 2-2\Re(s_1)+\varepsilon}\mathbf{1}_{|k_{21}|_v\leq \varepsilon^*\textbf{C}_v^{\, -1}}}{|(2-2s_1)L_v(2-2s_1,\omega_v\omega_v'^{-1})|},
\end{equation}
where the implied constant depends only on $\varepsilon$, $F_v$ and the fixed $h_v.$
\end{lemma}
\begin{proof}
By Tate's thesis and \eqref{J_v}, $\mathcal{J}_v^{\heartsuit}(s_1,s_2;k_v)/L_v(2s_1-1,\overline{\omega_v}\omega_v')$ converges absolutely in the region $\{(s_1,s_2)\in \mathbb{C}^2:\ 0<\Re(s_1)<1\}$.

According to the definition \eqref{eq6.28} and utilizing the uncertainty principle,  
\begin{align*}
\widehat{\beta_v}(t_v,k_v):= \int_{F_v}\int_{F_v}\Phi_v\left((u_v,b_v)
k_v\right)\psi_v(u_vt_v)db_vdu_v\ll \textbf{C}_v^{\, -1/2}\mathbf{1}_{|t_v|_v\ll \textbf{C}_v^{\, 1+\varepsilon}}.
\end{align*}

Substituting this into \eqref{f6.26} we obtain  
\begin{equation}\label{equ6.29}
Z_v(2-2s_1,\widehat{\beta_v}(\cdot,k_v),\omega_v\omega_v'^{-1})\ll |2-2s_1|^{-1}\cdot \textbf{C}_v^{\, (1+\varepsilon)(2-2\Re(s_1))-1/2},
\end{equation}
where the implied constant depends only on $F_v$ and the fixed $h_v$. In particular, it is independent of $k_v\in K_v$. 

By the definition of archimedean $\varepsilon$-factors (see \cite[(3.7) and (4.7)]{Kna94}) we obtain $|\varepsilon_v(2s_1-1,\omega_v^{-1}\omega_v',\psi_v)|=1$. Thus, by \eqref{J_v} and \eqref{equ6.29}, 
\begin{equation}\label{eq6.31}
\Bigg|\frac{\mathcal{J}_v^{\heartsuit}(s_1,s_2;k_v)}{L_v(2s_1-1,\overline{\omega_v}\omega_v')}\Bigg|\ll \Bigg|\frac{\textbf{C}_v^{\, (1+\varepsilon)(2-2\Re(s_1))-1/2}}{(2-2s_1)L_v(2-2s_1,\omega_v\omega_v'^{-1})}\Bigg|\cdot |\overline{h_v(k_v,\overline{s_2})}|.
\end{equation}

Write $k_v=\begin{pmatrix}
k_{11}& k_{12}\\
k_{21}& k_{22}
\end{pmatrix}\in K_v$. By the definition \eqref{2.3}, we have  
\begin{equation}\label{f6.32}
|\overline{h_v(k_v,\overline{s_2})}|\ll 4^{\Re(s_2)}\textbf{C}_v^{\, 1/2}\mathbf{1}_{|k_{21}|_v\leq \varepsilon^*\textbf{C}_v^{\, -1}}.
\end{equation}

Therefore, \eqref{f6.29} follows from \eqref{eq6.31} and \eqref{f6.32}. 
\end{proof}

\begin{lemma}\label{lem6.7}
Let notation be as before. Let $v\mid\infty$. Then 
\begin{equation}\label{equa6.27}
\mathcal{J}_v(s_1,s_2;k_v):=\int_{N(F_v)}h_v\left(wu_vk_v,s_1\right)du_v\overline{h_v(k_v,\overline{s_2})}
\end{equation}
converges absolutely in $\{(s_1,s_2)\in \mathbb{C}^2:\ \Re(s_1)>1/2\}$. Moroever, it admits the meromorphic continuation $\mathcal{J}_v^{\heartsuit}(s_1,s_2;k_v)$ to $\{(s_1,s_2)\in \mathbb{C}^2:\ \Re(s_1)>0\}$. Here $\mathcal{J}_v^{\heartsuit}(s_1,s_2;k_v)$ is defined by \eqref{J_v}.
\end{lemma}
\begin{proof}
Let $\Re(s_1)>1/2$. Using a change of variables we obtain  
\begin{align*}
\int_{N(F_v)}h_v\left(wu_vk_v,s_1\right)du_v=\int_{F_v^{\times}}\beta_v(t_v,k_v)\omega_v^{-1}\omega_v'(t_v)|t_v|_v^{2s_1-1}d^{\times}t_v,
\end{align*}
where $\beta_v$ is defined by \eqref{beta}. Consequently,
\begin{align*}
\int_{N(F_v)}h_v\left(wu_vk_v,s_1\right)du_v=Z_v(2s_1-1,\beta_v(\cdot,k_v),\omega_v^{-1}\omega_v').
\end{align*}
By Tate's thesis, this integral converges absolutely in $\Re(s_1)>1/2$. As a result, $\mathcal{J}_v(s_1,s_2;k_v)$ converges absolutely in $\Re(s_1)>1/2.$

For $1/2<\Re(s_1)<1$, we have, by \eqref{eq6.26}, that
\begin{align*}
\int_{N(F_v)}h_v\left(wu_vk_v,s_1\right)du_v=\frac{Z_v(2-2s_2,\widehat{\beta_v}(\cdot,k_v),\omega_v\omega_v'^{-1})}{\gamma_v(2s_1-1,\omega_v^{-1}\omega_v',\psi_v)}.
\end{align*}

Then $\mathcal{J}_v(s_1,s_2;k_v)=\mathcal{J}_v^{\heartsuit}(s_1,s_2;k_v)$ for $1/2<\Re(s_1)<1$, leveraging the relation
\begin{align*}
\gamma_v(2s_1-1,\omega_v^{-1}\omega_v',\psi_v)=\varepsilon_v(2s_1-1,\omega_v^{-1}\omega_v',\psi_v)\frac{L_v(2-2s_1,\omega_v\omega_v'^{-1})}{L_v(2s_1-1,\omega_v^{-1}\omega_v')}.
\end{align*}

Hence, the function $\mathcal{J}_v^{\heartsuit}(s_1,s_2;k_v)$ could be regarded as the meromorphic continuation of $\mathcal{J}_v(s_1,s_2;k_v)$ in the region $\{(s_1,s_2)\in \mathbb{C}^2:\ \Re(s_1)>0\}.$
\end{proof}

\subsubsection{The Contribution from $\mathcal{M}_{\const,2}^{\du,\heartsuit}(s_1,s_2)$}
By definition in \textsection\ref{sec3.2}, the function $\mathcal{M}_{\const,2}^{\du,\heartsuit}(s_1,s_2;\mathfrak{L})$ is defined by 
\begin{equation}\label{6.19}
\gamma\cdot\int_{B(F)Z(\mathbb{A}_F)\backslash G(\mathbb{A}_F)}\phi(x)\overline{\phi(x\boldsymbol{d}_{\mathfrak{L}})}\int_{N(\mathbb{A}_F)}h(wux,s_1)du\overline{h(x\boldsymbol{d}_{\mathfrak{L}},\overline{s_2})}dx,
\end{equation}
where $\gamma:=\prod_{v\in\mathcal{P}}\omega_v(\varpi_v^{\ell_v'})(q_v^{\ell_v'/2}+\ell_v'q_v^{-\ell_v'/2}|\mathbb{F}_v^{\times}|)$,  $\boldsymbol{d}_{\mathfrak{L}}=\otimes_{v\in \Sigma_F}\ \boldsymbol{d}_v\in G(\mathbb{A}_F)$, where $\boldsymbol{d}_v=\begin{pmatrix}
1&\\
&\varpi_v^{\ell_v'}
\end{pmatrix}$ with the extension that $\ell_v':=0$ if $v\not\in \mathcal{P}$.  

The following is a more concise version of Lemma \ref{lem3.7}. 
\begin{lemma}\label{lem6.5}
Let notation be as before. Let $0<\Re(s_1)<1$ and $s_2\in \mathbb{C}$. Then $\mathcal{M}_{\const,2}^{\du,\heartsuit}(s_1,s_2;\mathfrak{L})$ is equal to 
\begin{align*}
\frac{L(1+s_2-s_1,\pi'\times\overline{\pi}'\otimes \omega\overline{\omega'})L(2s_1-1,\overline{\omega}\omega')L(2s_2,\omega\overline{\omega'})}{L(2+2s_2-2s_1,\omega^2\omega'^{-2})}\prod_{v\in\Sigma_F}J_v^{\sharp}(s_1,s_2;\mathfrak{L}),
\end{align*}
where the local factors $J_v^{\sharp}(s_1,s_2;\mathfrak{L})$ are defined as follows:
\begin{itemize}
\item at $v\mid\mathfrak{L}$, $J_v^{\sharp}(s_1,s_2;\mathfrak{L})$ is defined by  
\begin{align*}
\frac{(q_v^{\ell_v'/2}+\ell_v'q_v^{-\ell_v'/2}|\mathbb{F}_v^{\times}|)\omega_v(\varpi_v^{\ell_v'})\mathbb{L}_1(s_2,s_2)\mathcal{RS}(1-s_1,s_2)}{q_v^{s_2\ell_v'}\mathbb{L}_2(s_2,s_2)}+\frac{\overline{\lambda_{\pi'}(\mathfrak{p}_v^{\ell})}\omega_v(\varpi_v^{\ell_v'})}{q_v^{s_2\ell_v'}},
\end{align*}
where $\mathbb{L}_1(s_2,s_2):=L_v(2s_1,\omega_v^{-1}\omega_v')L_v(2+2s_2-2s_1,\omega_v^2\omega_v'^{-2})$, $\mathbb{L}_2(s_2,s_2):=L_v(2s_2,\omega_v\omega_v'^{-1})L_v(2s_1-1,\omega_v^{-1}\omega_v')L_v(1+s_2-s_1,\pi_v'\times\overline{\pi}_v')$, and the function $\mathcal{RS}(s_1,s_2)$ is defined right before Lemma \ref{lemm6.5};
\item at $v\mid\mathfrak{M}\mathfrak{N}$, with $m_v:=\max\{r_{\pi_v},r_{\pi_v'}\}$, we define  
\begin{align*}
J_v^{\sharp}(s_1,s_2;\mathfrak{L}):=\frac{q_v^{m_v(1-2s_1)}\mathbf{1}_{\text{$\omega_v\omega_v'^{-1}$ is unramified}}}{\zeta_v(1)\Vol(K_v[m_v])L_v(2s_2,\omega_v\overline{\omega_v'})}.
\end{align*}

\item and at archimedean places $v\mid\infty$, the integral $J_v^{\sharp}(s_1,s_2;\mathfrak{L})$ is defined by 
\begin{align*}
\int_{K_v}\int_{F_v^{\times}}\Big|W_{\pi_v'}\left(\begin{pmatrix}
	a_v\\
	&1
\end{pmatrix}k_v\right)\Big|^2|a_v|_v^{s_2-s_1}\omega_v\omega_v'^{-1}(a_v)d^{\times}a_v\mathcal{J}_v^{\heartsuit}(s_1,s_2;k_v)dk_v.
\end{align*}
Here $\mathcal{J}_v^{\heartsuit}(s_1,s_2;k_v)$ is defined as in \eqref{J_v}.
\item at the remaining places, i.e., $v<\Sigma_{F,\fin}$ and $v\nmid\mathfrak{M}\mathfrak{N}$, $J_v^{\sharp}(s_1,s_2;\mathfrak{L})\equiv 1$. 
\end{itemize}
\end{lemma}
\begin{proof}
Let $\Re(s_1)\gg 1$ and $\Re(s_2)
\gg 1$. Executing Fourier expansion in \eqref{6.19} and using the Iwasawa decomposition we have
\begin{align*}
\mathcal{M}_{\const,2}^{\du,\heartsuit}(s_1,s_2;\mathfrak{L})=\prod_{v\in \Sigma_F}\Psi_{2,v}(s_1,s_2),
\end{align*}
where for each $v\in \Sigma_F$,, the factor $\Psi_{2,v}(s_1,s_2)$ is defined by 
\begin{align*}
&\omega_v(\varpi_v^{\ell_v'})(q_v^{\ell_v'/2}+\ell_v'q_v^{-\ell_v'/2}|\mathbb{F}_v^{\times}|)\int_{K_v}\int_{F_v^{\times}}W_{\pi_v'}\left(\begin{pmatrix}
	a_v\\
	&1
\end{pmatrix}k_v\right)\overline{W_{\pi_v'}\left(\begin{pmatrix}
	a_v\\
	&1
\end{pmatrix}k_v\boldsymbol{d}_v\right)}\\
&\int_{N(F_v)}h_v\left(wu_v\begin{pmatrix}
	a_v\\
	&1
\end{pmatrix}k_v,s_1\right)du_v\overline{h_v(k_v\boldsymbol{d}_v,\overline{s_2})}|a_v|_v^{s_2-1}d^{\times}a_vdk_v.
\end{align*}

By definition
\begin{align*}
h_v(x_v,s_1)=|\det x_v|_v^{s_1}\int_{F_v^{\times}}\Phi_v((0,t_v)x_v)\omega_v^{-1}\omega_v'(t_v)|t_v|_v^{2s_1}d^{\times}t_v.
\end{align*}
Therefore, the constant term 
\begin{align*}
\int_{N(F_v)}h_v\left(wu_v\begin{pmatrix}
	a_v\\
	&1
\end{pmatrix}k_v,s_1\right)du_v
\end{align*}
is equal to 
\begin{align*}
|a_v|_v^{1-s_1}\omega_v\omega_v'^{-1}(a_v)\int_{F_v}\int_{F_v^{\times}}\Phi_v\left((t_v,b_v)k_v\right)\omega_v^{-1}\omega_v'(t_v)|t_v|_v^{2s_1-1}d^{\times}t_vdb_v.
\end{align*}
Along with the definition \eqref{6.8} of $\overline{h_v(k_v\boldsymbol{d}_v,\overline{s_2})}$, we conclude that 
\begin{comment}
\begin{align*}
&\omega_v(\varpi_v^{\ell_v'})(q_v^{\ell_v'/2}+\ell_v'q_v^{-\ell_v'/2}|\mathbb{F}_v^{\times}|)\int_{K_v}\int_{F_v^{\times}}W_{\pi_v'}\left(\begin{pmatrix}
	a_v\\
	&1
\end{pmatrix}k_v\right)\overline{W_{\pi_v'}\left(\begin{pmatrix}
	a_v\\
	&1
\end{pmatrix}k_v\boldsymbol{d}_v\right)}\\
&|a_v|_v^{s_2-s_1}\omega_v\omega_v'^{-1}(a_v)\int_{F_v}\int_{F_v^{\times}}\Phi_v\left((t_v,b_v)k_v\right)\omega_v^{-1}\omega_v'(t_v)|t_v|_v^{2s_1-1}d^{\times}t_vdb_v\\
&q_v^{-s_2\ell_v'}\int_{F_v^{\times}}\overline{\Phi_v((0,t_v')k_v\boldsymbol{d}_v)}\omega_v\omega_v'^{-1}(t_v')|t_v'|_v^{2s_2}d^{\times}t_v'd^{\times}a_vdk_v
\end{align*}

\begin{align*}
&\omega_v(\varpi_v^{\ell_v'})(q_v^{\ell_v'/2}+\ell_v'q_v^{-\ell_v'/2}|\mathbb{F}_v^{\times}|)q_v^{-s_2\ell_v'}\int_{K_v}\int_{F_v}\int_{F_v^{\times}}\Phi_v\left((t_v,b_v)k_v\right)\omega_v^{-1}\omega_v'(t_v)|t_v|_v^{2s_1-1}d^{\times}t_vdb_v\\
&\int_{F_v^{\times}}\overline{\Phi_v((0,t_v')k_v\boldsymbol{d}_v)}\omega_v\omega_v'^{-1}(t_v')|t_v'|_v^{2s_2}d^{\times}t_v'dk_v
\end{align*}
\end{comment}
\begin{equation}\label{6.20.}
\Psi_{2,v}(s_1,s_2)=\frac{\omega_v(\varpi_v^{\ell_v'})(q_v^{\ell_v'/2}+\ell_v'q_v^{-\ell_v'/2}|\mathbb{F}_v^{\times}|)}{q_v^{s_2\ell_v'}}\int_{K_v}\mathcal{J}_1(k_v)\mathcal{J}_2(k_v)\mathcal{J}_3(k_v)dk_v,
\end{equation}
where 
\begin{align*}
\mathcal{J}_1(k_v):=&\int_{F_v^{\times}}\overline{\Phi_v((0,t_v')k_v\boldsymbol{d}_v)}\omega_v\omega_v'^{-1}(t_v')|t_v'|_v^{2s_2}d^{\times}t_v'\\
\mathcal{J}_2(k_v):=&\int_{F_v}\int_{F_v^{\times}}\Phi_v\left((t_v,b_v)k_v\right)\omega_v^{-1}\omega_v'(t_v)|t_v|_v^{2s_1-1}d^{\times}t_vdb_v,
\end{align*}
and the function $\mathcal{J}_3(k_v)$ is defined by 
\begin{align*}
\int_{F_v^{\times}}W_{\pi_v'}\left(\begin{pmatrix}
	a_v\\
	&1
\end{pmatrix}k_v\right)\overline{W_{\pi_v'}\left(\begin{pmatrix}
	a_v\\
	&1
\end{pmatrix}k_v\boldsymbol{d}_v\right)}|a_v|_v^{s_2-s_1}\omega_v\omega_v'^{-1}(a_v)d^{\times}a_v.
\end{align*}

Let $v<\infty$. Write $k_v=\begin{pmatrix}
k_{11}& k_{12}\\
k_{21}& k_{22}
\end{pmatrix}$, and $m_v:=\max\{r_{\pi_v},r_{\pi_v'}\}$.  We may assume $k_{22}\neq 0$. Making use of the definition of $\Phi_v(\cdot,\cdot)$ in \textsection\ref{sec4.4}, we have the following cases.
\begin{itemize}
\item Suppose $v\mid\mathfrak{L}$. Then by \eqref{4.10}, 
\begin{equation}\label{6.20}
\mathcal{J}_1(k_v)=L_v(2s_2,\omega_v\omega_v'^{-1})\mathbf{1}_{k_v\in K_v}+\sum_{j=1}^{\ell_v'}\omega_v^{-1}\omega_v'(\varpi_v^j)q_v^{2s_2j}\mathbf{1}_{k_v\in K_v[j]}.
\end{equation}

By definition \eqref{equ4.9} we obtain 
\begin{equation}\label{6.21}
\mathcal{J}_2(k_v)=\int_{\mathcal{O}_v}\int_{\mathcal{O}_v-\{0\}}\omega_v^{-1}\omega_v'(t_v)|t_v|_v^{2s_1-1}d^{\times}t_vdb_v=L_v(2s_1-1,\omega_v^{-1}\omega_v').
\end{equation}

Substituting \eqref{6.20} and \eqref{6.21} into \eqref{6.20.} we then derive that 
\begin{align*}
\Psi_{2,v}(s_1,s_2)=\delta_v^{(2)}\mathcal{RS}(1-s_1,s_2)+\delta_v^{(2)}L_v(2s_2,\omega_v\omega_v'^{-1})\Psi_{2,v}^{(0)}(s_2-s_1+1),
\end{align*}
where $\delta_v^{(2)}:=q_v^{-s_2\ell_v'}(q_v^{\ell_v'/2}+\ell_v'q_v^{-\ell_v'/2}|\mathbb{F}_v^{\times}|)\omega_v(\varpi_v^{\ell_v'})L_v(2s_1-1,\omega_v^{-1}\omega_v')$, the function $\mathcal{RS}(s_1,s_2)$ is defined right before Lemma \ref{lemm6.5}; and 
\begin{align*}
\Psi_{2,v}^{(0)}(s):=\int_{K_v}\int_{F_v^{\times}}W_{\pi_v'}\left(\begin{pmatrix}
	a_v\\
	&1
\end{pmatrix}k_v\right)\overline{W_{\pi_v'}\left(\begin{pmatrix}
	a_v\\
	&1
\end{pmatrix}k_v\boldsymbol{d}_v\right)}|a_v|_v^{s-1}d^{\times}a_vdk_v.
\end{align*}

Since $\Re(s_2-s_1)\gg 1$,  then $\Psi_{2,v}^{(0)}(s_2-s_1+1)$ converges absolutely. Moreover, by \eqref{5.2} and the standard calculation of Rankin--Selberg convolutions, we then obtain that $\delta_v^{(2)}L_v(2s_2,\omega_v\omega_v'^{-1})\Psi_{2,v}^{(0)}(s_2-s_1+1)$ becomes
\begin{align*}
\frac{\overline{\lambda_{\pi'}(\mathfrak{p}_v^{\ell})}L_v(1+s_2-s_1,\pi_v'\times\overline{\pi}_v'\otimes \omega_v\omega_v'^{-1})L_v(2s_1-1,\omega_v^{-1}\omega_v')L_v(2s_2,\omega_v\omega_v'^{-1})}{\omega_v^{-1}(\varpi_v^{\ell_v'})q_v^{s_2\ell_v'}L_v(2+2s_2-2s_1,\omega_v^2\omega_v'^{-2})}.
\end{align*}

Therefore, by Lemma \ref{lemm6.5} we derive 
\begin{align*}
\Psi_{2,v}(s_1,s_2)=&q_v^{-s_2\ell_v'}(q_v^{\ell_v'/2}+\ell_v'q_v^{-\ell_v'/2}|\mathbb{F}_v^{\times}|)\omega_v(\varpi_v^{\ell_v'})L_v(2s_1,\omega_v^{-1}\omega_v')\mathcal{RS}(1-s_1,s_2)\\
&+\frac{\overline{\lambda_{\pi'}(\mathfrak{p}_v^{\ell})}L_v(1+s_2-s_1,\pi_v'\times\overline{\pi}_v'\otimes \omega_v\omega_v'^{-1})L_v(2s_1-1,\omega_v^{-1}\omega_v')}{\omega_v^{-1}(\varpi_v^{\ell_v'})q_v^{s_2\ell_v'}L_v(2+2s_2-2s_1,\omega_v^2\omega_v'^{-2})L_v(2s_2,\omega_v\omega_v'^{-1})^{-1}}.
\end{align*}

As a result, $J_v^{\sharp}(s_1,s_2;\mathfrak{L})$ is given by 
\begin{align*}
\frac{(q_v^{\ell_v'/2}+\ell_v'q_v^{-\ell_v'/2}|\mathbb{F}_v^{\times}|)\omega_v(\varpi_v^{\ell_v'})\mathbb{L}_1(s_2,s_2)\mathcal{RS}(1-s_1,s_2)}{q_v^{s_2\ell_v'}\mathbb{L}_2(s_2,s_2)}+\frac{\overline{\lambda_{\pi'}(\mathfrak{p}_v^{\ell})}\omega_v(\varpi_v^{\ell_v'})}{q_v^{s_2\ell_v'}}.
\end{align*}

\item  Suppose $v\mid\mathfrak{M}\mathfrak{N}$, i.e.,  $m_v:=\max\{r_{\pi_v},r_{\pi_v'}\}\geq 1$. In this case, $\ell_v'=0$, i.e., $\boldsymbol{d}_v=I_2$. Recall the definition of $\Phi_v(\cdot,\cdot):$
\begin{align*}
\Phi_v(t_{1,v},t_{2,v}):=\Vol(K_v[m_v])^{-1}\omega_v\omega_v'^{-1}(t_{2,v})\mathbf{1}_{\mathfrak{p}_v^{m_v}}(t_{1,v})\mathbf{1}_{\mathcal{O}_v^{\times}}(t_{2,v}).\tag{\ref{equa4.5}}
\end{align*}

Plugging this into 
\begin{align*}
\mathcal{J}_1(k_v):=\int_{F_v^{\times}}\overline{\Phi_v((0,t_v')k_v)}\omega_v\omega_v'^{-1}(t_v')|t_v'|_v^{2s_2}d^{\times}t_v',
\end{align*}
we obtain that
\begin{equation}\label{6.24}
\mathcal{J}_1(k_v)=\Vol(K_v[m_v])^{-1}\omega_v^{-1}\omega_v'(k_{22})\mathbf{1}_{k_v\in K_v[m_v]}.
\end{equation}

Let $k_v\in K_v[m_v]$. By the definition \eqref{equa4.5}, we have
\begin{align*}
\int_{F_v}\Phi_v\left((t_v,b_v)k_v\right)db_v=\Vol(K_v[m_v])^{-1}\omega_v\omega_v'^{-1}(k_{22})\mathbf{1}_{t_v\in\mathfrak{p}_v^{m_v}}\int_{\mathcal{O}_v^{\times}}\omega_v\omega_v'^{-1}(b_v)db_v.
\end{align*}

Substituting this into the integral 
\begin{align*}
\mathcal{J}_2(k_v)=\int_{F_v^{\times}}\int_{F_v}\Phi_v\left((t_v,b_v)k_v\right)db_v\omega_v^{-1}\omega_v'(t_v)|t_v|_v^{2s_1-1}d^{\times}t_v,
\end{align*}
we thus deduce the following calculation 
\begin{equation}\label{6.25}
\mathcal{J}_2(k_v)=\frac{\mathbf{1}_{\text{$\omega_v\omega_v'^{-1}$ is unramified}}q_v^{m_v(1-2s_1)}L_v(2s_1-1,\omega_v^{-1}\omega_v')}{\zeta_v(1)\Vol(K_v[m_v])\omega_v^{-1}\omega_v'(k_{22})}.
\end{equation}

Let $k_v\in K_v[m_v]$. Then $\mathcal{J}_3(k_v)$ is equal to 
\begin{align*}
\sum_{m\geq 0}\Big|W_{\pi_v'}\left(\begin{pmatrix}
	\varpi_v^{m}\\
	&1
\end{pmatrix}\right)\Big|^2|a_v|_v^{s_2-s_1}\omega_v\overline{\omega_v'}(\varpi_v^{m})=\frac{L_v(1+s_2-s_1,\pi_v\times\overline{\pi}_v'\otimes\omega_v\overline{\omega_v'})}{L_v(2+2s_2-2s_1,\omega_v^2\omega_v'^{-2})}.
\end{align*}
In conjunction with \eqref{6.24} and \eqref{6.25} we then deduce that 
\begin{align*}
\Psi_{2,v}(s_1,s_2)=\frac{q_v^{m_v(1-2s_1)}L_v(1+s_2-s_1,\pi_v\times\overline{\pi}_v'\otimes\omega_v\overline{\omega_v'})\mathbf{1}_{\text{$\omega_v\omega_v'^{-1}$ is unramified}}}{\zeta_v(1)\Vol(K_v[m_v])L_v(2s_1-1,\omega_v^{-1}\omega_v')^{-1}L_v(2+2s_2-2s_1,\omega_v^2\omega_v'^{-2})}.
\end{align*}

\item  Suppose that $v<\infty$, and $v\nmid\mathfrak{M}\mathfrak{N}$. By the definition of $\Phi_v(t_{1,v},t_{2,v})$,
\begin{align*}
\mathcal{J}_1(k_v)=\int_{F_v^{\times}}\overline{\Phi_v((0,t_v')k_v)}\omega_v\omega_v'^{-1}(t_v')|t_v'|_v^{2s_2}d^{\times}t_v'=L_v(2s_2,\omega_v\omega_v'^{-1}),
\end{align*}
and $\mathcal{J}_2(k_v)$ is equal to 
\begin{align*}
\int_{F_v^{\times}}\int_{F_v}\Phi_v\left((t_v,b_v)k_v\right)db_v\omega_v^{-1}\omega_v'(t_v)|t_v|_v^{2s_1-1}d^{\times}t_v=L_v(2s_1-1,\omega_v^{-1}\omega_v').
\end{align*}
Since $v\nmid\mathfrak{M}\mathfrak{N}$, then $W_{\pi_v'}$ is spherical. Thus, $\mathcal{J}_3(k_v)$ is equal to 
\begin{align*}
\sum_{m\geq 0}\Big|W_{\pi_v'}\left(\begin{pmatrix}
	\varpi_v^{m}\\
	&1
\end{pmatrix}\right)\Big|^2|a_v|_v^{s_2-s_1}\omega_v\overline{\omega_v'}(\varpi_v^{m})=\frac{L_v(1+s_2-s_1,\pi_v\times\overline{\pi}_v'\otimes\omega_v\overline{\omega_v'})}{L_v(2+2s_2-2s_1,\omega_v^2\omega_v'^{-2})}.
\end{align*}

In conclusion we obtain 
\begin{align*}
\Psi_{2,v}(s_1,s_2)=\frac{L_v(1+s_2-s_1,\pi_v\times\overline{\pi}_v'\otimes\omega_v\overline{\omega_v'})L_v(2s_2,\omega_v\omega_v'^{-1})L_v(2s_1-1,\omega_v^{-1}\omega_v')}{L_v(2+2s_2-2s_1,\omega_v^2\omega_v'^{-2})}.
\end{align*}

\item Suppose $v\mid\infty$. Then $\Psi_{2,v}(s_1,s_2)$ is equal to 
\begin{equation}\label{6.35}
\int_{K_v}\int_{F_v^{\times}}\Big|W_{\pi_v'}\left(\begin{pmatrix}
	a_v\\
	&1
\end{pmatrix}k_v\right)\Big|^2|a_v|_v^{s_2-s_1}\omega_v\omega_v'^{-1}(a_v)d^{\times}a_v\mathcal{J}_v(s_1,s_2;k_v)dk_v,
\end{equation}
where $\mathcal{J}_v(s_1,s_2;k_v)$ is defined by \eqref{equa6.27} in Lemma \ref{lem6.7}. By construction, $W_{\pi_v'}\left(\begin{pmatrix}
	a_v\\
	&1
\end{pmatrix}\right)$ has compact support. Hence, \eqref{6.35} converges absolutely in $\Re(s_1)>1/2$. Making use of Lemma \ref{lemma6.5} we obtain the meromorphic continuation of $\Psi_{2,v}(s_1,s_2)$ in $\{(s_1,s_2)\in \mathbb{C}^2:\ 0<\Re(s_1)<1\}$ as
\begin{align*}
\int_{K_v}\int_{F_v^{\times}}\Big|W_{\pi_v'}\left(\begin{pmatrix}
	a_v\\
	&1
\end{pmatrix}k_v\right)\Big|^2|a_v|_v^{s_2-s_1}\omega_v\omega_v'^{-1}(a_v)d^{\times}a_v\mathcal{J}_v^{\heartsuit}(s_1,s_2;k_v)dk_v.
\end{align*}
\end{itemize}

Therefore, Lemma \ref{lem6.5} then follows from the above discussions. 
\end{proof}

\subsubsection{An upper bound  of $\mathcal{M}_{\const,2}^{\du,\heartsuit,\Reg}(1/2+it,1/2-it;\mathfrak{L})$}\label{sec6.1.3}
Let $C_{\fin}(\pi,\pi')$ be the conductor defined in \textsection\ref{sec4.1},  $t\in \mathbb{R}$ be the subconvexity parameter in \textsection\ref{sec4.2.}, and $L>10^{3}$ be the amplification parameter in \textsection\ref{sec4.2}. Let $(\log C_{\infty}(\pi')+\log C_{\fin}(\pi,\pi')+\log L)^{-10}\leq \varepsilon_0\leq 2(\log C_{\infty}(\pi')+\log C_{\fin}(\pi,\pi')+\log L)^{-10}$ be the parameter defined in \textsection\ref{sect4.4}. Recall that $\textbf{C}_v(\pi,t,\pi')$ is  defined by \eqref{f4.5} in \textsection\ref{sec4.4}. Define $\textbf{C}_{\infty}:=\prod_{v\mid\infty}\textbf{C}_v(\pi,t,\pi')$. Let $s_1=1/2+it-s$, and $s_2=1/2-it+2s$, where $|s|=\varepsilon_0$. 

\begin{prop}\label{prop6.6}
Let notation be as before. Let $\varepsilon>0$. Then   
\begin{align*}
\mathcal{M}_{\const,2}^{\du,\heartsuit,\Reg}(1/2+it,1/2-it;\mathfrak{L})\ll 
\frac{C_{\fin}(\pi,\pi')\Delta_2}{\sqrt{\|\mathfrak{L}\|}}\prod_{v\mid\mathfrak{L}}\Big[|\lambda_{\pi'}(\mathfrak{p}_v^{\ell_v'})|+1\Big],
\end{align*}
where $\Delta_2:=\textbf{C}_{\infty}^{\varepsilon}C_{\infty}(\pi')^{\varepsilon}C_{\infty}(\omega\omega'^{-1} |\cdot|^{2it})^{\varepsilon}C_{\infty}(\pi'\times\overline{\pi}')^{\varepsilon}C_{\fin}(\pi,\pi')^{\varepsilon}$. Here the implied constant depends only on $\varepsilon$ and the base field $F$.  
\end{prop} 
\begin{proof} 
For $v\in \Sigma_F$, let $J_v^{\sharp}(s_1,s_2;\mathfrak{L})$ be the local factor defined in Lemma \ref{lem6.5}. 

Let $v\mid\mathfrak{L}$. Then $J_v^{\sharp}(1/2+it-s,1/2-it+2s;\mathfrak{L})$ is equal to 
\begin{align*}
\frac{(q_v^{\ell_v'/2}+\ell_v'q_v^{-\ell_v'/2}|\mathbb{F}_v^{\times}|)\omega_v(\varpi_v^{\ell_v'})\mathcal{RS}(1/2-it+s,1/2-it+2s)}{q_v^{(1/2-it+2s)\ell_v'}\prod_{n=1}^2\mathbb{L}_n(1/2+it-s,1/2-it+2s)^{-1}}+\frac{\overline{\lambda_{\pi'}(\mathfrak{p}_v^{\ell})}\omega_v(\varpi_v^{\ell_v'})}{q_v^{(1/2-it+2s)\ell_v'}}.
\end{align*}

Similar to the estimate \eqref{6.11} we have 
\begin{equation}\label{6.26}
J_v^{\sharp}(1/2+it-s,1/2-it+2s;\mathfrak{L})\ll\ell_v'^2 q_v^{-\frac{\ell_v'}{2}+7\ell_v'\varepsilon_0}\sum_{0\leq i\leq \ell_v'}|\lambda_{\pi'}(\mathfrak{p}_v^{i})|,
\end{equation}
where the implied is absolute. In conjunction with \eqref{eq6.15}, we obtain 
\begin{equation}\label{6.28}
\prod_{v\mid\mathfrak{L}}J_v^{\sharp}(1/2+it-s,1/2-it+2s;\mathfrak{L})\ll \|\mathfrak{L}\|^{-\frac{1}{2}+\varepsilon}\prod_{v\mid\mathfrak{L}}\Big[|\lambda_{\pi'}(\mathfrak{p}_v^{\ell_v'})|+1\Big].
\end{equation}
Here the implied constant depends only on $\varepsilon.$ 

Let $v\mid\mathfrak{M}\mathfrak{N}$. Then $J_v^{\sharp}(1/2+it-s,1/2-it+2s;\mathfrak{L})$ is equal to 
\begin{equation}\label{6.29.}
\prod_{v\mid\mathfrak{M}\mathfrak{N}}\frac{q_v^{m_v(2s-2it)}\mathbf{1}_{\text{$\omega_v\omega_v'^{-1}$ is unramified}}}{\zeta_v(1)\Vol(K_v[m_v])L_v(1-2it+4s,\omega_v\overline{\omega_v'})}\ll C_{\fin}(\pi,\pi')^{1+\varepsilon}.
\end{equation}

Let $v\mid\infty$. Recall that $J_v^{\sharp}(s_1,s_2;\mathfrak{L})$ (see Lemma \ref{lem6.5}) is defined as  
\begin{align*}
\int_{K_v}\int_{F_v^{\times}}\Big|W_{\pi_v'}\left(\begin{pmatrix}
	a_v\\
	&1
\end{pmatrix}k_v\right)\Big|^2|a_v|_v^{s_2-s_1}\omega_v\omega_v'^{-1}(a_v)d^{\times}a_v\mathcal{J}_v^{\heartsuit}(s_1,s_2;k_v)dk_v.
\end{align*}
where $\mathcal{J}_v^{\heartsuit}(s_1,s_2;k_v)$ is defined as in \eqref{J_v}. 

Similar to the proof of \eqref{eq6.14}, we have
\begin{equation}\label{6.29}
\int_{K_v}\int_{F_v^{\times}}\Big|W_{\pi_v'}\left(\begin{pmatrix}
	a_v\\
	&1
\end{pmatrix}k_v\right)\Big|^2|a_v|_v^{3s-2it}\omega_v\omega_v'^{-1}(a_v)d^{\times}a_v\ll C_v(\pi')^{3\varepsilon_0},
\end{equation}
where the implied constant depends only on $\varepsilon$, $F_v$, and $\varphi_v$.

To bound $\mathcal{J}_v^{\heartsuit}(1/2+it-s,1/2-it+2s;k_v)$, we appeal to \eqref{f6.29} in Lemma \ref{lemma6.5}: 
\begin{equation}\label{6.41}
\mathcal{J}_v^{\heartsuit}(1/2+it-s,1/2-it+2s;k_v)\ll \frac{\textbf{C}_v^{\, 2\Re(s)+\varepsilon}|L_v(2it-2s,\overline{\omega_v}\omega_v')|}{|L_v(1+2s-2it,\omega_v\omega_v'^{-1})|}.
\end{equation}

By the construction of $\varepsilon_0$ in \textsection\ref{sect4.4}, 
\begin{align*}
\min_{v\mid\infty}|it+\gamma_v/2-s|\geq \frac{1}{4[F:\mathbb{Q}](\log C_{\infty}(\pi')+\log C_{\fin}(\pi,\pi')+\log L)^{10}}.\tag{\ref{equa4.2}}
\end{align*}
where $\gamma_v$ is the spectral parameter of $\overline{\omega_v}\omega_v'$. Hence, by \eqref{equa4.2} and Stirling formula, we obtain from \eqref{6.41} that 
\begin{equation}\label{6.42}
\mathcal{J}_v^{\heartsuit}(1/2+it-s,1/2-it+2s;k_v)\ll \frac{(\log C_{\infty}(\pi')+\log C_{\fin}(\pi,\pi')+\log L)^{10}}{\textbf{C}_v^{\, -2\Re(s)-\varepsilon}},
\end{equation}
where the implied constant depends on $\varepsilon$, $F_v$ and the fixed $h_v$ (see \textsection\ref{sec4.2.}).

Combining \eqref{6.28}, \eqref{6.29.}, \eqref{6.29}, and  \eqref{6.42} with Lemma \ref{lem6.5}, we thus obtain 
\begin{equation}\label{6.16}
\mathcal{M}_{\const,2}^{\du,\heartsuit}(\cdots)\ll \frac{|\mathcal{L}_2(s;\pi')|C_{\fin}(\pi,\pi')^{1+\varepsilon}}{\textbf{C}_{\infty}^{-\varepsilon}\sqrt{\|\mathfrak{L}\|}}\prod_{v\mid\mathfrak{L}}\Big[|\lambda_{\pi'}(\mathfrak{p}_v^{\ell_v'})|+1\Big],
\end{equation}
where $\mathcal{M}_{\const,2}^{\du,\heartsuit}(\cdots):=\mathcal{M}_{\const,2}^{\du,\heartsuit}(1/2+it-s,1/2-it+2s;\mathfrak{L})$, and 
\begin{align*}
\mathcal{L}_2(s;\pi'):=\frac{L(1+3s-2it,\pi'\times\overline{\pi}'\otimes \omega\overline{\omega'})L(2it-2s,\overline{\omega}\omega')L(1-2it+4s,\omega\overline{\omega'})}{L(2+6s-4it,\omega^2\omega'^{-2})}.
\end{align*}

Appealing the well known upper bound for $L$-functions in the zero-free region, 
\begin{equation}\label{6.17}
\mathcal{L}_2(s;\pi')\ll C_{\fin}(\pi,\pi')^{\varepsilon}C_{\infty}(\omega^{-1}\omega' |\cdot|^{2it})^{\varepsilon}C_{\infty}(\pi'\times\overline{\pi}')^{\varepsilon},
\end{equation}
where the implied constant depends on $\varepsilon$ and $F$.

By defintion \eqref{equ6.2} and the triangle inequality,
\begin{equation}\label{6.34}
\mathcal{M}_{\const,2}^{\du,\heartsuit,\Reg}(1/2+it,1/2-it;\mathfrak{L})\ll \max_{|s|=\varepsilon_0}|\mathcal{M}_{\const,2}^{\du,\heartsuit}(\cdots)|,
\end{equation}
where $\mathcal{M}_{\const,2}^{\du,\heartsuit}(\cdots):=\mathcal{M}_{\const,2}^{\du,\heartsuit}(1/2+it-s,1/2-it+2s;\mathfrak{L})$.

Therefore, Proposition \ref{prop6.6} follows from  \eqref{6.16}, \eqref{6.17} and \eqref{6.34}.
\end{proof}

\subsection{Proof of Proposition \ref{prop6.1}}\label{sec6.3}
Let $\varepsilon>0$ be a tiny parameter. We derive from the decomposition \eqref{equa6.1}, along with Propositions \ref{prop6.4}, and \ref{prop6.6}, that    
\begin{equation}\label{6.46}
\mathcal{M}_{\const}^{\du,\heartsuit,\Reg}(\cdots)\ll 
\frac{C_{\fin}(\pi,\pi')\Delta_3}{\sqrt{\|\mathfrak{L}\|}}\prod_{v\mid\mathfrak{L}}\Big[|\lambda_{\pi'}(\mathfrak{p}_v^{\ell_v'})|+1\Big],
\end{equation}
where $\mathcal{M}_{\const}^{\du,\heartsuit,\Reg}(\cdots)=\mathcal{M}_{\const}^{\du,\heartsuit,\Reg}(1/2+it,1/2-it;\mathfrak{L})$, $\mathfrak{L}=\prod_{v\in \mathcal{P}}\mathfrak{p}_v^{\ell_v'}$ for some $0\leq \ell_v'\leq 2\ell_v$, and 
\begin{align*}
\Delta_3:=C_{\infty}(\pi\otimes|\cdot|^{it})^{\varepsilon}C_{\infty}(\pi')^{\varepsilon}C_{\infty}(\omega\omega'^{-1} |\cdot|^{2it})^{\varepsilon}C_{\fin}(\pi,\pi')^{\varepsilon}.
\end{align*}
It's worth noting that the above constant $n_0$ may differ from that in Propositions \ref{prop6.4} and \ref{prop6.6}, as $\Delta_3$ is simplified from $\Delta_1$ and $\Delta_2$. Define 
\begin{align*}
\text{Diag}:=\sum_{\substack{v,v'\in \mathcal{P}\\ v\neq v'}}\alpha_v\overline{\alpha_{v'}}\omega_{v'}^{-1}(\varpi_{v'}^{\ell_{v'}})\mathcal{M}_{\const}^{\du,\heartsuit,\Reg}(1/2+it,1/2-it;\mathfrak{p}_v^{\ell_v}\mathfrak{p}_{v'}^{\ell_{v'}}).
\end{align*}

As a consequence of \eqref{6.46}, we have the following 
\begin{equation}\label{6.47}
\text{Diag}\ll C_{\fin}(\pi,\pi')\Delta_3\Bigg[\sum_{v\in \mathcal{P}}\frac{|\alpha_v|\cdot (|\lambda_{\pi'}(\mathfrak{p}_v^{\ell_v})|+1)}{q_v^{l_v/2}}\Bigg]^2.
\end{equation}

Executing Cauchy-Schwarz inequality, Iwaniec's estimate \cite[Lemma 1]{Iwa92} on the Rankin--Selberg convolution, alongside Kim-Sarnak's bound $0\leq \vartheta\leq 7/64$ (to handle the case that $\ell_v\geq 2$), we obtain from \eqref{6.47} that 
\begin{equation}\label{6.48}
\text{Diag}\ll C_{\fin}(\pi,\pi')\Delta_3^{1+\varepsilon}\sum_{\substack{v\in \mathcal{P}}}|\alpha_v|_v^2.
\end{equation}

Define the contribution from off-diagonal by 
\begin{align*}
\text{Off-diag}:=\sum_{v\in \mathcal{P}}\sum_{n_v=0}^{\ell_v}|\alpha_v|_v^2\omega_v^{-1}(\varpi_v^{\ell_v-n_v})\mathcal{M}_{\const}^{\du,\heartsuit,\Reg}(1/2+it,1/2-it;\mathfrak{p}_v^{2\ell_v-2n_v}).
\end{align*}

Taking advantage of the estimate \eqref{6.46},
\begin{equation}\label{6.53}
\text{Off-diag}\ll C_{\fin}(\pi,\pi')\Delta_3\sum_{v\in \mathcal{P}}\sum_{n_v=0}^{\ell_v}\frac{|\alpha_v|_v^2\cdot (|\lambda_{\pi'}(\mathfrak{p}_v^{2l_v-2n_v})|+1)}{q_v^{2l_v-2n_v}}.
\end{equation}

Notice that 
\begin{equation}\label{6.51}
\sum_{v\in \mathcal{P}}\sum_{n_v=0}^{\ell_v}\frac{|\alpha_v|_v^2\cdot (|\lambda_{\pi'}(\mathfrak{p}_v^{2l_v-2n_v})|+1)}{q_v^{2l_v-2n_v}}\ll\sum_{v\in \mathcal{P}}|\alpha_v|_v^2+ \sum_{j=1}^{2\ell_v}\sum_{v\in \mathcal{P}}\frac{|\alpha_v|_v^2|\lambda_{\pi'}(\mathfrak{p}_v^{2j})|}{q_v^{2j}}.
\end{equation}

For $j\geq 1$, $|\lambda_{\pi'}(\mathfrak{p}_v^{2j})|\ll q_v^{7j/32}$. So
\begin{equation}\label{6.52}
\sum_{j=1}^{2\ell_v}\sum_{v\in \mathcal{P}}\frac{|\alpha_v|_v^2|\lambda_{\pi'}(\mathfrak{p}_v^{2j})|}{q_v^{2j}}\ll \sum_{j=1}^{2\ell_v}\sum_{v\in \mathcal{P}}\frac{|\alpha_v|_v^2}{q_v^{57j/32}}\ll \sum_{v\in \mathcal{P}}|\alpha_v|_v^2.
\end{equation}

Substituting \eqref{6.51} and \eqref{6.52} into \eqref{6.53} we then derive that 
\begin{equation}\label{6.50}
\text{Off-diag}\ll C_{\fin}(\pi,\pi')\Delta_3^{1+\varepsilon}\sum_{\substack{v\in \mathcal{P}}}|\alpha_v|_v^2.
\end{equation}

Therefore, Proposition \ref{prop6.1} follows from \eqref{6.48} and \eqref{6.50}.

\section{Nonarchimedean Integrals on the Dual Side \RNum{1}}\label{sec8} 
In this section we evaluate the nonarchimedean local integrals appearing in the dual-side
contributions $\mathcal{M}_{\cusp}^{\du}(s_1,s_2)$ and $\mathcal{M}_{\Eis}^{\du}(s_1,s_2)$
at all finite places where the representation $\pi_v'$ is a principal series.  
We derive explicit formulas for the corresponding local factors.

\subsection{A Local Rankin--Selberg Integral}
\subsubsection{Notation}\label{7.1.1}
We fix some notation used in this section. Let $\sigma=\otimes_v\sigma_v$ be an automorphic generic representation of $\mathrm{PGL}(2)/F$.  For $v<\infty$, Let $r_v:=r_{\sigma_v}$ be the conductor exponent of $\sigma_v$ (see \textsection\ref{sec4.1}). Let $\chi_{i}$, $i=1,2$, be unitary characters of $F_v^{\times}$. Let $s_1, s_2\in \mathbb{C}$. Let $\sigma_v'=\pi_{\chi_1,\chi_2,s_1}$ be the principal series induced by $\chi_1|\cdot|_v^{s_1-1/2}\otimes \chi_2|\cdot|_v^{-s_1+1/2}$. Let $m_v=\max\{r_{\pi_v},r_{\pi_v'}\}$ be the exponent defined in  \textsection\ref{sec4.1}, and set 
\begin{equation}\label{m_v}
m_v^*:=r_{\overline{\chi_1}\chi_2}.
\end{equation}

For our purposes, we consider $\sigma_v'=\pi_v'$ if $\pi_v'$ is the principal series $\pi_{\chi_1,\chi_2,s'}$ for some $|\Re(s')-1/2|\leq \vartheta$ (and $s_1$ deforms to $s'$), or $\sigma_v'=\pi_{\mathbf{1},\omega_v\omega_v'^{-1}, s_1}$. Here $\pi'_v$ and $\pi_{\mathbf{1},\omega_v\omega_v'^{-1}, s_1}$ are the $v$-th components of the cuspidal representation $\pi'$ and the Eisenstein series $\pi_{\mathbf{1},\omega\omega'^{-1}, s_1}$ defined in \textsection\ref{sec4}, respectively. 

For $\Re(s)\gg 1$ and $x_v\in G(F_v)$, we define  
\begin{equation}\label{eq7.1}
	h_v(x_v,s):=\chi_1(\det x_v)|\det x_v|_v^{s}\int_{F_v^{\times}}\Phi_v^*\left((0,t)x_v\right)\chi_1\chi_2^{-1}(t)|t|_v^{2s}d^{\times}t,
\end{equation}
where $\Phi_v^*$ is defined by 
\begin{equation}\label{phi*}
	\Phi_v^*(t_{1,v},t_{2,v}):=\begin{cases}
		\frac{\overline{\chi_1}\chi_2(t_{2,v})\mathbf{1}_{\mathfrak{p}_v^{m_v}}(t_{1,v})\mathbf{1}_{\mathcal{O}_v^{\times}}(t_{2,v})}{\Vol(K_v[m_v])},\ &\ \text{if $m_v\geq 1$,}\\
		\mathbf{1}_{\mathcal{O}_v}(t_{1,v})\mathbf{1}_{\mathcal{O}_v}(t_{2,v}),\ &\ \text{if $m_v=0$}.
	\end{cases}
\end{equation}
Note that when $\chi_1=\mathbf{1}$ and $\chi_2=\omega_v\omega_v'^{-1}$, the function $\Phi_v^*$ boils down to $\Phi_v$, which is defined by \eqref{equa4.5} and \eqref{equ4.9}.  

For $v<\infty$, let $\ell_v'\geq 0$ be an integer such that $\ell_v'm_v=0$ for all $v<\infty$. For $\Re(s)\gg 1$ and $x_v\in G(F_v)$, we define 
\begin{align*}
	h_v(x_v,s;\ell_v'):=\chi_1(\det x_v\boldsymbol{d}_v)|\det x_v\boldsymbol{d}_v|_v^{s}\int_{F_v^{\times}}\Phi_v^*\left((0,t)x_v\boldsymbol{d}_v\right)\chi_1\chi_2^{-1}(t)|t|_v^{2s}d^{\times}t.
\end{align*}
Here $\boldsymbol{d}_v=\begin{pmatrix}
	1\\
	&\varpi_v^{\ell_v'}
\end{pmatrix}$ as defined in \textsection\ref{sec4.3}. Note that $h_v(x_v,s;0)=h_v(x_v,s)$, which is defined by \eqref{eq7.1}.

Let $L_v(s,\sigma_v)=\sum_{n\geq 0}\lambda_{\sigma}(\mathfrak{p}_v^n)q_v^{-ns}$ and $L_v(s,\sigma_v')=\sum_{n\geq 0}\lambda_{\sigma'}(\mathfrak{p}_v^n)q_v^{-ns}$ be the local $L$-factor of $\sigma_v$ and $\sigma_v'$, respectively.

%Then $\overline{h_v(x_v,\overline{s_2};\ell_v')}$ is a holomorphic function of $s_2$ in $\Re(s_2)>0$ and admits a meromorphic continuation to $s_2\in \mathbb{C}$. 

\subsubsection{A Local Period Integral} Let $v<\infty$. Let $W_v$ be a nonzero vector in the Whittaker models of $\sigma_v$  relative to $\theta_v$. Let $W_v^{\circ}$ be a local new vector in the Whittaker models of $\sigma_v$  relative to $\theta_v$. Let 
\begin{equation}\label{7.1}
	W_v'(x_v):=\int_{N(F_v)}h_v(wu_vx_v,s_1)\theta_v(u_v)du_v
\end{equation}
be the Whittake vector in the Whittaker model of $\sigma_v'$ relative to $\overline{\theta}_v$. In particular, if $\chi_1$ is unramified, then $W_v'$ is a new vector.

Let $\Re(s_2-s_1)\gg 1$ and $\Re(s_2)\gg 1$.  Define  
\begin{equation}\label{equa7.4}
	\Psi_v(s_1,s_2;\ell_v'):=\int_{N(F_v)Z(F_v)\backslash G(F_v)}W_v\left(x_v\right)W_v'(x_v)\overline{h_v(x_v,\overline{s_2};\ell_v')}dx_v.
\end{equation}

By definition, $\Psi_v(s_1,s_2)\propto L_v(s_2,\sigma_v\times\sigma_v'\otimes\overline{\chi_1})L_v(2s_1,\chi_1\chi_2^{-1})$, and thus admits a meromorphic continuation to $(s_1,s_2)\in \mathbb{C}^2$. 

\subsubsection{Explicit orthonormal basis}\label{7.1.3.}
Let $v<\infty$ and $l\geq 0$ be an integer. Let $\mathcal{W}(\sigma_v,\theta_v)$ be the Whittaker model of $\sigma_v$. Let $\mathcal{W}(\sigma_v,\theta_v)^{K_v[l]}$ be the subspace of right $K_v[l]$-invariant vectors in $\mathcal{W}(\sigma_v,\theta_v)$. By \cite[Lemma 9]{BM15}, the space $\mathcal{W}(\sigma_v,\theta_v)^{K_v[l]}$ has an orthonormal basis $\{W_{v,n_v}:\ 0\leq n_v \leq l-r_v\}$, where
\begin{equation}\label{eq7.5}
	W_{v,n_v}=\sum_{i=0}^{n_v}\xi_{\sigma_v}(\mathfrak{p}_v^{i},\mathfrak{p}_v^{n_v})q_v^{\frac{i-n_v}{2}}\sigma_v\left(\begin{pmatrix}
		1\\
		&\varpi_v^{i}
	\end{pmatrix}\right)W_v^{\circ}.
\end{equation}
Moreover, $\xi_{\sigma_v}(\mathfrak{p}_v^{i},\mathfrak{p}_v^{n_v})\ll 1$, where the implied constant is absolute. 

\subsubsection{The $K_v$-isotypical of $\sigma_v$}\label{7.1.3}
Analyzing the definition of $\Phi_v$, we obtain the following. 
\begin{itemize}
	\item At $v\nmid\mathfrak{MN}$, we have $\Psi_v(s_1,s_2;\ell_v')\equiv 0$ if $\int_{K_v[\ell_v']}\sigma_v(k_v)W_vdk_v\neq 0$.
	\item At $v\mid\mathfrak{MN}$,   we have $\Psi_v(s_1,s_2;\ell_v')\equiv 0$ if $\int_{K_v[m_v]}\sigma_v(k_v)W_vdk_v\neq 0$. 
\end{itemize}

Henceforth for each $v<\infty$ we may assume $W_v=W_{v,n_v}$ for some $0\leq n_v\leq \ell_v'+m_v-r_v$. Here $W_{v,n_v}$ is defined by \eqref{eq7.5}.

\subsection{Explicit Expression of the Whittaker Function $W_v'$}\label{sec7.1}
Recall 
\begin{align*}
	W_v'(x_v):=\int_{N(F_v)}h_v(wu_vx_v,s_1)\theta_v(u_v)du_v.\tag{\ref{7.1}}
\end{align*}

It is worth noting that $W_v'(\cdot)$ is also a function of $s_1$. For the sake of simplicity, we adopt the notation $W_v'(x_v)$ without explicitly involving $s_1.$

\subsubsection{Simplification of $W_v'(\diag(a_v,1)k_v)$}\label{7.1.1}
Let $\Re(s_1)>1$ and $x_v\in G(F_v)$. From the definition \eqref{eq7.1} we have  
\begin{align*}
	h_v(wu_vx_v,s_1)=\chi_1(-\det x_v)|\det x_v|_v^{s_1}\int_{F_v^{\times}}\Phi_v^*((0,t_v)wu_vx_v)\chi_1\chi_2^{-1}(t_v)|t_v|_v^{2s_1}d^{\times}t_v.
\end{align*}

Substituting this into \eqref{7.1} we obtain that $W_v'(x_v)$ is equal to 
\begin{align*}
	\frac{\chi_1(-1)\chi_1(\det x_v)}{|\det x_v|_v^{-s_1}}\int_{N(F_v)}\int_{F_v^{\times}}\Phi_v^*((0,t_v)wu_vg_v)\chi_1\chi_2^{-1}(t_v)|t_v|_v^{2s_1}\theta_v(u_v)d^{\times}t_vdu_v.
\end{align*}

Let $x_v=\diag(a_v,1)k_v$, where $a_v\in F_v^{\times}$ and $k_v\in K_v$. After a change of variable $W_v'(\diag(a_v,1)k_v)$ becomes  
\begin{equation}\label{eq7.2}
	\frac{\chi_1(-\det k_v)}{\overline{\chi_2}(a_v)|a_v|_v^{s_1-1}}\int_{F_v^{\times}}\int_{F_v}\Phi_v^*((t_v,b_v)
	k_v)\psi_v(a_vt_v^{-1}b_v)db_v\chi_1\overline{\chi_2}(t_v)|t_v|_v^{2s_1-1}d^{\times}t_v.
\end{equation}

\subsubsection{Unramified Whittaker functions} 
\begin{lemma}\label{lem7.2}
	Let notation be as before. Let $v<\infty$,  $v\nmid \mathfrak{MN}$. Let $a_v\in F_v^{\times}$, and $k_v\in K_v$. Suppose that $e_v(a_v)=j_v\geq 0$ (see \textsection\ref{2.1.1}). Then 
	\begin{equation}\label{eq7.3}
		W_v'(\diag(a_v,1)k_v)=\frac{\chi_1(-1)}{|a_v|_v^{s_1-1}}\cdot \frac{\chi_2(\varpi_v^{j_v})-\chi_1(\varpi_v^{j_1+1})\overline{\chi_2}(\varpi_v)q_v^{-(2s_1-1)(j_v+1)}}{1-\chi_1\overline{\chi_2}(\varpi_v)q_v^{-(2s_1-1)}},
	\end{equation}
\end{lemma}
\begin{proof}
	Substituting the definition of $\Phi_v^*$ (see \eqref{phi*}) into \eqref{eq7.2} we obtain 
	\begin{align*}
		W_v'(\diag(a_v,1)k_v)=\frac{\chi_1(-1)\chi_2(a_v)}{|a_v|_v^{s_1-1}}\int \int_{\mathcal{O}_v}\psi_v(a_vt_v^{-1}b_v)db_v\chi_1\overline{\chi_2}(t_v)|t_v|_v^{2s_1-1}d^{\times}t_v,
	\end{align*}
	where $t_v$ ranges over $\mathcal{O}_v-\{0\}$. Analyzing the integral relative to $b_v\in \mathcal{O}_v$ by orthogonality, we derive that 
	\begin{align*}
		W_v'(\diag(a_v,1)k_v)=\chi_1(-1)\chi_2(a_v)|a_v|_v^{1-s_1}\sum_{0\leq i\leq j_v}\chi_1\overline{\chi_2}(\varpi_v^i)q_v^{(1-2s_1)i}.
	\end{align*}
	Therefore, \eqref{eq7.3} follows from the above expression by computing the geometric series. 
\end{proof}

\subsubsection{Ramified Whittaker Functions}
\begin{lemma}\label{lem7.3}
	Let $v<\infty$ be such that $m_v\geq 1$. Let $a_v\in F_v^{\times}$, and  $k_v\in K_v[m_v]$. Suppose that $e_v(a_v)=j_v\geq 0$. Then 
	\begin{equation}\label{7.2}
		W_v'(\diag(a_v,1)k_v)=\frac{|a_v|_v^{1-s_1}\chi_1(-a_vk_{11})\chi_2(k_{22})}{\Vol(K_v[m_v])}\sum_{i=m_v}^{j_v+m_v}\frac{\chi_1\overline{\chi_2}(\varpi_v^{i})G_v(a_v,i)}{q_v^{(2s_1-1)i}},
	\end{equation}
	where $k_{ij}\in \mathcal{O}_v^{\times}$ is the $(i,j)$-th entry of $k_v$, and   
\begin{equation}\label{7.3}
G_v(a_v,i):=\int_{\mathcal{O}_v^{\times}}\overline{\chi_1}\chi_2(a_vb_v)\psi_v(\varpi_v^{-i}a_vb_v)db_v.
\end{equation}
Moreover, we have the following scenarios: 
\begin{itemize}
\item Suppose that $m_v^*:=r_{\overline{\chi_1}\chi_2}\geq 1$ (see \eqref{m_v}). Then \eqref{7.2} becomes 
\begin{equation}\label{f7.6}
W_v'(\diag(a_v,1)k_v)=\frac{\chi_1(-a_v\det k_v)\chi_1\overline{\chi_2}(\varpi_v^{m_v^*}k_{22}^{-1})\mathcal{G}_{\overline{\chi_1}\chi_2,\psi_v}\mathbf{1}_{j_v\geq m_v-m_v^*}}{\Vol(K_v[m_v])|a_v|_v^{s_1-1}q_v^{(2s_1-1)(j_v+m_v^*)+m_v^*}},
\end{equation}
where $\mathcal{G}_{\overline{\chi_1}\chi_2,\psi_v}$ is the Gauss sum with respect to $\overline{\chi_1}\chi_2$ and $\psi_v$.
		
\item Suppose that $m_v^*:=r_{\overline{\chi_1}\chi_2}=0$. Then \eqref{7.2} becomes 
\begin{equation}\label{f7.7}
W_v'(\diag(a_v,1)k_v)=W_v'^{(1)}(\diag(a_v,1)k_v)-W_v'^{(2)}(\diag(a_v,1)k_v),
\end{equation}
where $W_v'^{(1)}(\diag(a_v,1)k_v)$ is defined by 
\begin{align*}
\frac{\chi_1(-a_vk_{11})\chi_2(k_{22})\mathbf{1}_{j_v\geq m_v}}{\Vol(K_v[m_v])|a_v|_v^{s_1-1}\zeta_v(1)}\cdot \frac{\overline{\chi_1}\chi_2(\varpi_v^{j_v-m_v})q_v^{(1-2s_1)m_v}-\chi_1\overline{\chi_2}(\varpi_v)q_v^{(1-2s_1)(j_v+1)}}{1-\chi_1\overline{\chi_2}(\varpi_v)q_v^{-(2s_1-1)}};
\end{align*}
and 
\begin{align*}
W_v'^{(2)}(\diag(a_v,1)k_v):=\frac{|a_v|_v^{1-s_1}\chi_1(-a_vk_{11})\chi_2(k_{22})}{\Vol(K_v[m_v])}\cdot \frac{\chi_1\overline{\chi_2}(\varpi_v)\mathbf{1}_{j_v\geq m_v-1}}{q_v^{(2s_1-1)(j_v+1)+1}}.
\end{align*}
\end{itemize}
\end{lemma}
\begin{proof}
	Let $k_v\in K_v[m_v]$. Considering the construction of $\Phi_v^*$ in \eqref{phi*}, the integral 
	\begin{align*}
		\int_{F_v}\Phi_v^*((t_v,b_v)
		k_v)\psi_v(a_vt_v^{-1}b_v)db_v
	\end{align*}
	is equal to 
	\begin{align*}
		\Vol(K_v[m_v])^{-1}\overline{\chi_1}\chi_2(k_{22})\mathbf{1}_{t_v\in \mathfrak{p}_v^{m_v}}\int_{\mathcal{O}_v^{\times}}\overline{\chi_1}\chi_2(b_v)\psi_v(a_vt_v^{-1}b_v)db_v.
	\end{align*}
	
	Therefore, by \eqref{eq7.2},  $W_v'(\diag(a_v,1)k_v)$ is equal to 
	\begin{align*}
		&\frac{\chi_1(-\det k_v)\overline{\chi_1}\chi_2(k_{22})}{\Vol(K_v[m_v])\overline{\chi_2}(a_v)|a_v|_v^{s_1-1}}\int_{\mathfrak{p}_v^{m_v}-\{0\}}\int_{\mathcal{O}_v^{\times}}\overline{\chi_1}\chi_2(t_v^{-1}b_v)\psi_v(a_vt_v^{-1}b_v)db_v
		|t_v|_v^{2s_1-1}d^{\times}t_v.
	\end{align*}
	
	Parametrizing $t_v$ as $\varpi_v^{i}\gamma_v$, $\gamma_v\in \mathcal{O}_v^{\times}$, along with the changing variable $b_v\mapsto \gamma_vb_v$, we then obtain 
	\begin{align*}
		W_v'(\diag(a_v,1)k_v)=\frac{|a_v|_v^{1-s_1}\chi_1(-a_v\det k_v)\overline{\chi_1}\chi_2(k_{22})}{\Vol(K_v[m_v])}\sum_{i\geq m_v}\frac{\chi_1\overline{\chi_2}(\varpi_v^{i})G_v(a_v,i)}{q_v^{(2s_1-1)i}}.
	\end{align*}
	Moreover, by definition \eqref{7.3}, $G_v(a_v,i)\equiv 0$ unless $i\leq j_v+m_v$, where we have used the assumption that $m_v\geq 1$ as $v\mid\mathfrak{M}\mathfrak{N}$. Therefore, the formula \eqref{7.2} holds, and defines a holomorphic function of $s_1\in \mathbb{C}$. 
	
	We can compute $G_v(a_v,i)$ to establish a refined formula for $W_v'(\diag(a_v,1)k_v)$ according to the situations determined by $\chi_1$ and $\chi_2$  as follows.
	
	\begin{itemize}
		\item Suppose that $m_v^*\geq 1$,  namely, $\overline{\chi_1}\chi_2$ is ramified. Hence, 
		\begin{align*}
			G_v(a_v,i)=\begin{cases}
				\overline{\chi_1}\chi_2(\varpi_v^{j_v})q_v^{-m_v^*}\mathcal{G}_{\overline{\chi_1}\chi_2,\psi_v}, &\text{if $i=j_v+m_v^*$};\\
				0, &\text{if $i\neq j_v+m_v^*$},
			\end{cases}	
		\end{align*}
		where $\mathcal{G}_{\overline{\chi_1}\chi_2,\psi_v}$ is the Gauss sum. Therefore,
\begin{align*}
\sum_{i\geq m_v}\frac{\chi_1\overline{\chi_2}(\varpi_v^{i})G_v(a_v,i)}{q_v^{(2s_1-1)i}}=\chi_1\overline{\chi_2}(\varpi_v^{m_v^*})q_v^{(1-2s_1)(j_v+m_v^*)-m_v^*}\mathcal{G}_{\overline{\chi_1}\chi_2,\psi_v}\mathbf{1}_{j_v\geq m_v-m_v^*},
\end{align*}
from which the equality \eqref{f7.6} follows. 
		
\item Suppose that $m_v^*=0$. Then $\overline{\chi_1}\chi_2$ is unramified. So 
\begin{equation}\label{eq7.8}
G_v(a_v,i)=\begin{cases}
\overline{\chi_1}\chi_2(\varpi_v^{j_v})(1-q_v^{-1}),\ &\text{if $i\leq j_v$};\\
-\overline{\chi_1}\chi_2(\varpi_v^{j_v})q_v^{-1},\ &\text{if $i=j_v+1$};\\
0,\ &\text{if $i\geq j_v+2$}.
\end{cases}	
\end{equation}
		
Substituting this into \eqref{7.2} and computing the geometric series we then obtain  \eqref{f7.7}. 
	\end{itemize}
	
Therefore, Lemma \ref{lem7.3} holds from the above discussions. 
\end{proof}

\subsection{Calculation of Local Period Integrals}

\begin{lemma}\label{lemma7.4}
	Let notation be as before. Let $v<\infty$ be such that $m_v=\ell_v'=0$. Then $\Psi_v(s_1,s_2;\ell_v')$ converges absolutely in the region $\Re(s_2)-|\Re(s_1)|\gg 1$, and it admits a meromorphic continuation to $\mathbb{C}^2$ via the following 
	\begin{align*}
		\Psi_v(s_1,s_2;\ell_v')=W_v(I_2)W_v'(I_2)L_v(2s_1,\chi_1\chi_2^{-1})L_v(s_2,\sigma_v\times\sigma_v'\otimes\overline{\chi_1}).
	\end{align*}
\end{lemma}
\begin{proof}
	This follows from the unramified computation in \cite{JPSS79a}. 
\end{proof}

\begin{lemma}\label{lem7.6}
Let notation be as before. Let $v<\infty$ be such that $m_v^*\geq 1$. Let   
\begin{equation}\label{eq7.15}
W_v=W_{v,n_v}=\sum_{i=0}^{n_v}\xi_{\sigma_v}(\mathfrak{p}_v^{i},\mathfrak{p}_v^{n_v})q_v^{\frac{i-n_v}{2}}\sigma_v\left(\begin{pmatrix}
			1\\
			&\varpi_v^{i}
		\end{pmatrix}\right)W_v^{\circ}
\end{equation}
for some $0\leq n_v\leq m_v-r_v.$ Then $\Psi_v(s_1,s_2;\ell_v')$ converges absolutely in the region $\Re(s_2)-|\Re(s_1)|\gg 1$, and it admits a meromorphic continuation to $\mathbb{C}^2$ via the following   
\begin{equation}\label{equation7.15}
\frac{\Psi_v(s_1,s_2;\ell_v')}{L_v(s_2,\sigma_v\times\sigma_v'\otimes\overline{\chi_1})}=\frac{\chi_1\overline{\chi_2}(\varpi_v^{m_v^*})\mathcal{G}_{\overline{\chi_1}\chi_2,\psi_v}W_v^{\circ}(I_2)}{\overline{\chi_1}(-1)\Vol(K_v[m_v])q_v^{2s_1m_v^*}}\cdot \sum_{j=1}^2J_v^{(j)}(s_1,s_2),
\end{equation}
where 
\begin{align*}
J_v^{(1)}(s_1,s_2):=&\sum_{i=m_v-m_v^*+1}^{n_v}\xi_{\sigma_v}(\mathfrak{p}_v^{i},\mathfrak{p}_v^{n_v})q_v^{\frac{i-n_v}{2}}q_v^{-(s_1+s_2-1)i},
\end{align*}
and $J_v^{(2)}(s_1,s_2)$ is defined by 
\begin{align*} 
\sum_{i=0}^{\min\{m_v-m_v^*,n_v\}}\frac{\xi_{\sigma_v}(\mathfrak{p}_v^{i},\mathfrak{p}_v^{n_v})\big[\lambda_{\sigma}(\mathfrak{p}_v^{m_v-m_v^*-i})+\lambda_{\sigma}(\mathfrak{p}_v^{m_v-m_v^*-1-i})q_v^{-(s_2+s_1-1/2)}\big]}{q_v^{(m_v-m_v^*)(s_2+s_1-1/2)}q_v^{n_v/2-i}}.
\end{align*}
\end{lemma}
\begin{proof}
Recall the definition 
\begin{align*}
\Psi_v(s_1,s_2;\ell_v'):=\int_{N(F_v)Z(F_v)\backslash G(F_v)}W_v\left(x_v\right)W_v'(x_v)\overline{h_v(x_v,\overline{s_2};\ell_v')}dx_v.\tag{\ref{equa7.4}}
\end{align*}
	
Let $x_v=\diag(a_v,1)k_v$. By the definition \eqref{phi*} we have
\begin{equation}\label{equ7.12}
\overline{h_v(x_v,\overline{s_2};\ell_v')}=\Vol(K_v[m_v])^{-1}\overline{\chi_1}(a_v\det k_v)\chi_1\overline{\chi_2}(k_{22})|a_v|_v^{s_2}\mathbf{1}_{k_v\in K_v[m_v]}.
\end{equation}
	
Substituting \eqref{equ7.12} and \eqref{f7.6} into \eqref{equa7.4}, in conjunction with the Iwasawa decomposition, we derive that $\Psi_v(s_1,s_2;\ell_v')$ is equal to  
	\begin{equation}\label{f7.12}
		\frac{\chi_1(-1)\chi_1\overline{\chi_2}(\varpi_v^{m_v^*})\mathcal{G}_{\overline{\chi_1}\chi_2,\psi_v}}{\Vol(K_v[m_v])q_v^{2s_1m_v^*}}\cdot \sum_{j_v\geq m_v-m_v^*} W_v\left(\begin{pmatrix}
			\varpi_v^{j_v}\\
			&1
		\end{pmatrix}\right)q_v^{-(s_1+s_2-1)j_v},
	\end{equation}
	which converges absolutely in $\Re(s_1+s_2)\gg 1$. Here we also utilize the assumption that $W_v'$ is right $K_v[m_v]$-invariant, see \textsection\ref{7.1.3}. Let 
	\begin{equation}\label{eq7.18}
		\textbf{S}(m_v^*,i):=\sum_{j_v\geq m_v-m_v^*} W_v^{\circ}\left(\begin{pmatrix}
			\varpi_v^{j_v}\\
			&\varpi_v^i
		\end{pmatrix}\right)q_v^{-(s_1+s_2-1)j_v},
	\end{equation}
	where $\Re(s_1+s_2)\gg 1$. 
	Let $j^*:=\max\{m_v-m_v^*-i,0\}$. Then 
	\begin{align*}
		\textbf{S}(m_v^*,i)=q_v^{-(s_1+s_2-1)i}\sum_{j_v\geq j^*} W_v^{\circ}\left(\begin{pmatrix}
			\varpi_v^{j_v}\\
			&1
		\end{pmatrix}\right)q_v^{-j_v(s_1+s_2-1)}.
	\end{align*}
	
Substituting \eqref{eq7.15} into  \eqref{f7.12}, we obtain 
\begin{equation}\label{7.16}
\Psi_v(s_1,s_2;\ell_v')=\frac{\chi_1(-1)\chi_1\overline{\chi_2}(\varpi_v^{m_v^*})\mathcal{G}_{\overline{\chi_1}\chi_2,\psi_v}}{\Vol(K_v[m_v])q_v^{2s_1m_v^*}}\cdot \sum_{i=m_v-m_v^*}^{n_v}\widetilde{\xi}_{\sigma_v}(\mathfrak{p}_v^{i},\mathfrak{p}_v^{n_v})\textbf{S}(m_v^*,i),
\end{equation}
where $\widetilde{\xi}_{\sigma_v}(\mathfrak{p}_v^{i},\mathfrak{p}_v^{n_v}):=\xi_{\sigma_v}(\mathfrak{p}_v^{i},\mathfrak{p}_v^{n_v})q_v^{\frac{i-n_v}{2}}$.
	
	Since $\overline{\chi_1}\chi_2$ is ramified, then $L_v(s_2-s_1+1/2,\sigma_v\otimes \overline{\chi_1}\chi_2)\equiv 1$. It follows from  
	\begin{align*}
		L_v(s_2,\sigma_v\times\sigma_v'\otimes\overline{\chi_1})=L_v(s_2+s_1-1/2,\sigma_v)L_v(s_2-s_1+1/2,\sigma_v\otimes \overline{\chi_1}\chi_2)
	\end{align*}
	that $L_v(s_2,\sigma_v\times\sigma_v'\otimes\overline{\chi_1})=L_v(s_2+s_1-1/2,\sigma_v).$ We consider the following scenarios.

\begin{itemize}
\item Suppose that $i> m_v-m_v^*$. By Casselman-Shalika formula, 
\begin{equation}\label{7.17}
\textbf{S}(m_v^*,i)=\frac{W_v^{\circ}(I_2)L_v(s_2+s_1-1/2,\sigma_v)}{q_v^{(s_1+s_2-1)i}}=\frac{W_v^{\circ}(I_2)L_v(s_2,\sigma_v\times\sigma_v'\otimes\overline{\chi_1})}{q_v^{(s_1+s_2-1)i}}.
\end{equation}
In conjunction with \eqref{7.16} we obtain 
\begin{align*}
\frac{\Psi_v(s_1,s_2;\ell_v')}{L_v(s_2,\sigma_v\times\sigma_v'\otimes\overline{\chi_1})}=\frac{\chi_1(-1)\chi_1\overline{\chi_2}(\varpi_v^{m_v^*})\mathcal{G}_{\overline{\chi_1}\chi_2,\psi_v}W_v^{\circ}(I_2)}{\Vol(K_v[m_v])q_v^{2s_1m_v^*}} \sum_{i=m_v-m_v^*+1}^{n_v}\frac{\widetilde{\xi}_{\sigma_v}(\mathfrak{p}_v^{i},\mathfrak{p}_v^{n_v})}{q_v^{(s_1+s_2-1)i}}.
\end{align*}		
		
\item Suppose that $i\leq m_v-m_v^*$. Similar to \eqref{7.17}, using Hecke relations, we have  
\begin{equation}\label{7.18}
\textbf{S}(m_v^*,i)=\frac{W_v^{\circ}(I_2)\big[\lambda_{\sigma}(\mathfrak{p}_v^{j^*})+\lambda_{\sigma}(\mathfrak{p}_v^{j^*-1})q_v^{-(s_2+s_1-1/2)}\big]}{q_v^{j^*(s_2+s_1-1/2)}q_v^{(s_1+s_2-1)i}}L_v(s_2+s_1-1/2,\sigma_v).
\end{equation}
		
Combining \eqref{7.16} with \eqref{7.18} leads to 
\begin{align*}
\frac{\Psi_v(s_1,s_2;\ell_v')}{L_v(s_2,\sigma_v\times\sigma_v'\otimes\overline{\chi_1})}=&\frac{\chi_1(-1)\chi_1\overline{\chi_2}(\varpi_v^{m_v^*})\mathcal{G}_{\overline{\chi_1}\chi_2,\psi_v}W_v^{\circ}(I_2)}{\Vol(K_v[m_v])q_v^{2s_1m_v^*}}\\
&\sum_{i=0}^{m_v-m_v^*-1}\xi_{\sigma_v}(\mathfrak{p}_v^{i},\mathfrak{p}_v^{n_v})\frac{\big[\lambda_{\sigma}(\mathfrak{p}_v^{j^*})+\lambda_{\sigma}(\mathfrak{p}_v^{j^*-1})q_v^{-(s_2+s_1-1/2)}\big]}{q_v^{j^*(s_2+s_1-1/2)}q_v^{(s_1+s_2-1)i}},
\end{align*}
where $j^*=m_v-m_v^*-i$. 
\end{itemize}
	
Combining the above discussions, we then obtain \eqref{equation7.15}. 
\end{proof}

We proceed to consider the scenario that $m_v\geq 1$ but $m_v^*=0$. By \eqref{7.2} in Lemma \ref{lem7.3}, and the formula \eqref{equ7.12}, the function $\Psi_v(s_1,s_2;\ell_v')$ simplifies to 
\begin{align*}
	\int_{K_v[m_v]}\int_{F_v^{\times}}\overline{\chi_1}(a_v\det k_v)\chi_1\overline{\chi_2}(k_{22})W_v\left(x_v\right)W_v'\left(x_v\right)|a_v|_v^{s_2-1}d^{\times}a_vdk_v,
\end{align*}
where $x_v:=\diag(a_v,1)k_v$. Recall the decomposition  \eqref{f7.7} in Lemma \ref{lem7.3}:
\begin{align*}
	W_v'\left(\begin{pmatrix}
		a_v\\
		&1
	\end{pmatrix}k_v\right)=W_v'^{(1)}\left(\begin{pmatrix}
		a_v\\
		&1
	\end{pmatrix}k_v\right)-W_v'^{(2)}\left(\begin{pmatrix}
		a_v\\
		&1
	\end{pmatrix}k_v\right).
\end{align*}

As a consequence, we obtain the corresponding decomposition 
\begin{equation}\label{7.11}
	\Psi_v(s_1,s_2;\ell_v')=\Psi_v^{(1)}(s_1,s_2)-\Psi_v^{(2)}(s_1,s_2),
\end{equation}
where for $j=1,2$, the function $\Psi_v^{(j)}(s_1,s_2)$ is defined by 
\begin{align*}
	\int_{K_v[m_v]}\int_{F_v^{\times}}\overline{\chi_1}(a_v\det k_v)\chi_1\overline{\chi_2}(k_{22})W_v\left(x_v\right)W_v'^{(j)}\left(x_v\right)|a_v|_v^{s_2-1}d^{\times}a_vdk_v,
\end{align*}
where $x_v:=\diag(a_v,1)k_v$. Utilizing the calculation of $W_v'^{(j)}(\diag(a_v,1)k_v)$ in Lemma \ref{lem7.3}, we obtain 
\begin{align*}
	\Psi_v^{(1)}(s_1,s_2)=&\frac{\chi_1(-1)}{\zeta_v(1)\Vol(K_v[m_v])}\sum_{j_v\geq m_v}W_{v}\left(\begin{pmatrix}
		\varpi_v^{j_v}\\
		&1
	\end{pmatrix}\right)q_v^{(s_1-s_2)j_v}\\
	&\frac{\overline{\chi_1}\chi_2(\varpi_v^{j_v-m_v})q_v^{(1-2s_1)m_v}-\chi_1\overline{\chi_2}(\varpi_v)q_v^{(1-2s_1)(j_v+1)}}{1-\chi_1\overline{\chi_2}(\varpi_v)q_v^{-(2s_1-1)}},
\end{align*}
and 
\begin{equation}\label{7.23}
	\Psi_v^{(2)}(s_1,s_2)=\frac{\chi_1(-1)\chi_1\overline{\chi_2}(\varpi_v)}{q_v^{2s_1}\Vol(K_v[m_v])}\sum_{j_v\geq m_v-1}W_{v}\left(\begin{pmatrix}
		\varpi_v^{j_v}\\
		&1
	\end{pmatrix}\right)q_v^{(1-s_1-s_2)j_v}.
\end{equation}

\begin{lemma}\label{lem7.7}
	Let notation be as before. Let $v<\infty$ be such that $m_v\geq 1$ and  $m_v^*=0$. Let   
\begin{align*}
W_v=W_{v,n_v}=\sum_{i=0}^{n_v}\xi_{\sigma_v}(\mathfrak{p}_v^{i},\mathfrak{p}_v^{n_v})q_v^{\frac{i-n_v}{2}}\sigma_v\left(\begin{pmatrix}
1\\
&\varpi_v^{i}
		\end{pmatrix}\right)W_v^{\circ}
\end{align*}
for some $0\leq n_v\leq m_v-r_v.$ Then $\Psi_v^{(1)}(s_1,s_2)$ converge absolutely in the region $\Re(s_2)-|\Re(s_1)|\gg 1$. Moreover, $\Psi_v^{(1)}(s_1,s_2)$ admits a meromorphic continuation to $\mathbb{C}^2$ via the following identity 
\begin{equation}\label{7.22}
\frac{\Psi_v^{(1)}(s_1,s_2)}{L_v(s_2,\sigma_v\times\sigma_v'\otimes\overline{\chi_1})}=\frac{\chi_1(-1)W_v^{\circ}(I_2)}{\zeta_v(1)\Vol(K_v[m_v])q_v^{(s_2+s_1-1/2)m_v}}\cdot I_v^{(1)}(s_1,s_2),
\end{equation}
where 
\begin{align*}
I_v^{(1)}(s_1,s_2):=\sum_{i=0}^{n_v}q_v^{i-\frac{n_v}{2}}&\xi_{\sigma_v}(\mathfrak{p}_v^{i},\mathfrak{p}_v^{n_v})\Big[\lambda_{\sigma}(\mathfrak{p}_v^{m_v-i})+\lambda_{\sigma}(\mathfrak{p}_v^{m_v-2-i})\overline{\chi_1}\chi_2(\varpi_v)q_v^{-2s_2}\\
&-\lambda_{\sigma}(\mathfrak{p}_v^{m_v-i-1})(q_v^{1/2-s_1}+\overline{\chi_1}\chi_2(\varpi_v)q_v^{s_1-1/2})q_v^{-s_2}\Big].
\end{align*}
\end{lemma}
\begin{proof}
	By definition, we can decompose $\Psi_v^{(1)}(s_1,s_2)$ as 
\begin{equation}\label{eq7.12}
\Psi_v^{(1)}(s_1,s_2)=\frac{\chi_1(-1)}{\zeta_v(1)\Vol(K_v[m_v])}\cdot \Big[\Psi_v^{(11)}(s_1,s_2)-\Psi_v^{(12)}(s_1,s_2)\Big],
\end{equation}
where $\Psi_v^{(11)}(s_1,s_2)$ is defined by 
\begin{align*}
\frac{\chi_1\overline{\chi_2}(\varpi_v^{m_v})q_v^{(1-2s_1)m_v}}{1-\chi_1\overline{\chi_2}(\varpi_v)q_v^{-(2s_1-1)}}\sum_{i=0}^{n_v}\xi_{\sigma_v}(\mathfrak{p}_v^{i},\mathfrak{p}_v^{n_v})q_v^{\frac{i-n_v}{2}}\sum_{j_v\geq m_v}W_{v}^{\circ}\left(\begin{pmatrix}
\varpi_v^{j_v}\\
&\varpi_v^i
\end{pmatrix}\right)\frac{\overline{\chi_1}\chi_2(\varpi_v^{j_v})}{q_v^{(s_2-s_1)j_v}};
\end{align*}
and $\Psi_v^{(12)}(s_1,s_2)$ is defined by 
\begin{align*}
\frac{\chi_1\overline{\chi_2}(\varpi_v)q_v^{1-2s_1}}{1-\chi_1\overline{\chi_2}(\varpi_v)q_v^{-(2s_1-1)}}\sum_{i=0}^{n_v}\xi_{\sigma_v}(\mathfrak{p}_v^{i},\mathfrak{p}_v^{n_v})q_v^{\frac{i-n_v}{2}}
\sum_{j_v\geq m_v}W_{v}^{\circ}\left(\begin{pmatrix}
\varpi_v^{j_v}\\
&1
\end{pmatrix}\right)q_v^{(1-s_1-s_2)j_v}.
\end{align*}
	
Hence, $\Psi_v^{(11)}(s_1,s_2)$ and $\Psi_v^{(12)}(s_1,s_2)$ converge absolutely in the region $\Re(s_2+s_1)\gg 1$ and $\Re(s_2-s_1)\gg 1$. 
	
	We proceed to calculate $\Psi_v^{(11)}(s_1,s_2)$ and $\Psi_v^{(12)}(s_1,s_2)$ as follows.
	\begin{itemize}
		\item Through a brute force calculation using the Casselman-Shalika formula and Hecke relations, we obtain 
		\begin{align*}
			&\frac{1}{L(s_2-s_1+1/2,\sigma_v\otimes\overline{\chi_1}\chi_2)}\sum_{j_v\geq m_v}W_{v}^{\circ}\left(\begin{pmatrix}
				\varpi_v^{j_v}\\
				&\varpi_v^i
			\end{pmatrix}\right)\frac{\overline{\chi_1}\chi_2(\varpi_v^{j_v})}{q_v^{(s_2-s_1)j_v}}\\
			=&\frac{q_v^{\frac{i}{2}}W_v^{\circ}(I_2)\big[\lambda_{\sigma}(\mathfrak{p}_v^{m_v-i})-\lambda_{\sigma}(\mathfrak{p}_v^{m_v-1-i})\overline{\chi_1}\chi_2(\varpi_v)q_v^{-(s_2-s_1+1/2)}\big]}{\chi_1\overline{\chi_2}(\varpi_v^{m_v})q_v^{(s_2-s_1+1/2)m_v}},
		\end{align*} 
		where $\lambda_{\sigma}(\mathfrak{p}_v^{-1})=0$. In conjunction with the decomposition 
		\begin{align*}
			L_v(s_2,\sigma_v\times\sigma_v'\otimes\overline{\chi_1})=L_v(s_2+s_1-1/2,\sigma_v)L_v(s_2-s_1+1/2,\sigma_v\otimes \overline{\chi_1}\chi_2),
		\end{align*}
		we obtain that 
\begin{align*}
\frac{\Psi_v^{(11)}(s_1,s_2)}{L_v(s_2,\sigma_v\times\sigma_v'\otimes\overline{\chi_1})}=\frac{W_v^{\circ}(I_2)L_v(2s_1-1,\chi_1\overline{\chi_2})}{q_v^{(s_2+s_1-1/2)m_v}}\sum_{i=0}^{n_v}q_v^{\frac{i}{2}}\widetilde{\xi}_{\sigma_v}(\mathfrak{p}_v^{i},\mathfrak{p}_v^{n_v})J_v^{(11)}(s_1,s_2;i),
\end{align*}
where $\widetilde{\xi}_{\sigma_v}(\mathfrak{p}_v^{i},\mathfrak{p}_v^{n_v})=\xi_{\sigma_v}(\mathfrak{p}_v^{i},\mathfrak{p}_v^{n_v})q_v^{\frac{i-n_v}{2}}$, and 
\begin{align*}
J_v^{(11)}(s_1,s_2;i):=\frac{\lambda_{\sigma}(\mathfrak{p}_v^{m_v-i})-\lambda_{\sigma}(\mathfrak{p}_v^{m_v-1-i})\overline{\chi_1}\chi_2(\varpi_v)q_v^{-(s_2-s_1+1/2)}}{L_v(s_2+s_1-1/2,\sigma_v)}.
\end{align*}
		
\item Similar to the above computation, we obtain  
\begin{align*}
\frac{\Psi_v^{(12)}(s_1,s_2)}{L(s_2+s_1-1/2,\sigma_v)}=\frac{W_v^{\circ}(I_2)L_v(2s_1-1,\chi_1\overline{\chi_2})}{q_v^{(s_2+s_1-1/2)m_v}}\sum_{i=0}^{n_v}q_v^{\frac{i}{2}}\widetilde{\xi}_{\sigma_v}(\mathfrak{p}_v^{i},\mathfrak{p}_v^{n_v})J_v^{(12)}(s_1,s_2;i),
\end{align*}
where 
\begin{align*}
J_v^{(12)}(s_1,s_2;i):=\frac{\chi_1\overline{\chi_2}(\varpi_v)q_v^{1-2s_1}\cdot \big[\lambda_{\sigma}(\mathfrak{p}_v^{m_v-i})-\lambda_{\sigma}(\mathfrak{p}_v^{m_v-1-i})q_v^{-(s_2+s_1-1/2)}\big]}{L_v(s_2-s_1+1/2,\sigma_v\otimes \overline{\chi_1}\chi_2)}.
\end{align*}
\end{itemize}
	
After a detailed calculation, $J_v^{(11)}(s_1,s_2;i)-J_v^{(12)}(s_1,s_2;i)$ is equal to 
\begin{align*}
L_v(2s_1-1,\chi_1\overline{\chi_2})^{-1}\cdot \big[&\lambda_{\sigma}(\mathfrak{p}_v^{m_v-i})+\lambda_{\sigma}(\mathfrak{p}_v^{m_v-2-i})\overline{\chi_1}\chi_2(\varpi_v)q_v^{-2s_2}\\
&-\lambda_{\sigma}(\mathfrak{p}_v^{m_v-i-1})(q_v^{1/2-s_1}+\overline{\chi_1}\chi_2(\varpi_v)q_v^{s_1-1/2})q_v^{-s_2}\big].
\end{align*}
	
Substituting this into \eqref{eq7.12} we then obtain \eqref{7.22}. 
\end{proof}

\begin{lemma}\label{lem7.8}
Let notation be as before. Let $v<\infty$ be such that $m_v\geq 1$ and  $m_v^*=0$. Let   
\begin{align*}
W_v=W_{v,n_v}=\sum_{i=0}^{n_v}\xi_{\sigma_v}(\mathfrak{p}_v^{i},\mathfrak{p}_v^{n_v})q_v^{\frac{i-n_v}{2}}\sigma_v\left(\begin{pmatrix}
1\\
&\varpi_v^{i}
\end{pmatrix}\right)W_v^{\circ}
\end{align*}
for some $0\leq n_v\leq m_v-r_v.$ Then $\Psi_v^{(2)}(s_1,s_2)$ converge absolutely in the region $\Re(s_2)-|\Re(s_1)|\gg 1$. Moreover, $\Psi_v^{(2)}(s_1,s_2)$ admits a meromorphic continuation to $\mathbb{C}^2$ via the following identity 
\begin{align*}
\Psi_v^{(2)}(s_1,s_2)=&
\frac{\chi_1(-1)\chi_1\overline{\chi_2}(\varpi_v)W_v^{\circ}(I_2)L_v(s_2+s_1-1/2,\sigma_v)}{q_v^{2s_1}\Vol(K_v[m_v])}\sum_{i=0}^{n_v}\frac{\xi_{\sigma_v}(\mathfrak{p}_v^{i},\mathfrak{p}_v^{n_v})q_v^{\frac{i-n_v}{2}}}{q_v^{(s_1+s_2-1)i}}\\
&\bigg[\mathbf{1}_{i\geq m_v-1}+\frac{\lambda_{\sigma}(\mathfrak{p}_v^{m_v-1-i})+\lambda_{\sigma}(\mathfrak{p}_v^{m_v-2-i})q_v^{-(s_2+s_1-1/2)}}{q_v^{(m_v-1-i)(s_2+s_1-1/2)}}\mathbf{1}_{0\leq i\leq m_v-2}\bigg].
\end{align*}
\end{lemma}
\begin{proof}
By definition \eqref{7.23}, we have 
\begin{align*}
\Psi_v^{(2)}(s_1,s_2)=
\frac{\chi_1(-1)\chi_1\overline{\chi_2}(\varpi_v)}{q_v^{2s_1}\Vol(K_v[m_v])}\sum_{i=0}^{n_v}\xi_{\sigma_v}(\mathfrak{p}_v^{i},\mathfrak{p}_v^{n_v})q_v^{\frac{i-n_v}{2}}\textbf{S}(1,i),
\end{align*}
where 
\begin{align*}
\textbf{S}(1,i):=\sum_{j_v\geq m_v-1}W_{v}^{\circ}\left(\begin{pmatrix}
\varpi_v^{j_v}\\
&\varpi_v^i
\end{pmatrix}\right)q_v^{(1-s_1-s_2)j_v}.
\end{align*}
As a consequence, we conclude that $\Psi_v^{(2)}(s_1,s_2)$ converge absolutely in the region $\Re(s_2-s_1)\gg 1$. Moreover, we note that $\textbf{S}(1,i)$ is a special case of  \eqref{eq7.18}. Utilizing the calculations \eqref{7.17} and \eqref{7.18} in the proof of Lemma \ref{lem7.6}, $\textbf{S}(1,i)$ is equal to the product of $W_v^{\circ}(I_2)L_v(s_2+s_1-1/2,\sigma_v)q_v^{-(s_1+s_2-1)i}$ with 
\begin{align*}
\mathbf{1}_{i\geq m_v-1}+\frac{\lambda_{\sigma}(\mathfrak{p}_v^{m_v-1-i})+\lambda_{\sigma}(\mathfrak{p}_v^{m_v-2-i})q_v^{-(s_2+s_1-1/2)}}{q_v^{(m_v-1-i)(s_2+s_1-1/2)}}\mathbf{1}_{0\leq i\leq m_v-2}.
\end{align*}
	
Therefore, Lemma \ref{lem7.8} follows. 
\end{proof}

\begin{lemma}\label{lemma7.7}
Let notation be as before. Let $v<\infty$ be such that $\ell_v'\geq 1$. Let   
\begin{align*}
W_v=W_{v,n_v}=\sum_{i=0}^{n_v}\xi_{\sigma_v}(\mathfrak{p}_v^{i},\mathfrak{p}_v^{n_v})q_v^{\frac{i-n_v}{2}}\sigma_v\left(\begin{pmatrix}
1\\
&\varpi_v^{i}
\end{pmatrix}\right)W_v^{\circ}
\end{align*}
for some $0\leq n_v\leq l_v'-r_v.$ Then $\Psi_v(s_1,s_2;\ell_v')$ converges absolutely in the region $\Re(s_2)-|\Re(s_1)|\gg 1$, and it admits a meromorphic continuation to $\mathbb{C}^2$ via the following equality 
\begin{equation}\label{7.14}
\frac{\Psi_v(s_1,s_2;\ell_v')}{L_v(s_2,\sigma_v\times\sigma_v'\otimes\overline{\chi_1})}=\frac{W_v^{\circ}(I_2)W_v'(I_2)}{\chi_1(\varpi_v^{\ell_v'})q_v^{\ell_v's_2}}\cdot \Big[J_v^{(1)}(s_1,s_2;\ell_v')+J_v^{(2)}(s_1,s_2;\ell_v')\Big],
\end{equation}
where 
\begin{align*}
J_v^{(1)}(s_1,s_2;\ell_v'):=&\mathbf{1}_{\substack{r_v=0\\
n_v\leq l_v'}}\cdot \sum_{i=0}^{n_v}\lambda_{\sigma}(\mathfrak{p}_v^i)\xi_{\sigma_v}(\mathfrak{p}_v^{i},\mathfrak{p}_v^{n_v})q_v^{\frac{i-n_v}{2}}(q_v^{i/2}+\ell_v'q_v^{-i/2}|\mathbb{F}_v^{\times}|)^{-1},\\
J_v^{(2)}(s_1,s_2;\ell_v'):=&\sum_{l=n_v+r_v}^{\ell_v'}\frac{\chi_1\chi_2^{-1}(\varpi_v^l)q_v^{2s_2l}}{\Vol(K_v[l])^{-1}}\cdot \sum_{i=0}^{n_v}\frac{\xi_{\sigma_v}(\mathfrak{p}_v^{i},\mathfrak{p}_v^{n_v})\overline{\chi_1}(\varpi_v^{i})H_v^{(i)}(s_2)}{q_v^{i(s_2-1/2)}q_v^{(n_v-i)/2}}.
\end{align*}
Here the function $H_v^{(i)}(s_2)$ is defined by 
\begin{equation}\label{eq7.27}
H_v^{(i)}(s_2):=\lambda_{\sigma'}(\mathfrak{p}_v^{i})-\chi_2(\varpi_v)\lambda_{\sigma}(\mathfrak{p}_v)\lambda_{\sigma'}(\mathfrak{p}_v^{i-1})q_v^{-s_2}+\chi_2^2(\varpi_v)\lambda_{\sigma'}(\mathfrak{p}_v^{i-2})q_v^{-2s_2}.
\end{equation}
\end{lemma}
\begin{proof}
Recall the definition 
\begin{align*}
\Psi_v(s_1,s_2;\ell_v'):=\int_{N(F_v)Z(F_v)\backslash G(F_v)}W_v\left(x_v\right)W_v'(x_v)\overline{h_v(x_v,\overline{s_2};\ell_v')}dx_v.\tag{\ref{equa7.4}}
\end{align*}
	
By the assumption $\ell_v'\geq 1$, we have $m_v=0$. So $\chi_1$ and $\chi_2$ are unramified. By Iwasawa decomposition, we may write  $x_v=\diag(\varpi_v^{j_v},1)k_v$. Then 
\begin{align*}
\overline{h_v(x_v,\overline{s_2};\ell_v')}=\overline{\chi_1}(\varpi_v^{j_v+\ell_v'})q_v^{-s_2(j_v+\ell_v')}\int_{F_v^{\times}}\overline{\Phi_v^*\left((0,t)x_v\boldsymbol{d}_v\right)}\chi_1^{-1}\chi_2(t)|t|_v^{2s_2}d^{\times}t.
\end{align*}
	
According to the definition \eqref{phi*}, $\Phi^*(t_{1,v},t_{2,v})=\mathbf{1}_{\mathcal{O}_v}(t_{1,v})\mathbf{1}_{\mathcal{O}_v}(t_{2,v})$. Similar to Lemma \ref{lem4.1}, the integral $\overline{h_v(x_v,\overline{s_2};\ell_v')}$ is equal to 
\begin{align*}
\overline{\chi_1}(\varpi_v^{j_v+\ell_v'})q_v^{-(j_v+\ell_v')s_2}\Big[L_v(2s_2,\chi_1^{-1}\chi_2)\mathbf{1}_{k_v\in K_v}+\sum_{l=1}^{\ell_v'}\chi_1\chi_2^{-1}(\varpi_v^l)q_v^{2s_2l}\mathbf{1}_{k_v\in K_v[l]}\Big]. 
\end{align*}

Substituting this into \eqref{equa7.4} we obtain 
\begin{equation}\label{7.27}
\Psi_v(s_1,s_2;\ell_v')=\overline{\chi_1}(\varpi_v^{\ell_v'})q_v^{-\ell_v's_2}\Big[\Psi_v^{(1)}(s_1,s_2;\ell_v')+\Psi_v^{(2)}(s_1,s_2;\ell_v')\Big],
\end{equation}
where $\Psi_v^{(1)}(s_1,s_2;\ell_v'):=
	L_v(2s_2,\chi_1^{-1}\chi_2)\mathcal{I}(s_1,s_2;0)$, 
and
\begin{equation}\label{7.28}
\Psi_v^{(2)}(s_1,s_2;\ell_v'):=\sum_{l=1}^{\ell_v'}\chi_1\chi_2^{-1}(\varpi_v^l)q_v^{2s_2l}\cdot \mathcal{I}(s_1,s_2;l).
\end{equation}
Here, for $l\geq 0$ and $\Re(s_2)-\Re(|s_1|)\gg 1$, the function $\mathcal{I}(s_1,s_2;l)$ is defined by 
\begin{equation}\label{7.29}
\sum_{j_v\geq 0}\int_{K_v[l]}W_v\left(\begin{pmatrix}
\varpi_v^{j_v}\\
			&1
\end{pmatrix}k_v\right)dk_vW_v'\left(\begin{pmatrix}
\varpi_v^{j_v}\\
	&1
\end{pmatrix}\right)\overline{\chi_1}(\varpi_v^{j_v})q_v^{j_v(1-s_2)}.
\end{equation}
	
We proceed to calculate $\Psi_v^{(1)}(s_1,s_2;\ell_v')$ and $\Psi_v^{(2)}(s_1,s_2;\ell_v')$ as follows. 
\begin{itemize}
\item Let $\iota_i:=q_v^{i/2}+iq_v^{-i/2}|\mathbb{F}_v^{\times}|$. Note that $K_v[0]=K_v$. By \eqref{5.2} we obtain 
\begin{align*}
\iota_i\cdot\int_{K_v[0]}W_v^{\circ}\left(x_vk_v\begin{pmatrix}
1\\
&\varpi_v^i
\end{pmatrix}\right)dk_v=\lambda_{\sigma}(\mathfrak{p}_v^i)W_v^{\circ}\left(x_v\right)\mathbf{1}_{r_v=0}
\end{align*}
for all $x_v\in G(F_v)$. So  $\mathcal{I}(s_1,s_2;0)$ boils down to 
\begin{align*}
\mathbf{1}_{r_v=0}\sum_{i=0}^{n_v}\frac{\lambda_{\sigma}(\mathfrak{p}_v^i)\xi_{\sigma_v}(\mathfrak{p}_v^{i},\mathfrak{p}_v^{n_v})}{\iota_iq_v^{(n_v-i)/2}}\sum_{j_v\geq 0}W_v^{\circ}\left(\begin{pmatrix}
\varpi_v^{j_v}\\
&1
\end{pmatrix}\right)W_v'\left(\begin{pmatrix}
\varpi_v^{j_v}\\
&1
\end{pmatrix}\right)\frac{\overline{\chi_1}(\varpi_v^{j_v})}{q_v^{j_v(s_2-1)}},
\end{align*}
which, by Casselman-Shalika formula, is equal to 
\begin{align*}
W_v^{\circ}(I_2)W_v'(I_2)\mathbf{1}_{r_v=0}\sum_{i=0}^{n_v}\iota_i^{-1}\lambda_{\sigma}(\mathfrak{p}_v^i)\xi_{\sigma_v}(\mathfrak{p}_v^{i},\mathfrak{p}_v^{n_v})q_v^{\frac{i-n_v}{2}}\cdot \frac{L_v(s_2,\sigma_v\times\sigma_v'\otimes\overline{\chi_1})}{L_v(2s_2,\chi_1^{-1}\chi_2)}.
\end{align*}

As a consequence, we obtain  
\begin{equation}\label{eq7.30}
\frac{\Psi_v^{(1)}(s_1,s_2;\ell_v')}{L_v(s_2,\sigma_v\times\sigma_v'\otimes\overline{\chi_1})}=W_v^{\circ}(I_2)W_v'(I_2)\mathbf{1}_{r_v=0}\cdot \sum_{i=0}^{n_v}\frac{\lambda_{\sigma}(\mathfrak{p}_v^i)\xi_{\sigma_v}(\mathfrak{p}_v^{i},\mathfrak{p}_v^{n_v})q_v^{\frac{i-n_v}{2}}}{(q_v^{i/2}+\ell_v'q_v^{-i/2}|\mathbb{F}_v^{\times}|)},
\end{equation}
which yields a meromorphic continuation of $\Psi_v^{(1)}(s_1,s_2;\ell_v')$. 
		
\item Recall that $W_v=W_{v,n_v}$ for some $0\leq n_v\leq \ell_v'+m_v-r_v=\ell_v'-r_v$. Since $\{W_{v,i}:\ 0\leq i \leq l-r_v\}$ form an orthonormal basis $\mathcal{W}(\sigma_v,\theta_v)^{K_v[l]}$, then 
\begin{equation}\label{7.30}
\int_{K_v[l]}\sigma_v(k_v)W_vdk_v=\Vol(K_v[l])W_v\mathbf{1}_{n_v\leq l-r_v}.
\end{equation} 
		
Substituting \eqref{7.29} and \eqref{7.30} into \eqref{7.28} we obtain 
\begin{align*}
\Psi_v^{(2)}(s_1,s_2;\ell_v')=\sum_{l=n_v+r_v}^{\ell_v'}\chi_1\chi_2^{-1}(\varpi_v^l)q_v^{2s_2l}\Vol(K_v[l])\cdot \sum_{i=0}^{n_v}\xi_{\sigma_v}(\mathfrak{p}_v^{i},\mathfrak{p}_v^{n_v})q_v^{\frac{i-n_v}{2}}\mathcal{J}_i(s_1,s_2),
\end{align*}
where 
\begin{align*}
\mathcal{J}_i(s_1,s_2):=\sum_{j_v\geq 0}W_v^{\circ}\left(\begin{pmatrix}
\varpi_v^{j_v}\\
&\varpi_v^i
\end{pmatrix}\right)W_v'\left(\begin{pmatrix}
\varpi_v^{j_v}\\
&1
\end{pmatrix}\right)\overline{\chi_1}(\varpi_v^{j_v})q_v^{j_v(1-s_2)}.
\end{align*}
		
Similar to Lemma \ref{lem6.2}, we deduce by a direct calculation that 
\begin{align*}
\mathcal{J}_i(s_1,s_2)=\overline{\chi_1}(\varpi_v^{i})q_v^{i(1/2-s_2)}W_v^{\circ}(I_2)W_v'(I_2)H_v^{(i)}(s_2)L_v(s_2,\sigma_v\times\sigma_v'\otimes\overline{\chi_1}),
\end{align*}
where $H_v^{(i)}(s_2)$ is defined by \eqref{eq7.27}.
		
Therefore, $\Psi_v^{(2)}(s_1,s_2;\ell_v')L_v(s_2,\sigma_v\times\sigma_v'\otimes\overline{\chi_1})^{-1}$ is equal to 
\begin{equation}\label{7.32}
W_v^{\circ}(I_2)W_v'(I_2)\sum_{l=n_v+r_v}^{\ell_v'}\frac{\chi_1\chi_2^{-1}(\varpi_v^l)q_v^{2s_2l}}{\Vol(K_v[l])^{-1}}\cdot \sum_{i=0}^{n_v}\frac{\xi_{\sigma_v}(\mathfrak{p}_v^{i},\mathfrak{p}_v^{n_v})\overline{\chi_1}(\varpi_v^{i})H_v^{(i)}(s_2)}{q_v^{i(s_2-1/2)}q_v^{(n_v-i)/2}},
\end{equation}
which gives a meromorphic continuation of $\Psi_v^{(2)}(s_1,s_2;\ell_v')$ to $\mathbb{C}^2$. 
\end{itemize}
	
Therefore, \eqref{7.14} follows from \eqref{7.27}, \eqref{eq7.30} and \eqref{7.32}. 
\end{proof}

\begin{comment}

	\bigskip
	\bigskip
	\bigskip
	
	\begin{align*}
		\frac{\omega_v\overline{\omega_v'}(a_v)}{|a_v|_v^{s_1-1}}\int_{F_v^{\times}}\int_{F_v}\Phi_v((t_v,b_v)
		k_v)\psi_v(a_vt_v^{-1}b_v)db_v\overline{\omega_v}\omega_v'(t_v)|t_v|_v^{2s_1-1}d^{\times}t_v.
	\end{align*}
	
	\begin{align*}
		\frac{\omega_v\overline{\omega_v'}(a_v)}{|a_v|_v^{s_1-1}}\int_{F_v^{\times}}\int_{F_v}\Phi_v((t_v^{-1},b_v)
		k_v)\psi_v(a_vb_vt_v)db_v\omega_v\omega_v'^{-1}(t_v)|t_v|_v^{1-2s_1}d^{\times}t_v.
	\end{align*}
	
	\bigskip
	\bigskip
	\bigskip
	
	\begin{align*}
		|t_v|_v\gg\textbf{C}_v
	\end{align*}
	
	\bigskip
	\bigskip
	\bigskip

	Here $k_v=\begin{pmatrix}
		k_{11}& k_{12}\\
		k_{21}& k_{22}
	\end{pmatrix}$, with $|k_{12}|\asymp |k_{21}|\ll C_v^{-1}$. 
	
	\begin{align*}
		\frac{\omega_v\overline{\omega_v'}(a_v)}{|a_v|_v^{s_1-1}}\int_{F_v^{\times}}\mathcal{F}_2\Phi_v((t_v,a_vt_v^{-1})
		k_v)\overline{\omega_v}\omega_v'(t_v)|t_v|_v^{2s_1-1}d^{\times}t_v.
	\end{align*}
	
	\begin{align*}
		\int_{F_c^{\times}}W_{\varphi,v}\left(\begin{pmatrix}
			a_v&\\
			&1
		\end{pmatrix}k_v\right)\int_{F_v^{\times}}\mathcal{F}_2\Phi_v((t_v,a_vt_v^{-1})
		k_v)\overline{\omega_v}\omega_v'(t_v)|t_v|_v^{2s_1-1}d^{\times}t_v|a_v|_v^{s_2-s_1}d^{\times}a_v.
	\end{align*}
	
	\bigskip
	\bigskip
	\bigskip
\end{comment}

\section{Nonarchimedean Integrals on the Dual Side \RNum{2}}\label{sec9}
In this section we treat the supercuspidal case where the local representation
$\pi_v'$ is supercuspidal and the representation occurring in the $\mathrm{GL}_2$-spectrum
on the dual side is also supercuspidal.  
We compute the associated nonarchimedean local integrals in $\mathcal{M}_{\cusp}^{\du}(s_1,s_2)$
and $\mathcal{M}_{\Eis}^{\du}(s_1,s_2)$, obtaining explicit formulas and uniform bounds adapted
to this situation.

\subsection{Notation}\label{8.1.1}
Let $v$ be a finite place, and $\pi_v'$ be the $v$-th component of $\pi'$. Let $\sigma_v$ be an admissible irreducible generic representation of $\mathrm{PGL}(2)/F$. Suppose that $\sigma_v$ and $\pi_v'$ are  supercuspidal. Let $r_v$ and $r_v'=r_{\pi_v'}$ be the conductor exponent of $\sigma_v$ and $\pi_v'$ (see \textsection\ref{sec4.1}), respectively. Then $r_v\geq 2$ and $r_v'\geq 2$.

%With a slight abuse of notation, we denote by $\omega_v'$ the central character of $\pi_v'$. This notation is used only in this local section and should not be confused with the $\omega_v'$ used in \textsection\ref{sec4} for the central character of $\pi_v'$ as $\pi_v'$ may coincide with $\pi_v'$ under certain character twists. 

Let $W_v^{\circ}$ and $W_v'$ be local newforms of the Whittaker models of $\sigma_v$ and $\pi_v'$, respectively. As in \eqref{eq7.5} in \textsection\ref{7.1.3.}, we set 
\begin{equation}\label{eq8.1}
W_v=W_{v,n_v}=\sum_{i=0}^{n_v}\xi_{\sigma_v}(\mathfrak{p}_v^{i},\mathfrak{p}_v^{n_v})q_v^{\frac{i-n_v}{2}}\sigma_v\left(\begin{pmatrix}
1\\
&\varpi_v^{i}
\end{pmatrix}\right)W_v^{\circ},\ \ 0\leq n_v\leq r_v'-r_v. 
\end{equation}
 
Let notation be as above. We define the local triple integral by   
\begin{equation}\label{8.1}
\mathcal{I}_v^{\sharp}(W_v,W_v',\overline{W_v'}):=\int_{\overline{G'}(F_v)}\frac{\langle\sigma_v(x_v)W_v,W_v\rangle\cdot |\langle\pi_v'(x_v)W_v',W_v'\rangle|^2}{\langle W_v,W_v\rangle\cdot \langle W_v',W_v'\rangle^2}dx_v.
\end{equation}

Since $\pi_v'$ is supercuspidal, the matrix coefficient $\langle \pi_v'(x_v) W_v', W_v' \rangle$ has compact support. Consequently, the integral \eqref{8.1} converges absolutely, ensuring that the local integral $\mathcal{I}_v^{\sharp}(W_v, W_v', \overline{W_v'})$ is well-defined. The main result in this section is an explicit calculation of $\mathcal{I}_v^{\sharp}(W_v, W_v', \overline{W_v'})$. See Proposition \ref{prop9.6} in \textsection\ref{sec9.5}, which will be proved in \textsection\ref{sec9.6}--\textsection\ref{sec9.11}.

\begin{comment}
Moreover, for every unitary character $\chi_v$ of $F_v^{\times}$, we have
\begin{equation}\label{min}
\mathcal{I}_v^{\sharp}(W_v,W_v'\otimes \chi_v,\overline{W_v'\otimes \chi_v})=\mathcal{I}_v^{\sharp}(W_v,W_v',\overline{W_v'}).
\end{equation}
Henceforth, we assume that $\sigma_v$ is \textit{twist-minimal}, i.e., $r_v'\leq r_{\pi_v'\otimes\chi_v}$ for all characters $\chi_v$ of $F_v^{\times}$.	
\end{comment}

\subsection{Evaluation of Whittaker Functions}
\subsubsection{Bessel Functions}
Let $j_{\pi_v'}$ be the Bessel function associated with $\pi_v'$, e.g., see \cite{Cog14}. Then for all $g_v\in G(F_v)$, it is well known that 
\begin{equation}\label{8.2}
W_v'\left(\begin{pmatrix}
a_v\\
&1	
\end{pmatrix}wg_v\right)=\omega_v'(a_v)\int_{F_v^{\times}}j_{\pi_v'}(a_vy_v)W_v'\left(\begin{pmatrix}
y_v\\
&1
\end{pmatrix}g_v\right)d^{\times}y_v.
\end{equation}
This is the $p$-adic version of \eqref{eq5.7}. By functional equation, for $\Re(s)\ll 0$ we have 
\begin{align*}
\gamma_v(s,\pi_v',\psi_v)
=\int_{F_v^{\times}}j_{\pi_v'}(y_v)|y_v|_v^{1/2-s}d^{\times}y_v,
\end{align*}
where $\gamma_v(s,\pi_v',\psi_v)$ is the $\gamma$-factor associated with $\pi_v'$ relative to the unramified character $\psi_v$. 
\begin{comment}
\begin{align*}
W_v'\left(\begin{pmatrix}
a_v\\
&1	
\end{pmatrix}\right)=&\omega_v'(a_v)\int_{F_v^{\times}}j_{\pi_v'}(a_vy_v)W_v'\left(\begin{pmatrix}
y_v\\
&1
\end{pmatrix}w\right)d^{\times}y_v\\
=&\int_{F_v^{\times}}j_{\pi_v'}(a_vy_v)\omega_v'(a_vy_v)\int_{\mathcal{O}_v^{\times}}j_{\pi_v'}(y_vy_v')d^{\times}y_v'd^{\times}y_v\\
=&\gamma_v(s,\sigma_v,\psi_v)\int_{F_v^{\times}}j_{\pi_v'}(a_vy_v)\omega_v'(a_vy_v)\mathbf{1}_{e_v(y_v)=-r_v'}d^{\times}y_v\\
=&\gamma_v(1/2,\sigma_v,\psi_v)\int_{\varpi_v^{e_v(a_v)-r_v'}\mathcal{O}_v^{\times}}j_{\pi_v'}(a_vy_v)\omega_v'(a_vy_v)d^{\times}y_v\\
=&\gamma_v(1/2,\sigma_v,\psi_v)\omega_v'(-1)\gamma_v(1/2,\sigma_v\otimes\omega_v'^{-1},\psi_v)\mathbf{1}_{e_v(a_v)-r_v'=-r_{\pi_v'\otimes\omega_v'^{-1}}}\\
=&\mathbf{1}_{e_v(a_v)=0}.
\end{align*}
\end{comment}

Let $\chi_v$ be a general multiplicative character of $F_v^{\times}$. Upon replacing $W_v'(\cdot)$ with $W_v'(\cdot)\cdot \chi_v(\det(\cdot))$ in \eqref{8.2} we obtain $j_{\pi_v'\otimes\chi_v}(y_v)=\chi_v^{-1}(-y_v)j_{\pi_v'}(y_v)$. Thus,
\begin{comment}
$j_{\pi_v'\otimes\chi_v}(y_v)=\overline{\chi_v}(-y_v)j_{\pi_v'}(y_v)$.

\begin{align*}
\chi_v(-a_v\det g_v)W_v'\left(\begin{pmatrix}
a_v\\
&1	
\end{pmatrix}wg_v\right)=\omega_v'(a_v)\chi_v(a_v)^2\int_{F_v^{\times}}\overline{\chi_v}(-a_vy_v)j_{\pi_v'}(a_vy_v)\chi_v(y_v\det g_v)W_v'\left(\begin{pmatrix}
y_v\\
&1
\end{pmatrix}g_v\right)d^{\times}y_v.
\end{align*}
\end{comment} 
\begin{equation}\label{8.3}
\gamma_v(s,\pi_v'\otimes\chi_v,\psi_v)=\chi_v(-1)\int_{F_v^{\times}}\chi_v^{-1}(y_v)j_{\pi_v'}(y_v)|y_v|_v^{1/2-s}d^{\times}y_v,
\end{equation}
where $\gamma_v(s,\pi_v'\otimes\chi_v,\psi_v)$ is the $\gamma$-factor associated with $\pi_v'\otimes\chi_v$. 

The integral \eqref{8.3} converges in $\vartheta<\Re(s)<1-\vartheta$, where we obtain 
\begin{equation}\label{e8.3}
\gamma_v(s,\pi_v'\otimes\chi_v,\psi_v)=\chi_v(-1)\sum_{n\in\mathbb{Z}}q_v^{n(s-1/2)}\int_{\varpi_v^n\mathcal{O}_v^{\times}}\chi_v^{-1}(y_v)j_{\pi_v'}(y_v)d^{\times}y_v.	
\end{equation}

Since $\pi_v'$ is cuspercuspidal, so is $\pi_v'\otimes\chi_v$.  Consequently, the local $L$-factor of $\pi_v'\otimes\chi_v$ is trivial. Therefore,  
\begin{equation}\label{8.4}
\gamma_v(s,\pi_v'\otimes\chi_v,\psi_v)=\varepsilon(s,\pi_v'\otimes\chi_v,\psi_v)=\varepsilon(1/2,\pi_v'\otimes\chi_v,\psi_v) q_v^{(1/2-s)r_{\pi_v'\otimes\chi_v}},
\end{equation}
where $\varepsilon(s,\pi_v',\psi_v)$ is the $\varepsilon$ factor, and $\varepsilon(1/2,\pi_v',\psi_v)$ is the root number of $\pi_v'$ relative to the unramified character $\psi_v$. In particular, $|\varepsilon(1/2,\pi_v',\psi_v)|=1$. Here we recall that $r_{\pi_v'\otimes\chi_v}$ is the conductor exponent of $\pi_v'\otimes\chi_v$ (see \textsection\ref{sec4.1}). In particular, $r_{\pi_v'\otimes\chi_v}=r_v'=r_{\pi_v'}$ if $\chi_v$ is unramified. 

Combining \eqref{e8.3} and \eqref{8.4} (and changing variable $\chi_v\leftrightarrow\chi_v^{-1}$) we obtain 
\begin{equation}\label{8.5}
\int_{\varpi_v^{n}\mathcal{O}_v^{\times}}\chi_v(y_v)j_{\pi_v'}(y_v)d^{\times}y_v=\chi_v(-1)\varepsilon(1/2,\pi_v'\otimes\chi_v^{-1},\psi_v)\mathbf{1}_{n=-r_{\pi_v'\otimes\chi_v^{-1}}}.
\end{equation}

\subsubsection{Evaluation of Whittaker functions}

\begin{lemma}\label{lem8.1}
Let notation be as before. Then 
\begin{equation}\label{8.6}
W_v'\left(\begin{pmatrix}
a_v\\
&1	
\end{pmatrix}w\right)=\omega_v'(a_v)\varepsilon(1/2,\pi_v',\psi_v)W_v'(I_2)\mathbf{1}_{e_v(a_v)=-r_v'}.
\end{equation}
\end{lemma}
\begin{proof}
Since $W_v'$ is a local newform of the supercuspidal representation $\pi_v'$, then $W_v'\left(\begin{pmatrix}
y_v\\
&1
\end{pmatrix}\right)=W_v'(I_2)\mathbf{1}_{y_v\in \mathcal{O}_v^{\times}}$. Hence, it follows from \eqref{8.2} that
\begin{equation}\label{8.7}
W_v'\left(\begin{pmatrix}
a_v\\
&1	
\end{pmatrix}w\right)=W_v'(I_2)\omega_v'(a_v)\int_{\mathcal{O}_v^{\times}}j_{\pi_v'}(a_vy_v)d^{\times}y_v.
\end{equation}

Substituting \eqref{8.5} into \eqref{8.7} yields 
\begin{align*}
W_v'\left(\begin{pmatrix}
a_v\\
&1	
\end{pmatrix}w\right)=W_v'(I_2)\omega_v'(a_v)\varepsilon(1/2,\pi_v',\psi_v)\mathbf{1}_{e_v(a_v)=-r_v'},
\end{align*} 
which gives \eqref{8.6}.
\end{proof}

\begin{lemma}\label{lem8.2}
Let notation be as before. Let $j, l\in\mathbb{Z}$, and $\beta_v\in \mathcal{O}_v^{\times}$. 
\begin{itemize}
\item If $l':=r_v'-l\leq 0$, then 
\begin{equation}\label{8.8}
W_v'\left(\begin{pmatrix}
\varpi_v^j\\
\varpi_v^l\beta_v & 1
\end{pmatrix}\right)=W_v'(I_2)\mathbf{1}_{j=0}.
\end{equation}

\item If $l':=r_v'-l\geq 1$, then 
\begin{multline*}
W_v'\left(\begin{pmatrix}
\varpi_v^j\\
\varpi_v^l\beta_v & 1
\end{pmatrix}\right)=-\frac{\zeta_v(1)W_v'(I_2)}{q_v}\mathbf{1}_{\substack{j=0\\ l'=1}}+\omega_v'(-1)\varepsilon(1/2,\pi_v',\psi_v)q_v^{-\frac{l'}{2}}\zeta_v(1)W_v'(I_2)\\
\sum_{\substack{\chi_v\in \widehat{\mathcal{O}_v^{\times}/(1+\mathfrak{p}_v^{l'})}\\
r_{\chi_v}=l'\\ 
r_{\pi_v'\otimes\chi_v}=r_v'-j}}\chi_v(-\beta_v)\varepsilon(1/2,\pi_v'\otimes\omega_v'^{-1}\chi_v^{-1},\psi_v)\varepsilon(1/2,\chi_v).
\end{multline*}
\end{itemize}
\end{lemma}

\begin{proof}
Notice that $\begin{pmatrix}
1\\
\varpi_v^l\beta_v& 1
\end{pmatrix}=w\begin{pmatrix}
1& \varpi_v^l\beta_v\\
& 1
\end{pmatrix}w$. By \eqref{8.2},
\begin{align*}
W_v'\left(\begin{pmatrix}
\varpi_v^j\\
\varpi_v^l\beta_v & 1
\end{pmatrix}\right)=\omega_v'(\varpi_v^j)\int_{F_v^{\times}}j_{\pi_v'}(\varpi_v^jy_v)W_v'\left(\begin{pmatrix}
y_v\\
&1
\end{pmatrix}\begin{pmatrix}
1& \varpi_v^l\beta_v\\
& 1
\end{pmatrix}w\right)d^{\times}y_v,
\end{align*}
which further simplifies to 
\begin{equation}\label{8.9}
\omega_v'(\varpi_v^j)\int_{F_v^{\times}}j_{\pi_v'}(\varpi_v^jy_v)W_v'\left(\begin{pmatrix}
y_v\\
&1
\end{pmatrix}w\right){\psi_v(y_v\varpi_v^l\beta_v)}d^{\times}y_v.
\end{equation}

Substituting \eqref{8.6} into \eqref{8.9}, $W_v'\left(\begin{pmatrix}
\varpi_v^j\\
\varpi_v^l\beta_v & 1
\end{pmatrix}\right)$ boils down to  
\begin{equation}\label{8.10}
\varepsilon(1/2,\pi_v',\psi_v)W_v'(I_2)\int_{\mathcal{O}_v^{\times}}\omega_v'(\varpi_v^{j-r_v'}y_v)j_{\pi_v'}(\varpi_v^{j-r_v'}y_v){\psi_v(\varpi_v^{l-r_v'}y_v\beta_v)}d^{\times}y_v.
\end{equation}

We have the following two scenarios according to the sign of $l-r_v'.$
\begin{itemize}
\item Suppose that $l-r_v'\geq 0$. Then \eqref{8.10} becomes 
\begin{equation}\label{8.10.}
\varepsilon(1/2,\pi_v',\psi_v)W_v'(I_2)\int_{\mathcal{O}_v^{\times}}\omega_v'(\varpi_v^{j-r_v'}y_v)j_{\pi_v'}(\varpi_v^{j-r_v'}y_v)d^{\times}y_v,
\end{equation}
which, by \eqref{8.5} with $\chi_v=\mathbf{1}$, is equal to 
\begin{equation}\label{8.12}
\omega_v'(-1)\varepsilon(1/2,\pi_v',\psi_v)\varepsilon(1/2,\pi_v'\otimes\omega_v'^{-1},\psi_v)W_v'(I_2)\mathbf{1}_{j=0}
=W_v'(I_2)\mathbf{1}_{j=0}.
\end{equation}
Here we use the fact that 
\begin{align*}
\varepsilon(1/2,\pi_v',\psi_v)\varepsilon(1/2,\pi_v'\otimes\omega_v'^{-1},\psi_v)=\omega_v'(-1).
\end{align*}

\item Suppose that $l-r_v'<0$, namely, $l'=r_v'-l\geq 1$. Let $\beta_v'\in \mathcal{O}_v^{\times}$. By orthogonality of characters, 
\begin{equation}\label{8.11}
\psi_v(\varpi_v^{-l'}\beta_v')=\zeta_v(1)q_v^{-l'}\sum_{\chi_v\in \widehat{\mathcal{O}_v^{\times}/(1+\mathfrak{p}_v^{l'})}}\overline{\chi_v}(\beta_v')\mathcal{G}_{\chi_v,\psi_v},
\end{equation}
where $\widehat{\mathcal{O}_v^{\times}/(1+\mathfrak{p}_v^{l'})}$ is the Pontryagin dual of $\mathcal{O}_v^{\times}/(1+\mathfrak{p}_v^{l'})$, and $\mathcal{G}_{\chi_v,\psi_v}$ is the Gauss sum defined by 
\begin{align*}
\mathcal{G}_{\chi_v,\psi_v}:=\sum_{\alpha_v\in \mathcal{O}_v^{\times}/(1+\mathfrak{p}_v^{l'})}\chi_v(\alpha_v)\psi_v(\varpi_v^{-l'}\alpha_v).
\end{align*}

As a consequence of \eqref{8.11}, the integral 
\begin{align*}
\zeta_v(1)^{-1}q_v^{l'}\int_{\mathcal{O}_v^{\times}}\omega_v'(\varpi_v^{j-r_v'}y_v)j_{\pi_v'}(\varpi_v^{j-r_v'}y_v){\psi_v(\varpi_v^{l-r_v'}y_v\beta_v)}d^{\times}y_v.
\end{align*}
is equal to 
\begin{equation}\label{8.14}
\sum_{\chi_v\in \widehat{\mathcal{O}_v^{\times}/(1+\mathfrak{p}_v^{l'})}}\overline{\chi_v}(\beta_v)\mathcal{G}_{\chi_v,\psi_v}\int_{\mathcal{O}_v^{\times}}\omega_v'(\varpi_v^{j-r_v'}y_v)\overline{\chi_v}(y_v)j_{\pi_v'}(\varpi_v^{j-r_v'}y_v)d^{\times}y_v.
\end{equation}

Making use of \eqref{8.5}, the expression \eqref{8.14} reduces to
\begin{comment}
\begin{align*}
\sum_{\chi_v\in \widehat{\mathcal{O}_v^{\times}/(1+\mathfrak{p}_v^{l'})}}\overline{\chi_v}(\varpi_v^{r_v'-j}\beta_v)\mathcal{G}_{\chi_v,\psi_v}\omega_v'\chi_v(-1)\varepsilon(1/2,\pi_v'\otimes\omega_v'^{-1}\chi_v,\psi_v)\mathbf{1}_{r_{\pi_v'\otimes\omega_v'^{-1}\chi_v}=r_v'-j}
\end{align*}
\end{comment}
\begin{equation}\label{8.16}
\omega_v'(-1)\sum_{\chi_v}\overline{\chi_v}(-\beta_v)\varepsilon(1/2,\pi_v'\otimes\omega_v'^{-1}\chi_v,\psi_v)\mathcal{G}_{\chi_v,\psi_v}\mathbf{1}_{j=r_v'-r_{\pi_v'\otimes\omega_v'^{-1}\chi_v}},
\end{equation}
where $\chi_v\in \widehat{\mathcal{O}_v^{\times}/(1+\mathfrak{p}_v^{l'})}$. Substituting \eqref{8.11} and \eqref{8.16} into \eqref{8.10} we derive that 
\begin{multline*}
W_v'\left(\begin{pmatrix}
\varpi_v^j\\
\varpi_v^l\beta_v & 1
\end{pmatrix}\right)=\varepsilon(1/2,\pi_v',\psi_v)W_v'(I_2)\omega_v'(-1)\zeta_v(1)q_v^{-l'}\\
\sum_{\chi_v}\overline{\chi_v}(-\beta_v)\varepsilon(1/2,\pi_v'\otimes\omega_v'^{-1}\chi_v,\psi_v)\mathcal{G}_{\chi_v,\psi_v}\mathbf{1}_{j=r_v'-r_{\pi_v'\otimes\omega_v'^{-1}\chi_v}}.
\end{multline*}

By a brute force calculation, noting that $l'\geq 1$, we have  
\begin{align*}
\mathcal{G}_{\chi_v,\psi_v}=\begin{cases}
-1,\ &\text{if $\chi_v=\mathbf{1}_v$ and $l'=1$,}\\
q_v^{l'/2}\varepsilon(1/2,\chi_v^{-1})\chi_v(\varpi_v^{l'})\mathbf{1}_{r_{\chi_v}=l'},\ &\text{if $\chi_v\neq\mathbf{1}_v$, and $r_{\chi_v}=l'$,}\\
0,\ &\text{otherwise.}
\end{cases}
\end{align*} 

Therefore, using $\varepsilon(1/2,\pi_v',\psi_v)\varepsilon(1/2,\pi_v'\otimes\omega_v'^{-1},\psi_v)=\omega_v'(-1)$, we obtain 
\begin{multline*}
W_v'\left(\begin{pmatrix}
\varpi_v^j\\
\varpi_v^l\beta_v & 1
\end{pmatrix}\right)=-\frac{\zeta_v(1)W_v'(I_2)}{q_v}\mathbf{1}_{\substack{j=0\\ l'=1}}+\varepsilon(1/2,\pi_v',\psi_v)W_v'(I_2)\omega_v'(-1)\zeta_v(1)q_v^{-\frac{l'}{2}}\\
\sum_{\substack{\chi_v\in \widehat{\mathcal{O}_v^{\times}/(1+\mathfrak{p}_v^{l'})}\\
r_{\chi_v}=l'\\ 
r_{\pi_v'\otimes\omega_v'^{-1}\chi_v}=r_v'-j}}\chi_v^{-1}(-\beta_v)\varepsilon(1/2,\pi_v'\otimes\omega_v'^{-1}\chi_v,\psi_v)\varepsilon(1/2,\chi_v^{-1}).
\end{multline*}

Notice that $r_{\pi_v'\otimes\omega_v'^{-1}\chi_v}=r_{\pi_v'\otimes\chi_v^{-1}}$. We deduce the second part of Lemma \ref{lem8.2} by the replacement $\chi_v\leftrightarrow\chi_v^{-1}$.
\end{itemize}

As a result, Lemma \ref{lem8.2} follows from the above discussions. 
\end{proof}

By Lemmas \ref{lem8.1} and \ref{lem8.2}, we can evaluate $W_v'(g_v)$ explicitly (in terms of root numbers) for all $g_v\in G(F_v)$.

\subsection{Matrix Coefficients of $\pi_v'$}\label{sec9.3}
Since $W_v'$ is a local new form, then 
\begin{align*}
\langle\pi_v'(x_v)W_v',W_v'\rangle=\overline{W_v'(I_2)}\int_{\mathcal{O}_v^{\times}}W_v'\left(\begin{pmatrix}
\alpha_v\\
&1
\end{pmatrix}x_v\right)d^{\times}\alpha_v.
\end{align*}
\begin{comment}
Write $x_v=\begin{pmatrix}
1& u_v\\
&1
\end{pmatrix}\begin{pmatrix}
\varpi_v^j\\
&1
\end{pmatrix}k_v$ in the Iwasawa coordinates. Then 
\begin{equation}\label{8.17.}
\langle\pi_v'(x_v)W_v',W_v'\rangle=\overline{W_v'(I_2)}\int_{\mathcal{O}_v^{\times}}W_v'\left(\begin{pmatrix}
\varpi_v^j\alpha_v\\
&1
\end{pmatrix}k_v\right){\psi_v}(\alpha_vu_v)d^{\times}\alpha_v.
\end{equation}
\end{comment}

Since $\pi_v'$ is supercuspidal, $\langle\pi_v'(x_v)W_v',W_v'\rangle$ has compact support. The lemma below provides an explicit calculation of $\langle\pi_v'(x_v)W_v',W_v'\rangle$.
\begin{lemma}\label{lem8.3}
Let notation be as before. 
\begin{itemize}
\item Suppose $x_v=\begin{pmatrix}
1& \varpi_v^c\gamma_v\\
&1
\end{pmatrix}\begin{pmatrix}
\varpi_v^j\\
&1
\end{pmatrix}wk_v'$, where $c\in\mathbb{Z}$ and $\gamma_v\in \mathcal{O}_v^{\times}$, $ k_v'\in K_v[r_v']$. Then 
\begin{equation}\label{8.17}
\frac{\langle\pi_v'(x_v)W_v',W_v'\rangle}{|W_v'(I_2)|^2}=\frac{\varepsilon(1/2,\pi_v',\psi_v)\mathbf{1}_{j=-r_v'}}{\omega_v'(\gamma_v)\omega_v'(\varpi_v)^{r_v'}}\int_{\mathcal{O}_v^{\times}}\omega_v'(\alpha_v){\psi_v}(\varpi_v^c\alpha_v)d^{\times}\alpha_v.
\end{equation}
In particular, $\langle\pi_v'(x_v)W_v',W_v'\rangle=0$ unless $c\geq \min\{-r_{\omega_v'},-1\}$, and $j=-r_v'$. 
 
\item Suppose  $x_v=\begin{pmatrix}
1& \varpi_v^c\gamma_v\\
&1
\end{pmatrix}\begin{pmatrix}
\varpi_v^j\\
&1
\end{pmatrix}\begin{pmatrix}
1\\
\varpi_v^l\beta_v& 1
\end{pmatrix}k_v'$, where $c\in\mathbb{Z}$ and $\gamma_v\in \mathcal{O}_v^{\times}$, $l\leq r_v'$, $\beta_v\in \mathcal{O}_v^{\times}$, and $k_v'\in K_v[r_v']$. 
\begin{itemize}
\item If $l=r_v'$, then 
\begin{equation}\label{8.19}
\langle\pi_v'(x_v)W_v',W_v'\rangle=|W_v'(I_2)|^2\cdot \mathbf{1}_{j=0}\cdot \int_{\mathcal{O}_v^{\times}}\psi_v(\varpi_v^c\alpha_v)d^{\times}\alpha_v.
\end{equation}

\item If $l\leq r_v'-1$, then 
\begin{equation}\label{eq9.20}
\langle\pi_v'(x_v)W_v',W_v'\rangle=\mathcal{M}_{\mathrm{trivial}}'(x_v)+\mathcal{M}_{\mathrm{nontrivial}}'(x_v),
\end{equation}
where 
\begin{equation}\label{eq9.21}
\mathcal{M}_{\mathrm{trivial}}'(x_v):=\begin{cases}
-\zeta_v(1)|W_v'(I_2)|^2q_v^{-1}\mathbf{1}_{j=0,\ r_v'-l=1},\ & \text{if $c\geq 0$,}\\
\zeta_v(1)^2|W_v'(I_2)|^2q_v^{-2}\mathbf{1}_{j=0,\ r_v'-l=1},\ & \text{if $c=-1$,}\\
0,\ & \text{if $c\leq -2$,}
\end{cases}
\end{equation}
and 
\begin{multline*}
\mathcal{M}_{\mathrm{nontrivial}}'(x_v):=\omega_v'(-1)\varepsilon(1/2,\pi_v',\psi_v)q_v^{l-r_v'}\zeta_v(1)^2|W_v'(I_2)|^2\mathbf{1}_{c=l-r_v'}\\
\sum_{\substack{\chi_v\in \widehat{\mathcal{O}_v^{\times}/(1+\mathfrak{p}_v^{r_v'-l})}\\
r_{\chi_v}=r_v'-l\\ 
r_{\pi_v'\otimes\chi_v}=r_v'-j}}\chi_v(-\beta_v\gamma_v)\varepsilon(1/2,\pi_v'\otimes\omega_v'^{-1}\chi_v^{-1},\psi_v)\varepsilon(1/2,\chi_v)^2.
\end{multline*}
In particular, we have $\langle\pi_v'(x_v)W_v',W_v'\rangle=0$ unless $c=l-r_v'$ and $\min\{0,2l-r_v'\}\leq j\leq r_v'$.
\end{itemize}
\end{itemize}
\end{lemma}

\begin{comment}
\begin{align*}
r_v'-r_{\widetilde{\sigma}_v'\otimes\chi_v}=j=i\leq r_v'-r_v
\end{align*}
\begin{align*}
r_v'-l\geq 1,\ \ l-i\geq r_v',
\end{align*}

\begin{align*}
l-i<r_v',\ j+i-2l=-r_v,\ \Rightarrow j-l<r_v'-r_v
\end{align*}
\end{comment}
\begin{proof}
When $x_v=\begin{pmatrix}
1& \varpi_v^c\gamma_v\\
&1
\end{pmatrix}\begin{pmatrix}
\varpi_v^j\\
&1
\end{pmatrix}wk_v'$, we have 
\begin{align*}
\langle\pi_v'(x_v)W_v',W_v'\rangle=&\overline{W_v'(I_2)}\int_{\mathcal{O}_v^{\times}}W_v'\left(\begin{pmatrix}
\varpi_v^j\alpha_v\\
&1
\end{pmatrix}w\right){\psi_v}(\varpi_v^c\alpha_v\gamma_v)d^{\times}\alpha_v.
\end{align*}

Substituting Lemma \ref{lem8.1} into the above integral we obtain 
\begin{align*}
\langle\pi_v'(x_v)W_v',W_v'\rangle=&\frac{\varepsilon(1/2,\pi_v',\psi_v)|W_v'(I_2)|^2\mathbf{1}_{j=-r_v'}}{\omega_v(\varpi_v)^{r_v'}}\int_{\mathcal{O}_v^{\times}}\omega_v'(\alpha_v){\psi_v}(\varpi_v^c\alpha_v\gamma_v)d^{\times}\alpha_v.
\end{align*}

As a result, \eqref{8.17} follows from a change of variable $\alpha_v\mapsto \gamma_v^{-1}\alpha_v$. 

Suppose $k_v=\begin{pmatrix}
1\\
\varpi_v^l\beta_v& 1
\end{pmatrix}k_v'$. Then \eqref{8.19} follows from \eqref{8.8}. Now we assume that $l<r_v'$. By definition, we have 
\begin{equation}\label{eq8.19}
\langle\pi_v'(x_v)W_v',W_v'\rangle=\overline{W_v'(I_2)}\int_{\mathcal{O}_v^{\times}}W_v'\left(\begin{pmatrix}
\varpi_v^j\alpha_v\\
\varpi_v^l\beta_v &1
\end{pmatrix}\right){\psi_v}(\varpi_v^c\alpha_v\gamma_v)d^{\times}\alpha_v.
\end{equation}

Notice that $W_v'\left(\begin{pmatrix}
\varpi_v^j\alpha_v\\
\varpi_v^l\beta_v &1
\end{pmatrix}\right)=W_v'\left(\begin{pmatrix}
\varpi_v^j\\
\varpi_v^l\alpha_v^{-1}\beta_v &1
\end{pmatrix}\right)$. Then 
\begin{align*}
\langle\pi_v'(x_v)W_v',W_v'\rangle=&\overline{W_v'(I_2)}\int_{\mathcal{O}_v^{\times}}W_v'\left(\begin{pmatrix}
\varpi_v^j\\
\varpi_v^l\alpha_v^{-1}\beta_v &1
\end{pmatrix}\right){\psi_v}(\varpi_v^c\alpha_v\gamma_v)d^{\times}\alpha_v.
\end{align*}

Utilizing Lemma \ref{lem8.2}, where the contribution from $\chi_v=\mathbf{1}$ is  
\begin{equation}\label{8.20a}
\mathcal{M}_{\mathrm{trivial}}'(x_v)=-\zeta_v(1)|W_v'(I_2)|^2q_v^{-1}\mathbf{1}_{j=0,\ r_v'-l=1}\int_{\mathcal{O}_v^{\times}}\psi_v(\varpi_v^c\alpha_v)d^{\times}\alpha_v.
\end{equation}
Here the Ramanujan integral can be computed explicitly:  
\begin{align*}
\int_{\mathcal{O}_v^{\times}}\psi_v(\varpi_v^c\alpha_v)d^{\times}\alpha_v=\begin{cases}
1,\ & \text{if $c\geq 0$,}\\
-\zeta_v(1)q_v^{-1},\ & \text{if $c=-1$,}\\
0,\ & \text{if $c\leq -2$.}
\end{cases}
\end{align*}

As a consequence, appealing to Lemma \ref{lem8.2} we obtain 
\begin{align*}
\langle\pi_v'(x_v)W_v',W_v'\rangle=&\mathcal{M}_{\mathrm{trivial}}'(x_v)+\omega_v'(-1)\varepsilon(1/2,\pi_v',\psi_v)q_v^{-l'/2}\zeta_v(1)|W_v'(I_2)|^2\\
&\sum_{\substack{\chi_v\in \widehat{\mathcal{O}_v^{\times}/(1+\mathfrak{p}_v^{l'})}\\
r_{\chi_v}=r_v'-l\\ 
r_{\pi_v'\otimes\chi_v}=r_v'-j}}\boldsymbol{\varepsilon}(\chi_v)\cdot \int_{\mathcal{O}_v^{\times}}\chi_v^{-1}(\alpha_v)\psi_v(\varpi_v^c\alpha_v)d^{\times}\alpha_v,
\end{align*}
where $l'=r_v'-l$, and 
\begin{align*}
\boldsymbol{\varepsilon}(\chi_v):=\chi_v(-\beta_v\gamma_v)\varepsilon(1/2,\pi_v'\otimes\omega_v'^{-1}\chi_v^{-1},\psi_v)\varepsilon(1/2,\chi_v).
\end{align*}

Hence, the desired form for $\langle\pi_v'(x_v)W_v',W_v'\rangle$ follows from the calculation 
\begin{align*}
\int_{\mathcal{O}_v^{\times}}\chi_v^{-1}(\alpha_v)\psi_v(\varpi_v^c\alpha_v)d^{\times}\alpha_v=\zeta_v(1)q_v^{-\frac{r_v'-l}{2}}\varepsilon(1/2,\chi_v)\mathbf{1}_{c=l-r_v'}.
\end{align*}

Let $\chi_v\in \widehat{\mathcal{O}_v^{\times}/(1+\mathfrak{p}_v^{r_v'-l})}$. Combining the constraints $r_{\chi_v}=r_v'-l$ and $r_{\pi_v'\otimes\chi_v}=r_v'-j$ with the relation %(e.g., see \cite[Proposition 5]{RY21})
$r_{\pi_v'\otimes\chi_v}\leq \max\{2r_{\chi_v},r_v'\}$, we obtain  
\begin{align*}
r_v'-j\leq \max\{r_v',2r_v'-2l\},
\end{align*}
which implies that $j\geq \min\{0,2l-r_v'\}$. Since $r_{\pi_v'\otimes\chi_v}=r_v'-j$ and $r_{\pi_v'\otimes\chi_v}\geq 0$, then $j\leq r_v'$. Thus, $\min\{0,2l-r_v'\}\leq j\leq r_v'$.
\end{proof}

\subsection{Matrix Coefficients of $\sigma_v$: Old Vectors}
Let $\sigma_v$ and $\pi_v'$ be the representations  introduced in \textsection\ref{8.1.1}. Suppose that $\sigma_v$ is supercuspidal. Recall the Whittaker vector 
\begin{align*}
W_v=W_{v,n_v}=\sum_{i=0}^{n_v}\xi_{\sigma_v}(\mathfrak{p}_v^{i},\mathfrak{p}_v^{n_v})q_v^{\frac{i-n_v}{2}}\sigma_v\left(\begin{pmatrix}
1\\
&\varpi_v^{i}
\end{pmatrix}\right)W_v^{\circ},\tag{\ref{eq8.1}}
\end{align*}
where $n_v\leq r_v'-r_v$. %In this subsection we assume that $r_v\leq \min\{r_v',m_v^*\}$, where $m_v':=r_{\omega_v\omega_v'^{-1}}$ is defined as in  \textsection\ref{sec4.1}. Recall that $\pi_v'$ has central character $\omega_v'$. 

Let $x_v=\begin{pmatrix}
1& u_v\\
&1
\end{pmatrix}\begin{pmatrix}
\varpi_v^j\\
&1
\end{pmatrix}k_v$ be the Iwasawa decomposition. Let 
\begin{align*}
\mathcal{MC}_{i}(x_v):=\langle\sigma_v(x_v)\sigma_v(\diag(1,\varpi_v^i))W_v^{\circ},\sigma_v(\diag(1,\varpi_v^i)W_v^{\circ})\rangle
\end{align*}
be the matrix coefficient. By a change of variable, we obtain 
\begin{equation}\label{8.22}
\mathcal{MC}_{i}(x_v)=\int_{F_v^{\times}}W_v^{\circ}\left(\begin{pmatrix}
t_v\\
&\varpi_v^{-i}
\end{pmatrix}x_v\begin{pmatrix}
1\\
&\varpi_v^i
\end{pmatrix}\right)\overline{W_v^{\circ}\left(\begin{pmatrix}
t_v\\
&1
\end{pmatrix}\right)}d^{\times}t_v.
\end{equation}

\begin{comment}
\begin{align*}
\mathcal{MC}_{i}(x_v):=\int_{F_v^{\times}}W_v^{\circ}\left(\begin{pmatrix}
t_v\\
&1
\end{pmatrix}x_v\begin{pmatrix}
1\\
&\varpi_v^i
\end{pmatrix}\right)\overline{W_v^{\circ}\left(\begin{pmatrix}
t_v\\
&\varpi_v^{i'}
\end{pmatrix}\right)}d^{\times}t_v
\end{align*}
\end{comment}
\begin{lemma}\label{lem9.4}
Let notation be as before. Let $x_v=\begin{pmatrix}
1& \varpi_v^c\gamma_v\\
&1
\end{pmatrix}\begin{pmatrix}
\varpi_v^j\\
&1
\end{pmatrix}wk_v'$, where $c\in\mathbb{Z}$ and $\gamma_v\in \mathcal{O}_v^{\times}$, $ k_v'\in K_v[r_v']$. Then 
\begin{equation}\label{8.21}
\mathcal{MC}_{i}(x_v)=\varepsilon(1/2,\sigma_v,\psi_v)|W_v^{\circ}(I_2)|^2\mathbf{1}_{j+2i=-r_v}\int_{\mathcal{O}_v^{\times}}\psi_v(\varpi_v^{c+i}\alpha_v)d^{\times}\alpha_v.
\end{equation}

\end{lemma}
\begin{proof}
Since $W_v^{\circ}$ is right $K_v[r_v]$-invariant, and $i\leq n_v\leq r_v'-r_v$, then 
\begin{align*}
W_v^{\circ}\left(\begin{pmatrix}
t_v\\
&\varpi_v^{-i}
\end{pmatrix}x_v\begin{pmatrix}
1\\
&\varpi_v^i
\end{pmatrix}\right)=W_v^{\circ}\left(\begin{pmatrix}
t_v\varpi_v^{i}\\
&1
\end{pmatrix}\begin{pmatrix}
\varpi_v^{j+i}& \varpi_v^c\gamma_v\\
&1
\end{pmatrix}w\right).
\end{align*}
Here we also make use of the fact that $\sigma_v$ has trivial central character. Hence, 
\begin{align*}
\mathcal{MC}_{i}(x_v)=\overline{W_v^{\circ}(I_2)}\int_{\mathcal{O}_v^{\times}}W_v^{\circ}\left(\begin{pmatrix}
\varpi_v^{j+2i}\alpha_v\\
&1
\end{pmatrix}w\right)\psi_v(\varpi_v^{c+i}\alpha_v\gamma_v)d^{\times}\alpha_v.
\end{align*}

Then \eqref{8.21} follows from Lemma \ref{lem8.1} and the change of variable $\alpha_v\mapsto \alpha_v\gamma_v^{-1}$.
\end{proof}

\begin{lemma}\label{lem8.5}
Let notation be as before. Let 
\begin{align*}
x_v=\begin{pmatrix}
1& \varpi_v^c\gamma_v\\
&1
\end{pmatrix}\begin{pmatrix}
\varpi_v^j\\
&1
\end{pmatrix}\begin{pmatrix}
1\\
\varpi_v^l\beta_v& 1
\end{pmatrix}k_v',
\end{align*}
where $c\in\mathbb{Z}$ and $\gamma_v\in \mathcal{O}_v^{\times}$, $l\leq r_v'$, $\beta_v\in \mathcal{O}_v^{\times}$, and $k_v'\in K_v[r_v']$.  
\begin{itemize}
\item If $r_v+i\leq l\leq r_v'$, then 
\begin{align*}
\mathcal{MC}_{i}(x_v)=|W_v^{\circ}(I_2)|^2\cdot \mathbf{1}_{j=0}\cdot \int_{\mathcal{O}_v^{\times}}\psi_v(\varpi_v^{c+i}\alpha_v)d^{\times}\alpha_v.
\end{align*}
\item If $l\leq r_v+i-1$, then 
\begin{equation}\label{equ9.25}
\mathcal{MC}_{i}(x_v)=\mathcal{M}_{\mathrm{trivial}}(x_v)+\mathcal{M}_{\mathrm{nontrivial}}(x_v),
\end{equation}
where  
\begin{equation}\label{eq9.27}
\mathcal{M}_{\mathrm{trivial}}(x_v):=\begin{cases}
-\zeta_v(1)|W_v^{\circ}(I_2)|^2q_v^{-1}\mathbf{1}_{j=0,\ r_v-l+i=1},\ & \text{if $c\geq -i$,}\\
\zeta_v(1)^2|W_v^{\circ}(I_2)|^2q_v^{-2}\mathbf{1}_{j=0,\ r_v-l+i=1},\ & \text{if $c=-1-i$,}\\
0,\ & \text{if $c\leq -2-i$,}
\end{cases}
\end{equation}
and 
\begin{multline*}
\mathcal{M}_{\mathrm{nontrivial}}(x_v):=\varepsilon(1/2,\sigma_v,\psi_v)q_v^{l-i-r_v}\zeta_v(1)^2|W_v^{\circ}(I_2)|^2\mathbf{1}_{c=l-r_v-2i}\\
\sum_{\substack{\chi_v\in \widehat{\mathcal{O}_v^{\times}/(1+\mathfrak{p}_v^{r_v-l+i})}\\
r_{\chi_v}=r_v-l+i\\ 
r_{\sigma_v\otimes\chi_v}=r_v-j}}\chi_v(\varpi_v^{r_v-j})\chi_v(-\beta_v\gamma_v)\varepsilon(1/2,\sigma_v\otimes\chi_v^{-1},\psi_v)\varepsilon(1/2,\chi_v)^2.
\end{multline*}
\end{itemize}
\end{lemma}
\begin{proof}
By definition and $\diag(1,\varpi_v^{-i})K_v[r_v']\diag(1,\varpi_v^i)\in K_v[r_v],$ the function 
\begin{align*}
W_v^{\circ}\left(\begin{pmatrix}
t_v\\
&\varpi_v^{-i}
\end{pmatrix}x_v\begin{pmatrix}
1\\
&\varpi_v^i
\end{pmatrix}\right)
\end{align*}
is equal to 
\begin{align*}
W_v^{\circ}\left(\begin{pmatrix}
t_v\\
&\varpi_v^{-i}
\end{pmatrix}\begin{pmatrix}
1& \varpi_v^c\gamma_v\\
&1
\end{pmatrix}\begin{pmatrix}
\varpi_v^j\\
&\varpi_v^i
\end{pmatrix}\begin{pmatrix}
1\\
\varpi_v^{l-i}\beta_v& 1
\end{pmatrix}\right).
\end{align*}
As a consequence, we obtain 
\begin{align*}
W_v^{\circ}\left(\begin{pmatrix}
t_v\\
&\varpi_v^{-i}
\end{pmatrix}x_v\begin{pmatrix}
1\\
&\varpi_v^i
\end{pmatrix}\right)=\psi_v(\varpi_v^{c+i}t_v\gamma_v)W_v^{\circ}\left(\begin{pmatrix}
\varpi_v^{j}t_v\\
&1
\end{pmatrix}\begin{pmatrix}
1\\
\varpi_v^{l-i}\beta_v& 1
\end{pmatrix}\right).
\end{align*}
Substituting this into \eqref{8.22}, we derive that  
\begin{align*}
\mathcal{MC}_{i}(x_v)=\overline{W_v(I_2)}\int_{\mathcal{O}_v^{\times}}W_v^{\circ}\left(\begin{pmatrix}
\varpi_v^{j}\alpha_v\\
\varpi_v^{l-i}\beta_v& 1
\end{pmatrix}\right)\psi_v(\varpi_v^{c+i}\alpha_v\gamma_v)d^{\times}\alpha_v.
\end{align*}

As a consequence, Lemma \ref{lem8.5} follows from a similar proof as Lemma \ref{lem8.3}, with the following substitutions: $\pi_v'\mapsto\sigma_v$, $r_v'\mapsto r_v$, $j\mapsto j$, $l\mapsto l-i$, $c\mapsto c+i$, and $\omega_v'\mapsto\mathbf{1}_v$. 
\end{proof}

\subsection{The Local Triple Integral}\label{sec9.5}
Let notation be as before. Define the local triple integral (relative to $\sigma_v(\diag(1,\varpi_v^i))W_v^{\circ}$, $W_v'$ and $W_v'$) by 
\begin{equation}\label{eq9.28}
\mathcal{I}_{i}:=|W_v'(I_2)|^{-4}|W_v^{\circ}(I_2)|^{-2}\int_{\overline{G'}(F_v)}\mathcal{MC}_{i}(x_v)\cdot |\langle\pi_v'(x_v)W_v',W_v'\rangle|^2dx_v,
\end{equation}
where $\mathcal{MC}_{i}(x_v)$ is defined in \eqref{8.22}. The main result in this section is the following.

\begin{prop}\label{prop9.6}
Let notation be as before. Let $0\leq i\leq r_v'-r_v$. Then 
\begin{equation}
\mathcal{I}_{i}\ll r_v'\cdot\Big[q_v^{-r_v'}+q_v^{-\frac{r_v'}{2}}\cdot \mathbf{1}_{r_v=r_v',\ i=0}\Big],
\end{equation}
where the implied constant is absolute. 	
\end{prop}
\begin{remark}
The expected ``sharp'' bound is
\begin{align*}
\mathcal{I}_{i}\ll r_v'\cdot\Big[q_v^{-r_v'}+q_v^{-r_v'}\cdot \mathbf{1}_{r_v=r_v',\ i=0}\Big],
\end{align*}
where a nontrivial estimate of certain character sums (see $\mathcal{S}(j)$ in Lemma \ref{lem9.14} in \textsection\ref{sec9.9.4}) is required. 
\end{remark}

\subsection{Decomposition of $\mathcal{I}_{i}$}\label{sec9.6}
Making use of the Iwasawa decomposition 
\begin{align*}
x_v=\begin{pmatrix}
1& \varpi_v^c\gamma_v\\
&1
\end{pmatrix}\begin{pmatrix}
\varpi_v^j\\
&1
\end{pmatrix}k_v
\end{align*} 
we may decompose $\mathcal{I}_{i}$ as 
\begin{equation}\label{eq9.25}
\mathcal{I}_{i}=\mathcal{I}_{i}^{^{\mathrm{big}}}+\mathcal{I}_{i}^{^{\mathrm{small}}},
\end{equation}
where 
\begin{align*}
&\mathcal{I}_{i}^{^{\mathrm{big}}}:=|W_v'(I_2)|^{-4}|W_v^{\circ}(I_2)|^{-2}\int_{\overline{G'}(F_v)}\mathcal{MC}_{i}(x_v)\cdot |\langle\pi_v'(x_v)W_v',W_v'\rangle|^2\mathbf{1}_{k_v\in K_v-K_v[1]}dx_v,\\
&\mathcal{I}_{i}^{^{\mathrm{small}}}:=|W_v'(I_2)|^{-4}|W_v^{\circ}(I_2)|^{-2}\int_{\overline{G'}(F_v)}\mathcal{MC}_{i}(x_v)\cdot |\langle\pi_v'(x_v)W_v',W_v'\rangle|^2\mathbf{1}_{k_v\in K_v[1]}dx_v.
\end{align*}

We use the superscript ``$\mathrm{big}$'' to refer the contribution from $k_v\in K_v-K_v[1]$, the big cell in the Bruhat decomposition; and ``$\mathrm{small}$'' refers to the contribution from the small cell $k_v\in K_v[1]$ in the Bruhat decomposition.

\subsection{Calculation of $\mathcal{I}_{i}^{\mathrm{big}}$}

\begin{lemma}\label{lem9.6}
Let notation be as before. Let $0\leq i\leq r_v'-r_v$.  Then 
\begin{align*}
\mathcal{I}_{i}^{\mathrm{big}}=\varepsilon(1/2,\sigma_v,\psi_v)\Vol(K_v[r_v'])\mathbf{1}_{r_v'=r_v+2i}\cdot S(i;\omega_v'),
\end{align*}
where
\begin{align*}
S(i;\omega_v')=\begin{cases}
-\zeta_v(1)^2q_v^{-2}+\zeta_v(1)^{-1},\ &\text{if $r_{\omega_v'}=0$ and $i=0$,}\\
\zeta_v(1)q_v^{-1}+\zeta_v(1)^{-1},\ &\text{if $r_{\omega_v'}=0$ and $i\geq 1$,}\\
\zeta_v(1),\ &\text{if $r_{\omega_v'}\geq 1$ and $i\geq r_{\omega_v'}$,}\\
-\zeta_v(1)^2q_v^{-1},\ &\text{if $r_{\omega_v'}\geq 1$ and $i=r_{\omega_v'}-1$,}\\
0,\ &\text{otherwise.}
\end{cases}
\end{align*}
\end{lemma}
\begin{proof}
By Lemmas \ref{lem8.3}, and \ref{lem9.4}, we have 
\begin{align*}
\mathcal{I}_{i}^{\mathrm{big}}=\frac{\varepsilon(1/2,\sigma_v,\psi_v)\Vol(K_v[r_v'])\mathbf{1}_{r_v'=r_v+2i}}{\zeta_v(1)}\cdot S(i;\omega_v'), 
\end{align*}
where 
\begin{equation}\label{eq9.24}
S(i;\omega_v'):=\sum_{c\in \mathbb{Z}}q_v^{-c}\Big|\int_{\mathcal{O}_v^{\times}}\omega_v'(\beta_v){\psi_v}(\varpi_v^c\beta_v)d^{\times}\beta_v\Big|^2\int_{\mathcal{O}_v^{\times}}\psi_v(\varpi_v^{c+i}\alpha_v)d^{\times}\alpha_v.
\end{equation}

\begin{itemize}
\item Suppose that $r_{\omega_v'}=0$, i.e., $\omega_v'$ is unramified. Then 
\begin{equation}\label{9.24.}
\int_{\mathcal{O}_v^{\times}}\psi_v(\varpi_v^{c+i}\alpha_v)d^{\times}\alpha_v=\begin{cases}
1,\ &\text{if $c\geq -i$,}\\
-\zeta_v(1)q_v^{-1},\ &\text{if $c=-1-i$,}\\
0,\ &\text{if $c\leq -2-i$.}
\end{cases}
\end{equation}
Similarly we have 
\begin{equation}\label{eq9.32}
\Big|\int_{\mathcal{O}_v^{\times}}\omega_v'(\beta_v){\psi_v}(\varpi_v^c\alpha_v)d^{\times}\beta_v\Big|^2=\begin{cases}
1,\ &\text{if $c\geq 0$,}\\
\zeta_v(1)^2q_v^{-2},\ &\text{if $c=-1$,}\\
0,\ &\text{if $c\leq -2$.}
\end{cases}
\end{equation}

As a consequence, we obtain 
\begin{align*}
S(i;\omega_v')=\zeta_v(1)^2q_v^{-1}\int_{\mathcal{O}_v^{\times}}\psi_v(\varpi_v^{i-1}\alpha_v)d^{\times}\alpha_v+\sum_{c\geq 0}q_v^{-c}\int_{\mathcal{O}_v^{\times}}\psi_v(\varpi_v^{c+i}\alpha_v)d^{\times}\alpha_v,
\end{align*}
which simplifies to 
\begin{equation}\label{9.24}
S(i;\omega_v')=-\zeta_v(1)^3q_v^{-2}\mathbf{1}_{i=0}+\zeta_v(1)^2q_v^{-1}\mathbf{1}_{i\geq 1}+1.
\end{equation}

\item Suppose that $r_{\omega_v'}\geq 1$, i.e., $\omega_v'$ is ramified. Then 
\begin{equation}\label{eq9.34}
\Big|\int_{\mathcal{O}_v^{\times}}\omega_v'(\beta_v){\psi_v}(\varpi_v^c\alpha_v)d^{\times}\beta_v\Big|^2=\begin{cases}
\zeta_v(1)^2q_v^{-r_{\omega_v'}},\ &\text{if $c=-r_{\omega_v'}$,}\\
0,\ &\text{if $c\neq -r_{\omega_v'}$.}
\end{cases}
\end{equation}
In conjunction with \eqref{9.24.} we obtain 
\begin{equation}\label{9.25}
S(i;\omega_v')=\zeta_v(1)^2\int_{\mathcal{O}_v^{\times}}\psi_v(\varpi_v^{i-r_{\omega_v'}}\alpha_v)d^{\times}\alpha_v=\begin{cases}
\zeta_v(1)^2,\ &\text{if $i\geq r_{\omega_v'}$,}\\
-\frac{\zeta_v(1)^3}{q_v},\ &\text{if $i=r_{\omega_v'}-1$,}\\
0,\ &\text{if $i\leq r_{\omega_v'}-2$.}
\end{cases}
\end{equation}
\end{itemize}

Therefore, Lemma \ref{lem9.6} follows from \eqref{9.24} and \eqref{9.25}. 
\end{proof}

\subsection{Calculation of $\mathcal{I}_{i}^{\mathrm{small}}$: $r_v+i\leq l\leq r_v'$}
By Cramer's rule, we may write $k_v\in K_v[1]$ uniquely as $k_v=\begin{pmatrix}
1\\
\alpha_v&1
\end{pmatrix}K_v[r_v']$ for some $\alpha_v\in \mathfrak{p}_v/\mathfrak{p}_v^{r_v'}$.  Therefore, 
\begin{align*}
\mathcal{I}_{i}^{\mathrm{small}}=\sum_{l=1}^{r_{v}'-1}\mathcal{I}_{i,l}^{\mathrm{small}}+\mathcal{I}_{i,r_{v}'}^{\mathrm{small}},
\end{align*}
where for $1\leq l\leq r_v'-1$, 
\begin{align*}
\mathcal{I}_{i,l}^{\mathrm{small}}:=|W_v'(I_2)|^{-4}|W_v^{\circ}(I_2)|^{-2}\int_{\overline{G'}(F_v)}\mathcal{MC}_{i}(x_v)\cdot |\langle\pi_v'(x_v)W_v',W_v'\rangle|^2\mathbf{1}_{e_v(\alpha_v)=l}dx_v,
\end{align*}
and for $l=r_v'$, we define 
\begin{align*}
\mathcal{I}_{i,r_v'}^{\mathrm{small}}:=|W_v'(I_2)|^{-4}|W_v^{\circ}(I_2)|^{-2}\int_{\overline{G'}(F_v)}\mathcal{MC}_{i}(x_v)\cdot |\langle\pi_v'(x_v)W_v',W_v'\rangle|^2\mathbf{1}_{k_v\in K_v[r_v']}dx_v.
\end{align*}

\begin{lemma}\label{lemma9.8}
Let notation be as before. Then 
\begin{equation}\label{9.28}
\mathcal{I}_{i,r_v'}^{\mathrm{small}}=\Vol(K_v[r_v'])\cdot \big[-\zeta_v(1)^3q_v^{-2}\mathbf{1}_{i=0}+\zeta_v(1)^2q_v^{-1}\mathbf{1}_{i\geq 1}+1\big].
\end{equation}
\end{lemma}
\begin{proof}
By Lemmas \ref{lem8.3} and Lemma \ref{lem8.5}, we have
\begin{align*}
\mathcal{I}_{i,r_v'}^{\mathrm{small}}=\Vol(K_v[r_v'])S(i;\mathbf{1}),
\end{align*}
where $S(i;\mathbf{1})$ is the sum $S(i;\omega_v')$ defined as in \eqref{eq9.24} with $\omega_v'=\mathbf{1}$. Hence, \eqref{9.28} follows from \eqref{9.24}.
\end{proof}

\begin{lemma}\label{lem9.8}
Let notation be as before. Let $r_v+i\leq l\leq r_v'-1$. Then 
\begin{align*}
\mathcal{I}_{i,l}^{\mathrm{small}}=&\big[\zeta_v(1)^2q_v^{-2}\mathbf{1}_{i\geq 1}-\zeta_v(1)^3q_v^{-3}\mathbf{1}_{i=0}+\zeta_v(1)^{-1}q_v^{-1}\big]\cdot \Vol(K_v[r_v'])\cdot \mathbf{1}_{r_v'-l=1}\\
&+\frac{\Vol(K_v[r_v'])}{\zeta_v(1)^{2}}\Big[\mathbf{1}_{l\geq r_v'-i}-\zeta_v(1)q_v^{-1}\mathbf{1}_{l=r_v'-i-1}\Big]\cdot \sum_{\substack{\chi_v\in \widehat{\mathcal{O}_v^{\times}/(1+\mathfrak{p}_v^{r_v'-l})}\\
r_{\chi_v}=r_v'-l, \ r_{\pi_v'\otimes\chi_v}=r_v'}}1.
\end{align*}
In particular, we have
\begin{equation}\label{eq9.30}
\mathcal{I}_{i,l}^{\mathrm{small}}\ll q_v^{-r_v-i},
\end{equation}
where the implied constant is absolute. 
\end{lemma}
\begin{proof}
By definition and Iwasawa decomposition
\begin{align*}
x_v=\begin{pmatrix}
1& \varpi_v^c\gamma_v\\
&1
\end{pmatrix}\begin{pmatrix}
\varpi_v^j\\
&1
\end{pmatrix}\begin{pmatrix}
1\\
\varpi_v^l\beta_v& 1
\end{pmatrix}k_v',
\end{align*}
the integral $\mathcal{I}_{i,l}^{\mathrm{small}}$ is equal to 
\begin{align*}
\Vol(K_v[r_v'])\sum_jq_v^j\sum_cq_v^{-c}\int_{\mathcal{O}_v^{\times}}\sum_{\beta_v\in \varpi_v^l\mathcal{O}_v^{\times}/(1+\mathfrak{p}_v^{r_v'})}\frac{\mathcal{MC}_{i}(x_v)\cdot |\langle\pi_v'(x_v)W_v',W_v'\rangle|^2}{|W_v^{\circ}(I_2)|^{2}\cdot |W_v'(I_2)|^{4}}d\gamma_v.
\end{align*}

Taking advantage of Lemmas \ref{lem8.3} and Lemma \ref{lem8.5}, we obtain 
\begin{align*}
\frac{\mathcal{I}_{i,l}^{\mathrm{small}}}{\Vol(K_v[r_v'])}=\sum_cq_v^{-c}\int_{\mathcal{O}_v^{\times}}\sum_{\beta_v}\int_{\mathcal{O}_v^{\times}}\psi_v(\varpi_v^{c+i}\alpha_v)d^{\times}\alpha_v \Big|\mathcal{M}_{\mathrm{trivial}}'^{\sharp}+\mathcal{M}_{\mathrm{nontrivial}}'^{\sharp}\Big|^2d\gamma_v,
\end{align*}
where $\beta_v\in \varpi_v^l\mathcal{O}_v^{\times}/(1+\mathfrak{p}_v^{r_v'})$, $\mathcal{M}_{\mathrm{trivial}}'^{\sharp}:=\mathcal{M}_{\mathrm{trivial}}'(x_v)\cdot |W_v^{\circ}(I_2)|^{-2}$, and $\mathcal{M}_{\mathrm{nontrivial}}'^{\sharp}:=\mathcal{M}_{\mathrm{trivial}}'(x_v)\cdot |W_v'(I_2)|^{-2}$. Here $\mathcal{M}_{\mathrm{trivial}}'(x_v)$ and $\mathcal{M}_{\mathrm{nontrivial}}'(x_v)$ are defined in \eqref{eq9.20} (see Lemma \ref{lem8.3} in \textsection\ref{sec9.3}).

Recall that $\mathcal{M}_{\mathrm{trivial}}'(x_v)$ is independent of $\gamma_v$ and  
\begin{align*}
\mathcal{M}_{\mathrm{nontrivial}}'(x_v)=&\omega_v'(-1)\varepsilon(1/2,\pi_v',\psi_v)q_v^{l-r_v'}\zeta_v(1)^2|W_v'(I_2)|^2\mathbf{1}_{c=l-r_v'}\\
&\sum_{\substack{\chi_v\in \widehat{\mathcal{O}_v^{\times}/(1+\mathfrak{p}_v^{r_v'-l})}\\
r_{\chi_v}=r_v'-l\\ 
r_{\pi_v'\otimes\chi_v}=r_v'}}\chi_v(-\beta_v\gamma_v)\varepsilon(1/2,\pi_v'\otimes\omega_v'^{-1}\chi_v^{-1},\psi_v)\varepsilon(1/2,\chi_v)^2.
\end{align*}

Hence, after swapping sums, the integral 
\begin{align*}
\int_{\mathcal{O}_v^{\times}}\sum_{\beta_v}\int_{\mathcal{O}_v^{\times}}\psi_v(\varpi_v^{c+i}\alpha_v)d^{\times}\alpha_v \cdot \mathcal{M}_{\mathrm{trivial}}'^{\sharp}\cdot \overline{\mathcal{M}_{\mathrm{nontrivial}}'^{\sharp}}d\gamma_v
\end{align*}
factors through 
\begin{equation}\label{eq9.33}
\int_{\mathcal{O}_v^{\times}}\chi_v(\gamma_v)d\gamma_v=\zeta_v(1)^{-1}\int_{\mathcal{O}_v^{\times}}\chi_v(\gamma_v)d^{\times}\gamma_v=0
\end{equation}
as $r_{\chi_v}=r_v'-l\geq 1$, i.e., $\chi_v$ is nontrivial. Therefore, we obtain 
\begin{align*}
\mathcal{I}_{i,l}^{\mathrm{small}}=\mathcal{I}_{i,l}^{\mathrm{small}}(1)+\mathcal{I}_{i,l}^{\mathrm{small}}(2),
\end{align*}
where 
\begin{align*}
\mathcal{I}_{i,l}^{\mathrm{small}}(1):=&\Vol(K_v[r_v'])\sum_cq_v^{-c}\int_{\mathcal{O}_v^{\times}}\sum_{\beta_v}\int_{\mathcal{O}_v^{\times}}\psi_v(\varpi_v^{c+i}\alpha_v)d^{\times}\alpha_v \Big|\mathcal{M}_{\mathrm{trivial}}'^{\sharp}\Big|^2d\gamma_v,\\
\mathcal{I}_{i,l}^{\mathrm{small}}(2):=&\Vol(K_v[r_v'])\sum_cq_v^{-c}\int_{\mathcal{O}_v^{\times}}\sum_{\beta_v}\int_{\mathcal{O}_v^{\times}}\psi_v(\varpi_v^{c+i}\alpha_v)d^{\times}\alpha_v \Big|\mathcal{M}_{\mathrm{nontrivial}}'^{\sharp}\Big|^2d\gamma_v.
\end{align*}

\begin{itemize}
\item Combining the calculation of the Ramanujan sum 
\begin{align*}
\int_{\mathcal{O}_v^{\times}}\psi_v(\varpi_v^{c+i}\alpha_v)d^{\times}\alpha_v=\begin{cases}
1,\ &\text{if $c\geq -i$,}\\
-\zeta_v(1)q_v^{-1},\ &\text{if $c=-1-i$,}\\
0,\ &\text{if $c\leq -2-i$,}
\end{cases}\tag{\ref{9.24.}}
\end{align*}
and the formula 
\begin{align*}
\mathcal{M}_{\mathrm{trivial}}'(x_v)=\begin{cases}
-\zeta_v(1)|W_v'(I_2)|^2q_v^{-1}\mathbf{1}_{j=0,\ r_v'-l=1},\ & \text{if $c\geq 0$,}\\
\zeta_v(1)^2|W_v'(I_2)|^2q_v^{-2}\mathbf{1}_{j=0,\ r_v'-l=1},\ & \text{if $c=-1$,}\\
0,\ & \text{if $c\leq -2$}
\end{cases}
\end{align*}
into the integral $\mathcal{I}_{i,l}^{\mathrm{small}}(1)$, we obtain 
\begin{align*}
\frac{\mathcal{I}_{i,l}^{\mathrm{small}}(1)}{\Vol(K_v[r_v'])}=&\Big[\zeta_v(1)^2q_v^{r_v'-l-3}\int_{\mathcal{O}_v^{\times}}\psi_v(\varpi_v^{i-1}\alpha_v)d^{\times}\alpha_v+\zeta_v(1)^{-1}q_v^{r_v'-l-2}\big]\cdot \mathbf{1}_{r_v'-l=1}.
\end{align*}

Applying \eqref{9.24.} one more time we obtain 
\begin{equation}\label{9.30}
\frac{\mathcal{I}_{i,l}^{\mathrm{small}}(1)}{\Vol(K_v[r_v'])}=\big[\zeta_v(1)^2q_v^{-2}\mathbf{1}_{i\geq 1}-\zeta_v(1)^3q_v^{-3}\mathbf{1}_{i=0}+\zeta_v(1)^{-1}q_v^{-1}\big]\cdot \mathbf{1}_{r_v'-l=1}.
\end{equation}

\item Substituting the expression of $\mathcal{M}_{\mathrm{nontrivial}}'(x_v)$ into $\mathcal{I}_{i,l}^{\mathrm{small}}(2)$, along with the fact that $|\varepsilon(1/2,\pi_v',\psi_v)|=1$, we obtain 
\begin{equation}\label{9.31}
\mathcal{I}_{i,l}^{\mathrm{small}}(2)=\zeta_v(1)^{-1}\Vol(K_v[r_v'])q_v^{l-r_v'}\sum_{\beta_v}\int_{\mathcal{O}_v^{\times}}\psi_v(\varpi_v^{l-r_v'+i}\alpha_v)d^{\times}\alpha_v\cdot \mathcal{I}(\beta_v),
\end{equation}
where the auxiliary integral $\mathcal{I}(\beta_v)$ is defined by 
\begin{equation}\label{eq9.38}
\int_{\mathcal{O}_v^{\times}}\Big|\sum_{\substack{\chi_v\in \widehat{\mathcal{O}_v^{\times}/(1+\mathfrak{p}_v^{r_v'-l})}\\
r_{\chi_v}=r_v'-l,\ r_{\pi_v'\otimes\chi_v}=r_v'}}\frac{\varepsilon(1/2,\pi_v'\otimes\omega_v'^{-1}\chi_v^{-1},\psi_v)\varepsilon(1/2,\chi_v)^2}{\chi_v^{-1}(-\beta_v\gamma_v)}\Big|^2d^{\times}\gamma_v.
\end{equation}

Open the square and swap sums, utilizing the orthogonality and the fact that $|\varepsilon(1/2,\pi_v'\otimes\omega_v'^{-1}\chi_v^{-1},\psi_v)|=|\varepsilon(1/2,\chi_v)|=1$, leading to 
\begin{equation}\label{9.32}
\mathcal{I}(\beta_v)=\sum_{\substack{\chi_v\in \widehat{\mathcal{O}_v^{\times}/(1+\mathfrak{p}_v^{r_v'-l})}}}\mathbf{1}_{r_{\chi_v}=r_v'-l, \ r_{\pi_v'\otimes\chi_v}=r_v'}.
\end{equation}

Substituting \eqref{9.32} into \eqref{9.31} we then deduce that 
\begin{equation}\label{9.33}
\frac{\mathcal{I}_{i,l}^{\mathrm{small}}(2)}{\zeta_v(1)^{-2}\Vol(K_v[r_v'])}=\Big[\mathbf{1}_{l\geq r_v'-i}-\zeta_v(1)q_v^{-1}\mathbf{1}_{l=r_v'-i-1}\Big]\cdot \sum_{\substack{\chi_v\in \widehat{\mathcal{O}_v^{\times}/(1+\mathfrak{p}_v^{r_v'-l})}\\
r_{\chi_v}=r_v'-l, \ r_{\pi_v'\otimes\chi_v}=r_v'}}1.
\end{equation}
\end{itemize}

Therefore, the formula for $\mathcal{I}_{i,l}^{\mathrm{small}}$ follows from \eqref{9.30} and \eqref{9.33}. Moreover, the estimate \eqref{eq9.30} follows from the fact that $\Vol(K_v[r_v'])=(q_v-1)/(q_v^{r_v'+1}-1)$ and the trivial bound for the sum over the characters $\chi_v$'s. 
\end{proof}

\subsection{Calculation of $\mathcal{I}_{i}^{\mathrm{small}}$: $1\leq l\leq r_v+i-1$}
By definition, for $1\leq l\leq r_v+i-1$ we have the series expansion:
\begin{align*}
\frac{\mathcal{I}_{i,l}^{\mathrm{small}}}{\Vol(K_v[r_v'])}=\sum_jq_v^j\sum_cq_v^{-c}\int_{\mathcal{O}_v^{\times}}&\sum_{\beta_v\in \varpi_v^l\mathcal{O}_v^{\times}/(1+\mathfrak{p}_v^{r_v'})}\frac{\mathcal{M}_{\mathrm{trivial}}(x_v)+\mathcal{M}_{\mathrm{nontrivial}}(x_v)}{|W_v^{\circ}(I_2)|^2}\\
&\Big|\frac{\mathcal{M}_{\mathrm{trivial}}'(x_v)+\mathcal{M}_{\mathrm{nontrivial}}'(x_v)}{|W_v'(I_2)|^2}\Big|^2d\gamma_v.
\end{align*}
Here $\mathcal{M}_{\mathrm{trivial}}'(x_v)$ and $\mathcal{M}_{\mathrm{nontrivial}}'(x_v)$ are defined in \eqref{eq9.20} in Lemma \ref{lem8.3}; $\mathcal{M}_{\mathrm{trivial}}(x_v)$ and $\mathcal{M}_{\mathrm{nontrivial}}(x_v)$ are defined in \eqref{equ9.25} in Lemma \ref{lem8.5}. Thus, 
\begin{equation}\label{9.37}
\mathcal{I}_{i,l}^{\mathrm{small}}=\sum_{\textbf{j}=(j_1,j_2,j_3)\in \{\text{trivial}, \text{nontrivial}\}^{\otimes 3}}\mathcal{I}_{i,l}^{\mathrm{small}}(\textbf{j}),
\end{equation}
where for the subscripts $j_1, j_2, j_3\in \{\text{trivial}, \text{nontrivial}\}$, we define 
\begin{equation}\label{9.43}
\mathcal{I}_{i,l}^{\mathrm{small}}(\textbf{j})=\Vol(K_v[r_v'])\sum_jq_v^j\sum_cq_v^{-c}\int_{}\ \sum_{\beta_v}\frac{\mathcal{M}_{j_1}(x_v)\mathcal{M}_{j_2}'(x_v)\overline{\mathcal{M}_{j_3}'(x_v)}}{|W_v^{\circ}(I_2)|^2|W_v'(I_2)|^4}d\gamma_v.
\end{equation}
Here $\beta_v\in \varpi_v^l\mathcal{O}_v^{\times}/(1+\mathfrak{p}_v^{r_v'})$, and the integral relative to $\gamma_v$ is over $\mathcal{O}_v^{\times}$.

\subsubsection{Vanishing integrals}
\begin{lemma}\label{lem9.9}
Let notation be as before. Let $1\leq l\leq r_v+i-1$. Then 
\begin{equation}\label{9.38}
\mathcal{I}_{i,l}^{\mathrm{small}}(\textbf{j})=0
\end{equation}
if $\textbf{j}=(j_1,j_2,j_3)$ satisfies $\{j_1,j_2,j_3\}=\{\text{trivial},\text{trivial}, \text{nontrivial}\}$ as a set.
\end{lemma}
\begin{proof}
Suppose that $\textbf{j}=(j_1,j_2,j_3)$ with $\{j_1,j_2,j_3\}=\{\text{trivial},\text{trivial}, \text{nontrivial}\}$. Then \eqref{9.38} holds as a consequence of orthogonality of characters, as $\mathcal{I}_{i,l}^{\mathrm{small}}(\textbf{j})$ factors through \eqref{eq9.33}.
\end{proof}

\subsubsection{Triple trivial case}
\begin{lemma}\label{lemma9.11}
Let notation be as before. Let $1\leq l\leq r_v+i-1$. Let $\textbf{j}=(\text{trivial},\text{trivial}, \text{trivial})$. Then 
\begin{equation}\label{9.42}
\mathcal{I}_{i,l}^{\mathrm{small}}(\textbf{j})=\Bigg[(\zeta_v(1)^4q_v^{-4}-q_v^{-2})\mathbf{1}_{\substack{r_v'-l=1\\
r_v-l=1\\
i=0}}-\zeta_v(1)^3q_v^{-3}\mathbf{1}_{\substack{r_v'-l=1\\
r_v-l+i=1\\ i\geq 1}}\Bigg]\cdot \Vol(K_v[r_v']).
\end{equation}
\end{lemma}
\begin{proof}
ecall the definitions in Lemma \ref{lem8.3} and Lemma \ref{lem8.5}: 
\begin{align*}
\mathcal{M}_{\mathrm{trivial}}'(x_v)=\begin{cases}
-\zeta_v(1)|W_v'(I_2)|^2q_v^{-1}\mathbf{1}_{j=0,\ r_v'-l=1},\ & \text{if $c\geq 0$,}\\
\zeta_v(1)^2|W_v'(I_2)|^2q_v^{-2}\mathbf{1}_{j=0,\ r_v'-l=1},\ & \text{if $c=-1$,}\\
0,\ & \text{if $c\leq -2$,}
\end{cases}\tag{\ref{eq9.21}}
\end{align*}
and 
\begin{align*}
\mathcal{M}_{\mathrm{trivial}}(x_v)=\begin{cases}
-\zeta_v(1)|W_v^{\circ}(I_2)|^2q_v^{-1}\mathbf{1}_{j=0,\ r_v-l+i=1},\ & \text{if $c\geq -i$,}\\
\zeta_v(1)^2|W_v^{\circ}(I_2)|^2q_v^{-2}\mathbf{1}_{j=0,\ r_v-l+i=1},\ & \text{if $c=-1-i$,}\\
0,\ & \text{if $c\leq -2-i$.}\tag{\ref{eq9.27}}
\end{cases}
\end{align*}

Substituting \eqref{eq9.21} and \eqref{eq9.27} into \eqref{9.43}, then the formula 
\eqref{9.42} follows from a direct simplification.  
\end{proof}

\subsubsection{Mixed trivial case}

\begin{lemma}\label{lem9.11}
Let notation be as before. Let $1\leq l\leq r_v+i-1$. Let $\textbf{j}=(\text{trivial},\text{nontrivial}, \text{nontrivial})$. Then 
\begin{align*}
\mathcal{I}_{i,l}^{\mathrm{small}}(\textbf{j})=\frac{\zeta_v(1)^5\Vol(K_v[r_v'])}{q_v^{r_v'-l+1}}\Big[\frac{\zeta_v(1)}{q_v}\cdot \mathbf{1}_{\substack{l=r_v+i-1\\
l=r_v'-i-1}}-\mathbf{1}_{\substack{l=r_v+i-1\\ l\geq r_v'-i}}\Big]\sum_{\substack{\chi_v\in \widehat{\mathcal{O}_v^{\times}/(1+\mathfrak{p}_v^{r_v'-l})}\\
r_{\chi_v}=r_v'-l, \ r_{\pi_v'\otimes\chi_v}=r_v'}}1.
\end{align*}
In particular, we have
\begin{equation}\label{eq9.45}
\mathcal{I}_{i,l}^{\mathrm{small}}(\textbf{j})\ll q_v^{-r_v'-1}.
\end{equation}
\end{lemma}
\begin{proof}
By \eqref{eq9.27} (which implies that $j=0$) and \eqref{9.43} we have 
\begin{align*}
\mathcal{I}_{i,l}^{\mathrm{small}}(\textbf{j})=\Vol(K_v[r_v'])\sum_cq_v^{-c}\frac{\mathcal{M}_{\text{trivial}}(x_v)}{|W_v^{\circ}(I_2)|^2}\int_{\mathcal{O}_v^{\times}}\sum_{\beta_v}\frac{|\mathcal{M}_{\text{nontrivial}}'(x_v)|^2}{|W_v'(I_2)|^4}d\gamma_v.
\end{align*}

Notice that 
\begin{align*}
\int_{\mathcal{O}_v^{\times}}\sum_{\beta_v}\frac{|\mathcal{M}_{\text{nontrivial}}'(x_v)|^2}{|W_v'(I_2)|^4}d\gamma_v=\zeta_v(1)^4q_v^{2(l-r_v')}\mathcal{I}(\beta_v)\mathbf{1}_{c=l-r_v'},
\end{align*}
where $\mathcal{I}(\beta_v)$ is defined as in \eqref{eq9.38}. Hence, by \eqref{9.32}, we obtain 
\begin{equation}\label{9.45}
\mathcal{I}_{i,l}^{\mathrm{small}}(\textbf{j})=\frac{\zeta_v(1)^4\Vol(K_v[r_v'])}{q_v^{2(r_v'-l)}}\sum_{\substack{\chi_v\in \widehat{\mathcal{O}_v^{\times}/(1+\mathfrak{p}_v^{r_v'-l})}\\
r_{\chi_v}=r_v'-l, \ r_{\pi_v'\otimes\chi_v}=r_v'}}\sum_{c\in\mathbb{Z}}\frac{\mathcal{M}_{\text{trivial}}(x_v)\mathbf{1}_{c=l-r_v'}}{q_v^c\cdot |W_v^{\circ}(I_2)|^2}.
\end{equation}

By \eqref{eq9.27} we have
\begin{equation}\label{9.46}
q_v^{l-r_v'}\sum_{c\in\mathbb{Z}}\frac{\mathcal{M}_{\text{trivial}}(x_v)\mathbf{1}_{c=l-r_v'}}{q_v^c\cdot |W_v^{\circ}(I_2)|^2}=\zeta_v(1)^2q_v^{-2}\mathbf{1}_{\substack{l=r_v+i-1\\
l=r_v'-i-1}}-\zeta_v(1)q_v^{-1}\mathbf{1}_{\substack{l=r_v+i-1\\ l\geq r_v'-i}}.
\end{equation}

Therefore, Lemma \ref{lem9.11} follows from \eqref{9.45} and \eqref{9.46}. 
\end{proof}

\begin{lemma}\label{lem9.12}
Let notation be as before. Let $1\leq l\leq r_v+i-1$. Let $\textbf{j}=(\text{nontrivial}, \text{trivial},\text{nontrivial})$. Then 
\begin{align*}
\frac{\mathcal{I}_{i,l}^{\mathrm{small}}(\textbf{j})}{\Vol(K_v[r_v'])}=\frac{\varepsilon(1/2,\sigma_v,\psi_v)\zeta_v(1)^4q_v^{-2}}{\omega_v'(-1)\varepsilon(1/2,\pi_v',\psi_v)}\mathbf{1}_{\substack{r_v'-l=1\\
r_v'=r_v\\ i=0}}\sum_{\substack{r_{\chi_v}=1\\ 
r_{\sigma_v\otimes\chi_v}=r_v}}\frac{\varepsilon(1/2,\sigma_v\otimes\chi_v^{-1},\psi_v)}{\varepsilon(1/2,\pi_v'\otimes\omega_v'^{-1}\chi_v^{-1},\psi_v)}.
\end{align*}
In particular, we have $\mathcal{I}_{i,l}^{\mathrm{small}}(\textbf{j})\ll q_v^{-r_v'-1}$, where the implied constant is absolute.
\end{lemma}
\begin{proof}
By \eqref{eq9.21} (which implies that $j=0$) and \eqref{9.43} we have
\begin{equation}\label{eq9.48}
\mathcal{I}_{i,l}^{\mathrm{small}}(\textbf{j})=\Vol(K_v[r_v'])|W_v'(I_2)|^{-2}\sum_{c\in\mathbb{Z}}q_v^{-c}\cdot \mathcal{M}_{\text{trivial}}'(x_v)\cdot \mathcal{J}(c),
\end{equation}
where for $c\in \mathbb{Z}$, $\mathcal{J}(c)$ is defined by 
\begin{align*}
\frac{1}{|W_v^{\circ}(I_2)|^2|W_v'(I_2)|^2}\int_{\mathcal{O}_v^{\times}}\sum_{\beta_v\in \varpi_v^l\mathcal{O}_v^{\times}/(1+\mathfrak{p}_v^{r_v'})}\mathcal{M}_{\text{nontrivial}}(x_v)\cdot \overline{\mathcal{M}_{\text{nontrivial}}'(x_v)}d\gamma_v.
\end{align*}

By definition, $|W_v^{\circ}(I_2)|^{-2}|W_v'(I_2)|^{-2}\mathcal{M}_{\text{nontrivial}}(x_v)\overline{\mathcal{M}_{\text{nontrivial}}'(x_v)}$ is equal to 
\begin{align*}
&\overline{\omega_v'(-1)\varepsilon(1/2,\pi_v',\psi_v)}\varepsilon(1/2,\sigma_v,\psi_v)q_v^{2l-i-r_v-r_v'}\zeta_v(1)^4\mathbf{1}_{c=l-r_v'}\mathbf{1}_{c=l-r_v-2i}\\
&\sum_{\substack{\chi_v\in \widehat{\mathcal{O}_v^{\times}/(1+\mathfrak{p}_v^{r_v-l+i})}\\
r_{\chi_v}=r_v-l+i\\ 
r_{\sigma_v\otimes\chi_v}=r_v-j}}\sum_{\substack{\xi_v\in \widehat{\mathcal{O}_v^{\times}/(1+\mathfrak{p}_v^{r_v'-l})}\\
r_{\xi_v}=r_v'-l\\ 
r_{\pi_v'\otimes\xi_v}=r_v'-j}}\frac{\chi_v(-\beta_v\gamma_v)\varepsilon(1/2,\sigma_v\otimes\chi_v^{-1},\psi_v)\varepsilon(1/2,\chi_v)^2}{\xi_v(-\beta_v\gamma_v)\varepsilon(1/2,\pi_v'\otimes\omega_v'^{-1}\xi_v^{-1},\psi_v)\varepsilon(1/2,\xi_v)^2}.
\end{align*}
As a consequence of orthogonality, the integral $\mathcal{J}(c)$ is equal to 
\begin{equation}\label{9.48}
\frac{\varepsilon(1/2,\sigma_v,\psi_v)q_v^{l-r_v}\zeta_v(1)^2}{\omega_v'(-1)\varepsilon(1/2,\pi_v',\psi_v)}\mathbf{1}_{\substack{c=l-r_v'\\
r_v'=r_v\\ i=0}}\sum_{\substack{r_{\chi_v}=r_v-l\\ 
r_{\sigma_v\otimes\chi_v}=r_v}}\frac{\varepsilon(1/2,\sigma_v\otimes\chi_v^{-1},\psi_v)}{\varepsilon(1/2,\pi_v'\otimes\omega_v'^{-1}\chi_v^{-1},\psi_v)},
\end{equation}
where $\chi_v\in \widehat{\mathcal{O}_v^{\times}/(1+\mathfrak{p}_v^{r_v-l})}$. 

Substituting \eqref{eq9.21} and \eqref{9.48} into \eqref{eq9.48} we then derive that
\begin{align*}
\frac{\mathcal{I}_{i,l}^{\mathrm{small}}(\textbf{j})}{\Vol(K_v[r_v'])}=\mathbf{1}_{\substack{r_v'-l=1\\
r_v'=r_v\\ i=0}}\frac{\varepsilon(1/2,\sigma_v,\psi_v)\zeta_v(1)^4q_v^{-2}}{\omega_v'(-1)\varepsilon(1/2,\pi_v',\psi_v)}\sum_{\substack{r_{\chi_v}=1\\ 
r_{\sigma_v\otimes\chi_v}=r_v}}\frac{\varepsilon(1/2,\sigma_v\otimes\chi_v^{-1},\psi_v)}{\varepsilon(1/2,\pi_v'\otimes\omega_v'^{-1}\chi_v^{-1},\psi_v)}.
\end{align*}

Executing the trivial bound for the sum over $\chi_v$'s we then obtain $\mathcal{I}_{i,l}^{\mathrm{small}}(\textbf{j})\ll q_v^{-r_v'-1}$, where the implied constant is absolute. 
\end{proof}

\begin{lemma}\label{lem9.13}
Let notation be as before. Let $1\leq l\leq r_v+i-1$. Let $\textbf{j}=(\text{nontrivial},\text{nontrivial},\text{trivial})$. Then 
\begin{align*}
\frac{\mathcal{I}_{i,l}^{\mathrm{small}}(\textbf{j})}{\Vol(K_v[r_v'])}=\frac{\varepsilon(1/2,\sigma_v,\psi_v)\zeta_v(1)^4q_v^{-2}}{\overline{\omega_v'(-1)\varepsilon(1/2,\pi_v',\psi_v)}}\mathbf{1}_{\substack{r_v'-l=1\\
r_v'=r_v\\ i=0}}\sum_{\substack{r_{\chi_v}=1\\ 
r_{\sigma_v\otimes\chi_v}=r_v}}\frac{\varepsilon(1/2,\sigma_v\otimes\chi_v^{-1},\psi_v)}{\overline{\varepsilon(1/2,\pi_v'\otimes\omega_v'^{-1}\chi_v,\psi_v)}}.
\end{align*}
In particular, we have $\mathcal{I}_{i,l}^{\mathrm{small}}(\textbf{j})\ll q_v^{-r_v'-1}$, where the implied constant is absolute.
\end{lemma}
\begin{proof}
By \eqref{eq9.21} (which implies that $j=0$) and \eqref{9.43} we have
\begin{equation}\label{9.50}
\mathcal{I}_{i,l}^{\mathrm{small}}(\textbf{j})=\Vol(K_v[r_v'])|W_v'(I_2)|^{-2}\sum_{c\in\mathbb{Z}}q_v^{-c}\cdot \mathcal{M}_{\text{trivial}}'(x_v)\cdot \mathcal{J}'(c),
\end{equation}
where for $c\in \mathbb{Z}$, $\mathcal{J}'(c)$ is defined by 
\begin{align*}
\frac{1}{|W_v^{\circ}(I_2)|^2|W_v'(I_2)|^2}\int_{\mathcal{O}_v^{\times}}\sum_{\beta_v\in \varpi_v^l\mathcal{O}_v^{\times}/(1+\mathfrak{p}_v^{r_v'})}\mathcal{M}_{\text{nontrivial}}(x_v)\cdot \mathcal{M}_{\text{nontrivial}}'(x_v)d\gamma_v.
\end{align*}

By definition, $|W_v^{\circ}(I_2)|^{-2}|W_v'(I_2)|^{-2}\mathcal{M}_{\text{nontrivial}}(x_v)\mathcal{M}_{\text{nontrivial}}'(x_v)$ is equal to 
\begin{align*}
&\omega_v'(-1)\varepsilon(1/2,\pi_v',\psi_v)\varepsilon(1/2,\sigma_v,\psi_v)q_v^{2l-i-r_v-r_v'}\zeta_v(1)^4\mathbf{1}_{c=l-r_v'}\mathbf{1}_{c=l-r_v-2i}\\
&\sum_{\substack{\chi_v\in \widehat{\mathcal{O}_v^{\times}/(1+\mathfrak{p}_v^{r_v-l+i})}\\
r_{\chi_v}=r_v-l+i\\ 
r_{\sigma_v\otimes\chi_v}=r_v-j}}\chi_v(-\beta_v\gamma_v)\varepsilon(1/2,\sigma_v\otimes\chi_v^{-1},\psi_v)\varepsilon(1/2,\chi_v)^2\\
&\sum_{\substack{\xi_v\in \widehat{\mathcal{O}_v^{\times}/(1+\mathfrak{p}_v^{r_v'-l})}\\
r_{\xi_v}=r_v'-l\\ 
r_{\pi_v'\otimes\xi_v}=r_v'-j}}\xi_v(-\beta_v\gamma_v)\varepsilon(1/2,\pi_v'\otimes\omega_v'^{-1}\xi_v^{-1},\psi_v)\varepsilon(1/2,\xi_v)^2.
\end{align*}

As a consequence of orthogonality and the relation $\varepsilon(1/2,\chi_v)\varepsilon(1/2,\chi_v^{-1})=\chi_v(-1)$, we obtain 
\begin{align*}
\mathcal{J}'(c)=\omega_v'(-1)&\varepsilon(1/2,\pi_v',\psi_v)\varepsilon(1/2,\sigma_v,\psi_v)q_v^{l-r_v}\zeta_v(1)^2\mathbf{1}_{\substack{c=l-r_v',\
r_v'=r_v,\ i=0}}\\
&\sum_{\substack{\chi_v\in \widehat{\mathcal{O}_v^{\times}/(1+\mathfrak{p}_v^{r_v-l})}\\
r_{\chi_v}=r_v-l\\ 
r_{\sigma_v\otimes\chi_v}=r_v}}\varepsilon(1/2,\sigma_v\otimes\chi_v^{-1},\psi_v)\varepsilon(1/2,\pi_v'\otimes\omega_v'^{-1}\chi_v,\psi_v).
\end{align*}
In conjunction with  \eqref{eq9.21} we  derive from \eqref{eq9.48} that 
\begin{align*}
\frac{\mathcal{I}_{i,l}^{\mathrm{small}}(\textbf{j})}{\Vol(K_v[r_v'])}=\mathbf{1}_{\substack{r_v'-l=1\\
r_v'=r_v\\ i=0}}\frac{\varepsilon(1/2,\sigma_v,\psi_v)\zeta_v(1)^4q_v^{-2}}{\overline{\omega_v'(-1)\varepsilon(1/2,\pi_v',\psi_v)}}\sum_{\substack{r_{\chi_v}=1\\ 
r_{\sigma_v\otimes\chi_v}=r_v}}\frac{\varepsilon(1/2,\sigma_v\otimes\chi_v^{-1},\psi_v)}{\overline{\varepsilon(1/2,\pi_v'\otimes\omega_v'^{-1}\chi_v,\psi_v)}}.
\end{align*}

Executing the trivial bound for the sum over $\chi_v$'s we then obtain $\mathcal{I}_{i,l}^{\mathrm{small}}(\textbf{j})\ll q_v^{-r_v'-1}$, where the implied constant is absolute. 
\end{proof}

\subsubsection{Triple nontrivial case}\label{sec9.9.4}
\begin{lemma}\label{lem9.14}
Let notation be as before. Let $1\leq l\leq r_v+i-1$. Let $\textbf{j}=(\text{nontrivial}, \text{nontrivial},\text{nontrivial})$. Then 
\begin{equation}\label{eq9.51}
\mathcal{I}_{i,l}^{\mathrm{small}}(\textbf{j})=\frac{\zeta_v(1)^4\varepsilon(1/2,\sigma_v,\psi_v)\Vol(K_v[r_v'])\mathbf{1}_{2i=r_v'-r_v}}{q_v^{r_v+i-l}}\sum_{j\in \mathbb{Z}}q_v^j\cdot \mathcal{S}(j),
\end{equation}
where  
\begin{multline*}
\mathcal{S}(j):=\sum_{\substack{\chi_v,\xi_v\in \widehat{\mathcal{O}_v^{\times}/(1+\mathfrak{p}_v^{r_v'-l})}\\
r_{\chi_v}=r_{\xi_v}=r_v'-l\\ 
r_{\chi_v^{-1}\xi_v}=r_v-l+i\\ 
r_{\pi_v'\otimes\chi_v}=r_{\pi_v'\otimes\xi_v}=r_v'-j\\
r_{\sigma_v\otimes\chi_v^{-1}\xi_v}=r_v-j}}\varepsilon(1/2,\pi_v'\otimes\omega_v'^{-1}\chi_v^{-1},\psi_v)\varepsilon(1/2,\chi_v)^2\varepsilon(1/2, \chi_v^{-1} \xi_v)^2\\
\overline{\varepsilon(1/2,\pi_v'\otimes\omega_v'^{-1}\xi_v^{-1},\psi_v)\varepsilon(1/2,\xi_v)^2}\varepsilon(1/2, \sigma_v \otimes \chi_v \xi_v^{-1}, \psi_v).
\end{multline*}
Moroever, if $\pi_v'$ is twist-minimal, then $\mathcal{S}(j)=0$ unless $j=\min\{0,2l-r_v'\}$.  
\end{lemma}
\begin{proof}
By definition \eqref{9.43},  
\begin{align*}
\mathcal{I}_{i,l}^{\mathrm{small}}(\textbf{j})=\Vol(K_v[r_v'])\sum_jq_v^j\sum_cq_v^{-c}\int_{\mathcal{O}_v^{\times}}\sum_{\beta_v}\frac{\mathcal{M}_{\text{nontrivial}}(x_v)|\mathcal{M}_{\text{nontrivial}}'(x_v)|^2}{|W_v^{\circ}(I_2)|^2|W_v'(I_2)|^4}d\gamma_v,
\end{align*}
where $\beta_v\in \varpi_v^l\mathcal{O}_v^{\times}/(1+\mathfrak{p}_v^{r_v'})$. 

Utilizing the explicit formulas of $\mathcal{M}_{\text{nontrivial}}(x_v)$ and $\mathcal{M}_{\text{nontrivial}}'(x_v)$ in Lemmas \ref{lem8.3} and Lemma \ref{lem8.5}, respectively, along with orthogonality of characters, we derive 
\begin{multline*}
\int_{\mathcal{O}_v^{\times}}\frac{\mathcal{M}_{\text{nontrivial}}(x_v)|\mathcal{M}_{\text{nontrivial}}'(x_v)|^2}{|W_v^{\circ}(I_2)|^2|W_v'(I_2)|^4}d\gamma_v=\zeta_v(1)^5\varepsilon(1/2,\sigma_v,\psi_v)q_v^{3l-2r_v'-r_v-i}\\
\mathbf{1}_{c=l-r_v-2i}\mathbf{1}_{c=l-r_v'}\sum_{\substack{\chi_v,\xi_v\in \widehat{\mathcal{O}_v^{\times}/(1+\mathfrak{p}_v^{r_v'-l})}\\
r_{\chi_v}=r_{\xi_v}=r_v'-l\\ 
r_{\pi_v'\otimes\chi_v}=r_{\pi_v'\otimes\xi_v}=r_v'-j\\
r_{\chi_v^{-1}\xi_v}=r_v-l+i\\ 
r_{\sigma_v\otimes\chi_v^{-1}\xi_v}=r_v-j}}\varepsilon(1/2,\pi_v'\otimes\omega_v'^{-1}\chi_v^{-1},\psi_v)\varepsilon(1/2,\chi_v)^2\\
\overline{\varepsilon(1/2,\pi_v'\otimes\omega_v'^{-1}\xi_v^{-1},\psi_v)\varepsilon(1/2,\xi_v)^2}\varepsilon(1/2, \sigma_v \otimes \chi_v \xi_v^{-1}, \psi_v)
\varepsilon(1/2, \chi_v^{-1} \xi_v)^2.
\end{multline*}

Suppose $\pi_v'$ is twist-minimal, then $r_{\pi_v'\otimes\xi_v}=r_v'-j\geq r_v'$, implying that $j\leq 0$. 
\begin{itemize}
\item Suppose $j<\min\{0,2l-r_v'\}$. So $j\leq -1$. Then it follows from  
\begin{align*}
r_v'-j=r_{\pi_v'\otimes\xi_v}\leq \max\{r_v',2r_{\xi_v}\}=\max\{r_v',2r_v'-2l\}
\end{align*}
that $r_v'-j\leq 2r_v'-2l$, namely, $j\geq 2l-r_v'$. A contradiction.
\item Suppose $j> \min\{0,2l-r_v'\}$. Then $2l-r_v'<j\leq 0$. By $2(r_v'-l)=2r_{\xi_v}\leq \max\{r_v',r_{\pi_v'\otimes\xi_v}\}=\max\{r_v',r_v'-j\}\leq r_v'-j$, forcing that $j\leq 2l-r_v'$. Another contradiction. 
\end{itemize}
Therefore, $\mathcal{S}(j)=0$ unless $j=\min\{0,2l-r_v'\}$. Hence, Lemma \ref{lem9.14} follows. 
\end{proof}

\begin{comment}

\bigskip
\bigskip
\bigskip

Recall that $1\leq l\leq r_v+i-1$ and $j=\min\{0,2l-r_v'\}$. 

\begin{multline*}
\mathcal{I}_{i,l}^{\mathrm{small}}(\textbf{j})=\frac{\zeta_v(1)^4\varepsilon(1/2,\sigma_v,\psi_v)\Vol(K_v[r_v'])\mathbf{1}_{2i=r_v'-r_v}}{q_v^{r_v+i-l}}q_v^{\min\{0,2l-r_v'\}}\\
\sum_{\substack{\chi_v,\xi_v\in \widehat{\mathcal{O}_v^{\times}/(1+\mathfrak{p}_v^{r_v'-l})}\\
r_{\chi_v}=r_{\xi_v}=r_v'-l\\ 
r_{\pi_v'\otimes\chi_v}=r_{\pi_v'\otimes\xi_v}=r_v'-j\\
r_{\chi_v^{-1}\xi_v}=r_v-l+i\\ 
r_{\sigma_v\otimes\chi_v^{-1}\xi_v}=r_v-j}}\chi_v\xi_v^{-1} (\varpi_v^{r_v'-r_v})\varepsilon(1/2,\pi_v'\otimes\omega_v'^{-1}\chi_v^{-1},\psi_v)\\
\varepsilon(1/2,\chi_v)^2
\overline{\varepsilon(1/2,\pi_v'\otimes\omega_v'^{-1}\xi_v^{-1},\psi_v)\varepsilon(1/2,\xi_v)^2}\varepsilon(1/2, \sigma_v \otimes \chi_v \xi_v^{-1}, \psi_v)
\varepsilon(1/2, \chi_v^{-1} \xi_v)^2
\end{multline*}

The trivial bound is 
\begin{align*}
\mathcal{I}_{i,l}^{\mathrm{small}}(\textbf{j})\ll \frac{q_v^{-r_v'}}{q_v^{r_v+i-l}}\cdot q_v^{\min\{0,2l-r_v'\}}\cdot q_v^{2r_v'-2l}\ll q_v^{-r_v'}\cdot q_v^{2r_v'-2l}\cdot q_v^{l-r_v-i}\ll q_v^{\frac{r_v'-r_v}{2}-l}
\end{align*}

\bigskip
\bigskip
\bigskip
\bigskip

\end{comment}

\subsection{Expressions as Character Sums}
Suppose that $v$ is not above $2$. Then $\sigma_v$ and $\pi_v'$ are dihedral. For simplicity we denote by $F=F_v$ in this section. 

\subsubsection{Construction of supercuspidal representations}
We briefly recall the construction of $\sigma_v$ (see \cite[\textsection 1.1]{JL70}). There is a quadratic field extension $E/F$ and a character $\kappa$ of $E^{\times}$ that is not trivial on the kernel of the norm map $N_{E/F}$, such that $\sigma_v=\omega_{\kappa}$, which is the induction from the Weil representation of $\mathrm{SL}_2(F)$ associated with $\kappa$ and $\psi_v$. The central character of $\sigma_v\simeq \omega_{\kappa}$ is $\eta_{E/F}\kappa|_{F^{\times}}$, where $\eta_{E/F}$ is the quadratic character associated with $E/F$. Moreover, for any character $\mu$ of $F^{\times}$, we have the isomorphism $\omega_{\kappa}\otimes\mu\simeq \omega_{\kappa(\mu\circ N_{E/F})}$, and 
\begin{equation}\label{9.51}
\varepsilon(s,\omega_{\kappa}\otimes\mu,\psi_F)=\gamma_{E/F}\cdot\varepsilon(s,\kappa\cdot (\mu\circ N_{E/F}),\psi_E),
\end{equation}
where $\gamma_{E/F}$ is a fourth root of unity with the property $\gamma_{E/F}^2=\eta_{E/F}(-1)$, and $\psi_E=\psi_v\circ \Tr_{E/F}.$ Here $\Tr_{E/F}$ is the trace map. 

Let $\mathfrak{P}$ be a uniformizer in $E$. Let $\mathcal{O}_E$ be the ring of integral elements of $E$. Let $f=f_{E/F}$ be the  degree of the residue field extension. 

For any character $\chi_v$, we have (see \cite[Lemma 1.1.1]{Sch02})
\begin{equation}\label{9.52}
\varepsilon(1/2,\chi_v)=\zeta_v(1)^{-1}q_v^{\frac{r_{\chi_v}}{2}}\int_{\mathcal{O}_v^{\times}}\chi_v^{-1}(\alpha_v)\psi_v(\varpi_v^{-r_{\chi_v}}\alpha_v)d^{\times}\alpha_v.
\end{equation}

As a consequence of \eqref{9.51} and \eqref{9.52}, we obtain that 
\begin{align*}
\varepsilon(1/2,\sigma_v\otimes\mu,\psi_v)=\gamma_{E/F}q_v^{\frac{r_{\sigma_v\otimes\mu}}{2}}\int_{\mathcal{O}_{E}^{\times}}\kappa^{-1}(x)\mu^{-1}(N_{E/F}(x))\psi_E(\mathfrak{P}^{-r_{\sigma_v\otimes\mu}/f}x)dx,
\end{align*}
where $\mu$ us an arbitrary character of $E^{\times}$. 

Since $\pi_v'$ is dihedral, there is a quadratic field extension $E'/F$ and a character $\kappa'$ of $E'$ such that $\pi_v'\simeq \omega_{\kappa'}$. Let $f'=f_{E'/F}$ be the degree, and $\mathfrak{P}'$ be a uniformizer of $E'$.

\subsubsection{Character sums}
Let $\mathcal{S}(j)$ be defined as in \eqref{eq9.51} in Lemma \ref{lem9.14}, i.e., 
\begin{align*}
\sum_{\substack{\chi_v,\xi_v\in \widehat{\mathcal{O}_v^{\times}/(1+\mathfrak{p}_v^{r_v'-l})}\\
r_{\chi_v}=r_{\xi_v}=r_v'-l\\ 
r_{\pi_v'\otimes\chi_v}=r_{\pi_v'\otimes\xi_v}=r_v'-j\\
r_{\chi_v^{-1}\xi_v}=r_v-l+i\\ 
r_{\sigma_v\otimes\chi_v^{-1}\xi_v}=r_v-j}}\frac{\varepsilon(1/2,\pi_v'\otimes\omega_v'^{-1}\chi_v^{-1},\psi_v)\varepsilon(1/2,\sigma_v\otimes\chi_v\xi_v^{-1},\psi_v)\varepsilon(1/2,\chi_v)^2}{\varepsilon(1/2,\pi_v'\otimes\omega_v'^{-1}\xi_v^{-1},\psi_v)\varepsilon(1/2,\xi_v)^2\varepsilon(1/2,\chi_v^{-1}\xi_v)^{-2}}.
\end{align*}
In this section we aim to establish refined description for $\mathcal{S}(j)$.  

\begin{lemma}\label{lem9.15}
Let notation be as before. Then 
\begin{align*}
\mathcal{S}(j)=&\frac{q_v^{5(r_v-l)+\frac{3(r_v-j)}{2}}\mathbf{1}_{r_v=r_v'-i}}{\zeta_v(1)^{8}}\int_{\mathcal{O}_{E}^{\times}}\int_{(\mathcal{O}_{E'}^{\times})^{\otimes 2}}\kappa'^{-1}\omega_v'\circ N_{E'/F}(x_1x_2^{-1})\\
&\kappa^{-1}(x_3)\psi_{E'}(\mathfrak{P}'^{-\frac{r_v-j}{f'}}(x_1-x_2))\psi_E(\mathfrak{P}^{-\frac{r_v-j}{f}}x_3)\mathcal{K}(x_1,x_2,x_3)dx_1dx_2dx_3,
\end{align*}
where $\mathcal{K}(x_1,x_2,x_3)$ is defined by 
\begin{align*}
\int \sum_{\substack{\delta_v\in \mathcal{O}_v^{\times}/(1+\mathfrak{p}_v^{r_v-l})\\\alpha_v\alpha_v'\equiv \delta_vN_{E'/F}(x_1)\pmod{1+\mathfrak{p}_v^{r_v-l}}\\
\beta_v\beta_v'\equiv \delta_vN_{E'/F}(x_2)\pmod{1+\mathfrak{p}_v^{r_v-l}}\\
\gamma_v\gamma_v'\equiv \delta_vN_{E/F}(x_3)\pmod{1+\mathfrak{p}_v^{r_v-l}}}}\psi_v(\varpi_v^{l-r_v}(\alpha_v+\alpha_v'+\beta_v+\beta_v'+\gamma_v+\gamma_v')).
\end{align*}
Here the integral is over $(\mathcal{O}_v^{\times})^{\otimes 6}$ relative to $d^{\times}\alpha_v'd^{\times}\alpha_vd^{\times}\beta_v'd^{\times}\beta_vd^{\times}\gamma_v'd^{\times}\gamma_v$. 
\end{lemma}
\begin{proof}
Denote by $r_{\chi_v}=r_{\xi_v}=r_v'-l=:r_1$, and $r_{\chi_v^{-1}\xi_v}=r_v-l+i=:r_2$. Since 
\begin{align*}
\overline{\varepsilon(1/2,\xi_v)}=\varepsilon(1/2,\xi_v)^{-1},\ \ \varepsilon(1/2,\xi_v)\varepsilon(1/2,\xi_v^{-1})=\xi_v(-1),
\end{align*}
then $\overline{\varepsilon(1/2,\xi_v)^2}=\varepsilon(1/2,\xi_v^{-1})^2$. Hence, by \eqref{9.52}, 
\begin{align*}
\varepsilon(1/2,\chi_v)^2\overline{\varepsilon(1/2,\xi_v)^2}\varepsilon(1/2,\chi_v^{-1}\xi_v)^2=\varepsilon(1/2,\chi_v)^2\varepsilon(1/2,\xi_v^{-1})^2\varepsilon(1/2,\chi_v^{-1}\xi_v)^2
\end{align*}
is further equal to 
\begin{align*}
&\zeta_v(1)^{-6}q_v^{2r_1+r_2}\int_{\mathcal{O}_v^{\times}}\int_{\mathcal{O}_v^{\times}}\int_{\mathcal{O}_v^{\times}}\int_{\mathcal{O}_v^{\times}}\int_{\mathcal{O}_v^{\times}}\int_{\mathcal{O}_v^{\times}}\chi_v^{-1}(\alpha_v\alpha_v')\xi_v(\beta_v\beta_v')\chi_v\xi_v^{-1}(\gamma_v\gamma_v')\\
&\psi_v(\varpi_v^{-r_1}(\alpha_v+\alpha_v'+\beta_v+\beta_v')+\varpi_v^{-r_2}(\gamma_v+\gamma_v'))d^{\times}\alpha_v'd^{\times}\alpha_vd^{\times}\beta_v'd^{\times}\beta_vd^{\times}\gamma_v'd^{\times}\gamma_v.
\end{align*}

Let $r_{\pi_v'\otimes\chi_v}=r_{\pi_v'\otimes\xi_v}=r_v'-j$ and $r_{\sigma_v\otimes\chi_v^{-1}\xi_v}=r_v-j$ be as in Lemma \ref{lem9.14}. Then 
\begin{align*}
\varepsilon(1/2,\pi_v'\otimes\omega_v'^{-1}\chi_v^{-1},\psi_v)\varepsilon(1/2,\sigma_v\otimes\chi_v\xi_v^{-1},\psi_v)\overline{\varepsilon(1/2,\pi_v'\otimes\omega_v'^{-1}\xi_v^{-1},\psi_v)}
\end{align*}
is equal to 
\begin{align*}
&\gamma_{E/F}q_v^{r_v'-j}q_v^{\frac{r_v-j}{2}}\int_{\mathcal{O}_{E}^{\times}}\int_{\mathcal{O}_{E'}^{\times}}\int_{\mathcal{O}_{E'}^{\times}}\kappa'^{-1}(x_1)\kappa'(x_2)\kappa^{-1}(x_3)\psi_{E'}(\mathfrak{P}'^{-\frac{r_v'-j}{f'}}(x_1-x_2))\\
&\psi_E(\mathfrak{P}^{-\frac{r_v-j}{f}}x_3)\omega_v'\chi_v(N_{E'/F}(x_1))\omega_v'^{-1}\xi_v^{-1}(N_{E'/F}(x_2))\chi_v^{-1}\xi_v(N_{E/F}(x_3))dx_1dx_2dx_3.
\end{align*}

Combining the above calculations and properties of Gauss sums, we obtain 
\begin{align*}
\mathcal{S}(j)=&\frac{q_v^{2r_1+r_2}q_v^{r_v'-j}q_v^{\frac{r_v-j}{2}}}{\zeta_v(1)^{6}}\int_{(\mathcal{O}_v^{\times})^{\otimes 6}}\int_{\mathcal{O}_{E}^{\times}}\int_{(\mathcal{O}_{E'}^{\times})^{\otimes 2}}\kappa'^{-1}\omega_v'\circ N_{E'/F}(x_1x_2^{-1})\kappa^{-1}(x_3)\\
&\psi_{E'}(\mathfrak{P}'^{-\frac{r_v'-j}{f'}}(x_1-x_2))\psi_E(\mathfrak{P}^{-\frac{r_v-j}{f}}x_3)\psi_v(\cdots)\mathcal{E}dx_1dx_2dx_3d^{\times}\cdots ,
\end{align*}
where $d^{\times}\cdots:=d^{\times}\alpha_v'd^{\times}\alpha_vd^{\times}\beta_v'd^{\times}\beta_vd^{\times}\gamma_v'd^{\times}\gamma_v$, $\psi_v(\cdots):=\psi_v(\varpi_v^{-r_1}(\alpha_v+\alpha_v'+\beta_v+\beta_v')+\varpi_v^{-r_2}(\gamma_v+\gamma_v'))$, and  
\begin{align*}
\mathcal{E}:=\sum_{\substack{\chi_v,\xi_v\in \widehat{\mathcal{O}_v^{\times}/(1+\mathfrak{p}_v^{r_v'-l})}}}&\chi_v(N_{E'/F}(x_1))\xi_v^{-1}(N_{E'/F}(x_2))\chi_v^{-1}\xi_v(N_{E/F}(x_3))\\
&\chi_v^{-1}(\alpha_v\alpha_v')\xi_v(\beta_v\beta_v')\chi_v\xi_v^{-1}(\gamma_v\gamma_v').
\end{align*}

In particular, in the above expression for $\mathcal{S}(j)$, the constraints on $r_{\chi_v}$, $r_{\xi_v}$,  $r_{\chi_v^{-1}\xi_v}$,  $r_{\pi_v'\otimes\chi_v}$, $r_{\pi_v'\otimes\xi_v}$ and $r_{\sigma_v\otimes\chi_v^{-1}\xi_v}$ in the sum over $\chi_v$ and $\xi_v$ are removed. 

\begin{comment}
\begin{align*}
&\sum_{\substack{\chi_v,\xi_v\in \widehat{\mathcal{O}_v^{\times}/(1+\mathfrak{p}_v^{r_v'-l})}}}\int_{\mathcal{O}_{E}^{\times}}\int_{\mathcal{O}_{E'}^{\times}}\int_{\mathcal{O}_{E'}^{\times}}\kappa'^{-1}(x_1)\kappa'(x_2)\kappa^{-1}(x_3)\psi_{E'}(\mathfrak{P}'^{-\frac{r_v'-j}{f'}}(x_1-x_2))\\
&\psi_E(\mathfrak{P}^{-\frac{r_v-j}{f}}x_3)\omega_v'\chi_v(N_{E'/F}(x_1))\omega_v'^{-1}\xi_v^{-1}(N_{E'/F}(x_2))\chi_v^{-1}\xi_v(N_{E/F}(x_3))dx_1dx_2dx_3\\
&\int_{\mathcal{O}_v^{\times}}\int_{\mathcal{O}_v^{\times}}\int_{\mathcal{O}_v^{\times}}\int_{\mathcal{O}_v^{\times}}\int_{\mathcal{O}_v^{\times}}\int_{\mathcal{O}_v^{\times}}\chi_v^{-1}(\alpha_v\alpha_v')\xi_v(\beta_v\beta_v')\chi_v\xi_v^{-1}(\gamma_v\gamma_v')\\
&\psi_v(\varpi_v^{-r_1}(\alpha_v+\alpha_v'+\beta_v+\beta_v')+\varpi_v^{-r_2}(\gamma_v+\gamma_v'))d^{\times}\alpha_v'd^{\times}\alpha_vd^{\times}\beta_v'd^{\times}\beta_vd^{\times}\gamma_v'd^{\times}\gamma_v.
\end{align*}

\begin{align*}
&\sum_{\substack{\chi_v\in \widehat{\mathcal{O}_v^{\times}/(1+\mathfrak{p}_v^{r_v'-l})}}}\chi_v(N_{E'/F}(x_1)N_{E/F}(x_3^{-1}))\chi_v(\alpha_v^{-1}\alpha_v'^{-1}\gamma_v\gamma_v')\\
&\sum_{\substack{\xi_v\in \widehat{\mathcal{O}_v^{\times}/(1+\mathfrak{p}_v^{r_v'-l})}}}\xi_v^{-1}(N_{E'/F}(x_2)N_{E/F}(x_3^{-1}))\xi_v^{-1}(\beta_v^{-1}\beta_v'^{-1}\gamma_v\gamma_v')
\end{align*}

\end{comment}
 
By orthogonality of characters, we obtain  
\begin{align*}
\mathcal{E}=\zeta_v(1)^{-2}q_v^{2(r_v'-l)}\mathbf{1}_{N_{E'/F}(x_2^{-1})\beta_v\beta_v'\equiv N_{E'/F}(x_1^{-1})\alpha_v\alpha_v'\equiv N_{E/F}(x_3^{-1})\gamma_v\gamma_v'\pmod{1+\mathfrak{p}_v^{r_v'-l}}}.
\end{align*}

Substituting this into the expression of $\mathcal{S}(j)$ we then derive that 
\begin{align*}
\mathcal{S}(j)=&\frac{q_v^{4(r_v'-l)+r_v-l+i}q_v^{r_v'-j}q_v^{\frac{r_v-j}{2}}}{\zeta_v(1)^{8}}\int_{\mathcal{O}_{E}^{\times}}\int_{(\mathcal{O}_{E'}^{\times})^{\otimes 2}}\kappa'^{-1}\omega_v'\circ N_{E'/F}(x_1x_2^{-1})\kappa^{-1}(x_3)\\
&\psi_{E'}(\mathfrak{P}'^{-\frac{r_v'-j}{f'}}(x_1-x_2))\psi_E(\mathfrak{P}^{-\frac{r_v-j}{f}}x_3)\mathcal{K}^*(x_1,x_2,x_3)dx_1dx_2dx_3,
\end{align*}
where
\begin{multline*}
\mathcal{K}^*(x_1,x_2,x_3):=\int_{(\mathcal{O}_v^{\times})^{\otimes 6}}\psi_v(\varpi_v^{-r_1}(\alpha_v+\alpha_v'+\beta_v+\beta_v')+\varpi_v^{-r_2}(\gamma_v+\gamma_v'))\\
\mathbf{1}_{N_{E'/F}(x_2^{-1})\beta_v\beta_v'\equiv N_{E'/F}(x_1^{-1})\alpha_v\alpha_v'\equiv N_{E/F}(x_3^{-1})\gamma_v\gamma_v'\pmod{1+\mathfrak{p}_v^{r_v'-l}}}\\
d^{\times}\alpha_v'd^{\times}\alpha_vd^{\times}\beta_v'd^{\times}\beta_vd^{\times}\gamma_v'd^{\times}\gamma_v.
\end{multline*}

By a change of variable, $\mathcal{K}^*(x_1,x_2,x_3)\equiv 0$ unless  $r_1=r_2$, i.e., $r_v'-l=r_v-l+i$.
 which implies that  $i=r_v'-r_v$. Hence, Lemma \ref{lem9.15} follows.
\end{proof}

\begin{cor}\label{cor9.16}
Let notation be as before. Let $1\leq l\leq r_v+i-1$. Let $\textbf{j}=(\text{nontrivial}, \text{nontrivial},\text{nontrivial})$. Suppose $\pi_v'$ is twist-minimal. Then 
\begin{equation}\label{eq9.57}
\mathcal{I}_{i,l}^{\mathrm{small}}(\textbf{j})\ll q_v^{\min\{-l,l-r_v'\}}\cdot \mathbf{1}_{r_v=r_v',\ i=0},
\end{equation}
where the implied constant is absolute. 
\end{cor}
\begin{proof}
Substituting the constraint $\mathbf{1}_{r_v=r_v'-i}$ in Lemma \ref{lem9.15} into 
Lemma \ref{lem9.14} and utilizing the trivial bound for the sum over $\chi_v$ and $\xi_v$, we then obtain the estimate \eqref{eq9.57}.
\end{proof}
\begin{remark}
Let $\mu_{\mathrm{P}}$ be the Plancherel measure. 
By \cite[Proposition 18.8.5]{Dix77} we have $d\mu_{\mathrm{P}}(\sigma_v)=\mathrm{deg}(\sigma_v)$, where $\mathrm{deg}(\sigma_v)$ is the formal degree of $\sigma_v$. According to  \cite{Zin93} and \cite{SZ96}, $\mathrm{deg}(\sigma_v)=2^{-1}|\gamma_v(0,\sigma_v,\Ad,\psi_v)|$. It then follows from the Plancherel formula that 
\begin{align*}
\int_{\overline{G'}(F_v)}\frac{\mathcal{MC}_{i}(x_v)\cdot |\langle\pi_v'(x_v)W_v',W_v'\rangle|^2}{|W_v^{\circ}(I_2)|^{2}}dx_v\cdot |\gamma_v(0,\sigma_v,\Ad,\psi_v)|\ll |\langle W_v',W_v'\rangle|^2,
\end{align*}
i.e., $\mathcal{I}_{i}\ll |\gamma_v(0,\sigma_v,\Ad,\psi_v)|^{-1}$. Here $\mathcal{I}_i$ is defined by \eqref{eq9.28}. This amounts to the trivial bound \eqref{eq9.57} when $r_v=r_v'$. 
\end{remark}

\subsection{Proof of Proposition \ref{prop9.6}}\label{sec9.11}
Notice that replacing the Whittaker function $W_v'(x_v)$ with $W_v'(x_v)\chi_v(\det x_v)$ for any  unitary character 
$\chi_v$ of $F_v^{\times}$ does not affect  $\mathcal{I}_i$. Hence we may assume $\pi_v'$ is twist-minimal. Therefore, Proposition \ref{prop9.6} follows from gathering Lemmas \ref{lem9.6} and \ref{lemma9.8}, \eqref{eq9.30} in Lemma \ref{lem9.8}, \eqref{9.37}, Lemmas \ref{lem9.9}, \ref{lemma9.11}, \ref{lem9.11}, \ref{lem9.12}, \ref{lem9.13}, and Corollary \ref{cor9.16} into the decomposition \eqref{eq9.25}.

\section{Nonarchimedean Integrals on the Dual Side \RNum{3}}\label{sec10} 
In this section we treat the remaining supercuspidal case, where $\pi_v'$ is supercuspidal
but the representation appearing in the $\mathrm{GL}_2$-spectrum is not supercuspidal.
We evaluate the corresponding nonarchimedean local integrals in
$\mathcal{M}_{\cusp}^{\du}(s_1,s_2)$ and $\mathcal{M}_{\Eis}^{\du}(s_1,s_2)$,
and derive explicit formulas together with sharp upper bounds.

\subsection{Notation}
Let $v$ be a finite place, and $\pi_v'$ be the $v$-th component of $\pi'$. Suppose $\pi_v'$ is supercuspidal as in \textsection\ref{8.1.1}. Let $\sigma_v=\eta_v\boxplus\eta_v^{-1}$ be a principal series or $\sigma_v\subseteq \eta_v\boxplus\eta_v^{-1}$ be a twisted Steinberg representation (in which case $\eta_v$ is quadratic). We will continue to use the notation introduced in \textsection\ref{sec9}. Let $W_v^{\circ}$ and $W_v'$ be local newforms of the Whittaker models of $\sigma_v$ and $\pi_v'$, respectively. As in \textsection\ref{8.1.1}, we set 
\begin{align*}
W_v=W_{v,n_v}=\sum_{i=0}^{n_v}\xi_{\sigma_v}(\mathfrak{p}_v^{i},\mathfrak{p}_v^{n_v})q_v^{\frac{i-n_v}{2}}\sigma_v\left(\begin{pmatrix}
1\\
&\varpi_v^{i}
\end{pmatrix}\right)W_v^{\circ},\  0\leq n_v\leq r_v'-r_v. \tag{\ref{eq8.1}}
\end{align*}
 
For $0\leq i\leq n_v\leq r_v'-r_v$, we define the normalized local triple integral by   
\begin{equation}\label{eq9.28.}
\mathcal{I}_{i}:=|W_v'(I_2)|^{-4}\langle W_v^{\circ}, W_v^{\circ}\rangle^{-1}\int_{\overline{G'}(F_v)}\mathcal{MC}_{i}(x_v)\cdot |\langle\pi_v'(x_v)W_v',W_v'\rangle|^2dx_v,
\end{equation}
where $\mathcal{MC}_{i}(x_v)$ is the matrix coefficient relative to $\sigma_v(\diag(1,\varpi_v^i))W_v^{\circ}$ as defined in \eqref{8.22}. Explicitly, $\mathcal{MC}_{i}(x_v)=\langle \sigma_v(x_v)\sigma_v(\diag(1,\varpi_v^i))W_v^{\circ},\sigma_v(\diag(1,\varpi_v^i))W_v^{\circ}\rangle$. 

In this section we will compute $\mathcal{I}_{i}$ in the case that $\sigma_v$ is a principal series or a special representation, and $\pi_v'$ is supercuspidal.

\subsection{$\eta_v$ is Unramified}
Suppose that $r_{\eta_v}=0$, i.e., $\eta_v$ is unramified. By \cite[Corollary 1.2]{HS01}, $\mathcal{I}_{i}\equiv 0$ if $\sigma_v$ is special. As a consequence, we may assume that $\sigma_v$ is a principal series. 

For $\Re(s)\gg 1$ and $x_v\in G(F_v)$, we define  
\begin{align*}
h_v^{\sharp}(x_v,s;i):=\frac{\eta_v^2(\varpi_v^i)\eta_v(\det x_v)|\det x_v|_v^{s}}{q_v^{2is}L_v(2s,\eta_v^2)}\int_{F_v^{\times}}\Phi_v^{\sharp}\left((0,t)x_v\diag(1,\varpi_v^i)\right)\eta_v^2(t)|t|_v^{2s}d^{\times}t,
\end{align*}
where $\Phi_v^{\sharp}(t_{1,v},t_{2,v}):=\mathbf{1}_{\mathcal{O}_v}(t_{1,v})\mathbf{1}_{\mathcal{O}_v}(t_{2,v})$. Then $h_v^{\sharp}(\cdot,s;i)$ is a spherical section in $\sigma_v$ with $h_v^{\sharp}(I_2,s;i)\equiv 1$. For given $x_v\in G(F_v)$, $h_v^{\sharp}(\cdot,s;i)$ extends to a holomorphic function in $s$. Let $\Re(s)\gg 1$. Denote by 
\begin{align*}
\Psi_v^{\sharp}(s;i):=\int_{N(\mathbb{A}_F)\backslash\overline{G}(\mathbb{A}_F)}W_v'(x_v)\overline{W_v'(x_v)}h_v^{\sharp}(x_v,s;\diag(1,\varpi_v^i))dx_v.
\end{align*}
Let $\Re(s)\gg 1$. Define 
\begin{equation}\label{equ10.1}
\mathcal{I}_{i}(s):=\zeta_v(1)\langle W_v^{\circ}, W_v^{\circ}\rangle^{-1}|W_v'(I_2)|^{-4}\cdot |\Psi_v^{\sharp}(s;i)|^2,
\end{equation}
which converges absolutely in $\Re(s)\gg 1$ and admits a meromorphic continuation to $s\in \mathbb{C}$. As a consequence of Lemma 3.4.2 (and the remark immediately following it) in \cite{MV10} we have $\mathcal{I}_i(0)=\mathcal{I}_i$.

%(or its generalization \cite[Proposition 5.1]{Hsi21}, noting that the hypothesis (Hb) therein is satisfied if $\sigma_v$ is not special since $\pi_v'$ is supercuspidal) 

\begin{prop}\label{prop10.1}
Let notation be as before. Suppose  $\eta_v$ is unramified. Let $0\leq i\leq r_v'$. Then
\begin{align*}
\mathcal{I}_i\ll r_v'^2\cdot q_v^{-r_v'},
\end{align*}	
where the implied constant is absolute. 
\end{prop}

\subsubsection{Expansion of $\Psi_v^{\sharp}(\cdot;i)$}
Under the Iwasawa coordinates the function $\Psi_v^{\sharp}(\cdot;i)$ boils down to 
\begin{align*}
\sum_{j\in \mathbb{Z}}q_v^{j(1-s)}\int_{K_v}\Big|W_v'\left(\begin{pmatrix}
\varpi_v^j\\
&1
\end{pmatrix}k_v\right)\Big|^2\sum_{m\in\mathbb{Z}}\frac{\Phi_v^{\sharp}(\eta \varpi_v^mk_v\diag(1,\varpi_v^i))\eta_v(\varpi_v^{j+2m})}{\eta_v^{-2}(\varpi_v^i)q_v^{2(m+i)s}L_v(2s,\eta_v^2)}dk_v.
\end{align*}

Let $\Re(s)\gg 1$. Then 
\begin{equation}\label{eq10.1}
\Psi_v^{\sharp}(s;i)=\Psi_v^{\mathrm{small}}(s;i)+\Psi_v^{\mathrm{big}}(s;i),
\end{equation}
where $\Psi_v^{\mathrm{small}}(s;i)$ is defined by 
\begin{align*}
\sum_{j\in \mathbb{Z}}q_v^{j(1-s)}\int \Big|W_v'\left(\begin{pmatrix}
\varpi_v^j\\
&1
\end{pmatrix}k_v\right)\Big|^2\sum_{m\in\mathbb{Z}}\frac{\Phi_v^{\sharp}(\eta \varpi_v^mk_v\diag(1,\varpi_v^i))}{\eta_v^{-1}(\varpi_v^{j+2(m+i)})q_v^{2(m+i)s}L_v(2s,\eta_v^2)}dk_v;
\end{align*}
with $k_v$ ranges over $K_v[1]$; and $\Psi_v^{\mathrm{big}}(s;i)$ is defined by 
\begin{align*}
\sum_{j\in \mathbb{Z}}q_v^{j(1-s)}\int\Big|W_v'\left(\begin{pmatrix}
\varpi_v^j\\
&1
\end{pmatrix}k_v\right)\Big|^2\sum_{m\in\mathbb{Z}}\frac{\Phi_v^{\sharp}(\eta \varpi_v^mk_v\diag(1,\varpi_v^i))}{\eta_v^{-1}(\varpi_v^{j+2(m+i)})q_v^{2(m+i)s}L_v(2s,\eta_v^2)}dk_v,
\end{align*}
with $k_v$ ranges over $K_v-K_v[1]$. 

\subsubsection{Calculation of $\Psi_v^{\mathrm{big}}(s;i)$}
\begin{lemma}\label{lemma10.2}
Let notation be as above. Let $0\leq i\leq r_v'$.  Let $\Re(s)\gg 1$. Then 
\begin{equation}\label{equ10.2}
\Psi_v^{\mathrm{big}}(s;i)=\eta_v(\varpi_v^{2i-r_v'})\Vol(K_v[r_v'])|W_v'(I_2)|^2q_v^{(r_v'-2i)s}.
\end{equation}
In particular, $\Psi_v^{\mathrm{big}}(s;i)$ extends to a holomorphic function in $s\in \mathbb{C}$. 
\end{lemma}
\begin{proof}
Let $k_v\in K_v-K_v[1]$. We may write  $k_v=\begin{pmatrix}
1& b\\
&1
\end{pmatrix}
 wk_v'$, where $b\in \mathcal{O}_v/\mathfrak{p}_v^{r_v'}$, and $k_v'\in K_v[r_v']$. Since $0\leq i\leq r_v'$, then 
\begin{equation}\label{10.2}
\sum_{m\in\mathbb{Z}}\Phi_v^{\sharp}(\eta \varpi_v^mk_v\diag(1,\varpi_v^i))\eta_v(\varpi_v^{j+2m})q_v^{-2ms}=\eta_v(\varpi_v^j)L_v(s,\eta_v^2).
\end{equation}

By \eqref{8.6} in Lemma \ref{lem8.1} we have
\begin{equation}\label{10.3}
\Big|W_v'\left(\begin{pmatrix}
\varpi_v^j\\
&1
\end{pmatrix}k_v\right)\Big|^2=\Big|W_v'\left(\begin{pmatrix}
\varpi_v^j\\
&1
\end{pmatrix}w\right)\Big|^2=|W_v'(I_2)|^2\mathbf{1}_{j=-r_v'}.	
\end{equation}

Therefore, \eqref{equ10.2} follows from 
\eqref{10.2} and \eqref{10.3}. 
\end{proof}

\subsubsection{Calculation of $\Psi_v^{\mathrm{small}}(s;i)$}

\begin{lemma}\label{lem10.2}
Let notation be as above. Then 
\begin{align*}
\frac{\Psi_v^{\mathrm{small}}(s;i)}{|W_v'(I_2)|^2\Vol(K_v[r_v'])}=&\frac{q_v^{2\min\{0,r_v'-i\}s}}{\eta_v(\varpi_v^{2\min\{0,r_v'-i\}})}+\zeta_v(1)\sum_{l=1}^{r_v'-1}\frac{q_v^{2\min\{0,l-i\}s}}{\eta_v(\varpi_v^{2\min\{0,l-i\}})}\\
&\sum_{j\in \mathbb{Z}}\sum_{\substack{\chi_v\in \widehat{\mathcal{O}_v^{\times}/(1+\mathfrak{p}_v^{r_v'-l})}\\
r_{\chi_v}=r_v'-l,\ 
r_{\pi_v'\otimes\chi_v}=r_v'-j}}\frac{\eta_v(\varpi_v^j)}{q_v^{j(s-1)}}.
\end{align*}
\end{lemma}
\begin{proof}
Suppose that $k=\begin{pmatrix}
1\\
\varpi_v^l\beta_v& 1
\end{pmatrix}k_v'$, where $1\leq l\leq r_v'$, $\beta_v\in \mathcal{O}_v$, and $k_v'\in K_v[r_v']$. Then 
\begin{align*}
\sum_{m\in\mathbb{Z}}\frac{\Phi_v^{\sharp}(\eta \varpi_v^mk_v\diag(1,\varpi_v^i))\eta_v(\varpi_v^{j+2m})}{q_v^{2ms}L_v(s,\eta_v^2)}=\sum_{m\in\mathbb{Z}}\frac{\Phi_v^{\sharp}((\varpi_v^{m+l}\beta_v,\varpi_v^{m+i}))\eta_v(\varpi_v^{j+2m})}{q_v^{2ms}L_v(s,\eta_v^2)},
\end{align*}
which is equal to 
\begin{align*}
\sum_{m\geq -\min\{i,l\}}\eta_v(\varpi_v^{j+2m})q_v^{-2ms}L_v(s,\eta_v^2)^{-1}=\eta_v(\varpi_v^{j})\cdot \eta_v(\varpi_v^{-2\min\{i,l\}})q_v^{2\min\{i,l\}s}.
\end{align*}

Therefore, $\Psi_v^{\mathrm{small}}(s;i)=\Psi_v^{(1)}(s;i)+\Psi_v^{(2)}(s;i)$,  
where $\Psi_v^{(1)}(s;i)$ is defined by 
\begin{align*}
\sum_{j\in \mathbb{Z}}q_v^{j(1-s)}\int_{K_v[r_v']}\Big|W_v'\left(\begin{pmatrix}
\varpi_v^j\\
&1
\end{pmatrix}k_v\right)\Big|^2\cdot \eta_v(\varpi_v^{j-2\min\{0,r_v'-i\}})q_v^{2\min\{0,r_v'-i\}s}dk_v,
\end{align*}
and $\Psi_v^{(2)}(s;i)$ is defined by 
\begin{align*}
\Vol(K_v[r_v'])\sum_{j\in \mathbb{Z}}\frac{\eta_v(\varpi_v^j)}{q_v^{j(s-1)}}\sum_{l} \frac{\eta_v(\varpi_v^{-2\min\{0,l-i\}})}{q_v^{-2\min\{0,l-i\}s}}\sum_{\beta_v}\Big|W_v'\left(\begin{pmatrix}
\varpi_v^j\\
\varpi_v^l\beta_v& 1
\end{pmatrix}\right)\Big|^2.
\end{align*}
Here $1\leq l<r_v'$, and $\beta_v\in \varpi_v^l\mathcal{O}_v^{\times}/(1+\mathfrak{p}_v^{r_v'})$. 
\begin{itemize}
\item Since $W_v'$ is a local new vector, then 
\begin{align*}
\Psi_v^{(1)}(s;i)=\eta_v(\varpi_v^{-2\min\{0,r_v'-i\}})q_v^{2\min\{0,r_v'-i\}s}\Vol(K_v[r_v'])|W_v'(I_2)|^2.
\end{align*}

\item For $1\leq l<r_v'$,  by Lemma \ref{lem8.2},
\begin{align*}
\Psi_v^{(2)}(s;i)=\frac{|W_v'(I_2)|^2\Vol(K_v[r_v'])}{\zeta_v(1)^{-2}q_v^{r_v'-l}}\sum_{j\in \mathbb{Z}}\frac{\eta_v(\varpi_v^j)}{q_v^{j(s-1)}}\sum_{l=1}^{r_v'-1}\frac{\eta_v(\varpi_v^{-2\min\{0,l-i\}})}{q_v^{-2\min\{0,l-i\}s}}\cdot\mathcal{I}^{\sharp}(l),
\end{align*}
where $\mathcal{I}^{\sharp}(l)$ is defined by 
\begin{align*}
\sum_{\beta_v\in \varpi_v^l\mathcal{O}_v^{\times}/(1+\mathfrak{p}_v^{r_v'})}\bigg|\sum_{\substack{\chi_v\in \widehat{\mathcal{O}_v^{\times}/(1+\mathfrak{p}_v^{r_v'-l})}\\
r_{\chi_v}=r_v'-l\\ 
r_{\pi_v'\otimes\chi_v}=r_v'-j}}\chi_v(-\beta_v)\varepsilon(1/2,\pi_v'\otimes\omega_v'^{-1}\chi_v^{-1},\psi_v)\varepsilon(1/2,\chi_v)\bigg|^2.
\end{align*}
Similar to the calculation of \eqref{eq9.38}, opening the square and swap sums, along with the orthogonality and the fact that $|\varepsilon(1/2,\pi_v'\otimes\omega_v'^{-1}\chi_v^{-1},\psi_v)|=|\varepsilon(1/2,\chi_v)|=1$, we obtain 
\begin{align*}
\mathcal{I}^{\sharp}(l)=\zeta_v(1)^{-1}q_v^{r_v'-l}\sum_{\substack{\chi_v\in \widehat{\mathcal{O}_v^{\times}/(1+\mathfrak{p}_v^{r_v'-l})}\\
r_{\chi_v}=r_v'-l,\ 
r_{\pi_v'\otimes\chi_v}=r_v'-j}}1.
\end{align*}

Therefore, $\Psi_v^{(2)}(s;i)$ is equal to 
\begin{align*}
\frac{|W_v'(I_2)|^2\Vol(K_v[r_v'])}{\zeta_v(1)^{-1}}\sum_{l=1}^{r_v'-1}\frac{\eta_v(\varpi_v^{-2\min\{0,l-i\}})}{q_v^{-2\min\{0,l-i\}s}}\sum_{j\in \mathbb{Z}}\sum_{\substack{\chi_v\in \widehat{\mathcal{O}_v^{\times}/(1+\mathfrak{p}_v^{r_v'-l})}\\
r_{\chi_v}=r_v'-l,\ 
r_{\pi_v'\otimes\chi_v}=r_v'-j}}\frac{\eta_v(\varpi_v^j)}{q_v^{j(s-1)}}.
\end{align*}
\end{itemize}

Lemma \ref{lem10.2} follows from the above calculations of $\Psi_v^{(1)}(s;i)$ and $\Psi_v^{(2)}(s;i)$.
\end{proof}

\begin{cor}\label{cor10.3}
Let notation be as above. Then 
\begin{align*}
\frac{\Psi_v^{\mathrm{small}}(s;i)}{|W_v'(I_2)|^2\Vol(K_v[r_v'])}=&\frac{q_v^{2\min\{0,r_v'-i\}s}}{\eta_v(\varpi_v^{2\min\{0,r_v'-i\}})}+\sum_{l=1}^{r_v'-1}\frac{\zeta_v(1)\eta_v(\varpi_v^{\min\{0,2l-r_v'\}})}{q_v^{(s-1)\min\{0,2l-r_v'\}}}\\
&\frac{q_v^{2\min\{0,l-i\}s}}{\eta_v(\varpi_v^{2\min\{0,l-i\}})}\sum_{\substack{\chi_v\in \widehat{\mathcal{O}_v^{\times}/(1+\mathfrak{p}_v^{r_v'-l})}\\
r_{\chi_v}=r_v'-l,\ 
r_{\pi_v'\otimes\chi_v}=r_v'-\min\{0,2l-r_v'\}}}1.
\end{align*}	
In particular, for $0\leq i\leq r_v'$, we have the uniform bound 
\begin{equation}\label{10.5}
\Psi_v^{\mathrm{small}}(0;i)\ll r_v'\cdot q_v^{-\frac{r_v'}{2}}|W_v'(I_2)|^2,
\end{equation}
where the implied constant is absolute. 
\end{cor}
\begin{proof}
Since replacing $W_v'(x_v)$ with $W_v'(x_v)\chi_v(\det x_v)$ for any  unitary character 
$\chi_v$ of $F_v^{\times}$ does not affect  $\Psi_v^{\mathrm{small}}(s;i)$, we may assume $\pi_v'$ is twist-minimal. Arguing as in the proof of Lemma \ref{lem9.14}, we obtain from Lemma \ref{lem10.2} that 
\begin{align*}
\Psi_v^{(2)}(s;i)=&\zeta_v(1)\sum_{l=1}^{r_v'-1}\eta_v(\varpi_v^{\min\{0,2l-r_v'\}})|W_v'(I_2)|^2\Vol(K_v[r_v'])q_v^{(1-s)\min\{0,2l-r_v'\}}\\
&\eta_v(\varpi_v^{-2\min\{0,l-i\}})q_v^{2\min\{0,l-i\}s}\sum_{\substack{\chi_v\in \widehat{\mathcal{O}_v^{\times}/(1+\mathfrak{p}_v^{r_v'-l})}\\
r_{\chi_v}=r_v'-l,\ 
r_{\pi_v'\otimes\chi_v}=r_v'-\min\{0,2l-r_v'\}}}1.
\end{align*}

As a consequence, we obtain
\begin{align*}
\Psi_v^{(2)}(0;i)\ll &\sum_{l=1}^{r_v'-1}q_v^{\min\{0,2l-r_v'\}-l}\cdot |W_v'(I_2)|^2\ll r_v'\cdot q_v^{-\frac{r_v'}{2}}\cdot |W_v'(I_2)|^2.
\end{align*}
In conjunction with the trivial bound $\Psi_v^{(1)}(0;i)\ll q_v^{-r_v'}|W_v'(I_2)|^2$, the estimate 
 \eqref{10.5} follows. 
\end{proof}

\subsubsection{Proof of Proposition \ref{prop10.1}}
By \eqref{equ10.1} and meromorphic continuation of each $\Psi_v^{\sharp}(s;i)$, we obtain 
\begin{equation}\label{10.7}
\mathcal{I}_{i}=\mathcal{I}_{i}(0)=\zeta_v(1)\langle W_v^{\circ}, W_v^{\circ}\rangle^{-1}|W_v'(I_2)|^{-4}\cdot |\Psi_v^{\sharp}(0;i)|^2.
\end{equation}
Then Proposition \ref{prop10.1} follows from \eqref{equ10.2} in Lemm \ref{lemma10.2}, \eqref{10.5} in Corollary \ref{cor10.3}, and the above identity \eqref{10.7}. 

\subsection{$\eta_v$ is Ramified}
In this section we aim to prove an upper bound for $\mathcal{I}_i$ in the case that $r_{\eta_v}\geq 1$, i.e., $\eta_v$ is ramified. The main result is the following.  
\begin{prop}\label{prop10.5}
Let notation be as before. Suppose $\eta_v$ is ramified. Let $0\leq i\leq r_v'-r_v$. Then
\begin{equation}\label{eq10.10}
\mathcal{I}_i\ll q_v^{-r_v-i}, 
\end{equation}
where the implied constant is absolute. 
\end{prop}
 
\subsubsection{Sections}
For $n\in \mathbb{Z}_{\geq 1}$, we define 
\begin{align*}
K_v^{\circ}[n]:=\Big\{\begin{pmatrix}
k_{11}& k_{12}\\
k_{21}& k_{22}	
\end{pmatrix}\in K_v:\ k_{21}\in \varpi_v^n\mathcal{O}_v,\ k_{22}\in 1+\mathfrak{p}_v^n
\Big\}.
\end{align*}

Recall that $r_v=2r_{\eta_v}$. By \cite[Proposition 2.1.2]{Sch02}, the function 
\begin{align*}
h_v^{\flat}(g_v)=\begin{cases}
\eta_v(\varpi_v^{-r_{\eta_v}}a_vd_v^{-1})|a_vd_v^{-1}|_v^{1/2}, & \text{if $g_v\in \begin{pmatrix}
a_v& *\\
&d_v	
\end{pmatrix}\begin{pmatrix}
1&\\
\varpi_v^{r_{\eta_v}}&1
\end{pmatrix}K_v^{\circ}[r_v],$}\\
0, & \text{if $g_v\notin B(F_v)\begin{pmatrix}
1&\\
\varpi_v^{r_{\eta_v}}&1
\end{pmatrix}K_v^{\circ}[r_v]$.}
\end{cases}
\end{align*}
is a nonzero right-$K_v^{\circ}[r_{v}]$-invariant section in $\sigma_v$. The factor $\eta_v(\varpi_v^{-r_{\eta_v}})$ in the above definition makes $h_v^{\flat}$ independent of the choice of $\varpi_v$. 

Moreover, since $\eta_v$ is ramified, then $h_v^{\flat}$ is also a section when $\sigma_v=\mathrm{St}\otimes\eta_v$, where $\mathrm{St}$ is the Steinberg representation.

%$h_v^{\flat}(I_2)=\eta_v(\varpi_v^{-r_{\eta_v}})$.

\subsubsection{Auxiliary period integral}
Define the Rankin--Selberg period
\begin{equation}\label{eq9.10}
\Psi_v^{\flat}(i):=\int_{N(\mathbb{A}_F)\backslash\overline{G}(\mathbb{A}_F)}W_v'(x_v)\overline{W_v'(x_v)}h_v^{\flat}(x_v\diag(1,\varpi_v^i))dx_v.
\end{equation}

Let $\mathcal{I}_i$ be defined by \eqref{eq9.28.}. By \cite[Proposition 5.1]{Hsi21} we have 
\begin{equation}\label{10.8}
\mathcal{I}_i=\zeta_v(1)\langle h_v^{\flat}, h_v^{\flat}\rangle^{-1}|W_v'(I_2)|^{-4}\cdot |\Psi_v^{\flat}(i)|^2.
\end{equation}

Making use of Iwasawa decomposition we obtain  
\begin{align*}
\Psi_v^{\flat}(i)=\sum_{j\in \mathbb{Z}}q_v^j\int_{K_v}\Big|W_v'\left(\begin{pmatrix}
\varpi_v^j\\
&1
\end{pmatrix}k_v\begin{pmatrix}
\varpi_v^i\\
&1
\end{pmatrix}\right)\Big|^2h_v^{\flat}\left(\begin{pmatrix}
\varpi_v^j\\
&1
\end{pmatrix}k_v\right)dk_v.
\end{align*}

\begin{lemma}\label{lem10.6}
Let the notation be as above, and assume $0\le i\le r_v'-r_v$ and $r_{\eta_v}\ge 1$. Denote by $H(\mathcal{O}_v)=B(\mathcal{O}_v)\cap \begin{pmatrix}
1\\
\varpi_v^{r_{\eta_v}} & 1
\end{pmatrix}K_v^{\circ}[r_v]\begin{pmatrix}
1\\
-\varpi_v^{r_{\eta_v}} & 1
\end{pmatrix}$. Then 
\begin{multline}\label{10.12}
\Psi_v^{\flat}(i)=\eta_v(\varpi_v^{-r_{\eta_v}})\zeta_v(1)^2|W_v'(I_2)|^2\#(\mathcal{O}_v^{\times}/(1+\mathfrak{p}_v^{r_{\eta_v}}))^{-1}[B(\mathcal{O}_v):H(\mathcal{O}_v)]\\
\Vol(K_v[r_v'-i,i])
q_v^{i-r_{\eta_v}}\sum_{j\in \mathbb{Z}}\eta_v(\varpi_v^{j})q_v^{j/2}\sum_{\substack{\chi_v\in \widehat{\mathcal{O}_v^{\times}/(1+\mathfrak{p}_v^{r_v'-i-r_{\eta_v}})}\\
r_{\chi_v}=r_v'-i-r_{\eta_v}\\ 
r_{\pi_v'\otimes\chi_v}=r_v'-j-i}}\mathbf{1}_{\substack{r_{\chi_v\eta_v}=r_v'-i-r_{\eta_v}\\ r_{\pi_v'\otimes\chi_v\eta_v}=r_v'-j-i}}\\
\varepsilon(1/2,\pi_v'\otimes\omega_v'^{-1}\chi_v^{-1}\eta_v^{-1},\psi_v)\varepsilon(1/2,\chi_v\eta_v)\overline{\varepsilon(1/2,\pi_v'\otimes\omega_v'^{-1}\chi_v^{-1},\psi_v)\varepsilon(1/2,\chi_v)}.
\end{multline}
In particular, if $\pi_v'$ is twist-minimal, then 
\begin{equation}\label{fc9.12}
\Psi_v^{\flat}(i)\ll q_v^{-\frac{3r_{\eta_v}+i}{2}}|W_v'(I_2)|^2.
\end{equation}
\end{lemma}
\begin{proof}
By definition, $h_v^{\flat}\left(k_v\right)=0$ unless 
\begin{equation}\label{fc9.13}
k_v\in B(\mathcal{O}_v)\begin{pmatrix}
1&\\
\varpi_v^{r_{\eta_v}} &1
\end{pmatrix}K_v^{\circ}[r_v]=K_v[r_{\eta_v}].
\end{equation}

We obtain 
\begin{multline*}
\Psi_v^{\flat}(i)=\frac{\eta_v(\varpi_v^{-r_{\eta_v}})}{\Vol(H(\mathcal{O}_v))}\cdot \sum_{j\in \mathbb{Z}}\eta_v(\varpi_v^{j})q_v^{j/2}\int_{K_v^{\circ}[r_v]}\int_{B(\mathcal{O}_v)}h_v^{\flat}\left(b\begin{pmatrix}
1\\
\varpi_v^{r_{\eta_v}} & 1
\end{pmatrix}k_v'\right)\\
\Big|W_v'\left(\begin{pmatrix}
\varpi_v^j\\
&1
\end{pmatrix}b\begin{pmatrix}
1\\
\varpi_v^{r_{\eta_v}} & 1
\end{pmatrix}k_v\begin{pmatrix}
\varpi_v^i\\
&1
\end{pmatrix}\right)\Big|^2dbdk_v.
\end{multline*}

As a result, we obtain 
\begin{multline}\label{10.10}
\Psi_v^{\flat}(i)=\frac{\eta_v(\varpi_v^{-r_{\eta_v}})[B(\mathcal{O}_v):H(\mathcal{O}_v)]}{\#(\mathcal{O}_v^{\times}/(1+\mathfrak{p}_v^{r_{\eta_v}}))}\cdot \sum_{j\in \mathbb{Z}}\eta_v(\varpi_v^{j})q_v^{j/2}\int_{\mathcal{O}_v^{\times}}\overline{\eta}_v(\gamma)\int_{K_v[r_v]}\\
\Big|W_v'\left(\begin{pmatrix}
\varpi_v^j\\
\varpi_v^{r_{\eta_v}}\gamma & 1
\end{pmatrix}k_v\begin{pmatrix}
\varpi_v^i\\
&1
\end{pmatrix}\right)\Big|^2dk_vd^{\times}\gamma.
\end{multline}

Recall that $i\leq r_v'-r_v$. By Cramer's rule, we have the decomposition  
\begin{equation}\label{10.11}
K_v[r_v]=\bigsqcup_{\alpha\in \mathfrak{p}_v^{r_v}/\mathfrak{p}_v^{r_v'-i}}\bigsqcup_{\beta\in \mathcal{O}_v/\mathfrak{p}_v^i}\begin{pmatrix}
1&\\
\alpha & 1
\end{pmatrix}\begin{pmatrix}
1&\beta \\
& 1
\end{pmatrix}K_v[r_v'-i,i],
\end{equation}
where 
\begin{align*}
K_v[r_v'-i,i]:=\Big\{\begin{pmatrix}
k_{11}& k_{12}\\
k_{21} & k_{22}
\end{pmatrix}\in K_v:\ k_{21}\in \mathfrak{p}_v^{r_v'-i},\ k_{12}\in \mathfrak{p}_v^i\Big\}.
\end{align*}

Substituting \eqref{10.11} into \eqref{10.10} leads to 
\begin{multline}\label{eq10.12}
\Psi_v^{\flat}(i)=\eta_v(\varpi_v^{-r_{\eta_v}})\#(\mathcal{O}_v^{\times}/(1+\mathfrak{p}_v^{r_{\eta_v}}))^{-1}[B(\mathcal{O}_v):H(\mathcal{O}_v)]\\
\Vol(K_v[r_v'-i,i])\sum_{j\in \mathbb{Z}}\eta_v(\varpi_v^{j})q_v^{j/2}\sum_{\alpha\in \mathfrak{p}_v^{r_v}/\mathfrak{p}_v^{r_v'-i}}\sum_{\beta\in \mathcal{O}_v/\mathfrak{p}_v^i}\mathcal{J}^{\sharp},
\end{multline}
where 
\begin{align*}
\mathcal{J}^{\sharp}:=\int_{\mathcal{O}_v^{\times}}\overline{\eta}_v(\gamma)\Big|W_v'\left(\begin{pmatrix}
\varpi_v^j&\\
\varpi_v^{r_{\eta_v}}\gamma & 1
\end{pmatrix}\begin{pmatrix}
1&\\
\alpha & 1
\end{pmatrix}\begin{pmatrix}
1&\beta \\
& 1
\end{pmatrix}\begin{pmatrix}
\varpi_v^i\\
&1
\end{pmatrix}\right)\Big|^2d^{\times}\gamma.
\end{align*}

Let $\alpha'=\varpi_v^{r_{\eta_v}}\gamma+\alpha$. Using the identity 
\begin{align*}
\begin{pmatrix}
1&\\
\alpha' & 1
\end{pmatrix}\begin{pmatrix}
1&\beta \\
& 1
\end{pmatrix}=\begin{pmatrix}
(1+\alpha'\beta)^{-1}&\beta(1+\alpha'\beta)^{-1}\\
& 1
\end{pmatrix}\begin{pmatrix}
1&\\
\alpha' & 1+\alpha'\beta
\end{pmatrix}
\end{align*}
and the property of Whittaker functions we obtain 
\begin{align*}
\bigg|W_v'\left(\begin{pmatrix}
\varpi_v^j&\\
\gamma & 1
\end{pmatrix}\begin{pmatrix}
1&\\
\alpha & 1
\end{pmatrix}\begin{pmatrix}
1&\beta \\
& 1
\end{pmatrix}\begin{pmatrix}
\varpi_v^i\\
&1
\end{pmatrix}\right)\bigg|=\bigg|W_v'\left(\begin{pmatrix}
\varpi_v^{j+i}&\\
\varpi_v^i(\varpi_v^{r_{\eta_v}}\gamma+\alpha)& 1
\end{pmatrix}\right)\bigg|.
\end{align*}

Since $e_v(\alpha)\geq r_v=2r_{\eta_v}$ and $r_{\eta_v}\geq 1$, we may apply the change of variables $\gamma\mapsto \gamma-\varpi_v^{-r_{\eta_v}}\alpha$ to derive that  
\begin{align*}
\mathcal{J}^{\sharp}=\int_{\mathcal{O}_v^{\times}}\overline{\eta}_v(\gamma)\Big|W_v'\left(\begin{pmatrix}
\varpi_v^{j+i}&\\
\varpi_v^{i+r_{\eta_v}}\gamma & 1
\end{pmatrix}\right)\Big|^2d^{\times}\gamma.
\end{align*}

By Lemma \ref{lem8.2}, and orthogonality of characters (so that the cross terms when open the square vanish as in Lemma \ref{lem9.9}), we obtain 
\begin{multline}\label{eq10.13}
\mathcal{J}^{\sharp}=q_v^{i+r_{\eta_v}-r_v'}\zeta_v(1)^2|W_v'(I_2)|^2\int_{\mathcal{O}_v^{\times}}\overline{\eta}_v(\gamma)\\
\Bigg|\sum_{\substack{\chi_v\in \widehat{\mathcal{O}_v^{\times}/(1+\mathfrak{p}_v^{r_v'-i-r_{\eta_v}})}\\
r_{\chi_v}=r_v'-i-r_{\eta_v}\\ 
r_{\pi_v'\otimes\chi_v}=r_v'-j-i}}\chi_v(-\gamma)\varepsilon(1/2,\pi_v'\otimes\omega_v'^{-1}\chi_v^{-1},\psi_v)\varepsilon(1/2,\chi_v)\Bigg|^2d^{\times}\gamma.
\end{multline}

Opening the square and taking advantage of orthogonality we derive that 
\begin{multline}\label{10.13}
\mathcal{J}^{\sharp}=q_v^{i+r_{\eta_v}-r_v'}\zeta_v(1)^2|W_v'(I_2)|^2\sum_{\substack{\chi_v\in \widehat{\mathcal{O}_v^{\times}/(1+\mathfrak{p}_v^{r_v'-i-r_{\eta_v}})}\\
r_{\chi_v}=r_v'-i-r_{\eta_v}\\ 
r_{\pi_v'\otimes\chi_v}=r_v'-j-i}}\mathbf{1}_{\substack{r_{\chi_v\eta_v}=r_v'-i-r_{\eta_v}\\ r_{\pi_v'\otimes\chi_v\eta_v}=r_v'-j-i}}\\
\varepsilon(1/2,\pi_v'\otimes\omega_v'^{-1}\chi_v^{-1}\eta_v^{-1},\psi_v)\varepsilon(1/2,\chi_v\eta_v)\overline{\varepsilon(1/2,\pi_v'\otimes\omega_v'^{-1}\chi_v^{-1},\psi_v)\varepsilon(1/2,\chi_v)}.
\end{multline}

Therefore, substituting \eqref{eq10.13} and  \eqref{10.13} into \eqref{eq10.12} we  obtain \eqref{10.12}. If moreover $\pi_v'$ is twist-minimal, then 
\begin{align*}
r_{\pi_v'\otimes\chi_v}=r_v'-j-i\geq r_v',
\end{align*}
which forces $j\leq -i\leq 0$. Hence  \eqref{fc9.12} follows from \eqref{10.12}. 
\end{proof}

\subsubsection{Proof of Proposition \ref{prop10.5}}
By \eqref{fc9.13} we obtain 
\begin{equation}\label{9.19}
\langle h_v^{\flat}, h_v^{\flat}\rangle\gg \int_{K_v[r_{\eta_v}]}dk_v\gg q_v^{-r_{\eta_v}}.
\end{equation}

Observe that in the definition of \eqref{eq9.10} we may replace the function 
$W_v'(x_v)$ with $W_v'(x_v)\chi_v(\det x_v)$ for any unitary character 
$\chi_v$ of $F_v^{\times}$. 
Hence, without loss of generality, we may assume that $\pi_v'$ is twist-minimal 
in the definition of $\mathcal{I}_i$. 
Therefore, \eqref{eq10.10} follows from \eqref{10.8}, \eqref{fc9.12}, and \eqref{9.19}.

\section{Archimedean Integrals on the Dual Side \RNum{1}}\label{sec11} 

In this section we analyze the local integrals on the dual side associated with the 
cusp forms $\phi_j$, $j=1,2$, at real places where the $\mathrm{GL}_2$-spectrum lies 
in the principal series.  
We obtain uniform upper bounds for the resulting integral transforms.  
The remaining Archimedean cases are treated in \textsection\ref{sect11} and \textsection\ref{sect12}.

\subsection{Representations of $\mathrm{PGL}_2(\mathbb{R})$}
Let $F_v\simeq \mathbb{R}$ and  $K_v=\mathrm{SO}_2(\mathbb{R})$. The irreducible representations of $K_v$ are parametrized by nonnegative integers. Denote by such representations $\tau_n$, $n\in\mathbb{Z}_{n\geq 0}$. Then the restriction of $\sigma_v$ to $K_v$ decomposes into a direct sum of $\tau_n$'s. Let $n_{\sigma_v}$ be the smallest $n\geq 0$ such that $0\neq \tau_n\subset \sigma_v|_{K_v}$. 

Then $\sigma_v$ is isomorphic either to a principal series representation $\Ind\mu\otimes\mu^{-1}$ for some characters $\mu=|\cdot|_v^{\nu}\cdot \sgn^{\delta}$, $\nu\in i\mathbb{R}_+\sqcup (-\vartheta, \vartheta)$ and $\delta\in \{0,1\}$; or a discrete series $\sigma(\mu,\mu^{-1})$, which is the infinite dimensional factor in the composition series of the reducible representation $\Ind\mu\otimes\mu^{-1}$, where $\mu=|\cdot|_v^{\nu}\cdot \sgn^{\delta}$, with $2\nu$ being a positive odd integer. Explicitly we have 
\begin{align*}
n_{\sigma_v}=\begin{cases}
0,\ & \text{if $\sigma_v=\Ind\mu\otimes\mu^{-1}$ is a principal series}\\
2\nu+1,\ & \text{if $\sigma_v=\sigma(\mu,\mu^{-1})$ is a discrete series}.
\end{cases}
\end{align*}

Let $n\in 2\mathbb{Z}$ and $\Re(s)>1/2$. Consider the Godement section in the induced representation $\Ind\mu|\cdot |_v^{s}\otimes\mu^{-1}|\cdot |_v^{-s}$: 
\begin{align*}
f_n^{\sharp}(g_v,s):=\mu(\det g_v)|\det g_v|_v^{s+1/2}\int_{F_v^{\times}}\Phi_v^{\sharp}((0,t)g_v)\mu^2(t)|t|^{2s+1}d^{\times}t,
\end{align*}
where $g_v\in \mathrm{GL}_2(F_v)$, and $\Phi_v^{\sharp}(t_1,t_2)=e^{-\pi(t_1^2+t_2^2)}(t_2-i\sgn(n)t_1)^{|n|}$. Then the function  
\begin{align*}
f_n(g_v,s):=f_n^{\sharp}(I_2,s)^{-1}f_n^{\sharp}(g_v,s)
\end{align*} 
admits an analytic continuation to $s\in \mathbb{C}$. Let $f_n=f_n(\cdot,0)$ if $n\in 2\mathbb{Z}$ and $f_n=0$ if $n$ is odd. Then $f_n$ is the section in $\Ind\mu\otimes\mu^{-1}$ such that 
\begin{equation}\label{equ11.1}
f_n\left(\begin{pmatrix}
a_v\\
&1
\end{pmatrix}\begin{pmatrix}
\cos\theta &\sin\theta\\
-\sin\theta & \cos\theta
\end{pmatrix}\right)=|a_v|_v^{1/2+\nu}\sgn(a_v)^{\delta}e^{in\theta}\mathbf{1}_{n\in 2\mathbb{Z}}
\end{equation}
for all $a_v\in F_v^{\times}$ and $\theta\in [0,2\pi)$. The corresponding Whittaker vector associated with $f_n$ is given by  
\begin{equation}\label{eq11.1}
W_n(g_v):=\int_{N(F_v)}f_n(wu_vg_v,0)\overline{\theta_v(u_v)}du_v.
\end{equation}

\begin{comment}

\bigskip
\bigskip
\bigskip
\bigskip
\bigskip
\bigskip
\bigskip

\begin{align*}
|a|_v^{\frac{1}{2}+\nu}\int_{N(F_v)}\int_{F_v^{\times}}\Phi_v^{\sharp}((ta,tu))|t|^{1+2\nu}d^{\times}t\psi_v(-u)du
\end{align*}

\begin{align*}
\Phi_v^{\sharp}((ta,tu))=t^{|n|}e^{-\pi t^2(a^2+u^2)}(u-i\sgn(n)a)^{|n|}
\end{align*}

\begin{align*}
&|a|_v^{\frac{1}{2}+\nu}\int_{N(F_v)}\int_{F_v^{\times}}|t|^{1+2\nu+n}e^{-\pi t^2(a^2+u^2)}(u-ia)^{n}d^{\times}t\psi_v(-u)du\\
=&2|a|_v^{\frac{1}{2}+\nu}\int_{N(F_v)}\int_0^{\infty}t^{m+n}e^{-\pi t^2(a^2+u^2)}(u-ia)^{n}d^{\times}t\psi_v(-u)du\\
=&2|a|_v^{\frac{m}{2}}\int_0^{\infty}t^{m+n}e^{-\pi t^2}d^{\times}t\int_{\mathbb{R}}\frac{(u-ia)^{n}e^{-2\pi iu}}{(a^2+u^2)^{\frac{m+n}{2}}}du\\
\end{align*}

\bigskip
\bigskip
\bigskip
\bigskip
\bigskip
\bigskip
\bigskip

\begin{align*}
\Phi_v^{\sharp}((0,t))=t^{n}e^{-\pi t^2}
\end{align*}

\begin{align*}
f_n^{\sharp}(\begin{pmatrix}
a\\
&1	
\end{pmatrix}
,s):=|a|_v^{\frac{1}{2}+\nu}\int_{F_v^{\times}}e^{-\pi t^2}|t|^{1+2\nu+n}d^{\times}t=2|a|_v^{\frac{m}{2}}\int_0^{\infty}e^{-\pi t^2}t^{m+n}d^{\times}t
\end{align*}

\bigskip
\bigskip
\bigskip
\bigskip
\bigskip
\bigskip
\bigskip
\end{comment}

This also holds for the discrete series case. A brute force calculation yields that 
\begin{equation}\label{whi11.2}
W_n(\diag(a,1))=|a|_v^{1/2-\nu}\sgn(a)^{\delta}\int_{F_v}\frac{(u-i)^ne^{-2\pi iua}}{(u^2+1)^{\nu+\frac{n+1}{2}}}du,
\end{equation}
which converges (at least conditionally) when $\Re(\nu)>-1/2$. 
\begin{comment}
W_n(\diag(a,1))=|a|^{\nu+1/2}\sgn(a)^{\delta}\int_{F_v}\frac{(u-ia)^ne^{-2\pi iu}}{(u^2+a^2)^{\nu+\frac{n+1}{2}}}du.
\end{comment}

In particular, when $n=n_{\sigma_v}$, $W_n$ is the lowest weight vector in the Whittaker model of $\sigma_v$, satisfying the following: 
\begin{itemize}
\item if $\sigma_v$ is not a discrete series, then 
\begin{align*}
W_{n_{\sigma_v}}(\diag(a,1))\propto \sgn(a)^{\delta}|a|_v^{1/2}K_{\nu}(2\pi |a|).
\end{align*}
\item if $\sigma_v$ is a discrete series, then 
\begin{align*}
W_{n_{\sigma_v}}(\diag(a,1))\propto  (1+\sgn(a))|a|_v^{n_{\sigma_v}/2}e^{-2\pi a}.
\end{align*}
\end{itemize}

Therefore, 
\begin{equation}\label{eq10.4}
\{\langle W_n,W_n\rangle^{-1/2}W_n:\ n\equiv n_{\sigma_v}\pmod{2}\}
\end{equation}
forms an orthonormal basis of $\sigma_v$ if it is not a discrete series; and 
\begin{equation}\label{eq10.5}
\{\langle W_n,W_n\rangle^{-1/2}W_n:\ |n|\geq n_{\sigma_v},\ n\equiv n_{\sigma_v}\pmod{2}\}
\end{equation}
forms an orthonormal basis of $\sigma_v$ if it is a discrete series. 

In particular, when $\sigma_v$ is a unitary principal series, by \cite[\textsection 3.1.6--\textsection 3.1.10]{MV10}, 
\begin{align*}
\langle W_n,W_n\rangle=\frac{\zeta_v(1)}{\zeta_v(2)}\int_{K_v}|f_n(k_v)|^2dk_v=\frac{\zeta_v(1)\mathbf{1}_{n\in 2\mathbb{Z}}}{\zeta_v(2)}=\pi\cdot \mathbf{1}_{n\in 2\mathbb{Z}}.
\end{align*}

\subsection{Triple Integrals: Setup}
Let $W_n$ be defined as in \eqref{eq11.1}. Let $W_{\pi_v'}$ be the vector in the Kirillov model of $\pi_v'$ as defined in \textsection\ref{sec4.3}. Define 
\begin{equation}\label{f11.3}
|\mathcal{P}_v(W_n,W_{\pi_v'},\overline{W_{\pi_v'}})|^2=\int_{\overline{G}(F_v)}\frac{\langle\sigma_v(g_v)W_n,W_n\rangle\cdot |\langle\pi_v'(g_v)W_{\pi_v'},W_{\pi_v'}\rangle|^2}{\langle W_n,W_n\rangle\cdot |\langle W_{\pi_v'},W_{\pi_v'}\rangle|^2} dg_v.  
\end{equation}

If $\sigma_v$ is a principal series representation $\Ind\mu\otimes\mu^{-1}$ for some characters $\mu=|\cdot|_v^{\nu}\cdot \sgn^{\delta}$, $\nu\in i\mathbb{R}_+\sqcup (-\vartheta, \vartheta)$ and $\delta\in \{0,1\}$, it follows from \cite[Lemma 3.4.2 and the accompanying Remark]{MV10} that $|\mathcal{P}_v(W_n,W_{\pi_v'},\overline{W_{\pi_v'}})|=|\Psi_v(f_n,W_{\pi_v'},\overline{W_{\pi_v'}})|$, where 
\begin{align*}
|\Psi_v(f_n,W_{\pi_v'},\overline{W_{\pi_v'}})|=\frac{1}{|\langle W_{\pi_v'},W_{\pi_v'}\rangle|^2}\bigg|\int_{N(F_v)\backslash G(F_v)}|W_{\pi_v'}(g_v)|^2f_n(g_v)dg_v\bigg|.
\end{align*}

\begin{lemma}\label{lemma11.1}
Let the notation remain as before. Suppose $0<\varepsilon<10^{-1}$, and let $\beta:(-\pi,\pi]\to\mathbb{R}$ be a function with compact support in $\{|\theta|<\pi:\ \varepsilon<|\theta|<\pi-\varepsilon\}$, and $\beta(-\pi/2)=\beta(\pi/2)=1$. Then
\begin{equation}\label{equ11.2}
\Psi_v(f_n,W_{\pi_v'},\overline{W_{\pi_v'}})=|\langle W_{\pi_v'},W_{\pi_v'}\rangle|^{-2}\cdot\Big[\Psi_n^{(1)}(W_{\pi_v'};\beta_1)+\Psi_n^{(2)}(W_{\pi_v'};\beta_2)\Big],
\end{equation}
where 
\begin{multline*}
\Psi_n^{(1)}(W_{\pi_v'};\beta_1):=\frac{1}{2\pi}\int_{F_v^{\times}}\int_{-\pi}^{\pi}\bigg|W_{\pi_v'}\left(\begin{pmatrix}
a_v\\
&1
\end{pmatrix}\begin{pmatrix}
1 &\\
c & 1
\end{pmatrix}\right)\bigg|^2\\
e^{in\theta}\beta_1(\theta)(1+c^2)^{\frac{1}{2}-\nu}d\theta |a_v|_v^{-\frac{1}{2}+\nu}\sgn(a_v)^{\delta}d^{\times}a_v,
\end{multline*} 
with $c=-\tan\theta\in \mathbb{R}$ and $\beta_1(\theta)=1-\beta(\theta)$; and 
\begin{multline*}
\Psi_n^{(2)}(W_{\pi_v'};\beta_2):=\frac{1}{2\pi}\int_{F_v^{\times}}\int_{-\pi}^{\pi}\bigg|W_{\pi_v'}\left(\begin{pmatrix}
a_v\\
&1
\end{pmatrix}\begin{pmatrix}
& 1\\
1& c'
\end{pmatrix}\right)\bigg|^2\\
e^{in\theta}\beta_2(\theta)(1+c'^2)^{\frac{1}{2}-\nu}d\theta |a_v|_v^{-\frac{1}{2}+\nu}\sgn(a_v)^{\delta}d^{\times}a_v,
\end{multline*}
with $c'=-\cot\theta\in \mathbb{R}$ and $\beta_2(\theta)=\beta(\theta)$. 
\end{lemma}
\begin{proof}
For simplicity we denote by 
\begin{align*}
\Psi_n(W_{\pi_v'}):=\int_{N(F_v)\backslash G(F_v)}|W_{\pi_v'}(g_v)|^2f_n(g_v)dg_v.
\end{align*}	
	
By Iwasawa decomposition, in conjunction with $\beta_1+\beta_2=1$, 
\begin{align*}
\Psi_n(W_{\pi_v'})=\int_{F_v^{\times}}\int_{-\pi}^{\pi}\big|W_{\pi_v'}\left(g_v\right)\big|^2f_n\left(g_v\right)(\beta_1(\theta)+\beta_2(\theta))d\theta |a_v|_v^{-1}d^{\times}a_v,
\end{align*}
where $g_v:=\begin{pmatrix}
a_v\\
&1
\end{pmatrix}\begin{pmatrix}
\cos\theta & \sin\theta\\
-\sin\theta & \cos\theta
\end{pmatrix}$.	

Let $\theta\in (-\pi,\pi)$. We have the following scenarios. 
\begin{itemize}
\item Suppose $\theta\neq \pm\pi/2$. Let $c=-\tan\theta\in \mathbb{R}$. We have 
\begin{equation}\label{f11.4}
\begin{pmatrix}
\cos\theta & \sin\theta\\
-\sin\theta & \cos\theta
\end{pmatrix}=\sgn(\cos\theta)\begin{pmatrix}
1 & -c\\
& 1\end{pmatrix}\begin{pmatrix}
\sqrt{1+c^2}& \\
& \frac{1}{\sqrt{1+c^2}}\end{pmatrix}\begin{pmatrix}
1\\
c & 1
\end{pmatrix}.
\end{equation}

\item Suppose $\theta\neq 0, \pm\pi$. Let $c'=-\cot\theta\in \mathbb{R}$. We have  
\begin{equation}\label{f11.5}
\begin{pmatrix}
\cos\theta & \sin\theta\\
-\sin\theta & \cos\theta
\end{pmatrix}=\sgn(\sin\theta)\begin{pmatrix}
1& c'\\
&1
\end{pmatrix}\begin{pmatrix}
\sqrt{1+c'^2}& \\
&-\frac{1}{\sqrt{1+c'^2}}
\end{pmatrix}\begin{pmatrix}
& 1\\
1& c'
\end{pmatrix}.
\end{equation}
\end{itemize}

Consequently, for $g_v=\begin{pmatrix}
a_v\\
&1
\end{pmatrix}\begin{pmatrix}
\cos\theta & \sin\theta\\
-\sin\theta & \cos\theta
\end{pmatrix}$, $W_{\pi_v'}\left(g_v\right)f_n\left(g_v\right)$ is equal to  
\begin{equation}\label{equat11.4}
\omega_v'^{-1}(\sgn(\cos\theta)\sqrt{1+c^2})\psi_v(-a_vc)W_{\pi_v'}\left(\begin{pmatrix}
a_v(1+c^2)\\
&1
\end{pmatrix}\begin{pmatrix}
1 & \\
c & 1
\end{pmatrix}\right)
\end{equation}
if $\theta\neq \pm\pi/2$; and $W_{\pi_v'}\left(\begin{pmatrix}
a_v\\
&1
\end{pmatrix}\begin{pmatrix}
\cos\theta & \sin\theta\\
-\sin\theta & \cos\theta
\end{pmatrix}\right)$ is equal to
\begin{equation}\label{equa11.3}
\omega_v'^{-1}(-\sgn(\sin\theta)\sqrt{1+c'^2})\psi_v(a_vc')W_{\pi_v'}\left(\begin{pmatrix}
a_v(1+c'^2)\\
&1
\end{pmatrix}\begin{pmatrix}
& 1\\
1& c'
\end{pmatrix}\right)
\end{equation}
if $\theta\neq 0, \pm\pi$. Therefore, \eqref{equ11.2} follows from \eqref{equa11.3} and \eqref{equat11.4}. 
\end{proof}

\begin{comment}
\marginpar{This is the case $\Re(\nu)\neq 0$}
We need to handle 
\begin{align*}
\int_0^{\infty}J_m(ay_1)J_m(ay_2)a^{m-1}da
\end{align*}
\end{comment}

\subsection{Upper Bounds for the Local Integrals}\label{sec11.3}
Suppose $\sigma_v$ is a principal series, which may be reducible or non-tempered. Specifically, we assume that $\sigma_v$ is isomorphic to a principal series representation $\Ind\mu\otimes\mu^{-1}$ for some characters $\mu=|\cdot|_v^{\nu}\cdot \sgn^{\delta}$, $\nu\in i\mathbb{R}_+\sqcup (-\vartheta, \vartheta)$ and $\delta\in \{0,1\}$. 

In this subsection, we establish an upper bound for $\Psi_v(f_n,W_{\pi_v'},\overline{W_{\pi_v'}})$ when $\sigma_v$ is a principal series. The main result is the following.
\begin{prop}\label{prop11.4}
Let notation be as before. Let $n\in \mathbb{Z}$. Suppose $\sigma_v$ is a principal series. Let $l\geq 0$. Then 
\begin{align*}
\Psi_v(f_n,W_{\pi_v'},\overline{W_{\pi_v'}})\ll C_v(\pi')^{1+\varepsilon}\cdot \min \bigg\{\frac{C_v(\pi')^l}{(1+|\nu|)^l},\frac{(1+|\nu|)^m}{|n|^m}\bigg\},
\end{align*}
where the implied constant depends only on $\varepsilon$, $l$, $m$, and the function $\varphi_v$ (see \textsection\ref{sec4.3}).
\end{prop}
\begin{proof}
The above upper bound for $\Psi_v(f_n,W_{\pi_v'},\overline{W_{\pi_v'}})$ follows readily from Lemmas \ref{lemma11.1}, \ref{lem11.2}, \ref{lem11.4}, \ref{lem11.6} and \ref{lem11.9}, together with the assumption that $\langle W_{\pi_v'},W_{\pi_v'}\rangle\asymp 1$ (see \textsection\ref{sec4.3}). 
\end{proof}

Let $\beta:(-\pi,\pi]\to [0,1]$ be a smooth function with compact support in $J:=\{|\theta|<\pi:\ \pi/5<|\theta|<4\pi/5\}$ such that $\beta(\theta)\equiv 1$ for all $\pi/4<|\theta|<3\pi/4$. Let $\beta_2=\beta$ and $\beta_1=1-\beta$ be as in the statement of Lemma \ref{lemma11.1}. 

\begin{lemma}\label{lem11.2}
Let notation be as before. Let $\varepsilon>0$. Let $n\in\mathbb{Z}$, and $l\geq 0$. Let $\beta$ be defined as above. Suppose $\sigma_v$ is a principal series. Then 
\begin{equation}\label{eq11.3}
\Psi_n^{(2)}(W_{\pi_v'};\beta_2)\ll \frac{C_v(\pi')^{-\frac{1}{2}+\Re(\nu)+\varepsilon}}{(1+|\nu|)^l},
\end{equation}
where the implied constant depends only on $\varepsilon$, $l$, and the function $\varphi_v$ (see \textsection\ref{sec4.3}). 
\end{lemma}
\begin{proof}
Let $c'=-\cot\theta\in \mathbb{R}$. By definition, $\Psi_n^{(2)}(W_{\pi_v'};\beta_2)$ is equal to 
\begin{align*}
\iint\bigg|W_{\pi_v'}\left(\begin{pmatrix}
a_v\\
&1
\end{pmatrix}\begin{pmatrix}
& 1\\
1& c'
\end{pmatrix}\right)\bigg|^2e^{in\theta}(1+c'^2)^{\frac{1}{2}-\nu}\beta_2(\theta)d\theta |a_v|_v^{-\frac{1}{2}+\nu}\sgn(a_v)^{\delta}d^{\times}a_v,
\end{align*}
where $a_v\in F_v^{\times}$, $\theta\in (-\pi,\pi]$. In particular, analyzing the support of $\beta_2$, we have 
$|c'|\leq 10$ unless $\beta_2(\theta)=0$. 

Specifying $g_v=\begin{pmatrix}
	1& c'\\
	&1
\end{pmatrix}$ in \eqref{eq5.7} leads to 
\begin{align*}
W_{\pi_v'}\left(\begin{pmatrix}
a_v\\
&1
\end{pmatrix}\begin{pmatrix}
& 1\\
1&c'
\end{pmatrix}\right)=\int_{F_v^{\times}}\psi_v(y_vc')j_{\pi_v'}(a_vy_v)W_{\pi_v'}\left(\begin{pmatrix}
y_v\\
&1
\end{pmatrix}\right)d^{\times}y_v.
\end{align*}

Substituting \eqref{5.17} into the above integral we then obtain (parallel to \eqref{5.11}):
\begin{equation}\label{equa11.4}
W_{\pi_v'}\left(\cdots\right)=\sum_{\delta_1}\frac{\sgn(-a_v)^{\delta_1}}{4\pi i}\int_{(\alpha)}\frac{\gamma_v(1/2-s,\pi_v'\otimes\sgn^{\delta_1},\psi_v)G_{\delta_1}(c',s)}{|a_v|_v^s}ds,
\end{equation}
where $\delta_1\in \mathbb{Z}/2\mathbb{Z}$, $\alpha>-1/2+\vartheta$,  $W_{\pi_v'}\left(\cdots\right)=W_{\pi_v'}\left(\begin{pmatrix}
a_v\\
&1
\end{pmatrix}\begin{pmatrix}
& 1\\
1&c'
\end{pmatrix}\right)$, and $G_{\delta_1}(c',s)$ is defined by 
\begin{equation}\label{equa11.5}
G_{\delta_1}(c',s):=\int_{F_v^{\times}}\psi_v(y_vc')\varphi_v(y_v)\sgn(y_v)^{\delta_1}|y_v|_v^{-s}d^{\times}y_v.
\end{equation}
Here $\varphi_v\in C_c^{\infty}(F_v^{\times})$ is defined as in \textsection\ref{sec4.3}. Let 
\begin{align*}
\Upsilon_v(s_1,s_2,\delta_1,\delta_2):=\gamma_v(1/2-s_1,\pi_v'\otimes\sgn^{\delta_1},\psi_v)\overline{\gamma_v(1/2-s_2,\pi_v'\otimes\sgn^{\delta_2},\psi_v)},
\end{align*}
and 
\begin{equation}\label{e11.11}
\Omega_2(a_v):=\int_{(\alpha)}\int_{(\alpha)}\frac{\Upsilon_v(s_1,s_2,\delta_1,\delta_2)G_{\delta_1}(c',s_1)\overline{G_{\delta_2}(c',s_2)}}{|a_v|_v^{s_1+\overline{s_2}}}ds_1d\overline{s_2}.
\end{equation}

According to \cite[Lemma 3.1.14]{MV10}, for $A<\Re(s)<B$ and $m\geq 0$, we have 
\begin{equation}\label{11.19}
G_{\delta_j}(c',s_j)\ll \min\left\{\left(\frac{1+|c'|_v}{1+|s_j|}\right)^m,\left(\frac{1+|s_j|}{|c'|_v}\right)^m,(1+|s_j|)^{-1/2+\varepsilon}\right\},
\end{equation}
where $j\in\{1,2\}$, and the implied constant relies on $A$, $B$, $\varepsilon$, $m$, and $\varphi_v$. %Hence, the integral in \eqref{e11.11} converges absolutely. 

Let $l\geq 0$, $m\geq 0$, $\varepsilon>0$, and $|c'|\leq 10$. By applying \eqref{eq5.21} and \eqref{11.19} with the choices $\alpha=-1/2+\Re(\nu)+\varepsilon$ and $\alpha=m$, respectively, we obtain  
\begin{equation}\label{e11.13.}
a_v^l\frac{\partial^l\Omega_2(a_v)}{\partial a_v^l}\ll \min\bigg\{C_v(\pi')^{-1+2\Re(\nu)+2\varepsilon}|a_v|_v^{1-2\Re(\nu)-2\varepsilon},\frac{C_v(\pi')^{2m+\varepsilon}}{|a_v|_v^{2m}}\bigg\},
\end{equation}
where the implied constant depends on $l$, $m$, $\varepsilon$ and $\varphi_v$. 

Substituting \eqref{equa11.4} into the definition of $\Psi_n^{(2)}(W_{\pi_v'};\beta_2)$ yields 
\begin{equation}\label{11.18}
\Psi_n^{(2)}(W_{\pi_v'};\beta_2)=\sum_{\delta_1,\delta_2}\frac{(-1)^{\delta_1+\delta_2+\delta}}{16\pi^2}\int_{-\pi}^{\pi}\mathcal{J}_v(\theta)e^{in\theta}(1+c'^2)^{\frac{1}{2}-\nu}\beta_2(\theta)d\theta,
\end{equation}
where  
\begin{equation}\label{e11.15}
\mathcal{J}_v(\theta):=\int_{F_v^{\times}}\sgn(a_v)^{\delta+\delta_1+\delta_2}|a_v|_v^{-1/2+\nu}\Omega_2(a_v)d^{\times}a_v.
\end{equation}

Since $|\Re(\nu)|\leq \vartheta\leq 7/64$, the the right hand side of \eqref{e11.13.} tends to $0$ as $|a_v|_v\to\infty$ or $|a_v|_v\to 0$. Therefore, we may preform a change of variable in \eqref{e11.15}, obtaining
\begin{equation}\label{e11.16}
\mathcal{J}_v(\theta)\ll \prod_{j=0}^{l-1}\frac{1}{|1/2-\nu+j|}\int_{F_v^{\times}}a_v^{-\frac{1}{2}+\Re(\nu)}\bigg|a_v^l\frac{\partial^l\Omega_2(a_v)}{\partial a_v^l}\bigg|d^{\times}a_v.	
\end{equation}

Substituting \eqref{e11.13.} into \eqref{e11.15}, we obtain 
\begin{equation}\label{e11.17}
\int_{0}^{C_v(\pi')^{1+\varepsilon}}a_v^{-\frac{1}{2}+\Re(\nu)}\bigg|a_v^l\frac{\partial^l\Omega_2(a_v)}{\partial a_v^l}\bigg|d^{\times}a_v\ll C_v(\pi')^{-\frac{1}{2}+\Re(\nu)+10\varepsilon},
\end{equation}
and 
\begin{equation}\label{e11.18}
\int_{C_v(\pi')^{1+\varepsilon}}^{\infty}a_v^{-\frac{1}{2}+\Re(\nu)}\bigg|a_v^l\frac{\partial^l\Omega_2(a_v)}{\partial a_v^l}\bigg|d^{\times}a_v\ll \frac{C_v(\pi')^{2m+\varepsilon}|c'|_v^{\varepsilon}}{C_v(\pi')^{(2m+1/2-\Re(\nu))(1+\varepsilon)}}.
\end{equation}

Therefore, \eqref{eq11.3}  follows from \eqref{11.18}, \eqref{e11.16}, \eqref{e11.17}, and  \eqref{e11.18} (with $m=\floor{100\varepsilon^{-1}}$). 
\end{proof}

\begin{lemma}\label{lem11.4}
Let notation be as before. Let $\varepsilon>0$. Let $n\in\mathbb{Z}$, and $l\geq 0$. Let $\beta$ be defined as above. Suppose $\sigma_v$ is a principal series. Then 
\begin{equation}\label{eq11.3}
\Psi_n^{(1)}(W_{\pi_v'};\beta_1)\ll \frac{C_v(\pi')^{1+l+\varepsilon}}{(1+|\nu|)^l},
\end{equation}
where the implied constant depends only on $\varepsilon$, $l$, and the function $\varphi_v$ (see \textsection\ref{sec4.3}). 
\end{lemma}
\begin{proof}
Let $c=-\tan\theta\in \mathbb{R}$. By a change of variable, $\Psi_n^{(1)}(W_{\pi_v'};\beta_1)$ is equal to 
\begin{align*}
\iint\bigg|W_{\pi_v'}\left(\begin{pmatrix}
a_v\\
&1
\end{pmatrix}\begin{pmatrix}
1& \\
c& 1
\end{pmatrix}\right)\bigg|^2e^{in\theta}(1+c^2)^{\frac{1}{2}-\nu}\beta_1(\theta)d\theta |a_v|_v^{-\frac{1}{2}+\nu}\sgn(a_v)^{\delta}d^{\times}a_v,
\end{align*}
where $a_v\in F_v^{\times}$, $\theta\in (-\pi,\pi]$. In particular, analyzing the support of $\beta_1$, we have 
$|c|\leq 10$ unless $\beta_1(\theta)=0$. 

Parallel to \eqref{equa11.4}, we obtain 
\begin{equation}\label{equ11.25}
W_{\pi_v'}\left(\cdots\right)=\sum_{\delta_1}\frac{\sgn(-a_v)^{\delta_1}}{4\pi i}\int_{(\alpha)}\frac{\gamma_v(1/2-s,\pi_v'\otimes\sgn^{\delta_1},\psi_v)H_{\delta_1}(c',s)}{|a_v|_v^s}ds,
\end{equation}
where $\delta_1\in \mathbb{Z}/2\mathbb{Z}$, $\alpha>-1/2+\vartheta$,  $W_{\pi_v'}\left(\cdots\right)=W_{\pi_v'}\left(\begin{pmatrix}
a_v\\
&1
\end{pmatrix}\begin{pmatrix}
1& \\
c&1
\end{pmatrix}\right)$, and  
\begin{equation}\label{e11.26}
H_{\delta_1}(c,s):=\int_{F_v^{\times}}\psi_v(y_vc)W_{\pi_v'}\left(\begin{pmatrix}
y_v&\\
& 1
\end{pmatrix}w\right)\sgn(y_v)^{\delta_1}|y_v|_v^{-s}d^{\times}y_v.
\end{equation}

Let $\beta_0\geq 0$ be a fixed smooth function such that $\beta_0(t)\equiv 1$ if $|t|\leq 1$ and $\beta_0(t)\equiv 0$ if $|t|>2$. Let 
\begin{align*}
H_{\delta_1}(c,s;\beta_0):=\int_{F_v^{\times}}\psi_v(y_vc)W_{\pi_v'}\left(\begin{pmatrix}
y_v&\\
& 1
\end{pmatrix}w\right)\sgn(y_v)^{\delta_1}\beta(y_vC_v(\pi')^{-1-\varepsilon})|y_v|_v^{-s}d^{\times}y_v.
\end{align*}

Since $1-\beta(y_vC_v(\pi')^{-1-\varepsilon})\equiv 0$ unless $|y_v|_v\geq C_v(\pi')^{1+\varepsilon}$. By Lemma \ref{lem5.4} we have
\begin{equation}\label{11.24}
H_{\delta_1}(c,s)-H_{\delta_1}(c,s;\beta_0)\ll C_v(\pi')^{-(1+\varepsilon)\Re(s)}C_v(\pi')^{-\varepsilon m},
\end{equation}
where $m>-\Re(s)$ and the implied constant depends on $m$ and $\varepsilon$. 

Substituting \eqref{5.11} into the definition of $H_{\delta_1}(c,s;\beta_0)$ yields 
\begin{multline*}
H_{\delta_1}(c,s;\beta_0)
= \int_{(0)}\mathcal{M}\varphi_v(\lambda)\gamma_v(1/2+\lambda,\pi_v',\psi_v)C_v(\pi')^{-(1+\varepsilon)(s-\lambda)}\\
\int_{F_v^{\times}}\beta(y_v)\psi_v(y_vC_v(\pi')^{1+\varepsilon}c)\sgn(y_v)^{\delta_1}|y_v|_v^{\lambda-s}d^{\times}y_vd\lambda.
\end{multline*}

Therefore, it follows from \eqref{11.19} that 
\begin{equation}\label{11.25}
H_{\delta_1}(c,s;\beta_0)\ll \frac{1}{C_v(\pi')^{(1+\varepsilon)\Re(s)}}
\min\left\{\left(\frac{C_v(\pi')^{1+\varepsilon}}{1+|s|}\right)^m,(1+|s|)^{-1/2+\varepsilon}\right\},
\end{equation}
where the implied constant depends on $m$. 

Combining \eqref{11.24} with \eqref{11.25} leads to 
\begin{multline}\label{eq11.27}
H_{\delta_1}(c,s)\ll C_v(\pi')^{-(1+\varepsilon)\Re(s)}\bigg[C_v(\pi')^{-\varepsilon m}\\
+\min\left\{\left(\frac{C_v(\pi')^{1+\varepsilon}}{1+|s|}\right)^m,(1+|s|)^{-1/2+\varepsilon}\right\}\bigg].
\end{multline}

Let $\Upsilon_v(s_1,s_2,\delta_1,\delta_2):=\gamma_v(1/2-s_1,\pi_v'\otimes\sgn^{\delta_1},\psi_v)\overline{\gamma_v(1/2-s_2,\pi_v'\otimes\sgn^{\delta_2},\psi_v)}$ be as before, and define 
\begin{equation}\label{e11.22}
\Omega_1(a_v):=\int_{(\alpha)}\int_{(\alpha)}\frac{\Upsilon_v(s_1,s_2,\delta_1,\delta_2)H_{\delta_1}(c',s_1)\overline{H_{\delta_2}(c',s_2)}}{|a_v|_v^{s_1+\overline{s_2}}}ds_1d\overline{s_2}.
\end{equation}

Substituting \eqref{equ11.25} into the definition of $\Psi_n^{(1)}(W_{\pi_v'};\beta_1)$ leads to 
\begin{equation}\label{equ11.26}
\Psi_n^{(1)}(W_{\pi_v'};\beta_1)=\sum_{\delta_1,\delta_2}\frac{(-1)^{\delta_1+\delta_2+\delta}}{16\pi^2}\int_{-\pi}^{\pi}\mathcal{I}_v(\theta)e^{in\theta}(1+c^2)^{\frac{1}{2}-\nu}\beta_1(\theta)d\theta,
\end{equation} 
where $\mathcal{I}_v(\theta)$ is defined by 
\begin{align*}
\mathcal{I}_v(\theta):=\int_{F_v^{\times}}\sgn(a_v)^{\delta+\delta_1+\delta_2}|a_v|_v^{-1/2+\nu}\Omega_1(a_v)d^{\times}a_v.
\end{align*}

\begin{comment}
Let $-1/2+\vartheta+\varepsilon\leq \Re(s_j)<0$, $j=1,2$. By Lemma \ref{lem5.4}, we may integrate by parts to deduce 
\begin{align*}
H_{\delta_j}(c,s_j)\ll \left(\frac{1+|c|_v}{1+|s_j|}\right)^m\int_{F_v^{\times}}|y_v|_v^{-\Re(s_j)+m}\sum_{j=0}^m\bigg|\frac{\partial^j}{\partial y_v^j}W_{\pi_v'}\left(\begin{pmatrix}
y_v&\\
& 1
\end{pmatrix}w\right)\bigg|d^{\times}y_v.
\end{align*}

Splitting the integral into the regions $|y_v|_v\geq C_v(\pi')^{1+\varepsilon}$ and $|y_v|_v\leq C_v(\pi')^{1+\varepsilon}$, and applying \eqref{5.12}, together with the lower bound  $\Re(s_j)\geq -1/2$, we obtain 
\begin{equation}
H_{\delta_j}(c,s_j)\ll \Big[C_v(\pi')^{-\Re(s_j)+\varepsilon}+C_v(\pi')^{-m\varepsilon+1}\Big]\cdot \left(\frac{1+|c|_v}{1+|s_j|}\right)^m.
\end{equation}

\end{comment}

Let $l\geq 0$, $m\geq 0$, $\varepsilon>0$, and $|c|\leq 10$. Applying  \eqref{eq5.21} and \eqref{eq11.27} with the choices $\alpha=-1/2+\Re(\nu)-\varepsilon$ and $\alpha=m$, respectively, we obtain  
\begin{equation}\label{e11.13}
a_v^l\frac{\partial^l\Omega_1(a_v)}{\partial a_v^l}\ll C_v(\pi')^{l+1+\varepsilon}\min\bigg\{|a_v|_v^{1-2\Re(\nu)+10\varepsilon},\frac{C_v(\pi')^{m\varepsilon}}{|a_v|_v^{2m}}\bigg\},
\end{equation}
where the implied constant depends on $l$, $m$, $\varepsilon$ and $\varphi_v$. 

Therefore, \eqref{eq11.3} follows by applying a similar argument as in the proof of Lemma \ref{lem11.2}. 
\end{proof}

\begin{remark}
The key distinction between $\Psi_n^{(2)}(W_{\pi_v'};\beta_2)$ and 
$\Psi_n^{(1)}(W_{\pi_v'};\beta_1)$ lies in the behaviors of $G_{\delta}(c',s)$ and 
$H_{\delta}(c,s)$, described in \eqref{11.19} and \eqref{eq11.27}.  
In particular, $G_{\delta}(c',s)$ imposes a sharp truncation 
$|s|\ll |c'|\ll 1$, while $H_{\delta}(c,s)$ yields only the weaker cutoff 
$|s|\ll C_v(\pi')^{1+\varepsilon}$.
\end{remark}

\begin{comment}
\begin{remark}
In the proof of Lemma \ref{lem11.2}, the bound $G_{\delta_j}(c', s_j) \ll (1+|s_j|)^{-1/2+\varepsilon}$ in \eqref{11.19} is crucial for estimating the $s_1$- and $s_2$-integrals (e.g., see \eqref{eq11.22}) since $|s_j| \asymp |c'|_v$, which can be as large as $T^{-1} \asymp C_v(\pi')^{\frac{1-\varepsilon}{2}}$.

In contrast, in the proof of Lemma \ref{lem11.4}, we do not have access to the refined bound for $H_{\delta_j}(c, s_j)$ (i.e., \eqref{eq11.27} versus \eqref{11.19}). However, in this case, $|s_j| \asymp 1 + |c| \ll 1 + \tan T\ll 1$, so the $s_1$- and $s_2$-integrals contribute $O(1)$.   
\end{remark}
\end{comment}

\begin{lemma}\label{lem11.6}
Let notation be as before. Let $\varepsilon>0$. Let $|n|\geq 1$, $m\geq 1$, and $l\geq 1$. Let $\beta$ be defined as above. Suppose $\sigma_v$ is a principal series. Then 
\begin{equation}\label{eq11.30}
\Psi_n^{(2)}(W_{\pi_v'};\beta_2)\ll  \frac{(1+|\nu|)^m}{|n|^m}\cdot C_v(\pi')^{-\frac{1}{2}+\Re(\nu)+\varepsilon},
\end{equation}
where the implied constant depends only on $\varepsilon$, $l$, $m$ and the function $\varphi_v$ (see \textsection\ref{sec4.3}). 
\end{lemma}
\begin{proof}
Recall the definition of $\Psi_n^{(2)}(W_{\pi_v'};\beta_2)$:  
\begin{align*}
\iint\bigg|W_{\pi_v'}\left(\begin{pmatrix}
a_v\\
&1
\end{pmatrix}\begin{pmatrix}
& 1\\
1& c'
\end{pmatrix}\right)\bigg|^2e^{in\theta}(1+c'^2)^{\frac{1}{2}-\nu}\beta_2(\theta)d\theta |a_v|_v^{-\frac{1}{2}+\nu}\sgn(a_v)^{\delta}d^{\times}a_v,
\end{align*}
where $a_v\in F_v^{\times}$, $\theta\in (-\pi,\pi]$, and $c'=-\cot\theta$. Similar to \eqref{11.18}, we have
\begin{align*}
\Psi_n^{(2)}(W_{\pi_v'};\beta_2)=\sum_{\delta_1,\delta_2}\frac{(-1)^{\delta_1+\delta_2+\delta}}{16\pi^2}\int_{-\pi}^{\pi}\mathcal{J}_v(\theta)(1+c'^2)^{\frac{1}{2}-\nu}\beta_2(\theta)e^{in\theta}d\theta,	
\end{align*}
where the function $\mathcal{J}_v(\theta)$ is defined by \eqref{e11.15} in the proof of Lemma \ref{lem11.2}. 

Since $\beta_2$ is smooth with compact support in $J':=\{|\theta|<\pi:\ \pi/5<|\theta|<4\pi/5\}$, all of its derivatives are supported in $J$. By integrating by parts, for $m\geq 1$, we can express $\Psi_n^{(2)}(W_{\pi_v'};\beta_2)$ as 
\begin{equation}\label{e11.28}
\sum_{\delta_1,\delta_2}\frac{(-1)^{\delta_1+\delta_2+\delta}}{16\pi^2}\cdot \frac{1}{(in)^m}\int_{J'}e^{in\theta}\frac{\partial^m}{\partial\theta^m}\left(\mathcal{J}_v(\theta)(1+c'^2)^{\frac{1}{2}-\nu}\beta_2(\theta)\right)d\theta.
\end{equation}

Recall that $\mathcal{J}_v(\theta)$ is defined by 
\begin{align*}
\int_{F_v^{\times}}\frac{\sgn(a_v)^{\delta+\delta_1+\delta_2}}{|a_v|_v^{1/2-\nu}}\int_{(\alpha)}\int_{(\alpha)}\frac{\Upsilon_v(s_1,s_2,\delta_1,\delta_2)G_{\delta_1}(c',s_1)\overline{G_{\delta_2}(c',s_2)}}{|a_v|_v^{s_1+\overline{s_2}}}ds_1d\overline{s_2}d^{\times}a_v,
\end{align*}
where $G_{\delta_1}(c',s_1)$ and $G_{\delta_2}(c',s_2)$ are defined by \eqref{equa11.5}. Notice that 
\begin{equation}\label{eq11.32}
\partial G_{\delta_1}(c',s_1)/\partial c'=2\pi i G_{\delta_1}(c',s_1-1),\ \text{and}\ \ dc'/d\theta=-(1+c'^2).
\end{equation}

Hence, the partial derivative $\partial(G_{\delta_1}(c',s_1)\overline{G_{\delta_2}(c',s_2)})/\partial \theta$ is equal to 
\begin{equation}\label{equation11.10}
-2\pi i\cdot \big[G_{\delta_1}(c',s_1-1)\overline{G_{\delta_2}(c',s_2)}-G_{\delta_1}(c',s_1)\overline{G_{\delta_2}(c',s_2-1)}\big]\cdot (1+c'^2).
\end{equation}

As a consequence, for $0\leq m'\leq m$, the function $\partial^{m'}(G_{\delta_1}(c',s_1)\overline{G_{\delta_2}(c',s_2)})/\partial \theta^{m'}$ can be expressed as a linear combination of functions of the form  $G_{\delta_1}(c',s_1+m_1)\overline{G_{\delta_2}(c',s_2+m_2)}P(c')$ for some $m_1, m_2\in \mathbb{Z}$ with $|m_1|\leq m'$ and $|m_2|\leq m'$, and $P$ is a polynomial whose degree and coefficients are all  $O((2m')!)$.   

Moreover, the partial derivative $\partial^{m'}((1+c'^2)^{\frac{1}{2}-\nu}\beta_2(\theta))/\partial \theta^{m'}$ can be expressed as a linear combination of functions of the form  $P_1(\nu)(1+c'^2)^{\frac{1}{2}-\nu}P_2(c')\beta^{(m_3)}(\theta)$, where $|m_3|\leq m'$, $P_1$ and $P_2$ are polynomials whose degree and coefficients are all  $O((2m')!)$, and $\beta^{(m_3)}$ is the $m_3$-th derivative of $\beta$.   

Incorporating this analysis into $\partial^m(\mathcal{J}_v(\theta)(1+c'^2)^{\frac{1}{2}-\nu}\beta_2(\theta))/\partial\theta^m$, together with the arguments in the proof of \eqref{e11.13} and \eqref{e11.16}, we derive, for $|c'|\leq 10$, that 
\begin{equation}\label{e11.30}
\frac{\partial^m(\mathcal{J}_v(\theta)(1+c'^2)^{\frac{1}{2}-\nu}\beta_2(\theta))}{\partial\theta^m}\ll (1+|\nu|)^m\cdot C_v(\pi')^{-\frac{1}{2}+\Re(\nu)+10\varepsilon}.	
\end{equation}

Hence, \eqref{eq11.30} follows from \eqref{e11.28} and \eqref{e11.30}.  
\end{proof}

\begin{lemma}\label{lem11.9}
Let notation be as before. Let $\varepsilon>0$. Let $|n|>0$, $m\geq 1$, and $l\geq 1$. Let $\beta$ be defined as above. Suppose $\sigma_v$ is a principal series. Then 
\begin{equation}\label{eq11.34}
\Psi_n^{(1)}(W_{\pi_v'};\beta_1)\ll  \frac{(1+|\nu|)^m}{|n|^m}\cdot C_v(\pi')^{1+\varepsilon},
\end{equation}
where the implied constant depends only on $\varepsilon$, $l$, $m$, and the function $\varphi_v$ (see \textsection\ref{sec4.3}). 
\end{lemma}
\begin{proof}
According to the proof of Lemma \ref{lem11.4}, we have 
\begin{align*}
\Psi_n^{(1)}(W_{\pi_v'};\beta_1)=\sum_{\delta_1,\delta_2}\frac{(-1)^{\delta_1+\delta_2+\delta}}{16\pi^2}\int_{J}\mathcal{I}_v(\theta)(1+c^2)^{\frac{1}{2}-\nu}\beta_1(\theta)e^{in\theta}d\theta,
\end{align*}
where $J=\supp \beta_1=(-\pi,\pi]-\supp\beta$, and $\mathcal{I}_v(\theta)$ is defined by 
\begin{align*}
\int_{F_v^{\times}}\frac{\sgn(a_v)^{\delta+\delta_1+\delta_2}}{|a_v|_v^{1/2-\nu}}\int_{(\alpha)}\int_{(\alpha)}\frac{\Upsilon_v(s_1,s_2,\delta_1,\delta_2)H_{\delta_1}(c,s_1)\overline{H_{\delta_2}(c,s_2)}}{|a_v|_v^{s_1+\overline{s_2}}}ds_1d\overline{s_2}d^{\times}a_v.
\end{align*}  
Here $H_{\delta_1}(c,s_1)$ and $H_{\delta_2}(c,s_2)$ are defined by \eqref{e11.26}.

Since $\mathcal{I}_v(\theta)$ is $\pi$-periodic and $\beta_1$ has compact support modulo $\pi$, the expression for $\Psi_n^{(1)}(W_{\pi_v'};\beta_1)$ simplifies for $m\geq 1$ by integrating by parts:
\begin{equation}\label{e11.34}
\sum_{\delta_1,\delta_2}\frac{(-1)^{\delta_1+\delta_2+\delta}}{16\pi^2}\cdot \frac{1}{(in)^m}\int_{J}e^{in\theta}\frac{\partial^m}{\partial\theta^m}\left(\mathcal{I}_v(\theta)(1+c^2)^{\frac{1}{2}-\nu}\beta_1(\theta)\right)d\theta. 
\end{equation}

By the definition \eqref{e11.26}, the function $H_{\delta_1}(c,s_1)\overline{H_{\delta_2}(c,s_2)}$ admits similar properties as $G_{\delta_1}(c',s_1)\overline{G_{\delta_2}(c',s_2)}$ in  \eqref{eq11.32} and \eqref{equation11.10}. Consequently, for $0\leq m'\leq m$, the partial derivative $\partial^{m'}(H_{\delta_1}(c,s_1)\overline{H_{\delta_2}(c,s_2)})/\partial \theta^{m'}$ can be expressed as a linear combination of functions of the form  $H_{\delta_1}(c,s_1+m_1)\overline{H_{\delta_2}(c,s_2+m_2)}P(c)$ for some $m_1, m_2\in \mathbb{Z}$ with $|m_1|\leq m'$ and $|m_2|\leq m'$, and $P$ is a polynomial whose degree and coefficients are all  $O((2m)!)$. 

Incorporating this analysis into $\partial^m(H_{\delta_1}(c,s_1)\overline{H_{\delta_2}(c,s_2)})/\partial\theta^m$, together with the arguments in the proof of \eqref{e11.13} and \eqref{e11.16}, we derive, for $|c|\leq 10$, that 
\begin{equation}\label{e11.33}
\frac{\partial^m(\mathcal{I}_v(\theta)(1+c^2)^{\frac{1}{2}-\nu}\beta_1(\theta))}{\partial\theta^m}\ll (1+|\nu|)^m\cdot C_v(\pi')^{10\varepsilon}.	
\end{equation}

Hence, \eqref{eq11.34} follows from \eqref{e11.13} (with $l=0$ therein), \eqref{e11.34} and \eqref{e11.33}. 
\end{proof}

\section{Archimedean Integrals on the Dual Side \RNum{2}}\label{sect11}
In this section we analyze the local period 
$\mathcal{P}_v(W_n,W_{\pi_v'},\overline{W_{\pi_v'}})$ (see \eqref{f11.3}) in the 
case where $\sigma_v$ is a \textit{discrete series}.  
By \cite[Theorem 9.5]{Pra90},  
\begin{align*}
\dim_{\mathbb{C}}
\mathrm{Hom}_{\mathrm{GL}_2(\mathbb{R})}
(\sigma_v \times \pi_v' \times \widetilde{\pi}_v',\,\mathbb{C})
= 1
\end{align*}
whenever $\pi_v'$ is a discrete series of weight $m'$ satisfying 
$2\max\{m', -m',m\}\ge m$.
The goal of this section is to obtain a uniform estimate for  
$\mathcal{P}_v(W_n,W_{\pi_v'},\overline{W_{\pi_v'}})$ that applies to all possible 
types of $\pi_v'$. 
\begin{prop}\label{prop11.11}
Let notation be as before. Let $n\in\mathbb{Z}$. Suppose $\sigma_v$ is a discrete series of weight $m$ and trivial central character. Let $W_n$ be defined by \eqref{eq11.1}, where $|n|\geq m$ with $|n|\equiv m\pmod{2}$. Let $l\geq 0$. Then 
\begin{equation}\label{11.39}
\big|\mathcal{P}_v(W_n,W_{\pi_v'},\overline{W_{\pi_v'}})\big|\ll |n|^{-l}C_v(\pi')^{\frac{5}{4}+l-\vartheta+(1+l)\varepsilon},
\end{equation}	
where the implied constant depends only on $l$, $\varepsilon$ and the function $\varphi_v$ (see \textsection\ref{sec4.3}).  
\end{prop}
\begin{proof}
The above upper bound \eqref{11.39} follows readily from Lemmas \ref{lem11.5} and \ref{lem11.17} below, together with the fact that $\langle W_{\pi_v'},W_{\pi_v'}\rangle\asymp 1$. 
\end{proof}

Notice that the estimate \eqref{11.39} roughly reduces the contribution from the full basis $\{\langle W_n,W_n\rangle^{-1/2}W_n:\ |n|\geq n_{\sigma_v},\ n\equiv n_{\sigma_v}\pmod{2}\}$ on the dual side to those terms where $|n|\ll C_v(\pi')^{1+2\varepsilon}$.

\begin{comment}

\begin{align*}
\widetilde{f}_n^{\sharp}(g_v):=\mu(\det g_v)|\det g_v|_v^{-1/2}\int_{F_v^{\times}}\widehat{\Phi}_v^{\sharp}((0,t)wg_v^{\iota})\mu^{-2}(t)|t|d^{\times}t,
\end{align*}
where $g_v^{\iota}$ refers to the inverse transpose of $g_v$, and $\widehat{\Phi}_v^{\sharp}$ is the Fourier transform of $\Phi_v^{\sharp}$. Let $\widetilde{f}_n(g_v):=f_n^{\sharp}(I_2,s)^{-1}\widetilde{f}_n^{\sharp}(g_v)$. By \cite[Lemma 6]{BJN23}, we have
\begin{align*}
|\mathcal{P}_v(W_n,W_{\pi_v'},\overline{W_{\pi_v'}})|^2=C_{\nu}\cdot |\Psi_v(f_n,W_{\pi_v'},\overline{W_{\pi_v'}})|\cdot |\Psi_v(\widetilde{f}_n,W_{\pi_v'},\overline{W_{\pi_v'}})|,
\end{align*}
where 
\begin{align*}
C_{\nu}:=
\end{align*}

\begin{align*}
\int_{N(F_v)\backslash G(F_v)}|W_{\pi_v'}(g_v)|^2f_n(g_v)dg_v=
\end{align*}

\end{comment}
\subsection{Matrix coefficients of $K$-isotypical vectors}
\begin{lemma}\label{lem11.10}
Let notation be as before. Let $\mu=|\cdot|_v^{\nu}\cdot \sgn^{\delta}$, with $2\nu$ being a positive odd integer. Let $\sigma_v=\sigma(\mu,\mu^{-1})$ be a discrete series. Let $n\geq 2\nu+1=:m$ be an even integer. Let $W_n$ be defined as in \eqref{whi11.2}. Then 
\begin{equation}\label{11.30}
W_n(\diag(a,1))=(-2\pi i)^{m}a^{\frac{m}{2}}e^{-2\pi a}{}_2F_1\left( -\frac{n - m}{2},\; 1;\; m;\; 4\pi a \right)\mathbf{1}_{a>0},
\end{equation}
where ${}_2F_1(a,b;c;z)$ is the Gaussian hypergeometric function. 
\end{lemma}
\begin{proof}
Substituting $m=2\nu+1$ (which is the weight $n_{\sigma_v}$) into \eqref{whi11.2} yields 
\begin{align*}
W_n(\diag(a,1))=|a|^{1-\frac{m}{2}}\sgn(a)^{\delta}\int_{F_v}\frac{(u-i)^ne^{-2\pi iua}}{(u^2+1)^{\frac{n+m}{2}}}du,
\end{align*}
where $|n|\geq m$ is an even integer. Notice that 
\begin{align*}
\int_{F_v}\frac{(u-i)^ne^{-2\pi iua}}{(u^2+1)^{\frac{n+m}{2}}}du=\int_{\mathbb{R}}\frac{(u-i)^{\frac{n-m}{2}}e^{-2\pi iua}}{(u+i)^{\frac{n+m}{2}}}du.
\end{align*}

\begin{itemize}
\item Suppose $a<0$. Shifting contour from $\mathbb{R}\to \mathbb{R}+i\infty$ leads to
\begin{equation}\label{11.31}
W_n(\diag(a,1))\equiv 0.
\end{equation}

\item Suppose $a>0$. We then obtain by shifting contour from $\mathbb{R}\to \mathbb{R}-i\infty$ that 
\begin{align*}
W_n(\diag(a,1))=\frac{-2\pi i\cdot a^{1-\frac{m}{2}}}{((n+m)/2-1)!}\cdot \frac{\partial^{\frac{n+m}{2}-1}((u-i)^{\frac{n-m}{2}}e^{-2\pi iua})}{\partial u^{\frac{n+m}{2}-1}}\bigg|_{u=-i}.
\end{align*}

\begin{comment}
\bigskip
\bigskip
\bigskip

\begin{align*}
\frac{a^{1-\frac{m}{2}}}{\left(\frac{n+m}{2}-1\right)!}\sum_{j=0}^{\frac{n+m}{2}-1}\frac{\left(\frac{n+m}{2}-1\right)!}{j!\left(\frac{n+m}{2}-j-1\right)!}\cdot \frac{\left(\frac{n-m}{2}\right)!(-2i)^{\frac{n-m}{2}-j}}{\left(\frac{n-m}{2}-j\right)!}\cdot (-2\pi i a)^{\frac{n+m}{2}-1-j}e^{-2\pi a}
\end{align*}

\begin{align*}
a^{1-\frac{m}{2}}\sum_{j=0}^{\frac{n-m}{2}}\binom{\frac{n-m}{2}}{j}\frac{(-2i)^{\frac{n-m}{2}-j}\cdot (-2\pi i a)^{\frac{n+m}{2}-1-j}}{\left(\frac{n+m}{2}-j-1\right)!}\cdot e^{-2\pi a}
\end{align*}

change variable $j\mapsto \frac{n-m}{2}-j$:
\begin{align*}
(-2\pi i)^{m-1}\sum_{j=0}^{\frac{n-m}{2}}\binom{\frac{n-m}{2}}{j}\frac{(-1)^{j}(4\pi)^j}{\left(m-1+j\right)!}\cdot a^{j+\frac{m}{2}}e^{-2\pi a}
\end{align*}

\bigskip
\bigskip
\bigskip
\end{comment}

A straightforward calculation yields 
\begin{equation}\label{11.32}
W_n\left(\begin{pmatrix}
a\\
&1
\end{pmatrix}\right)=(-2\pi i)^{m}\sum_{j=0}^{\frac{n-m}{2}}\binom{\frac{n-m}{2}}{j}\frac{(-1)^{j}(4\pi)^j}{\left(m-1+j\right)!}\cdot a^{j+\frac{m}{2}}e^{-2\pi a}\mathbf{1}_{a>0}.
\end{equation}
\end{itemize}
\begin{comment}
\begin{align*}
(-2i)^{\frac{n-m}{2}}
\sum_{j=\frac{m}{2}}^{\frac{n}{2}}\binom{\frac{n-m}{2}}{j-\frac{m}{2}}\frac{(-2\pi i )^{\frac{m}{2}-1+j}}{(j+m/2-1)!}a^{j}e^{-2\pi a}
\end{align*}
\end{comment}

By the definition of the hypergeometric function, we have 
\begin{equation}\label{e11.39}
\sum_{j=0}^{\frac{n - m}{2}} \binom{\frac{n - m}{2}}{j} \frac{(-1)^{j} (4\pi)^{j}a^{j}}{(m - 1 + j)!} = {}_2F_1\left( -\frac{n - m}{2},\; 1;\; m;\; 4\pi a \right).
\end{equation}

Therefore, \eqref{11.30} follows from \eqref{11.31},  \eqref{11.32}, and \eqref{e11.39}. 
\end{proof}

\begin{comment}
\begin{align*}
{}_2F_1\left(-\frac{n-m}{2}, 1; m; 4\pi a\right) = (m-1)\int_0^1 (1-t)^{m-2}(1-4\pi a t)^{\frac{n-m}{2}} dt
\end{align*}
\end{comment}

\begin{lemma}\label{lem11.11}
Let notation be as in Lemma \ref{lem11.10}. Let $W_n$ be the $K$-isotypical Whittaker function defined as in \eqref{whi11.2}. Let 
\begin{align*}
g_v=\begin{pmatrix}
1& b\\
& 1
\end{pmatrix}\begin{pmatrix}
a\\
& 1
\end{pmatrix}\begin{pmatrix}
\cos\theta & \sin\theta\\
-\sin\theta & \cos\theta
\end{pmatrix},\ \ \ \theta\in [-\pi,\pi).
\end{align*}
Then the function $\langle \sigma_v(g_v)W_n,W_n\rangle\langle W_n,W_n\rangle^{-1}$ can be expressed as 
\begin{equation}\label{eq11.40}
\frac{2^ma^{\frac{m}{2}}(1+a+ib)^{\frac{n-m}{2}}e^{in\theta}\mathbf{1}_{a>0}}{(1+a-ib)^{\frac{n+m}{2}}}\cdot (-1)^{\frac{n-m}{2}}P_{\frac{n-m}{2}}^{(m-1,0)}(2|c|^2-1),
\end{equation}
where $c=(1-a+ib)(1+a-ib)^{-1}$. Here $P_{\frac{n-m}{2}}^{(m-1,0)}(2|c|^2-1)$ is the Jacobi polynomial defined by  \eqref{e11.3} (see \textsection\ref{sec11.1.2}).
\end{lemma}

\begin{comment}
\begin{align*}
\langle \sigma_v(g_v)W_n,W_n\rangle=\sum_{j=0}^{\frac{n-m}{2}}\sum_{j'=0}^{\frac{n-m}{2}}\binom{\frac{n-m}{2}}{j}\binom{\frac{n-m}{2}}{j'}\frac{(-4\pi)^{j+j'}(2\pi)^{2m-2}e^{in\theta}}{\left(m-1+j\right)!\left(m-1+j'\right)!}\cdot \frac{a^{j+\frac{m}{2}}\Gamma(j+j'+m)\mathbf{1}_{a>0}}{[2\pi (a+1-ib)]^{j+j'+m}}.
\end{align*}

\begin{proof}
By \eqref{11.30} and the $K$-type of $W_n$, $\sigma_v(g_v)W_n\left(\begin{pmatrix}
y\\
&1
\end{pmatrix}\right)$ is equal to  
\begin{equation}\label{11.37}
(-2\pi i)^{m-1}e^{in\theta}\psi_v(by)\mathbf{1}_{ay>0}\sum_{j=0}^{\frac{n-m}{2}}\binom{\frac{n-m}{2}}{j}\frac{(-4\pi)^j}{\left(m-1+j\right)!}\cdot y^{j+\frac{m}{2}}e^{-2\pi y}.
\end{equation}

Using contour integral we obtain 
\begin{equation}\label{11.38}
\int_{0}^{\infty}(ay)^{j+\frac{m}{2}}e^{-2\pi ay}y^{j'+\frac{m}{2}}e^{-2\pi y}\psi_v(by)d^{\times}y=\frac{a^{j+\frac{m}{2}}\Gamma(j+j'+m)}{[2\pi (a+1-ib)]^{j+j'+m}}.
\end{equation}

Therefore, Lemma \ref{lem11.11} follows from \eqref{11.30}, \eqref{11.37} and \eqref{11.38}.
\end{proof}
\end{comment}

\begin{proof}
By \eqref{11.32} and the $K$-type of $W_n$, $\sigma_v(g_v)W_n\left(\begin{pmatrix}
y\\
&1
\end{pmatrix}\right)$ is equal to  
\begin{equation}\label{11.37}
(-2\pi i)^{m}e^{in\theta}\psi_v(by)\mathbf{1}_{ay>0}\sum_{j=0}^{\frac{n-m}{2}}\binom{\frac{n-m}{2}}{j}\frac{(-4\pi)^j}{\left(m-1+j\right)!}\cdot (ay)^{j+\frac{m}{2}}e^{-2\pi ay}.
\end{equation}

Since $W_n\left(\begin{pmatrix}
y\\
&1
\end{pmatrix}\right)\equiv 0$ unless $y>0$. Hence, it follows from \eqref{11.37} and the definition of $\langle \sigma_v(g_v)W_n,W_n\rangle$ that $\langle \sigma_v(g_v)W_n,W_n\rangle\equiv 0$ unless $a>0$.

Define the modified intertwining operator
\begin{align*}
\widetilde{f}_n^{\sharp}(g_v):=\mu(\det g_v)|\det g_v|_v^{-1/2}\int_{F_v^{\times}}\widehat{\Phi}_v^{\sharp}((0,t)wg_v^{\iota})\mu^{-2}(t)|t|d^{\times}t\in \Ind (\mu^{-1}\otimes\mu),
\end{align*}
where $g_v^{\iota}$ refers to the inverse transpose of $g_v$, and $\widehat{\Phi}_v^{\sharp}$ is the Fourier transform of $\Phi_v^{\sharp}$. Let $\widetilde{f}_n(g_v):=f_n^{\sharp}(I_2,0)^{-1}\widetilde{f}_n^{\sharp}(g_v)$. 
\begin{comment}
Then $\widetilde{f}_n(I_2)$ is equal to 
\begin{equation}\label{11.43}
\frac{\int_{F_v^{\times}}\Phi_v^{\sharp}((0,t))\mu^{-2}(t)|t|d^{\times}t}{\int_{F_v^{\times}}\Phi_v^{\sharp}((0,t))\mu^{2}(t)|t|d^{\times}t}=\frac{\int_{0}^{\infty}e^{-\pi t^2}|t|^{n-2\nu+1}d^{\times}t}{\int_{0}^{\infty}e^{-\pi t^2}|t|^{n+2\nu+1}d^{\times}t}=\frac{\pi^{m-1}\Gamma\bigl(\frac{n-m+2}{2}\bigr)}{\Gamma\bigl(\frac{n + m}{2}\bigr)}.
\end{equation}
\end{comment}

By comparing the norms in different models, there is a constant $c_m$ depending only on $\sigma_v$ such that 
\begin{equation}\label{11.42.}
\langle W_n,W_n\rangle=c_m\int_{K_v}f_n(k_v)\overline{\widetilde{f}_n(k_v)}dk_v,\ \ n\geq m:=2\nu+1,\ n\equiv m\pmod{2}.
\end{equation}

It the follows from \eqref{11.42.} that 
\begin{equation}\label{e11.42}
\frac{\langle \sigma_v(g_v)W_n,W_n\rangle}{\langle W_n,W_n\rangle}=\int_{K_v}f_n(k_vg_v)e^{-in\theta'}dk_v,\ \ k_v=\begin{pmatrix}
\cos\theta' & \sin\theta'\\
-\sin\theta' & \cos\theta'
\end{pmatrix}.
\end{equation}

\begin{comment}
\begin{align*}
\frac{\langle \sigma_v(g_v)W_n,W_n\rangle}{\langle W_n,W_n\rangle}=\int_{K_v}f_n(k_vg_v)e^{-in\theta}dk_v=\frac{a^{m/2}}{2\pi}\int_{-\pi}^{\pi}\frac{(\cos\theta+ia\sin\theta)^{\frac{n-m}{2}}e^{-in\theta}}{(\cos\theta-ia\sin\theta)^{\frac{n+m}{2}}}d\theta
\end{align*}

Suppose $0<a<1$. Then

\begin{align*}
&\cos\theta+ia\sin\theta=\frac{z+z^{-1}}{2}+a\frac{z-z^{-1}}{2}=\frac{1+a}{2}z+\frac{1-a}{2}z^{-1},\\
&\cos\theta-ia\sin\theta=\frac{z+z^{-1}}{2}-a\frac{z-z^{-1}}{2}=\frac{1-a}{2}z+\frac{1+a}{2}z^{-1}.
\end{align*}

Let $c=\frac{1-a}{1+a}$. Then 
\begin{align*}
&\frac{2^{m}}{(a+1)^{m}}\frac{1}{2\pi i}\int_{|z|=1}\frac{\left(cz^{-1}+z\right)^{\frac{n-m}{2}}z^{-n}}{\left(z^{-1}+cz\right)^{\frac{n+m}{2}}}\frac{dz}{z}\\
=&\frac{2^{m}}{(a+1)^{m}}\frac{1}{2\pi i}\int_{|z|=1}\frac{\left(1+cz^{-2}\right)^{\frac{n-m}{2}}}{\left(1+cz^2\right)^{\frac{n+m}{2}}}\frac{dz}{z}\\
\end{align*}

Now we only need to find the constant terms of $\left(1+cz^{-2}\right)^{\frac{n-m}{2}}\left(1+cz^2\right)^{-\frac{n+m}{2}}$.

\begin{align*}
\left(1+cz^{-2}\right)^{\frac{n-m}{2}}\left(1+cz^2\right)^{-\frac{n+m}{2}}=z^{m-n}\left(c+z^2\right)^{\frac{n-m}{2}}\left(1+cz^2\right)^{-\frac{n+m}{2}}
\end{align*}

Now we need to find the coefficient of $z^{n-m}$: 
\begin{align*}
\sum_{\substack{l+j=\frac{n-m}{2}\\ 0\leq l\leq \frac{n-m}{2}}}(-1)^j\binom{\frac{n-m}{2}}{l}\binom{\frac{n+m}{2}+j-1}{j}c^{2j}=\sum_{j=0}^{\frac{n-m}{2}}(-1)^j\binom{\frac{n-m}{2}}{j}\binom{\frac{n+m}{2}+j-1}{j}c^{2j}
\end{align*}

\bigskip
\bigskip
\bigskip
\bigskip

\end{comment}
Write $g_v'=\begin{pmatrix}
a& b\\
& 1
\end{pmatrix}$. By the  Iwasawa decomposition we have
\begin{equation}\label{f11.45}
k_vg_v'=\begin{pmatrix}
ay & *\\
& y^{-1}
\end{pmatrix}\begin{pmatrix}
\cos\gamma & \sin\gamma\\
-\sin\gamma & \cos\gamma
\end{pmatrix},
\end{equation}
where $*\in F_v$, $\tan\gamma=a\sin\theta'\cdot(\cos\theta'-b\sin\theta')^{-1}$, and 
\begin{align*}
y^2=\frac{1}{(\cos\theta'-b\sin\theta')^2+(a\sin\theta')^2}.
\end{align*}

Substituting  \eqref{equ11.1} and \eqref{f11.45} into 
\eqref{e11.42} yields  
\begin{equation}\label{e11.43}
\frac{\langle \sigma_v(g_v')W_n,W_n\rangle}{\langle W_n,W_n\rangle}=\frac{a^{m/2}}{2\pi}\int_{-\pi}^{\pi}\frac{(\cos\theta'-b\sin\theta'+ia\sin\theta')^{\frac{n-m}{2}}e^{-in\theta'}}{(\cos\theta'-b\sin\theta'-ia\sin\theta')^{\frac{n+m}{2}}}d\theta'.
\end{equation}

Let $z=e^{i\theta'}$. By a straightforward calculation, 
\begin{equation}\label{e11.44}
\begin{cases}
\cos\theta'-b\sin\theta'+ia\sin\theta'=\frac{1+a+ib}{2}z+\frac{1-a-ib}{2}z^{-1},\\
\cos\theta'-b\sin\theta'-ia\sin\theta'=\frac{1-a+ib}{2}z+\frac{1+a-ib}{2}z^{-1}.
\end{cases}	
\end{equation}

Let $c=(1-a+ib)(1+a-ib)^{-1}$. Substituting \eqref{e11.44} into \eqref{e11.43} leads to
\begin{equation}\label{e11.45}
\frac{\langle \sigma_v(g_v)W_n,W_n\rangle}{\langle W_n,W_n\rangle}=\frac{2^m(1+a+ib)^{\frac{n-m}{2}}\mathbf{1}_{a>0}}{2\pi ia^{-\frac{m}{2}}(1+a-ib)^{\frac{n+m}{2}}}\int_{|z|=1}\frac{z^{m-n}(\overline{c}+z^2)^{\frac{n-m}{2}}}{(1+cz^2)^{\frac{n+m}{2}}}\frac{dz}{z}.	
\end{equation}

\begin{comment}
\begin{align*}
c=\frac{1-a+ib}{1+a-ib}=\frac{2}{1+a-ib}-1
\end{align*}

\begin{align*}
&\frac{1}{2\pi}\int_{-\pi}^{\pi}\frac{(\cos\theta'-b\sin\theta'+ia\sin\theta')^{\frac{n-m}{2}}e^{-in\theta}}{(\cos\theta'-b\sin\theta'-ia\sin\theta')^{\frac{n+m}{2}}}d\theta\\
=&\frac{1}{2\pi i}\int_{|z|=1}\frac{z^{m-n}(\frac{1+a+ib}{2}z^2+\frac{1-a-ib}{2})^{\frac{n-m}{2}}}{(1+cz^2)^{\frac{n+m}{2}}(\frac{1+a-ib}{2})^{\frac{n+m}{2}}}\frac{dz}{z}\\
=&\frac{1}{2\pi i}\int_{|z|=1}\frac{z^{m-n}(z^2+c')^{\frac{n-m}{2}}(\frac{1+a+ib}{2})^{\frac{n-m}{2}}}{(1+cz^2)^{\frac{n+m}{2}}(\frac{1+a-ib}{2})^{\frac{n+m}{2}}}\frac{dz}{z}\\
=&\frac{2^m(1+a+ib)^{\frac{n-m}{2}}}{(1+a-ib)^{\frac{n+m}{2}}}\frac{1}{2\pi i}\int_{|z|=1}\frac{z^{m-n}(c'+z^2)^{\frac{n-m}{2}}}{(1+cz^2)^{\frac{n+m}{2}}}\frac{dz}{z}\\
\end{align*}
\end{comment}

The coefficient of $z^{n-m}$ in the Taylor expansion of $(\overline{c}+z^2)^{\frac{n-m}{2}}(1+cz^2)^{-\frac{n+m}{2}}$ is 
\begin{equation}\label{eq11.46}
\sum_{j=0}^{\frac{n-m}{2}}(-1)^j\binom{\frac{n-m}{2}}{j}\binom{\frac{n+m}{2}+j-1}{j}|c|^{2j}={}_2F_1\left(-\frac{n-m}{2}, \frac{n+m}{2}; 1; |c|^2\right),
\end{equation}
where ${}_2F_1(a,b;c;z)$ is the Gaussian hypergeometric function. 

As a consequence of \eqref{e11.45} and \eqref{eq11.46}, in conjunction with the identity 
\begin{align*}
{}_2F_1\left(\frac{m-n}{2}, \frac{n+m}{2}; 1; |c|^2\right)=P_{\frac{n-m}{2}}^{(0,m-1)}(1-2|c|^2)=(-1)^{\frac{n-m}{2}}P_{\frac{n-m}{2}}^{(m-1,0)}(2|c|^2-1),
\end{align*}
we derive the formula \eqref{eq11.40}. Here the above equalities can be found in  \cite[p.210-212]{MOS66}. 
\end{proof}

Furthermore, we also establish the following coarser bound, which holds uniformly in $n$.
\begin{lemma}\label{lem11.14}
Let notation be as in Lemma \ref{lem11.10}. Let $W_n$ be defined as in \eqref{whi11.2}. Let $g_v=\begin{pmatrix}
1& b\\
& 1
\end{pmatrix}\begin{pmatrix}
a\\
& 1
\end{pmatrix}\begin{pmatrix}
\cos\theta & \sin\theta\\
-\sin\theta & \cos\theta
\end{pmatrix}$, $a>0$, $b\in\mathbb{R}$, and $\theta\in [-\pi,\pi)$. Then
\begin{equation}\label{eq11.48}
\frac{|\langle \sigma_v(g_v)W_n,W_n\rangle|}{\langle W_n,W_n\rangle}\ll \frac{1}{\sqrt{a+a^{-1}+a^{-1}b^2}},
\end{equation} 
where the implied constant is absolute. 
\end{lemma}
\begin{proof}
By a straightforward calculation, we have 
\begin{align*}
g_v\in \begin{pmatrix}
1& b\\
& 1
\end{pmatrix}\begin{pmatrix}
a\\
& 1
\end{pmatrix}K_v\subseteq  Z(F_v)K_v\begin{pmatrix}
y\\
& 1
\end{pmatrix}K_v,
\end{align*}
where $y=\frac{a^2+b^2+1+\sqrt{(a^2+b^2+1)^2-4a^2}}{2a}$. 

Let $\theta'\in[0,\pi/2]$ be such that $\cos\theta'=|1-y|/(1+y)$.
Taking $b=0$ and $a=y$ in \eqref{eq11.40}, we obtain 
\begin{equation}\label{f11.51}
\frac{|\langle \sigma_v(g_v)W_n,W_n\rangle|}{\langle W_n,W_n\rangle}=\Big|(\sin\theta')^mP_{\frac{n-m}{2}}^{(m-1,0)}(\cos2\theta')\Big|.
\end{equation}

Since $m\geq 2$, we may  substitute the bound (see \cite[(20) on p. 234]{HS14})
\begin{align*}
\max_{\theta'\in[0,\pi/2]}\Big|(\sin\theta')^{m-1}P_{\frac{n-m}{2}}^{(m-1,0)}(\cos2\theta')\Big|\leq 1
\end{align*}
into the formula \eqref{f11.51}, obtaining 
\begin{align*}
\frac{|\langle \sigma_v(g_v)W_n,W_n\rangle|}{\langle W_n,W_n\rangle}\leq \sin\theta'=\frac{2\sqrt{y}}{1+y}\ll y^{-1/2},
\end{align*}
which leads to \eqref{eq11.48}. 
\end{proof}

\begin{remark}
It follows from \cite{CHH88} that 
\begin{equation}\label{eq11.49}
|\langle \sigma_v(g_v)W_n,W_n\rangle|\leq \langle W_n,W_n\rangle\cdot |\Xi_v(g_v)|,
\end{equation}
where $\Xi_v$ is the Harish-Chandra spherical function on $\mathrm{PGL}_2(F_v)=\mathrm{PGL}_2(\mathbb{R})$. Utilizing the standard bound for $|\Xi_v(g_v)|$ would give the bound 
\begin{align*}
\frac{|\langle \sigma_v(g_v)W_n,W_n\rangle|}{\langle W_n,W_n\rangle}\ll \frac{\log(a+a^{-1}+a^{-1}b^2)}{\sqrt{a+a^{-1}+a^{-1}b^2}},
\end{align*}
which is slightly weaker than \eqref{eq11.48}. 
\end{remark}

\subsection{Coordinates and partitions}
For $g_v\in \overline{G}(F_v)$, we may write 
\begin{align*}
g_v=\begin{pmatrix}
1& b\\
& 1
\end{pmatrix}\begin{pmatrix}
a\\
& 1
\end{pmatrix}\begin{pmatrix}
\cos\theta & \sin\theta\\
-\sin\theta & \cos\theta
\end{pmatrix}
\end{align*} 
by Iwasawa decomposition, where $b\in F_v$, $a\in F_v^{\times}$, and $-\pi\leq \theta<\pi$.

\begin{itemize}
\item Suppose $\theta\neq \pm\pi/2$. Let $c=-\tan\theta\in \mathbb{R}$. By \eqref{f11.4}, we can parametrize $g_v$ by 
\begin{equation}\label{11.40}
g_v=\sgn(\cos\theta)\begin{pmatrix}
1& b\\
& 1
\end{pmatrix}\begin{pmatrix}
a\\
& 1
\end{pmatrix}\begin{pmatrix}
1 & -c\\
& 1\end{pmatrix}\begin{pmatrix}
\sqrt{1+c^2}& \\
& \frac{1}{\sqrt{1+c^2}}\end{pmatrix}\begin{pmatrix}
1\\
c & 1
\end{pmatrix}.
\end{equation}

\item Suppose $\theta\neq 0, \pm\pi$. Let $c'=-\cot\theta\in \mathbb{R}$. By \eqref{f11.5}, we have
\begin{equation}\label{11.41}
g_v=\sgn(\sin\theta)\begin{pmatrix}
1& b\\
& 1
\end{pmatrix}\begin{pmatrix}
a\\
& 1
\end{pmatrix}\begin{pmatrix}
1& c'\\
&1
\end{pmatrix}\begin{pmatrix}
\sqrt{1+c'^2}& \\
&-\frac{1}{\sqrt{1+c'^2}}
\end{pmatrix}\begin{pmatrix}
& 1\\
1& c'
\end{pmatrix}.
\end{equation}
\end{itemize}

Let $\beta:(-\pi,\pi]\to [0,1]$ be a smooth function with compact support in $J:=\{|\theta|<\pi:\ \pi/5<|\theta|<4\pi/5\}$ such that $\beta(\theta)\equiv 1$ for all $\pi/4<|\theta|<3\pi/4$. Let $\beta_2=\beta$ and $\beta_1=1-\beta$ be as in the statement of Lemma \ref{lemma11.1}.

Recall the definition 
\begin{align*}
|\mathcal{P}_v(W_n,W_{\pi_v'},\overline{W_{\pi_v'}})|^2=\int_{\overline{G}(F_v)}\frac{\langle\sigma_v(g_v)W_n,W_n\rangle\cdot |\langle\pi_v'(g_v)W_{\pi_v'},W_{\pi_v'}\rangle|^2}{\langle W_n,W_n\rangle\cdot |\langle W_{\pi_v'},W_{\pi_v'}\rangle|^2} dg_v.\tag{\ref{f11.3}}  
\end{align*}

By \eqref{11.40} and \eqref{11.41}, we have
\begin{equation}\label{11.42}
|\mathcal{P}_v(W_n,W_{\pi_v'},\overline{W_{\pi_v'}})|^2=\mathcal{P}_v^{(1)}(W_{\pi_v'};\beta)+\mathcal{P}_v^{(2)}(W_{\pi_v'};\beta) ,
\end{equation}
where 
\begin{multline}\label{11.45}
\mathcal{P}_v^{(1)}(W_{\pi_v'};\beta):=\frac{|\langle W_{\pi_v'},W_{\pi_v'}\rangle|^{-2}}{\langle W_n,W_n\rangle}\int_{F_v}\int_{F_v^{\times}}\int_{-\pi}^{\pi}\langle\sigma_v\left(\begin{pmatrix}
a& b\\
& 1
\end{pmatrix}\right)W_n,W_n\rangle\\
\Big|\int_{F_v^{\times}}\overline{\varphi_v(y)}W_{\pi_v'}\left(\begin{pmatrix}
ay(1+c^2)\\
c & 1
\end{pmatrix}\right)\psi_v((b-ac)y)d^{\times}y\Big|^2e^{in\theta}\beta_1(\theta)|a|_v^{-1}d\theta d^{\times}adb,
\end{multline}
and 
\begin{multline}\label{11.47}
\mathcal{P}_v^{(2)}(W_{\pi_v'};\beta):=\frac{|\langle W_{\pi_v'},W_{\pi_v'}\rangle|^{-2}}{\langle W_n,W_n\rangle}\int_{F_v}\int_{F_v^{\times}}\int_{-\pi}^{\pi}\langle\sigma_v\left(\begin{pmatrix}
a& b\\
& 1
\end{pmatrix}\right)W_n,W_n\rangle\\
\Big|\int_{F_v^{\times}}\overline{\varphi_v(y)}W_{\pi_v'}\left(\begin{pmatrix}
& ay(1+c'^2)\\
1 & c'
\end{pmatrix}\right)\psi_v((b+ac')y)d^{\times}y\Big|^2e^{in\theta}\beta_2(\theta)|a|_v^{-1}d\theta d^{\times}adb.
\end{multline}

\subsection{Estimates of $\mathcal{P}_v^{(1)}(W_{\pi_v'};\beta)$ and $\mathcal{P}_v^{(2)}(W_{\pi_v'};\beta)$}
By \eqref{equa11.4},  
\begin{multline}\label{11.48}
W_{\pi_v'}\left(\begin{pmatrix}
&ay(1+c'^2)\\
1 & c'
\end{pmatrix}\right)=\sum_{\delta_1}\frac{(-1)^{\delta_1}}{4\pi i}\\
\int_{(\alpha)}\frac{\gamma_v(1/2-s,\pi_v'\otimes\sgn^{\delta_1},\psi_v)G_{\delta_1}(c',s)}{|ay(1+c'^2)|_v^s}ds,
\end{multline}
where $G_{\delta_1}(c',s)$ is defined by 
\begin{align*}
G_{\delta_1}(c',s):=\int_{F_v^{\times}}\psi_v(y_vc')\varphi_v(y_v)\sgn(y_v)^{\delta_1}|y_v|_v^{-s}d^{\times}y_v.\tag{\ref{equa11.5}}
\end{align*}

Substituting \eqref{11.48} into \eqref{11.47}, we derive 
\begin{multline*}
\mathcal{P}_v^{(2)}(W_{\pi_v'};\beta)=-\frac{|\langle W_{\pi_v'},W_{\pi_v'}\rangle|^{-2}}{\langle W_n,W_n\rangle}\int_{\mathbb{R}}\int_{0}^{\infty}\int_{-\pi}^{\pi}\langle\sigma_v\left(\begin{pmatrix}
a& b\\
& 1
\end{pmatrix}\right)W_n,W_n\rangle e^{in\theta}\\
\sum_{\delta_1,\delta_2}\frac{(-1)^{\delta_1+\delta_2}}{16\pi^2}\int_{(\alpha)}\int_{(\alpha)}\frac{\Upsilon_v(s_1,s_2,\delta_1,\delta_2)G_{\delta_1}(c',s_1)\overline{G_{\delta_2}(c',s_2)}}{(a(1+c'^2))^{s_1+\overline{s_2}}}\int_{0}^{\infty}\overline{\varphi_v(y_1)}|y_1|_v^{-s_1} \\
\psi_v((b+ac')y_1)d^{\times}y_1\int_{0}^{\infty}\varphi_v(y_2)\overline{\psi}_v((b+ac')y_2)|y_2|_v^{-\overline{s_2}}d^{\times}y_2ds_1ds_2\beta_2(\theta)d\theta
\frac{d^{\times}a}{|a|_v}db,
\end{multline*}
where $\alpha>-1/2+\vartheta$, and 
\begin{align*}
\Upsilon_v(s_1,s_2,\delta_1,\delta_2):=\gamma_v(1/2-s_1,\pi_v'\otimes\sgn^{\delta_1},\psi_v)\overline{\gamma_v(1/2-s_2,\pi_v'\otimes\sgn^{\delta_2},\psi_v)}.
\end{align*}

Similarly, we have
we have 
\begin{multline*}
\mathcal{P}_v^{(1)}(W_{\pi_v'};\beta)= -\frac{|\langle W_{\pi_v'},W_{\pi_v'}\rangle|^{-2}}{\langle W_n,W_n\rangle}\int_{\mathbb{R}}\int_{0}^{\infty}\int_{-\pi}^{\pi}\langle\sigma_v\left(\begin{pmatrix}
a& b\\
& 1
\end{pmatrix}\right)W_n,W_n\rangle e^{in\theta}\\
\sum_{\delta_1,\delta_2}\frac{(-1)^{\delta_1+\delta_2}}{16\pi^2}\int_{(\alpha)}\int_{(\alpha)}\frac{\Upsilon_v(s_1,s_2,\delta_1,\delta_2)H_{\delta_1}(c',s_1)\overline{H_{\delta_2}(c',s_2)}}{(a(1+c'^2))^{s_1+\overline{s_2}}}\int_{0}^{\infty}\overline{\varphi_v(y_1)}|y_1|_v^{-s_1} \\
\psi_v((b-ac)y_1)d^{\times}y_1\int_{0}^{\infty}\varphi_v(y_2)\overline{\psi}_v((b-ac)y_2)|y_2|_v^{-\overline{s_2}}d^{\times}y_2ds_1ds_2\beta_1(\theta)d\theta
a^{-1}d^{\times}adb,
\end{multline*} 
where $H_{\delta_j}(c,s_j)$ is defined by \eqref{e11.26}.

For convenience, we define the auxiliary function 
\begin{align*}
h(a,b,s_1,s_2):=\int_{-\pi}^{\pi}\frac{G_{\delta_1}(c',s_1)\overline{G_{\delta_2}(c',s_2)}}{(1+c'^2)^{s_1+\overline{s_2}}}\kappa(\theta;a,b,s_1,s_2)\beta_2(\theta)e^{in\theta}d\theta,
\end{align*}
where $\kappa(\theta;a,b,s_1,s_2)$ is defined by 
\begin{align*}
\int_{0}^{\infty}\overline{\varphi_v(y_1)}|y_1|_v^{-s_1}
\psi_v((b+ac')y_1)d^{\times}y_1\int_{0}^{\infty}\varphi_v(y_2)\overline{\psi}_v((b+ac')y_2)|y_2|_v^{-\overline{s_2}}d^{\times}y_2.
\end{align*}

Let $l'\geq 1$ be an \textit{odd} integer  and $\varepsilon>0$. We have the following decomposition: 
\begin{equation}\label{eq11.58}
\mathcal{P}_v^{(2)}(W_{\pi_v'};\beta)=-\mathcal{P}_v^{(2,1)}(W_{\pi_v'};\beta)-\mathcal{P}_v^{(2,2)}(W_{\pi_v'};\beta),
\end{equation}
where the term $\mathcal{P}_v^{(2,1)}(W_{\pi_v'};\beta)$ is defined by 
\begin{multline*}
\frac{1}{|\langle W_{\pi_v'},W_{\pi_v'}\rangle|^{2}}\sum_{\delta_1,\delta_2}\frac{(-1)^{\delta_1+\delta_2}}{16\pi^2}\int_{0}^{C_v(\pi')^{1+\varepsilon}}\int_{\mathbb{R}}\langle\sigma_v\left(\begin{pmatrix}
1& b\\
& 1
\end{pmatrix}\begin{pmatrix}
a& \\
& 1
\end{pmatrix}\right)W_n,W_n\rangle \\
\langle W_n,W_n\rangle^{-1}\int_{(-1/2+\vartheta+\varepsilon)}\int_{(-1/2+\vartheta+\varepsilon)}\frac{\Upsilon_v(s_1,s_2,\delta_1,\delta_2)}{a^{1+s_1+\overline{s_2}}}
h(a,b,s_1,s_2)ds_1ds_2
dbd^{\times}a,
\end{multline*}
and the term $\mathcal{P}_v^{(2,2)}(W_{\pi_v'};\beta)$ is defined by 
\begin{multline*}
\frac{1}{|\langle W_{\pi_v'},W_{\pi_v'}\rangle|^{2}}\sum_{\delta_1,\delta_2}\frac{(-1)^{\delta_1+\delta_2}}{16\pi^2}\int_{C_v(\pi')^{1+\varepsilon}}^{\infty}\int_{\mathbb{R}}\langle\sigma_v\left(\begin{pmatrix}
1& b\\
& 1
\end{pmatrix}\begin{pmatrix}
a\\
& 1
\end{pmatrix}\right)W_n,W_n\rangle \\
\langle W_n,W_n\rangle^{-1}\int_{(l'/2)}\int_{(l'/2)}\frac{\Upsilon_v(s_1,s_2,\delta_1,\delta_2)}{a^{1+s_1+\overline{s_2}}}
h(a,b,s_1,s_2)ds_1ds_2
dbd^{\times}a.
\end{multline*}

Similarly, we have the decomposition
\begin{align*}
\mathcal{P}_v^{(1)}(W_{\pi_v'};\beta)=-\mathcal{P}_v^{(1,1)}(W_{\pi_v'};\beta)-\mathcal{P}_v^{(1,2)}(W_{\pi_v'};\beta)
\end{align*}
as in \eqref{eq11.58}, where each $\mathcal{P}_v^{(1,j)}(W_{\pi_v'};\beta)$, $1\leq j\leq 2$, is defined in the same manner as 
$\mathcal{P}_v^{(2,1)}(W_{\pi_v'};\beta)$, except that $G_{\delta_j}(c,s_j)$ is replaced by $H_{\delta_j}(c,s_j)$.

%Building on the arguments in the proof of \eqref{11.63}, while replacing $G_{\delta_j}(c,s_j)$ with $H_{\delta_j}(c,s_j)$, and substituting the bound \eqref{11.19} with \eqref{eq11.27}, and further utilizing the manipulations employed in the proofs of Lemmas \ref{lem11.4} and \ref{lem11.9}, we obtain the following result:

%\subsubsection{Majorization of $\mathcal{P}_v^{(i,1)}(W_{\pi_v'};\beta)$: $i=1, 2$}
%We may integrate by parts to bound $\mathcal{P}_v^{(2,1)}(W_{\pi_v'};\beta)$ and $\mathcal{P}_v^{(2,2)}(W_{\pi_v'};\beta)$, as in the proof of Proposition \ref{prop11.4}. However, the integral $\mathcal{P}_v^{(2,2)}(W_{\pi_v'};\beta)$ is more delicate and requires extra manipulation. We will take advantage of Lemma \ref{lem11.15} and   

\begin{lemma}\label{lem11.5}
Let $l\geq 0$ be an integer and $0<\varepsilon<10^{-3}$. Then 
\begin{equation}\label{11.21}
\mathcal{P}_v^{(1,1)}(W_{\pi_v'};\beta)+\mathcal{P}_v^{(2,1)}(W_{\pi_v'};\beta)\ll n^{-l}\cdot C_v(\pi')^{\frac{5}{2}+l-2\vartheta+10(l+1)\varepsilon},
\end{equation}	
where the implied constant depends only on $l$, $\varepsilon$, and the function $\varphi_v$ (see \textsection\ref{sec4.3}).  
\end{lemma}
\begin{proof}
Note that $\beta_2(\theta)\equiv 0$ unless $c'=-\cot\theta\in [-10,10]$. Integrating by parts we can express $h(a,b,s_1,s_2)$ as 
\begin{equation}\label{c11.82}
\frac{1}{(in)^l}\int_{-\pi}^{\pi}e^{in\theta}\frac{\partial^l}{\partial\theta^l}\bigg[\frac{G_{\delta_1}(c',s_1)\overline{G_{\delta_2}(c',s_2)}}{(1+c'^2)^{s_1+\overline{s_2}}}
\kappa(\theta;a,b,s_1,s_2)\beta_2(\theta)\bigg]d\theta.
\end{equation}

By \eqref{eq5.21}, Lemma \ref{lem11.14} and the triangle inequality, we obtain 
\begin{multline}\label{fc11.82}
\mathcal{P}_v^{(2,1)}(W_{\pi_v'};\beta)\ll \frac{1}{n^l}\int_{0}^{C_v(\pi')^{1+\varepsilon}}\frac{C_v(\pi')^{-1+2\vartheta+2\varepsilon}}{a^{-1+2\vartheta+2\varepsilon}\sqrt{a+a^{-1}}}\int_{(\alpha)}\int_{(\alpha)}\int_{\mathbb{R}}\\
\bigg|\int_{-\pi}^{\pi}e^{in\theta}\frac{\partial^l}{\partial\theta^l}\bigg[\frac{G_{\delta_1}(c',s_1)\overline{G_{\delta_2}(c',s_2)}}{(1+c'^2)^{s_1+\overline{s_2}}}\kappa(\theta;a,b,s_1,s_2)\beta_2(\theta)\bigg]d\theta\bigg|dbds_1ds_2
d^{\times}a,
\end{multline}
where $\alpha=-1/2+\vartheta+\varepsilon$. Here, we recall that 
\begin{multline*}
\kappa(\theta;a,b,s_1,s_2):=\int_{0}^{\infty}\overline{\varphi_v(y_1)}|y_1|_v^{-s_1}
\psi_v((b+ac')y_1)d^{\times}y_1\\
\int_{0}^{\infty}\varphi_v(y_2)\overline{\psi}_v((b+ac')y_2)|y_2|_v^{-\overline{s_2}}d^{\times}y_2,
\end{multline*}
where $\varphi_v$ has a compact support and $|c'|\leq 10$. Hence, for $l_1, l_2\geq 0$, integrating by parts, we have (in parallel to \eqref{11.19}) that 
\begin{multline}\label{fc11.84}
\frac{\partial^{l_2}}{\partial\theta^{l_2}}\kappa(\theta;a,b,s_1,s_2)\ll_{l_1,l_2} (1+a)^{l_2}(1+c'^2)^{l_2}\\
\prod_{j=1}^2\min\left\{\left(\frac{1+|b+ac'|_v}{1+|s_j|}\right)^{l_1},\left(\frac{1+|s_j|}{|b+ac'|_v}\right)^{l_1}, (1+|s_j|)^{-1/2+\varepsilon}\right\}.
\end{multline}

Observe that, for $l$ sufficiently large, \eqref{fc11.84} restricts the variable $b$ to the window 
\begin{equation}\label{11.87.}
|b+ac'|\ll 1+|s_j|+C_v(\pi')^{\varepsilon}.
\end{equation}

Substituting \eqref{11.19},  \eqref{equation11.10}, \eqref{fc11.84} and \eqref{11.87.} into \eqref{fc11.82},  we thus obtain 
\begin{multline*}
\mathcal{P}_v^{(2,1)}(W_{\pi_v'};\beta)\ll \frac{1}{n^l}\int_{0}^{C_v(\pi')^{1+\varepsilon}}\frac{C_v(\pi')^{-1+2\vartheta+2\varepsilon}}{a^{-1+2\vartheta+2\varepsilon}\sqrt{a+a^{-1}}}\int_{(\alpha)}\int_{(\alpha)}\\
(1+|s_1|+|s_2|)^{l+5\varepsilon}\mathbf{1}_{|s_1|\ll C_v(\pi')^{\varepsilon}}\mathbf{1}_{|s_1|\ll C_v(\pi')^{\varepsilon}}ds_1ds_2
d^{\times}a+n^{-l}C_v(\pi')^{-l},
\end{multline*}
which leads to the estimate 
\begin{equation}\label{c11.86}
\mathcal{P}_v^{(2,1)}(W_{\pi_v'};\beta)\ll n^{-l}\cdot C_v(\pi')^{-1/2+10(l+1)\varepsilon}.
\end{equation}
Here we have used the fact that $2\vartheta<1/2$ to compute the $a$-integral. 

Replacing the estimates \eqref{11.19} for $G_{\delta_j}(c',s_j)$ with the bound \eqref{eq11.27} for 
$H_{\delta_j}(c',s_j)$, we note that $|s_j|$ may now be as large as 
$C_v(\pi')^{1+\varepsilon}$ (rather than $C_v(\pi')^{\varepsilon}$ as in the estimate for 
$\mathcal{P}_v^{(2,1)}(W_{\pi_v'};\beta)$).  
This leads to
\begin{multline*}
\mathcal{P}_v^{(2,1)}(W_{\pi_v'};\beta)\ll \frac{1}{n^l}\int_{0}^{C_v(\pi')^{1+\varepsilon}}\frac{C_v(\pi')^{-1+2\vartheta+2\varepsilon}}{a^{-1+2\vartheta+2\varepsilon}\sqrt{a+a^{-1}}}\int_{(\alpha)}\int_{(\alpha)}C_v(\pi')^{1-2\vartheta}\\
(1+|s_1|+|s_2|)^{l+5\varepsilon}\mathbf{1}_{|s_1|\ll C_v(\pi')^{1+\varepsilon}}\mathbf{1}_{|s_1|\ll C_v(\pi')^{1+\varepsilon}}ds_1ds_2
d^{\times}a+n^{-l}C_v(\pi')^{-l}.
\end{multline*}
As a consequence,  
\begin{equation}\label{cf11.87}
\mathcal{P}_v^{(1,1)}(W_{\pi_v'};\beta)\ll n^{-l}\cdot C_v(\pi')^{\frac{1}{2}-2\vartheta+2\varepsilon}\cdot C_v(\pi')^{2+l+10(l+1)\varepsilon}.
\end{equation}

Therefore, \eqref{11.21} follows from \eqref{c11.86} and \eqref{cf11.87}. 
\end{proof}

\begin{lemma}\label{lem11.17}
Let $l\geq 0$ be an integer. Then 
\begin{equation}\label{11.66}
\mathcal{P}_v^{(1,2)}(W_{\pi_v'};\beta)+\mathcal{P}_v^{(2,2)}(W_{\pi_v'};\beta)\ll n^{-l}\cdot C_v(\pi')^{-l}, 
\end{equation}	
where the implied constant depends only on $l$ and the function $\varphi_v$ (see \textsection\ref{sec4.3}).  
\end{lemma}
\begin{proof}
Recall that $\mathcal{P}_v^{(2,2)}(W_{\pi_v'};\beta)$ is defined by 
\begin{multline}\label{fc11.87}
\frac{1}{|\langle W_{\pi_v'},W_{\pi_v'}\rangle|^{2}}\sum_{\delta_1,\delta_2}\frac{(-1)^{\delta_1+\delta_2}}{16\pi^2}\int_{C_v(\pi')^{1+\varepsilon}}^{\infty}\int_{\mathbb{R}}\langle\sigma_v\left(\begin{pmatrix}
a& b\\
& 1
\end{pmatrix}\right)W_n,W_n\rangle \\
\langle W_n,W_n\rangle^{-1}\int_{(l'/2)}\int_{(l'/2)}\frac{\Upsilon_v(s_1,s_2,\delta_1,\delta_2)}{a^{1+s_1+\overline{s_2}}}
h(a,b,s_1,s_2)ds_1ds_2
dbd^{\times}a.
\end{multline}

Notice that the integrand $a^{-1-s_1-\overline{s_2}}\Upsilon_v(s_1,s_2,\delta_1,\delta_2)
h(a,b,s_1,s_2)$ is holomorphic in $\Re(s_1)\geq 0$ and $\Re(s_2)\geq 0$. Therefore, the formula in \eqref{fc11.87} is independent of $l'$ as long as $l'\geq 0$. Thus, we may take $l'=\floor{100\varepsilon^{-1}l}$. 

By \eqref{eq5.21}, Lemma \ref{lem11.14}, \eqref{c11.82} and the triangle inequality, we obtain
\begin{multline}\label{c11.87}
\mathcal{P}_v^{(2,2)}(W_{\pi_v'};\beta)\ll C_v(\pi')^{l'}\int_{C_v(\pi')^{1+\varepsilon}}^{\infty}\int_{\mathbb{R}}\frac{1}{a^{l'+3/2}} \int_{(l'/2)}\int_{(l'/2)}\\
\frac{1}{n^l}\int_{-\pi}^{\pi}\bigg|\frac{\partial^l}{\partial\theta^l}\bigg[\frac{G_{\delta_1}(c',s_1)\overline{G_{\delta_2}(c',s_2)}}{(1+c'^2)^{s_1+\overline{s_2}}}
\kappa(\theta;a,b,s_1,s_2)\beta_2(\theta)\bigg]\bigg|d\theta ds_1ds_2
dbd^{\times}a.
\end{multline}

Substituting \eqref{11.19},  \eqref{equation11.10} and \eqref{fc11.84} into \eqref{c11.87}, we derive 
\begin{equation}\label{c11.88}
\mathcal{P}_v^{(2,2)}(W_{\pi_v'};\beta)\ll \frac{C_v(\pi')^{l'+\varepsilon}}{n^l}\int_{C_v(\pi')^{1+\varepsilon}}^{\infty}\frac{(1+a)^{l}}{a^{l'+3/2}} 
d^{\times}a\ll \frac{n^{-l}\cdot C_v(\pi')^{2l}}{C_v(\pi')^{l'\varepsilon}}.
\end{equation}

Replacing the estimates \eqref{11.19} for $G_{\delta_j}(c',s_j)$ with the bound \eqref{eq11.27} for 
$H_{\delta_j}(c',s_j)$, we obtain, in parallel with \eqref{c11.88}, that 
\begin{equation}\label{c11.89}
\mathcal{P}_v^{(1,2)}(W_{\pi_v'};\beta)\ll \frac{n^{-l}\cdot C_v(\pi')^{100+2l}}{C_v(\pi')^{l'\varepsilon}}.
\end{equation}

Since $l'=\floor{100\varepsilon^{-1}l}$, then it follows from \eqref{c11.88} and \eqref{c11.89} that 
\begin{align*}
\mathcal{P}_v^{(1,2)}(W_{\pi_v'};\beta)+\mathcal{P}_v^{(2,2)}(W_{\pi_v'};\beta)\ll n^{-l}\cdot C_v(\pi')^{-l}.
\end{align*}

Therefore, \eqref{11.66} holds. 
\end{proof}

\section{Archimedean Integrals on the Dual Side \RNum{3}}\label{sect12}

In this section we analyze the local integrals of 
$\mathcal{M}_{\cusp}^{\du}(s_1,s_2)$ and $\mathcal{M}_{\Eis}^{\du}(s_1,s_2)$ 
at complex places and  obtain uniform upper bounds for the associated integral transforms.

\subsection{Representations of $\mathrm{PGL}_2(\mathbb{C})$}\label{sec11.1.2}
Let $F_v\simeq\mathbb{C}$ and $K_v=\mathrm{SU}_2(\mathbb{C})$. The representations $\sigma_v$ of $\mathrm{PGL}_2(\mathbb{C})$ are of the form $\Ind\mu\otimes\mu^{-1}$, where $\mu$ is the character given by $\mu(z)=|z|_v^{\nu}[z]^{n_0}$. Here $z=re^{i\theta}\in F_v=\mathbb{C}$ with  $r>0$, $0\leq\theta<2\pi$, and $[z]:=e^{i\theta}$. The parameters satisfy either $\nu\in i\mathbb{R}_+$ and $n_0\in \mathbb{Z}$, or $\nu\in (-\vartheta, \vartheta)$ and $n_0=0$. 

Upon exchanging $\mu$ and $\mu^{-1}$, we may assume $n_0\geq 0$. Let $n\geq 2n_0$ be an even integer. Let $V_n$ be the space of homogeneous polynomials in $\mathbb{C}[z_1,z_2]$ of degree $n$, equipped with the inner product
\begin{align*}
\langle P(z_1,z_2),Q(z_1,z_2)\rangle_{\mathrm{SU}_2(\mathbb{C})}=\int_{\mathrm{SU}_2(\mathbb{C})}P(z_1,z_2)\overline{Q(z_1,z_2)}dz_1dz_2,
\end{align*}
where the action of $\mathrm{SU}_2(\mathbb{C})$ is given by 
\begin{align*}
k_v. P(z_1,z_2)=P((z_1,z_2)g),\ \ k_v\in \mathrm{SU}_2(\mathbb{C}). 
\end{align*}

Then $\{z_1^{n-l}z_2^l:\ 0\leq l\leq n\}$ forms an orthogonal basis of $V_n$. Let $n\geq 2n_0$ be an even integer. Let $f_{n,l}$ be a section in $\sigma_v$ defined by 
\begin{align*}
f_{n,l}\left(\begin{pmatrix}
a& b\\
&d
\end{pmatrix}g\right)=\mu(a)\mu^{-1}(d)|a/d|_v^{1/2}f_{n,l}(g),\ \ g\in G(F_v),
\end{align*}
and for $k_v\in K_v=\mathrm{SU}_2(\mathbb{C})$, 
\begin{align*}
f_{n,l}(k_v)=\frac{\langle k_v.z_1^{n-l}z_2^{l}, z_1^{n-l_0}z_2^{l_0}\rangle_{\mathrm{SU}_2(\mathbb{C})}}{\langle z_1^{n-l_0}z_2^{l_0}, z_1^{n-l_0}z_2^{l_0}\rangle_{\mathrm{SU}_2(\mathbb{C})}},\ \ \ l_0=n/2-n_0.
\end{align*}

An integral representation  for $f_{n,l}$ as a Godement section is given explicitly in \cite[\textsection 6]{JL70}. Let $l_0=n/2-n_0$, $\alpha=l-l_0,$
 and $\beta=n-l-l_0
$. A straightforward calculation yields 
\begin{align*}
&f_{n,l}\left(\begin{pmatrix}
e^{i\gamma_1}& \\
&e^{-i\gamma_1}
\end{pmatrix}k_v\begin{pmatrix}
e^{i\gamma_2}& \\
&e^{-i\gamma_2}
\end{pmatrix}\right)=e^{2in_0\gamma_1}e^{i(n-2l)\gamma_2}f_{n,l}(k_v),\ \ k_v\in \mathrm{SU}_2(\mathbb{C}),\\
&f_{n,l}\left(\begin{pmatrix}
\cos\theta & \sin\theta \\
-\sin\theta  & 
\cos\theta 
\end{pmatrix}\right)=(\sin\theta)^{\alpha} (\cos\theta)^{\beta}P_{l_0}^{(\alpha,\beta)}(\cos2\theta),\ \ 0\leq \theta<2\pi,
\end{align*}
where $P_{l_0}^{(\alpha,\beta)}(x)$ is the Jacobi polynomial defined by 
\begin{equation}\label{e11.3}
P_{l_0}^{(\alpha,\beta)}(x)=\sum_{i=0}^{l_0}\binom{l_0+\alpha}{i}\binom{l_0+\beta}{l_0-i}\left(\frac{x-1}{2}\right)^{l_0-i}\left(\frac{x+1}{2}\right)^{i}.
\end{equation}
\begin{comment}
$P_{l_0}^{(\alpha,\beta)}(\cos2\theta)$ is the Jacobi polynomial defined by 
\begin{align*}
P_{l_0}^{(\alpha,\beta)}(\cos2\theta)=\sum_{i=0}^{l_0}\binom{l_0+\alpha}{i}\binom{l_0+\beta}{l_0-i}(-1)^{l_0-i}\cos^{2i}\theta\sin^{2l_0-2i}\theta.
\end{align*} 
\end{comment}

Utilizing polar coordinates and orthogonality condition of Jacobi polynomials (e.g., see \cite[p.212]{MOS66}), we have 
\begin{align*}
\langle f_{n,l},f_{n,l}\rangle=\int_{K_v}|f_{n,l}(k_v)|^2dk_v=\frac{1}{2l_0+\alpha+\beta+1}\cdot \frac{\Gamma(l_0+\alpha+1)\Gamma(l_0+\beta+1)}{\Gamma(l_0+\alpha+\beta+1)\Gamma(l_0+1)}.
\end{align*}

Let $n\geq 2n_0$ be an even integer and $0\leq l\leq n$. Let $l_0=n/2-n_0$, $\alpha=l-l_0,$
 and $\beta=n-l-l_0
$. It is known that 
\begin{equation}\label{eq12.2}
\big\{\langle f_{n,l},f_{n,l}\rangle^{-1/2}f_{n,l}:\ 0\leq l\leq n,\ n\geq 2n_0,\ n\equiv 0\pmod{2}\big\}
\end{equation}
forms an orthonormal basis of $\sigma_v$.

\begin{prop}\label{prop11.18}
Let notation be as before. Let $n\geq 2n_0$ be an even integer and $0\leq l\leq n$. Suppose $\sigma_v$ is a principal series. Let $l_1, l_2\geq 0$. Then 
\begin{multline*}
\frac{\Psi_v(f_{n,l},W_{\pi_v'},\overline{W_{\pi_v'}})}{\sqrt{\langle W_{n,l},W_{n,l}\rangle}}\ll C_v(\pi')^{1+\varepsilon}\mathbf{1}_{l=n/2}\cdot \min\bigg\{\frac{C_v(\pi')^{l_1+l_2}}{(1+|\nu|)^{l_1}(1+|n_0|)^{l_2}},\\
\frac{\log n\cdot C_v(\pi')^{l_1+l_2}\mathbf{1}_{n_0\geq 3}}{n^{3/2}(1+|\nu|)^{l_1-1}(1+|n_0|)^{l_2}},\ \frac{C_v(\pi')^{l_1}\mathbf{1}_{n_0=2, n\geq 1000}}{n^{3/2}(1+|\nu|)^{l_1}},\ \frac{C_v(\pi')^{l_1}\mathbf{1}_{n_0\in\{0,1\}, n\geq 10}}{n^{3/2}(1+|\nu|)^{l_1}}\bigg\},
\end{multline*}
where the implied constant depends on $l_1$, $l_2$, $\varepsilon$, and the function $\varphi_v$ (see \textsection\ref{sec4.3}).
\end{prop}
\begin{proof}
Proposition \ref{prop11.18} is a direct consequence of the decomposition presented in \eqref{11.76} in \textsection\ref{sec11.5.1}, combined with the results of Propositions \ref{prop11.19} and \ref{prop11.24} in \textsection\ref{sec11.5.3} and \textsection\ref{sec11.5.4}, respectively.
\end{proof}  

The proof of Proposition \ref{prop11.18} follows a structure analogous to that of Proposition \ref{prop11.4}, albeit with certain technical modifications. The primary challenge arises from the significantly more intricate representation theory of $\mathrm{SU}(2)$ compared to that of $\mathrm{SO}(2)$. This complexity manifests in the appearance of Jacobi polynomials in the complex setting, contrasting with the straightforward exponential functions $e^{in\theta}$ that arise in the real case.

\subsection{Decomposition of $\Psi_v(f_{n,l},W_{\pi_v'},\overline{W_{\pi_v'}})$}\label{sec11.5.1}
Recall that 
\begin{align*}
\mathrm{SU}(2)=\bigg\{\begin{pmatrix}
e^{i\gamma_1} & 0 \\
0 & e^{-i\gamma_1}
\end{pmatrix}
\begin{pmatrix}
\cos\theta & \sin\theta \\
-\sin\theta & \cos\theta
\end{pmatrix}
\begin{pmatrix}
e^{i\gamma_2}& 0 \\
0 & e^{-i\gamma_2}
\end{pmatrix}:\ \theta, \gamma_1,\gamma_2\in [0,2\pi)\bigg\}.
\end{align*}
The measure on $\mathrm{SU}(2)$ is taken as $\pi^{-2}\sin\theta\cos\theta d\theta d\gamma_1 d\gamma_2$ so that $\Vol(\mathrm{SU}(2))=1$. 

Let $k_v=\begin{pmatrix}
e^{i\gamma_1} & 0 \\
0 & e^{-i\gamma_1}
\end{pmatrix}
\begin{pmatrix}
\cos\theta & \sin\theta \\
-\sin\theta & \cos\theta
\end{pmatrix}
\begin{pmatrix}
e^{i\gamma_2}& 0 \\
0 & e^{-i\gamma_2}
\end{pmatrix}\in \mathrm{SU}(2)$. 

\begin{itemize}
\item Suppose $\theta\neq \pm\pi/2$. Let $c=-\tan\theta\in \mathbb{R}$. By \eqref{f11.4}, we can write $k_v$ as 
\begin{align*}
\sgn(\cos\theta)\begin{pmatrix}
e^{i\gamma_1} & 0 \\
0 & e^{-i\gamma_1}
\end{pmatrix}
\begin{pmatrix}
1 & -c\\
& 1\end{pmatrix}
\begin{pmatrix}
\sqrt{1+c^2}& \\
& \frac{1}{\sqrt{1+c^2}}\end{pmatrix}\begin{pmatrix}
1\\
c& 1
\end{pmatrix}\begin{pmatrix}
e^{i\gamma_2}& 0 \\
0 & e^{-i\gamma_2}
\end{pmatrix}.
\end{align*}

\item Suppose $\theta\neq 0, \pm\pi$. Let $c'=-\cot\theta\in \mathbb{R}$. By \eqref{f11.5}, we can write $k_v$ as  
\begin{align*}
\sgn(\sin\theta)\begin{pmatrix}
e^{i\gamma_1} & 0 \\
0 & e^{-i\gamma_1}
\end{pmatrix}\begin{pmatrix}
1& c'\\
&1
\end{pmatrix}\begin{pmatrix}
\sqrt{1+c'^2}& \\
&-\frac{1}{\sqrt{1+c'^2}}
\end{pmatrix}\begin{pmatrix}
&1\\
1& c'
\end{pmatrix}\begin{pmatrix}
e^{i\gamma_2}& 0 \\
0 & e^{-i\gamma_2}
\end{pmatrix}.
\end{align*}
\end{itemize}

Let $\beta:(-\pi,\pi]\to [0,1]$ be a smooth function with compact support in 
\begin{align*}
J:=\{|\theta|<\pi:\ \pi/5<|\theta|<4\pi/5\}
\end{align*}
such that $\beta(\theta)\equiv 1$ for all $\pi/4<|\theta|<3\pi/4$. Let $\beta_2=\beta$ and $\beta_1=1-\beta$. 

Similar to \eqref{equ11.2} we have the decomposition 
\begin{equation}\label{11.76}
\frac{\Psi_v(f_{n,l},W_{\pi_v'},\overline{W_{\pi_v'}})}{\langle W_{n,l},W_{n,l}\rangle^{1/2}}=|\langle W_{\pi_v'},W_{\pi_v'}\rangle|^{-2}\cdot\Big[\Psi_{n,l}^{(1)}(W_{\pi_v'};\beta_1)+\Psi_{n,l}^{(2)}(W_{\pi_v'};\beta_2)\Big],
\end{equation}
where 
\begin{multline}\label{11.77}
\Psi_{n,l}^{(1)}(W_{\pi_v'};\beta_1):=\frac{\varsigma}{\sqrt{\langle f_{n,l},f_{n,l}\rangle}}\int_{F_v^{\times}}\int_{-\pi}^{\pi}\\
\int_{-\pi}^{\pi}\bigg|W_{\pi_v'}\left(\begin{pmatrix}
a_v\\
&1
\end{pmatrix}\begin{pmatrix}
1 & \\
c & 1
\end{pmatrix}\begin{pmatrix}
e^{i\gamma_2}& 0 \\
0 & e^{-i\gamma_2}
\end{pmatrix}\right)\bigg|^2e^{i(n-2l)\gamma_2}d\gamma_2\\
(\sin\theta)^{\alpha+1} (\cos\theta)^{\beta+1}P_{l_0}^{(\alpha,\beta)}(\cos2\theta)(1+c^2)^{\frac{1}{2}-\nu}\beta_1(\theta)d\theta |a_v|_v^{-\frac{1}{2}+\nu}[a_v]^{n_0}d^{\times}a_v,
\end{multline}
with $c=-\tan\theta\in \mathbb{R}$ and $\beta_1(\theta)=1-\beta(\theta)$; and 
\begin{multline}\label{11.78}
\Psi_{n,l}^{(2)}(W_{\pi_v'};\beta_2):=\frac{\varsigma}{\sqrt{\langle f_{n,l},f_{n,l}\rangle}}\int_{F_v^{\times}}\int_{-\pi}^{\pi}\\
\int_{-\pi}^{\pi}\bigg|W_{\pi_v'}\left(\begin{pmatrix}
a_v\\
&1
\end{pmatrix}\begin{pmatrix}
&1\\
1& c'
\end{pmatrix}\begin{pmatrix}
e^{i\gamma_2}& 0 \\
0 & e^{-i\gamma_2}
\end{pmatrix}\right)\bigg|^2e^{i(n-2l)\gamma_2}d\gamma_2\\
(\sin\theta)^{\alpha+1} (\cos\theta)^{\beta+1}P_{l_0}^{(\alpha,\beta)}(\cos2\theta)(1+c'^2)^{\frac{1}{2}-\nu}\beta_2(\theta)d\theta  |a_v|_v^{-\frac{1}{2}+\nu}[a_v]^{n_0}d^{\times}a_v,
\end{multline}
with $c'=-\cot\theta\in \mathbb{R}$ and $\beta_2(\theta)=\beta(\theta)$. Here $\varsigma$ is an absolute constant depending only on the normalization of the Haar measure on $\mathrm{SU}(2)$. 

\subsection{Uniform bounds for Jacobi polynomials}\label{sec11.5.2}
To handle the Jacobi polynomials \( P_{l_0}^{(\alpha,\beta)}(\cos 2\theta) \) in \eqref{11.77} and \eqref{11.78}, we will take advantage of the asymptotic Hilb-type formula (see Szeg\"{o} \cite[Theorem 8.21.12]{Sze59}):  
\begin{multline}\label{c11.77}
(\sin\theta)^{\alpha}(\cos\theta)^{\beta}P_{n}^{(\alpha,\beta)}(\cos 2\theta) = \frac{\Gamma(n+\alpha+1)}{n!}(\theta/\sin\theta)^{1/2}N^{-\alpha}J_{\alpha}(N\theta)\\
+\begin{cases}
\theta^{1/2} O(n^{-3/2}), & \quad \text{if } cn^{-1} \leq \theta \leq \pi - \varepsilon, \\
\theta^{\alpha+2} O(n^{\alpha}), & \quad \text{if } 0<\theta \leq cn^{-1},
\end{cases}
\end{multline} 
where \( N = n+(\alpha+\beta+1)/2 \), and \( c \) and \( \varepsilon \) are fixed positive numbers. Here, \( J_{\alpha}(z) \) denotes the Bessel function of the first kind of order \( \alpha \), and the implied constants in \eqref{c11.77} depend on \( \alpha \), \( \beta \), \( c \), and \( \varepsilon \). 

Let $\alpha = \beta = n_0 - 1 \geq 0$, and let $N = l_0 + n_0 + 1/2 = (n + 1)/2$. For $0 < \theta < \pi/2$, \cite[Main Theorem]{FW85} refines \eqref{c11.77} by establishing an asymptotic Hilb-type formula for Jacobi polynomials:
\begin{multline}\label{f11.77}
(\sin\theta)^{n_0 - 1} (\cos\theta)^{n_0 - 1} P_{l_0+1}^{(n_0 - 1, n_0 - 1)}(\cos 2\theta) 
= \frac{\Gamma(l_0 + n_0 + 1)}{\Gamma(l_0 + 2)} \sqrt{\frac{\theta}{\sin\theta \cos\theta}} \\
\cdot \left[ \frac{J_{n_0 - 1}(N\theta)}{N^{n_0 - 1}} + \sum_{l = 1}^{m - 1} \frac{A_l(\theta) J_{n_0 - 1 + l}(N\theta)}{N^{n_0 - 1 + l}} + E_m(\theta; N) \right],
\end{multline}
where the coefficients $A_l(\theta)$ and the error term $E_m(\theta; N)$ are explicitly defined in \cite[(3.15)--(3.16)]{FW85}.

In this section we estimate $A_l(\theta)$ and $E_m(\theta;N)$ directly and derive an effective version of \eqref{f11.77}. The main result it the following.
\begin{prop}\label{prop11.16}
Let notation be as before. Let $\alpha=n_0-1\geq 0$, $0<\theta<\pi/2$, and  $m\geq 1$. Suppose $N\geq 100m$. Then we have the expansion 
\begin{multline}\label{f11.79}
(\sin\theta)^{n_0-1}(\cos\theta)^{n_0-1}P_{l_0+1}^{(n_0-1,n_0-1)}(\cos2\theta)=\frac{\Gamma(l_0+n_0+1)}{\Gamma(l_0+2)}\sqrt{\frac{\theta}{\sin\theta\cos\theta}}\\
\cdot \bigg[\frac{J_{n_0-1}(N\theta)}{N^{n_0-1}}+O\left(\frac{2^{\alpha}e^{\alpha}m\sqrt{\theta}}{N^{n_0+1/2}}\right)+O\left(\frac{e^{\alpha}2^mm!(m+\alpha)^m}{2^{\alpha}\Gamma(\alpha+1/2)}\cdot \frac{\theta^{\alpha}}{N^m}\right)\bigg],
\end{multline}
and 
\begin{multline}\label{f11.80.}
(\sin\theta)^{n_0-1}(\cos\theta)^{n_0-1}P_{l_0+1}^{(n_0-1,n_0-1)}(\cos2\theta)=\frac{\Gamma(l_0+n_0+1)}{\Gamma(l_0+2)}\sqrt{\frac{\theta}{\sin\theta\cos\theta}}\\
\cdot \bigg[\frac{J_{n_0-1}(N\theta)}{N^{n_0-1}}+\frac{1-\theta\cot\theta}{2\theta}\cdot\frac{(n_0-3/2)(n_0-1/2)J_{n_0}(N\theta)}{N^{n_0}}+O\left(\frac{n_0^5\theta^{n_0+1}}{2^{n_0}\Gamma(n_0)N^{2}}\right)\bigg].
\end{multline}
Here the implied constants in \eqref{f11.79}
and \eqref{f11.80.} are absolute.

\end{prop}

The estimates \eqref{f11.79} and \eqref{f11.80.} follow from Lemmas \ref{lem11.16.} and \ref{lem11.17.}, which are proved below. The proof of Proposition \ref{prop11.16} will be presented immediately after that of Lemma \ref{lem11.17.}.

\begin{lemma}\label{lem11.16.}
Let notation be as before. Let $\alpha=n_0-1\geq 0$ and $0<\theta<\pi/2$.  Then
\begin{equation}\label{11.105}
E_m(\theta;N)\ll\frac{e^{\alpha}m!(2m+2\alpha)^m}{2^{\alpha}\Gamma(\alpha+1/2)}\cdot \frac{\theta^{\alpha}}{N^m},	
\end{equation}
where the implied constant is absolute. 
\end{lemma}
\begin{proof}
For simplicity, we will follow the notation in \cite{FW85} and use it freely throughout. Recall that $\alpha=\beta$, i.e., $b=\alpha-\beta=0$. So $b_k(\theta)\equiv 0$ (see (3.2) in loc. cit.) unless $k=0$. Let $\Delta_{0,m}(\theta,\phi)$ be defined as in (2.29) in loc. cit.. Let 
\begin{equation}\label{eq11.79}
\phi_j(\theta):=\frac{1}{2^{j-1}j!}\cdot \frac{J_{j-1/2}(\theta)}{\theta^{j-1}J_{1/2}(\theta)},\ \ \ j\geq 2.
\end{equation}

Define the $\psi_{0,i}(\theta)$ by 
\begin{equation}\label{f11.80}
\bigg[1+\sum_{j=2}^{\infty}x^{j-1}\phi_j(\theta)\bigg]^{\alpha-1/2}=\sum_{i=0}^{\infty}\psi_{0,i}(\theta)x^i.
\end{equation}

By Cauchy theorem, we derive 
\begin{equation}\label{eq11.80}
\psi_{0,i}(\theta)=\frac{1}{2\pi i}\int_{|z|=100}\bigg[1+\sum_{j=2}^{\infty}z^{j-1}\phi_j(\theta)\bigg]^{\alpha-1/2}\frac{dz}{z^{i+1}}.
\end{equation}

Substituting \eqref{eq11.79} into \eqref{eq11.80}, along with the bounds for the Bessel $J$-functions, we obtain 
\begin{equation}\label{f11.82}
|\psi_{0,i}(\theta)|\ll  \frac{1}{100^i}\bigg[\sum_{j=1}^{\infty}\frac{50^{j-1}}{j!}\cdot \left(\frac{e\theta}{2j}\right)^{j}\bigg]^{\alpha-1/2}\ll 100^{-i}e^{\alpha}, 
\end{equation}
where the implied constant is absolute. 

Incorporating  \eqref{f11.80} into the definition of $\Delta_{0,m}(\theta,\phi)$ leads to 
\begin{align*}
\Delta_{0,m}(\theta,\phi)=\sum_{i=m}^{\infty}\psi_{0,i}(\theta)(\theta^2-\phi^2)^{i-m}=\sum_{j=0}^{\infty}\psi_{j+m}(\theta)(\theta^2-\phi^2)^{j}.
\end{align*}

As a consequence, we obtain 
\begin{equation}\label{f11.83}
\frac{d^{m-i}}{d\phi^{m-i}}\Delta_{0,m}(\theta,\phi)=\sum_{j=0}^{\infty}\psi_{j+m}(\theta)\frac{d^{m-i}}{d\phi^{m-i}}(\theta^2-\phi^2)^{j}.
\end{equation}

By a straightforward calculation, we have 
\begin{equation}\label{f11.84}
\frac{d^{m-i}}{d\phi^{m-i}}(\theta^2-\phi^2)^{j}=(m-i)!\sum_{u=0}^{j}\binom{j}{u}\binom{2u}{m-i}(-1)^u\theta^{2(j-u)}\phi^{2u-m+i}.
\end{equation}

\begin{comment}
\begin{align*}
\max_{0\leq \phi\leq \theta}\bigg|\frac{d^{m-i}}{d\phi^{m-i}}(\theta^2-\phi^2)^{j}\bigg|\leq (m-i)!\theta^{2j-m+i}2^j\binom{2j}{m-i} 
\end{align*}
\end{comment}

It then follows from \eqref{f11.82}, \eqref{f11.83} and \eqref{f11.84} that 
\begin{equation}\label{f11.85}
\max_{0\leq \phi\leq \theta}\bigg|\frac{d^{m-i}}{d\phi^{m-i}}\Delta_{0,m}(\theta,\phi)\bigg|\ll \sum_{j=0}^{\infty}\frac{e^{\alpha} (m-i)!\theta^{2j-m+i}2^j}{100^{j+m}}\cdot\binom{2j}{m-i}.
\end{equation}

Let $M_{0,m}(\theta)$ be defined by \cite[(2.39)]{FW85}, namely,   
\begin{equation}\label{f11.86}
M_{0,m}(\theta)=\sum_{i=0}^m\binom{m}{i}(2m+2\alpha-1)^i\theta^{m-i}\max_{0\leq \phi\leq \theta}\bigg|\frac{d^{m-i}}{d\phi^{m-i}}\Delta_{0,m}(\theta,\phi)\bigg|
\end{equation}

Combining \eqref{f11.85} and \eqref{f11.86} yields 
\begin{equation}\label{f11.87}
M_{0,m}(\theta)\ll \frac{e^{\alpha}m!}{100^{m}}\sum_{i=0}^m\binom{m}{i}(2m+2\alpha-1)^i=\frac{e^{\alpha}m!(2m+2\alpha)^m}{100^{m}}.
\end{equation}

According to the definition of $E_m(\theta;N)$ in \cite[(3.16)]{FW85}, we have 
\begin{equation}\label{f11.88}
E_m(\theta;N)\ll\frac{\theta^{-\alpha}}{2^{\alpha}\Gamma(\alpha+1/2)}\cdot \frac{M_{0,m}(\theta)\theta^m}{N^m}\int_0^{\theta}(\theta^2-\phi^2)^{\alpha-\frac{1}{2}}d\phi.
\end{equation}

Therefore, \eqref{11.105} follows from \eqref{f11.87} and \eqref{f11.88}.
\end{proof}

\begin{lemma}\label{lem11.17.}
Let notation be as before. Let $\alpha=n_0-1\geq 0$. Then
\begin{equation}\label{f11.89}
A_l(\theta)\ll \frac{\theta^l\cdot \Gamma(l+\alpha+1/2)e^{\alpha}}{50^l\Gamma(\alpha+1/2)}. 	
\end{equation}
\end{lemma}
\begin{proof}
By the definition \cite[(3.13) \& (3.15)]{FW85}, we have
\begin{equation}\label{f11.90}
A_l(\theta)\ll \frac{\theta^l\cdot 2^{l}\Gamma(l+\alpha+1/2)\psi_{0,l}(\theta)}{\Gamma(\alpha+1/2)},
\end{equation}
where $\psi_{0,l}(\theta)$ is the coefficient defined by \eqref{f11.80}. Therefore, \eqref{f11.89} follows from \eqref{f11.90} and \eqref{f11.82}.
\end{proof}

\begin{proof}[Proof of Proposition \ref{prop11.16}]
By Lemma \ref{lem11.17.} and uniform bounds for the Bessel $J$-functions, we obtain 
\begin{equation}\label{f11.93}
\sum_{l=1}^{m-1}\frac{A_l(\theta)J_{n_0-1+l}(N\theta)}{N^{n_0-1+l}}\ll \frac{\sqrt{\theta}e^{\alpha}}{N^{n_0-3/2}}\sum_{l=1}^{m-1}\frac{\Gamma(l+\alpha+1/2)}{10^l N^l\Gamma(\alpha+1/2)}.
\end{equation}
	
Since $\Gamma(l+\alpha+1/2)/\Gamma(\alpha+1/2)\leq 2^{l+\alpha-1/2}l!\leq 2^{l+\alpha-1/2}(m-1)^l$, it follows from \eqref{f11.93} that 
\begin{equation}\label{f11.94}
\sum_{l=1}^{m-1}\frac{A_l(\theta)J_{n_0-1+l}(N\theta)}{N^{n_0-1+l}}\ll \frac{2^{\alpha}e^{\alpha}\sqrt{\theta}}{N^{n_0-1/2}}\sum_{l=1}^{m-1}\frac{(m-1)^l}{N^l}\ll \frac{2^{\alpha}e^{\alpha}m\sqrt{\theta}}{N^{n_0+1/2}}.
\end{equation}
	
Hence, \eqref{f11.79}  follows from \eqref{11.105} and \eqref{f11.94}. The formula \eqref{f11.80.} follows from \cite[Corollary 1]{FW85}, with an explicit upper bound for the coefficient $E_2$ defined in \cite[(4.27)]{FW85}.
\end{proof}

\subsection{Estimates of $\Psi_{n,l}^{(2)}(W_{\pi_v'};\beta_2)$}\label{sec11.5.3}
Recall that $\Psi_{n,l}^{(2)}(W_{\pi_v'};\beta_2)$ (see \eqref{11.78} in \textsection\ref{sec11.5.1}) is defined by 
\begin{multline*}
\frac{\varsigma}{\sqrt{\langle f_{n,l},f_{n,l}\rangle}}\int_{F_v^{\times}}\int_{-\pi}^{\pi}\int_{-\pi}^{\pi}\bigg|W_{\pi_v'}\left(\begin{pmatrix}
a_v\\
&1
\end{pmatrix}\begin{pmatrix}
&1\\
1& c'
\end{pmatrix}\begin{pmatrix}
e^{i\gamma_2}& 0 \\
0 & e^{-i\gamma_2}
\end{pmatrix}\right)\bigg|^2e^{i(n-2l)\gamma_2}\\
(\sin\theta)^{\alpha+1} (\cos\theta)^{\beta+1}P_{l_0}^{(\alpha,\beta)}(\cos2\theta)(1+c'^2)^{\frac{1}{2}-\nu}\beta_2(\theta)d\theta d\gamma_2 |a_v|_v^{-\frac{1}{2}+\nu}[a_v]^{n_0}d^{\times}a_v,
\end{multline*}
with $c'=-\cot\theta\in \mathbb{R}$ and $\beta_2(\theta)=\beta(\theta)$. Here $\varsigma$ is an absolute constant depending only on the normalization of the Haar measure on $\mathrm{SU}(2)$. 

The main result in this subsection is the following.
\begin{prop}\label{prop11.19}
Let notation be as in \textsection\ref{sec11.5.1}. Let $l'\geq 0$. Then 
\begin{multline}\label{e11.96}
\Psi_{n,l}^{(2)}(W_{\pi_v'};\beta_2)\ll C_v(\pi')^{-\frac{1}{2}+\Re(\nu)+\varepsilon}\mathbf{1}_{l=n/2}\cdot \min\bigg\{\frac{1}{(1+|\nu|)^{l'}(1+|n_0|)^{l'}},\\
\frac{\log n\cdot \mathbf{1}_{n_0\geq 3}}{n^{3/2}(1+|\nu|)^{l'-1}(1+|n_0|)^{l'}},\ \frac{\mathbf{1}_{n_0=2, n\geq 1000}}{n^{3/2}(1+|\nu|)^{l'}},\ \frac{\mathbf{1}_{n_0\in\{0,1\}, n\geq 10}}{n^{3/2}(1+|\nu|)^{l'}}\bigg\},
\end{multline}
where the implied constant depends on $l'$, $\varepsilon$, and the function $\varphi_v$ (see \textsection\ref{sec4.3}).
\end{prop}
\begin{proof}
The estimate \eqref{e11.96} follows readily from Lemmas \ref{lem11.16}, \ref{lemma11.17}, \ref{lem11.21}, and \ref{lem11.22} below, together with the assumption that $\langle W_{\pi_v'},W_{\pi_v'}\rangle\asymp 1$ (see \textsection\ref{sec4.3}). 
\end{proof}

\begin{lemma}\label{lem11.16}
Let notation be as in \textsection\ref{sec11.5.1}. Let $l'\geq 0$. Then 
\begin{equation}\label{11.79}
\Psi_{n,l}^{(2)}(W_{\pi_v'};\beta_2)\ll \frac{C_v(\pi')^{-\frac{1}{2}+\Re(\nu)+\varepsilon}\cdot \mathbf{1}_{l=n/2}}{(1+|\nu|)^{l'-1}(1+|n_0|)^{l'}},
\end{equation}
where the implied constant depends on $l'$, $\varepsilon$, and the function $\varphi_v$ (see \textsection\ref{sec4.3}).
\end{lemma}
\begin{proof}
\begin{comment}
\begin{align*}
W_{\pi_v'}\left(\begin{pmatrix}
a_v\\
&1
\end{pmatrix}w\begin{pmatrix}
1& c'\\
&1
\end{pmatrix}\right)=\int_{F_v^{\times}}\psi_v(y_vc')j_{\pi_v'}(a_vy_v)W_{\pi_v'}\left(\begin{pmatrix}
y_v\\
&1
\end{pmatrix}\right)d^{\times}y_v.
\end{align*}

\begin{align*}
j_{\pi_v'}(y_v)=\frac{1}{8\pi^2 i}\sum_{m\in \mathbb{Z}}[-y_v]^{-m}\int_{(\alpha)}\gamma_v(1/2-s,\pi_v'\otimes[\cdot]^{-m},\psi_v)|y_v|_v^{-s}ds.
\end{align*}

\begin{align*}
j_{\pi_v'}(a_vy_v)=\frac{1}{8\pi^2 i}\sum_{m\in \mathbb{Z}}[-a_vy_v]^{-m}\int_{(\alpha)}\gamma_v(1/2-s,\pi_v'\otimes[\cdot]^{-m},\psi_v)|a_vy_v|_v^{-s}ds.
\end{align*}
\end{comment}
By \eqref{eq5.7} and \eqref{5.17.}, $W_{\pi_v'}\left(\begin{pmatrix}
a_v\\
&1
\end{pmatrix}\begin{pmatrix}
&1\\
1& c'
\end{pmatrix}\begin{pmatrix}
e^{i\gamma_2}& 0 \\
0 & e^{-i\gamma_2}
\end{pmatrix}\right)$ is equal to  
\begin{equation}\label{e11.80}
\frac{1}{8\pi^2 i}\sum_{m\in \mathbb{Z}}[-a_v]^{-m}
\int_{(\alpha)}\frac{\gamma_v(1/2-s,\pi_v'\otimes[\cdot]^{-m},\psi_v)\widetilde{G}_m(c',s)}{|a_v|_v^{s}}ds,
\end{equation}
where 
\begin{align*}
\widetilde{G}_m(c',s):=\int_{F_v^{\times}}W_{\pi_v'}\left(\begin{pmatrix}
y_v\\
& 1
\end{pmatrix}\begin{pmatrix}
e^{i\gamma_2}& 0 \\
0 & e^{-i\gamma_2}
\end{pmatrix}\right)\psi_v(y_vc')[y_v]^{-m}|y_v|_v^{-s}d^{\times}y_v.
\end{align*}

Recall that $W_{\pi_v'}(\diag(y_v,1))=\varphi_v(y_v)=\varphi_v(|y_v|_v^{1/2})$ (see \textsection\ref{sec4.3}). We obtain $\widetilde{G}_m(c',s)=e^{-im_0\gamma_2}G_m(c',s)$, where 
\begin{equation}\label{11.80}
G_m(c',s):=\int_{F_v^{\times}}\varphi_v(y_v)\psi_v(y_vc')[y_v]^{-m}|y_v|_v^{-s}d^{\times}y_v.
\end{equation}
Here the factor $e^{-im_0\gamma_2}$ comes from the central character of $\pi_v'$ (see \textsection\ref{sec11.1.2}). 

\begin{comment}
\begin{multline*}
\bigg|W_{\pi_v'}\left(\begin{pmatrix}
a_v\\
&1
\end{pmatrix}w\begin{pmatrix}
1& c'\\
&1
\end{pmatrix}\right)\bigg|^2=\frac{1}{64\pi^4}\sum_{m_1\in \mathbb{Z}}[-a_v]^{-m_1}\sum_{m_2\in \mathbb{Z}}[-a_v]^{m_2}\\ 
\int_{(\alpha)}\int_{(\alpha)}\frac{\Upsilon_v(s_1,s_2,m_1,m_2)G_{m_1}(c',s_1)\overline{G_{m_2}(c',s_2)}}{|a_v|_v^{s_1+\overline{s_2}}}ds_1ds_2.
\end{multline*}
\end{comment}

Substituting \eqref{e11.80} into the definition \eqref{11.78} yields 
\begin{multline}\label{11.82}
\Psi_{n,l}^{(2)}(W_{\pi_v'};\beta_2)=\frac{\varsigma}{64\pi^4\sqrt{\langle f_{n,l},f_{n,l}\rangle}}\int_{F_v^{\times}}\int_{-\pi}^{\pi}\int_{-\pi}^{\pi}\sum_{m_1\in \mathbb{Z}}\sum_{m_2\in \mathbb{Z}}[-a_v]^{-m_1+m_2}\\ 
\int_{(\alpha)}\int_{(\alpha)}\frac{\Upsilon_v(s_1,s_2,m_1,m_2)G_{m_1}(c',s_1)\overline{G_{m_2}(c',s_2)}}{|a_v|_v^{s_1+\overline{s_2}}}ds_1ds_2e^{i(n-2l)\gamma_2}\\
(\sin\theta)^{\alpha+1} (\cos\theta)^{\beta+1}P_{l_0}^{(\alpha,\beta)}(\cos2\theta)(1+c'^2)^{\frac{1}{2}-\nu}\beta_2(\theta)d\theta d\gamma_2 |a_v|_v^{-\frac{1}{2}+\nu}[a_v]^{n_0}d^{\times}a_v,
\end{multline}
where 
\begin{align*}
\Upsilon_v(s_1,s_2,m_1,m_2):=\gamma_v(1/2-s_1,\pi_v'\otimes[\cdot]^{-m_1},\psi_v)\overline{\gamma_v(1/2-s_2,\pi_v'\otimes[\cdot]^{-m_2},\psi_v)}.
\end{align*}

Write $a_v=re^{i\phi}$. 
Substituting \eqref{11.80} into \eqref{11.82}, we derive  
\begin{multline}\label{11.83}
\Psi_{n,l}^{(2)}(W_{\pi_v'};\beta_2)=\frac{\varsigma}{64\pi^4\sqrt{\langle f_{n,l},f_{n,l}\rangle}}\sum_{m_1\in \mathbb{Z}}\sum_{m_2\in \mathbb{Z}}(-1)^{-m_1+m_2}\\ 
\int_{-\pi}^{\pi}\int_{0}^{\infty}r^{-1+2\nu}\int_{(\alpha)}\int_{(\alpha)}\frac{\Upsilon_v(s_1,s_2,m_1,m_2)G_{m_1}(c',s_1)\overline{G_{m_2}(c',s_2)}}{r^{2s_1+2\overline{s_2}}}ds_1ds_2\\
(\sin\theta)^{\alpha+1} (\cos\theta)^{\beta+1}P_{l_0}^{(\alpha,\beta)}(\cos2\theta)(1+c'^2)^{\frac{1}{2}-\nu}\beta_2(\theta)d^{\times}rd\theta   \\
\int_{-\pi}^{\pi}e^{i(n-2l)\gamma_2}d\gamma_2\int_{-\pi}^{\pi}e^{(n_0-m_1+m_2)i\phi}d\phi.
\end{multline}

Analyzing the integrals relative to $\gamma_2$ and $\phi$ in \eqref{11.83}, we deduce $\Psi_{n,l}^{(2)}(W_{\pi_v'};\beta_2)=0$ unless 
\begin{equation}\label{11.84}
\begin{cases}
n-2l=0\\
n_0-m_1+m_2=0
\end{cases}\ \Leftrightarrow\ \ \ \begin{cases}
l=n/2\\
m_1=n_0+m_2.
\end{cases}
\end{equation}

Assume the constraints \eqref{11.84} holds henceforth. From \eqref{11.83}, we obtain the following expression: 
\begin{multline}\label{11.85}
\Psi_{n,l}^{(2)}(W_{\pi_v'};\beta_2)=\frac{(-1)^{n_0}\varsigma}{16\pi^2\sqrt{\langle f_{n,l},f_{n,l}\rangle}}\int_{-\pi}^{\pi}\int_{0}^{\infty}r^{-1+2\nu}\Omega_2(r,c')d^{\times}r\\
(\sin\theta)^{\alpha+1} (\cos\theta)^{\beta+1}P_{l_0}^{(\alpha,\beta)}(\cos2\theta)(1+c'^2)^{\frac{1}{2}-\nu}\beta_2(\theta)d\theta   ,
\end{multline}
where 
\begin{align*}
\Omega_2(r,c'):=\sum_{m\in \mathbb{Z}}\int_{(\alpha')}\int_{(\alpha')}\frac{\Upsilon_v(s_1,s_2,n_0+m,m)G_{n_0+m}(c',s_1)\overline{G_{m}(c',s_2)}}{r^{2s_1+2\overline{s_2}}}ds_1ds_2.
\end{align*}

We will present two approaches to apply \eqref{11.85} in proving \eqref{11.79}, following arguments similar to those in Lemma \ref{lem11.2}.

Recall \cite[Lemma 3.1.14]{MV10}: for $A<\Re(s)<B$ and $l\geq 0$, we have 
\begin{equation}\label{11.86}
G_m(c',s)\ll \min\bigg\{\frac{(1+c'^2)^l}{(1+|s|+|m|)^l},\frac{(1+|s|+|m|)^l}{c'^{2l}},(1+|s|+|m|)^{-1+\varepsilon}\bigg\},
\end{equation}
where the implied constant relies on $A$, $B$, $\varepsilon$, $l$, and $\varphi_v$.

Let $l'\geq 0$, $m', m''\geq 0$, $\varepsilon>0$, and $|c'|\leq 10$. By applying \eqref{eq5.21} and \eqref{11.86} with the choices $\alpha'=-1/2+\Re(\nu)+\varepsilon$ and $\alpha'=m'$, respectively, together with the analysis in the proof of Lemma \ref{lem5.7} to handle the sum over $m$, we obtain  
\begin{equation}\label{11.87}
r^{l'}\frac{\partial^{l'}\Omega_2(r,c')}{\partial r^{l'}}\ll \frac{C_v(\pi')^{10\varepsilon}}{(1+|n_0|)^{m''}}\cdot\min\bigg\{\frac{C_v(\pi')^{-1+2\Re(\nu)}}{r^{-2+4\Re(\nu)+4\varepsilon}},\frac{C_v(\pi')^{2m'}}{r^{4m'}}\bigg\},
\end{equation}
where the implied constant depends on $l'$, $m'$, $m''$, $\varepsilon$ and $\varphi_v$. In particular, 
\begin{equation}\label{11.88}
\lim_{r\to 0^{+}}r^{l'}\frac{\partial^{l'}\Omega_2(r,c')}{\partial r^{l'}}=\lim_{r\to \infty}r^{l'}\frac{\partial^{l'}\Omega_2(r,c')}{\partial r^{l'}}=0.
\end{equation}

By \eqref{11.87} and \eqref{11.88}, and upon applying integration by parts to \eqref{11.85},  
\begin{multline}\label{11.89}
\Psi_{n,l}^{(2)}(W_{\pi_v'};\beta_2)\ll \frac{1}{\sqrt{\langle f_{n,l},f_{n,l}\rangle}}\cdot \frac{C_v(\pi')^{-\frac{1}{2}+\Re(\nu)+\varepsilon}}{(1+|\nu|)^{l'}(1+|n_0|)^{l'}}\\
\int_{-\pi}^{\pi}\big|(\sin\theta)^{\alpha+1} (\cos\theta)^{\beta+1}P_{l_0}^{(\alpha,\beta)}(\cos2\theta)\big|(1+c'^2)^{\frac{1}{2}-\nu}\beta_2(\theta)d\theta   .
\end{multline}

Executing Cauchy-Schwarz inequality, we obtain 
\begin{align*}
\int_{-\pi}^{\pi}\big|(\sin\theta)^{\alpha+1} (\cos\theta)^{\beta+1}P_{l_0}^{(\alpha,\beta)}(\cos2\theta)\big|(1+c'^2)^{\frac{1}{2}-\nu}\beta_2(\theta)d\theta\ll \sqrt{\langle f_{n,l},f_{n,l}\rangle}.
\end{align*}
In conjunction with \eqref{11.89} we deduce that 
\begin{equation}\label{11.90}
\Psi_{n,l}^{(2)}(W_{\pi_v'};\beta_2)\ll \frac{C_v(\pi')^{-\frac{1}{2}+\Re(\nu)+\varepsilon}}{(1+|\nu|)^{l'}(1+|n_0|)^{l'}}.
\end{equation}

Therefore, \eqref{11.79} follows from \eqref{11.90}. 
\end{proof}

\begin{lemma}\label{lemma11.17}
Let notation be as in \textsection\ref{sec11.5.1}. Suppose $n_0\geq 3$ and $l'\geq 1$. Then 
\begin{equation}\label{equ11.93}
\Psi_{n,l}^{(2)}(W_{\pi_v'};\beta_2)\ll \frac{C_v(\pi')^{-\frac{1}{2}+\Re(\nu)+\varepsilon}\log n\cdot \mathbf{1}_{l=n/2}}{n^{3/2}(1+|\nu|)^{l'-1}(1+|n_0|)^{l'}},
\end{equation}
where the implied constant depends on $l'$, $\varepsilon$, and the function $\varphi_v$ (see \textsection\ref{sec4.3}).
\end{lemma}
\begin{proof}
Let notation be as in the proof of Lemma \ref{lem11.16}. Swapping the integrals in \eqref{11.85} and integrating by parts we derive 
\begin{equation}\label{11.91}
\Psi_{n,l}^{(2)}(W_{\pi_v'};\beta_2)=\frac{(-1)^{n_0}\varsigma}{16\pi^2\sqrt{\langle f_{n,l},f_{n,l}\rangle}}\prod_{j=0}^{l'-1}\frac{1}{2\nu-1+j}\int_{0}^{\infty}r^{-1+2\nu}\Omega_2^*(r)   d^{\times}r,
\end{equation}
where $\Omega_2^*(r)$ is defined by 
\begin{align*}
\int_{-\pi}^{\pi}r^{l'}\frac{\partial^{l'
}\Omega_2(r,c')}{\partial r^{l'}}(\sin\theta)^{\alpha+1} (\cos\theta)^{\beta+1}P_{l_0}^{(\alpha,\beta)}(\cos2\theta)(1+c'^2)^{\frac{1}{2}-\nu}\beta_2(\theta)d\theta.
\end{align*}

\begin{comment}
By \cite[p.210]{MOS66} we have $P_{l_0}^{(\alpha,\beta)}(-\cos2\theta)=(-1)^{l_0}P_{l_0}^{(\beta,\alpha)}(\cos2\theta)$ and 
\begin{align*}
P_{l_0}^{(-\alpha,\beta)}(\cos2\theta)=\binom{l_0}{\alpha}^{-1}\binom{l_0+\beta}{\alpha}(-1)^{\alpha}\sin^{2\alpha}\theta P_{l_0-\alpha}^{(\alpha,\beta)}(\cos2\theta).
\end{align*}
Hence, we may assume $\beta\geq 0$ in $\Omega^*(r)$.
\end{comment} 

Under the constraints \eqref{11.84} we have $\alpha=l-l_0=n/2-l_0=n_0$, and $\beta=n-l-l_0=n/2-l_0=n_0$ (see \textsection\ref{sec11.1.2}).  

By the differentiation formula in \cite[p.213]{MOS66}, we obtain 
\begin{equation}\label{11.111}
P_{l_0}^{(n_0,n_0)}(\cos2\theta)=-\frac{1}{2(n-l_0+1)\sin\theta\cos\theta}\frac{d}{d\theta}P_{l_0+1}^{(n_0-1,n_0-1)}(\cos2\theta).
\end{equation}

Write $\Omega_{2,l'}(r,c'):=r^{l'}\frac{\partial^{l'
}\Omega_2(r,c')}{\partial r^{l'}}$. Integrating by parts, 
\begin{multline}\label{e11.90}
2(n-l_0+1)\Omega_2^*(r)=\int_{-\pi}^{\pi}\frac{P_{l_0+1}^{(n_0-1,n_0-1)}(\cos2\theta)}{{2^{n_0}}}\\
\frac{\partial}{\partial\theta}\Big[\Omega_{2,l'}(r,c')(\sin2\theta)^{n_0} (1+c'^2)^{\frac{1}{2}-\nu}\beta_2(\theta)\Big]d\theta.
\end{multline}

Therefore, it follows from \eqref{e11.90} that 
\begin{equation}\label{eq11.93}
2(n-l_0+1)\Omega_2^*(r)=S_1+S_2,	
\end{equation}
where
\begin{align*}
S_1:=&\frac{1}{2^{n_0}}\int_{-\pi}^{\pi}(\sin2\theta)^{n_0} P_{l_0+1}^{(n_0-1,n_0-1)}(\cos2\theta)\frac{\partial}{\partial\theta}\big[\Omega_{2,l'}(r,c')(1+c'^2)^{\frac{1}{2}-\nu}\beta_2(\theta)\big]d\theta,\\
S_2:=&\frac{2n_0}{2^{n_0}}\int_{-\pi}^{\pi}(\sin2\theta)^{n_0-1} P_{l_0+1}^{(n_0-1,n_0-1)}(\cos2\theta)\Omega_{2,l'}(r,c')(1+c'^2)^{\frac{1}{2}-\nu}\beta_2(\theta)\cos2\theta d\theta.
\end{align*}

Since $n_0 \geq 2$, it follows, similarly to \eqref{11.111}, that
\begin{equation}\label{11.114}
P_{l_0+1}^{(n_0-1,n_0-1)}(\cos2\theta)=-\frac{1}{2(n-l_0)\sin\theta\cos\theta}\frac{d}{d\theta}P_{l_0+2}^{(n_0-2,n_0-2)}(\cos2\theta).
\end{equation}

Substituting \eqref{11.114} into the definition of $S_1$ and $S_2$, we obtain from integrating by parts that 
\begin{equation}\label{11.115}
S_1=S_{11}+S_{12},\ \ \ S_2=S_{21}+S_{22},
\end{equation}
where 
\begin{multline*}
S_{11}:=\frac{1}{2^{n_0}(n-l_0)}\int_{-\pi}^{\pi}(\sin2\theta)^{n_0-1} P_{l_0+2}^{(n_0-2,n_0-2)}(\cos2\theta)\\
\frac{\partial^2}{\partial\theta^2}\big[\Omega_{2,l'}(r,c')(1+c'^2)^{\frac{1}{2}-\nu}\beta_2(\theta)\big]d\theta,
\end{multline*}
\begin{multline*}
S_{12}:=\frac{n_0-1}{2^{n_0-1}(n-l_0)}\int_{-\pi}^{\pi}(\sin2\theta)^{n_0-2} P_{l_0+2}^{(n_0-2,n_0-2)}(\cos2\theta)\\
\cos2\theta\frac{\partial}{\partial\theta}\big[\Omega_{2,l'}(r,c')(1+c'^2)^{\frac{1}{2}-\nu}\beta_2(\theta)\big]d\theta,
\end{multline*}
\begin{multline*}
S_{21}=\frac{n_0}{2^{n_0-1}(n-l_0)}\int_{-\pi}^{\pi}(\sin2\theta)^{n_0-2}P_{l_0+2}^{(n_0-2,n_0-2)}(\cos2\theta)\\
\frac{\partial}{\partial\theta} \big[\Omega_{2,l'}(r,c')(1+c'^2)^{\frac{1}{2}-\nu}\beta_2(\theta)\cos2\theta\big]d\theta,
\end{multline*}
\begin{multline*}
S_{22}=\frac{n_0(n_0-2)}{2^{n_0-2}(n-l_0)}\int_{-\pi}^{\pi}(\sin2\theta)^{n_0-3}P_{l_0+2}^{(n_0-2,n_0-2)}(\cos2\theta)\\
\Omega_{2,l'}(r,c')(1+c'^2)^{\frac{1}{2}-\nu}\beta_2(\theta)\cos^22\theta d\theta.
\end{multline*}

We proceed to bound $S_{ij}$, $1\leq i,j\leq 2$, as follows. 
\begin{itemize}
\item By Cauchy inequality, we have
\begin{multline}\label{11.93}
|S_{11}|\ll \bigg[\int_{-\pi}^{\pi}(\sin\theta)^{2n_0-3}(\cos\theta)^{2n_0-3} \big[P_{l_0+2}^{(n_0-2,n_0-2)}(\cos2\theta)\big]^2d\theta\bigg]^{1/2}\\
\frac{1}{n-l_0}\bigg[\int_{-\pi}^{\pi}\sin\theta\cos\theta\left(\frac{\partial^2}{\partial\theta^2}\big[\Omega_{2,l'}(r,c')(1+c'^2)^{\frac{1}{2}-\nu}\beta_2(\theta)\big]\right)^2d\theta\bigg]^{1/2}.
\end{multline}

By \cite[p.212]{MOS66} we derive 
\begin{equation}\label{e11.94}
\int_{-\pi}^{\pi}(\sin\theta)^{2n_0-3}(\cos\theta)^{2n_0-3} \big[P_{l_0+2}^{(n_0-2,n_0-2)}(\cos2\theta)\big]^2d\theta\ll \langle f_{n,l},f_{n,l}\rangle.	
\end{equation}

Similar to \eqref{11.87} we have 
\begin{multline}\label{11.95}
\bigg[\int_{-\pi}^{\pi}\sin\theta\cos\theta\left(\frac{\partial^2}{\partial\theta^2}\big[\Omega_{2,l'}(r,c')(1+c'^2)^{\frac{1}{2}-\nu}\beta_2(\theta)\big]\right)^2d\theta\bigg]^{1/2}\ll (1+|\nu|)\\
\cdot \min\bigg\{\frac{C_v(\pi')^{-1+2\Re(\nu)+2\varepsilon}}{(1+|n_0|)^{m''}r^{-2+4\Re(\nu)+4\varepsilon}},\frac{C_v(\pi')^{2m'+\varepsilon}}{(1+|n_0|)^{m''}r^{4m'}}\bigg\}.
\end{multline}

Therefore, it follows from \eqref{11.93}, \eqref{e11.94} and \eqref{11.95} that 
\begin{multline}\label{11.96}
\frac{(-1)^{n_0}\varsigma}{16\pi^2\sqrt{\langle f_{n,l},f_{n,l}\rangle}}\int_{0}^{\infty}r^{-1+2\nu}S_{11}d^{\times}r\\
\ll \frac{\log n\cdot \mathbf{1}_{l=n/2}}{(n-l_0)(1+|\nu|)^{l'-1}(1+|n_0|)^{l'}C_v(\pi')^{1/2-\Re(\nu)-\varepsilon}}.
\end{multline}

\item By \cite[inequality (20) on p. 234]{HS14} we have, for $0\leq\theta\leq \pi$, that 
\begin{equation}\label{11.97}
2^{-n_0}(\sin2\theta)^{n_0-2} P_{l_0+2}^{(n_0-2,n_0-2)}(\cos2\theta)\ll \bigg[\frac{\Gamma(l_0+n_0)\Gamma(l_0+n_0)}{\Gamma(l_0+2n_0-2)\Gamma(l_0+2)}\bigg]^{\frac{1}{2}},
\end{equation}
where the implied constant is absolute. It follows from \eqref{11.97} and the definition of $\langle f_{n,l},f_{n,l}\rangle$ (see \textsection\ref{sec11.1.2}) that 
\begin{equation}\label{11.121}
2^{-n_0}(\sin2\theta)^{n_0-2} P_{l_0+2}^{(n_0-2,n_0-2)}(\cos2\theta)\ll n^{1/2}\langle f_{n,l},f_{n,l}\rangle^{1/2}.
\end{equation}

Hence, it follows from \eqref{11.121} and \eqref{11.87} that 
\begin{multline}\label{11.98}
\frac{(-1)^{n_0}\varsigma}{16\pi^2\sqrt{\langle f_{n,l},f_{n,l}\rangle}}\int_{0}^{\infty}r^{-1+2\nu}(|S_{12}|+|S_{21}|)d^{\times}r\\
\ll \frac{n^{1/2}}{n-l_0}\cdot \frac{(1+|\nu|)C_v(\pi')^{-1/2+\Re(\nu)+\varepsilon}}{(1+|n_0|)^{l'}}.
\end{multline}

\item Let $I_n:=\{\theta\in[-\pi,\pi):\ |\sin 2\theta|\leq  n^{-10-10n_0}\}$. Then we have
\begin{equation}\label{e11.123}
S_{22}=S_{22}^{(1)}+S_{22}^{(2)},
\end{equation}
where 
\begin{multline}\label{eq11.124}
S_{22}^{(1)}:=\frac{n_0(n_0-2)}{2^{n_0-2}(n-l_0)}\int_{I_n}(\sin2\theta)^{n_0-3}P_{l_0+2}^{(n_0-2,n_0-2)}(\cos2\theta)\\
\Omega_{2,l'}(r,c')(1+c'^2)^{\frac{1}{2}-\nu}\beta_2(\theta)\cos^22\theta d\theta,
\end{multline}
and 
\begin{multline}\label{eq11.125}
S_{22}^{(2)}:=\frac{n_0(n_0-2)}{2^{n_0-2}(n-l_0)}\int_{[-\pi,\pi]-I_n}(\sin2\theta)^{n_0-3}P_{l_0+2}^{(n_0-2,n_0-2)}(\cos2\theta)\\
\Omega_{2,l'}(r,c')(1+c'^2)^{\frac{1}{2}-\nu}\beta_2(\theta)\cos^22\theta d\theta.
\end{multline}

By \cite[Theorem 7.32.1]{Sze59} and the Stirling formula, we have 
\begin{equation}\label{e11.124}
\max_{\theta\in [0,\pi]}\big|P_{l_0+2}^{(n_0-2,n_0-2)}(\cos2\theta)\big|\leq \binom{l_0+n_0}{l_0+2}\ll (l_0+2)^{10n_0}\ll n^{10n_0}.
\end{equation}

Note that $n_0\geq 3$. Substituting \eqref{e11.124} into \eqref{eq11.124} leads to 
\begin{align*}
S_{22}^{(1)}\ll \frac{n_0(n_0-2)n^{10n_0}}{n-l_0}\int_{I_n}\big|\Omega_{2,l'}(r,c')\big|(1+c'^2)^{\frac{1}{2}-\nu}\beta_2(\theta) d\theta.
\end{align*}
In conjunction with the estimate \eqref{11.87} we derive 
\begin{multline}\label{e11.125}
\frac{(-1)^{n_0}\varsigma}{16\pi^2\sqrt{\langle f_{n,l},f_{n,l}\rangle}}\int_{0}^{\infty}r^{-1+2\nu}S_{22}^{(1)}d^{\times}r\\
\ll \frac{n^{-10}}{n-l_0}\cdot \frac{(1+|\nu|)C_v(\pi')^{-1/2+\Re(\nu)+\varepsilon}}{(1+|n_0|)^{l'}}.
\end{multline}

Furthermore, plugging \eqref{11.121} into \eqref{eq11.125} yields 
\begin{multline*}
S_{22}^{(2)}\ll \frac{n_0(n_0-2)n^{1/2}\langle f_{n,l},f_{n,l}\rangle^{1/2}}{n-l_0}\int_{[-\pi,\pi]-I_n}|\sin\theta\cos\theta|^{-1}\\
\big|\Omega_{2,l'}(r,c')\big|(1+c'^2)^{\frac{1}{2}-\nu}\beta_2(\theta)\cos^22\theta d\theta.
\end{multline*}

Therefore, it follows from \eqref{11.87} that  
\begin{multline}\label{e11.126}
\frac{(-1)^{n_0}\varsigma}{16\pi^2\sqrt{\langle f_{n,l},f_{n,l}\rangle}}\int_{0}^{\infty}r^{-1+2\nu}S_{22}^{(2)}d^{\times}r\\
\ll \frac{n^{1/2}\log n^{10+10n_0}}{n-l_0}\cdot \frac{(1+|\nu|)C_v(\pi')^{-1/2+\Re(\nu)+\varepsilon}}{(1+|n_0|)^{l'}}.
\end{multline}
\end{itemize}

Therefore, \eqref{equ11.93} follows from \eqref{11.91}, \eqref{eq11.93}, \eqref{11.115}, \eqref{11.96}, \eqref{11.98}, \eqref{e11.123}, \eqref{e11.125}, and \eqref{e11.126}.  
\end{proof}

An additional integration by parts does not yield an improvement of \eqref{equ11.93} unless the estimate \eqref{11.121} is sharpened. Moreover, the condition $n_0 \geq 3$ is essential in the preceding arguments. For the case $n_0 \leq 2$, we will adopt a different approach based on the asymptotic Hilb-type formula for Jacobi polynomials given in Proposition \ref{prop11.16}. 

\begin{lemma}\label{lem11.21}
Let notation be as in \textsection\ref{sec11.5.1}. Suppose $n_0=2$ and $n\geq 1000$. Let $l'\geq 1$. Then 
\begin{equation}\label{11.104}
\Psi_{n,l}^{(2)}(W_{\pi_v'};\beta_2)\ll \frac{C_v(\pi')^{-\frac{1}{2}+\Re(\nu)+\varepsilon}\cdot \mathbf{1}_{l=n/2}}{n^{3/2}(1+|\nu|)^{l'}},
\end{equation}
where the implied constant depends on $l'$, $\varepsilon$, and the function $\varphi_v$ (see \textsection\ref{sec4.3}).
\end{lemma}
\begin{proof}
Let notation be as in the proof of Lemma \ref{lemma11.17}. Let $N=(n+1)/2$. Let $\langle f_{n,l},f_{n,l}\rangle$ be defined as in \textsection\ref{sec11.1.2}. Suppose $l=n/2$ throughout this proof. Then 
\begin{equation}\label{11.122}
	\langle f_{n,l},f_{n,l}\rangle\asymp n^{-1/2},
\end{equation}
where the implied constant is absolute.   

By Proposition \ref{prop11.16} with $m=10$, we obtain  
\begin{multline}\label{11.125}
(\sin\theta)^{n_0-1}(\cos\theta)^{n_0-1}P_{l_0+1}^{(n_0-1,n_0-1)}(\cos2\theta)\\
=(l_0+2)\sqrt{\frac{\theta}{\sin\theta\cos\theta}}
\cdot \frac{J_{1}(N\theta)}{N}+O\left(\frac{\theta}{N^{3/2}\sqrt{\sin\theta\cos\theta}}\right).
\end{multline}	

Substituting \eqref{11.125} into the definition of $S_2$ (see \eqref{eq11.93}) we derives 
\begin{equation}\label{11.126}
|S_2-S_2'|\ll \int_{-\pi}^{\pi}\frac{\big|\Omega_{2,l'}(r,c')\big|(1+c'^2)^{\frac{1}{2}-\nu}\beta_2(\theta)}{N^{3/2}\sqrt{\sin\theta\cos\theta}}d\theta,
\end{equation}
where 
\begin{equation}\label{e11.133}
S_2':=\frac{2(l_0+2)}{N}\int_{-\pi}^{\pi}\frac{\sqrt{\theta}\cos2\theta}{\sqrt{\sin\theta\cos\theta}}
\cdot J_1(N\theta)\Omega_{2,l'}(r,c')(1+c'^2)^{\frac{1}{2}-\nu}\beta_2(\theta)d\theta.
\end{equation}

As a consequence of \eqref{11.91},  \eqref{eq11.93} and \eqref{11.126}, we obtain 
\begin{equation}\label{eq11.134}
\Psi_{n,l}^{(2)}(W_{\pi_v'};\beta_2)	=\Psi_{n,l}^{(2,1)}(W_{\pi_v'};\beta_2)+\Psi_{n,l}^{(2,2)}(W_{\pi_v'};\beta_2)+\Psi_{n,l}^{(2,3)}(W_{\pi_v'};\beta_2),
\end{equation}
where 
\begin{align*}
\Psi_{n,l}^{(2,1)}(W_{\pi_v'};\beta_2):=\frac{(-1)^{n_0}\varsigma}{32\pi^2(n-l_0+1)\sqrt{\langle f_{n,l},f_{n,l}\rangle}}\prod_{j=0}^{l'-1}\frac{1}{2\nu-1+j}\int_{0}^{\infty}r^{-1+2\nu}S_1d^{\times}r,
\end{align*}
and 
\begin{multline*}
\Psi_{n,l}^{(2,2)}(W_{\pi_v'};\beta_2):=\frac{(-1)^{n_0}\varsigma}{32\pi^2(n-l_0+1)\sqrt{\langle f_{n,l},f_{n,l}\rangle}}\cdot \frac{2(l_0+2)}{N}\prod_{j=0}^{l'-1}\frac{1}{2\nu-1+j}\\
\int_{0}^{\infty}r^{-1+2\nu}\int_{-\pi}^{\pi}\frac{\sqrt{\theta}\cos2\theta}{\sqrt{\sin\theta\cos\theta}}
\cdot J_1(N\theta)\Omega_{2,l'}(r,c')(1+c'^2)^{\frac{1}{2}-\nu}\beta_2(\theta)d\theta d^{\times}r,
\end{multline*}
and the term $\Psi_{n,l}^{(2,3)}(W_{\pi_v'};\beta_2)$ satisfies 
\begin{multline}\label{equ11.135}
\Psi_{n,l}^{(2,3)}(W_{\pi_v'};\beta_2)\ll \frac{1}{(n-l_0+1)\sqrt{\langle f_{n,l},f_{n,l}\rangle}}\prod_{j=0}^{l'-1}\frac{1}{2\nu-1+j}\\
\int_{0}^{\infty}r^{-1+2\nu}\int_{-\pi}^{\pi}\frac{\big|\Omega_{2,l'}(r,c')\big|(1+c'^2)^{\frac{1}{2}-\nu}\beta_2(\theta)}{N^{3/2}\sqrt{\sin\theta\cos\theta}}d\theta d^{\times}r.
\end{multline}

By \eqref{11.115}, \eqref{11.96} and \eqref{11.98} we have
\begin{equation}\label{e11.135}
\Psi_{n,l}^{(2,1)}(W_{\pi_v'};\beta_2)\ll \frac{C_v(\pi')^{-\frac{1}{2}+\Re(\nu)+\varepsilon}\cdot \mathbf{1}_{l=n/2}}{n^{3/2}(1+|\nu|)^{l'}}. 
\end{equation}

Moreover, substuting \eqref{11.87} and \eqref{11.122} into \eqref{equ11.135} leads to 
\begin{equation}\label{e11.137}
\Psi_{n,l}^{(2,3)}(W_{\pi_v'};\beta_2)\ll \frac{C_v(\pi')^{-\frac{1}{2}+\Re(\nu)+\varepsilon}\cdot \mathbf{1}_{l=n/2}}{n^2(1+|\nu|)^{l'}}. 
\end{equation}

Now we proceed to bound $\Psi_{n,l}^{(2,2)}(W_{\pi_v'};\beta_2)$, which is more delicate due to the singular term $\sqrt{\frac{\theta}{\sin\theta\cos\theta}}$. We will make use of the properties of the Bessel function $J_1(N\theta)$ together with a decomposition as follows. 

Define $J:=\{\theta\in[-\pi,\pi):\ |\sin 2\theta|\leq  N^{-1}\}$, and $J^c:=[-\pi,\pi)-J$. Then 
\begin{equation}\label{eq11.138}
\Psi_{n,l}^{(2,2)}(W_{\pi_v'};\beta_2):=\Psi_{n,l}^{(2,4)}(W_{\pi_v'};\beta_2)+\Psi_{n,l}^{(2,5)}(W_{\pi_v'};\beta_2),
\end{equation}
where 
\begin{multline*}
\Psi_{n,l}^{(2,4)}(W_{\pi_v'};\beta_2):=\frac{(-1)^{n_0}\varsigma}{32\pi^2(n-l_0+1)\sqrt{\langle f_{n,l},f_{n,l}\rangle}}\cdot \frac{2(l_0+2)}{N}\prod_{j=0}^{l'-1}\frac{1}{2\nu-1+j}\\
\int_{0}^{\infty}r^{-1+2\nu}\int_{J}\frac{\sqrt{\theta}\cos2\theta}{\sqrt{\sin\theta\cos\theta}}
\cdot J_1(N\theta)\Omega_{2,l'}(r,c')(1+c'^2)^{\frac{1}{2}-\nu}\beta_2(\theta)d\theta d^{\times}r,
\end{multline*}
and 
\begin{multline*}
\Psi_{n,l}^{(2,5)}(W_{\pi_v'};\beta_2):=\frac{(-1)^{n_0}\varsigma}{32\pi^2(n-l_0+1)\sqrt{\langle f_{n,l},f_{n,l}\rangle}}\cdot \frac{2(l_0+2)}{N}\prod_{j=0}^{l'-1}\frac{1}{2\nu-1+j}\\
\int_{0}^{\infty}r^{-1+2\nu}\int_{J^c}\frac{\sqrt{\theta}\cos2\theta}{\sqrt{\sin\theta\cos\theta}}
\cdot J_1(N\theta)\Omega_{2,l'}(r,c')(1+c'^2)^{\frac{1}{2}-\nu}\beta_2(\theta)d\theta d^{\times}r.
\end{multline*}

Utilizing the bound $J_1(N\theta)\ll (N\theta)^{-1/2}$,  \eqref{11.87} and \eqref{11.122}, we deduce that
\begin{align*}
\Psi_{n,l}^{(2,4)}(W_{\pi_v'};\beta_2)\ll \frac{C_v(\pi')^{-\frac{1}{2}+\Re(\nu)+\varepsilon}}{n(1+|\nu|)^{l'}(1+|n_0|)^{l'}}
\int_{I}(\sin\theta\cos\theta)^{-\frac{1}{2}}(1+c'^2)^{\frac{1}{2}-\nu}\beta_2(\theta)d\theta. 
\end{align*}

Notice that $|c'|\leq 10$ if $\beta_2(\theta)\neq 0$. Therefore, the above estimate gives 
\begin{equation}\label{11.132}
\Psi_{n,l}^{(2,4)}(W_{\pi_v'};\beta_2)\ll\frac{C_v(\pi')^{-\frac{1}{2}+\Re(\nu)+\varepsilon}}{n^{3/2}(1+|\nu|)^{l'}}. 
\end{equation}

On the other hand, by \eqref{11.122} we have 
\begin{equation}\label{11.133}
\Psi_{n,l}^{(2,5)}(W_{\pi_v'};\beta_2)\ll \frac{1}{n^{1/2}(1+|\nu|)^{l'}}\bigg|\int_{J^c}f(\theta)J_1(N\theta)d\theta\bigg|,
\end{equation}
where 
\begin{align*}
f(\theta):=\frac{\sqrt{\theta}\cos2\theta}{\sqrt{\sin\theta\cos\theta}}\int_{0}^{\infty}r^{-1+2\nu}\Omega_{2,l'}(r,c')d^{\times}r(1+c'^2)^{\frac{1}{2}-\nu}\beta_2(\theta),\ \ c'=\cot\theta. 
\end{align*}

For $m\geq 1$, by taking advantage of the differential relation
\begin{align*}
\frac{d}{dx}(x^mJ_m(x))=x^mJ_{m-1}(x),
\end{align*}
we may integrate by parts: 
\begin{multline}\label{11.130}
\int_{a}^{b} f(\theta)J_{m-1}(N\theta)d\theta=N^{-1}\int_{a}^{b} \frac{df(\theta)}{d\theta}J_{m-1}(N\theta)d\theta\\
-\frac{m}{N}\int_{a}^{b}  f(\theta)J_{m-1}(N\theta)\frac{d\theta}{\theta}+\frac{f(b)J_m(Nb)-f(a)J_m(Na)}{N},
\end{multline}
where $(a,b)\subseteq (-\pi,\pi)$, and $f$ is a differentiable function. 

Substituting \eqref{11.130} into \eqref{11.133}, along with the bound $J_0(N\theta)\ll (N\theta)^{-1/2}$, 
we derive that
\begin{equation}\label{11.135}
\Psi_{n,l}^{(2,5)}(W_{\pi_v'};\beta_2)\ll \frac{C_v(\pi')^{-\frac{1}{2}+\Re(\nu)+\varepsilon}}{n^{3/2}(1+|\nu|)^{l'}}.
\end{equation}

Therefore, the estimate \eqref{11.104} follows from substituting \eqref{e11.135}, \eqref{e11.137}, \eqref{eq11.138}, \eqref{11.132} and \eqref{11.135} into \eqref{eq11.134}. 
\end{proof}

\begin{lemma}\label{lem11.22}
Let notation be as in \textsection\ref{sec11.5.1}. Suppose $n_0\in\{0,1\}$ and $n\geq 10$. Let $l'\geq 0$. Then 
\begin{equation}\label{11.136}
\Psi_{n,l}^{(2)}(W_{\pi_v'};\beta_2)\ll \frac{C_v(\pi')^{-\frac{1}{2}+\Re(\nu)+\varepsilon}\cdot \mathbf{1}_{l=n/2}}{n^{3/2}(1+|\nu|)^{l'-1}},
\end{equation}
where the implied constant depends on $l'$, $\varepsilon$, and the function $\varphi_v$ (see \textsection\ref{sec4.3}).
\end{lemma}
\begin{proof}
Following the analysis in Lemma \ref{lem11.16}, we may assume $l=n/2$. Consider the following scenarios.
\begin{itemize}
\item Suppose $n_0=0$. The constraints \eqref{11.84} implies $\alpha=\beta=0$. Therefore, \eqref{11.91} boils down to 
\begin{multline}\label{11.137}
\Psi_{n,l}^{(2)}(W_{\pi_v'};\beta_2)=\frac{\varsigma}{16\pi^2\sqrt{\langle f_{n,l},f_{n,l}\rangle}}\prod_{j=0}^{l'-1}\frac{1}{2\nu-1+j}\int_{-\pi}^{\pi}\sin\theta\cos\theta \\P_{l_0}(\cos2\theta)\int_{0}^{\infty}r^{-1+2\nu}\Omega_{2,l'}(r,c')  d^{\times}r(1+c'^2)^{\frac{1}{2}-\nu}\beta_2(\theta)d\theta, 
\end{multline}
where $P_{l_0}(\cos2\theta)=P_{l_0}^{(0,0)}(\cos2\theta)$ is the Legendre polynomial. 

By \cite[p.231]{MOS66} we have 
\begin{equation}\label{11.138}
\sin\theta\cos\theta P_{l_0}(\cos2\theta)=-\frac{1}{4(2l_0+1)}\cdot \frac{d}{d\theta}\big[P_{l_0+1}(\cos2\theta)-P_{l_0-1}(\cos2\theta)\big].
\end{equation}

Substituting \eqref{11.138} into \eqref{11.137} yields 
\begin{multline}\label{11.139}
\Psi_{n,l}^{(2)}(W_{\pi_v'};\beta_2)=\frac{\varsigma}{64\pi^2(2l_0+1)\sqrt{\langle f_{n,l},f_{n,l}\rangle}}\prod_{j=0}^{l'-1}\frac{1}{2\nu-1+j}\int_{-\pi}^{\pi}\\ \big[P_{l_0+1}(\cos2\theta)-P_{l_0-1}(\cos2\theta)\big]\frac{d}{d\theta}\bigg[\int_{0}^{\infty}r^{-1+2\nu}\Omega_{2,l'}(r,c')  d^{\times}r(1+c'^2)^{\frac{1}{2}-\nu}\beta_2(\theta)\bigg]d\theta. 
\end{multline} 

By \cite[the equation (1.8)]{FW85}, we have
\begin{equation}\label{11.140}
(2\sin\theta\cos\theta)^{1/2}P_{l_0}(\cos2\theta)=\sqrt{\theta}J_0(N\theta)+\frac{\theta\cot\theta-1}{8\theta^{1/2}N}J_1(N\theta)+O(N^{-2}\theta^2),
\end{equation}
where $N=l_0+1/2$, and the implied constant is absolute. 

Therefore, \eqref{11.136} follows from substituting \eqref{11.140} into \eqref{11.139}, along with the integration by parts \eqref{11.130} following the arguments in the proof of Lemma \ref{lem11.21}.

\item Suppose $n_0=1$. By \eqref{11.91} and  \eqref{eq11.93}, we obtain 
\begin{equation}\label{11.148}
\Psi_{n,l}^{(2)}(W_{\pi_v'};\beta_2)=\Psi_{n,l}^{(2,6)}(W_{\pi_v'};\beta_2)+\Psi_{n,l}^{(2,7)}(W_{\pi_v'};\beta_2),	
\end{equation}
where
\begin{multline*}
\Psi_{n,l}^{(2,6)}(W_{\pi_v'};\beta_2)=\frac{(-1)^{n_0}\varsigma}{32\pi^2(n-l_0+1)\sqrt{\langle f_{n,l},f_{n,l}\rangle}}\prod_{j=0}^{l'-1}\frac{1}{2\nu-1+j}\\
\int_{-\pi}^{\pi}\sin\theta\cos\theta P_{l_0+1}(\cos2\theta)\frac{d}{d\theta}\bigg[\int_{0}^{\infty}r^{-1+2\nu}\Omega_{2,l'}(r,c')d^{\times}r(1+c'^2)^{\frac{1}{2}-\nu}\beta_2(\theta)\bigg]d\theta ,
\end{multline*}
and 
\begin{multline*}
\Psi_{n,l}^{(2,7)}(W_{\pi_v'};\beta_2)=\frac{(-1)^{n_0}\varsigma}{32\pi^2(n-l_0+1)\sqrt{\langle f_{n,l},f_{n,l}\rangle}}\prod_{j=0}^{l'-1}\frac{1}{2\nu-1+j}\\
\int_{-\pi}^{\pi} P_{l_0+1}(\cos2\theta)\int_{0}^{\infty}r^{-1+2\nu}\Omega_{2,l'}(r,c')d^{\times}r(1+c'^2)^{\frac{1}{2}-\nu}\beta_2(\theta)\cos2\theta d\theta .
\end{multline*}

Notice that $\Psi_{n,l}^{(2,7)}(W_{\pi_v'}; \beta_2)$ is of a similar form to \eqref{11.139}. Moreover, by \eqref{11.138}, $\Psi_{n,l}^{(2,6)}(W_{\pi_v'}; \beta_2)$ also reduces to the form in \eqref{11.139}. Therefore, by applying the same strategy as in the treatment of \eqref{11.139}, we obtain
\begin{equation}\label{11.149}
\Psi_{n,l}^{(2,6)}(W_{\pi_v'}; \beta_2) + \Psi_{n,l}^{(2,7)}(W_{\pi_v'}; \beta_2) 
\ll \frac{C_v(\pi')^{-\frac{1}{2} + \Re(\nu) + \varepsilon}\cdot \mathbf{1}_{l=n/2}}{n^{3/2}(1+|\nu|)^{l'-1}}.
\end{equation}

Therefore, \eqref{11.136} follows from \eqref{11.148} and \eqref{11.149}. 
\end{itemize}

Combining the above discussions, we conclude the proof of Lemma \ref{lem11.22}. 
\end{proof}

\subsection{Estimates of $\Psi_{n,l}^{(1)}(W_{\pi_v'};\beta_1)$}\label{sec11.5.4}
Recall that $\Psi_{n,l}^{(1)}(W_{\pi_v'};\beta_1)$ (see \eqref{11.77} in \textsection\ref{sec11.5.1}) is defined by 
\begin{multline*}\frac{\varsigma}{\sqrt{\langle f_{n,l},f_{n,l}\rangle}}\int_{F_v^{\times}}\int_{-\pi}^{\pi}\int_{-\pi}^{\pi}\bigg|W_{\pi_v'}\left(\begin{pmatrix}
a_v\\
&1
\end{pmatrix}\begin{pmatrix}
1 & \\
c & 1
\end{pmatrix}\begin{pmatrix}
e^{i\gamma_2}& 0 \\
0 & e^{-i\gamma_2}
\end{pmatrix}\right)\bigg|^2e^{i(n-2l)\gamma_2}\\
(\sin\theta)^{\alpha+1} (\cos\theta)^{\beta+1}P_{l_0}^{(\alpha,\beta)}(\cos2\theta)(1+c^2)^{\frac{1}{2}-\nu}\beta_1(\theta)d\theta d\gamma_2|a_v|_v^{-\frac{1}{2}+\nu}[a_v]^{n_0}d^{\times}a_v,
\end{multline*}
with $c=-\tan\theta\in \mathbb{R}$ and $\beta_1(\theta)=1-\beta(\theta)$.

The main result in this subsection is the following.
\begin{prop}\label{prop11.24}
Let notation be as in \textsection\ref{sec11.5.1}. Let $l_1, l_2\geq 0$. Then 
\begin{multline*}
\Psi_{n,l}^{(1)}(W_{\pi_v'};\beta_2)\ll C_v(\pi')^{1+\varepsilon}\mathbf{1}_{l=n/2}\cdot \min\bigg\{\frac{C_v(\pi')^{l_1+l_2}}{(1+|\nu|)^{l_1}(1+|n_0|)^{l_2}},\\
\frac{\log n\cdot C_v(\pi')^{l_1+l_2}\mathbf{1}_{n_0\geq 3}}{n^{3/2}(1+|\nu|)^{l_1-1}(1+|n_0|)^{l_2}},\ \frac{C_v(\pi')^{l_1}\mathbf{1}_{n_0=2, n\geq 1000}}{n^{3/2}(1+|\nu|)^{l_1}},\ \frac{C_v(\pi')^{l_1}\mathbf{1}_{n_0\in\{0,1\}, n\geq 10}}{n^{3/2}(1+|\nu|)^{l_1}}\bigg\},
\end{multline*}
where the implied constant depends on $l_1$, $l_2$, $\varepsilon$, and the function $\varphi_v$ (see \textsection\ref{sec4.3}).
\end{prop}
\begin{proof}
By formulas \eqref{eq5.7},  \eqref{5.17.}, in conjunction with the construction of $W_{\pi_v'}$ in \textsection\ref{sec4.3},  the function  $W_{\pi_v'}\left(\begin{pmatrix}
a_v\\
&1
\end{pmatrix}\begin{pmatrix}
1& \\
c& 1
\end{pmatrix}\begin{pmatrix}
e^{i\gamma_2}& 0 \\
0 & e^{-i\gamma_2}
\end{pmatrix}\right)$ is equal to  
\begin{align*}
\frac{1}{8\pi^2 i}\sum_{m\in \mathbb{Z}}[-a_v]^{-m}
\int_{(\alpha)}\frac{\gamma_v(1/2-s,\pi_v'\otimes[\cdot]^{-m},\psi_v)H_m(c,s)}{|a_v|_v^{s}}ds,
\end{align*}
where 
\begin{align*}
H_m(c,s):=e^{-im_0\gamma_2}\int_{F_v^{\times}}W_{\pi_v'}\left(\begin{pmatrix}
y_v\\
& 1
\end{pmatrix}w\right)\psi_v(y_vc)[y_v]^{-m}|y_v|_v^{-s}d^{\times}y_v.
\end{align*}

Parallel to \eqref{11.85} we have
\begin{multline}\label{11.142}
\Psi_{n,l}^{(1)}(W_{\pi_v'};\beta_2)=\frac{(-1)^{n_0}\varsigma}{16\pi^2\sqrt{\langle f_{n,l},f_{n,l}\rangle}}\int_{-\pi}^{\pi}\int_{0}^{\infty}r^{-1+2\nu}\Omega_1(r,c)d^{\times}r\\
(\sin\theta)^{\alpha+1} (\cos\theta)^{\beta+1}P_{l_0}^{(\alpha,\beta)}(\cos2\theta)(1+c^2)^{\frac{1}{2}-\nu}\beta_1(\theta)d\theta   ,
\end{multline}
where 
\begin{align*}
\Omega_1(r,c):=\sum_{m\in \mathbb{Z}}\int_{(\alpha')}\int_{(\alpha')}\frac{\Upsilon_v(s_1,s_2,n_0+m,m)H_{n_0+m}(c,s_1)\overline{H_{m}(c,s_2)}}{r^{2s_1+2\overline{s_2}}}ds_1ds_2.
\end{align*}

Parallel to the estimate \eqref{11.86}, we obtain from the proof of \eqref{eq11.27} that 
\begin{multline}\label{11.143}
H_m(c,s)\ll C_v(\pi')^{-(1+\varepsilon)\Re(s)}\Bigg[C_v(\pi')^{-\varepsilon l}\\
+\min\left\{\left(\frac{C_v(\pi')^{1+\varepsilon}}{1+|s|+|m|}\right)^l,(1+|s|+|m|)^{-1+\varepsilon}\right\}\Bigg],	
\end{multline}
where $l\ge 0$ is an integer, and the implied constant depends only on $l$, $\varepsilon$, and the choice of $\varphi_v$.  
In particular, this yields the truncation
\begin{equation}\label{e11.152}
|s|+|m|\ll C_v(\pi')^{1+2\varepsilon}.	
\end{equation}

Consequently, \eqref{e11.152} implies that 
$H_{n_0+m}(c,s_1)\,\overline{H_m(c,s_2)}$ enforces
\begin{align*}
|s_1|+|s_2|+|m+n_0|+|m|\ll C_v(\pi')^{1+2\varepsilon},
\end{align*}
which in particular leads to $|n_0|\ll C_v(\pi')^{1+2\varepsilon}$. 

Let $l'\geq 0$, $m', m''\geq 0$, $\varepsilon>0$, and $|c|\leq 10$. By applying \eqref{eq5.21} and \eqref{11.143} with the choices $\alpha'=-1/2+\Re(\nu)-\varepsilon$ and $\alpha'=m'$, respectively, together with the analysis in the proof of Lemma \ref{lem5.7} to handle the sum over $m$, we obtain the following counterpart of \eqref{11.87}: 
\begin{equation}\label{11.144}
r^{l'}\frac{\partial^{l'}\Omega_1(r,c)}{\partial r^{l'}}\ll \frac{C_v(\pi')^{l'+m''+1+10m''\varepsilon}}{(1+|n_0|)^{m''}}\min\bigg\{r^{2-4\Re(\nu)+8\varepsilon},\frac{C_v(\pi')^{m'\varepsilon}}{r^{4m'}}\bigg\},
\end{equation}
where the implied constant depends on $l'$, $m'$, $m''$, $\varepsilon$ and $\varphi_v$.

Therefore, Proposition \ref{prop11.24} follows from the same arguments used in the proof of 
Proposition \ref{prop11.19}.  
One proceeds exactly as in the proofs of Lemmas \ref{lem11.16}, \ref{lemma11.17}, 
\ref{lem11.21}, and \ref{lem11.22}, except that the bound \eqref{11.87} is replaced by 
\eqref{11.144}.  
The main differences are that the factor 
$C_v(\pi')^{-\frac12+\Re(\nu)+\varepsilon}$ occurring in those lemmas is now replaced by 
$C_v(\pi')^{1+\varepsilon}$, and the ranges of $m$, $n_0$, $|s_1|$, and $|s_2|$ increase from 
$C_v(\pi')^{\varepsilon}$ to $C_v(\pi')^{1+2\varepsilon}$.  
An analogous phenomenon also appears in Lemmas~\ref{lem11.4} and \ref{lem11.9} in 
\textsection\ref{sec11.3}.
\end{proof}

\section{Archimedean Integrals on the Dual Side \RNum{4}}\label{sect13}

In this section we establish uniform bounds for the local period integrals on the dual side 
associated with the Eisenstein series $E_j(\cdot,s_j)$, $j=1,2$, at all Archimedean places, 
for every $W_v$ arising in the preceding sections.

\subsection{Bounds for $\Psi_v(W_v,W_{1,v},\overline{h_{2,v}})$: the Real Case \RNum{1}}
Suppose $F_v\simeq\mathbb{R}$ and $\sigma_v$ is a holomorphic discrete series of trivial central character. Let $W_v=W_n$ be the $K$-isotypical Whittaker function defined as in \eqref{whi11.2}. 
\begin{lemma}\label{lemma11.25}
Suppose $F_v\simeq\mathbb{R}$. Let $a_v\in F_v^{\times}$. Then 
\begin{equation}\label{f11.156}
\Big|W_v\left(\begin{pmatrix}
a_v\\
& 1
\end{pmatrix}\right)\Big|\ll |a_v|_v^{1/2}\langle W_v,W_v\rangle^{1/2}\log n. 	
\end{equation}	
\end{lemma}
\begin{proof}
Let $g_v=\begin{pmatrix}
a_v\\
& 1
\end{pmatrix}$. By the Fourier inversion formula (see \cite[Proof of Lemma 3.4.2]{MV10}), we have
\begin{equation}\label{e11.156}
|W_v(g_v)|^2=\zeta_v(1)^{-1}\int_{N(F_v)}\langle \sigma_v(u_vg_v)W_v,\sigma_v(g_v)W_v\rangle\theta_v(u_v)du_v.
\end{equation}

Note that $\langle \sigma_v(u_vg_v)W_v,\sigma_v(g_v)W_v=\langle \sigma_v(g_v^{-1}u_vg_v)W_v,W_v$. Making use of the change of variable $u_v\mapsto g_vu_vg_v^{-1}$, we obtain 
\begin{equation}\label{e11.157}
\frac{|W_v(g_v)|^2}{\langle W_v,W_v\rangle}=\frac{|a_v|_v}{\zeta_v(1)}\int_{F_v}\frac{\langle \sigma_v(u_v)W_v,W_v\rangle\theta_v(a_vb_v)}{\langle W_v,W_v\rangle}db_v,\ \ u_v=\begin{pmatrix}
1& b_v\\
& 1
\end{pmatrix}.
\end{equation}

Let $C>0$. We have the decomposition 
\begin{equation}\label{e11.158}
\int_{F_v}\frac{\langle \sigma_v(u_v)W_v,W_v\rangle\theta_v(a_vb_v)}{\langle W_v,W_v\rangle}db_v=\mathcal{I}_1+\mathcal{I}_2,
\end{equation}
where 
\begin{align*}
&\mathcal{I}_1:=\langle W_v,W_v\rangle^{-1}\int_{|b_v|_v\leq C}\langle \sigma_v(u_v)W_v,W_v\rangle\theta_v(a_vb_v)db_v,\\
&\mathcal{I}_2:=\langle W_v,W_v\rangle^{-1}\int_{|b_v|_v>C}\langle \sigma_v(u_v)W_v,W_v\rangle\theta_v(a_vb_v)db_v.
\end{align*}

By Lemma \ref{lem11.14} and the triangle inequality, we obtain 
\begin{equation}\label{e11.159}
\mathcal{I}_1\ll \int_{|b_v|_v\leq C}\frac{1}{\sqrt{2+b_v^2}}db_v\ll \log (1+C).
\end{equation}

Moreover, by Lemma \ref{lem11.11}, 
\begin{equation}\label{e11.160}
\mathcal{I}_2\ll \int_{|b_v|_v>C}\frac{2^m}{(b_v^2+4)^{\frac{m}{2}}}\cdot \big|P_{\frac{n-m}{2}}^{(m-1,0)}((b_v^2-4)(b_v^2+4)^{-1})\big|db_v. 
\end{equation}

For $b_v\in F_v$, we have the bound (see \cite[(7.32.2) on p. 168]{Sze59})
\begin{align*}
\big|P_{\frac{n-m}{2}}^{(m-1,0)}((b_v^2-4)(b_v^2+4)^{-1})\big|\leq \binom{\frac{n+m}{2}-1}{\frac{n-m}{2}}\ll \left(\frac{n-m}{2}\right)^{m-1}\ll \frac{n^{m-1}}{2^m}.
\end{align*}

\begin{comment}
\begin{align*}
\log\frac{(n+q)!}{n!q!}& \doteq (n+q)\log (n+q)-n\log n-q\log q+\frac{1}{2}\log\frac{n+q}{nq}\\
& \doteq n\log (1+q/n)+q\log (1+n/q)-\frac{1}{2}\log\frac{nq}{n+q}\\
&\ll q\log n
\end{align*}
\end{comment}

As a consequence, it follows from  \eqref{e11.160} that 
\begin{equation}\label{e11.161}
\mathcal{I}_2\ll \int_{|b_v|_v>C}\frac{n^{m-1}}{(b_v^2+4)^{\frac{m}{2}}}db_v\ll \frac{n^{m-1}}{C^{m-1}}.	
\end{equation}

Therefore, by taking $C=n^2$, the estimate \eqref{f11.156} follows from \eqref{e11.157}, \eqref{e11.158}, \eqref{e11.159}, and \eqref{e11.161}. 
\end{proof}

\begin{comment}

\bigskip
\bigskip
\bigskip
\bigskip
\bigskip
\bigskip
\bigskip
\bigskip
\bigskip

As a result of the triangle inequality and the $G(F_v)$-equivariance, we derive 
\begin{equation}
|W_v(g_v)|^2\leq \zeta_v(1)^{-1}\int_{N(F_v)}\big|\langle \sigma_v(g_v^{-1}u_vg_v)W_v,W_v\rangle\big|du_v.	
\end{equation}

Notice that $\dim(K_vW_v)=1$ when $F_v=\mathbb{R}$. We have the following bounds.  
\begin{itemize}
\item Suppose $F_v \simeq \mathbb{R}$ and $\sigma_v$ is tempered. Then, by \cite{CHH88}, we have 
\begin{align*}
\frac{\big|\langle \sigma_v(g_v^{-1}u_vg_v)W_v, W_v \rangle\big|}{\big|\langle W_v, W_v \rangle\big|} 
\leq \Xi(g_v^{-1}u_vg_v).
\end{align*}
\item Suppose $F_v \simeq \mathbb{R}$ and $\sigma_v$ belongs to the complementary series. For any $\varepsilon > 0$, by \cite[Theorem 2.11]{S00}, we have 
\begin{align*}
\frac{\big|\langle \sigma_v(g_v^{-1}u_vg_v)W_v, W_v \rangle\big|}{\big|\langle W_v, W_v \rangle\big|} 
\ll_{\varepsilon} \Xi(g_v^{-1}u_vg_v)^{1-\vartheta-\varepsilon}.
\end{align*}
\end{itemize}

In conjunction with Lemma \ref{lem11.14}, we derive, for $k_v\in K_v$, that 
\begin{equation}\label{11.146}
\frac{W_v(\diag(a_v,1)k_v)}{\langle W_v,W_v\rangle^{1/2}}\ll \sqrt{\int_{F_v}\left(\frac{\log(2+a_v^{-2}b_v^2)}{\sqrt{2+a_v^{-2}b_v^2}}\right)^{1-\vartheta-\varepsilon}db_v}\ll |a_v|_v^{\frac{1}{2}}.
\end{equation}
\end{comment}

\begin{lemma}\label{lem11.25}
Let $\varepsilon>0$. With the above notation, we have  
\begin{equation}\label{11.147}
\langle W_v,W_v\rangle^{-1/2}\Psi_v(W_v,W_{1,v},\overline{h_{2,v}})\ll\mathbf{C}_v^{-\frac{1}{2}+\varepsilon}\log n.  	
\end{equation} 
\end{lemma}
\begin{proof}
By the Iwasawa decomposition, 
\begin{multline}\label{11.152}
\Psi_v(W_v,W_{1,v},\overline{h_{2,v}})=\int_{F_v^{\times}}\int_{F_v^{\times}}\int_{K_v}W_{v}(\diag(a_v,1)k_v)\\
W_{1,v}(\diag(a_v,1)k_v;\overline{\theta}_v)
\overline{\Phi_{v}(\mathbf{e}_2z_vk_v)}|a|_v^{-1/2}|z_v|_vdk_vd^{\times}ad^{\times}z,
\end{multline}
where $\Phi_v$ is defined in \textsection\ref{sec4.4}. 

According to the construction of $\Phi_v(\cdot, \cdot)$ in \textsection\ref{sec4.4}, we obtain from \cite[Lemma 3.7.2]{MV10} that for \( m \geq 0 \) and \( \varepsilon > 0 \),
\begin{equation}\label{11.153}
W_{1,v}(\diag(a_v,1)k_v; \overline{\theta}_v) \ll 
\frac{\mathbf{C}_v^{1/2+\varepsilon} |a_v|_v^{1/2-\varepsilon}}{C_v(\omega\omega'^{-1})^{1/2}} \cdot 
\left[1 + \frac{|a_v|_v}{C_v(\omega\omega'^{-1})\mathbf{C}_v^{-1}}\right]^{-m},
\end{equation}
where the implied constant depends on $\varepsilon$ and $m$.
	 
Substituting \eqref{11.153} into \eqref{11.152}, along with \eqref{f11.156} and the fact that 
\begin{align*}
|W_v(\diag(a_v,1)k_v)|=|W_v(\diag(a_v,1))|,
\end{align*}  
we derive 
\begin{align*}
\frac{\Psi_v(W_v,W_{1,v},\overline{h_{2,v}})}{\langle W_v,W_v\rangle^{1/2}}\ll  \frac{\mathbf{C}_v^{-1/2+\varepsilon}\cdot\mathbf{C}_v^{1/2+\varepsilon} }{C_v(\omega\omega'^{-1})^{1/2}}\int_{|a_v|_v\ll C_v(\omega\omega'^{-1})^{1-\varepsilon}\mathbf{C}_v^{-1-\varepsilon}}|a_v|_v^{\frac{1}{2}-\varepsilon}d^{\times}a_v,
\end{align*}
which implies \eqref{11.147}. 
\end{proof}

\subsection{Bounds for $\Psi_v(W_v,W_{1,v},\overline{h_{2,v}})$: the Real Case \RNum{2}}
Suppose $F_v\simeq\mathbb{R}$ and $\sigma_v=\Ind\mu\otimes\mu^{-1}$ is a principal series. Let $W_v=W_n$ be the $K$-isotypical Whittaker function defined as in \eqref{whi11.2}.

Define the unramified character $\mu^{\circ}$ by $\mu^{\circ}(a_v)=|\mu(a_v)|$, $a_v\in F_v^{\times}$. Let $\sigma_v^{\circ}=\Ind\mu^{\circ}\otimes(\mu^{\circ})^{-1}$. Then $\sigma_v^{\circ}$ is unramified with bounded spectral parameters. Let $W_v^{\circ}$ be the unit spherical vector in the Whittaker model of $\sigma_v^{\circ}$. 

\begin{lemma}\label{lemm11.25}
Let $\varepsilon>0$. With the above notation, we have  
\begin{equation}\label{f11.166}
\langle W_v,W_v\rangle^{-1/2}\Psi_v(W_v,W_{1,v},\overline{h_{2,v}})\ll \mathbf{C}_v^{-\frac{1}{2}+\vartheta+\varepsilon}C_v(\omega\omega'^{-1})^{-\vartheta}.  	
\end{equation} 
\end{lemma}
\begin{proof}
By the Iwasawa decomposition and the triangle inequality, along with the fact that $W_{1,v}=W_{2,v}$, we derive that 
\begin{align*}
|\Psi_v(W_{n},W_{1,v},\overline{h_{2,v}})|\ll \int_{K_v}\int_{F_v^{\times}}|W_{1,v}(\diag(a_v,1)k_v)|^2\mu^{\circ}(a_v)|a_v|_v^{-\frac{1}{2}}d^{\times}a_vdk_v.	
\end{align*}

Making use of \cite[Lemma 3.4.2]{MV10} we obtain 
\begin{multline}\label{f11.157}
|\Psi_v(W_{n},W_{1,v},\overline{h_{2,v}})|\ll \int_{K_v}\int_{F_v^{\times}}\big|W_v^{\circ}(\diag(a_v,1)k_v)\big|\\
\big|W_{1,v}(\diag(a_v,1)k_v)\overline{h_{2,v}(\diag(a_v,1)k_v)}\big|\cdot |a_v|_v^{-1}d^{\times}y_vdk_v.	
\end{multline}

By \cite[Proposition 3.2.3]{MV10}, for $l\geq 0$, we have  
\begin{equation}\label{f11.168}
W_v^{\circ}(\diag(a_v,1)k_v)\ll_l \begin{cases}
|a_v|_v^{1/2-\vartheta},\ & \text{if $|a_v|_v\leq 1$},\\
|a_v|_v^{-l},\ & \text{if $|a_v|_v>1$}.
\end{cases}	
\end{equation}

Therefore, \eqref{f11.166} follows from substituting \eqref{f11.168} into \eqref{f11.157}, together with \cite[Lemma 3.7.2]{MV10} (see \eqref{11.153}).
\end{proof}

\subsection{Bounds for $\Psi_v(W_{n,l},W_{1,v},\overline{h_{2,v}})$: the Complex Case}
Suppose $F_v \simeq \mathbb{C}$. Let $f_{n,l}$ be the  section in $\sigma_v = \Ind \mu \otimes \mu^{-1}$, where $\mu$ is the character given by $\mu(z) = |z|_v^{\nu}[z]^{n_0}$ (see \textsection\ref{sec11.1.2}), and \begin{align*}
W_{n,l}(g_v) = \int_{N(F_v)} f_{n,l}(w u_v g_v) \, \overline{\theta_v(u_v)} \, du_v
\end{align*}
be the associated Whittaker function. Then 
\begin{align*}
\langle W_{n,l},W_{n,l}\rangle:=\int_{F_v^{\times}}\bigg|W_{n,l}\left(\begin{pmatrix}
y_v\\
& 1
\end{pmatrix}\right)\bigg|^2d^{\times}y_v\asymp \langle f_{n,l},f_{n,l}\rangle,
\end{align*}
where the implied constant depends only on $F_v$. 

In this subsection, we aim to establish an upper bound for $\Psi_v(W_{n,l}, W_{1,v},\overline{h_{2,v}})$. The main result is the following.
\begin{prop}\label{prop11.26}
Let $\varepsilon>0$. Then 
\begin{align*}
\frac{\Psi_v(W_{n,l}, W_{1,v},\overline{h_{2,v}})}{\langle W_{n,l},W_{n,l}\rangle^{1/2}}\ll \min\big\{n^{1/2}\mathbf{C}_v^{-1/2+\vartheta+\varepsilon}C_v(\omega\omega'^{-1})^{-\vartheta},n^{-1/2} \textbf{C}_v^{7/2+\vartheta+\varepsilon}\big\}.
\end{align*}
\end{prop}
\begin{proof}
Proposition \ref{prop11.26} follows readily from Lemmas \ref{lem11.26} and \ref{lem11.27} below.
\end{proof}

\begin{lemma}\label{lem11.26}
Let notation be as above. Then 
\begin{equation}\label{11.158}
\Psi_v(W_{n,l}, W_{1,v},\overline{h_{2,v}})\ll n^{1/2}\cdot \mathbf{C}_v^{-\frac{1}{2}+\vartheta+\varepsilon}C_v(\omega\omega'^{-1})^{-\vartheta}\langle W_{n,l},W_{n,l}\rangle^{1/2}.
\end{equation}
\end{lemma}
\begin{proof}
By \cite[Lemma 3.4.2]{MV10} and the triangle inequality, along with the fact that $W_{1,v}=W_{2,v}$, we derive that 
\begin{multline}\label{e11.154}
|\Psi_v(W_{n,l},W_{1,v},\overline{h_{2,v}})|\ll \int_{K_v}\int_{F_v^{\times}}|W_{1,v}(\diag(y_v,1)k_v)|^2\\
|y_v|_v^{-\frac{1}{2}+\Re(\nu)}d^{\times}y_v|\sin\theta\cos\theta|^{n_0} \big|P_{l_0}^{(n_0,n_0)}(\cos2\theta)\big|dk_v.
\end{multline}

Substituting \eqref{11.121} into \eqref{e11.154} leads to 
\begin{equation}\label{11.161}
\frac{|\Psi_v(W_{n,l},W_{1,v},\overline{h_{2,v}})|}{\langle W_v,W_v\rangle^{1/2}}\ll n^{1/2}\int_{K_v}\int_{F_v^{\times}}\frac{|W_{1,v}(\diag(y_v,1)k_v)|^2}{|y_v|_v^{1/2-\Re(\nu)}}d^{\times}y_vdk_v.
\end{equation}

Let $\sigma_v^{\circ}= \Ind \mu' \otimes \mu'^{-1}$, where $\mu'$ is the character given by $\mu'(z) = |z|_v^{\Re(\nu)}$. Let $W_v^{\circ}$ be the unit spherical vector in the Whittaker function in $\sigma_v^{\circ}$.  

By \cite[Lemma 3.4.2]{MV10}, we have 
\begin{equation}\label{11.162}
\int_{K_v}\int_{F_v^{\times}}\frac{|W_{1,v}(\diag(y_v,1)k_v)|^2}{|y_v|_v^{1/2-\Re(\nu)}}d^{\times}y_vdk_v\ll |\Psi_v(W_v^{\circ},W_{1,v},\overline{h_{2,v}})|.
\end{equation}

Making use of \eqref{f11.168} as in the proof of Lemma \ref{lemm11.25}, we have 
\begin{equation}\label{11.163}
|\Psi_v(W_v^{\circ},W_{1,v},\overline{h_{2,v}})|\ll \mathbf{C}_v^{-\frac{1}{2}+\vartheta+\varepsilon}C_v(\omega\omega'^{-1})^{-\vartheta}. 	
\end{equation}

Therefore, \eqref{11.158} follows from \eqref{11.161}, \eqref{11.162} and \eqref{11.163}.
\end{proof}

\begin{lemma}\label{lemma11.27}
Let $a_v\in F_v^{\times}$, $|\theta|\leq \pi$, and $|\gamma_2|\leq \pi$. Let $j\in\{0,1\}$,  $\varepsilon>0$, and $m\in \mathbb{Z}_{\geq 1}$. Then 
\begin{multline}\label{eq11.163}
\frac{\partial^j}{\partial\theta^j}W_{1,v}\left(\begin{pmatrix}
a_v\\
&1
\end{pmatrix}\begin{pmatrix}
\cos\theta & \sin\theta \\
-\sin\theta & \cos\theta
\end{pmatrix}
\begin{pmatrix}
e^{i\gamma_2}& 0 \\
0 & e^{-i\gamma_2}
\end{pmatrix}\right)\\
\ll_{m,\varepsilon} \textbf{C}_v^{\, \frac{1}{2}+j}|a_v|_v^{\frac{1}{2}}\mathbf{1}_{|a_v|_v\leq 2\mathbf{C}_v^{1+\varepsilon}}\log (2\mathbf{C}_v^{1+\varepsilon}|a_v|_v^{-1})+\frac{|a_v|_v^{\frac{1}{2}}\mathbf{C}_v^{m+\frac{1}{2}+j}}{\max\{\mathbf{C}_v,|a_v|_v\}^{(1+\varepsilon)m}}.
\end{multline}
\end{lemma}
\begin{proof}
By the definitions \eqref{f2.3} and \eqref{eq6.28}, we have 
\begin{multline}\label{11.167}
W_{1,v}\left(\begin{pmatrix}
a_v\\
&1
\end{pmatrix}\begin{pmatrix}
\cos\theta & \sin\theta \\
-\sin\theta & \cos\theta
\end{pmatrix}
\begin{pmatrix}
e^{i\gamma_2}& 0 \\
0 & e^{-i\gamma_2}
\end{pmatrix}\right)\\
=\textbf{C}_v^{\, \frac{1}{2}}|a_v|_v^{\frac{1}{2}}
\int_{F_v^{\times}}I(t_v)\overline{\omega_v'}\omega_v(t_v)d^{\times}t_v,
\end{multline}
where $I(t_v)$ refers to the integral 
\begin{align*}
\int_{F_v}
h_v(\textbf{C}_v|a_vt_v^{-1}\cos\theta-b_v\sin\theta|_v)
h_v(|a_vt_v^{-1}\sin\theta+b_v\cos\theta|_v-1)\overline{\psi}_v(b_vt_v)db_v.
\end{align*}

Recall that $h_v$ is a fixed nonnegative smooth  function on $\mathbb{R}$, supported in the interval $|t|\leq \varepsilon^*$, satisfying $h_v(t_v)\equiv 1$ when $|t|\leq \varepsilon^*/2$. Hence, 
\begin{align*}
h_v(\textbf{C}_v|a_vt_v^{-1}\cos\theta-b_v\sin\theta|_v)h_v(|a_vt_v^{-1}\sin\theta+b_v\cos\theta|_v-1)\equiv 0
\end{align*}
unless $|a_vt_v^{-1}|_v+|b_v|_v\leq 1+\varepsilon^*$, which implies that 
\begin{equation}\label{e11.165}
|t_v|_v\geq (1+\varepsilon^*)^{-1}|a_v|_v>0,\ \ |b_v|_v\leq 1+\varepsilon^*. 
\end{equation}

Let $m\geq 1$ be an integer. 
 Integrating by parts $m$ times leads to 
\begin{multline}\label{eq11.167}
I(t_v)\ll \max_{|b_v|_v\ll 1}\big|\frac{\partial^m}{\partial b_v^m}h_v(\textbf{C}_v|a_vt_v^{-1}\cos\theta-b_v\sin\theta|_v)\\
h_v(|a_vt_v^{-1}\sin\theta+b_v\cos\theta|_v-1)\big|\ll \mathbf{C}_v^{m}|t_v|_v^{-m},
\end{multline} 
where the implied constant depends on $m$ and $h_v$. 

\begin{itemize}
\item It follows from \eqref{e11.165} (which implies that $|a_v|_v\leq (1+\varepsilon^*)|t_v|_v< 2|t_v|_v$) and $I(t_v)\ll 1$  that 
\begin{equation}\label{e11.166}
\int_{|t_v|_v\leq \mathbf{C}_v^{1+\varepsilon}}I(t_v)\overline{\omega_v'}\omega_v(t_v)d^{\times}t_v\ll \mathbf{1}_{|a_v|_v\leq 2\mathbf{C}_v^{1+\varepsilon}}\log (2\mathbf{C}_v^{1+\varepsilon}|a_v|_v^{-1}). 
\end{equation}

\item As a result of \eqref{e11.165} and \eqref{eq11.167} we derive 
\begin{multline}\label{e11.167}
\int_{|t_v|_v>\mathbf{C}_v^{1+\varepsilon}}I(t_v)\overline{\omega_v'}\omega_v(t_v)d^{\times}t_v\\
\ll \int_{|t_v|_v>\max\{\mathbf{C}_v,|a_v|_v\}^{1+\varepsilon}}\mathbf{C}_v^{m}|t_v|_v^{-m}d^{\times}t_v\ll_m\mathbf{C}_v^{m}\max\{\mathbf{C}_v,|a_v|_v\}^{-(1+\varepsilon)m}.
\end{multline}
\end{itemize}
	
Combining \eqref{e11.166} with \eqref{e11.167} we obtain 
\begin{align*}
\int_{F_v^{\times}}I(t_v)\overline{\omega_v'}\omega_v(t_v)d^{\times}t_v\ll \mathbf{1}_{|a_v|_v\leq 2\mathbf{C}_v^{1+\varepsilon}}\log (2\mathbf{C}_v^{1+\varepsilon}|a_v|_v^{-1})+\frac{\mathbf{C}_v^{m}}{\max\{\mathbf{C}_v,|a_v|_v\}^{(1+\varepsilon)m}},
\end{align*}
from which the $j=0$ case of the estimate 
\eqref{eq11.163}  follows. 

Now we consider the derivative scenario. By \eqref{11.167}, 
\begin{multline}\label{11.170}
\frac{\partial}{\partial\theta}W_{1,v}\left(\begin{pmatrix}
a_v\\
&1
\end{pmatrix}\begin{pmatrix}
\cos\theta & \sin\theta \\
-\sin\theta & \cos\theta
\end{pmatrix}
\begin{pmatrix}
e^{i\gamma_2}& 0 \\
0 & e^{-i\gamma_2}
\end{pmatrix}\right)\\
=\textbf{C}_v^{\, \frac{1}{2}}|a_v|_v^{\frac{1}{2}}
\int_{F_v^{\times}}J(t_v)\overline{\omega_v'}\omega_v(t_v)d^{\times}t_v,
\end{multline}
where  
\begin{multline*}
J(t_v):=\int_{F_v}\overline{\psi}_v(b_vt_v)
\frac{\partial}{\partial\theta}\big[h_v(\textbf{C}_v|a_vt_v^{-1}\cos\theta-b_v\sin\theta|_v)\\ 
h_v(|a_vt_v^{-1}\sin\theta+b_v\cos\theta|_v-1)\big]db_v.
\end{multline*}

By \eqref{e11.165} and the triangle inequality, we obtain $J(t_v)\ll \textbf{C}_v$, where the implied constant depends on $h_v$. Moreover,  integrating by parts $m$ times leads to
\begin{multline}\label{e11.170}
J(t_v)\ll \max_{|b_v|_v\ll 1}\big|\frac{\partial^m}{\partial b_v^m}\frac{\partial}{\partial\theta}\big[h_v(\textbf{C}_v|a_vt_v^{-1}\cos\theta-b_v\sin\theta|_v)\\
h_v(|a_vt_v^{-1}\sin\theta+b_v\cos\theta|_v-1)\big]\big|\ll \mathbf{C}_v^{m+1}|t_v|_v^{-m},
\end{multline} 
where the implied constant depends on $m$ and $h_v$. 

Making use of the estimates $J(t_v)\ll \textbf{C}_v$ and \eqref{e11.170},  in parallel to \eqref{e11.166} and \eqref{e11.167}, we derive that 
\begin{align*}
\int_{F_v^{\times}}J(t_v)\overline{\omega_v'}\omega_v(t_v)d^{\times}t_v\ll\mathbf{C}_v \mathbf{1}_{|a_v|_v\leq 2\mathbf{C}_v^{1+\varepsilon}}\log\frac{2\mathbf{C}_v^{1+\varepsilon}}{|a_v|_v}+\frac{\mathbf{C}_v^{m+1}}{\max\{\mathbf{C}_v,|a_v|_v\}^{(1+\varepsilon)m}}.
\end{align*}

Substituting this into \eqref{11.170} yields the estimate 
\eqref{eq11.163} in the  $j=1$ case. 
\end{proof}

\begin{lemma}\label{lem11.27}
Suppose $n>0$, $\varepsilon>0$, and $n=2l$. Then 
\begin{equation}\label{e11.163}
\Psi_v(W_{n,l}, W_{1,v},\overline{h_{2,v}})\ll_{\varepsilon} n^{-1/2}\cdot \textbf{C}_v^{\, 7/2+\vartheta+\varepsilon}\langle W_{n,l},W_{n,l}\rangle^{1/2}.
\end{equation}
\end{lemma}
\begin{proof}
According to \cite[Lemma 3.4.2]{MV10}, 
\begin{equation}\label{11.164}
|\Psi_v(W_{n,l}, W_{1,v},\overline{h_{2,v}})|=|\Psi_v(f_{n,l},W_{1,v},\overline{h_{2,v}})|.
\end{equation} 

Utilizing the Iwasawa decomposition, we obtain 
\begin{multline}\label{11.165}
\langle W_{n,l},W_{n,l}\rangle^{-1/2}\Psi_v(f_{n,l},W_{1,v},\overline{h_{2,v}})=\frac{\varsigma}{\sqrt{\langle f_{n,l},f_{n,l}\rangle}}\int_{F_v^{\times}}\int_{-\pi}^{\pi}\\
\int_{-\pi}^{\pi}\bigg|W_{1,v}\left(\begin{pmatrix}
a_v\\
&1
\end{pmatrix}\begin{pmatrix}
\cos\theta & \sin\theta \\
-\sin\theta & \cos\theta
\end{pmatrix}
\begin{pmatrix}
e^{i\gamma_2}& 0 \\
0 & e^{-i\gamma_2}
\end{pmatrix}\right)\bigg|^2d\gamma_2\\
(\sin\theta)^{n_0+1} (\cos\theta)^{n_0+1}P_{l_0}^{(n_0,n_0)}(\cos2\theta)(1+c^2)^{\frac{1}{2}-\nu}\beta_1(\theta)d\theta |a_v|_v^{-\frac{1}{2}+\nu}[a_v]^{n_0}d^{\times}a_v.
\end{multline}

Substituting \eqref{11.111} into \eqref{11.165}, and applying integration by parts, we obtain   
\begin{equation}\label{11.166}
\langle W_{n,l},W_{n,l}\rangle^{-1/2}\Psi_v(f_{n,l},W_{1,v},\overline{h_{2,v}})=\Sigma_1+\Sigma_2,
\end{equation}
where  
\begin{multline*}
\Sigma_1:=\frac{\varsigma n_0}{4(n-l_0+1)\sqrt{\langle f_{n,l},f_{n,l}\rangle}}\int_{F_v^{\times}}\int_{-\pi}^{\pi}\int_{-\pi}^{\pi}\\
\bigg|W_{1,v}\left(\begin{pmatrix}
a_v\\
&1
\end{pmatrix}\begin{pmatrix}
\cos\theta & \sin\theta \\
-\sin\theta & \cos\theta
\end{pmatrix}
\begin{pmatrix}
e^{i\gamma_2}& 0 \\
0 & e^{-i\gamma_2}
\end{pmatrix}\right)\bigg|^2d\gamma_2\\
(\sin\theta)^{n_0-1} (\cos\theta)^{n_0-1}P_{l_0+1}^{(n_0-1,n_0-1)}(\cos2\theta)\cos2\theta d\theta |a_v|_v^{-\frac{1}{2}+\nu}[a_v]^{n_0}d^{\times}a_v,
\end{multline*}
and 
\begin{multline*}
\Sigma_2:=\frac{\varsigma}{2(n-l_0+1)\sqrt{\langle f_{n,l},f_{n,l}\rangle}}\int_{F_v^{\times}}\int_{-\pi}^{\pi}\int_{-\pi}^{\pi}\\
\frac{\partial}{\partial\theta}\bigg|W_{1,v}\left(\begin{pmatrix}
a_v\\
&1
\end{pmatrix}\begin{pmatrix}
\cos\theta & \sin\theta \\
-\sin\theta & \cos\theta
\end{pmatrix}
\begin{pmatrix}
e^{i\gamma_2}& 0 \\
0 & e^{-i\gamma_2}
\end{pmatrix}\right)\bigg|^2d\gamma_2\\
(\sin\theta)^{n_0} (\cos\theta)^{n_0}P_{l_0+1}^{(n_0-1,n_0-1)}(\cos2\theta)d\theta |a_v|_v^{-\frac{1}{2}+\nu}[a_v]^{n_0}d^{\times}a_v.
\end{multline*}

By \cite[inequality (20) on p. 234]{HS14} we have
\begin{equation}\label{11.175}
(\sin\theta)^{n_0-1} (\cos\theta)^{n_0-1}P_{l_0+1}^{(n_0-1,n_0-1)}(\cos2\theta)\ll 1.
\end{equation}

Hence, combining \eqref{11.175} with Lemma \ref{lemma11.27}, we obtain 
\begin{multline}\label{11.176}
\Sigma_1+\Sigma_2\ll \frac{1}{(n-l_0+1)\sqrt{\langle f_{n,l},f_{n,l}\rangle}}\int_{F_v^{\times}}|a_v|_v^{-\frac{1}{2}+\Re(\nu)}\\
\bigg[\textbf{C}_v^{\, \frac{3}{2}}|a_v|_v^{\frac{1}{2}}\mathbf{1}_{|a_v|_v\leq 2\mathbf{C}_v^{1+\varepsilon}}\log (2\mathbf{C}_v^{1+\varepsilon}|a_v|_v^{-1})+\frac{|a_v|_v^{\frac{1}{2}}\mathbf{C}_v^{m+\frac{3}{2}}}{\max\{\mathbf{C}_v,|a_v|_v\}^{(1+\varepsilon)m}}\bigg]^2d^{\times}a_v.
\end{multline}

Since $\langle f_{n,l},f_{n,l}\rangle\asymp n^{-1}$, it follows from \eqref{11.176} (with $m=\floor{100\varepsilon^{-1}}$) that 
\begin{equation}\label{11.177}
\Sigma_1+\Sigma_2\ll n^{-\frac{1}{2}}\bigg[\textbf{C}_v^{\, 3}\int_{|a_v|_v\leq 2\mathbf{C}_v^{1+\varepsilon}}|a_v|_v^{\frac{1}{2}+\Re(\nu)+\varepsilon}d^{\times}a_v+1\bigg]\ll \frac{\textbf{C}_v^{\, 7/2+\vartheta+\varepsilon}}{n^{1/2}}.
\end{equation}

Therefore, \eqref{e11.163} follows from \eqref{11.166} and \eqref{11.177}.
\end{proof}

\section{Majorization of the Dual Side}\label{sec14}
In this section we assemble the local results proved in 
\textsection\ref{sec8}--\textsection\ref{sect13} to obtain 
upper bounds for the amplified dual side  at $s_1=\overline{s_2}=\frac{1}{2}+it$. 

\subsection{Reduction to Moments of Central $L$-values}
Let the notation be as in \textsection\ref{sec4}, and let 
$\mathcal{M}_{\cusp}^{\du}(s_1,s_2;\mathfrak{L})$, 
$\mathcal{M}_{\Eis}^{\du,\heartsuit}(s_1,s_2;\mathfrak{L})$, and 
$\mathcal{M}_{\Res}^{\du,\heartsuit}(s_1,s_2;\mathfrak{L})$ denote the 
meromorphic functions defined in \textsection\ref{sec4.5}.

\begin{thm}\label{thm14.1}
Let $0<\varepsilon<10^{-3}$, $s_1=1/2+it$, and $s_2=1/2-it$. Then 
\begin{multline}\label{14.1}
\mathcal{M}_{\cusp}^{\du}(s_1,s_2;\mathfrak{L})+\mathcal{M}_{\Eis}^{\du,\heartsuit}(s_1,s_2;\mathfrak{L})\ll 
C_{\infty}(\pi')^{9/4+\varepsilon}\mathbf{C}_{\infty}^{-1/2+\vartheta+\varepsilon}C_{\infty}(\omega\overline{\omega}')^{-\vartheta}\\
[N_F(\mathfrak{N}),N_F(\mathfrak{M})]^{1/2+\varepsilon}N_F(\mathfrak{L})^{-1/2+\varepsilon}\sum_{\mathfrak{a}\supseteq \mathfrak{L}}|\lambda_{\pi'}(\mathfrak{a})|\int_{\substack{\sigma\in \mathcal{F}(\mathfrak{L}\mathfrak{M};\mathbf{1})\\ C_v(\sigma)\ll C_v(\pi')^{2+\varepsilon},\ v\mid\infty}}\\
L(1/2,\pi'\times\widetilde{\pi}'\times\sigma)^{1/2}\big|L(1/2,\sigma)L(1/2-2it,\sigma \times  \omega\overline{\omega}')\big|d\mu_{\sigma},
\end{multline}	
where the implied constant depends only on $F$ and $\varepsilon$. 
\end{thm}
\begin{proof}
Recall the definition: 
\begin{align*}
\mathcal{M}_{\cusp}^{\du}(s_1,s_2;\mathfrak{L})=\sum_{\sigma\in \mathcal{A}_0([G],\mathbf{1})}\sum_{\varphi\in\mathfrak{B}(\sigma)}\langle \phi_1\overline{\phi_2},\varphi\rangle\Psi(s_2,\varphi,E_1(\cdot,s_1),\overline{h_2}).\tag{\ref{f3.1}}
\end{align*}

By Watson-Ichino formula, we have
\begin{equation}\label{e7.1}
\frac{|\langle \phi_1\overline{\phi_2},\varphi\rangle|^2}{\langle \phi_1,\phi_1\rangle\langle \phi_2,\phi_2\rangle\langle \varphi,\varphi\rangle}=\Lambda(1/2,\pi'\times\widetilde{\pi}'\times\sigma)\prod_{v\leq\infty}\frac{|\mathcal{P}_v(W_v,W_{\phi_1,v},\overline{W_{\phi_2,v}})|^2}{L_v(1/2,\pi_v'\times\widetilde{\pi}_v'\times\sigma_v)},  
\end{equation}
where the local factor $|\mathcal{P}_v(W_v,W_{\phi_1,v},\overline{W_{\phi_2,v}})|^2$ is defined by 
\begin{align*}
\int_{\overline{G}(F_v)}\frac{\langle\sigma_v(g_v)W_v,W_v\rangle \langle\pi_v'(g_v)W_{\phi_1,v},W_{\phi_1,v}\rangle\overline{\langle\pi_v'(g_v)W_{\phi_2,v},W_{\phi_2,v}\rangle}}{\langle W_v,W_v\rangle\langle W_{\phi_1,v},W_{\phi_1,v}\rangle\overline{\langle W_{\phi_2,v},W_{\phi_2,v}\rangle}} dg_v.  
\end{align*}

Furthermore, we have 
\begin{align*}
|\mathcal{P}_v(W_v,W_{\phi_1,v},\overline{W_{\phi_2,v}})|^2=L_v(1/2,\pi_v'\times\widetilde{\pi}_v'\times\sigma_v)
\end{align*}
for all but finitely many places $v$, and thereby the right hand side of \eqref{e7.1} converges absolutely. 

Notice that $E_1(\cdot,s_1)$ lies in the induced representation 
$\sigma':= |\cdot|^{s_1-1/2}\boxplus \omega\overline{\omega}'|\cdot|^{1/2-s_1}$.  
Hence, by the Rankin--Selberg convolution,
\begin{equation}\label{eq14.2}
\Psi(s_2,\varphi,E_1(\cdot,s_1),\overline{h_2})=\Lambda(s_2,\sigma\times \sigma')\prod_v\frac{\Psi_v(W_v,W_{1,v},\overline{h_{2,v}})}{L_v(s_2,\sigma_v\times \sigma'_v)},	
\end{equation}

By the bound for the Ramanujan parameters (see \eqref{eq6.6} with 
$\vartheta<1/2$), each local integral 
$\Psi_v(W_v,W_{1,v},\overline{h_{2,v}})$ converges absolutely.  Moreover, 
\begin{align*}
\Psi_v(W_v,W_{1,v},\overline{h_{2,v}})=W_v(I_2)W_{1,v}(I_2)h_{2,v}(I_2)L_v(s_2,\sigma_v\times \sigma_v')
\end{align*}
for all $v<\infty$ with $v\nmid \mathfrak{N}\mathfrak{M}\mathfrak{L}$,

As a consequence of \eqref{e7.1} and \eqref{eq14.2}, we may choose a 
factorizable orthonormal basis to deduce that  
\begin{multline}\label{eq14.3}
\sum_{\varphi\in\mathfrak{B}(\sigma)}\big|\langle \phi_1\overline{\phi_2},\varphi\rangle\Psi(s_2,\varphi,E_1(\cdot,s_1),\overline{h_2})\big|=\langle \phi_1,\phi_1\rangle\langle \phi_2,\phi_2\rangle\prod_{v\mid\infty}\mathcal{Q}_v(\sigma)\\
\sqrt{L(1/2,\pi'\times\widetilde{\pi}'\times\sigma)}\cdot |L_v(s_2,\sigma_v\times (|\cdot|_v^{s_1-1/2}\boxplus \omega_v\overline{\omega}_v'|\cdot|_v^{1/2-s_1}))|\prod_{v<\infty}\mathcal{Q}_v^{\sharp}(\sigma),
\end{multline}
where for $v<\infty$,
\begin{align*}
\mathcal{Q}_v^{\sharp}(\sigma):=\sum_{W_v}\frac{|\mathcal{P}_v(W_v,W_{\phi_1,v},\overline{W_{\phi_2,v}})|}{\sqrt{|L_v(1/2,\pi_v'\times\widetilde{\pi}_v'\times\sigma_v)|}}\cdot \frac{|\Psi_v(W_v,W_{1,v},\overline{h_{2,v}})|}{|L_v(1/2,\sigma_v\times (\mathbf{1}\boxplus \omega_v\overline{\omega}_v'))|}, 
\end{align*}
and for $v\mid\infty$, 
\begin{align*}
\mathcal{Q}_v(\sigma):=\sum_{W_v}\, |\mathcal{P}_v(W_v,W_{\phi_1,v},\overline{W_{\phi_2,v}})|\cdot |\Psi_v(W_v,W_{1,v},\overline{h_{2,v}})|.
\end{align*}
Here $W_v$ ranges over an orthonormal basis of the Whittaker model of 
$\sigma_v$.  An explicit description is given in \textsection\ref{7.1.3} for 
non-Archimedean $v$, in \eqref{eq10.4} and \eqref{eq10.5} for real $v$, and 
via the Jacquet integral of the sections in \eqref{eq12.2} for complex $v$.

Let $v<\infty$ be such that $m_v>0$.  For $0\le n_v\le m_v-r_v$, where 
$r_v=r_{\sigma_v}$, write 
\begin{align*}
W_v=W_{v,n_v}=\sum_{i=0}^{n_v}\xi_{\sigma_v}(\mathfrak{p}_v^{i},\mathfrak{p}_v^{n_v})q_v^{\frac{i-n_v}{2}}\sigma_v\left(\begin{pmatrix}
			1\\
			&\varpi_v^{i}
		\end{pmatrix}\right)W_v^{\circ}.
\end{align*}
By Lemmas \ref{lem7.6}--\ref{lem7.8}, we obtain
\begin{equation}\label{14.6}
\frac{|\Psi_v(W_v,W_{1,v},\overline{h_{2,v}})|}{|L_v(1/2,\sigma_v\times (\mathbf{1}\boxplus \omega_v\overline{\omega}_v'))|}\ll q_v^{m_v/2+\varepsilon+n_v/2}\sum_{j=0}^{m_v}|\lambda_{\sigma}(\mathfrak{p}_v^j)|.
\end{equation}

Moreover, 
$\mathcal{P}_v(W_v,W_{\phi_1,v},\overline{W_{\phi_2,v}})\equiv 0$
unless $0\leq n_v\le r_v'-r_v$.  In particular, $n_v=0$ whenever $r_v=r_v'$.  
Combining \eqref{14.6} with Propositions \ref{prop9.6}, \ref{prop10.1}, and  
\ref{prop10.5}, we conclude that for $v<\infty$ with 
$v\mid\mathfrak{M}\mathfrak{N}$,
\begin{equation}\label{14.7}
\mathcal{Q}_v^{\sharp}(\sigma)\ll q_v^{(1/2+\vartheta)m_v+\varepsilon}\mathbf{1}_{r_v\leq r_v'}.
\end{equation} 

Let $v\mid\mathfrak{L}$. By Lemma \ref{lemma7.7}, noting the normalization in \eqref{eq4.1}, we obtain 
\begin{equation}\label{14.8}
\mathcal{Q}_v^{\sharp}(\sigma)\ll q_v^{-l_v'/2+\varepsilon}\mathbf{1}_{r_v\leq l_v'}\sum_{j=0}^{l_v'-r_v}|\lambda_{\pi'}(\mathfrak{p}_v^j)|. 
\end{equation}  

Let $v\mid\infty$. By definition, we have $\mathcal{P}_v(W_v,W_{\phi_1,v},\overline{W_{\phi_2,v}})=\mathcal{P}_v(W_v,W_{\pi_v'},\overline{W_{\pi_v'}})$. Consider the following scenarios. 
\begin{itemize}
\item Suppose $F_v\simeq \mathbb{R}$ and $\sigma_v$ is a principal series. By Proposition \ref{prop11.4} and Lemma \ref{lemm11.25} we obtain 
\begin{equation}\label{14.9}
\mathcal{Q}_v(\sigma)\ll C_v(\pi')^{2+\varepsilon}\mathbf{C}_v^{-\frac{1}{2}+\vartheta+\varepsilon}C_v(\omega\omega'^{-1})^{-\vartheta}\mathbf{1}_{C_v(\sigma)\ll C_v(\pi')^{2+\varepsilon}}.
\end{equation}

\item Suppose $F_v\simeq \mathbb{R}$ and $\sigma_v$ is a discrete series. By Proposition \ref{prop11.11} and Lemma \ref{lem11.25} we obtain 
\begin{equation}\label{14.10}
\mathcal{Q}_v(\sigma)\ll C_v(\pi')^{\frac{9}{4}+\varepsilon}\mathbf{C}_v^{-\frac{1}{2}+\varepsilon}\mathbf{1}_{C_v(\sigma)\ll C_v(\pi')^{2+\varepsilon}}.
\end{equation}

\item Suppose $F_v\simeq \mathbb{C}$. By Propositions \ref{prop11.18} and \ref{prop11.26} we obtain 
\begin{equation}\label{14.11}
\mathcal{Q}_v(\sigma)\ll C_v(\pi')^{1+\varepsilon}\mathbf{C}_v^{-\frac{1}{2}+\vartheta+\varepsilon}C_v(\omega\omega'^{-1})^{-\vartheta}\mathbf{1}_{C_v(\sigma)\ll C_v(\pi')^{2+\varepsilon}}.
\end{equation}
\end{itemize}

Notice that $\mathbf{C}_v \ge C_v(\omega\omega'^{-1})$.  
Therefore, combining \eqref{14.9}, \eqref{14.10}, and \eqref{14.11}, we obtain for any 
Archimedean place $v$ that
\begin{equation}\label{14.12}
\mathcal{Q}_v(\sigma)\ll C_v(\pi')^{\frac{9}{4}+\varepsilon}\mathbf{C}_v^{-\frac{1}{2}+\vartheta+\varepsilon}C_v(\omega\omega'^{-1})^{-\vartheta}\mathbf{1}_{C_v(\sigma)\ll C_v(\pi')^{2+\varepsilon}}.
\end{equation}

Substituting \eqref{14.7}, \eqref{14.8} and  \eqref{14.12} into \eqref{f3.1} yields 
\begin{multline}\label{eq14.11}
\mathcal{M}_{\cusp}^{\du}(s_1,s_2;\mathfrak{L})\ll \frac{[N_F(\mathfrak{N}),N_F(\mathfrak{M})]^{1/2+\varepsilon}}{N_F(\mathfrak{L})^{1/2-\varepsilon}}
\left(\frac{\mathbf{C}_{\infty}}{C_{\infty}(\omega\overline{\omega}')}\right)^{\vartheta}\frac{C_{\infty}(\pi')^{9/4+\varepsilon}}{\mathbf{C}_{\infty}^{1/2-\varepsilon}}\\
\sum_{\mathfrak{a}\supseteq \mathfrak{L}}|\lambda_{\pi'}(\mathfrak{a})|\sum_{\substack{\sigma\in \mathcal{F}(\mathfrak{L}\mathfrak{M};\mathbf{1})\\ C_v(\sigma)\ll C_v(\pi')^{2+\varepsilon},\ v\mid\infty}}L(1/2,\pi'\times\widetilde{\pi}'\times\sigma)^{1/2}\big|L(s_2,\sigma\times \sigma')\big|.
\end{multline}	

Likewise, following the above local analysis (in which $\sigma_v$ cannot be a 
discrete series), we obtain 
\begin{multline}\label{eq14.12}
\mathcal{M}_{\Eis}^{\du,\heartsuit}(s_1,s_2;\mathfrak{L})\ll \frac{[N_F(\mathfrak{N}),N_F(\mathfrak{M})]^{1/2+\varepsilon}}{N_F(\mathfrak{L})^{1/2-\varepsilon}}
\left(\frac{\mathbf{C}_{\infty}}{C_{\infty}(\omega\overline{\omega}')}\right)^{\vartheta}\frac{C_{\infty}(\pi')^{2+\varepsilon}}{\mathbf{C}_{\infty}^{1/2-\varepsilon}}\\
\sum_{\mathfrak{a}\supseteq \mathfrak{L}}|\lambda_{\pi'}(\mathfrak{a})|\int_{\substack{\sigma\in \mathcal{F}(\mathfrak{L}\mathfrak{M};\mathbf{1})\\ C_v(\sigma)\ll C_v(\pi')^{2+\varepsilon},\ v\mid\infty}}L(1/2,\pi'\times\widetilde{\pi}'\times\sigma)^{1/2}\big|L(s_2,\sigma\times \sigma')\big|d\mu_{\sigma}.	
\end{multline}

Therefore, \eqref{14.1} follows from \eqref{eq14.11} and \eqref{eq14.12}.
\end{proof}

\subsection{Moments Estimates and Explicit Hybrid Subconvexity}
Let $\mathfrak{q}\subseteq \mathcal{O}_F$ be an integral ideal and $C_{\infty}:=\prod_{v\mid\infty}C_v$ with each $C_v>10$. Define the integral 
\begin{equation}\label{fc14.14}
\mathcal{I}_{\pi'}(\mathfrak{q};C_{\infty}):=\int_{\substack{\sigma\in \mathcal{F}(\mathfrak{q};\mathbf{1})\\ C_v(\sigma)\leq C_v,\ v\mid\infty}}
L(1/2,\Ad\pi'\times\sigma)^{\frac{1}{2}}\big|L(1/2,\sigma)\big|^{\frac{3}{2}}d\mu_{\sigma},
\end{equation}
where $\Ad\pi'$ is the adjoint lift of $\pi'$. 

The main result of this subsection is the following bound for 
$\mathcal{I}_{\pi'}(\mathfrak{q};C_{\infty})$, which is a key input in 
controlling the right-hand side of \eqref{14.1}.
\begin{thm}\label{theorem14.2}
Suppose $\mathfrak{M}\supseteq \mathfrak{q}$ and $C_v\gg C_v(\pi')^{2+\varepsilon}$ for all $v\mid\infty$. Then  
\begin{equation}\label{14.15}
\mathcal{I}_{\pi'}(\mathfrak{q};C_{\infty})
\ll C_{\infty}^{1+\varepsilon}N_F(\mathfrak{q})^{\frac{1}{2}+\varepsilon}\cdot [N_F(\mathfrak{q}),C_{\fin}(\Ad\pi')]^{\frac{1}{2}},
\end{equation}
where the implied constant depends only on $F$ and $\varepsilon$.	
\end{thm}

\begin{remark}
Notice that \eqref{14.15} gives the average Lindel\"{o}f hypothesis when $C_{\fin}(\pi')$ is fixed.  
\end{remark}

Theorem \ref{theorem14.2} follows from several independent auxiliary results on 
moment estimates, given in \textsection\ref{sec14.2.1} and 
\textsection\ref{sec14.2.3} below.

\subsubsection{A Hybrid First Moment for $\mathrm{GL}_3\times\mathrm{GL}_2$}\label{sec14.2.1}
A special case of \cite[Theorem C]{Yan25c} yields the following result. 
\begin{thm}\label{thm14.2}
Let $F$ be a number field and $\mathfrak{q}$ an integral ideal. For each Archimedean place $v\mid\infty$, let $C_v \geq 10$, and set $C_{\infty}:=\prod_{v\mid\infty}C_v$. Let $\pi$ be a unitary self-dual cuspidal automorphic representation of $\mathrm{PGL}_3/F$. Then
\begin{multline}\label{fc1.1}
\int_{\substack{\sigma\in\mathcal{F}(\mathfrak{q};\mathbf{1})\\ C(\sigma_v)\leq C_v,\ v\mid\infty}}L(1/2,\pi\times\sigma)d\mu_{\sigma}\ll C_{\infty}^{1+\varepsilon}\cdot [N_F(\mathfrak{q}),C_{\fin}(\pi)]^{1+\varepsilon}\\
+C_{\infty}^{\frac{1}{2}+\varepsilon}[N_F(\mathfrak{q}),C_{\fin}(\pi)]^{\frac{1}{2}+\varepsilon}
C(\pi)^{\frac{1}{4}+\varepsilon}\prod_{v\mid\infty}\Big[1+ (C_v^{-1}C(\pi_v))^2\Big],
\end{multline}
where the implied constant depends only on $F$ and $\varepsilon$. Here $[\cdot,\cdot]$ refers to the least common multiple.	
\end{thm}

\subsubsection{The Fourth Moment for $\mathrm{GL}_2$}\label{sec14.2.3}
We also recall the fourth-moment estimate established in 
\cite[Theorem D \& E]{Yan26b}.
\begin{thm}\label{thm14.4}
Let $\mathfrak{q}\subseteq \mathcal{O}_F$ be an integral ideal, and $\omega$ be a unitary character of $F^{\times}\backslash\mathbb{A}_F^{\times}$. For $v\mid\infty$, let $C_v>1$ be a constant, and $C_{\infty}:=\prod_{v\mid\infty}C_v$. Then 
\begin{equation}\label{10.1}
\int_{\substack{\pi\in \mathcal{F}(\mathfrak{q},\omega)\\
C_v(\pi)\leq C_v,\ v\mid\infty}}|L(1/2,\pi)|^4d\mu_{\pi}\ll C_{\infty}^{1+\varepsilon}N_F(\mathfrak{q})^{1+\varepsilon},
\end{equation}
where the implied constant depends only on $\varepsilon$ and $F$. 
\end{thm}

\subsubsection{Proof of Theorem \ref{theorem14.2}}
By H\"{o}lder's inequality, we obtain 
\begin{multline}\label{14.17}
\mathcal{I}_{\pi'}(\mathfrak{q};C_{\infty})\ll \Bigg[\int_{\substack{\sigma\in \mathcal{F}(\mathfrak{q};\mathbf{1})\\ C_v(\sigma)\leq C_v,\ v\mid\infty}}
L(1/2,\Ad\pi'\times\sigma)d\mu_{\sigma}\Bigg]^{\frac{1}{2}}\\
\Bigg[\int_{\substack{\sigma\in \mathcal{F}(\mathfrak{q};\mathbf{1})\\ C_v(\sigma)\leq C_v,\ v\mid\infty}}
\big|L(1/2,\sigma)\big|^{4}d\mu_{\sigma}\Bigg]^{\frac{3}{8}}\Bigg[\int_{\substack{\sigma\in \mathcal{F}(\mathfrak{1};\mathbf{1})\\ C_v(\sigma)\leq C_v,\ v\mid\infty}}
1d\mu_{\sigma}\Bigg]^{\frac{1}{8}}.
\end{multline}

Notice that $C_v(\Ad\pi_v')\ll C_v(\pi_v')^{2+\varepsilon}\ll C_v$ for all 
$v\mid\infty$.  
Substituting Theorems \ref{thm14.2} and \ref{thm14.4} into \eqref{14.17}, we 
therefore obtain \eqref{14.15}.

\subsubsection{Hybrid Subconvexity for $\mathrm{GL}_2\times \mathrm{GL}_1$}\label{14.2.4}
Finally, we recall the hybrid subconvexity bound established in \cite[Corollary 1.10]{Yan26}.
\begin{thm}\label{thm14.3}
Let $\pi$ be a unitary generic automorphic representation of 
$\mathrm{GL}_2/F$ with central character $\omega$, and let $\chi$ be a unitary 
Hecke character. Then  
\begin{multline*}
L(1/2,\pi\times\chi)\ll  [C_{\fin}(\pi),C_{\fin}(\omega)C_{\fin}(\chi)]^{\frac{1}{4}+\varepsilon}C_{\fin}(\chi)^{\frac{1}{8}+\varepsilon}C_{\infty}(\pi)^{\frac{3}{8}+\varepsilon}C_{\infty}(\pi\otimes\chi)^{\frac{3}{16}+\varepsilon}\\
+C(\pi)^{\frac{1}{2}+\varepsilon}C_{\fin}(\chi)^{\frac{1}{4}+\varepsilon}C_{\infty}(\pi\otimes\chi)^{\frac{1}{8}+\varepsilon},
\end{multline*}
where the implied constant depends on $\varepsilon$ and $F$.  
\end{thm}

\subsection{Bounds for the Amplified Dual Side}
Let $\boldsymbol{\alpha}$, $\boldsymbol{\ell}$, and $\mathfrak{L}$ be the 
automorphic data defined in \textsection\ref{sec4.2}, and let 
$\mathcal{M}_{\cusp}^{\du}(s_1,s_2;\boldsymbol{\alpha},\boldsymbol{\ell})$ and 
$\mathcal{M}_{\Eis}^{\du,\heartsuit}(s_1,s_2;\boldsymbol{\alpha},\boldsymbol{\ell})$ 
be as defined in \textsection\ref{sec4.6}.  
We then have the following estimate.
  
\begin{prop}\label{prop14.7}
Let notation be as above. Then 
\begin{multline}\label{14.25}
\mathcal{M}_{\cusp}^{\du}(s_1,s_2;\boldsymbol{\alpha},\boldsymbol{\ell})+\mathcal{M}_{\Eis}^{\du,\heartsuit}(s_1,s_2;\boldsymbol{\alpha},\boldsymbol{\ell})\ll \mathbf{C}_{\infty}^{-\frac{1}{2}+\vartheta+\varepsilon}C_{\infty}(\omega\overline{\omega}')^{-\vartheta}\\
[N_F(\mathfrak{N}),N_F(\mathfrak{M})]^{\frac{1}{2}+\varepsilon}L^{4+\varepsilon}N_F(\mathfrak{M})C_{\infty}(\pi')^{\frac{43}{8}+\varepsilon}C_{\fin}(\Ad\pi')^{\frac{1}{2}}\\
\Big[L^2N_F(\mathfrak{M})^{\frac{1}{2}}C(\omega\overline{\omega}'|\cdot|^{-2i t})^{\frac{1}{4}}C_{\infty}(\pi')^{\frac{1}{8}}
+LN_F(\mathfrak{M})^{\frac{1}{4}}C(\omega\overline{\omega}'|\cdot|^{-2i t})^{\frac{3}{8}}\Big].
\end{multline}	
\end{prop}
\begin{proof}
Taking $\chi=\omega\overline{\omega}'|\cdot|^{-2i t}$ into  Theorem \ref{thm14.3} gives 
\begin{multline}\label{14.16}
L(1/2-2it,\sigma \times  \omega\overline{\omega}')\ll 
C(\sigma)^{\frac{1}{2}+\varepsilon}C_{\fin}(\omega\overline{\omega}')^{\frac{1}{4}+\varepsilon}C_{\infty}(\sigma\otimes\omega\overline{\omega}'|\cdot|^{-2i t})^{\frac{1}{8}+\varepsilon}\\
+[C_{\fin}(\sigma),C_{\fin}(\omega\overline{\omega}')]^{\frac{1}{4}+\varepsilon}C_{\fin}(\omega\overline{\omega}')^{\frac{1}{8}+\varepsilon}C_{\infty}(\sigma)^{\frac{3}{8}+\varepsilon}C_{\infty}(\sigma\otimes\omega\overline{\omega}'|\cdot|^{-2i t})^{\frac{3}{16}+\varepsilon}.
\end{multline}

By definition we have
\begin{equation}\label{equ14.17}
C_{\infty}(\sigma\otimes\omega\overline{\omega}'|\cdot|^{-2i t})\ll C_{\infty}(\sigma\otimes|\cdot|^{-2i t})C_{\infty}(\omega\overline{\omega}')^2. 
\end{equation}

Let $\sigma\in \mathcal{F}(\mathfrak{L}\mathfrak{M};\mathbf{1})$ satisfy $ C_v(\sigma)\ll C_v(\pi')^{2+\varepsilon}$ for all Archimedean places $ v\mid\infty$. Substituting \eqref{equ14.17} into \eqref{14.16} then gives 
\begin{multline}\label{eq14.18}
L(1/2-2it,\sigma \times  \omega\overline{\omega}')\ll 
N_F(\mathfrak{L}\mathfrak{M})^{\frac{1}{2}+\varepsilon}C(\omega\overline{\omega}'|\cdot|^{-2i t})^{\frac{1}{4}+\varepsilon}C_{\infty}(\pi')^{\frac{5}{4}+\varepsilon}\\
+N_F(\mathfrak{L}\mathfrak{M})^{\frac{1}{4}+\varepsilon}C(\omega\overline{\omega}'|\cdot|^{-2i t})^{\frac{3}{8}+\varepsilon}C_{\infty}(\pi')^{\frac{9}{8}+\varepsilon}.
\end{multline}

\begin{comment}

\begin{multline*}
L(1/2-2it,\sigma \times  \omega\overline{\omega}')\ll 
C(\sigma)^{\frac{1}{2}+\varepsilon}C_{\fin}(\omega\overline{\omega}')^{\frac{1}{4}+\varepsilon}C_{\infty}(\sigma\otimes\omega\overline{\omega}'|\cdot|^{-2i t})^{\frac{1}{8}+\varepsilon}\\
+[C_{\fin}(\sigma),C_{\fin}(\omega\overline{\omega}')]^{\frac{1}{4}+\varepsilon}C_{\fin}(\omega\overline{\omega}')^{\frac{1}{8}+\varepsilon}C_{\infty}(\sigma)^{\frac{3}{8}+\varepsilon}C_{\infty}(\sigma\otimes\omega\overline{\omega}'|\cdot|^{-2i t})^{\frac{3}{16}+\varepsilon}.
\end{multline*}

\begin{multline*}
L(1/2-2it,\sigma \times  \omega\overline{\omega}')\ll 
N_F(\mathfrak{L}\mathfrak{M})^{\frac{1}{2}+\varepsilon}C_{\fin}(\omega\overline{\omega}')^{\frac{1}{4}+\varepsilon}C_{\infty}(\pi')^{1+\varepsilon}C_{\infty}(\sigma\otimes\omega\overline{\omega}'|\cdot|^{-2i t})^{\frac{1}{8}+\varepsilon}\\
+N_F(\mathfrak{L}\mathfrak{M})^{\frac{1}{4}+\varepsilon}C_{\fin}(\omega\overline{\omega}')^{\frac{3}{8}+\varepsilon}C_{\infty}(\pi')^{\frac{3}{4}+\varepsilon}C_{\infty}(\sigma\otimes\omega\overline{\omega}'|\cdot|^{-2i t})^{\frac{3}{16}+\varepsilon}.
\end{multline*}

\end{comment}

According to the definition in \eqref{equ4.2}, we obtain the bound 
\begin{equation}\label{eq14.23}
\sum_{\mathfrak{a}\supseteq \mathfrak{L}}|\lambda_{\pi'}(\mathfrak{a})|\ll L^{100\varepsilon}.	
\end{equation}

Inserting \eqref{eq14.23} and \eqref{eq14.18} into Theorem \ref{thm14.1} yields
\begin{multline}\label{eq14.17}
\mathcal{M}_{\cusp}^{\du}(s_1,s_2;\mathfrak{L})+\mathcal{M}_{\Eis}^{\du,\heartsuit}(s_1,s_2;\mathfrak{L})\ll 
C_{\infty}(\pi')^{\frac{27}{8}+\varepsilon}\mathbf{C}_{\infty}^{-\frac{1}{2}+\vartheta+\varepsilon}C_{\infty}(\omega\overline{\omega}')^{-\vartheta}\\
[N_F(\mathfrak{N}),N_F(\mathfrak{M})]^{\frac{1}{2}+\varepsilon}N_F(\mathfrak{L})^{-\frac{1}{2}+\varepsilon}
\Big[N_F(\mathfrak{L}\mathfrak{M})^{\frac{1}{2}}C(\omega\overline{\omega}'|\cdot|^{-2i t})^{\frac{1}{4}}C_{\infty}(\pi')^{\frac{1}{8}}\\
+N_F(\mathfrak{L}\mathfrak{M})^{\frac{1}{4}}C(\omega\overline{\omega}'|\cdot|^{-2i t})^{\frac{3}{8}}\Big]\cdot \mathcal{I}_{\pi'}(\mathfrak{L}\mathfrak{M};C_{\infty}(\pi')^{2+\varepsilon}),
\end{multline}
where $\mathcal{I}_{\pi'}(\mathfrak{L}\mathfrak{M}; C_{\infty}(\pi')^{2+\varepsilon})$ is defined in \eqref{fc14.14}. 
We also note that the exponent $\varepsilon$ appearing in \eqref{eq14.17} may be taken as $1000\varepsilon'$, where $\varepsilon'$ denotes the parameter in \eqref{eq14.23}. 
For notational simplicity, we continue to write both simply as $\varepsilon$, indicating an arbitrarily small constant. 
\begin{comment}
\mathcal{M}_{\cusp}^{\du}(s_1,s_2;\mathfrak{L})+\mathcal{M}_{\Eis}^{\du,\heartsuit}(s_1,s_2;\mathfrak{L})\ll 
C_{\infty}(\pi')^{9/4+\varepsilon}\mathbf{C}_{\infty}^{-1/2+\vartheta+\varepsilon}C_{\infty}(\omega\overline{\omega}')^{-\vartheta}\\
[N_F(\mathfrak{N}),N_F(\mathfrak{M})]^{1/2+\varepsilon}N_F(\mathfrak{L})^{-1/2+\varepsilon}\\
\Big[N_F(\mathfrak{L}\mathfrak{M})^{\frac{1}{2}+\varepsilon}C(\omega\overline{\omega}'|\cdot|^{-2i t})^{\frac{1}{4}+\varepsilon}C_{\infty}(\pi')^{\frac{5}{4}+\varepsilon}\\
+N_F(\mathfrak{L}\mathfrak{M})^{\frac{1}{4}+\varepsilon}C(\omega\overline{\omega}'|\cdot|^{-2i t})^{\frac{3}{8}+\varepsilon}C_{\infty}(\pi')^{\frac{9}{8}+\varepsilon}\Big]\cdot \mathcal{I}_{\pi'}(\mathfrak{L}\mathfrak{M};C_{\infty}(\pi')^{2+\varepsilon}).
\end{comment}

By Theorem \ref{theorem14.2}, we obtain 
\begin{equation}\label{14.18}
\mathcal{I}_{\pi'}(\mathfrak{L}\mathfrak{M};C_{\infty}(\pi')^{2+\varepsilon})
\ll C_{\infty}(\pi')^{2+\varepsilon}N_F(\mathfrak{L}\mathfrak{M})^{\frac{1}{2}+\varepsilon}\cdot [N_F(\mathfrak{L}\mathfrak{M}),C_{\fin}(\Ad\pi')]^{\frac{1}{2}},
\end{equation}
where the implied constant depends only on $F$ and $\varepsilon$. 

Therefore, it follows from \eqref{eq14.17} and \eqref{14.18} that 
\begin{multline}\label{14.24}
\mathcal{M}_{\cusp}^{\du}(s_1,s_2;\mathfrak{L})+\mathcal{M}_{\Eis}^{\du,\heartsuit}(s_1,s_2;\mathfrak{L})\ll 
C_{\infty}(\pi')^{\frac{43}{8}+\varepsilon}\mathbf{C}_{\infty}^{-\frac{1}{2}+\vartheta+\varepsilon}C_{\infty}(\omega\overline{\omega}')^{-\vartheta}\\
[N_F(\mathfrak{N}),N_F(\mathfrak{M})]^{\frac{1}{2}+\varepsilon}N_F(\mathfrak{L})^{\frac{1}{2}+\varepsilon}N_F(\mathfrak{M})C_{\fin}(\Ad\pi')^{\frac{1}{2}}\\
\Big[N_F(\mathfrak{L}\mathfrak{M})^{\frac{1}{2}}C(\omega\overline{\omega}'|\cdot|^{-2i t})^{\frac{1}{4}}C_{\infty}(\pi')^{\frac{1}{8}}
+N_F(\mathfrak{L}\mathfrak{M})^{\frac{1}{4}}C(\omega\overline{\omega}'|\cdot|^{-2i t})^{\frac{3}{8}}\Big].
\end{multline}

Therefore, \eqref{14.25} follows from the definition (see \textsection\ref{sec4.6}).
\end{proof}

\subsection{Bounds for the Amplified Residual Side}
Let $\varepsilon_0$ be the small parameter introduced in 
\textsection\ref{sect4.4}. Recall the definition \eqref{eq4.5}:
\begin{align*}
\mathcal{M}_{\Res}^{\du,\heartsuit,\Reg}(s_1,s_2;\boldsymbol{\alpha},\boldsymbol{\ell}):=\frac{1}{2\pi i}\int_{|s|=\varepsilon_0}\frac{\mathcal{M}_{\Res}^{\du,\heartsuit}(s_1-s,s_2+2s;\boldsymbol{\alpha},\boldsymbol{\ell})}{s}ds.
\end{align*}

\begin{prop}\label{prop14.8}
Let $0<\varepsilon<10^{-3}$, $s_1=1/2+it$, and $s_2=1/2-it$. Then
\begin{align*}
\mathcal{M}_{\Res}^{\du,\heartsuit,\Reg}(s_1,s_2;\boldsymbol{\alpha},\boldsymbol{\ell})\ll C_{\infty}(\pi')^{\frac{19}{4}+\varepsilon}\mathbf{C}_{\infty}^{\varepsilon}L^{1+\varepsilon}N_F(\mathfrak{M})^{\frac{3}{2}+\varepsilon}[N_F(\mathfrak{N}),N_F(\mathfrak{M})]^{1+\varepsilon}.
\end{align*}
\end{prop}

\subsubsection{Local Estimates}
Let $\xi\in \widehat{F^{\times}\backslash\mathbb{A}_F^{(1)}}$ and $f\in \mathfrak{B}(\xi,\xi^{-1})$. Recall that 
\begin{align*}
\mathfrak{F}_{\xi}(s,s_1,s_2):=\frac{1}{2}\sum_{f\in \mathfrak{B}(\xi,\xi^{-1})}\langle \phi_1\overline{\phi_2},E(\cdot,f,-\overline{s})\rangle \Psi(s_2,E(\cdot,f,s),E_1(\cdot,s_1),\overline{h_2})\tag{\ref{f}}.
\end{align*}

Define the $L$-function
\begin{multline*}
\mathcal{L}(s,s_1,s_2;\xi):=\frac{L(1/2-s,\pi'\times\overline{\pi}'\otimes \overline{\xi})}{L(1-2s,\xi^{-2})L(1+2s,\xi^2)}\cdot L(s_2+s_1+s-1/2,\xi)\\
L(s_2-s_1+s+1/2,\xi\omega\overline{\omega}')L(s_2+s_1-s-1/2, \overline{\xi})L(s_2-s_1-s+1/2, \overline{\xi}\omega\overline{\omega}').
\end{multline*}

For a non-Archimedean place $v$, we let $\mathcal{L}_v(s,s_1,s_2;\xi)$ be the $v$-th component of $\mathcal{L}(s,s_1,s_2;\xi)$. 

Suppose that $\Re(s)<-1$. In this region the Eisenstein series
$E(\cdot,f,-\overline{s})$ converges absolutely. Unfolding 
$E(\cdot,f,-\overline{s})$ in the inner product 
$\langle \phi_1\overline{\phi_2}, E(\cdot,f,-\overline{s})\rangle$, we obtain
\begin{align*}
\langle \phi_1\overline{\phi_2},E(\cdot,f,-\overline{s})\rangle=\prod_{v\leq\infty}\mathcal{P}_v(f,W_{1,v},\overline{W_{2,v}}),
\end{align*}
where for each $v\leq\infty$, 
\begin{align*}
\mathcal{P}_v(f,W_{1,v},\overline{W_{2,v}}):=\int_{N(F_v)Z(F_v)\backslash G(F_v)}W_{\phi_1,v}(x_v;\theta)\overline{W_{\phi_2,v}(x_v;\theta)}\overline{S(f_v)(x_v,-\overline{s})}dx_v.
\end{align*}

In parallel with \eqref{eq14.2}, we have the identity
\begin{align*}
\Psi(s_2,E(\cdot,f,s),E_1(\cdot,s_1),\overline{h_2})=\Lambda(s_2,\sigma\times \sigma')\prod_{v\leq\infty}\frac{\Psi_v(W_v,W_{1,v},\overline{h_{2,v}})}{L_v(s_2,\sigma_v\times \sigma'_v)},
\end{align*}
where $W_v$ is the Jacquet integral associated with the section 
$f\in \mathfrak{B}(\xi,\xi^{-1})$.

Let $\Re(s)<-1$. We have the decomposition 
\begin{equation}\label{14.27}
\mathfrak{F}_{\xi}(s,s_1,s_2)=2^{-1}\mathcal{L}(s,s_1,s_2;\xi)\mathcal{Q}(s),
\end{equation}
where $\mathcal{Q}(s):=\prod_{v\mid\infty}\mathcal{Q}_v(s)\; \prod_{v<\infty}\mathcal{Q}_v^{\sharp}(s)$, with 
\begin{align*}
&\mathcal{Q}_v^{\sharp}(s):=\mathcal{L}_v(s,s_1,s_2;\xi)^{-1}\sum_{W_v}\mathcal{P}_v(f,W_{1,v},\overline{W_{2,v}})\Psi_v(W_v,W_{1,v},\overline{h_{2,v}}),\\
&\mathcal{Q}_v(s):=\sum_{W_v}\mathcal{P}_v(f,W_{1,v},\overline{W_{2,v}})\Psi_v(W_v,W_{1,v},\overline{h_{2,v}}).
\end{align*}

By Lemma \ref{lem3.6}, for every finite place $v<\infty$ the factor 
$\mathcal{Q}_v^{\sharp}(s)$ admits a holomorphic continuation to all 
$s\in\mathbb{C}$. Moreover, for all $v<\infty$ with $v\nmid \mathfrak{N}\mathfrak{M}\mathfrak{L}$,
\begin{align*}
\mathcal{Q}_v^{\sharp}(s)=W_{\phi_1,v}(I_2)W_{1,v}(I_2)\overline{W_{\phi_2,v}(I_2)}\overline{W_{2,v}(I_2)}.
\end{align*}

Furthermore, by the estimates established in 
\textsection\ref{sec11}--\textsection\ref{sect13}, the local factor 
$\mathcal{Q}_v(s)$ converges absolutely for each $v\mid\infty$ whenever 
$|\Re(s)|\leq 2$. It follows that $\mathcal{Q}(s)$ 
converges absolutely in the region $|\Re(s)|\leq 2$, and therefore defines a 
holomorphic function there. Consequently, \eqref{14.27} extends to an identity 
of meromorphic functions throughout the region $|\Re(s)|\leq 2$.

\begin{lemma}\label{lem14.9}
Let $|\Re(s)|\leq 2$ and $m\geq 0$. Then 
\begin{multline}\label{eq14.29}
\mathcal{Q}(s)\ll C_{\infty}(\pi')^{\frac{9}{4}+\varepsilon}\mathbf{C}_{\infty}^{-\frac{1}{2}+|\Re(s)|+\varepsilon}C_{\infty}(\omega\omega'^{-1})^{-|\Re(s)|}N_F(\mathfrak{L})^{-\frac{1}{2}+\varepsilon}\\
[N_F(\mathfrak{N}),N_F(\mathfrak{M})]^{\frac{1}{2}+|\Re(s)|+\varepsilon}\sum_{\mathfrak{a}\supseteq \mathfrak{L}}|\lambda_{\pi'}(\mathfrak{a})|\prod_{v\mid\infty}\bigg[\frac{C_v(\pi')^{1+\varepsilon}}{C_v(\xi|\cdot|^s)}\bigg]^m\prod_{v\mid\mathfrak{L}\mathfrak{M}}\mathbf{1}_{2r_{\xi_v}\leq r_v'+l_v'}.
\end{multline}
\end{lemma}
\begin{proof}
Let $|\Re(s_1)-1/2|<10^{-1}$ and $|\Re(s_2)-1/2|<10^{-1}$. Throughout this subsection we suppose that $|\Re(s)|\leq 2$.  
In parallel with \eqref{14.7}, we obtain the local bound
\begin{equation}\label{14.28}
\mathcal{Q}_v^{\sharp}(s)\ll q_v^{(1/2+|\Re(s)|)m_v+\varepsilon}\mathbf{1}_{2r_{\xi_v}\leq r_v'}.	
\end{equation}

For places $v\mid\mathfrak{L}$, Lemma \ref{lemma7.7} yields the analogue of 
\eqref{14.8}, namely 
\begin{equation}\label{14.29}
\mathcal{Q}_v^{\sharp}(s)\ll q_v^{-l_v'/2+\varepsilon}\mathbf{1}_{2r_{\xi_v}\leq l_v'}\sum_{j=0}^{l_v'-r_v}|\lambda_{\pi'}(\mathfrak{p}_v^j)|. 
\end{equation}  

Finally, for each archimedean place $v\mid\infty$, an argument parallel to the 
one leading to \eqref{14.12} yields  
\begin{equation}\label{14.30}
\mathcal{Q}_v(s)\ll C_v(\pi')^{\frac{9}{4}+\varepsilon}\mathbf{C}_v^{-\frac{1}{2}+|\Re(s)|+\varepsilon}C_v(\omega\omega'^{-1})^{-|\Re(s)|}\bigg[\frac{C_v(\pi')^{1+\varepsilon}}{C_v(\xi|\cdot|^s)}\bigg]^m,
\end{equation}
where the implied constant depends on $m$.  
Observe that, for each $v\mid\infty$, the factor $(C_v(\pi')^{1+\varepsilon}C_v(\xi|\cdot|^s)^{-1})^m$ serves effectively as a truncation 
$\mathbf{1}_{C_v(\xi|\cdot|^s)\ll C_v(\pi')^{1+\varepsilon}}$.

Therefore, \eqref{eq14.29} follows from \eqref{14.28}, \eqref{14.29} and \eqref{14.30}. 
\end{proof}

\subsubsection{Proof of Proposition \ref{prop14.8}}
According to the definition \eqref{3.27} we have 
\begin{multline}\label{e14.32}
\mathcal{M}_{\Res}^{\du,\heartsuit}(s_1,s_2;\mathfrak{L}):=\underset{s=s_1-s_2+\frac{1}{2}}{\Res}\  \mathfrak{F}_{\overline{\omega}\omega'}(s,s_1,s_2)-\underset{s=s_2-s_1-\frac{1}{2}}{\Res}\  \mathfrak{F}_{\omega\overline{\omega}'}(s,s_1,s_2)\\
+\underset{s=\frac{3}{2}-s_1-s_2}{\Res}\  \mathfrak{F}_{\mathbf{1}}(s,s_1,s_2)-\underset{s=s_1+s_2-\frac{3}{2}}{\Res}\  \mathfrak{F}_{\mathbf{1}}(s,s_1,s_2).
\end{multline}

Let $s_1=1/2+it$ and $s_2=1/2-it$. By \eqref{f}, \eqref{14.27} and Lemma \ref{lem14.9}, together with the convexity bound, we obtain  
\begin{multline}\label{14.32}
\underset{s=\frac{3}{2}-s_1-s_2}{\Res}\  \mathfrak{F}_{\mathbf{1}}(s,s_1,s_2)-\underset{s=s_1+s_2-\frac{3}{2}}{\Res}\  \mathfrak{F}_{\mathbf{1}}(s,s_1,s_2)\\
\ll C_{\infty}(\pi')^{\frac{9}{4}+\varepsilon}\mathbf{C}_{\infty}^{\varepsilon}N_F(\mathfrak{L})^{-\frac{1}{2}+\varepsilon} 
[N_F(\mathfrak{N}),N_F(\mathfrak{M})]^{1+\varepsilon},
\end{multline}
and 
\begin{multline}\label{14.33}
\underset{s=s_1-s_2+\frac{1}{2}}{\Res}\  \mathfrak{F}_{\overline{\omega}\omega'}(s,s_1,s_2)-\underset{s=s_2-s_1-\frac{1}{2}}{\Res}\  \mathfrak{F}_{\omega\overline{\omega}'}(s,s_1,s_2)\\
\ll C_{\infty}(\pi')^{\frac{9}{4}+\varepsilon}\mathbf{C}_{\infty}^{\varepsilon}C_{\infty}(\omega\omega'^{-1})^{-\frac{1}{2}}N_F(\mathfrak{L})^{-\frac{1}{2}+\varepsilon}[N_F(\mathfrak{N}),N_F(\mathfrak{M})]^{1+\varepsilon}C(\omega\overline{\omega}'|\cdot|^{2it})^{1+\varepsilon}\\
C(\Ad\pi'\otimes\omega\overline{\omega}'|\cdot|^{2it})^{\frac{1}{2}+\varepsilon}\prod_{v\mid\infty}\mathbf{1}_{C_v(\omega\overline{\omega}'|\cdot|^{2it})\ll C_v(\pi')^{1+\varepsilon}}\prod_{v\mid\mathfrak{L}\mathfrak{M}}\mathbf{1}_{2r_{\omega_v\overline{\omega}_v'}\leq r_v'+l_v'}.
\end{multline}

Utilizing the constraints $C_v(\omega\overline{\omega}'|\cdot|^{2it})\ll C_v(\pi')^{1+\varepsilon}$ at $v\mid\infty$ and $2r_{\omega_v\overline{\omega}_v'}\leq r_v'+l_v'$ at $v<\infty$ we derive that 
\begin{align*}
C(\omega\overline{\omega}'|\cdot|^{2it})^{1+\varepsilon}
\ll C_{\infty}(\pi')^{1+\varepsilon}C_{\fin}(\omega\overline{\omega}')^{1+\varepsilon}\ll C_{\infty}(\pi')^{1+\varepsilon}C_{\fin}(\pi')^{\frac{1}{2}},
\end{align*}
and 
\begin{align*}
C(\Ad\pi'\otimes\omega\overline{\omega}'|\cdot|^{2it})\ll C(\Ad\pi')^{1+\varepsilon}C_{\infty}(\pi')^{1+\varepsilon}\ll C_{\infty}(\pi')^{1+\varepsilon}C(\pi')^{2+\varepsilon}.
\end{align*}

Substituting these into \eqref{14.33} yields 
\begin{multline}\label{14.34}
\underset{s=s_1-s_2+\frac{1}{2}}{\Res}\  \mathfrak{F}_{\overline{\omega}\omega'}(s,s_1,s_2)-\underset{s=s_2-s_1-\frac{1}{2}}{\Res}\  \mathfrak{F}_{\omega\overline{\omega}'}(s,s_1,s_2)\\
\ll C_{\infty}(\pi')^{\frac{19}{4}+\varepsilon}\mathbf{C}_{\infty}^{\varepsilon}C_{\infty}(\omega\omega'^{-1})^{-\frac{1}{2}}N_F(\mathfrak{L})^{-\frac{1}{2}+\varepsilon}N_F(\mathfrak{M})^{\frac{3}{2}+\varepsilon}[N_F(\mathfrak{N}),N_F(\mathfrak{M})]^{1+\varepsilon}.
\end{multline}

It follows from \eqref{e14.32}, \eqref{14.32} and \eqref{14.34} that 
\begin{equation}\label{14.36}
\mathcal{M}_{\Res}^{\du,\heartsuit}(s_1,s_2;\mathfrak{L})\ll \frac{C_{\infty}(\pi')^{\frac{19}{4}+\varepsilon}\mathbf{C}_{\infty}^{\varepsilon}N_F(\mathfrak{M})^{\frac{3}{2}+\varepsilon}[N_F(\mathfrak{N}),N_F(\mathfrak{M})]^{1+\varepsilon}}{N_F(\mathfrak{L})^{\frac{1}{2}-\varepsilon}}.
\end{equation}

Therefore, Proposition \ref{prop14.8} follows from the definitions \eqref{equ4.13} and \eqref{eq4.5}, and the above estimate \eqref{14.36}.

\section{Proof of Theorems \ref{A} and \ref{thmD}}\label{sec15}
\subsection{Proof of Theorem \ref{A}}
Let $s_1=1/2+it$ and $s_2=1/2-it$. By Theorem \ref{thm4.4} we have
\begin{multline}\label{15.1}
\mathcal{M}_{\cusp}^{\tw}(s_1,s_2;\boldsymbol{\alpha},\boldsymbol{\ell})+\mathcal{M}_{\Eis}^{\tw}(s_1,s_2;\boldsymbol{\alpha},\boldsymbol{\ell})=\mathcal{M}_{\const}^{\du,\heartsuit,\Reg}(s_1,s_2;\boldsymbol{\alpha},\boldsymbol{\ell})\\
+\mathcal{M}_{\cusp}^{\du}(s_1,s_2;\boldsymbol{\alpha},\boldsymbol{\ell})+\mathcal{M}_{\Eis}^{\du,\heartsuit}(s_1,s_2;\boldsymbol{\alpha},\boldsymbol{\ell})+\mathcal{M}_{\Res}^{\du,\heartsuit,\Reg}(s_1,s_2;\boldsymbol{\alpha},\boldsymbol{\ell}).
\end{multline}

Let $\boldsymbol{\ell}=(\ell_v)_{v\in \mathcal{P}}$, where $\ell_v\in\{1,2\}$ satisfies $|\lambda_{\pi}(\mathfrak{p}_v^{\ell_v})|\gg 1$. The existence of such an $\ell_v$ follows from the Hecke relation. Let $\boldsymbol{\alpha}=(\alpha_v)_{v\in \mathcal{P}}$, where $\alpha_v=\lambda_{\pi}(\mathfrak{p}_v^{\ell_v})^{-1}|\lambda_{\pi}(\mathfrak{p}_v^{\ell_v})|$. As a consequence of \eqref{equ4.3}, we obtain 
\begin{equation}\label{15.2}
\big|\sum_{v\in \mathcal{P}}\alpha_v\lambda_{\pi}(\mathfrak{p}_v^{\ell_v})\big|\gg L^{1-10\varepsilon}.
\end{equation} 

Therefore, it follows from Proposition \ref{prop5.3} that  
\begin{equation}\label{15.3}
\mathcal{M}_{\cusp}^{\tw}(s_1,s_2;\boldsymbol{\alpha},\boldsymbol{\ell})+\mathcal{M}_{\Eis}^{\tw}(s_1,s_2;\boldsymbol{\alpha},\boldsymbol{\ell})\gg \frac{|L(1/2+it,\overline{\pi}\times\pi')|^2L^{2-\varepsilon}}{\textbf{C}_{\infty}^{1+\varepsilon}C_{\fin}(\pi,\pi')^{\varepsilon}},
\end{equation}
where $\textbf{C}_{\infty}:=\textbf{C}_{\infty}(\pi,t,\pi')$ is defined in \textsection\ref{sec4.4}. 

By Proposition \ref{prop6.1} we obtain 
\begin{equation}\label{15.4}
\mathcal{M}_{\const}^{\du,\heartsuit,\Reg}(s_1,s_2;\boldsymbol{\alpha},\boldsymbol{\ell})\ll 
\textbf{C}_{\infty}^{\varepsilon}C_{\fin}(\pi,\pi')L^{1+\varepsilon}.
\end{equation}

Substituting \eqref{15.3} and \eqref{15.4}, and Propositions \ref{prop14.7} and \ref{prop14.8} into \eqref{15.1} yields
\begin{multline}\label{15.5}
|L(1/2+it,\overline{\pi}\times\pi')|^2\ll \frac{C_{\infty}(\pi')^{\frac{19}{4}+\varepsilon}\mathbf{C}_{\infty}^{1+\varepsilon}N_F(\mathfrak{M})^{\frac{3}{2}+\varepsilon}[N_F(\mathfrak{N}),N_F(\mathfrak{M})]^{1+\varepsilon}}{L^{1+\varepsilon}}\\
+\mathbf{C}_{\infty}^{\frac{1}{2}+\vartheta+\varepsilon}C_{\infty}(\omega\overline{\omega}')^{-\vartheta}
[N_F(\mathfrak{N}),N_F(\mathfrak{M})]^{\frac{1}{2}+\varepsilon}N_F(\mathfrak{M})C_{\infty}(\pi')^{\frac{43}{8}+\varepsilon}C_{\fin}(\Ad\pi')^{\frac{1}{2}}\\
\Big[L^{4+\varepsilon}N_F(\mathfrak{M})^{\frac{1}{2}}C(\omega\overline{\omega}'|\cdot|^{-2i t})^{\frac{1}{4}}C_{\infty}(\pi')^{\frac{1}{8}}
+L^{3+\varepsilon}N_F(\mathfrak{M})^{\frac{1}{4}}C(\omega\overline{\omega}'|\cdot|^{-2i t})^{\frac{3}{8}}\Big].
\end{multline}

We then proceed to optimize $L$ in \eqref{15.5} as follows. 
\begin{itemize}
\item Assuming the following constraint:
\begin{equation}\label{15.6}
\mathbf{C}_{\infty}^{\frac{1}{2}-\vartheta}[N_F(\mathfrak{N}),N_F(\mathfrak{M})]^{\frac{1}{2}}\leq 
N_F(\mathfrak{M})^{-\frac{1}{4}}C_{\infty}(\pi')^{\frac{1}{8}}
C_{\infty}(\omega\overline{\omega}')^{-\vartheta}C(\omega\overline{\omega}'|\cdot|^{-2i t})^{\frac{7}{8}}.
\end{equation}

Let $L_0$ be the paramerer defined by 
\begin{align*}
\mathbf{C}_{\infty}^{\frac{1}{2}-\vartheta+\varepsilon}C_{\infty}(\omega\overline{\omega}')^{\vartheta}[N_F(\mathfrak{N}),N_F(\mathfrak{M})]^{\frac{1}{2}+\varepsilon}=
N_F(\mathfrak{M})^{\frac{3}{4}}C_{\infty}(\pi')^{\frac{5}{8}+\varepsilon}
L_0^4C(\omega\overline{\omega}'|\cdot|^{-2i t})^{\frac{3}{8}}.
\end{align*}

It follows from \eqref{15.6} that 
\begin{equation}\label{15.7}
L_0^{1-\varepsilon}N_F(\mathfrak{M})^{\frac{1}{4}-\varepsilon}C_{\infty}(\pi')^{\frac{1}{8}-\varepsilon}
\ll C(\omega\overline{\omega}'|\cdot|^{-2i t})^{\frac{1}{8}+\varepsilon}.	
\end{equation}

Let $L=\max\{L_0, 100\}$.  Substituting \eqref{15.7} into \eqref{15.5} leads to 
\begin{multline}\label{15.8}
|L(1/2+it,\overline{\pi}\times\pi')|^8\ll \\
\frac{C_{\infty}(\pi')^{20+\varepsilon}\mathbf{C}_{\infty}^{\frac{7}{2}+\vartheta+\varepsilon}N_F(\mathfrak{M})^{7+\varepsilon}[N_F(\mathfrak{N}),N_F(\mathfrak{M})]^{\frac{7}{2}+\varepsilon}
C(\omega\overline{\omega}'|\cdot|^{-2i t})^{\frac{3}{8}}}{C_{\infty}(\omega\overline{\omega}')^{\vartheta}}.
\end{multline}

Notice that $C(\omega\overline{\omega}'|\cdot|^{-2i t})\ll \mathbf{C}_{\infty}[N_F(\mathfrak{N}),N_F(\mathfrak{M})]$. We thus derive from the estimate \eqref{15.8} that 
\begin{equation}\label{e15.9}
L(1/2+it,\overline{\pi}\times\pi')\ll 
C_{\infty}(\pi')^{\frac{5}{2}+\varepsilon}N_F(\mathfrak{M})^{\frac{21}{16}+\varepsilon}C(\pi\otimes |\cdot|^{it})^{\frac{1}{2}-\frac{1}{64}+\varepsilon}.
\end{equation}

\item Assuming the oppsite constraint of \eqref{15.6}: 
\begin{equation}\label{15.9}
\mathbf{C}_{\infty}^{\frac{1}{2}-\vartheta}[N_F(\mathfrak{N}),N_F(\mathfrak{M})]^{\frac{1}{2}}>
\frac{N_F(\mathfrak{M})^{-\frac{1}{4}}C_{\infty}(\pi')^{\frac{1}{8}}
C(\omega\overline{\omega}'|\cdot|^{-2i t})^{\frac{7}{8}}}{C_{\infty}(\omega\overline{\omega}')^{\vartheta}}.
\end{equation}

Let $L_0$ be the paramerer defined by 
\begin{multline*}
\mathbf{C}_{\infty}^{1+\varepsilon}[N_F(\mathfrak{N}),N_F(\mathfrak{M})]^{\frac{1}{2}+\varepsilon}\\
=\mathbf{C}_{\infty}^{\frac{1}{2}+\vartheta+\varepsilon}C_{\infty}(\omega\overline{\omega}')^{-\vartheta}
C_{\infty}(\pi')^{\frac{3}{4}+\varepsilon}N_F(\mathfrak{M})
C(\omega\overline{\omega}'|\cdot|^{-2i t})^{\frac{1}{4}}L_0^5.
\end{multline*}

It follows from \eqref{15.9} that 
\begin{equation}\label{15.10}
L_0^{1+\varepsilon}N_F(\mathfrak{M})^{\frac{1}{4}+\varepsilon}C_{\infty}(\pi')^{\frac{1}{8}+\varepsilon}
\gg C(\omega\overline{\omega}'|\cdot|^{-2i t})^{\frac{1}{8}-\varepsilon}.	
\end{equation}

Let $L=\max\{L_0, 100\}$.  Substituting \eqref{15.10} into \eqref{15.5} leads to
\begin{multline}\label{15.11}
|L(1/2+it,\overline{\pi}\times\pi')|^{10}\ll\\
\frac{C_{\infty}(\pi')^{25+\varepsilon}\mathbf{C}_{\infty}^{\frac{9}{2}+\vartheta+\varepsilon}N_F(\mathfrak{M})^{\frac{17}{2}+\varepsilon}[N_F(\mathfrak{N}),N_F(\mathfrak{M})]^{\frac{9}{2}+\varepsilon}
C(\omega\overline{\omega}'|\cdot|^{-2i t})^{\frac{1}{4}}}{C_{\infty}(\omega\overline{\omega}')^{\vartheta}}.
\end{multline} 

Similar to \eqref{e15.9}, the estimate \eqref{15.11} implies 
\begin{equation}\label{15.13}
L(1/2+it,\overline{\pi}\times\pi')\ll\\
C_{\infty}(\pi')^{\frac{5}{2}+\varepsilon}N_F(\mathfrak{M})^{\frac{13}{10}+\varepsilon}
C(\pi\otimes |\cdot|^{it})^{\frac{1}{2}-\frac{1}{40}+\varepsilon}.
\end{equation}
\end{itemize}

Therefore, Theorem \ref{A} follows from \eqref{e15.9} and \eqref{15.13}.

\subsection{Proof of Theorem \ref{thmD}}\label{sec15.2}
\subsubsection{Nice Rankin--Selberg Test Vectors}
Let $v \mid \infty$ and let $\pi_v$ be a unitary, irreducible, generic 
representation of $\mathrm{GL}_2(F_v)$. We say that a vector $W_v$ in the 
Whittaker model of $\pi_v$ is \emph{nice} if, for all $s$ with 
$\lvert \Re(s)-\tfrac12 \rvert < \tfrac{1}{10}$,
\begin{equation}\label{e15.14}
\int_{N(F_v)\backslash \overline{G}(F_v)} |W_v(x_v)|^2|\det x_v|_v^sd^{\times}x_v
= 
\frac{L_v(s,\pi_v \times \widetilde{\pi}_v)}{L_v(1,\pi_v \times \widetilde{\pi}_v)}.
\end{equation}

When $F_v \simeq \mathbb{R}$, such a \emph{nice} vector $W_v$ is a scalar multiple of a lowest--weight vector. When $F_v \simeq \mathbb{C}$, one may construct 
such a $W_v$ using \cite[Theorem~6.1]{Miy18} (after assuming by symmetry that 
$d_2 \ge d_1$, and then taking $m = d_2 - d_1$ and $m' = 0$ therein).

\begin{defn}\label{defn1.7}
Let $\pi = \otimes_v' \pi_v$ be a unitary cuspidal automorphic representation 
of $\mathrm{GL}_2/F$.  
A vector $\phi \in \pi$ is called a \emph{nice} cusp form if its associated 
Whittaker vector $W = \otimes_v' W_v$ satisfies:
\begin{itemize}
    \item for $v < \infty$, $W_v$ is a local newvector with $\langle W_v, W_v\rangle=1$;
    \item for $v \mid \infty$, $W_v$ is a \emph{nice} Whittaker function with $\langle W_v, W_v\rangle=1$.
\end{itemize}
\end{defn}

\subsubsection{Estimate of Trilinear Forms}
Let $\sigma=\otimes_v'\sigma_v$ be a unitary automorphic representation of $\mathrm{PGL}_2/F$. Suppose $\sigma$ is everywhere unramified, and $\varphi\in\sigma$ is a spherical vector, i.e., $\varphi$ is right-$K$-invariant. 

\begin{prop}\label{prop15.2}
Let $\pi$ be a unitary cuspidal automorphic  representation of $\mathrm{GL}_2/F$ and $\phi\in\pi$ be a nice cusp form. Then 
\begin{equation}\label{eq15.15}
\frac{|\langle \varphi, |\phi|^2\rangle|^2}{\langle \varphi,\varphi\rangle}\ll C(\sigma)^{10}C(\pi)^{-1+\varepsilon}L(1/2,\sigma\times \pi\times\widetilde{\pi}),
\end{equation}	
where the implied constant depends on $F$ and $\varepsilon$. 
\end{prop}
\begin{proof}
Let $W_{\varphi}=\otimes_v'W_{\varphi,v}$ and $W=\otimes_v'W_{v}$ be the Whittaker function associated with $\varphi$ and $\phi$, respectively. For each place $v\leq\infty$, define 
\begin{align*}
|\mathcal{P}_v(W_{\varphi,v},W_{v},\overline{W_{v}})|^2:=\int_{\overline{G}(F_v)}\frac{\langle\sigma_v(g_v)W_{\varphi,v},W_{\varphi,v}\rangle|\langle\pi_v(g_v)W_v,W_v\rangle|^2}{\langle W_{\varphi,v},W_{\varphi,v}\rangle\langle W_{v},W_{v}\rangle^2} dg_v.
\end{align*}

By definition $\langle \phi,\phi\rangle=1$. As a consequence of  Watson-Ichino formula, we have
\begin{equation}\label{15.14}
\frac{|\langle \varphi, |\phi|^2\rangle|^2}{\langle \varphi,\varphi\rangle}=L(1/2,\sigma\times\pi\times\widetilde{\pi})\prod_{v\mid\infty}|\mathcal{P}_v(W_{\varphi,v},W_{v},\overline{W_{v}})|^2\prod_{v<\infty}\mathcal{I}_v^{\sharp}(\sigma),
\end{equation}
where, for $v<\infty$, 
\begin{align*}
\mathcal{I}_v^{\sharp}(\sigma):=L_v(1/2,\sigma_v\times\pi_v\times\widetilde{\pi}_v)^{-1}|\mathcal{P}_v(W_{\varphi,v},W_{v},\overline{W_{v}})|^2.
\end{align*}

By Lemmas \ref{lemma7.4}--\ref{lemma7.7} and Proposition \ref{prop10.1} we have 
\begin{equation}\label{15.15}
\prod_{v<\infty}\mathcal{I}_v^{\sharp}(\sigma)\ll C_{\fin}(\pi)^{-1+\varepsilon}.	
\end{equation}

Let $v\mid\infty$. Since $\sigma_v$ is unramified, we may assume $\sigma_v=|\cdot|_v^{\nu}\boxplus |\cdot|_v^{-\nu}$, where $|\Re(\nu)|\leq 7/64$. Therefore, it follows from \cite[Lemma 3.4.2]{MV10} that 
\begin{align*}
|\mathcal{P}_v(W_{\varphi,v},W_{v},\overline{W_{v}})|^2=\frac{\zeta_v(1)}{\langle W_{\varphi,v},W_{\varphi,v}\rangle}\prod_{\epsilon\in\{\pm 1\}}	\int_{X_v} |W_v(x_v)|^2|\det x_v|_v^{\frac{1}{2}+\epsilon\nu}d^{\times}x_v,
\end{align*}
where $X_v:=N(F_v)\backslash \overline{G}(F_v)$. Substituting  \eqref{e15.14} into the above integral yields 
\begin{equation}\label{15.18}
|\mathcal{P}_v(W_{\varphi,v},W_{v},\overline{W_{v}})|^2=\zeta_v(1)\cdot  \frac{L_v(1/2,\sigma_v\times \pi_v \times \widetilde{\pi}_v)}{L_v(1,\pi_v \times \widetilde{\pi}_v)^2}.
\end{equation}

Therefore, \eqref{eq15.15} follows from \eqref{15.14}, \eqref{15.15}, \eqref{15.18}, and the Stirling's formula. 
\end{proof}

\begin{cor}
Suppose $\pi$ is a unitary dihedral cuspidal automorphic representation of $\mathrm{GL}_2/F$ and $\phi\in\pi$ be a nice cusp form. Then 
\begin{equation}\label{15.19}
\frac{|\langle \varphi, |\phi|^2\rangle|}{\sqrt{\langle \varphi,\varphi\rangle}}\ll C(\sigma)^{10}C(\pi)^{-\frac{1}{128}+\varepsilon}C_{\fin}(\pi)^{-\frac{1}{16}},
\end{equation}	
where the implied constant depends on $F$ and $\varepsilon$.
\end{cor}
\begin{proof}
By definition, there exists a quadratic extension $E/F$ and a unitary Hecke character $\chi$ of $E^{\times}\backslash\mathbb{A}_E^{\times}$, such that  $\pi=\mathrm{AI}_{E/F}(\chi)$ is the automorphic induction of $\chi$.

Let $\iota$ be the nontrivial element in the Galois group $\mathrm{Gal}(E/F)$. Define 
\begin{align*}
\chi^{\iota}(t):=\chi(\iota(t)),\ \ t\in \mathbb{A}_E^{\times}. 
\end{align*}

Let $\eta_{E/F}$ be the quadratic character associated with $E/F$. Then 
\begin{equation}\label{15.20}
\pi\times\widetilde{\pi}=\mathbf{1}\boxplus \eta_{E/F}\boxplus \mathrm{AI}_{E/F}(\chi/\chi^{\iota}). 
\end{equation} 

Therefore, it follows from \eqref{15.20} that 
\begin{equation}\label{15.21}
L(1/2,\sigma\times \pi\times\widetilde{\pi})=L(1/2,\sigma)L(1/2,\sigma\times \eta_{E/F})L(1/2,\sigma\times \mathrm{AI}_{E/F}(\chi/\chi^{\iota})).
\end{equation}

Notice the relation between analytic conductors 
\begin{equation}\label{15.22}
C(\mathrm{AI}_{E/F}(\chi/\chi^{\iota}))\ll C(\mathrm{AI}_{E/F}(\chi))=C(\pi).
\end{equation}

Substituting \eqref{15.22}, Theorem \ref{thm14.3}, and Theorem \ref{A}, together with its non-cuspidal analogue \cite{Yan26b}, into \eqref{15.21}, we obtain 
\begin{equation}\label{15.23}
L(1/2,\sigma\times \pi\times\widetilde{\pi})\ll C(\sigma)^{10}	C(\eta_{E/F})^{\frac{3}{8}+\varepsilon}C(\pi)^{\frac{1}{2}-\frac{1}{64}+\varepsilon}.
\end{equation} 

Notice that $C(\eta_{E/F})\leq C_{\fin}(\pi)$. Therefore, \eqref{15.19} follows from \eqref{15.23}. 
\end{proof}

\subsubsection{Proof of Theorem \ref{thmD}}
Let $f \in C_c(X)$. We have the spectral decomposition
\begin{equation}\label{15.24}
f(x)=\int_X f(g)dg
+\int_{\varphi}
\frac{\langle f,\varphi\rangle\varphi(x)}{\langle \varphi,\varphi\rangle}
d\mu_{\varphi},
\end{equation}
where the integral over $\varphi$ runs over unitary generic automorphic forms on $X$.

Let $\sigma$ be a unitary automorphic representation of $\mathrm{PGL}_2/F$ that is unramified at all places, and let $\varphi \in \sigma$. By repeated integration by parts, we obtain
\begin{equation}\label{15.25}
\langle f,\varphi\rangle
=
C(\sigma)^{-100}\cdot 
\langle \mathcal{D}f,\varphi\rangle,
\end{equation}
for a suitable differential operator $\mathcal{D}$ depending only on $X$.

As a consequence of \eqref{15.24} and \eqref{15.25}, we deduce
\begin{multline}\label{15.26}
\bigg|
\int_{Z(\mathbb{A}_F)G(F)\backslash G(\mathbb{A}_F)}
    f(x)|\phi(x)|^2dx-\int_X f(x)dx\bigg|\\
\ll C(\pi)^{-\frac{1}{128}+\varepsilon}
C_{\fin}(\pi)^{-\frac{1}{16}}\int_{\varphi}
\frac{| 
\langle \mathcal{D}f,\varphi\rangle|}{\sqrt{\langle \varphi,\varphi\rangle}}\cdot 
C(\sigma)^{-50}
d\mu_{\varphi}.
\end{multline}

By Cauchy--Schwarz inequality, 
\begin{equation}\label{15.27}
\int_{\varphi}
\frac{| 
\langle \mathcal{D}f,\varphi\rangle|}{\sqrt{\langle \varphi,\varphi\rangle}}\cdot 
C(\sigma)^{-50}
d\mu_{\varphi}\ll \bigg[\int_{\varphi}
\frac{| 
\langle \mathcal{D}f,\varphi\rangle|^2}{\langle \varphi,\varphi\rangle}
d\mu_{\varphi}\bigg]^{\frac{1}{2}}\leq \langle \mathcal{D}f,\mathcal{D}f\rangle^{\frac{1}{2}}.
\end{equation}

Therefore, \eqref{1.3} follows from \eqref{15.26} and \eqref{15.27}.

\bibliographystyle{alpha}
\bibliography{HHY}

\end{document}